\documentclass[twoside,a4paper,12pt,centertags]{amsart}
\usepackage[T1]{fontenc}
\usepackage[utf8]{inputenc}
\usepackage{subfiles}
\usepackage{multicol}
\makeatletter
\usepackage[section]{placeins}
\AtBeginDocument{%
	\expandafter\renewcommand\expandafter\subsection\expandafter{%
		\expandafter\@fb@secFB\subsection
	}%
}
\makeatother
\makeatletter
\@namedef{subjclassname@2020}{%
	\textup{2020} Mathematics Subject Classification}
\makeatother
\usepackage{amsfonts,amssymb,amsmath,amsthm}
\usepackage{mathrsfs}
\usepackage{nicefrac}
\usepackage{latexsym}
\usepackage{mathtools}
\usepackage[noadjust]{cite}
\usepackage{paralist}
\usepackage{longtable}
\usepackage{float}
\usepackage{pdfsync}
\usepackage{tikz, tikz-cd}
\usetikzlibrary{patterns}
\usepackage{hhline}
\usepackage{import}
\usepackage{enumitem}

\newtheorem{thm}{Theorem}[section]
\newtheorem{lem}[thm]{Lemma}
\newtheorem{prop}[thm]{Proposition}
\newtheorem{cor}[thm]{Corollary}

\theoremstyle{definition}
\newtheorem{defn}[thm]{Definition}
\newtheorem{rem}[thm]{Remark}

\numberwithin{equation}{section}
\usepackage[plainpages=false,pdfpagelabels,backref=page,
citecolor=red]{hyperref}
\usepackage{xcolor}
\hypersetup{
 colorlinks,
 linkcolor={blue!80!black},
 citecolor={red!80!black},
 urlcolor={green!40!black}
} 

\newcommand{\IB}{\mathbb{B}}
\newcommand{\IC}{\mathbb{C}}

\newcommand{\IH}{\mathbb{H}}
\newcommand{\IN}{\mathbb{N}}
\newcommand{\IP}{\mathbb{P}}
\newcommand{\IQ}{\mathbb{Q}}
\newcommand{\R}{\mathbb{R}}
\newcommand{\IS}{\mathbb{S}}
\newcommand{\IX}{\mathbb{X}}
\newcommand{\IZ}{\mathbb{Z}}
\newcommand{\cD}{\mathcal{D}} 
\newcommand{\cF}{\mathcal{F}} 
\newcommand{\cH}{\mathcal{H}}
\newcommand{\cI}{\mathcal{I}}
\newcommand{\cJ}{\mathcal{J}}
\newcommand{\cL}{\mathcal{L}}
\newcommand{\cM}{\mathcal{M}}
\newcommand{\cN}{\mathcal{N}}
\newcommand{\cP}{\mathcal{P}}
\newcommand{\cS}{\mathcal{S}} 

\newcommand{\cZ}{\mathcal{Z}}
\newcommand{\sol}{\operatorname{sol}}
\newcommand{\loc}{\operatorname{loc}}
\newcommand{\glob}{\operatorname{glob}}
\renewcommand{\L}{\operatorname{L}} 
\newcommand{\Lloc}{\L_{\operatorname{loc}}} 
\newcommand{\C}{\operatorname{C}} 
\renewcommand{\H}{\operatorname{H}} 
\newcommand{\W}{\operatorname{W}}
\newcommand{\T}{\operatorname{T}}
\newcommand{\B}{\operatorname{B}}
\newcommand{\F}{\operatorname{F}}
\DeclareRobustCommand{\Hdot}{\dot{\H}\protect{\vphantom{H}}} 
\DeclareRobustCommand{\Wdot}{\dot{\W}\protect{\vphantom{W}}} 
\DeclareRobustCommand{\Fdot}{\dot{\F}\protect{\vphantom{F}}} 
\DeclareRobustCommand{\Bdot}{\dot{\B}\protect{\vphantom{B}}} 
\DeclareRobustCommand{\Xdot}{\dot{\X}\protect{\vphantom{X}}} 
\DeclareRobustCommand{\Lamdot}{\dot{\Lambda}\protect{\vphantom{\Lambda}}} 
\DeclareRobustCommand{\BMOdot}{\dot{\mathrm{BMO}}\protect{\vphantom{\mathrm{BMO}}}} 
\newcommand{\psiH}{\psi  \IH} 
\newcommand{\psiB}{\psi \IB} 
\newcommand{\psiX}{\psi  \IX} 
\newcommand{\varphiX}{\varphi  \IX} 
\newcommand{\Wloc}{\W_{\operatorname{loc}}} 
\newcommand{\X}{\operatorname{X}} 
\newcommand{\Y}{\operatorname{Y}} 
\newcommand{\Z}{\operatorname{Z}} 
\let\SS\S	
\renewcommand{\S}{\mathrm{S}} 
\newcommand{\BMO}{\mathrm{BMO}}
\newcommand{\I}{\mathrm{\, I \, }}
\newcommand{\II}{\mathrm{\, II \,}}
\newcommand{\III}{\mathrm{\, III \,}}
\newcommand{\reu}{{\mathbb{R}^{1+n}_+}}
\newcommand{\ree}{{\mathbb{R}^{1+n}}}
\newcommand{\dual}[2]{\langle #1,#2 \rangle}

\newcommand{\tdd}[2]{\tfrac{\partial #1}{\partial #2}}
\newcommand{\wt}{\widetilde}

\newcommand{\ta}{{\scriptscriptstyle \parallel}}
\newcommand{\no}{{\scriptscriptstyle\perp}}
\newcommand{\pd}{\partial}
\newcommand{\tL}{\widetilde{L}}
\newcommand{\tM}{\widetilde{M}}
\newcommand{\gradx}{\nabla_x}
\renewcommand{\div}{\operatorname{div}}
\newcommand{\curl}{\operatorname{curl}}
\newcommand{\dnuA}{\partial_{\nu_A}} 
\newcommand{\NT}{\widetilde{N}_*} 
\newcommand{\NTq}{\widetilde{N}_{*,q}} 
\newcommand{\NTr}{\widetilde{N}_{*,r}} 
\newcommand{\NTone}{\widetilde{N}_{*,1}} 
\newcommand{\NTsharpalpha}{\widetilde{N}_{\sharp,\alpha}} 
\newcommand{\wtD}{\smash{\wt{D}}\vphantom{D^2}}
\newcommand{\e}{\mathrm{e}} 

\renewcommand{\i}{\mathrm{i}} 
\renewcommand{\d}{\mathrm{d}} 
\newcommand{\eps}{\varepsilon} 
\renewcommand\Re{\operatorname{Re}}
\renewcommand\Im{\operatorname{Im}}
\def\angle#1#2{\langle #1,#2 \rangle} 
\newcommand{\Lop}{\mathcal{L}} 
\newcommand{\Le}{\mathcal{L}}
\newcommand{\cl}[1]{\overline{#1}} 
\DeclareMathOperator{\supp}{supp} 
\DeclareMathOperator{\ran}{\mathsf{R}} 
\DeclareMathOperator{\nul}{\mathsf{N}} 
\DeclareMathOperator{\dom}{\mathsf{D}} 
\DeclareMathOperator\Max{\mathcal{M}} 

\newcommand{\sgn}{\operatorname{sgn}}
\newcommand{\ind}{\mathbf{1}}
\newcommand{\dist}{\mathrm{d}}
\newcommand{\bd}{\partial}
\newcommand{\eo}{\Lambda}
\def\Xint#1{\mathchoice
{\XXint\displaystyle\textstyle{#1}}%
{\XXint\textstyle\scriptstyle{#1}}%
{\XXint\scriptstyle\scriptscriptstyle{#1}}%
{\XXint\scriptscriptstyle%
\scriptscriptstyle{#1}}%
\!\int}
\def\XXint#1#2#3{{\setbox0=\hbox{$#1{#2#3}{%
\int}$ }
\vcenter{\hbox{$#2#3$ }}\kern-.6\wd0}}
\def\barint{\,\Xint -} 
\def\bariint{\barint_{} \kern-.4em \barint}
\def\bariiint{\bariint_{} \kern-.4em \barint}
\renewcommand{\iint}{\int_{}\kern-.34em \int} 
\renewcommand{\iiint}{\iint_{}\kern-.34em \int} 
\title[Boundary value problems and Hardy spaces]{Boundary value problems and Hardy spaces for elliptic systems with block structure}
\author{Pascal Auscher}
\address{Universit\'e Paris-Saclay, CNRS, Laboratoire de Math\'{e}matiques d'Orsay, 91405 Orsay, France}
\email{pascal.auscher@universite-paris-saclay.fr}
\author{Moritz Egert}
\address{TU Darmstadt, Fachbereich Mathematik, Schlossgartenstr.\ 7, 64289 Darmstadt, Germany}
\address{Universit\'e Paris-Saclay, CNRS, Laboratoire de Math\'{e}matiques d'Orsay, 91405 Orsay, France}
\email{egert@mathematik.tu-darmstadt.de}
\thanks{This material is based upon work that got started under NSF Grant DMS-1440140 while Auscher was in residence at the MSRI in Berkeley, California, during the Spring 2017 semester. Egert also thanks this institute for hospitality. The authors were supported by the ANR projects HAB ANR-12-BS01-0013  and RAGE ANR-18-CE40-0012. The authors would like to thank Tim B\"ohnlein for a very careful reading of earlier versions of their manuscript and Hans Triebel and Winfried Sickel for helping them with the literature on interpolation spaces. This is a preprint of the following work: P.~Auscher and M.~Egert, Boundary value problems and Hardy spaces for elliptic systems with block structure, 2023, Birkh\"auser reproduced with permission of Birkh\"auser. The final authenticated version is available online at: \url{https://doi.org/10.1007/978-3-031-29973-5}}
\date{August 3, 2022}
\dedicatory{}
\allowdisplaybreaks
\makeindex
\begin{document}
\begin{abstract}
For elliptic systems with block structure in the upper half-space and $t$-independent coefficients, we settle the study of boundary value problems by proving compatible well-posedness of Dirichlet, regularity and Neumann problems in optimal ranges of exponents. Prior to this work, only the two-dimensional situation was fully understood. In higher dimensions, partial results for existence in smaller ranges of exponents and for a subclass of such systems had been established. 
The presented uniqueness results are completely new. We also elucidate optimal ranges for problems with fractional regularity data.

The first part of the monograph, which can be read independently, provides optimal ranges of exponents for functional calculus and adapted Hardy spaces for the associated boundary operator.

Methods use and improve, with new results, all the machinery developed over the last two decades to study such problems:  the Kato square root estimates and Riesz transforms, Hardy spaces associated to operators,  off-diagonal estimates, non-tangen\-tial estimates and square functions, and abstract layer potentials to replace fundamental solutions in the absence of local regularity of solutions.

This mostly self-contained monograph provides a comprehensive overview on the field and unifies many earlier results that have been obtained by a variety of methods.  
\end{abstract}

\maketitle

\newpage

\tableofcontents

\newpage

\section{Introduction and main results}
\label{sec: Introduction}
\subsection{Objective of the monograph}
\label{intro:objective}

Consider the elliptic system of $m$ equations in $(1+n)$ dimensions\index{dimension!$n$ (boundary dimension)}, $n\ge 1$, given by
\begin{align*}  
\sum_{i,j=0}^n\sum_{\beta= 1}^m \pd_i\big( A_{i,j}^{\alpha, \beta}(x) \pd_j u^{\beta}(t,x)\big) =0 \quad (\alpha=1,\ldots, m, \  t>0, \ x\in \R^n),
\end{align*}
where $\pd_0 \coloneqq \tdd{}{t}$ and $\pd_i \coloneqq \tdd{}{x_{i}}$ if $i=1,\ldots,n$. Note that the coefficients do not depend on the normal variable $t>0$. Ellipticity will be described below, but when $m=1$, the uniformly elliptic equations are included. 

Boundary value problems for such systems have been extensively studied since the pioneering work of Dahlberg~\cite{Da} in the late 1970s. The upper half-space situation is prototypical for Lipschitz graph domains. The case of $t$-independent coefficients is already challenging and meaningful since $t$-dependent coefficients are usually treated via perturbation techniques.\footnote{The reader can refer to Kenig's excellent survey~\cite{K} for background on these topics. They lie beyond the scope of our monograph.}
As usual in the harmonic analysis treatment of elliptic  boundary  value problems, solutions are taken in the weak sense, interior estimates involve non-tangential maximal functions and/or conical square functions and convergence at the boundary is to be understood in an appropriate non-tangential sense. 

In this monograph, we consider the class of systems in \emph{block form}\index{block form!system in}, that is, when there are no mixed $\frac{\partial}{\partial t} \frac{\partial}{\partial {x_{i}}}$-derivatives. In short notation, the system can be written as 
\begin{equation}
\label{eq:block}
\partial_{t}(a \partial_{t}u) + \div_{x}(d\gradx  u)=0 
\end{equation}
where the matrix $A=(A_{i,j}^{\alpha, \beta}(x))$ above is block diagonal with diagonal (matrix) entries $a=a(x)$ and $d=d(x)$, hence the name. These systems enjoy the additional feature that one can always produce strong solutions using the Poisson semigroup $\e^{-t L^{1/2}}$ associated with the sectorial operator $L\coloneqq -a^{-1}\div_x d \nabla_x$ on the boundary.\footnote{We identify the boundary of the upper half-space with $\R^n$.} 
Existence and uniqueness of solutions to the boundary value problems are therefore inseparably tied to operator theoretic properties of $L$.

Our goal is to identify all spaces of boundary data of Hardy, Lebesgue and homogeneous H\"older-type, for which the Dirichlet and Neumann boundary value problems have weak solutions, and then prove uniqueness in these cases. Thus, we aim at proving well-posedness results for the largest possible class of boundary spaces. 

To this end, we unify and improve, with several new results along the way,  all the machinery developed over the last two decades to study such problems: the Kato square root estimates and Riesz transforms, Hardy spaces associated to operators, off-diagonal estimates, non-tangential estimates and square functions, abstract layer potentials replacing fundamental solutions in the absence of local regularity of solutions, \ldots

Prior to this work, only the two-dimensional situation was fully understood for the boundary value problems. In higher dimensions, partial results for existence in smaller ranges of exponents and for a subclass of such systems had been established. The uniqueness results are completely new. We essentially close this topic by obtaining well-posedness in ranges of boundary spaces likely to be sharp in all dimensions. 

For Dirichlet-type problems these ranges go beyond the semigroup theory for $\e^{-tL^{1/2}}$ on Lebesgue or Sobolev spaces. The global picture is that for the regularity problem, one can go one Sobolev exponent down from the semigroup range and for the Dirichlet problem, one can go one Sobolev exponent up. In particular, we exhibit for the first time the possibility of solving Dirichlet problems for H\"older and $\BMO$-data without relying on any sort of duality with an adjoint problem with data in a Hardy space. For the Neumann problem, we shall provide a missing link to the existing literature, so that well-posedness in the optimal range of boundary spaces follows from earlier results. This range is the one provided by the semigroup theory. 

Natural extensions of the results above are the Dirichlet and Neumann problems for data with fractional regularity between $0$ and $1$, for which we also provide well-posedness results. This concerns data in Besov and even Hardy--Sobolev spaces. We believe they are optimal in the formulation of the problem as well as in the ranges of spaces.

Most recent results in the field rely on one of two opposing strategies, sometimes referred to as \emph{second-} and \emph{first}-order approaches. None of these two approaches can be used `off-the-shelf' in order to cover the full range of results that we are aiming at here. Indeed, in the former, the Poisson semigroup $\e^{-tL^{1/2}}$ is usually treated by comparison with the heat semigroup $\e^{-t^2 L}$, which offers better decay properties.\footnote{References for these techniques are \cite{May, AHM, CMP1, CMP2}.} 
When $a \neq 1$, it may happen that $L$ is sectorial of angle larger than $\nicefrac{\pi}{2}$, and hence $-L$ does not generate a heat semigroup. This forces us to rely on resolvents $(1+t^2L)^{-1}$ instead, which offer sufficient off-diagonal decay but introduce new and partly unsuspected technicalities. In the first-order approach, the elliptic equation is rewritten as an equivalent first-order system of Cauchy--Riemann-type for the variable $F = [a\partial_t u, \nabla_x u]^\top$ called the \emph{conormal gradient}\index{conormal gradient}.\footnote{In this context the idea is pioneered in \cite{AAMc, AA1}.}
The approach is genuinely built on the use of resolvents of a first-order operator, but the range of admissible data spaces is limited since it treats the interior estimates for Dirichlet and Neumann problems simultaneously.

Most of our arguments are carried out at the second-order level, but whenever convenient, we employ first-order methods to give more efficient proofs and novel results, even when $a=1$. Readers, who are not familiar with the first-order approach, may find in this monograph a light introduction to some important features of the theory, while keeping technicalities at the absolute minimum. We also characterize all ranges of boundary spaces that have previously been obtained through first-order methods, using only the second-order operator $L$. We believe that this helps in rendering accessible the cornerstones of the first-order method to the broader audience that they deserve. At the same time, the block structure will reveal interesting new phenomena that could not be captured by the first-order method.
\subsection{The elliptic equation}
\label{intro:elliptic}
\noindent Consider again the elliptic equation \eqref{eq:block}. The value of $m$ (the number of equations)\index{dimension!$m$ (number of equations)} is irrelevant to everything that follows and the reader may assume $m=1$ when it comes to differential operators such as gradient and divergence.\footnote{Notation in the case $m>1$ looks exactly the same and is explained in Section~\ref{subsec: Notation}.}  We write \eqref{eq:block} as
\begin{align*}
\Le u \coloneqq - \div(A \nabla u) = - \partial_t (a \partial_t u) - \div_x d \nabla_x u = 0,
\end{align*}
where 
\begin{align*}
A = \begin{bmatrix} a & 0 \\ 0 & d \end{bmatrix}: \R^n \to \Lop(\IC^m \times \IC^{mn})
\end{align*}
is the  coefficient matrix of dimension $m(1+n)$ in block form. The equation is understood in the weak sense: By $\Le u = 0$ we mean that $u \in \Wloc^{1,2}(\reu; \IC^m)$ satisfies 
\begin{align*}
\iint_{\reu}  A \nabla u \cdot \cl{\nabla \phi} \, \d t \d x = 0 \quad (\phi \in \C_0^\infty(\reu; \IC^m)).
\end{align*}
We assume that $A$ is measurable and that there is a constant $\lambda \in (0,\infty)$ called \emph{ellipticity constant}\index{ellipticity!constant}, such that the following hold. First, $A$ is bounded from above:
\begin{align*}
\|A\|_\infty \leq \lambda^{-1}.
\end{align*}
Second, $A$ is bounded from below on the subspace $\cH$\index{H@$\cH$ (closure of $\ran(D)$)} of vector fields $f=[f_0,\ldots f_n]^\top$ in $\L^2(\R^n;(\IC^m)^{1+n})$ that satisfy the curl-free condition $\partial_j f_k = \partial_k f_j$ whenever $1 \leq j,k \leq n$:
\begin{align}
\label{eq: accretivity A}
\Re \langle A f, f \rangle \geq \lambda \|f\|_2^2 \quad (f \in \cH),
\end{align}
where the angular brackets denote the inner product on $\L^2$. Due to the block form, this lower bound can be written equivalently as two separate conditions\footnote{This follows since in the definition of $\cH$ the first component $f_0$ is arbitrary and the curl-free condition is equivalent to $[f_1,\ldots,f_n]^\top = \nabla_x h$ for some distribution $h$, see \cite[p.~59]{Schwartz-Dist}. Then use that $\C_0^\infty$ is dense in the homogeneous Sobolev space $\Wdot^{1,2}$, see \cite[Thm.~1]{Sohr-Specovius}.} : 
\emph{Strict ellipticity}\footnote{The term \emph{strict accretivity} is also common.}\index{ellipticity!strict} of $a$,
\begin{align}
\label{eq: accretivity a}
\Re \langle a(x) \xi, \xi \rangle \geq \lambda |\xi|^2 \quad (x \in \R^n, \, \xi \in \IC^m),
\end{align}
so that $a$ is also invertible in $\L^\infty(\R^n ;\IC^m)$, and the \emph{G\aa{}rding inequality}\index{G\aa{}rding inequality} for $d$,
\begin{align}
\label{eq: Garding d}
\Re \langle d \nabla_x v, \nabla_x v \rangle \geq \lambda \|\nabla_x v\|_2^2 \quad (v \in \C_0^\infty(\R^n; \IC^m)),
\end{align}
which in general is weaker than strict ellipticity even when $m=1$.\footnote{Take for instance $d=\begin{bmatrix} 42 & \i \\ -\i & 42 \end{bmatrix}$ and calculate.}

\subsection{The critical numbers}

\noindent We use Hardy and homogeneous Hardy--Sobolev spaces $\H^p$ and $\Hdot^{1,p}$ in the range $p\in (1_{*}, \infty)$ with the convention that for $p \in (1,\infty)$ they coincide with Lebesgue and homogeneous Sobolev spaces $\L^p$ and $\Wdot^{1,p}$, respectively. We denote by $p_*$ and $p^*$  the lower and upper Sobolev conjugates of $p$. In particular, $1_* \coloneqq \nicefrac{n}{(n+1)}$. 

We keep on denoting by 
\begin{align*}
L = -a^{-1} \div_x d \nabla_x
\end{align*}
the boundary operator associated with \eqref{eq:block}, defined as a sectorial operator in $\L^2$ with maximal domain in $\W^{1,2}$. 

The applications to boundary value problems require understanding the functional properties of the Poisson semigroup $(\e^{-tL^{1/2}})_{t>0}$,  which comes as the natural solution operator, on Hardy and Hardy--Sobolev spaces. The existence of the Poisson semigroup operator $\e^{-t L^{1/2}}$ is granted from the functional calculus for $L$ on $\L^2$.\footnote{This is a classical construction. We give the necessary background in Section~\ref{sec: L2}.}

Two intervals will rule our entire theory:
\begin{itemize}
	\item  $(p_{-}(L), p_{+}(L))$ is the maximal open set within $(1_*, \infty)$ for which the family $(a(1+t^2L)^{-1}a^{-1})_{t>0}$ is uniformly bounded on $\H^p$. \\[-8pt]
	\item  $(q_{-}(L),q_{+}(L))$ is the maximal open set within $(1_*, \infty)$ for which $(t\nabla_{x}(1+t^2L)^{-1}a^{-1})_{t>0}$ is uniformly bounded on $\H^p$. 
\end{itemize}
The endpoints $p_\pm(L), q_\pm(L)$ are called \emph{critical numbers}\index{critical numbers} associated with $L$.\footnote{The idea to use critical numbers for the sake of a flexible theory that applies to any given operator originates in \cite{A}. Therein, they have been defined for $a=1$ through $\L^p$-boundedness of the heat semigroup. We shall prove that when $a=1$ our intervals coincide with the ones of \cite{A} in the range $(1,\infty)$, see Section~\ref{sec: Poisson}.}
They have various characterizations proved throughout the monograph. For example, replacing $(1+t^2L)^{-1}$ by $\e^{-tL^{1/2}}$ leads to the same intervals, which shows that the critical numbers capture sharp uniform boundedness properties of the Poisson semigroup for $L$ in Hardy spaces.\footnote{This is proved in Section~\ref{sec: Poisson}.} We give a systematic study of these numbers, their inner relationship and their values depending on the dimension for the class of all $L$. In particular, we shall show that they are independent of $a$.\footnote{This is proved in Section \ref{sec: J and N}.} Of course, that does not mean that we can assume $a=1$ in general.

For now, all one needs to know is that the best conclusion for the critical numbers for the class of all $L$ is
\begin{align*}
	(p_{-}(L), p_{+}(L)) 
	&\supseteq 
	{\begin{cases}
			(\frac{1}{2}, \infty) & \text{if $n=1$} \\
			[1, \infty) & \text{if $n = 2$} \\
			[\frac{2n}{n+2}, \frac{2n}{n-2}] & \text{if $n \geq 3$}
	\end{cases}}
	\intertext{and}
	(q_{-}(L), q_{+}(L)) 
	&\supseteq
	{\begin{cases}
			(\frac{1}{2}, \infty) & \text{if $n=1$} \\
			[\tfrac{2n}{n+2}, 2] & \text{if $n \ge 2$} \\
	\end{cases}}
\end{align*}
and that in general $p_-(L) = q_-(L)$ and $p_+(L) \geq (q_+(L))^*$.
Including systematically exponents $p \in (1_*,1]$ is a novelty of our approach for both the functional properties of $L$ for their own sake\footnote{Section~\ref{sec: consequences identification} is about consequences for the functional calculus and Section~\ref{sec: Critical numbers and kernels} provides a connection to kernel estimates.} and the applications to boundary value problems.
\subsection{Square root problem and Hardy spaces}
\label{intro:square root and Hardy}

One may wonder how we determine the spaces of data for the boundary value problems. Typically, they should include Lebesgue spaces, Sobolev spaces in the range $p>1$ and also Hardy and Hardy--Sobolev spaces in the range $p\le 1$, as well as their intermediate fractional spaces.  Indeed, it is natural from the point of view of regularity theory to incorporate the possibility of having estimates for $p\leq 1$, as is the case for instance for equations with real coefficients.  
The limitation to $p>1_*$ can be understood from Sobolev embeddings and duality:  The best one can hope for in absence of smoothness of the coefficients is regularity theory in H\"older spaces of exponents less than $1$.

The whole theory is built from the case $p=2$.  For the regularity problem\footnote{More precisely, the problem $(R)_{2}^\Le$ defined in Section~\ref{intro:main}.}, it was Kenig\footnote{See \cite[Rem.~2.5.6]{K}.} who observed that the required interior estimates are linked to the Kato conjecture for $L$, that is, the homogeneous estimate
\begin{align*}
\|aL^{1/2} f\|_2 \simeq \|\nabla_x f\|_2, 
\end{align*}
which identifies the domain of $L^{1/2}$ as the Sobolev space $\W^{1,2}$ since $a$ is invertible in $\L^\infty$. This conjecture is now solved.\footnote{In the case $a=1$, these are the results in \cite{CMMc} when $n=1$, \cite{HMc} when $n=2$ and \cite{AHLMcT} in all dimensions. When $a\ne 1$, this is proved in \cite{KM} when $n=1$ and then \cite{AKMc} in all dimensions.}

The $\H^p$-theory for the square root of $L$ consists in comparing $aL^{1/2}$ and $\nabla_x$ in $\H^p$. One estimate is the $\H^p$-boundedness of the Riesz transform $\nabla_x L^{-1/2} a^{-1}$, namely $\|\nabla_x f\|_{\H^p}\lesssim \|aL^{1/2}f\|_{\H^p}$, and then there is the reverse estimate. Of course, the left multiplication with the strictly elliptic function $a$ can be omitted when $p>1$. The conclusion is\footnote{This is proved in Section~\ref{sec: Riesz 2}. In the Lebesgue range $(1,\infty)$ it was first done in \cite{AT} when $n=1$ and $1<p<\infty$, and reproved in \cite{AMcN1}. For all dimensions, when $a=1$, the optimal range of $p$ within $(1,\infty)$ was settled in \cite{A} after earlier works of \cite{BK, Hofmann-Martell}. For discussions in the Hardy range $p\le 1$ when $a=1$, see \cite{HMMc}. Smaller intervals within the Lebesgue and Hardy range when $a\ne 1$ have been obtained in \cite{HMcP, AusSta, FMcP}.}
\begin{align*}
\|\nabla_x f\|_{\H^p} \lesssim \|aL^{1/2}f\|_{\H^p} \quad \text{if and only if $q_{-}(L) < p < q_+(L)$}
\end{align*}
 for the Riesz transform and that the reverse  estimate holds in a larger range, namely
\begin{align*}
\|aL^{1/2}f\|_{\H^p} \lesssim \|\nabla_x f\|_{\H^p} \quad \text{if $(q_{-}(L)_* \vee 1_*) < p < p_+(L)$}.
\end{align*}
What allows us to push the discussion to the range of exponents $1_*<p\le 1$ is the systematic use of Hardy and Hardy--Sobolev spaces $\IH_L^p$ and $\IH_L^{1,p}$ associated to $L$ that are defined using square functions involving the functional calculus of $L$. 

This foreshadows the main operator theoretic result of the monograph. Indeed, our approach to obtaining square function bounds and non-tangential maximal function bounds as in Theorem~\ref{thm: blockdir} and Theorem~\ref{thm: blockreg} below is to determine the ranges of exponents for which abstract Hardy and Hardy--Sobolev spaces associated to $L$ coincide with concrete spaces.\footnote{This approach is of course not new and the very reason why these spaces have been introduced. The latest development and exposition can be found in \cite{AA}. Elaborations on Hardy--Sobolev spaces  associated to $L$ were previously considered in \cite{HMMc} when $a=1$ and then in \cite{AusSta, FMcP} for general Dirac operators.} The upshot is that up to equivalent $p$-quasinorms, we are able to show
\begin{align}
\label{eq: HpL identification intro}
\IH_L^p = a^{-1}(\H^p \cap \L^2) \quad \text{if and only if $p_{-}(L) < p < p_+(L)$}
\end{align}
and
\begin{align}
\label{eq: H1pL identification intro}
\IH_L^{1,p} = \Hdot^{1,p} \cap \L^2 \quad \text{if $(q_{-}(L)_* \vee 1_*) < p < q_+(L)$},
\end{align}
where identification fails at the upper endpoint.\footnote{This is proved in Section~\ref{sec: Identification of adapted spaces}, except for the openness of $\cH(L)$ and $\cH^1(L)$ at the upper endpoint, which are obtained in Section~\ref{sec: Riesz 2} as a consequence of the results for the Riesz transform. When $a=1$ and $m=1$, results are obtained in \cite{HMMc} with  a different definition for the Hardy--Sobolev space and $p\le 2$, and limitations to $p>1$ for the identification for the Hardy space.} Even for the functional calculus \emph{per se} these identifications yield interesting new results.\footnote{See Section~\ref{sec: consequences identification}.} We now come to the boundary value problems.
\subsection{Main results on Dirichlet problems}
\label{intro:main}

\noindent Since for general systems the solutions might not be regular, we use the Whitney average variants of the non-tangential maximal function in order to pose our boundary value problems. Also we formulate the approach to the boundary in a non-tangential fashion using Whitney averages. When we get back to systems where solutions have pointwise values, these variants turn out to be equivalent to the usual  non-tangential pointwise control and  limits. More precisely, we let
\begin{align*}
\NT(F)(x) \coloneqq \sup_{t>0}  \bigg(\bariint_{W(t,x)} |F(s,y)|^2\, {\d s \d y}\bigg)^{1/2} \quad (x\in \R^n),
\end{align*}
with $W(t,x) \coloneqq (\nicefrac{t}{2}, 2t) \times B(x,t)$. 

For  $1<p<\infty$, the \emph{$\L^p$ Dirichlet problem}\index{Dirichlet problem!with $\L^p$-data} with non-tangential maximal control  and data $f\in \L^p(\R^{n}; \IC^m)$ consists in solving 
\begin{equation*}
(D)_{p}^\Le  \quad\quad
\begin{cases}
\Le u=0   & (\text{in }\reu), \\
\NT (u)\in \L^p(\R^n),   \\
\lim_{t \to 0} \bariint_{W(t,x)} |u(s,y)-f(x)| \, \d s \d y = 0 & (\text{a.e. } x\in \R^n).\end{cases}
\end{equation*}
For the endpoint problem $(D)_{1}^\Le$ the natural data class turns out to be a subspace of $\L^1$, namely the image of $\H^1$ under multiplication with the bounded function $a^{-1}$.\index{Dirichlet problem!with $a^{-1}\H^1$-data}

As usual, well-posedness means existence, uniqueness and continuous dependence on the data. Compatible well-posedness\index{well-posedness!compatible} means well-posedness together with the fact that the solution agrees with the energy solution that can be constructed via the Lax--Milgram lemma if the data $f$ also belongs to the boundary space $\Hdot^{\nicefrac{1}{2},2}(\R^n;\IC^m)$ for energy solutions.

Let us formulate our principal result on the Dirichlet problem, where we denote by $S$ the standard conical square function 
\begin{align*}
	S(F)(x) \coloneqq \bigg(\iint_{|x-y|<s} |F(s,y)|^2 \, \frac{\d s \d y}{s^{1+n}} \bigg)^{\frac{1}{2}} \quad (x \in \R^n).
\end{align*}

\begin{thm}[Dirichlet problem\index{Dirichlet problem!main result with $\L^p$/$a^{-1}\H^1$-data}]
\label{thm: blockdir} 
Let  $p \geq 1$ be such that $p_{-}(L)< p <  p_+(L)^*$. Given $f\in \L^p(\R^n; \IC^m)$ when $p>1$ and $f\in a^{-1} \H^1(\R^n; \IC^m)$ when $p=1$, the Dirichlet problem  $(D)_{p}^{\Le}$ is compatibly well-posed. The solution has the following additional properties.  
\begin{enumerate}
	\item There is comparability
	\begin{align*} 
	\|\NT(u)\|_{p} \simeq \|af\|_{\H^p} \simeq \|S (t\nabla u)\|_{p}.
	\end{align*}
	\item The non-tangential convergence improves to $\L^2$-averages
	\begin{align*}
	\qquad \lim_{t \to 0} \bariint_{W(t,x)} |u(s,y)-f(x)|^2 \, \d s \d y = 0  \quad (\text{a.e. } x\in \R^n).
	\end{align*}
	\item When $p<p_+(L)$, then   $au$ is of class\footnote{As usual, the notation $\C_0([0,\infty))$ means continuity and limit $0$ at infinity.\index{C@$\C_0([0,\infty))$}} $\C_0([0,\infty); \H^p(\R^n; \IC^m)) \cap \C^\infty((0,\infty); \H^p(\R^n; \IC^m))$ with $u(0, \cdot) = f$ and 
	\begin{align*}
	\sup_{t>0} \|au(t,\cdot)\|_{\H^p} \simeq \|af\|_{\H^p}.
		\end{align*}
	\item When $p \geq p_+(L)$, then for all $T>0$ and compact $K \subseteq \R^n$, $u$ is of class $\C([0,T]; \L^2(K; \IC^m))$ with $u(0, \cdot) = f$ and there is a constant $c=c(T,K)$ such that
	\begin{align*}
	\sup_{0<t \leq T} \|u(t,\cdot)\|_{\L^2(K)} \lesssim c \|f\|_p.
	\end{align*}	
\end{enumerate}
\end{thm}

As expected, the solution above is given by $u(t,x) =  \e^{-tL^{1/2}} f(x)$ if in addition we have $f \in \L^2$ and by an extension by density of this expression for  the respective topologies for general $f$. In the range $p<p_+(L)$ we can use the extension to a proper $\C_0$-semigroup on the data space, which explains the regularity result (iii).  However, and this was never observed before, the range of exponents in the statement exceeds by one Sobolev exponent the range provided by the semigroup theory.\footnote{When $a=1$, Mayboroda~\cite{May} dealt with variants where the $\L^2$-averages in the maximal functions are replaced with $\L^p$-averages. Her range of exponents is not the same and indeed, she shows that well-posedness is limited to the semigroup range.} 
This means that in this case $u$ is understood as a function of both variables $t$ and $x$ simultaneously that does not come from a semigroup action. 

Parts (i) and (iii) in the theorem remain true for the Poisson semigroup extension $u(t,x) =  \e^{-tL^{1/2}} f(x)$ of data $f \in a^{-1}(\H^p \cap \L^2)$, even when $p_-(L) < p < 1$. This is why we have systematically incorporated multiplication by $a$ in our estimates, although it can be omitted when $p>1$.\footnote{These estimates can be extended to $f$ in a closure of the data class for the quasinorm $\|a \cdot\|_{\H^p}$. However, since $\H^p$ does not embed into $\Lloc^1$ for $p<1$ and $a$ is not smooth, it is unclear whether this abstract extension has any reasonable (e.g.\ distributional) interpretation on the level of the boundary value problem. Even if $a=1$, (ii) has no meaning for us.}

For $1_*<p<\infty$, the \emph{$\H^p$ regularity problem}\index{Regularity problem} consists in solving, given $f\in \Hdot^{1,p}(\R^n;\IC^m)$,
\begin{equation*}
	(R)_{p}^\Le  \quad\quad
	\begin{cases}
		\Le u=0   & (\text{in } \reu), \\
		\NT(\nabla u)\in \L^p(\R^n),    \\
		\lim_{t \to 0} \bariint_{W(t,x)} |u(s,y)-f(x)| \, \d s \d y = 0 & (\text{a.e. } x\in \R^n).
	\end{cases}
\end{equation*}
As a (quasi-)Banach space  $\Hdot^{1,p}$ is a space of tempered distributions modulo constants but this point of view is not appropriate for the regularity problem. What we mean here is that the data $f$ is a tempered distribution such that $\nabla_x f \in \H^p$. By Hardy--Sobolev embeddings any such distribution is a locally integrable function and this gives a meaning to the boundary condition.\footnote{In fact, the condition $\NT(\nabla u)\in \L^p(\R^n)$ guarantees existence of a trace in $\Hdot^{1,p}$ in the sense of this limit at the boundary. See Appendix \ref{sec:technical}.}

Our principal result exhibits again an extended range of compatible well-posedness.\footnote{The fact that there is an extended range related to a Sobolev exponent down was observed by Mayboroda~\cite{May} when $a=1$, who establishes  $\|\NT(\nabla u)\|_{p}  \lesssim  \|\nabla_x f\|_{p}$  for $p\in (p_-(L)_*\vee 1, 2+\varepsilon]$ inspired from  the estimate  $\|L^{1/2}f\|_{\H^p} \lesssim \|\nabla_x f\|_{\H^p}$ in a similar  range from \cite{A}. We point out that Step~V in the proof of \cite[Thm.~4.1]{May} has a flaw that can be fixed (personal communication of S. Hofmann) or treated differently, see the argument in \cite{CMP2}.}
The solution is given by the Poisson semigroup if the data also belongs to $\L^2$ and appropriate extensions thereof in the general case.

\begin{thm}[Regularity problem\index{regularity problem!main result}]
\label{thm: blockreg} 
Let  $(q_{-}(L)_*\vee 1_*) < p < q_+(L)$. Then the regularity problem $(R)_{p}^{\Le}$ is compatibly well-posed. Given $f\in \Hdot^{1,p}(\R^n; \IC^m)$, the unique solution $u$ has the following additional properties.
\begin{enumerate}
	\item There are estimates
	\begin{align*}
	\hspace{1.5cm} \|\NT(\nabla u)\|_{p} \simeq  \|S (t\nabla \partial_{t}u)\|_{p} \simeq  \|\nabla_x f\|_{\H^p} \gtrsim \|g\|_{\H^p}
	\end{align*}
	with $g = -a{L^{1/2}}f$ being the conormal derivative of $u$, where the square root extends from $\Hdot^{1,p}(\R^n; \IC^m) \cap \W^{1,2}(\R^n; \IC^m)$ by density. 
	\item For a.e.\ $x \in \R^n$ and all $t>0$,
	\begin{align*}
		\bigg(\bariint_{W(t,x)} |u(s,y) - f(x)|^2 \, \d s \d y \bigg)^{\frac{1}{2}} \lesssim t \NT(\nabla u)(x).
	\end{align*}
	In particular, the non-tangential convergence improves to $\L^2$-averages. {Moreover, $\lim_{t \to 0} u(t,\cdot) = f$ in $\cD'(\R^n)$.}

	\item If $p \geq 1$, then for a.e.\ $x \in \R^n$,
	\begin{align*}
	\lim_{t \to 0} \bariint_{W(t,x)} \bigg|\begin{bmatrix} a\partial_t u \\ \nabla_x u \end{bmatrix} - \begin{bmatrix} g(x) \\ \nabla_x f(x) \end{bmatrix}\bigg|^2 \, \d s \d y = 0,
	\end{align*}
	where $g$  is as in (i).
	\item {$\nabla_ x u$ is of class $\C_{0}([0,\infty); \H^p(\R^n; \IC^m)) \cap \C^\infty((0,\infty); \H^p(\R^n; \IC^m))$ with $\nabla_x u(0,\cdot) = \nabla_x f$} and
	\begin{align*}
	\|\nabla_x f\|_{\H^p} \simeq \sup_{t > 0} \|\nabla_{x}u(t, \cdot )\|_{\H^p}.
	\end{align*}
	If $p<n$, then up to a constant\footnote{The constant is chosen via Hardy--Sobolev embeddings such that $f\in\L^{p^*}$.} $u \in \C_{0}([0,\infty); \L^{p^*}(\R^n; \IC^m)) \cap \C^\infty((0,\infty);  \L^{p^*}(\R^n; \IC^m))$ with $u(0,\cdot) = f$ and
	\begin{align*}
	\|f\|_{p^*} \leq \sup_{t > 0} \|u(t, \cdot )\|_{p^*} \lesssim \|\nabla_x f\|_{\H^p} + \|f\|_{p^*}.
	\end{align*}
	\item If $p>p_{-}(L)$, then $a\partial_{t}u$ is of class $ \C_{0}([0,\infty); \H^p(\R^n; \IC^m))$ and, with $g$ as in (i),
	\begin{align*}
	\qquad \|\NT(\partial_{t} u)\|_{p} \simeq \sup_{t\ge 0} \|a\partial_{t}u(t, \cdot )\|_{\H^p} \simeq \|g\|_{\H^p}\simeq \|\nabla_x f\|_{\H^p}.
	\end{align*}
\end{enumerate}
\end{thm}

As mentioned earlier, prior to these two results the situation was fully understood only in the case of boundary dimension $n=1$.\footnote{This is due to \cite{AT}, where existence and uniqueness are shown in the largest possible range $1<p<\infty$ as well as existence for a Dirichlet problem in the Hardy range $1_* = \nicefrac{1}{2}<p\leq 1$.  When $n \geq 2$ and $a=1$,  non-tangential maximal functions estimates pertaining to the Dirichlet and regularity problems first appeared in~\cite{May} and some related square functions estimates are in~\cite{AHM}. Uniqueness has not been considered in general, except for systems having regular solutions~\cite{HMaMo, HMiMo}. A possible strategy for general elliptic systems has been developed in \cite{AE}, but it only covers some smaller range of exponents when it comes to the block situation.}

One may wonder whether in the case $p_+(L) > n$ there are results for the Dirichlet problem with exponents `beyond $\infty$', which, in view of Sobolev embeddings, we think of corresponding to the homogeneous Hölder spaces $\Lamdot^\alpha(\R^n;\IC^m)$, $0 \leq \alpha < 1$, with the endpoint case $\Lamdot^0 \coloneqq \BMO$. We define the Carleson functional
\begin{align*}
C_\alpha (F)(x) \coloneqq \sup_{t>0} \frac{1}{t^\alpha} \bigg(\frac{1}{t^n} \int_0^{t} \int_{B(x,t)} |F(s,y)|^2 \frac{\d y \d s}{s} \bigg)^{1/2} \quad (x \in \R^n).
\end{align*}

For $\alpha\in (0,1)$, one formulation of the Dirichlet problem with data  $f\in \Lamdot^\alpha(\R^n;\IC^m)$\index{Dirichlet problem!with $\Lamdot^\alpha$-data} consists in solving
\begin{equation*}
(D)_{\Lamdot^\alpha}^\Le  \quad\quad
\begin{cases}
\Le u=0   & (\text{in } \reu), \\
C_\alpha (t\nabla u)\in \L^\infty(\R^n),   \\
\lim_{t \to 0} \bariint_{W(t,x)} |u(s,y)-f(x)| \, \d s \d y = 0 & (\text{a.e. } x\in \R^n).\end{cases}
\end{equation*}
The interior control from the Carleson functional alone implies existence of a non-tangential trace $f\in \Lamdot^\alpha(\R^n;\IC^m)$ as in the third line\footnote{We include a proof of the trace theorem in Appendix~\ref{sec:technical}.}, so that this is the weakest possible formulation of the boundary behavior. Again, we regard $\Lamdot^\alpha$ as a space of functions to make sense of the limit condition. This non-tangential trace also satisfies $\|\NTsharpalpha (u-f)\|_\infty \lesssim \|C_\alpha (t\nabla u)\|_\infty$, where on the left-hand side we use the sharp functional on Whitney averages\index{Whitney average functional!sharp ($\NTsharpalpha$)}
\begin{align*}
	\NTsharpalpha (u-f)(x) \coloneqq \sup_{t>0} \frac{1}{t^\alpha} \bigg(\bariint_{W(t,x)} |u(s,y)-f(y)|^2 \, \d s \d y \bigg)^{1/2} \ (x \in \R^n).
\end{align*}
Such a trace result is not available for $\alpha = 0$ and we formulate the boundary behavior for the endpoint problem differently, using convergence of Cesàro averages, which is natural from the point of view of both our construction and our approach to uniqueness theorems:\index{Dirichlet problem!with $\BMO$-data}
\begin{equation*}
	(D)_{\Lamdot^0}^\Le  \quad\quad
	\begin{cases}
		\Le u=0   & (\text{in } \reu), \\
		C_0 (t\nabla u)\in \L^\infty(\R^n),   \\
		\lim_{t \to 0} \barint_{t}^{2t} |u(s,\cdot)-f| \, \d s = 0 & (\text{in } \Lloc^2(\R^n; \IC^m)).
	\end{cases}
\end{equation*}

The discussion of non-tangential traces naturally leads us to formulating a modified $\Lamdot^\alpha$ Dirichlet problem
\begin{equation*}
	(\wt D)_{\Lamdot^\alpha}^\Le  \quad\quad
	\begin{cases}
		\Le u=0   & (\text{in } \reu), \\
		\NTsharpalpha(u-f)\in \L^\infty(\R^n),   \\
		\lim_{t \to 0} \bariint_{W(t,x)} |u(s,y)-f(x)| \, \d s \d y = 0 & (\text{a.e. } x\in \R^n).\end{cases}
\end{equation*}
As we have seen above, this second problem is \emph{a priori} comparable to the first one when $\alpha>0$.\footnote{Uniqueness for the $\BMO$ Dirichlet problem with interior Carleson control and Whitney average convergence at the boundary appears to be out of reach. See \cite{MaMiMiMi-Survey, MaMiMiMi} for a very recent account on such Fatou-type theorems in the case of elliptic systems with \emph{constant} coefficients.}

We obtain compatible well-posedness for both problems in the same range of exponents. In order to formulate the theorem, and systematically throughout this book, we denote by $L^\sharp$ the boundary operator for the adjoint equation $\Le^* u = 0$, that is $L^\sharp= -(a^*)^{-1} \div_x d^* \nabla_x$.

\begin{thm}[$\Lamdot^\alpha$ Dirichlet problem\index{Dirichlet problem!main result with $\Lamdot^\alpha$-data}]
\label{thm: Holder-dir}
 Suppose that $p_+(L) > n$ and that $0 \leq \alpha < 1 - \nicefrac{n}{p_+(L)}$. Then the Dirichlet problems $(D)_{\Lamdot^\alpha}^{\Le}$ and  $(\wt D\smash)_{\Lamdot^\alpha}^{\Le}$ are compatibly well-posed. Given $f\in \Lamdot^\alpha(\R^n; \IC^m)$, the unique solution $u$ is the same for both problems and has the following additional properties.
\begin{enumerate}
	\item There is comparability
	\begin{align*}
	\|C_\alpha(t \nabla u)\|_\infty \simeq \|f\|_{\Lamdot^\alpha}.
	\end{align*}
	\item One has the upper bound  
	\begin{align*}
	\| \NTsharpalpha(u-f)\|_\infty \lesssim \|f\|_{\Lamdot^\alpha}
	\end{align*}
	and convergence
	\begin{align*}
	\qquad \qquad \lim_{t \to 0} \bariint_{W(t,x)} |u(s,y)-f(x)|^2 \, \d s \d y = 0 \quad (\text{a.e. } x\in \R^n).
	\end{align*}
	In addition, $u$ is of class $\C([0,T]; \Lloc^2(\R^n; \IC^m))$ with $u(0, \cdot) = f$ for every $T>0$.
	\item If, moreover, $p_{-}(L^\sharp) < 1$ and $\alpha < n(\nicefrac{1}{p_-(L^\sharp)} -1)$, then $u$ is of class $\C_0([0,\infty); \Lamdot^\alpha_{\text{weak}^*}(\R^n; \IC^m)) \cap \C^\infty((0,\infty); \Lamdot^\alpha_{\text{weak}^*}(\R^n; \IC^m))$ and
	\begin{align*}
	\sup_{t>0} \|u(t,\cdot)\|_{\Lamdot^\alpha} \simeq \|f\|_{\Lamdot^\alpha}.
	\end{align*}
	In addition, $u$ is of class $\Lamdot^\alpha(\overline{\reu};\IC^m)$, with 
	\begin{align*}
	\|u\|_{\Lamdot^\alpha(\overline{\reu})}\lesssim \|f\|_{\Lamdot^\alpha}.
	\end{align*}
\end{enumerate}
\end{thm}

Since $\Lamdot^\alpha \cap \L^2$ is \emph{not} dense in $\Lamdot^\alpha$, we cannot extend the Poisson semigroup to the boundary space by density. In (iii), $\Lamdot^\alpha$ is considered as the dual space of $\H^p$, where $\alpha = n(\nicefrac{1}{p} -1)$, with the weak$^*$ topology. The assumption in (iii) implies $p_+(L)=\infty$ and that the solution can be constructed by duality, using the extension of the Poisson semigroup for $L^* = a^*L^\sharp(a^*)^{-1}$ to $\H^p$. Therefore the solution keeps the $\Lamdot^\alpha$-regularity in the interior. This construction has appeared earlier.\footnote{References are \cite{AusSta, AM, HMiMo}.}

The construction of the solution under the mere assumption that $p_+(L) > n$ is much more general and we have
\begin{align*}
	u(t,x) = \lim_{j \to \infty} \e^{-t L^{1/2}} (\ind_{\{|\,\cdot\,| < 2^j\}} f)(x),
\end{align*}
where $p_+(L) > n$ is used already to prove convergence of the right-hand side in $\Lloc^2(\reu; \IC^m)$. This opens the possibility of uniquely solving Dirichlet problems for H\"older continuous (or $\BMO$) data, while producing solutions that have no reason to be in the same class in the interior of the domain. To the best of our knowledge this phenomenon is observed for the first time. Note also that $p_+(L) > n$ always holds in dimension $n \leq4$, so that in these dimensions both $\BMO$ Dirichlet problems are compatibly well-posed.

\subsection{Dirichlet problems with fractional spaces of data}
\label{intro:fractionalDP}

\noindent  If we think of the Dirichlet problem $(D)_p^\Le$ as a boundary value problem with regularity $s=0$ for the data and the regularity problem $(R)_p^\Le$ as a Dirichlet problem with regularity $s=1$,  we can depict the exponents for both problems simultaneously in an $(\nicefrac{1}{p},s)$-diagram. There are two classical scales of data spaces to fill the intermediate area of points with $0<s<1$: The homogeneous Hardy--Sobolev spaces $\Hdot^{s,p}$ that can be obtained from the endpoints by complex interpolation and the homogeneous Besov spaces $\Bdot^{s,p}$ that result from real interpolation.\footnote{We give a detailed account on all sorts of relevant function spaces in Section~\ref{sec: function spaces}.}

For $0 < p \leq \infty$ and $0<s<1$ satisfying $\nicefrac{1}{p} < 1 + \nicefrac{s}{n}$\footnote{\label{Footnote: Sobolev}When $0<p<\infty$, this Sobolev-type condition characterizes the spaces that can be obtained by interpolation between data spaces for the Dirichlet problem ($\L^p$ with $p > 1$) and the regularity problem ($\H^p$ with $p>1_*$), see Section~\ref{subsec: complex interpolation}. In particular, it is the natural restriction guaranteeing that all distributions in $\Xdot^{s,p}$ are locally integrable functions. The spaces $\Bdot^{s,\infty}$ and $\Hdot^{s,\infty}$ also have this property, see Section~\ref{subsec: homogeneous smoothness spaces}.}, the Dirichlet problem with data $f \in \Bdot^{s,p}(\R^n; \IC^m)$\index{Dirichlet problem!with $\Bdot^{s,p}$-data} consists in solving 
\begin{equation*}
	(D)_{\Bdot^{s,p}}^\Le  \quad\quad
	\begin{cases}
		\Le u=0   & (\text{in } \reu), \\
		W (t^{1-s} \nabla u)\in \L^p(\reu; \frac{\d t \d x}{t}),   \\
		\lim_{t \to 0} \bariint_{W(t,x)} |u(s,y)-f(x)| \, \d s \d y = 0 & (\text{a.e. } x\in \R^n),
	\end{cases}
\end{equation*}
where $W(F)$ is the Whitney average functional
\begin{align*}
W(F)(t,x) \coloneqq \bigg(\bariint_{W(t,x)} |F(s,y)|^2 \, \d s \d y \bigg)^{\frac{1}{2}} \quad ((t,x) \in \reu).
\end{align*}
For $0 < p < \infty$ and $0<s<1$  satisfying $\nicefrac{1}{p} < 1 + \nicefrac{s}{n}$, the Dirichlet problem with data $f \in \Hdot^{s,p}(\R^n; \IC^m)$\index{Dirichlet problem!with $\H^{s,p}$-data} consists in solving 
\begin{equation*}
(D)_{\Hdot^{s,p}}^\Le  \quad\quad
\begin{cases}
	\Le u=0   & (\text{in } \reu), \\
	S (t^{1-s} \nabla u)\in \L^p(\R^n),   \\
	\lim_{t \to 0} \bariint_{W(t,x)} |u(s,y)-f(x)| \, \d s \d y = 0 & (\text{a.e. } x\in \R^n),
\end{cases}
\end{equation*}
where $S$ is the same conical square function as before.\footnote{Boundary value problems for general elliptic equations ($m=1$)  with data of fractional regularity have been pioneered by Barton--Mayboroda~\cite{BM}. They treat $\Bdot^{s,p}$-data for equations  with the de~Giorgi--Nash--Moser property. This assumption was then removed in the first-order approach by Amenta along with the first author~\cite{AA} and their approach includes the problems with $\Hdot^{s,p}$-data. Thanks to the block structure we do not have to include a limiting condition for $u$ as $t \to \infty$ in the formulation of our fractional Dirichlet problems. Such a condition appears in the general framework of  \cite{AA} but not in \cite{BM}.}

For $p=\infty$ we can identify $\Bdot^{s,\infty} = \Lamdot^s$, so that $(D)_{\Bdot^{s,\infty}}^\Le$ is a third formulation of a Dirichlet problem for that space of data. The endpoint problems for the Hardy--Sobolev scale are formulated for data in Strichartz' $\BMO$-Sobolev spaces $\Hdot^{s,\infty} = \BMOdot^s$\index{Dirichlet problem!with $\BMOdot^s$-data} and consist in solving 
\begin{equation*}
(D)_{\Hdot^{s,\infty}}^\Le  \quad\quad
\begin{cases}
	\Le u=0   & (\text{in } \reu), \\
	\C_0 (t^{1-s} \nabla u)\in \L^\infty(\R^n),   \\
	\lim_{t \to 0} \bariint_{W(t,x)} |u(s,y)-f(x)| \, \d s \d y = 0 & (\text{a.e. } x\in \R^n).
\end{cases}
\end{equation*}

We note that the approach to the boundary in these problems is not in the sense of the usual trace theory, that is by extension of the restriction map to the boundary defined on smooth functions. In fact, this approach would work for Besov spaces\footnote{This is the point of view taken in \cite{BM}. See also \cite{AA}.} but not for Hardy--Sobolev spaces, which are not trace spaces in this sense. Our choice of a non-tangential convergence of Whitney averages has  one main advantage, valid for all situations: each interior control implies existence of a unique measurable function $f$, called non-tangential trace (in the sense of Whitney averages), such that the third condition holds, whether or not $u$ is a weak solution to $\Le u = 0$. In this sense, we prescribe the boundary limit in the weakest possible form. If, via a trace operator, $\lim_{t \to 0} u(t,\cdot)$ also exists in the sense of distributions (modulo constants), then the two notions of boundary trace coincide (modulo constants).\footnote{All this is shown in Appendix~\ref{sec:technical}. Similar trace theorems appear in \cite[Sec.~6.6]{AA}, where they are used to derive non-tangential convergence of the solution at the boundary \emph{a posteriori}.} The same limit condition was taken in the boundary value problems from the previous section (except for one of the Dirichlet problems with $\BMO$-data). We stress again that we consider the data spaces as classes of measurable functions and not as distributions (modulo constants) and that this is possible due to the assumption  $\nicefrac{1}{p} < 1 + \nicefrac{s}{n}$.

In the Figures~\ref{fig: diagram general}, \ref{fig: diagram p+large}, \ref{fig: diagram p-small} below we are going to unite \emph{all} results of the monograph on (compatible) well-posedness of boundary value problems. First, we collect compatible well-posedness results from the previous section on thick horizontal boundary segments at $s=0$ and $s=1$. For $p=\infty$, we can represent these results also on a thick vertical segment at $\nicefrac{1}{p}=0$.  Empty circles indicate boundary points that are not contained in a segment of well-posedness. This allows us to create a map  $f \mapsto u$ for different values of $(\nicefrac{1}{p},s)$ on these lines and, roughly speaking, we can interpolate to fill in a shaded region for compatible solvability of both fractional problems.\footnote{The fact that not only the data spaces but also the interior control from the functionals $S$ and $W$ interpolate, shows again that these are natural classes of solutions from our perspective.} 

Of course, interpolation does not preserve uniqueness. Still, we shall be able to show uniqueness (and hence compatible well-posedness) even in a possibly larger region than for existence of a solution.\footnote{{The corresponding regions and the proofs of the uniqueness theorems in the non-fractional and fractional cases can all be found in Section~\ref{sec: uniqueness}, together with illustrations of uniqueness regions (Figures~\ref{fig: diagram-uniqueness-general} and \ref{fig: diagram-uniqueness-p+large}). That section can be read almost independently of the rest of the monograph and we use quite original techniques for establishing uniqueness.}}

\begin{thm}
\label{thm: uniqueness fractional}
Let $0<s<1$ and $1_*<p\le \infty$.   If $(\nicefrac{1}{p},s)$ belongs  to the region displayed in Figure~\ref{fig: diagram general}, Figure~\ref{fig: diagram p+large} or Figure~\ref{fig: diagram p-small} (including the thick vertical segment), then $(D)_{\Xdot^{s,p}}^\Le$ is compatibly well-posed. 
\end{thm}

As customary, we obtain continuous dependence on the data: the interior control is bounded by the data in the boundary space. For the problems corresponding to all thick segments we have also seen the reverse estimates in the previous section.    
{Various additional regularity properties in the spirit of Theorems~\ref{thm: blockdir} - \ref{thm: Holder-dir} hold depending on the particular boundary space.\footnote{Precise results are stated and proved in Section~\ref{sec: fractional}.}}

A color code allows us to distinguish different zones {that explain the relation of the corresponding well-posedness results with the first- and second-order operator theory that we develop in parallel. A reader who is not familiar with these tools (yet) might ignore the different colors for the time being and focus only on the shape of the regions.}
\begin{itemize}
	\item Gray corresponds to what can be obtained from the theory of first-order $DB$-adapted spaces in \cite{AA}, although we shall argue independently of \cite{AA} when it comes to solving boundary value problems.
	\item Blue shows additional results obtained from the theory of $L$-adapted spaces.
	\item Red indicates results outside of the theory of operator-adapted spaces.
\end{itemize}
All shaded regions in the strip $0<s<1$ capture a situation that is common to Hardy--Sobolev and Besov data and we set $\Xdot$ to designate $\Hdot$ or $\Bdot$.   
They depict three different cases: first $p_+(L) \leq n$, next $p_+(L) > n$ but $p_-(L^\sharp) \geq 1$ and eventually  $p_-(L^\sharp) < 1$, {which turns out to imply $p_+(L) = \infty$ by duality}.\footnote{Let us mention that  the diagrams are up to scale when $p_-(L)\ge1$ but not when  $p_-(L) < 1$. In this latter case, the top right vertex of the blue region is always situated at $(\nicefrac{1}{p_-(L)}, 1)$, while the bottom right vertex would be $(1,0)$. 
}

We begin by illustrating the situation when $p_+(L) \leq n$. In this case we obtain the segment on the bottom line for $s=0$ and the top line for $s=1$ from Theorems~\ref{thm: blockdir} and \ref{thm: blockreg}, respectively. This leads to Figure~\ref{fig: diagram general}. In all such figures we shall write $p_+^L$ instead of $p_+(L)$ and so on for the sake of a clearer typeset.

\begin{figure}[h]
\begin{center}
\begin{tikzpicture}[scale=2.4]

	\newcommand\fracspace{\vphantom{\frac{1}{1}}};
	\def\dimension{6};
	\def\xlength{3.5};
	\pgfmathsetmacro\xstretch{\xlength/(1+1/\dimension)}; 
	
	\pgfmathsetmacro\UpSob{\xstretch*(1/2-1/\dimension)-0.2};
	\pgfmathsetmacro\UpSobStar{\xstretch*(\UpSob/\xstretch-1/\dimension)};
	\pgfmathsetmacro\LowSob{\xstretch*(1/2+1/\dimension)+0.5};
	\pgfmathsetmacro\LowSobStar{\xstretch*(\LowSob/\xstretch+1/\dimension)};
	\pgfmathsetmacro\UpGradient{\xstretch*(1/2)-0.2};
	\pgfmathsetmacro\UpGradientDual{\xstretch*(1/2)+0.3};
	\pgfmathsetmacro\Half{\xstretch*0.5};
	\pgfmathsetmacro\One{\xstretch*1};
	\coordinate (P00) at (1+\UpGradient,2);
	\coordinate (P01) at (1+\LowSob,2);
	\coordinate (P10) at (1+\UpSob,0);
	\coordinate (P11) at (1+\UpGradientDual,0);
	\coordinate (P12) at (1+\LowSob,0);
	\coordinate (DirExtra) at (1+\UpSobStar,0);
	\coordinate (RegExtra) at (1+\LowSobStar,2);

	\draw [thin] (1,2) -- (1+\xlength,2); 
	\draw [thin] (1,0) -- (1+\xlength,0); 
	
	\draw [thick,->] (1,-0.5) -- (1+\xlength+0.2,-0.5);
	\draw [fill=black] (1+\Half,-0.5) circle [radius = .5pt];
	\node [below] at (1+\Half,-0.5) {$\frac{1}{2 \fracspace}$};
	\draw [fill=black] (1+\UpGradient,-0.5) circle [radius = .5pt];
	\node [below] at (1+\UpGradient,-0.5) {$\frac{1}{q_+^L\fracspace}$};
	\draw [fill=black] (1+\UpGradientDual,-0.5) circle [radius = .5pt];
	\node [below] at (1+\UpGradientDual,-0.5) {$\frac{1}{(q_+^{L^\sharp})'\fracspace}$};
	\draw [fill=black] (1+\UpSob,-0.5) circle [radius = .5pt];
	\node [below] at (1+\UpSob,-0.5) {$\frac{1}{p_+^L\fracspace}$};
	\draw [fill=black] (1+\LowSob,-0.5) circle [radius = .5pt];
	\node [below] at (1+\LowSob,-0.5) {$\frac{1}{p_-^L \vee 1 \fracspace}$};
	\draw [fill=black] (1+\UpSobStar,-0.5) circle [radius = .5pt];
	\node [below] at (1+\UpSobStar,-0.5) {$\frac{1}{(p_+^L)^*\fracspace}$};
	\draw [fill=black] (1+\LowSobStar,-0.5) circle [radius = .5pt];
	\node [below] at (1+\LowSobStar,-0.5) {$\frac{1}{{(p_-^L)_* \vee 1_*\fracspace}}$};
	\draw [fill=black] (1+\xlength,-0.5) circle [radius = .5pt];
	\node [below] at (1+\xlength,-0.5) {$\frac{n+1}{n \fracspace}$};
	\node [right] at (1+\xlength+0.2,-0.5) {$\frac{1}{p \fracspace}$};
	\draw [fill=black] (1,-0.5) circle [radius = .5pt];
	\node [below] at (1,-0.5) {$0$};
	
	\draw [thick,->] (0.7,0) -- (0.7,2.2);
	\node [above] at (0.7,2.2) {$s$};
	\draw [fill=black] (0.7,2) circle [radius = .5pt];
	\node [left] at (0.7,2) {$1$};
	\draw [fill=black] (0.7,0) circle [radius = .5pt];
	\node [left] at (0.7,0) {$0$};

	\draw [thin,dotted] (1,2) -- (1,0); 
	\draw [thin,dotted] (1+\LowSob,2) -- (1+\LowSob,0); 
	\draw [thin,dotted] (1+\LowSobStar,2) -- (1+\LowSobStar,0); 
	\draw [thin,dotted] (1+\UpGradient,2) -- (1+\UpGradient,0);
	\draw [thin,dotted] (1+\xlength,2) -- (1+\xlength,0); 
	
	\path [fill=lightgray, opacity = 0.6] (P00)--(P01)--(P11)--(P10)--(P00);
	\path [fill=blue!80!black, opacity = 0.4] (P11)--(P12)--(RegExtra)--(P01)--(P11);
	\path [fill=red!80!black, opacity = 0.4] (DirExtra)--(P10)--(P00)--(DirExtra);
	\draw [ultra thick, gray] (P00) -- (P01); 
	\draw [ultra thick, blue!80!black] (P01) -- (RegExtra); 
	\draw [ultra thick, gray] (P10) -- (P11); 
	\draw [ultra thick, red!80!black] (DirExtra) -- (P10); 
	\draw [ultra thick, blue!80!black] (P11) -- (P12); 
	\draw [fill=white] (P00) circle [radius = .75pt];
	\draw [fill=blue!80!black] (P01) circle [radius = .75pt];
	\draw [fill=red!80!black] (P10) circle [radius = .75pt];
	\draw [fill=blue!80!black] (P11) circle [radius = .75pt];
	\draw [fill=white] (P12) circle [radius = .75pt];
	\draw [fill=white] (DirExtra) circle [radius = .75pt];
	\draw [fill=white] (RegExtra) circle [radius = .75pt];
\end{tikzpicture}
\end{center}
\caption{Compatible well-posedness region for Besov and Hardy--Sobolev data when $p_+(L) \leq n$.
}
\label{fig: diagram general}
\end{figure}

In the case $p_+(L)>n$ we can extend the  bottom line to exponents `beyond infinity', using  Theorem~\ref{thm: Holder-dir}. The point corresponding to compatible well-posedness of $(D)_{\Lamdot^\alpha}^\Le$ is $(-\nicefrac{\alpha}{n},0)$. We shall see that this also leads to compatible well-posedness of $(D)_{\Bdot^{\alpha,\infty}}^\Le$ at   $(0,\alpha)$ as stated. A similar result holds for $(D)_{\Hdot^{\alpha,\infty}}^\Le$ at   $(0,\alpha)$.\footnote{See Proposition~\ref{prop: BMO-Besov-dir}.} Figure~\ref{fig: diagram p+large} illustrates this extension in the case that $p_+(L) > n$ but $p_-(L^\sharp) \geq 1$. This is the generic situation in dimensions $n=3,4$.

\begin{figure}[h]
\begin{center}
\begin{tikzpicture}[scale=2.4]
	\newcommand\fracspace{\vphantom{\frac{1}{1}}};
	\def\dimension{6};
	\def\xlength{3.5};
	\pgfmathsetmacro\xstretch{\xlength/(1+1/\dimension)}; 
	
	\pgfmathsetmacro\UpSob{\xstretch*(1/\dimension)-0.3};
	\pgfmathsetmacro\UpSobStar{\xstretch*(\UpSob/\xstretch-1/\dimension)};
	\pgfmathsetmacro\UpSobStarBesov{(-2)*\UpSobStar*\dimension/\xlength};
	\pgfmathsetmacro\LowSob{\xstretch*(1/2+1/\dimension)+0.5};
	\pgfmathsetmacro\LowSobStar{\xstretch*(\LowSob/\xstretch+1/\dimension)};
	\pgfmathsetmacro\UpGradient{\xstretch*(1/2)-0.2};
	\pgfmathsetmacro\UpGradientDual{\xstretch*(1/2)+0.3};
	\pgfmathsetmacro\Half{\xstretch*0.5};
	\pgfmathsetmacro\One{\xstretch*1};
	\coordinate (P00) at (1+\UpGradient,2);
	\coordinate (P01) at (1+\LowSob,2);
	\coordinate (P10) at (1+\UpSob,0);
	\coordinate (P11) at (1+\UpGradientDual,0);
	\coordinate (P12) at (1+\LowSob,0);
	\coordinate (DirExtra) at (1,0);
	\coordinate (DirHolder) at (1+\UpSobStar,0);
	\coordinate (DirHolderBesov) at (1, \UpSobStarBesov);
	\coordinate (RegExtra) at (1+\LowSobStar,2);
	\coordinate (P000) at (1-\xstretch/\dimension,0);

	\draw [thin] (0.3,2) -- (1+\xlength,2); 
	\draw [thin] (0.3,0) -- (1+\xlength,0); 
	
	\draw [thick,->] (0.3,-0.5) -- (1+\xlength+0.2,-0.5);
	\draw [fill=black] (1+\Half,-0.5) circle [radius = .5pt];
	\node [below] at (1+\Half,-0.5) {$\frac{1}{2 \fracspace}$};
	\draw [fill=black] (1+\UpGradient,-0.5) circle [radius = .5pt];
	\node [below] at (1+\UpGradient,-0.5) {$\frac{1}{q_+^L\fracspace}$};
	\draw [fill=black] (1+\UpGradientDual,-0.5) circle [radius = .5pt];
	\node [below] at (1+\UpGradientDual,-0.5) {$\frac{1}{(q_+^{L^\sharp})'\fracspace}$};
	\draw [fill=black] (1+\UpSob,-0.5) circle [radius = .5pt];
	\node [below] at (1+\UpSob,-0.5) {$\frac{1}{p_+^L\fracspace}$};
	\draw [fill=black] (1+\LowSob,-0.5) circle [radius = .5pt];
	\node [below] at (1+\LowSob,-0.5) {$\frac{1}{p_-^L \vee 1 \fracspace}$};
	\draw [fill=black] (1+\UpSobStar,-0.5) circle [radius = .5pt];
	\node [below] at (1+\UpSobStar-0.05,-0.5) {$\frac{1}{p_+^L\fracspace} \!- \!\frac{1}{n\fracspace}$};
	\draw [fill=black] (1+\LowSobStar,-0.5) circle [radius = .5pt];
	\node [below] at (1+\LowSobStar,-0.5) {$\frac{1}{{(p_-^L)_* \vee 1_*\fracspace}}$};
	\draw [fill=black] (1+\xlength,-0.5) circle [radius = .5pt];
	\node [below] at (1+\xlength,-0.5) {$\frac{n+1}{n \fracspace}$};
	\node [right] at (1+\xlength+0.2,-0.5) {$\frac{1}{p \fracspace}$};
	\draw [fill=black] (1,-0.5) circle [radius = .5pt];
	\node [below] at (1,-0.5) {$0$};
	\draw [fill=black] (0.3,-0.5) circle [radius = .5pt];
	\node [below] at (0.2,-0.5) {$-\frac{1}{n \fracspace}$};
	
	\draw [thick,->] (-0,0) -- (0,2.2);
	\node [above] at (0,2.2) {$s$};
	\draw [fill=black] (0,2) circle [radius = .5pt];
	\node [left] at (0, \UpSobStarBesov) {$\hspace{-30pt}1\!-\! \frac{n}{p_+^L}$};
	\draw [fill=black] (0,  \UpSobStarBesov) circle [radius = .5pt];	
	\node [left] at (0,2) {$1$};
	\draw [fill=black] (0,0) circle [radius = .5pt];
	\node [left] at (0,0) {$0$};

	\draw [thin,dotted] (1,2) -- (1,0); 
	\draw [thin,dotted] (0.3,2) -- (0.3,0); 
	\draw [thin,dotted] (1+\UpGradient,2) -- (1+\UpGradient,0);
	\draw [thin,dotted] (1+\LowSob,2) -- (1+\LowSob,0);
	\draw [thin,dotted] (1+\xlength,2) -- (1+\xlength,0); 
	
	\draw [thin,dotted] (0.3, \UpSobStarBesov) -- (1,\UpSobStarBesov); 
	
	\path [fill=lightgray, opacity = 0.6] (P00)--(P01)--(P11)--(P10)--(P00);
	\path [fill=blue!80!black, opacity = 0.4] (P11)--(P12)--(RegExtra)--(P01)--(P11);
	\path [fill=red!80!black, opacity = 0.4] (DirExtra)--(P10)--(P00)--(DirExtra);
	\path [fill=red!80!black, opacity = 0.4]
	(DirExtra) -- (DirHolderBesov) -- (P00) -- (DirExtra);
	\draw [ultra thick, gray] (P00) -- (P01); 
	\draw [ultra thick, blue!80!black] (P01) -- (RegExtra); 
	\draw [ultra thick, gray] (P10) -- (P11); 
	\draw [ultra thick, red!80!black] (DirExtra) -- (P10); 
	\draw [ultra thick, red!80!black] (DirHolder) -- (DirExtra);
	\draw [ultra thick, red!80!black] (DirExtra) -- (DirHolderBesov); 
	\draw [ultra thick, blue!80!black] (P11) -- (P12); 
	\draw [fill=white] (P00) circle [radius = .75pt];
	\draw [fill=blue!80!black] (P01) circle [radius = .75pt];
	\draw [fill=red!80!black] (P10) circle [radius = .75pt];
	\draw [fill=blue!80!black] (P11) circle [radius = .75pt];	
	\draw [fill=white] (P12) circle [radius = .75pt];
	\draw [fill=red!80!black] (DirExtra) circle [radius = .75pt];
	\draw [fill=white] (RegExtra) circle [radius = .75pt];
	\draw [fill=white] (DirHolder) circle [radius = .75pt];
	\draw [fill=white] (DirHolderBesov) circle [radius = .75pt];
\end{tikzpicture}
\end{center}
\caption{Compatible well-posedness region for Besov and Hardy--Sobolev data when $p_+(L) > n$ but $p_-(L^\sharp) \geq 1$.}
\label{fig: diagram p+large}
\end{figure}

Figure~\ref{fig: diagram p-small} describes the case when $p_-(L^\sharp) < 1$, which happens for instance when $n=1,2$ or for special classes of systems such as equations ($m=1$) with real-valued coefficients $d$.\footnote{More examples are given in Section~\ref{subsec: Stability and examples}.} The caption of Figure~\ref{fig: diagram p-small} contains a more specific discussion of this particular case.

\begin{figure}[h]
\begin{center}
\begin{tikzpicture}[scale=2.4]

	\newcommand\fracspace{\vphantom{\frac{1}{1}}};
	\def\dimension{6};
	\def\xlength{3.5};
	\pgfmathsetmacro\xstretch{\xlength/(1+1/\dimension)}; 
	
	\pgfmathsetmacro\UpSob{-0.3};
	\pgfmathsetmacro\LowSob{\xstretch*(1/2+1/\dimension)+0.5};
	\pgfmathsetmacro\LowSobStar{\xstretch*(\LowSob/\xstretch+1/\dimension)};
	\pgfmathsetmacro\HolderMoritz{(-2)*\UpSob*\dimension/\xlength};
	\pgfmathsetmacro\UpGradient{\xstretch*(1/2)-0.2};
	\pgfmathsetmacro\UpGradientDual{\xstretch*(1/2)+0.3};
	\pgfmathsetmacro\Half{\xstretch*0.5};
	\pgfmathsetmacro\One{\xstretch*1};
	\coordinate (P00) at (1+\UpGradient,2);
	\coordinate (P01) at (1+\LowSob,2);
	\coordinate (P10) at (1+\UpSob,0);
	\coordinate (P11) at (1+\UpGradientDual,0);
	\coordinate (P12) at (1+\LowSob,0);
	\coordinate (DirExtra) at (1.2,0);
	\coordinate (DirHolder) at (0.3,0);
	\coordinate (DirHolderBesov) at (1.2, 2);
	\coordinate (HolderMoritz) at (1.2, \HolderMoritz);
	\coordinate (HolderAlex) at (1.2, 0.3);
	\coordinate (RegExtra) at (1+\LowSobStar,2);
	\coordinate (P000) at (1-\xstretch/\dimension,0);

	\draw [thin] (0.3,2) -- (1+\xlength,2); 
	\draw [thin] (0.3,0) -- (1+\xlength,0); 
	
	\draw [thick,->] (0.3,-0.5) -- (1+\xlength+0.2,-0.5);
	\draw [fill=black] (1+\Half,-0.5) circle [radius = .5pt];
	\node [below] at (1+\Half,-0.5) {$\frac{1}{2 \fracspace}$};
	\draw [fill=black] (1+\UpGradient,-0.5) circle [radius = .5pt];
	\node [below] at (1+\UpGradient,-0.5) {$\frac{1}{q_+^L\fracspace}$};
	\draw [fill=black] (1+\UpGradientDual,-0.5) circle [radius = .5pt];
	\node [below] at (1+\UpGradientDual,-0.5) {$\frac{1}{(q_+^{L^\sharp})'\fracspace}$};
	\draw [fill=black] (1+\UpSob,-0.5) circle [radius = .5pt];
	\node [below] at (1+\UpSob,-0.5) {$1\!-\!\frac{1}{p_-^{L^\sharp}\fracspace}$};
	\draw [fill=black] (1+\LowSob,-0.5) circle [radius = .5pt];
	\node [below] at (1+\LowSob,-0.5) {$\frac{1}{p_-^L \vee 1 \fracspace}$};
	\draw [fill=black] (1+\LowSobStar,-0.5) circle [radius = .5pt];
	\node [below] at (1+\LowSobStar,-0.5) {$\frac{1}{{(p_-^L)_* \vee 1_*\fracspace}}$};
	\draw [fill=black] (1+\xlength,-0.5) circle [radius = .5pt];
	\node [below] at (1+\xlength,-0.5) {$\frac{n+1}{n \fracspace}$};
	\node [right] at (1+\xlength+0.2,-0.5) {$\frac{1}{p \fracspace}$};
	\draw [fill=black] (1.2,-0.5) circle [radius = .5pt];
	\node [below] at (1.2,-0.5) {$0$};
	\draw [fill=black] (0.3,-0.5) circle [radius = .5pt];
	\node [below] at (0.2,-0.5) {$-\frac{1}{n \fracspace}$};
	
	\draw [thick,->] (0,0) -- (0,2.2);
	\node [above] at (0,2.2) {$s$};
	\draw [fill=black] (0,2) circle [radius = .5pt];
	\node [left] at (0,2) {$1$};
	\draw [fill=black] (0,0) circle [radius = .5pt];
	\node [left] at (0,0) {$0$};
	\node [left] at (0, \HolderMoritz) {$\hspace{-30pt}\frac{n}{p_-^{L^\sharp}} \!-\!n$};
	\draw [fill=black] (0,  \HolderMoritz) circle [radius = .5pt];	
	\node [left] at (0, 0.3) {$\theta$};
	\draw [fill=black] (0, 0.3) circle [radius = .5pt];

	\draw [thin,dotted] (1.2,2) -- (1.2,0); 
	\draw [thin,dotted] (0.3,2) -- (0.3,0); 
	\draw [thin,dotted] (1+\UpGradient,2) -- (1+\UpGradient,0);
	\draw [thin,dotted] (1+\LowSob,2) -- (1+\LowSob,0);
	\draw [thin,dotted] (1+\xlength,2) -- (1+\xlength,0); 
	
	\draw [thin,dotted] (0.3, \HolderMoritz) -- (1.2,\HolderMoritz);
	\draw [thin,dotted] (0.3, 0.3) -- (1.2, 0.3);
	
	\path [fill=red!80!black, opacity = 0.4](HolderAlex) -- (DirHolderBesov) -- (P00);
	\path [fill=lightgray, opacity = 0.6] (P00)--(P01)--(P11)--(P10)--(DirExtra)--(HolderAlex);
	\path [fill=blue!80!black, opacity = 0.4] (P11)--(P12)--(RegExtra)--(P01)--(P11);
	\draw [ultra thick, gray] (P00) -- (P01); 
	\draw [ultra thick, blue!80!black] (P01) -- (RegExtra);
	\draw [ultra thick, gray] (DirExtra) -- (P11); 
	\draw [ultra thick, blue!80!black] (P11) -- (P12); 
	\draw [ultra thick, red!80!black] (DirHolder) -- (P10);
	\draw [ultra thick, gray] (DirExtra) -- (P10);
	\draw [ultra thick, red!80!black] (DirExtra) -- (DirHolderBesov); 
	\draw [ultra thick, gray] (DirExtra) -- (HolderAlex); 
	\draw [fill=white] (P00) circle [radius = .75pt];
	\draw [fill=blue!80!black] (P01) circle [radius = .75pt];
	\draw [fill=red!80!black] (P10) circle [radius = .75pt];
	\draw [fill=blue!80!black] (P11) circle [radius = .75pt];	
	\draw [fill=white] (P12) circle [radius = .75pt];
	\draw [fill=gray] (DirExtra) circle [radius = .75pt];
	\draw [fill=red!80!black] (HolderMoritz) circle [radius = 0.75pt];
	\draw [fill=red!80!black] (HolderAlex) circle [radius = 0.75pt];
	\draw [fill=white] (RegExtra) circle [radius = .75pt];
	\draw [fill=white] (DirHolder) circle [radius = .75pt];
	\draw [fill=white] (DirHolderBesov) circle [radius = .75pt];
\end{tikzpicture}
\end{center}
\caption{Compatible well-posedness region for Besov and Hardy--Sobolev data when $p_-(L^\sharp) < 1$. This implies $p_+(L) = \infty$. Hence, at the bottom there is no horizontal thick red line as in Figure~\ref{fig: diagram p+large} for $\nicefrac{1}{p} > 0$ and Theorem~\ref{thm: Holder-dir} yields compatible well-posedness for the full horizontal segment with $\nicefrac{1}{p} \leq 0$. The number $\theta$ comes from the first-order approach in \cite{AA}. It has a specific meaning, see Proposition~\ref{prop: blue region large}, and is not larger than $n(\nicefrac{1}{p_-(L^\sharp)} -1)$, which is the limitation of part (iii) in Theorem~\ref{thm: Holder-dir} for $\Lamdot^\alpha$-data. When $1-\nicefrac{1}{p_-(L^\sharp)} < -\nicefrac{\alpha}{n} \leq 0$, well-posedness of the $\Lamdot^\alpha$ Dirichlet problem for $\cL$ can also be obtained by duality from well-posedness of the $\H^{\nicefrac{n}{(\alpha+n)}}$ regularity problem for $\cL^*$, using $DB^*$-adapted spaces~\cite[Thm.~7.11]{AA}. This is why the corresponding horizontal segment for $\nicefrac{1}{p} = \nicefrac{-\alpha}{n} \leq 0$ has been colored in gray.}
\label{fig: diagram p-small}
\end{figure}

\subsection{Neumann problems}
\label{intro:Neumann}

\noindent Although this is not central to our monograph, we complete the discussion with results on the Neumann problem. For $1_* < p < \infty$, the Neumann problem with data $g\in \H^p(\R^n; \IC^m)$\index{Neumann problem} consists in solving (modulo constants)
\begin{equation*}
(N)_{p}^\Le  \quad\quad
\begin{cases}
\Le u=0   & \text{in }\reu, \\
\NT (\nabla u)\in \L^p(\R^n),   \\
\lim_{t \to 0} a \partial_t u(t, \cdot) = g & (\text{in } \cD'(\R^n; \IC^m)).\end{cases}
\end{equation*}
Note that due to the block structure $a \partial_t u$ is indeed the conormal derivative $\dnuA u = e_0 \cdot A \nabla u$. Here, constants are solutions which do not change the Neumann data, so we must argue modulo constants. Once again there is a construction of energy solutions via the Lax--Milgram lemma, using the data space $\Hdot^{-\nicefrac{1}{2},2}$ for $g$.\footnote{We recall the construction in Section~\ref{sec: Neumann problems}.}

In order to understand how our results help in deducing a range of exponents for which the Neumann problem is compatibly well-posed from existing literature, we recall the first-order approach. For block systems it simply begins by writing \eqref{eq:block} in the equivalent form
\begin{align}
\label{eq: first-order reduction}
\partial_t \begin{bmatrix} a \partial_t u \\ \nabla_x u \end{bmatrix}
+ \begin{bmatrix} 0& \div_x \\ -\nabla_x & 0 \end{bmatrix} 
\begin{bmatrix} a^{-1} & 0  \\ 0 & d \end{bmatrix}
\begin{bmatrix} a \partial_t u \\ \nabla_x u \end{bmatrix}
= \begin{bmatrix}0 \\ 0 \end{bmatrix},
\end{align}
where the second line is a dummy equation, or in short notation
\begin{align}
\label{eq: DB equation intro}
\partial_t F + DB F = 0,
\end{align}
where $F = \nabla_A u \coloneqq [a \partial_t u, \nabla_x u]^\top$ is the \emph{conormal gradient}\index{conormal gradient} and $DB$ is called \emph{perturbed Dirac operator}. This operator is bisectorial and there are associated abstract Hardy spaces $\IH_{DB}^p$. The idea then is to work backwards from that: first classify all weak solutions to \eqref{eq: DB equation intro} in the usual classes and then try to reconstruct $u$ from its conormal gradient.

The principal thesis in the work of the first author with Stahlhut~\cite{AusSta} and Mourgoglou~\cite{AM} is that there is an open interval $I_L \subseteq (1_*, \infty)$ such that if $p \in I_L$, then
\begin{itemize}
	\item the conormal gradient of every weak solution to $(N)_{p}^\Le$ has an \emph{a priori} representation via the semigroup associated with $|DB| \coloneqq ((DB)^2)^{1/2}$,
	\item compatible well-posedness follows if the estimate $\|\NT(\nabla u)\|_p \lesssim \|g\|_{\H^p}$ holds for all $g\in \H^p(\R^n; \IC^m) \cap \H^{-\nicefrac{1}{2},2}(\R^n; \IC^m)$ and $u$ the energy solution with Neumann data $g$.\footnote{This is Theorem~1.8 in \cite{AM}.}
\end{itemize}
The interval $I_L$ corresponds to identification $\IH_{DB}^p = \IH_D^p$ of abstract and concrete Hardy spaces up to equivalent $p$-quasinorms and a certain $\L^p$-coercivity assumption of $B$ when $p>2$.

In the block case one can produce a formal solution to the Neumann problem by $u(t,x)\coloneqq - L^{-1/2} \e^{-t L^{1/2}} (a^{-1} g)(x)$, so that once this is made rigorous, compatible well-posedness of $(N)_p^{\Le}$ follows in the range $p \in I_L$. This being said, our main contribution for the Neumann problem lies in proving the equality\footnote{The proof is in Section~\ref{sec: Comparison with Auscher-Stahlhut interval}, Corollary~\ref{cor: AE interval and AM interval} and the principal issue is to \emph{prove} the $\L^p$-coercivity for $p>2$. Before it was only known that when $a=1$, $I_L$ cannot be larger than $(q_-(L),q_+(L))$ and that its upper endpoint is $q_+(L)$ if in addition $d$ is strictly elliptic, see \cite[Sec.~12.4.1]{AusSta}.} 
\begin{align}
\label{eq: IL identification intro}
I_L = (q_-(L), q_+(L))
\end{align}
and then we conclude the following result.

\begin{thm}[Neumann problem\index{Neumann problem!main result}]
	\label{thm: blockneu}
	Let  $q_-(L) < p < q_+(L)$. Then the Neumann problem $(N)_p^{\Le}$ is compatibly well-posed (modulo constants).
\end{thm}

With the determination of $I_L$ at hand, one can write down all further implications from \cite{AM} for solutions with the \emph{a priori} representation of $\nabla_A u$. This would lead us too far from the objective of our monograph. Let us just mention that there are additional regularity properties for solutions to $(N)_p^{\Le}$ in Theorem~\ref{thm: blockneu}, similar to Theorem~\ref{thm: blockreg}, and that well-posedness of an adjoint `rough' Neumann problem follows by duality.\footnote{For further regularity in the Neumann problem, see Corollary~1.2 in \cite{AM}. Therein, the Dirichlet data is given by $f=-{L^{-1/2}} (a^{-1}g)\in \Hdot^{1,p}$ using a suitable extension of the square root. For the duality with the rough Neumann problem see Theorem~1.6 and then Theorem~1.3 and Corollary~1.4 in \cite{AM} for the \emph{a priori} representation and regularity for its solutions.}
Finally, in the spirit of Section~\ref{intro:fractionalDP}, there are fractional Neumann problems in between for which ranges of compatible well-posedness have also been described via $I_L$.\footnote{See \cite{AA} for an introduction to and results on these problems.} In fact, this is the gray region in the diagrams above.
\subsection{Synthesis}
\label{subsec: Synthesis}
We close the introduction with a comment further explaining the color code in the diagrams in Section~\ref{intro:fractionalDP}.  Heuristically, the $\H^p$-theory for $DB$ comprises the theory for $L$ at both smoothness scales $s=0$ and $s=1$. On the level of Hardy spaces, this becomes apparent in the fact that the interval in \eqref{eq: IL identification intro} is the intersection of intervals of identification for $\IH_L^p$ and $\IH_L^{1,p}$, compare with \eqref{eq: HpL identification intro} and \eqref{eq: H1pL identification intro}. On the level of boundary value problems, the first-order approach via $DB$ yields ranges of exponents in which problems with Neumann and Dirichlet data are simultaneously well-posed --- this is the gray region. The $L$-adapted theory allows us to separate issues and obtain significantly larger ranges for the problems with Dirichlet data --- gray and blue regions. Finally, there is a new phenomenon ---  solving Dirichlet problems for one Sobolev conjugate above the limitation of the Hardy space theory in the red region.
\subsection{Notation}
\label{subsec: Notation}

\noindent The following notation will be used throughout the monograph.

\medskip

\noindent \emph{Geometry and measure.} We let $B(x,r) \subseteq \R^n$ the open ball of radius $r> 0$ around $x \in \R^n$. Given a ball $B \subseteq \R^n$ of radius $r(B)$, we write $cB$ for the concentric ball of radius $cr(B)$ and define the annular regions $C_j(B)$, $j \in \IN$, by
\begin{align*}
 C_1(B) \coloneqq 4B, \quad C_j(B) \coloneqq 2^{j+1} B \setminus 2^j B \quad (j \geq 2).
\end{align*}
The same type of notation will be used for cubes instead of balls. In this case, $\ell(Q)$ denotes the sidelength of $Q$. In order to avoid even the slightest confusion, let us explicitly state that for us $\IN \coloneqq \{1,2,3,\ldots\}$.

We write the Euclidean distance on finite-dimensional vector spaces as $\dist(x,y) \coloneqq |x-y|$ and extend the notation to sets $E, F \subseteq \R^n$ via 
\begin{align*}
\dist(E,F) \coloneqq \inf \{\dist(x,y) : x \in E, \, y \in F \}.
\end{align*}
The characteristic function of $E$ is $\ind_E$\index{$\ind_E$ (characteristic function)}. In $\ree$ we denote points by $(t,x)$ and define the open upper halfspace 
\begin{align*}
\reu \coloneqq \{(t,x) : t>0, \, x \in \R^n \}.
\end{align*}
We write $|\cdot|$ for the Lebesgue measure if the underlying Euclidean space is clear form the context. For integral averages we use $\barint$ and $\bariint$ in $\R^n$ and $\ree$, respectively, as well as the notation $(f)_E \coloneqq \barint_E f$. We use the (uncentered) Hardy--Littlewood maximal operator\index{maximal operator ($\Max$)} defined for measurable functions on $\R^n$ via
\begin{align*}
\Max(f)(x) \coloneqq \sup_{B \ni x} \barint_B |f| \, \d y \quad (x \in \R^n),
\end{align*}
where the supremum runs over all balls $B$ that contain $x$. Occasionally, we also use cubes instead of balls.
\medskip

\noindent \emph{Gradient and divergence of vector-valued functions}. Partial derivatives of $\IC^m$-valued functions are taken componentwise. If $f$ is a $\IC^m$-valued function on a subset of $\R^n$ or $\ree$, then 
\begin{align*}
\nabla_x f = [\partial_{x_1} f, \ldots, \partial_{x_n}f]^\top
\end{align*}
is a function valued in $\IC^{mn} \cong (\IC^m)^n$. In the opposite direction, if $F=[F_1,\ldots,F_n]^\top$ is $\IC^{mn}$-valued, then we let 
\begin{align*}
\div_x F = \partial_{x_1} F_1 + \ldots +\partial_{x_n} F_n.
\end{align*}
Gradient and divergence with respect to all variables in $\ree$ are defined as $\nabla f = [\partial_t f, \nabla_xf]^\top$ and $\div = \partial_t F_\no + \div_x F_\ta$ if $F = [F_\no, F_\ta]^\top$ is valued in $\IC^m \times \IC^{mn}$.

\medskip

\noindent \emph{Exponents}. We let
\begin{alignat*}{2}
\frac{1}{p'} &= 1 - \frac{1}{p} &&\qquad (p \in [1,\infty], \text{H\"older conjugate}\index{H\"older conjugate}), \\
\frac{1}{p_*} &= \frac{1}{p} + \frac{1}{n} &&\qquad (p \in (0,\infty], \text{lower Sobolev conjugate}\index{Sobolev conjugate!upper}), \\
\frac{1}{p^*} &= \frac{1}{p} - \frac{1}{n} &&\qquad (p \in (0,n), \text{upper Sobolev conjugate}\index{Sobolev conjugate!lower}), \\
\frac{1}{[p_0,p_1]_\theta} &= \frac{1-\theta}{p_0} + \frac{\theta}{p_1} &&\qquad (p_i \in (0,\infty], \, \theta \in [0,1], \text{interpolating index}\index{interpolating index}).
\end{alignat*}
The underlying dimension for Sobolev conjugates is usually $n$ and will always be clear from the context. We also agree on $p^* \coloneqq \infty$ for $p \geq n$.
\medskip

\noindent \emph{Constants}. Given $a,b \in [0,\infty]$, we write $a \lesssim b$ to mean $a \leq Cb$ for some $C \in (0,\infty)$ (oftentimes called `implicit constant'\index{implicit constant}) that is independent of $a$ and $b$. We write $a \simeq b$ to mean $a \lesssim b$ and $b \lesssim a$. In this case one of $a,b$ is equal to $\infty$ (or $0$) precisely when both are. Unless stated otherwise, estimates in this monograph are \emph{quantitative} in the sense that constants in estimates depend only on constants quantified in the relevant hypotheses. Such dependence will usually be clear.\index{quantitative estimates}

\medskip

\noindent \emph{Index}. This monograph has an index. For the sake of readability we shall occasionally refer to results by their name listed in the index instead of a number in the text.
\section{Preliminaries on function spaces}
\label{sec: function spaces}

\noindent This chapter contains all necessary background on function spaces that will be used later on. Throughout, we consider $\IC^k$-valued functions for some fixed $k \in \IN$. For simplicity we often drop the dependence of $k$ in the notation and write $\L^2(\R^n) = \L^2(\R^n; \IC^k)$, and so on. On $\R^n$ we abbreviate further $\L^2 = \L^2(\R^n)$. Concerning the dilemma that parts of the literature only treat scalar-valued functions, we agree on using such results for $k>1$ without further notice in the following cases:
\begin{itemize}
	\item splitting into components is immediately clear from the definition (e.g.\ $\L^2(\R^n; \IC^k) \cong \bigotimes_{j=1}^k \L^2(\R^n; \IC)$),
	
	\item proofs are exactly the same except for a systematic replacement of absolute values by Euclidean norms (e.g.\ Calder\'on--Zygmund decompositions or atomic decompositions).
\end{itemize}
\subsection{Lebesgue spaces and distributions}
\label{subsec: Lebesgue spaces}
On a (Lebesgue) measurable set $E \subseteq \R^n$ we let $\L^p(E)$, $p \in (0,\infty]$, be the (quasi-)Banach space of function classes with finite (quasi)norm
\begin{align*}
\|f\|_{\L^p(E)} \coloneqq \bigg(\int_E |f|^p \, \d x \bigg)^{\frac{1}{p}}.
\end{align*}
The right-hand side is interpreted as the essential supremum when $p = \infty$. We abbreviate $\|\cdot\|_p \coloneqq \|\cdot\|_{\L^p(\R^n)}$. The classes of functions that are $p$-integrable on compact subsets of $E$ are denoted by $\Lloc^p(E)$ and carry the natural Fréchet topology.

We write $\C_0^\infty(O)$, where $O \subseteq \R^n$ is open, and $\cS(\R^n)$ for the test functions with compact support and of Schwartz-type, respectively. Their topological duals are the distribution spaces $\cD'(O)$ and $\cS'(\R^n)$. The subspace $\cZ(\R^n) \subseteq \cS(\R^n)$ is the space of Schwartz functions $f$ whose Fourier transform\index{F@$\cF$ (Fourier transform)} $\cF f$ satisfies $D^\alpha \cF f(0) = 0$ for all multi-indices $\alpha \in \IN_0^n$. The dual $\cZ'(\R^n)$ can be identified with the quotient $\cS'(\R^n)/\cP(\R^n)$, where $\cP(\R^n)$ is the space of polynomials on $\R^n$, see~\cite[Sec.~5.2.1]{Triebel-TheoryOf}. 

For $p \in [1,\infty]$ the Sobolev spaces $\W^{1,p}(O)$\index{Sobolev space} is the collection of those $f \in \L^p(O)$ that satisfy $\nabla_x f \in \L^p(O)$ in the sense of distributions. Again, there are local versions denoted by $\Wloc^{1,p}(O)$.
\subsection{Tent spaces}
\label{subsec: tent spaces}

\noindent Tent spaces have been introduced by Coifman--Meyer--Stein in~\cite{CMS}. Good sources for detailed proofs are~\cite{Amenta-Tent, Amenta-Interpol}.

For $x \in \R^n$ we introduce the cone with vertex $x$,
\begin{align*}
\Gamma (x) \coloneqq \{(s,y) \in \reu : |x-y| < s \},
\end{align*}
and define the corresponding (conical) \emph{square function}\index{square function!$S$} for measurable functions $F : \reu \to \IC^k$ by
\begin{align}
\label{eq: S}
(S F)(x) \coloneqq \bigg(\iint_{\Gamma(x)} |F(s,y)|^2 \, \frac{\d s \d y}{s^{1+n}} \bigg)^{\frac{1}{2}} \quad (x \in \R^n).
\end{align}
For $\alpha \geq 0$ the Carleson functional\index{Carleson functional} is defined as
\begin{align}
\label{eq: C}
C_\alpha F(x) \coloneqq \sup_{r>0} \frac{1}{r^\alpha} \bigg(\frac{1}{r^n} \int_0^{r} \int_{B(x,r)} |F(t,y)|^2 \frac{\d y \d t}{t} \bigg)^{\frac{1}{2}}.
\end{align}
With a slight abuse of notation, we denote by $t^{-s} F$ the function $(t,y) \mapsto t^{-s} F(t,y)$.

\begin{defn}
\label{def: tent spaces}
Let $s \in \R$, $\alpha \geq 0$ and $p \in (0, \infty]$. For finite $p$ the \emph{tent space}\index{tent space} $\T^{s,p}$ consists of all functions $F \in \Lloc^2(\reu)$ with finite quasi-norm
\begin{align*}
\|F\|_{\T^{s,p}} \coloneqq \|S(t^{-s} F)\|_p.
\end{align*}
For $p=\infty$ the tent space $\T^{s,\infty;\alpha}$ consists of all functions $F \in \Lloc^2(\reu)$ with finite norm
\begin{align*}
\|F\|_{\T^{s,\infty; \alpha}} \coloneqq \|C_\alpha(t^{-s} F)\|_\infty.
\end{align*}
\end{defn}

\begin{rem}
\label{rem: tent spaces}
For brevity we set $\T^p \coloneqq \T^{0,p}$ for finite $p$ and we abbreviate and $\T^{s, \infty} = \T^{s,\infty;0}$ with the special case $\T^\infty \coloneqq \T^{0,\infty;0}$. We also note that $F \mapsto t^s F$ is an isometric isomorphism from $\T^p$ onto $\T^{s,p}$ and from $\T^{0,\infty;\alpha}$ onto $\T^{s,\infty;\alpha}$.
\end{rem}

All tent spaces are quasi-Banach spaces (Banach when $p \geq 1$) and their topology is finer than the one on $\Lloc^2(\reu)$. Both statements follow directly from the bounds 
\begin{align*}
S F(x) &\geq (2t)^{-\frac{1+n}{2}}\|F\|_{\L^2((t,2t) \times B(x,t))}, \\
C_\alpha F(x) &\geq t^{-\frac{1+n}{2}-\alpha}\|F\|_{\L^2((\frac{t}{2},t) \times B(x,t))}, 
\end{align*}
for $t>0$ and  $x \in \R^n$ and Fatou's lemma. Moreover, for $p<\infty$ there is a \emph{universal approximation technique}\index{universal approximation technique!for tent spaces} by functions in $\L^2(\reu)$ with compact support~\cite[Prop.~1.4]{Amenta-Interpol}: 
\begin{align*}
\forall F \in \T^{s,p}: \quad \lim_{j \to \infty} \ind_{(j^{-1},j) \times B(0,j)} F = F \quad (\text{in } \T^{s,p}).
\end{align*}
`Universal' refers to the fact that the same approximating sequence can be used in all tent spaces that $F$ belongs to. Results of this type will be important for us since we shall often work with intersections of spaces. We could also change the cones $\Gamma(x)$ to 
\begin{align*}
\Gamma_\alpha (x) \coloneqq \{(s,y) \in \reu : |x-y| < \alpha s \}
\end{align*}
for any fixed $\alpha>0$. This \emph{change of angle}\index{change of angle/aperture} yields equivalent tent space norms~\cite[Prop.~4]{CMS}. 

If $p \in (0,\infty)$, then the (anti-)dual space of $\T^{s,p}$ can be identified through the $\L^2$ duality pairing\index{duality!for $\T^{s,p}$}
\begin{align}
\label{eq: duality paring T}
\langle F, G \rangle = \iint_{\reu} F(s,y) \cdot \cl{G(s,y)} \, \frac{\d s \d y}{s},
\end{align}
see \cite[Prop.~1.9 \& Thm.~1.11]{Amenta-Interpol}. We have
\begin{align*}
	(\T^{s,p})^*
	= \begin{cases}
		\T^{-s,p'} & \text{ if } p>1 \\
		 \T^{-s,\infty; n(\frac{1}{p}-1)} & \text{ if } p \leq 1
	\end{cases}.
\end{align*}
In particular, $\T^2 = \L^2(\reu, \frac{\d t \d x}{t})$ with equivalent norms, which can also be seen directly by Fubini's theorem:
\begin{align*}
\int_{\R^n} \iint_{|x-y| < s} |F(s,y)|^2 \, \frac{\d s \d y}{s^{1+n}} \d x
= \omega_n \iint_{\reu} |F(s,y)|^2 \, \frac{\d s \d y}{s},
\end{align*}
where $\omega_n$ is the measure of the unit ball in $\R^n$. This technique is called \emph{averaging trick}\index{averaging trick} in the following.

We shall need one more tent space that is related to the (modified) \emph{non-tangential maximal function}\index{non-tangential maximal function} 
\begin{align}
\label{eq: NT}
\NT F(x) 
\coloneqq \sup_{t>0} \bigg(\bariint_{W(t,x)} |F(s,y)|^2 \, \d s \d y \bigg)^{\frac{1}{2}},
\end{align}
where $x \in \R^n$ and $W(t,x) \coloneqq (\nicefrac{t}{2}, 2t) \times B(x,t)$ is called \emph{Whitney box}\index{Whitney box}. 

\begin{defn}
\label{def: C and NT spaces}
Let $p\in (0,\infty)$. The tent spaces $\T^{0,p}_\infty$\index{tent space} consists of all functions $F \in \Lloc^2(\reu)$ with finite (quasi-)norm
\begin{align*}
\|F\|_{\T^{0,p}_\infty} \coloneqq \|\NT F\|_p.
\end{align*}
\end{defn}

As before, these are quasi-Banach (Banach when $p\geq 1$) spaces with a topology that is stronger than $\Lloc^2(\reu)$. Moreover, a \emph{change of Whitney parameters}\index{change of Whitney parameters!for $\NT$} to $W(t,x) = (c_0^{-1}t,c_0t) \times B(x,c_1t)$ with $c_0 > 1$ and $c_1 > 0$ leads to an equivalent $\T^{0,p}_\infty$-norm. For the reader's convenience we reprove this fact in Appendix~\ref{sec:technical} together with further auxiliary properties of non-tangential maximal functions.
\subsection{\texorpdfstring{$\boldsymbol{\Z}$}{Z}-spaces}
\label{subsec: Z spaces}

In the context of boundary value problems these spaces emerged from the work of Barton--Mayboroda~\cite{BM} under a different name. Their relation to tent spaces has been noted by Amenta~\cite{Amenta-Interpol}. 

For measurable functions $F$ on  $\reu$ we introduce the \emph{Whitney average functional}\index{Whitney average functional!usual ($W$)} 
\begin{align*}
	W(F)(t,x) = \bigg(\bariint_{W(t,x)} |F(s,y)|^2 \, \d s \d y \bigg)^{\frac{1}{2}} \quad ((t,x) \in \reu).
\end{align*}

\begin{defn}
\label{def: Z spaces}
Let $s \in \R$ and $p \in (0, \infty]$. The \emph{$\Z$-space}\index{Z@$\Z$-space} $\Z^{s,p}$ consists of all functions $F \in \Lloc^2(\reu)$ with finite quasi-norm
\begin{align*}
	\|F\|_{\Z^{s,p}} \coloneqq \|W(t^{-s} F)\|_{\L^p(\reu, \frac{\d t \d x}{t})}.
\end{align*}
\end{defn}

All $\Z$-spaces are quasi-Banach spaces (Banach when $p \geq 1$), their topology is finer than the one on $\Lloc^2(\reu)$ and for $p<\infty$ they have the same \emph{universal approximation technique}\index{universal approximation technique!for $\Z$-spaces} as the tent spaces. This can simply be checked by hand or deduced by real interpolation since $\Z$-spaces are the real interpolants of tent spaces, see Section~\ref{subsec: complex interpolation} below. Many properties of tent spaces have a $\Z$-space analog: A change of Whitney parameters\index{change of Whitney parameters!for $\Z$-spaces}  leads to equivalent quasi-norms (Remark~\ref{rem: independence of Whitney parameters}), the averaging trick reveals $\Z^{0,2} = \L^2(\reu, \frac{\d t \d x}{t}) = \T^{0,2}$ and the $\L^2$ duality pairing \eqref{eq: duality paring T} gives rise to \index{duality!for $\Z^{s,p}$}
\begin{align*}
	(\Z^{s,p})^*
	= \begin{cases}
		\Z^{-s,p'} & \text{ if } p>1 \\
		\Z^{-s + n(\frac{1}{p}-1), \infty} & \text{ if } p \leq 1
	\end{cases},
\end{align*}
see \cite[Prop.~2.22 \& Thm.~2.28]{AA}.
\subsection{Hardy spaces}
\label{subsec: Hardy spaces}

For $1<p<\infty$ we set $\H^p \coloneqq \L^p$ and for $p \leq 1$ we denote by $\H^p$ the real Hardy space\index{Hardy space} of Fefferman--Stein~\cite{Fefferman-Stein, Stein93}. For $p=1$ we have the continuous inclusion $\H^1 \subseteq \L^1$. 

We shall exclusively work in the range $p > 1_*$ and for most of our applications it will be convenient to think of $\H^p$-spaces in terms of atoms.

\begin{defn}
	\label{def: Hp atom}
	Let $p \in (1_*,1]$ and $q \in (1,\infty]$. An \emph{$\L^q$-atom for $\H^p$}\index{atom!for $\H^p$} is a function $a$ supported in a cube $Q \subseteq \R^n$ such that $\|a\|_q \leq \ell(Q)^{\frac{n}{q}-\frac{n}{p}}$ and $\int_{\R^n} a \d x = 0$.
\end{defn}

Of course we could also use balls instead of cubes in the definition. The \emph{atomic decomposition}\index{atomic decomposition!for $\H^p$}~\cite[Sec.~III.3.2]{Stein93} states that every $f \in \H^p$ can be written as $f = \sum_{i=1}^\infty \lambda_i a_i$, where the sum converges unconditionally in $\H^p$, the $a_i$ are $\L^\infty$-atoms for $\H^p$ and the scalars $\lambda_i$ satisfy
\begin{align}
	\label{eq: good atomic decomposition}
	\|(\lambda_i)\|_{\ell^p} \lesssim \|f\|_{\H^p}.
\end{align}
Moreover 
\begin{align*}
	\|f\|_{\H^p} \simeq \inf_{f = \sum_{i=1}^\infty \lambda_i a_i} \|(\lambda_i)\|_{\ell^p}.
\end{align*}
When working with operators that are defined on some space $\L^s$, $s \in (1,\infty)$, but not on distributions, the following compatibility property will be important: If $f \in \H^p \cap \L^s$, then the series that realizes \eqref{eq: good atomic decomposition} can be taken such that it also converges in $\L^s$. In fact, the explicit construction in \cite{Stein93} has this property, as has been carefully verified in \cite{Rocha}.

Occasionally, we shall need that for $p \in (1_*,1]$ smooth functions with compact support and integral zero are dense in $\H^p$. This follows, for example, by mollification of $\L^\infty$-atoms for $\H^p$ with a smooth kernel~\cite[Thm.~3.33]{Folland-Stein}.
\subsection{Homogeneous smoothness spaces}
\label{subsec: homogeneous smoothness spaces}

\noindent Good textbooks for further background are \cite{Triebel-TheoryOf, Sawano, Peetre-Besov, Grafakos2}. An operator-theoretic perspective on these spaces will emerge later on in Section~\ref{subsec: D-adapted spaces}. All function spaces will be on $\R^n$ and for the sake of a clear exposition we omit this from our notation.

Let $\psi \in \C_0^\infty$ be supported in the annulus $\frac{1}{2} \leq |\xi| \leq 2$ and normalized to
\begin{align*}
\sum_{j \in \IZ} \psi(2^j \xi) = 1 \quad (\xi \in \R^n \setminus \{0\})
\end{align*}
and introduce for $j \in \IZ$ the associated Littlewood--Paley operators\index{Littlewood--Paley!operator} $\Delta_j f \coloneqq \cF^{-1}(\psi(2^j \cdot) \cF f)$. Here $\cF$ denotes the Fourier transform on $\R^n$. Whenever $f \in \cZ'$, then
\begin{align}
\label{eq: L-P series in Z'}
\sum_{j \in \IZ} \Delta_j f = f \quad (\text{in } \cZ'),
\end{align}
see~\cite[Prop.~2.11]{Sawano}. The Paley--Wiener--Schwartz theorem~\cite[Thm.~1.7.7]{Hoermander}\index{Theorem!Paley--Wiener--Schwartz} asserts that the packets $\Delta_j f$ are smooth functions of moderate growth and the general idea behind the following homogeneous smoothness spaces is to measure them in Lebesgue-type norms.

\begin{defn}
\label{def: H-zero-dot}
Let $s \in \R$ and $p \in (0, \infty]$. The homogeneous \emph{Hardy--Sobolev space}\index{Hardy--Sobolev space!$\Hdot^{s,p}$} $\Hdot^{s,p}$ when $p<\infty$ is the set of those $f \in \cZ'$ with finite (quasi)norm
\begin{align*}
\|f\|_{\Hdot^{s,p}} &\coloneqq \| \|j \mapsto 2^{js} \Delta_j f(\cdot) \|_{\ell^2(\IZ)} \|_p.
\intertext{The endpoint space $\Hdot^{s,\infty}$ is determined by the norm}
\|f\|_{\Hdot^{s,\infty}} &\coloneqq \inf_{f = \sum_{j \in \IZ} \Delta_j f_j} \| \|j \mapsto 2^{js} f_j(\cdot) \|_{\ell^2(\IZ)} \|_\infty.
\intertext{The homogeneous \emph{Besov space}\index{Besov space!$\Bdot^{s,p}$} $\Bdot^{s,p}$ is the set of those $f \in \cZ'$ with finite (quasi)norm}
\|f\|_{\Bdot^{s,p}} &\coloneqq \| \|j \mapsto 2^{js} \|\Delta_j f\|_p \|_{\ell^p(\IZ)}.
\end{align*}
\end{defn}

\begin{rem}
\label{rem: H-zero-dot}
Within the full scale of homogeneous Besov--Triebel--Lizorkin spaces the common notation for $\Hdot^{s,p}$ and $\Bdot^{s,p}$ is $\Fdot^s_{p,2}$ and $\Bdot^s_{p,p}$, respectively.
\end{rem}

In the following let $\X$\index{X@$\X$ (one of $\B, \H$)} denote either $\B$ or $\H$. Then $\Xdot^{s,p}$ is a quasi-Banach space (Banach when $p \geq 1$), different choices of $\psi$ lead to equivalent (quasi)norms and there are continuous inclusions
\begin{align*}
	\cZ \subseteq \Xdot^{s,p} \subseteq \cZ'.
\end{align*}
Moreover, $\cZ$ is dense in $\Xdot^{s,p}$ when $p<\infty$ via a \emph{universal approximation technique}~\cite[Sec.~5.1.5]{Triebel-TheoryOf}\index{universal approximation technique!for homogeneous smoothness spaces}: If $\varphi \in \cS$ is such that $\varphi(0) =1$ and $\cF \varphi$ is supported in $|\xi| \leq 1$, then
\begin{align*}
	\forall f \in \Xdot^{s,p}: \quad \lim_{N \to \infty} \lim_{\delta \to 0} \bigg(\varphi(\delta \, \cdot) \sum_{|j| \leq N} \Delta_j f\bigg) = f \quad (\text{in } \Xdot^{s,p}).
\end{align*}
`Universal' has the same meaning and purpose as for the tent spaces and the approximants are in $\cZ$ provided that $\delta < 2^{-N-1}$.

While the ambient space $\cZ'$ is well-suited for general considerations, applications to boundary value problems require more concrete `realizations' of $\Xdot^{s,p}$. This issue can be resolved thanks to an observation due to Peetre~\cite[pp.~52-56]{Peetre-Besov}, see also \cite[Sec.~2.4.3]{Sawano}. Suppose that $L \in \IN_0$ is such that $L > s- \nicefrac{n}{p}$ and let $\cP_{L-1}$ be the space of polynomials of degree at most $L-1$. Then for any $f \in \Xdot^{s,p}$ the series in \eqref{eq: L-P series in Z'} converges in $\cS'/\cP_{L-1}$ and identifying $f$ with the limit yields an isometric copy of $\Xdot^{s,p}$ that is continuously embedded into the ambient space $\cS' /\cP_{L-1}$. In particular, the spaces of smoothness $s < 1$ can be viewed as subspaces of $\cS' / \IC$ and even of $\cS'$ if $s \leq 0$ and $p<\infty$. 

Within these smaller ambient spaces, $\Xdot^{s,p}$ can often be given an equivalent and more familiar quasinorm that does not make sense modulo \emph{all} polynomials. For example, we have the Littlewood--Paley theorem
\begin{alignat*}{2}
	\Hdot^{0,p} &= \H^p = \L^p &&\quad (1<p<\infty),\\
	\Hdot^{0,p} &= \H^p &&\quad (p \leq 1),
\end{alignat*} 
see \cite[Sec~6.2 \& 6.4]{Grafakos2} and in accordance with the observation above $\L^p$ and $\H^p$ do not contain any polynomials besides $0$. For $p=\infty$ we have
\begin{align}
\label{eq: BMO identification}
	\Hdot^{0,\infty} = \BMO \eqqcolon \Lamdot^0,
\end{align}
see \cite[Sec.~5.2.4]{Triebel-TheoryOf} and references therein. Here, $\BMO$ is the John-Nirenberg space of functions modulo constants with bounded mean oscillation\index{bounded mean oscillation ($\BMO$)}
\begin{align*}
	\|f\|_{\BMO} \coloneqq \sup_{B} \barint_B |f(x)-(f)_B| \, \d x,
\end{align*}
where the supremum is taken over all balls in $\R^n$. For $0<s<1$ we denote by $\Lamdot^s$ the H\"older space\index{H\"older space ($\Lamdot^\alpha$)} of functions modulo constants with finite norm
\begin{align*}
	\|f\|_{\Lamdot^s} \coloneqq \sup_{x\neq y} \frac{|f(x)-f(y)|}{|x-y|^s},
\end{align*}
which can be identified with
\begin{align*}
 \Lamdot^s = \Bdot^{s,\infty}  \quad (0<s<1),
\end{align*}
see \cite[Thm.~5.2.3.2]{Triebel-TheoryOf}. 

Next, we recall relevant duality results in the case of finite exponents $p$. Since in this case $\cZ$ is dense in $\Xdot^{s,p}$, we can view the (anti-)dual space $(\Xdot^{s,p})^*$ as a subspace of $\cZ'$ by restricting functionals to $\cZ$. In this sense we have\index{duality!for $\Xdot^{s,p}$}
\begin{align*}
	(\Xdot^{s,p})^* = \Xdot^{-s,p'} \quad (1 \leq p < \infty).
\end{align*}
A direct proof for inhomogeneous spaces that applies \emph{mutatis mutandis} in our homogeneous setting is given in~\cite[Sec.~2.11.2]{Triebel-TheoryOf}, see also \cite[Sec.~5.2.5]{Triebel-TheoryOf}. For $s=0$ and $p=1$ this is the famous $\H^1-\BMO$ duality of Fefferman--Stein~\cite{Fefferman-Stein}. In the case $p<1$ we shall only need the duality
\begin{align}
\label{eq: duality Hardy-Holder}
	(\Hdot^{0,p})^* &= \Lamdot^{n(\frac{1}{p} - 1)} \quad(1_* < p \leq 1),
\end{align}
see~\cite[Thm.~4.2]{Jawerth77} or again \cite[Sec.~2.11.2]{Triebel-TheoryOf}. An alternative proof is given in \cite[Rem.~5.14]{Jawerth-Frazier}.

Spaces for different smoothness parameters are related via a \emph{lifting property}\index{lifting property!for $\Xdot^{s,p}$}. The Riesz potential $I_{\sigma} \coloneqq \cF^{-1}(|\cdot|^{-\sigma} \cF f)$ is an isomorphism $\Xdot^{s,p} \to \Xdot^{s+\sigma,p}$. This is proved in~\cite[Sec.~5.2.3]{Triebel-TheoryOf} for $p<\infty$ and follows by duality for $p=\infty$, see \cite[Rem.~2.3.8.2]{Triebel-TheoryOf}. On the basis of \eqref{eq: BMO identification} we find that
\begin{align*}
	\Hdot^{s, \infty} = I_s (\BMO) \eqqcolon \BMOdot^s \quad (0<s<1)
\end{align*}
agree up to equivalent norms with Strichartz' $\BMOdot^s$-spaces\index{Strichartz' $\BMO$ Sobolev spaces}. We have $\Hdot^{s,\infty} \subseteq \Bdot^{s,\infty}$ with continuous inclusion as a mere consequence of the definitions and the inclusion $\ell^2(\IZ) \subseteq \ell^\infty(\IZ)$. In particular, $\Hdot^{s, \infty}$ is a space of H\"older continuous functions of exponent $s$. An equivalent, more concrete norm is given by
\begin{align}
\label{eq: Strichartz BMO norm}
\|f\|_{\BMOdot^s} \coloneqq \sup_{Q} \bigg(\frac{1}{|Q|} \int_Q \int_Q \frac{|f(x) - f(y)|^2}{|x-y|^{n + 2s}} \, \d x \d y \bigg)^{\frac{1}{2}},
\end{align}
where the supremum is taken over all cubes $Q \subseteq \R^n$, see~\cite[Thm.~3.3]{StrichartzBMO}\index{Sobolev space!$\BMO$}. 

Together with the Mihlin multiplier theorem\index{Theorem!Mihlin multiplier}~\cite[Sec.~5.2.2/3]{Triebel-TheoryOf} the lifting property also yields for $p<\infty$ that
\begin{align*}
\Hdot^{1,p} &= \{f \in \cS'/\IC : \nabla_x f \in \Hdot^{p} \} \\
\|f\|_{\Hdot^{1,p}} &\simeq \|\nabla_x f\|_{\Hdot^{p}}.
\end{align*}
For $p>1$ these are the more common homogeneous Sobolev spaces~\index{Sobolev space!homogeneous} and we write
\begin{align*}
	\Wdot^{1,p} \coloneqq \Hdot^{1,p} \quad  \& \quad \Wdot^{-1,p} = (\Wdot^{1,p'})^* = \Hdot^{-1,p} \quad (1<p<\infty).
\end{align*}
In our usual range of exponents $p \in (1_*, \infty)$ any distribution $f \in \Hdot^{1,p}$ can be identified with a locally integrable function. This follows by density of $\cZ$ in $\Hdot^{1,p}$ and the extended Sobolev embedding theorem\index{Sobolev embedding theorem} that we recall for later reference.

\begin{prop}
\label{prop: Sobolev}
There are continuous embeddings
\begin{alignat*}{2}
\Hdot^{1,p} &\subseteq \Hdot^{p^*} &&\quad (1_* < p < n), \\
\Hdot^{1,p} &\subseteq \Lamdot^{1-\frac{n}{p}} &&\quad(n \leq p < \infty).
\end{alignat*}
\end{prop}

The second part is the classical Morrey inequality~\cite[Thm.~7.17]{GT}. The first part is a special case of the general embedding theorem\index{Sobolev embedding theorem!for $\Xdot^{s,p}$} 
\begin{align*}
	\Xdot^{s_0,p_0} \subseteq \Xdot^{s_1,p_1} \quad (0 < p_0 < p_1 < \infty,\, s_0 - \nicefrac{n}{p_0} = s_1  - \nicefrac{n}{p_1}),
\end{align*}
see \cite[Thm.~2.1]{Jawerth77}. 
\subsection{Interpolation functors}
\label{subsec: complex interpolation}

Here, and throughout the monograph, `complex interpolation'\index{interpolation!complex} refers to the Kalton--Mitrea complex interpolation method~\cite[\SS 3]{Kalton-Mitrea}, which is well-defined for quasi-Banach spaces and agrees with the classical Calder\'on complex interpolation method on couples of Banach spaces. As usual, we write $[\cdot\,,\cdot]_\theta$, $\theta \in (0,1)$, for the complex interpolation bracket. `Real interpolation'\index{interpolation!real} refers to the classical $K$-method~\cite[Sec.~3.10]{BL} and the corresponding interpolation bracket is denoted by $(\cdot\,,\cdot)_{\theta,p}$, $\theta \in (0,1)$, $p \in (0,\infty]$.

We gather the standard interpolation formul\ae\, that will be needed in the further course. To this end we let $0<p_0,p_1 \leq \infty$, $s_0, s_1 \in \R$, $\theta \in (0,1)$ and set $p \coloneqq [p_0, p_1]_\theta$, $s \coloneqq (1-\theta)s_0 + \theta s_1$. 

As for tent and $\Z$ spaces, we have\index{interpolation! of tent spaces}\index{interpolation! of $\Z$-spaces} up to equivalent quasi-norms,
\begin{alignat*}{2}
\big[\T^{s_0,p_0}, \T^{s_1,p_1}\big]_\theta 
&= \T^{s, p} \quad &&(\text{one $p_i$ finite}), \\
\big(\T^{s_0,p_0}, \T^{s_1,p_1}\big)_{\theta, p} 
&= \Z^{s, p} \quad &&(\text{$s_0 \neq s_1$}), \\
\big(\Z^{s_0,p_0}, \Z^{s_1,p_1}\big)_{\theta, p} 
&= \Z^{s, p} \quad &&(\text{$s_0 \neq s_1$}),
\end{alignat*}
see \cite[Thm.~2.12 \& Thm.~2.30 \& Prop.~2.31]{AA}. A different proof for complex interpolation of tent spaces can be found in \cite[Thm.~4.3]{Huang-Tent}.

The required interpolation identities for $\Xdot^{s,p}$ have been proved in \cite[Thm.~4.28 \& Thm.~5.2]{AA} via an approach based on tent and $\Z$ spaces. Their proof uses the language of operator-adapted spaces that will be introduced in Section~\ref{sec: Hardy intro}. We have\index{interpolation! of Hardy--Sobolev spaces}\index{interpolation! of Besov spaces} up to equivalent quasinorms
\begin{alignat*}{2}
\big[\Hdot^{s_0,p_0}, \Hdot^{s_1,p_1}\big]_\theta 
&= \Hdot^{s,p} \quad &&(\text{one $p_i$ finite}), \\
\big(\Xdot^{s_0,p_0}, \Xdot^{s_1,p_1}\big)_{\theta,p} 
&= \Bdot^{s,p} \quad &&(\text{$s_0 \neq s_1$}).
\end{alignat*}
Different proofs for some of the identities have been given in many earlier references. A complete reference for complex interpolation is \cite{KMM}, see also \cite{Triebel-TheoryOf, Jawerth-Frazier}. For real interpolation, the techniques in \cite[Chap.~6]{BL} apply in the range $p_0, p_1 > 1$ whereas in the range $p_0, p_1 \leq 2$ all results can be found in \cite[Sect.~3.3]{Qui}. Finally, Wolff reiteration~\cite{Wolff-original} gives the full picture.
\section{Preliminaries on operator theory}
\label{sec: L2}

\noindent In this chapter, we introduce the elliptic operators used in this monograph and recall their main properties in the $\L^2$-setting. We also recall material on (bi)sectorial operators and their holomorphic functional calculus. A particularly useful reference for our purpose is Haase's book~\cite{Haase} and the reader is advised to refer thereto whenever necessary.

\subsection{Definition of the elliptic operators}
\label{subsec: definition elliptic operators}

We let $a$ and $d$ be the coefficients of $\Le$ as in \eqref{eq:block}. The bounded multiplication operator $B$ and the first-order \emph{Dirac operator}\index{Dirac operator} $D$ are defined with maximal domain\index{D@$\dom$ (domain)} in $\L^2(\R^n; \IC^{m} \times \IC^{mn})$ by
\begin{align*}
B\coloneqq\begin{bmatrix} a^{-1} & 0  \\ 
0 & d \end{bmatrix},
\quad 
D\coloneqq\begin{bmatrix} 0 & \div_x  \\ 
-\nabla_x & 0 \end{bmatrix}.
\end{align*}
We note that $D$ is self-adjoint. Hence, it splits $\L^2$ into an orthogonal sum $\nul(D) \oplus \cl{\ran(D)}$. The null space $\nul(D)$\index{N@$\nul$ (null space)} consists of all $f = [f_\no, f_\ta]^\top$\index{$\no$-$\ta$ notation} with $f_\no = 0$ and $\div_x f_\ta = 0$ and the closure of the range $\cl{\ran(D)} =\cH$\index{R@$\ran$ (range)} is the space in our ellipticity assumption \eqref{eq: accretivity A}\index{H@$\cH$ (closure of $\ran(D)$)}. Consequently, \eqref{eq: accretivity A} is equivalent to
\begin{align*}
\Re \int_{\R^n} Bf \cdot \cl{f} \, \d x \geq \lambda \int_{\R^n} |a^{-1}f_\no|^2 + |f_\ta|^2 \, \d x \quad (f \in \cl{\ran(D)}),
\end{align*}
or again, using angular brackets to denote inner products, equivalent to
\begin{align}
\label{eq: B accretive}
\Re \langle BDu,Du \rangle \gtrsim \|Du\|_2^2 \quad (u \in \dom(D)).
\end{align} 
Because of this, we say that $B$ is accretive (or elliptic)~\index{ellipticity!of $B$} on the range of $D$. 

The \emph{perturbed Dirac operators}\index{Dirac operator!perturbed}
\begin{align}
\label{eq: BD and BD}
BD \coloneqq \begin{bmatrix} 0& a^{-1}\div_x  \\ -d\nabla_x  & 0 \end{bmatrix},
\quad
DB \coloneqq \begin{bmatrix} 0& \div_x d \\ -\nabla_x a^{-1} & 0 \end{bmatrix}
\end{align}
are again considered with maximal domain in $\L^2$. Since $B$ is bounded, $DB$ is closed and as consequence of \eqref{eq: B accretive} also $BD$ is closed. Their squares contain the following second-order operators\index{second-order operator!$L, \tL, M, \tM$}:
\begin{align}
\label{eq: L and M}
\begin{bmatrix} L \vphantom{\tL} & 0 \\ 0 & M \vphantom{\tM} \end{bmatrix}
&\coloneqq 
\begin{bmatrix} -a^{-1}\div_x d \nabla_x \vphantom{\tL} & 0 \\ 0 & -d\nabla_x a^{-1} \div_x \vphantom{\tM}\end{bmatrix}
= (BD)^2,\\
\label{eq: tL and tM}
\begin{bmatrix} \tL& 0 \\ 0 & \tM\end{bmatrix}
&\coloneqq
\begin{bmatrix} -\div_x d \nabla_x a^{-1 \vphantom{\tL}}& 0 \\ 0 & -\nabla_x a^{-1} \div_x d \vphantom{\tM}\end{bmatrix}
=(DB)^2.
\end{align}
The definition of $L$ coincides with the more traditional variational approach to defining second-order operators. Indeed, the Lax--Milgram lemma provides an isomorphism, 
\begin{align}
\label{eq: Lax-Milgram operator}
\eo: \Wdot^{1,2} \to \Wdot^{-1,2}, \quad \langle \eo u, v \rangle = \int_{\R^n} d \nabla_x u \cdot \cl{\nabla_x v} \, \d x.
\end{align}
We have $\eo u \coloneqq -\div_x d \nabla_xu$ in the sense of distributions and one sees that $u \in \dom(L)$ means that $u \in \L^2 \cap \Wdot^{1,2}$ with $\eo u \in \L^2$ and $Lu = a^{-1}\eo u$. Note that the domain of $L$ does not depend on $a$. Occasionally, we will write \index{second-order operator!$L_0$}
\begin{align*}
	L_0 \coloneqq -\div_x d \nabla_x
\end{align*}
for the divergence form operator $L$ in the special case $a=1$, that is to say, the maximal restriction of $\eo$ to an operator in $\L^2$.
\subsection{(Bi)sectorial operators}
\label{subsec: bisectorial operators}

Statements and proofs for sectorial and bisectorial operators usually go \emph{mutadis mutandis}. Most authors have decided to showcase  sectorial operators. In case of doubt the reader can consult \cite[Ch.~3]{E}, which goes the other way round.

Let $\omega \in (0, \pi)$. We define the sector $\S_\omega^+ \coloneqq \{z\in\IC : | \arg z|<\omega\}$\index{sector ($\S_\mu^+$)} and agree on $\S_0^+ \coloneqq (0,\infty)$. \index{arg@$\arg$ (argument function with range $(-\pi, \pi]$.} A linear operator $T$ on a reflexive Banach space $X$ is \emph{sectorial}\index{operator!sectorial} of angle $\omega \in [0, \pi)$ if its spectrum is contained in $\cl{\S_\omega^+}$ and if for every $\mu \in (\omega, \pi)$,
\begin{align}
\label{eq: constant for resolvent bound}
M_{T,\mu} \coloneqq  \sup_{z \in \IC \setminus \cl{\S_\mu^+}}  \|z(z-T)^{-1}||_{X \to X} < \infty.
\end{align}
Usually, $\omega_T$ denotes the smallest angle $\omega$ with this property. A sectorial operator is densely defined, induces a topological kernel/range splitting
\begin{align}
\label{eq: kernel range splitting}
X = \nul(T) \oplus \cl{\ran(T)},
\end{align}
and the restriction of $T$ to $\cl{\ran(T)}$ is sectorial, injective and has dense range~\cite[Prop.~2.1.1]{Haase}. 

\emph{Bisectorial}\index{operator!bisectorial} operators of angle $\omega \in [0,\nicefrac{\pi}{2})$ are defined analogously upon replacing sectors with bisectors $\S_\omega \coloneqq \S_\omega^+ \cup (-\S_\omega^+)$\index{bisector ($\S_\mu$)} and share the same properties. If $T$ is bisectorial of angle $\omega$, then writing
\begin{align*}
(z^2-T^2)^{-1} = -(z-T)^{-1}(-z-T)^{-1},
\end{align*}
we see that $T^2$ is sectorial of angle $2 \omega$. Moreover, $\nul(T^2) = \nul(T)$ and hence $\cl{\ran(T^2)} = \cl{\ran(T)}$, see~\cite[Prop.~2.1.1e)]{Haase}.

As prototypical examples, $BD$ and $DB$ are bisectorial of the same angle $\omega_{BD} = \omega_{DB}$ with
\begin{align}
\label{eq: ran BD and DB}
\ran(BD) = B \ran(D), \quad \ran(DB) = \ran(D),
\end{align}
see \cite[Prop.~3.3]{AAMc}. From \eqref{eq: L and M} and \eqref{eq: tL and tM} we obtain that $L$, $M$, $\tL$, $\tM$ are sectorial of angle not larger than $2\omega_{BD}$, but possibly exceeding $\nicefrac{\pi}{2}$, with
\begin{align}
\label{eq: ran L and M}
 \cl{\ran(L)} \times \cl{\ran(M)}
 &= (\L^2(\R^n; \IC^m)) \times d \cl{\ran(\nabla_x)}\\
\label{eq: ran tL and tM}
 \cl{\ran(\tL)} \times \cl{\ran(\tM)}
&= (\L^2(\R^n; \IC^m)) \times \cl{\ran(\nabla_x)}.
\end{align}
In particular, $L$ and $\tL$ have dense range and hence they are injective. 
\subsection{Classes of holomorphic functions}
\label{subsec: classes of holomorphic functions}

Let $\mu \in (0,\pi)$. The classes $\Psi_\sigma^\tau(\S_\mu^+)$\index{$\Psi_\sigma^\tau$, $\Psi^+_+$, $\H^\infty$ (classes of holomorphic functions)}, $\sigma,\tau \in \R$, consist of those holomorphic functions $\varphi: \S_\mu^+ \to \IC$ that satisfy
\begin{align*}
|\varphi(z)| \lesssim |z|^\sigma \wedge |z|^{-\tau} \quad (z \in \S_\mu^+).
\end{align*}
We write $\H^\infty(\S_\mu^+) \coloneqq \Psi_0^0(\S_\mu^+)$ for the bounded holomorphic functions on $\S_\mu^+$. The classes of functions with some decay and arbitrarily large polynomial decay  at $0$ and $\infty$ are 
\begin{align*}
	\Psi_+^+(\S_\mu^+) \coloneqq \bigcup_{\sigma, \tau>0} \Psi_\sigma^\tau(\S_\mu^+) \quad \text{and} \quad \Psi_\infty^\infty(\S_\mu^+) \coloneqq \bigcap_{\sigma, \tau>0} \Psi_\sigma^\tau(\S_\mu^+),
\end{align*}
respectively. We suppress reference to $\S_\mu^+$ in the notation when the relevant sector is clear from the context.

On bisectors we use the same notation and call a function \emph{non-degenerate}\index{non-degenerate function} if it does not identically vanish on one of the two connected components. An example of a degenerate function is $z+[z]$, where\index{$[z]$ (holomorphic function)}
\begin{align*}
[z] \coloneqq \sqrt{z^2} \quad (z \in \IC \setminus \i \R)
\end{align*}
is defined via the principal branch of the logarithm.

\subsection{Holomorphic functional calculi}
\label{subsec: functional calculi}

For the same reason as before, we can focus on the sectorial case. So, let $T$ be sectorial and let $\mu \in (\omega_T,\pi)$. If $\psi$ is of the form $\psi(z) = \alpha + \beta(1+z)^{-1} + \varphi(z)$ for some $\alpha,\beta \in \IC$ and $\varphi \in \Psi_+^+(\S^+_\mu)$, then $\psi(T)$ is defined as a bounded operator on $X$ via
\begin{align}
\label{eq: def basic calculus}
\psi(T) = \alpha + \beta(1+T)^{-1} + \frac{1}{2 \pi \i} \int_{\bd \S^+_\nu} \varphi(z)(z-T)^{-1} \, \d z, 
\end{align}
where $\nu \in (\omega_T, \mu)$, the choice of which does not matter in view of Cauchy's theorem, and $\bd \S^+_\nu$ is oriented such that it surrounds the spectrum of $T$ counter-clockwise in the extended complex plane. The definition extends to larger classes of functions by \emph{regularization}: If $e(T)$ and $(e\psi)(T)$ are already defined by the procedure above and if $e(T)$ is injective, then 
\begin{align*}
\psi(T) \coloneqq e(T)^{-1}(e\psi)(T)
\end{align*}
is defined as a closed operator and can be shown not  to depend on the choice of $e$. The expected relations 
\begin{align*}
\psi(T) + \phi(T) &\subseteq (\psi + \phi)(T)\\
\psi(T)\phi(T) &\subseteq (\psi \phi)(T)
\end{align*}
hold and there is equality if $\psi(T)$ is bounded.

Since the restriction of $T$ to $\cl{\ran(T)}$ is an injective sectorial operator, $e(z)=z(1+z)^{-2}$ regularizes any bounded holomorphic function in $\H^\infty(\S_\mu^+)$. The \emph{convergence lemma}\index{convergence lemma} states that if $(\psi_j)_j$ is a bounded sequence in $\H^\infty(\S_\mu^+)$ that converges pointwise to $\psi$ and if $\sup_j \|\psi_j(T)\|_{\cl{\ran(T)} \to \cl{\ran(T)}} < \infty$, then
\begin{align}
\label{eq: convergence lemma}
\psi(T) = \lim_{j \to \infty} \psi_j(T) 
\end{align}
in the sense of strong convergence on $\cl{\ran(T)}$. In our applications the additional assumption on $\|\psi_j(T)\|_{\cl{\ran(T)} \to \cl{\ran(T)}}$ is automatically satisfied because of the following property.

We say that $T$ has a \emph{bounded $\H^\infty$-calculus on $\cl{\ran(T)}$}\index{Hinfty calculus@$\H^\infty$-calculus} (of angle $\mu \geq \omega_T$) if for all $\nu \in (\mu,\pi)$ there is a constant $M_{T,\nu}^\infty$ such that 
\begin{align}
\label{eq: Hoo bound}
\|\psi(T)\|_{\cl{\ran(T)} \to \cl{\ran(T)}} \leq M_{T,\nu}^\infty \|\psi\|_{\L^\infty(\S^+_\nu)} \quad (\psi \in \H^\infty(\S^+_\nu)).
\end{align}
In fact, by the convergence lemma, it suffices to have the bound for all $\psi \in \Psi_+^+(\S_\nu^+)$. In Hilbert spaces, these properties are independent of the angle $\mu$. This is one of the statements of the following fundamental result due to McIntosh~\cite{Mc}, see also \cite[Thm.~7.3.1]{Haase} or \cite[Thm.~3.4.11]{E}. 

\begin{thm}[McIntosh\index{Theorem!McIntosh's}]
\label{thm: McIntosh}
Let $T$ be a (bi)sectorial operator in a Hilbert space $H$. Then $T$ has a bounded $\H^\infty$-calculus of some angle  on $\cl{\ran(T)}$ (equivalently, of angle $\omega_T$) if and only if the quadratic estimate
\begin{align*}
\|f\|_H \simeq \bigg(\int_0^\infty \|\varphi(tT)f\|_H^2 \, \frac{\d t}{t} \bigg)^{1/2}
\end{align*}
holds for all $f \in \cl{\ran(T)}$ and some (equivalently, all) admissible and non-degenerate $\varphi \in \Psi_+^+$. 
\end{thm}

\begin{rem}
The following dependence of the implicit constants easily follows from the proof in \cite{E} and is also explicitly stated in \cite[Thm.~10.4.16/19]{HNVW2}.
For fixed angle $\nu$, a bound $M_{T,\nu}^\infty$ for the $\H^\infty$-calculus depends on $M_{T,\mu}$ for some $\mu \in (\omega_T,\nu)$ and implicit constants in the quadratic estimates. Conversely, for $\varphi \in \Psi_\sigma^\tau$ with $\sigma,\tau> 0$ on some (bi)sector of angle $\nu$, the quadratic estimate 
\begin{align*}
	\bigg(\int_0^\infty \|\varphi(tT)f\|_H^2 \, \frac{\d t}{t} \bigg)^{1/2} \leq C \|f\|_H \qquad (f \in \cl{\ran(T)}) 
\end{align*}
holds with $C$ depending on $M^\infty_{T,\nu}$, $\nu$, $\sigma$, $\tau$ and $\sup_z \nicefrac{|f(z)|}{(|z|^\sigma \wedge |z|^{-\tau})}$.
\end{rem}

We also recall the important reproducing formula for sectorial operators~\cite[Thm.~5.2.6]{Haase} and remark that up to the usual modifications there is a bisectorial version~\cite[Prop.~4.2]{AusSta}. 

\begin{lem}[Calderón reproducing formula\index{Calderón reproducing formula}]
\label{lem: Calderon reproducing}
Let $T$ be a sectorial operator in a reflexive Banach space $X$ and let $\varphi \in \Psi_+^+$ on a suitable sector be such that $\int_0^\infty \varphi(t) \frac{\d t}{t} = 1$. Then 
\begin{align*}
\int_0^\infty \varphi(tT)f \, \frac{\d t}{t} = f \quad (f \in \cl{\ran(T)})
\end{align*}
as an improper strong Riemann integral. 
\end{lem}

\begin{rem}
\label{rem: Calderon reproducing}
For any non-zero $\phi \in \H^\infty$ there is $\psi \in \Psi_\infty^\infty$ on the same sector such that $\varphi \coloneqq \phi \psi$ satisfies the Calderón reproducing formula, for example $\psi(z) \coloneqq c \cl{\phi(\cl{z})} \e^{-\sqrt{z}-1/\sqrt{z}}$, where $c^{-1} = \int_0^\infty |\phi(t)|^2 \e^{-\sqrt{t}-1/\sqrt{t}}\, \frac{\d t}{t}$. 
\end{rem}

Coming back to concrete operators, quadratic estimates (and hence bounded functional calculi) for $BD$ and $DB$ is a deep result due to Axelsson--Keith--McIntosh~\cite{AKMc}. For a condensed proof, see also~\cite[Thm.~1.1]{AAMc2}

\begin{thm}[Axelsson--Keith--McIntosh\index{Theorem!Axelsson--Keith--McIntosh's}]
\label{thm: QE for DB}
The operators $BD$ and $DB$ have bounded $\H^\infty$-calculi on the closure of their ranges.
\end{thm}

Let now $\mu \in (2\omega_{BD}, \pi)$ and $\psi \in \Psi_+^+(\S_\mu^+)$. Then $\varphi$ defined by $\varphi(z) \coloneqq \psi(z^2)$ belongs to $\Psi_+^+(\S_{\mu/2})$. From \eqref{eq: L and M} we obtain
\begin{align}
\label{eq: calculus for L and M from BD}
\begin{bmatrix}
\psi(L) & 0 \\ 0 & \psi(M)
\end{bmatrix}
=\psi((BD)^2)
=\varphi(BD).
\end{align}
The same argument works for $\tL$, $\tM$ by referring to $DB$ instead. McIntosh's theorem implies the following

\begin{cor}
\label{cor: QE for 2nd order ops}
The operators $L$ and $\tL$ have bounded $\H^\infty$-calculi on $\L^2$. Likewise, $M$ and $\tM$ have bounded $\H^\infty$-calculi on the closure of their ranges.
\end{cor}

Since $B$ is accretive on $\cH = \cl{\ran(DB)}$ and maps this space onto $B \cH = \cl{\ran(BD)}$, it follows that $B|_\cH: \cl{\ran(DB)} \to \cl{\ran(BD)}$ is invertible and that the restrictions of $BD$ and $DB$ to the closure of their ranges are similar under conjugation with $B|_{\cH}$. Therefore
\begin{align}
\label{eq: similarity BD and DB}
\varphi(BD)B = B \varphi(DB)
\end{align}
holds as unbounded operators from $\cl{\ran(DB)}$ into $\cl{\ran(BD)}$, whenever one side is defined by the respective functional calculus. Elaborating further along these line, we obtain 

\begin{lem}[Intertwining relations\index{intertwining relations}]
\label{lem: intertwining}
Let $\varphi \in \H^\infty$ on a suitable bisector and $\psi \in \H^\infty$ on a suitable sector. Then 
\begin{align*}
D \varphi(BD)f = \varphi(DB) Df \quad (f \in \dom(D))
\end{align*}
and 
\begin{align*}
\div_x \psi(M\vphantom{\tM})f_\ta &= \psi(\tL)\div_x f_\ta, \quad (f_\ta \in \dom(\div_x)),\\
\nabla_x \psi(L\vphantom{\tL})f_\no &= \psi(\tM) \nabla_x f_\no, \quad (f_\no \in \W^{1,2}).
\end{align*}
\end{lem}
	
\begin{proof}
For the first identity we note that $Df \in \ran(D) = \ran(DB)$ by \eqref{eq: ran BD and DB}. Hence, we can apply \eqref{eq: similarity BD and DB} to $Df$ in order to obtain 
\begin{align*}
BD \varphi(BD)f = B \varphi(DB) Df
\end{align*}
and the claim follows since $B$ is accretive on $\cl{\ran(D)}$. By means of \eqref{eq: calculus for L and M from BD} and the analogous identity for $DB$ the identities for $L$ and $M$ follow.
\end{proof}	
\subsection{Adjoints}
\label{subsec: adjoints}

We note that the adjoint of a (bi)sectorial operator in a Hilbert space is again bisectorial of the same angle~\cite[Prop.~2.1.1]{Haase} and that $B^*$ has the same properties as $B$. Since  $B$ is bounded, we have $(BD)^* = DB^*$ and likewise $(B^*D)^* = DB$, which yields $(DB)^* = B^*D$ because $B^*D$ is closed. Since all these operators are bisectorial, we obtain $((BD)^2)^* = (DB^*)^2$, which in matrix form reads
\begin{align}
\label{eq: adjoints second order}
\begin{bmatrix}
L^* & 0 \\ 0 & M^*
\end{bmatrix}
= 
\begin{bmatrix}
-\div_x d^* \nabla_x (a^*)^{-1} & 0 \\ 0 & -\nabla_x (a^*)^{-1} \div_x d^*
\end{bmatrix}.
\end{align}
The $\Psi_+^+$-calculus of any (bi)sectorial operator dualizes in the expected manner $\psi(T)^* = \psi^*(T^*)$, where $\psi^*(z) = \cl{\psi(\cl{z})}$~\index{$\psi^*$ (holomorphic conjugate $z \mapsto \cl{\psi(\cl{z})}$)}. If $T$ has dense range, for example $T=L$, then this relation also holds for all $\psi \in \cup_{\sigma,\tau \in \R} \Psi_\sigma^\tau$, see \cite[Prop.~7.0.1(d)]{Haase}. When $a=1$, the operator $L^*$ is in the same class as $L$. When $a\ne 1$, the operator $L^*$ is not in the same class as $L$ but is similar to such an operator under conjugation with $a^*$. This is why instead of $L^*$ we usually work with\index{second-order operator!$L^\sharp$}
\begin{align*}
L^\sharp \coloneqq -(a^*)^{-1} \div_x d^* \nabla_x = (a^*)^{-1} L^* a^*
\end{align*}
when it comes to duality arguments. 
\subsection{Kato problem and Riesz transform}
\label{subsec: Kato and Riesz}

Since $z \mapsto \nicefrac{[z]}{z}$ and its inverse are bounded and holomorphic on any bisector, the bounded $\H^\infty$-calculus for $BD$ entails that $BD$ and $[BD]$ share the same domain along with comparability
\begin{align*}
\big \|BD f \big\|_2 \simeq \big \|[BD]f \big \|_2 \quad (f \in \dom(BD)).
\end{align*}
The left-hand side is also comparable to $\|Df\|_2$ by ellipticity. Looking at the first component and using the specific form of $(BD)^2$, see \eqref{eq: L and M}, we obtain the resolution of the \emph{Kato conjecture}.

\begin{thm}[Resolution of the Kato conjecture\index{Kato problem/conjecture}]
\label{thm: Kato}
It follows that $\dom({L^{1/2}}) = \W^{1,2}$ with the homogeneous estimate $\|{L^{1/2}}f\|_2 \simeq \|\nabla_x f\|_2$.
\end{thm}

As a consequence, we obtain a bounded extension ${L^{1/2}} : \Wdot^{1,2} \to \L^2$ by density that is injective with closed range. It is an isomorphism since its range contains $\ran(L)$, which is dense in $\L^2$ by \eqref{eq: ran L and M}. We denote its inverse by $L^{-1/2}$\index{L@$L^{-1/2}$ (extension by density)}. In particular, the \emph{Riesz transform}\index{Riesz transform} $\nabla_x L^{-1/2}$ is a bounded operator on $\L^2$.

The domains of fractional powers of exponent $\alpha \in (0,\nicefrac{1}{2})$ can be determined by complex interpolation.

\begin{cor}
\label{cor: Kato}
If $\alpha \in (0,\nicefrac{1}{2})$, then $\dom({L^{\alpha}}) = \Hdot^{2\alpha,2} \cap \L^2$ with the homogeneous estimate $\|{L^{\alpha}}f\|_2 \simeq \|f\|_{\Hdot^{2\alpha,2}}$.
\end{cor}

\begin{proof}
By \cite[Thm.~5.1]{AMcN} we have $\dom({L^{\alpha}}) = [\L^2, \Wdot^{1,2}]_{2\alpha} \cap \L^2$ with the homogeneous estimate $\|L^\alpha f \|_2 \simeq \|f\|_{[\L^2, \Wdot^{1,2}]_{2\alpha}}$ and from Section~\ref{subsec: complex interpolation} we know that $[\L^2, \Wdot^{1,2}]_{2\alpha} = \Hdot^{2\alpha,2}$.
\end{proof}

\subsection{Off-diagonal estimates}
\label{subsec: off-diagonal estimates}

We develop on these estimates in Section~\ref{sec: Hp Hq boundedness} below. Here we only gather the well-known $\L^2$-bounds for our standard operators from Section~\ref{subsec: bisectorial operators}.

\begin{defn}
\label{def: off-diagonal}
Let $\Omega \subseteq \IC \setminus \{0\}$ and let $V_1, V_2$ be finite-dimensional Hilbert spaces. A family $(T(z))_{z \in \Omega}$ of linear operators $\L^2(\R^n; V_1) \to \L^2(\R^n; V_2)$ satisfies \emph{$\L^2$ off-diagonal estimates of order $\gamma>0$}\index{off-diagonal estimates!$\L^2$ of order $\gamma$} if there exists a constant $C$ such that 
\begin{align*}
\|\ind_F T(z) \ind_E f \|_2 \leq C \bigg(1+ \frac{\dist(E,F)}{|z|} \bigg)^{-\gamma} \|\ind_E f\|_2
\end{align*}
holds for all measurable subsets $E,F \subseteq \R^n$, all $z \in \Omega$ and all $f \in \L^2(\R^n; V_1)$. It there are constants $C, c>0$ such that the stronger estimate
\begin{align*}
	\|\ind_F T(z) \ind_E f \|_2 \leq C \e^{-c \frac{\dist(E,F)}{|z|}} \|\ind_E f\|_2
\end{align*}
holds, then the family is said to satisfies \emph{off-diagonal estimates of exponential order}\index{off-diagonal estimates!$\L^2$ of exponential order}.
\end{defn}

While decay of polynomial order is most suitable for the abstract theory that we develop in the upcoming sections, our prototypes actually satisfy the exponential estimate. For completeness, we include the argument from~\cite[Prop.~5.1]{AAMc2}.\index{off-diagonal estimates!for $DB, BD$}

\begin{prop}
\label{prop: OD for Dirac}
The resolvent families $((1+\i t BD)^{-1})_{t \in \R \setminus \{0\}}$ and $((1+\i t DB)^{-1})_{t \in \R \setminus \{0\}}$ satisfy $\L^2$ off-diagonal estimates of exponential order.
\end{prop}

\begin{proof}
We begin with the resolvents $T(t) \coloneqq (1+\i t BD)^{-1}$. Fix $t$, $E$, $F$ and set $\dist \coloneqq \d(E,F)$. The family $(T(t))_{t \in \R \setminus \{0\}}$ is uniformly bounded in $\L^2$ since $BD$ is bisectorial. Hence, it suffices to obtain the exponential estimate for $|t| \leq \alpha \dist$, where $\alpha > 0$ will be chosen later on in dependence of dimensions and ellipticity.

We introduce $G \coloneqq \{x \in \R^n: \d(x,F) \leq \nicefrac{\dist}{2}\}$. As $\d(F, {}^c G) \geq \nicefrac{\dist}{2}$, we can pick a smooth  function $\varphi$ that satisfies $\ind_F \leq \varphi \leq \ind_G$ and $\|\nabla_x \varphi\|_\infty \leq \nicefrac{C}{\dist}$ for some dimensional constant $C$. Let  $\eta \coloneqq \e^{(\alpha \dist/|t|)\varphi} - 1$ and observe that 
\begin{align*}
	\eta = 0 \quad  (\text{on } E) \quad \text{and} \quad \eta = \e^{\frac{\alpha \dist}{|t|}} -1 \geq \frac{1}{2} \e^{\frac{\alpha \dist}{|t|}} \quad (\text{on } F).
\end{align*}
Thus, we obtain for all $f \in \L^2$ that
\begin{align}
\label{eq1: OD for Dirac}
	\frac{1}{2} \e^{\frac{\alpha \dist}{|t|}} \|\ind_F T(t) \ind_E f\|_2
	\leq \|\eta T(t) \ind_E f\|_2 = \|[\eta, T(t)] \ind_E f\|_2,
\end{align}
where $[\eta, T(t)] = \eta T(t) - T(t)(\eta \, \cdot)$ is the commutator between $T(t)$ and multiplication with $\eta$. Next, we expand
\begin{align}
\label{eq2: OD for Dirac}
	[\eta, T(t)] = T(t) [1+ \i t BD, \eta] T(t) = \i t T(t) B [D, \eta] T(t).
\end{align}
By the product rule we find that $[D,\eta]$ acts via multiplication by a function $\theta  \e^{(\alpha \dist/|t|)\varphi}$, where $\theta$ is supported in $G$ and uniformly bounded by a dimensional multiple of  $\nicefrac{\alpha \dist}{|t|} \|\nabla \varphi\|_\infty \leq \nicefrac{C \alpha}{|t|}$. Since $T(t)$ and $B$ are (uniformly) bounded on $\L^2$, we conclude that  
\begin{align*}
	\|\eta T(t) \ind_E f\|_2 
	& \leq C \alpha   \big \|  \e^{ \frac{ \alpha \dist}{|t|}\varphi} T(t) \ind_E f \big \|_2 \\
	& \leq C \alpha \big(\|\eta T(t) \ind_E f\|_2 + \|T(t) \ind_E f\|_2 \big),
\end{align*}
where $C$ depends on ellipticity and dimension and the second step merely follows from $\eta = \e^{(\alpha \dist/|t|)\varphi} - 1$. Setting $\alpha \coloneqq \nicefrac{1}{2C}$, we can absorb the first term on the right-hand side back into the left-hand side and we are left with 
\begin{align*}
		\|\eta T(t) \ind_E f\|_2 \leq \|T(t) \ind_E f\|_2.
\end{align*}
Using \eqref{eq1: OD for Dirac} on the left and uniform boundedness of $T(t)$ on the right completes the proof for the resolvents of $BD$.

For $DB$ the only modification in the argument concerns \eqref{eq2: OD for Dirac}, where $B$ appears on the right of $[D, \eta]$.
\end{proof}

\begin{rem}
\label{rem: OD for Dirac}
The off-diagonal estimates extend to complex parameters $t=z \in \S_\mu$ for any $\mu \in (0, \nicefrac{\pi}{2} - \omega_{BD})$. The proof is literally the same but it is also instructive to remark that one can use Stein interpolation against the uniform resolvent bounds. This argument appears in greater generality in Lemma~\ref{lem: OD extrapolation to sectors} below.
\end{rem}

\begin{cor}
\label{cor: OD for second order}
The following families $(T(z))_{z \in \S_\mu^+}$ satisfy off-diagonal estimates of exponential order:
\begin{enumerate}
	\item $T(z) = (1+z^2 T)^{-1}$ if $\mu \in (0, \nicefrac{(\pi - \omega_T)}{2})$ and $T \in \{L, \tL, M, \tM\}$.
	\item $T(z) = z\nabla_x(1+z^2L)^{-1}$ if $\mu \in (0, \nicefrac{(\pi - \omega_L)}{2})$.
\end{enumerate}
In particular, these families satisfy $\L^2$ off-diagonal estimates of arbitrarily large order.
\end{cor}

\begin{proof}
By Stein interpolation, see the preceding remark, it suffices to argue for $z \in \R \setminus \{0\}$. Thanks to Proposition~\ref{prop: OD for Dirac} we have off-diagonal estimates of exponential order for
\begin{align*}
\frac{1}{2} \bigg((1+\i z BD)^{-1} + (1-\i z BD)^{-1} \bigg) = (1+z^2 (BD)^2)^{-1},
\end{align*}
as well as for the corresponding family with $DB$ replacing $BD$. Thus, (i) follows from \eqref{eq: L and M} and \eqref{eq: tL and tM}. Similarly, we have
\begin{align*}
\frac{1}{2} \bigg((1+\i z DB)^{-1} - (1-\i z DB)^{-1} \bigg)
= \begin{bmatrix}
- \i z \div_x d (1+z^2 \tM)^{-1} \\ \i z \nabla_x a^{-1} (1+z^2 \tL)^{-1}
\end{bmatrix}
\end{align*}
and we obtain the required off-diagonal estimates for
\begin{align*}
z \nabla_x a^{-1} (1+z^2 \tL)^{-1} = z \nabla_x (1+z^2L)^{-1} a^{-1}
\end{align*}
as stated in (ii).\index{off-diagonal estimates!for $L, \tL, M, \tM$}
\end{proof}
\section{\texorpdfstring{$\H^p - \H^q$}{Hp - Hq} bounded families}
\label{sec: Hp Hq boundedness}

\noindent In this section we discuss general principles for $\H^p - \H^q$-bounded operator families. We provide a toolbox that will allow us to manipulate resolvent families associated with our first and second-order operators efficiently on an abstract level.
\subsection{Abstract principles}
\label{subsec: Hp Hq abstract}

Throughout we work under the following assumption unless stated otherwise:
\begin{align}
\label{eq: standard assumptions Hp-Hq}
\begin{minipage}{0.89\linewidth}
\begin{itemize}
\item $(T(z))_{z \in \Omega}$ is a family of bounded operators $\L^2(\R^n; V_1) \to \L^2(\R^n; V_2)$ indexed over some set $\Omega \subseteq \IC \setminus \{0\}$, where the $V_i$ are finite-dimensional Hilbert spaces,
\item $a_i \in \L^\infty(\R^n; \Lop(V_i))$, $i=1,2$, are such that $a_i(x)$ is invertible for a.e.\ $x$ and $a_i^{-1} \in \L^\infty(\R^n; \Lop(V_i))$.
\end{itemize}
\end{minipage}
\end{align}

\begin{defn}
\label{def: p-q boundedness}
Let $(T(z))_{z \in \Omega}$ be an operator family as in \eqref{eq: standard assumptions Hp-Hq} and let $0<p\leq q<\infty$. This family is \emph{$a_1\H^p - a_2\H^q$-bounded}\index{boundedness!$a_1\H^p - a_2\H^q$} if 
\begin{align}
\label{eq: p-q boundedness}
\|a_2^{-1} T(z)a_1 f\|_{\H^q} \lesssim |z|^{\frac{n}{q}- \frac{n}{p}} \|f\|_{\H^p} \quad (z \in \Omega, f \in \H^p \cap \L^2).
\end{align}
\end{defn}

Usually, $\Omega$ is a half-line, a sector or a bisector in our application, hence the follow-up on the scaling in \eqref{eq: p-q boundedness}.

\begin{rem}
\label{rem: p-q boundedness}
\begin{enumerate}
\item We omit $\Omega$ and simply write $(T(z))$ when the context is clear. We speak of \emph{$a\H^p$-boundedness} when $a_1 = a_2 = a$ and $p=q$. If $q> 1$, then multiplication by $a_2$ is an automorphism of $\H^q = \L^q$ and hence $a_2$ may be dropped on the left-hand side of \eqref{eq: p-q boundedness}. We simply speak of \emph{$a_1\H^p-\L^q$-boundedness}. If also $p>1$, then $a_1$ may be dropped as well and we speak of \emph{$\L^p-\L^q$-boundedness}\index{boundedness!$\L^p-\L^q$} (\emph{$\L^p$-boundedness} if $p=q$).

\item Occasionally, we shall use the following extensions to the notions above. First, we can include endpoint Lebesgue spaces for $a_1\H^p - \L^q$, $q\in \{1,\infty\}$,  and $\L^p-\L^q$-boundedness, $p,q \in \{1,\infty\}$.  Second, when $0<p<\infty$ and $0 \leq \alpha <1$, we speak of $a_1 \H^p - a_2 \Lamdot^\alpha$-boundedness if 
\begin{align*}
\qquad \qquad \|a_2^{-1}T(z)a_1f\|_{\Lamdot^\alpha } \lesssim |z|^{-\alpha - \frac{n}{p}} \|f\|_{\H^p} \quad (z \in \Omega, f \in \H^p \cap \L^2)
\end{align*}
and make the same kind of notational abbreviations and extensions as before.\index{boundedness!$\H^p-\Lamdot^\alpha$}
\end{enumerate}
\end{rem}

Since the Hardy spaces interpolate by the complex method and have a universal approximation technique, the notion of $a_1\H^p-a_2\H^q$-boundedness interpolates as well. Moreover, the notions `dualize' in the expected way as the next lemma shows.\index{duality principle (for $p-q$-boundedness)}

\begin{lem}
\label{lem: p-q boundedness duality}
Let $(T(z))$ be as in \eqref{eq: standard assumptions Hp-Hq}. 
\begin{enumerate}
	\item If $1 \leq p \leq q \leq \infty$, then $(T(z))$ is $\L^p-\L^q$-bounded if and only if $(T(z)^*)$ is $\L^{q'}-\L^{p'}$-bounded.
	\item If $1_* < p \leq 1 \leq q \leq \infty$, then $(T(z))$ is $a_1\H^p-\L^q$-bounded if and only if $(T(z)^*)$ is $\L^{q'}-(a_1^*)^{-1}\Lamdot^{n(\frac{1}{p}-1)}$-bounded.
\end{enumerate}
\end{lem}

\begin{proof}
We can assume $a_1 = 1$ and $a_2 =1$ --- otherwise we replace $(T(z))$ by $(a_2^{-1} T(z) a_1)$. All of the claims take the abstract form that one of $(T(z))$ and $(T(z)^*)$ is $X_1-X_2$-bounded and the other one should be $X_3-X_4$-bounded. As
\begin{align*}
\langle T(z) f , g \rangle = \langle f, T(z)^* g \rangle \quad (z \in \Omega, \, f,g \in \L^2),
\end{align*}
it suffices to know that the $X_4$-norm can be computed by testing against functions in $X_1 \cap \L^2$. Above, either $X_1$ is a Hardy or Lebesgue space and $X_4$ is its dual (so the claim follows since $X_1 \cap \L^2$ is dense in $X_1$) or $X_1 = \L^\infty$ and $X_4 = \L^1$ (and the claim follows by testing against characteristic functions of bounded sets).
\end{proof}

The next lemma provides us with a useful criterion for a family to map a given $\H^q$-space back into $\H^2 = \L^2$.

\begin{lem}
\label{lem: extra} 
Let $(T(z))$ be a family as in \eqref{eq: standard assumptions Hp-Hq} with $V_1 = V_2 \eqqcolon V$ and $a_1 = a_2 \eqqcolon a$. Suppose that $(T(z))$ is $\L^2$-bounded and there exist $p,\varrho \in (0,2)$ for which $(T(z))$ is $a\H^p - a\H^p$ and $a\H^\varrho - \L^2$-bounded. 
Then, for each $q \in (p,2)$, there exists an integer $\beta = \beta(p,q,\varrho)$ such that $(T^{\beta}(z))$ is $a\H^q - \L^2$-bounded.
\end{lem}

\begin{proof} 
If $\varrho \leq p$, then we can simply interpolate and take $\beta = 1$. Henceforth, we assume $p < \varrho$.

Consider a $(\nicefrac{1}{s},\nicefrac{1}{t})$-plane as in Figure~\ref{fig: extra} where $(\nicefrac{1}{s},\nicefrac{1}{t})$ is marked provided $(T(z))$ is $a\H^s -a\H^t$-bounded. The initial configuration are the vertices $A=(\nicefrac{1}{p}, \nicefrac{1}{p})$, $B=(\nicefrac{1}{2}, \nicefrac{1}{2})$ and $C=(\nicefrac{1}{\varrho}, \nicefrac{1}{2})$. By interpolation, we obtain their convex hull, that is to say, the closed triangle $ABC$. 

Boundedness properties for $(T^2(z))$ are visualized in Figure~\ref{fig: extra} as follows: Take a point $X = (\nicefrac{1}{s},\nicefrac{1}{t})$ on $AC$, move to $AB$ on a horizontal line, then move to $AC$ on a vertical line and call that point $X' = (\nicefrac{1}{t}, \nicefrac{1}{t'})$. Then $(T^2(z))$ is $a\H^{s} - a\H^{t'}$-bounded.

\begin{figure}[ht]
	\centering
	\begin{tikzpicture}[scale=5]
	\draw [<->,thick] (0.2,1.2) node (yaxis) [above] {$\frac{1}{t}$}
	|- (1.2,0.2) node (xaxis) [right] {$\frac{1}{s}$};
	\coordinate (B) at (1/2, 1/2);
	\draw (B) node[left] {$B$};
	\fill (B) circle (0.3pt);
	\coordinate (A) at (1,1);
	\draw (A) node[above] {$A$};
	\fill (A) circle (0.3pt);
	\coordinate (C) at (0.8,1/2);
	\draw (C) node[right] {$C$};
	\fill (C) circle (0.3pt);
	\draw (A) -- (B);
	\draw (B) -- (C);
	\draw (C) -- (A);
	\coordinate (X) at (14/15, 5/6);
	\coordinate (Y) at (5/6,5/6);
	\coordinate (Y_1) at (5/6,0);
	\draw  (X) node[right] {$X$};
	\fill (X) circle (0.3pt);
	\draw[dashed] (X) -- (Y);
	\coordinate (X_1) at (intersection of A--C and Y--Y_1);
	\draw[dashed] (Y) -- (X_1);
	\draw (X_1) node[right] {$X'$};
	\fill (X_1) circle (0.3pt);
	\end{tikzpicture}
	\caption{Visualization of the proof of Lemma~\ref{lem: extra}.}
	\label{fig: extra}
\end{figure}

If $\nicefrac{1}{q} \leq \nicefrac{1}{\varrho}$, then $ABC$ contains the point $(\nicefrac{1}{q},\nicefrac{1}{2})$ and we can take $\beta = 1$. Otherwise, the segment $AC$ contains at least one point $X_0$ with abscissa $\nicefrac{1}{q}$. Starting from there, we construct $X_\beta \coloneqq (X_{\beta-1})'$ as above. After a finite number $\beta(p,q,\varrho)$ of steps, $X_\beta$ lies on the segment $BC$ with constant ordinate $\nicefrac{1}{2}$. Hence $(T^\beta(z))$ is $a\H^q - a\H^2$-bounded, that is, $a\H^q - \L^2$-bounded.
\end{proof}
\subsection{Off-diagonal estimates}
\label{subsec: Hp Hq OD}

For Lebesgue spaces we shall make extensive use of off-diagonal estimates.

\begin{defn}
\label{def: p-q OD}
Let $1\leq p \leq q \leq \infty$. An operator family $(T(z))_{z \in \Omega}$ as in \eqref{eq: standard assumptions Hp-Hq} satisfies \emph{$\L^p - \L^q$ off-diagonal estimates}\index{off-diagonal estimates!$\L^p - \L^q$ of order $\gamma$} of order $\gamma>0$ if 
\begin{align*}
\|\ind_F T(z) \ind_E f \|_q \lesssim |z|^{\frac{n}{q}- \frac{n}{p}} \bigg(1+ \frac{\dist(E,F)}{|z|} \bigg)^{-\gamma} \|\ind_E f\|_p
\end{align*}
for all measurable subsets $E,F \subseteq \R^n$, all $z \in \Omega$ and all $f \in \L^p \cap \L^2$. If there are is a constant $c>0$ such that the stronger estimate
\begin{align*}
	\|\ind_F T(z) \ind_E f \|_q \lesssim  |z|^{\frac{n}{q}- \frac{n}{p}} \e^{-c \frac{\dist(E,F)}{|z|}} \|\ind_E f\|_p
\end{align*}
holds, then the family is said to satisfies \emph{off-diagonal estimates of exponential order}\index{off-diagonal estimates!$\L^p - \L^q$ of exponential order}.
\end{defn}

As usual, we shall simply speak of \emph{$\L^p$ off-diagonal estimates} when $p=q$. For $p=q=2$ this notion is consistent with Definition~\ref{def: off-diagonal}. Duality for Lebesgue spaces yields the principle that $(T(z))$ satisfies $\L^p - \L^q$ off-diagonal estimates of order $\gamma$ (resp.\ of exponential order)  if and only if $(T(z)^*)$ satisfies $\L^{q'} - \L^{p'}$ off-diagonal estimates of order $\gamma$ (resp.\ of exponential order). As for composition of off-diagonal estimates, we have the following rule.\index{off-diagonal estimates!composition}

\begin{lem}
\label{lem: OD composition p,q}
Let $1\leq p \leq q \leq \infty$. Let $(T(z))$ and $(S(z))$ be families as in \eqref{eq: standard assumptions Hp-Hq} that are compatible in the sense that $(S(z)T(z))$ is defined. Suppose that they satisfy $\L^p - \L^q$ and $\L^q - \L^r$ off-diagonal estimates of orders $\gamma_T$ and $\gamma_S$, respectively. Then $(S(z)T(z))$ satisfies $\L^p - \L^r$ off-diagonal estimates of order $\gamma_S \wedge \gamma_T$. If the order is exponential for both families, then the same is true for the composition.
\end{lem}

\begin{proof}
Given $E, F \subseteq \R^n$, we put $\dist \coloneqq \dist(E,F)$ and define $G \coloneqq \{x \in \R^n: \dist(x,E) \leq \nicefrac{\dist}{2}\}$. Since we have  $\dist(E,{}^cG) \geq \nicefrac{\dist}{2}$ and $\dist(F,G) \geq \nicefrac{\dist}{2}$, the claim follows on splitting
\begin{align*}
\ind_F S(z)T(z) \ind_E = \ind_F S(z) \ind_G T(z) \ind_E + \ind_F S(z) \ind_{{}^cG}T(z) \ind_E
\end{align*}
and applying $\L^p - \L^q$ and $\L^q - \L^r$ off-diagonal estimates.
\end{proof}

Taking $E=F= \R^n$, we see that $\L^p - \L^q$ off-diagonal estimates are a stronger notion than $\L^p - \L^q$-boundedness, but more is true. This is well-known but we include a proof for convenience.

\begin{lem}
\label{lem: OD implies boundedness}
Let $1 \leq p \leq q \leq \infty$. If an operator family $(T(z))$ as in \eqref{eq: standard assumptions Hp-Hq} satisfies $\L^p - \L^q$ off-diagonal estimates of order $\gamma > n$, then it is $\L^q$-bounded and $\L^p$-bounded.
\end{lem}

\begin{proof}
If $p = \infty$, then $q=\infty$, and $\L^\infty$ off-diagonal estimates imply $\L^\infty$-boundedness. From now on we may assume $p<\infty$.

Let $f \in \L^p$. For fixed $z$, we partition $\R^n$ into closed, axis-parallel cubes $\{Q_k\}_{k \in \IZ^n}$ of sidelength $|z|$ with center $|z|k$. From H\"older's inequality and the assumption we obtain
\begin{align*}
\|T(z) f\|_p^p
&= \sum_{k \in \IZ^n} \|\ind_{Q_k} T(z)f\|_p^p \\
&\leq |z|^{n-\frac{np}{q}} \sum_{k \in \IZ^n} \|\ind_{Q_k} T(z)f\|_q^p \\
&\leq |z|^{n-\frac{np}{q}} \sum_{k \in \IZ^n} \bigg(\sum_{j \in \IZ^n} \|\ind_{Q_k} T(z)\ind_{Q_j}f\|_q\bigg)^p \\
&\leq \sum_{k \in \IZ^n} \bigg(\sum_{j \in \IZ^n} C \Big(1+\frac{\dist(Q_j,Q_k)}{|z|}\bigg)^{-\gamma}\|\ind_{Q_j}f\|_p\bigg)^p.
\end{align*}
Let $|\cdot|_\infty$ be the $\ell^\infty$-norm on $\R^n$ and $\dist_\infty$ the corresponding distance. Then
\begin{align*}
\dist(Q_j,Q_k) 
\geq \dist_\infty(Q_j,Q_k) 
= |z|\max\{|j-k|_\infty-1,0\}.
\end{align*}
Young's convolution inequality yields
\begin{align*}
\|T(z) f\|_p^p
&\leq \sum_{k \in \IZ^n} \bigg(\sum_{j \in \IZ^n} C \big(1+|j-k|_\infty \big)^{-\gamma} \|\ind_{Q_j}f\|_p\bigg)^p \\
&\leq C \bigg(\sum_{k \in \IZ^n} \big(1+|k|_\infty \big)^{-\gamma} \bigg)^p \bigg(\sum_{j \in \IZ^n} \|\ind_{Q_j}f\|_p^p \bigg).
\end{align*}
The sum in $k$ converges since for fixed $m \in \IN$ there are $(2m+1)^n- (2m-1)^n = \mathcal{O}(m^{n-1})$ points $k \in \IZ^n$ with $|k|_\infty = m$. The sum in $j$ equals $\|f\|_p^p$. This proves the $\L^p$-boundedness of $(T(z))$. 

The same argument applies to the dual family, which satisfies $\L^{q'}-\L^{p'}$ off-diagonal estimates of order $\gamma$. This yields $\L^{q}$-boundedness of $(T(z))$.
\end{proof}

\begin{rem}
 \label{rem:OD} 
Re-examining the above proof reveals that $(T(z))$ even satisfies $\L^p$ and $\L^q$ off-diagonal estimates, both of order $\gamma-n$, and that the order is exponential provided that this is the case for the $\L^p-\L^q$ off-diagonal estimates. 

Indeed, assume that $f$ is supported in a set $E$ and that the $\L^p$-norm is taken on a set $F$ with $\dist \coloneqq \dist_\infty(E,F) \geq 4 |z|$. All cubes $Q_k$ and $Q_j$ that are necessary to cover $E$ and $F$, respectively, satisfy $2|k-j|_{\infty}|z| \ge \d$. Consequently, we only need to sum over $k \in \IZ^n$ with $|k| \geq \nicefrac{\dist}{2|z|}$ in the final estimate. This sum is dominated by a multiple of $(1+ \nicefrac{\dist}{|z|})^{-\gamma+n}$. If the order for the $\L^p-\L^q$ off-diagonal estimates is exponential, then we would sum over $\e^{-c|k|_\infty}$ and get control by $\e^{-\frac{c \dist}{2|z|}}$. By duality, the same conclusions are true on $\L^q$.
\end{rem}

The previous lemma provides a means to obtain uniform boundedness in one space from sufficient decay between different spaces. We also need a result of this type for $p<1$.

\begin{lem}
\label{lem: OD Hardy implies boundedness}
Let $(T(z))$ be an operator family as in \eqref{eq: standard assumptions Hp-Hq}. Suppose that $\varrho \in (1_*,1)$ and $q \in (1,\infty)$ are such that $(T(z))$ is $a_1 \H^\varrho - \L^q$-bounded and satisfies $\L^q$ off-diagonal estimates of arbitrarily large order. In addition assume $\int_{\R^n} ((a_2)^{-1} T(z)a_1f)(x) \d x = 0$ for all $z$ and all compactly supported $f \in \L^2$ with integral zero. Then $(T(z))$ is $a_1\H^p - a_2 \H^p$-bounded for every $p \in (\varrho,1]$.
\end{lem}

\begin{proof}
We can assume $a_1=1$ and $a_2 = 1$. Otherwise we replace $T(z)$ with $(a_2)^{-1}T(z)a_1$. Relying on the ($\L^2$-convergent) atomic decomposition for $\H^p\cap \L^2$ (see Section~\ref{subsec: Hardy spaces}) it suffices to show that there is a constant $C$ such that $\|T(z)a\|_{\H^p} \leq C$ for all $z$ and all $\L^\infty$-atoms $a$ for $\H^p$.

\medskip

\noindent \emph{Step~1: Molecular decay}. We show that there exist $C, \eps$ independently of $a$, $z$ and $j \geq 1$ such that
\begin{align}
\label{eq0: OD Hardy implies boundedness}
\|T(z)a\|_{\L^q(C_j(Q))} \leq C (2^{j} \ell(Q))^{\frac{n}{q}-\frac{n}{p}} 2^{-\eps j},
\end{align}
where $Q$ is the cube associated with $a$. For $j=1$ we can simply use $\L^q$-boundedness and H\"older's inequality:  
\begin{align*}
\|T(z)a\|_q 
\leq C\|a\|_q 
\leq C \ell(Q)^{\frac{n}{q}} \|a\|_\infty
\leq C \ell(Q)^{\frac{n}{q}-\frac{n}{p}}.
\end{align*}
For $j \geq 2$ the off-diagonal assumption yields 
\begin{align}
\label{eq1: OD Hardy implies boundedness}
\begin{split}
\|T(z)a\|_{\L^q(C_j(Q))}
&\leq C_\gamma \bigg(\frac{2^{j}\ell(Q)}{|z|} \bigg)^{-\gamma} \|a\|_{\L^q(Q)}\\
&\leq C_\gamma \bigg(\frac{2^{j}\ell(Q)}{|z|} \bigg)^{-\gamma} \ell(Q)^{\frac{n}{q}-\frac{n}{p}}
\end{split}
\end{align}
with $\gamma>0$ at our disposal. Likewise, $\H^\varrho - \L^q$-boundedness yields
\begin{align}
\label{eq2: OD Hardy implies boundedness}
\|T(z)a\|_{\L^q(C_j(Q))}
\leq C |z|^{\frac{n}{q}-\frac{n}{\varrho}} \|a\|_{\H^\varrho}
\leq C |z|^{\frac{n}{q}-\frac{n}{\varrho}}  \ell(Q)^{\frac{n}{\varrho}-\frac{n}{p}},
\end{align}
where in the second step we have used that $\ell(Q)^{n/p-n/\varrho}a$ is an $\L^\infty$-atom for $\H^\varrho$. Now, fix $\delta > 0$ such that $\nicefrac{1}{p}-\nicefrac{1}{q} = (1-2\delta)(\nicefrac{1}{\varrho}-\nicefrac{1}{q})$. This is possible since we have $p>\varrho$. For $|z| \geq 2^{j(1-\delta)} \ell(Q)$ we use \eqref{eq2: OD Hardy implies boundedness} to get
\begin{align*}
\|T(z)a\|_{\L^q(C_j(Q))}
&\leq C 2^{j(1-\delta)(\frac{n}{q}-\frac{n}{\varrho})} \ell(Q)^{\frac{n}{q}-\frac{n}{p}} \\
&= C (2^{j} \ell(Q))^{\frac{n}{q}-\frac{n}{p}} 2^{-\delta(\frac{n}{\varrho}-\frac{n}{q})j},
\end{align*}
whereas for $|z| \leq 2^{j(1-\delta)} \ell(Q)$ we employ \eqref{eq1: OD Hardy implies boundedness} and find
\begin{align*}
\|T(z)a\|_{\L^q(C_j(Q))} \leq C_\gamma 2^{-j \gamma \delta} \ell(Q)^{\frac{n}{q}-\frac{n}{p}}.
\end{align*}
We take $\gamma > \delta^{-1}(\nicefrac{n}{p}-\nicefrac{n}{q})$ to make sure that these bounds take the form \eqref{eq0: OD Hardy implies boundedness}.

\medskip

\noindent \emph{Step~2: Conclusion}. Since $f \coloneqq T(z)a$ has integral zero by assumption, \eqref{eq0: OD Hardy implies boundedness} implies $\|f\|_{\H^p} \leq C$ for some constant independent of $f$. Indeed, this is due to the molecular theory of Taibleson--Weiss~\cite[Thm.~2.9]{Taibleson-Weiss} but we include their argument in our special case in the subsequent lemma.\index{molecule!for $\H^p$}
\end{proof}

\begin{lem}
\label{lem: Taibleson-Weiss}
Let $p \in (1_*,1]$ and $q \in (1,\infty)$. Suppose $f \in \L^2$ has integral zero and satisfies for some $C, \eps > 0$, some cube $Q \subseteq \R^n$ and all $j \geq 1$,
\begin{align*}
\|f\|_{\L^q(C_j(Q))} \leq C (2^{j} \ell(Q))^{\frac{n}{q}-\frac{n}{p}} 2^{-\eps j}.
\end{align*}
Then, there exists a constant $C'$ depending on $C,\eps$ and dimensions, and $\L^q$-atoms $a_j$ for $\H^p$ with support in $C_{j+1}(Q) \cup C_j(Q)$, such that
\begin{align*}
f = \sum_{j=1}^\infty C' 2^{-\eps j} a_j
\end{align*}
with unconditional convergence in $\Lloc^1$. In particular, the sum also converges in $\H^p$ and $\|f\|_{\H^p} \leq \frac{C'}{2^{\eps}-1}$.
\end{lem}

\begin{proof}
The final statement follows from the atomic representation, using the maximal function characterization of $\H^p$ to push the sum in $j$ through the $p$-th power of the $\H^p$-quasinorm, see also \cite[p.106]{Stein93}. To prove the rest, we set 
\begin{align*}
f_j \coloneqq \ind_{C_j(Q)} f, \quad p_j \coloneqq (f)_{C_j(Q)} \ind_{C_j(Q)}.
\end{align*}
Then $f_j - p_j$ has mean value zero and satisfies 
\begin{align*}
\|f_j - p_j\|_q \leq 2 \|f_j\|_q \leq 2C (2^{j} \ell(Q))^{\frac{n}{q}-\frac{n}{p}} 2^{-\eps j}.
\end{align*}
This means that $2^{\eps j} 2^{n/q-n/p} (2C)^{-1} (f_j-p_j)$ is an $\L^q$-atom for $\H^p$. Next, letting $c_j \coloneqq \sum_{k = j}^\infty \int_{C_k(Q)} f \d x$, summation by parts gives a pointwise identity
\begin{equation}
\label{eq1: Taibleson-Weiss}
\sum_{j=1}^\infty p_j 
= \sum_{j=1}^\infty (c_j-c_{j+1}) \frac{\ind_{C_j(Q)}}{|C_j(Q)|}=  \sum_{j=1}^\infty b_j  
\end{equation}
with 
\begin{equation*} b_j \coloneqq c_{j+1} \bigg(\frac{\ind_{C_{j+1}(Q)}}{|C_{j+1}(Q)|} - \frac{\ind_{C_j(Q)}}{|C_j(Q)|} \bigg).
\end{equation*}
There are no boundary terms since we have $c_1 = 0$ and, from the assumption, 
\begin{align*}
|c_j|
\leq \sum_{k=j}^\infty |C_k(Q)|^{1-\frac{1}{q}} \|\ind_{C_k(Q)}f\|_q
&\lesssim  \sum_{k=j}^\infty (2^k \ell(Q))^{n-\frac{n}{p}} 2^{-\eps k}\\
&\simeq (2^j \ell(Q))^{n-\frac{n}{p}} 2^{-\eps j},
\end{align*}
so $\nicefrac{|c_j|}{|C_j(Q)|}$ tends to $0$ as $j \to \infty$. Identity \eqref{eq1: Taibleson-Weiss} holds in $\Lloc^1$ with unconditional convergence because the sums are locally finite since $b_j$ has support in $C_{j+1}(Q) \cup C_j(Q)$. Moreover, $b_j$ has mean value zero and satisfies 
\begin{align*}
\|b_j\|_q \leq |c_{j+1}|\bigg(\frac{|C_{j+1}(Q)|^{\frac{1}{q}}}{|C_{j+1}(Q)|} + \frac{|C_j(Q)|^{\frac{1}{q}}}{|C_j(Q)|} \bigg)\leq C'' (2^{j} \ell(Q))^{\frac{n}{q}-\frac{n}{p}} 2^{-\eps j}.
\end{align*} 
Hence, $2^{\eps j} 4^{n/q-n/p} (C'')^{-1} b_j$ is an $\L^q$-atom for $\H^p$ and $f = \sum_{j=1}^\infty f_j - p_j + b_j$ is the representation we are looking for.
\end{proof}
\subsection{Interpolation principles}
\label{subsec: Hp Hq interpolation}

We continue with interpolation properties. Our main tool will be the Stein interpolation theorem\index{Theorem!Stein interpolation}, which we state in an abstract version due to Voigt~\cite{Voigt92}. In the following we work on the strip $S \coloneqq \{z \in \IC: 0 \leq \Re z \leq 1\}$.

\begin{prop}[{\cite{Voigt92}}]
\label{prop: Stein interpolation}
Let $(X_0,X_1)$ and $(Y_0,Y_1)$ be two interpolation couples of Banach spaces and let $Z$ be a dense subspace of $X_0 \cap X_1$. Let $(T(z))_{z \in S}$ be a family of linear mappings $Z \to Y_0 + Y_1$ with the following properties for all $f \in Z$:
\begin{enumerate}
	\item The function $T(\cdot)f: S \to Y_0 + Y_1$ is continuous, bounded and holomorphic in the interior of $S$.
	\item For $j=0,1$ the restriction $T(\cdot)f: j + \i \R \to Y_j$ is continuous and there is a constant $M_j$ that does not depend on $f$ such that
	\begin{align*}
	\sup_{t \in \R} \|T(j + \i t)f\|_{Y_j} \leq M_j \|f\|_{X_j}.
	\end{align*}
\end{enumerate}
Then for all $\theta \in (0,1)$ and all $f \in Z$,
\begin{align*}
 \|T(\theta)f\|_{[Y_0,Y_1]_\theta} \leq M_0^{1-\theta} M_1^\theta \|f\|_{[X_0,X_1]_\theta}.
\end{align*}
\end{prop}

\begin{rem}
\label{rem: Stein interpolation}
The classical version of the theorem is when $X_j$ and $Y_j$ are $\L^p$-spaces, $1 \leq p \leq \infty$, and $Z$ is the space of step functions, \cite[Thm.~1.3.4]{Grafakos}. Then continuity is not required in (ii) and in (i) it suffices to assume that for all $f, g \in Z$ and all $z \in S$ the integral $\int_{\R^n} T(z)f \cdot \cl{g} \, \d x$ converges absolutely and defines a bounded and continuous function of $z$ that is holomorphic in the interior of $S$. For example, it suffices that $T(\cdot): S \to \L^2$ is bounded, continuous and holomorphic in the interior. Such weakening of assumptions is not possible for general interpolation couples~\cite{Cwikel-Janson84}.
\end{rem}

As a first application we prove the following\index{off-diagonal estimates!interpolation}
 
\begin{lem}
\label{lem: OD extrapolation to sectors}
Let $p_0,q_0, p_1, q_1 \in [1,\infty]$, $p_i \leq q_i$ and $\omega \in (0,\pi)$. Let $(T(z))_{z \in \S_\omega^+}$ be a uniformly bounded family on $\L^2$ as in \eqref{eq: standard assumptions Hp-Hq} that depends holomorphically on $z$. Let $\theta \in (0,1)$ and set 
\begin{align*}
	p_\theta \coloneqq [p_0,p_1]_\theta \quad \text{and} \quad q_\theta \coloneqq [q_0,q_1]_\theta.
\end{align*}
\begin{enumerate}
	\item If $(T(z))_{z \in (0,\infty)}$ is $\L^{p_0} - \L^{q_0}$-bounded and $(T(z))_{z \in \S_\omega^+}$ is $\L^{p_1} - \L^{q_1}$-bounded, then $(T(z))_{z \in \S_{\theta \omega}^+}$ is $\L^{p_\theta} - \L^{q_\theta}$-bounded 
	\item If $(T(z))_{z \in (0,\infty)}$ is $\L^{p_0} - \L^{q_0}$-bounded and $(T(z))_{z \in \S_\omega^+}$ satisfies $\L^{p_1} - \L^{q_1}$ off-diagonal estimates of order $\gamma$, then $(T(z))_{z \in \S_{\theta \omega}^+}$ satisfies $\L^{p_\theta} - \L^{q_\theta}$ off-diagonal estimates of order $\theta \gamma$.
	\item If $(T(z))_{z \in (0,\infty)}$ satisfies $\L^{p_0} - \L^{q_0}$ off-diagonal estimates of order $\gamma$ and $(T(z))_{z \in \S_\omega^+}$ is $\L^{p_1} - \L^{q_1}$-bounded, then $(T(z))_{z \in \S_{\theta \omega}^+}$ satisfies $\L^{p_\theta} - \L^{q_\theta}$ off-diagonal estimates of order $(1-\theta) \gamma$.
\end{enumerate}
Exponential order in the assumptions leads to exponential order in the conclusion with the decay parameter $c$ changed accordingly.
\end{lem}

\begin{proof}
We prove the three statements first and indicate the necessary modifications in the case of exponential order afterwards.

\medskip

\noindent \emph{Proof of (i),(ii),(iii).}
We begin with part (ii). We fix $\nu \in (-\omega, \omega)$ and $r>0$. Then we fix measurable sets $E, F \subseteq \R^n$, set $\d \coloneqq \d(E,F)$ and consider the family $S(z) \coloneqq \e^{(z-\theta)^2} \ind_F T(r \e^{\i \nu z}) \ind_E$. This family is uniformly bounded on $\L^2$ and holomorphic in an open neighborhood of the strip $0 \leq \Re z \leq 1$. By assumption we have for all $t \in \R$ and all step functions $f$,
\begin{align*}
\|S(\i t) f\|_{q_0} 
&\leq  C_0 |r \e^{- \nu t}|^{\frac{n}{q_0}-\frac{n}{p_0}} \e^{1-t^2} \|f\|_{p_0} \\
&\leq  C_0 (r \e^{- \omega |t|})^{\frac{n}{q_0}-\frac{n}{p_0}} \e^{1-t^2} \|f\|_{p_0} \\
\intertext{and}
\|S(1+\i t) f\|_{q_1} 
&\leq  C_1 |r \e^{- \nu t}|^{\frac{n}{q_1}-\frac{n}{p_1}}  \Big(1+ \frac{\dist}{|r \e^{- \nu t}|} \Big)^{-\gamma}  \e^{1-t^2} \|f\|_{p_1} \\
&\leq C_1 (r \e^{- \omega |t|})^{\frac{n}{q_1}-\frac{n}{p_1}}  \Big(1+ \frac{\dist}{r \e^{\omega |t|}} \Big)^{-\gamma}  \e^{1-t^2} \|f\|_{p_1}.
\end{align*}
We use $1+\nicefrac{\d}{r \e^{\omega |t|}} \geq \e^{-\omega |t|}(1+\nicefrac{\d}{r})$ in the second line and that the additional factor $\e^{-t^2}$ acts in our favor, in order to give
\begin{align*}
\|S(\i t) f\|_{q_0} &
\leq  M_0 r^{\frac{n}{q_0}-\frac{n}{p_0}} \|f\|_{p_0}, \\
\|S(1+\i t) f\|_{q_1} &
\leq  M_1 r^{\frac{n}{q_1}-\frac{n}{p_1}}  \Big(1+ \frac{\dist}{r} \Big)^{-\gamma} \|f\|_{p_1},
\end{align*}
where the $M_j$ are still also independent of $r$, $\nu$ and $E,F$. Stein interpolation yields
\begin{align*}
\|S(\theta) f\|_{q_\theta} &\leq  M_0^{1-\theta} M_1^\theta r^{\frac{n}{q_\theta}-\frac{n}{p_\theta}} \Big(1+ \frac{\dist}{r} \Big)^{-\theta \gamma} \|f\|_{p_\theta}.
\end{align*}
Since we have $S(\theta)f = \ind_F T(r \e^{\i \nu \theta}) \ind_E f$, this estimates means that $T(z)$ satisfies $\L^{p_\theta} - \L^{q_\theta}$ off-diagonal estimates of order $\theta \gamma$ for $z \in \S_{\theta \omega}^+$. 

The proof of part (iii) is exactly the same except that now the estimate for $S(\i t)$ comes with decay.
 
The proof of part (i) does not need the sets $E,F$ and uses the same interpolation argument for $z \mapsto \e^{(z-\theta)^2} T(r \e^{\i \nu z})$.
 
\medskip

\noindent \emph{Proof of (ii), (iii) in the exponential setting.} The proof of (ii) in the case of exponential order follows the same idea but is slightly more technical. We fix again $E,F \subseteq \R^n$ and prove the estimate at $z = r \e^{\i \nu \theta}$. By (i) we already know that $(T(z))_{z \in \S_{\theta \omega}^+}$ is $\L^{p_\theta}-\L^{q_\theta}$-bounded. Hence, it suffices to prove the exponential estimate for $\nicefrac{\d}{r}$ large, say $\nicefrac{\d}{r} \geq \nicefrac{\omega}{\eps}(\nicefrac{n}{p_i} - \nicefrac{n}{q_i})$ for $i=0,1$, where $\eps>0$ will be chosen small later on in dependence of $\omega$ and $\theta$. Now, we set
\begin{align*}
	S(z) \coloneqq \e^{\eps \frac{\d}{r} (z-\theta)^2} \ind_F T(r \e^{\i \nu z}) \ind_E
\end{align*}
and obtain as before
\begin{align}
\label{eq1: OD extrapolation to sectors}
\begin{split}
\|S(\i t)f\|_{q_0} 
&\leq C_0 (r \e^{-\omega |t|})^{\frac{n}{q_0} - \frac{n}{p_0}} \e^{\eps \frac{\d}{r} (1-t^2)} \|f\|_{p_0}\\
&\leq C_0 r^{\frac{n}{q_0} - \frac{n}{p_0}} \e^{2 \eps  \frac{\d}{r}} \|f\|_{p_0},
\end{split}
\end{align}
where we have used the lower bound on $\nicefrac{\d}{r}$ and then that $1+|t|-t^2 \leq 2$ in the second step. Similarly, but now using the off-diagonal estimates of exponential order, we get
\begin{align*}
\|S(1+\i t)f\|_{q_1} 
&\leq C_1 (r \e^{-\omega |t|})^{\frac{n}{q_1} - \frac{n}{p_1}} \e^{-c_1 \frac{\d}{r \e^{\omega |t|}}} \e^{\eps \frac{\d}{r} (1-t^2)} \|f\|_{p_1}\\
&\leq C_1 r^{\frac{n}{q_1} - \frac{n}{p_1}} \e^{\frac{\d}{r}(-c_1 \e^{-\omega |t|} + \eps (1+|t|-t^2))} \|f\|_{p_1}.
\end{align*}
At this point, we claim that we can pick $\eps > 0$ such that 
\begin{align*}
	-c_1 \e^{-\omega |t|} + \eps (1+|t|-t^2)) \leq - \frac{2}{\theta} \eps \qquad (t \in \R).
\end{align*}
Indeed, we have $1+|t|-t^2 \leq \nicefrac{3-t^2}{2} \leq - \nicefrac{2}{\theta}$ if $|t| \geq \sqrt{3+ \nicefrac{4}{\theta}} \eqqcolon t_\theta$ and if $|t| \leq t_\theta$, then the left-hand side above is controlled by  $-c_1 \e^{-\omega t_\theta} + \eps(1+t_\theta) \leq - \nicefrac{2\eps }{\theta}$ provided $\eps$ is small enough. Altogether, with $\eps$ now being fixed, we conclude
\begin{align}
\label{eq2: OD extrapolation to sectors}
	\|S(1+\i t)f\|_{q_1} 
	\leq C_1 r^{\frac{n}{q_1} - \frac{n}{p_1}} \e^{- \frac{2}{\theta} \eps \frac{\d}{r}} \|f\|_{p_1}.
\end{align}
Using \eqref{eq1: OD extrapolation to sectors} and \eqref{eq2: OD extrapolation to sectors}, we get by Stein interpolation that
\begin{align*}
	\|S(\theta)f\|_{q_\theta} 
	\leq C_0^{1-\theta} C_1^\theta r^{\frac{n}{q_\theta} - \frac{n}{p_\theta}} \e^{-2 \theta \eps \frac{\d}{r}} \|f\|_{p_\theta},
\end{align*}
which is the required estimate since $S(\theta) = \ind_F T(r \e^{\i \nu \theta}) \ind_E$.

The proof of (iii) in the exponential setting is again the same except that $\theta$ should be replaced by $1-\theta$ in the construction of $\eps$ since now the estimate for $S(\i t)$ comes with decay.
\end{proof}

For families with constant domain in $z$ the argument is much simpler.\index{off-diagonal estimates!interpolation} 

\begin{lem}
\label{lem: OD interpolation}
Let $p_0,q_0, p_1, q_1 \in [1,\infty]$ with $p_i \leq q_i$ and suppose that a family as in \eqref{eq: standard assumptions Hp-Hq} is $\L^{p_0} - \L^{q_0}$-bounded and satisfies $\L^{p_1} - \L^{q_1}$ off-diagonal estimates of order $\gamma$ (of exponential order). Then for each $\theta \in (0,1)$ it satisfies  $\L^{[p_0,p_1]_\theta} - \L^{[q_0,q_1]_\theta}$ off-diagonal estimates of order $\theta \gamma$ (of exponential order).
\end{lem}

\begin{proof}
Apply the Riesz--Thorin theorem to the operator $\ind_F T(z) \ind_E$, where $z$ and the measurable sets $E,F \subseteq \R^n$ are fixed.
\end{proof}
\subsection{Applications to the functional calculus}
\label{subsec: Hp Hq resolvent}

We turn to the more specific setting that the family $(T(z))$ is modeled after the resolvents of a sectorial operator. In this section, we assume that 
\begin{align}
\label{eq: ass Hp Hq resolvent}
\begin{minipage}{0.89\linewidth}
\begin{itemize}
	\item $T$ is a sectorial operator on $\L^2(\R^n; V)$ of some angle $\omega \in [0,\pi)$, where $V$ is a finite-dimensional Hilbert space,
	\item $((1+t^2 T)^{-1})_{t >0}$ satisfies $\L^2$ off-diagonal estimates of arbitrarily large order.
\end{itemize}
\end{minipage}
\end{align}

\begin{lem}
\label{lem: Lp bounds for resolvents on sector}
Let $p \in (1,\infty)$ be such that $((1+t^2T)^{-1})_{t>0}$ is $\L^p$-bounded. Let $\theta \in (0,1]$. Then for every $\mu \in (0, \nicefrac{\theta (\pi - \omega)}{2})$ the family $((1+z^2T)^{-1})_{z \in \S_{\mu}^+}$ satisfies $\L^{[p,2]_\theta}$ off-diagonal estimates of arbitrarily large order. 
\end{lem}

\begin{proof}
The resolvent $z \mapsto (1+z^2 T)^{-1}$ on $\L^2$ is a bounded holomorphic function on $\S_\mu^+$ for any $\mu \in (0, \nicefrac{(\pi - \omega)}{2})$. We apply Lemma~\ref{lem: OD extrapolation to sectors} twice. 

First, interpolation between the $\L^2$-bounds on sectors and the $\L^2$ off-diagonal estimates on the positive real axis yields $\L^2$ off-diagonal estimates of arbitrarily large order on $\S_{\mu}^+$ for any $\mu \in (0, \nicefrac{(\pi - \omega)}{2})$. Second, interpolation between the $\L^2$ off-diagonal estimates on sectors and the $\L^p$-bounds on the positive real axis yields the claim.
\end{proof}

We obtain off-diagonal estimates for the functional calculus similar to \cite[Part~II]{AusSta}. In applications we usually work with holomorphic functions that are in the respective classes on any sector and the technical conditions on the angles can be ignored. On the other hand, the order of off-diagonal decay is of utmost importance: It is mainly the decay of $\psi$ at $z=0$, quantified by the classes $\Psi_\sigma^\tau$ from Section~\ref{subsec: functional calculi}, that limits the available off-diagonal decay for $(\psi(t^2T))_{t>0}$.\index{off-diagonal estimates!for the functional calculus}

\begin{lem}
\label{lem: functional calculus bounds from J(L) abstract}
Let $p \in (1,\infty)$ be such that $((1+t^2T)^{-1})_{t>0}$ is $\L^p$-bounded. Let $\theta \in (0,1]$, put $q \coloneqq [p,2]_\theta$ and fix an angle $\mu \in (0, \nicefrac{\theta(\pi - \omega)}{2})$. Let $\sigma, \tau > 0$ and $\psi \in \Psi_\sigma^\tau(\S_{\pi - 2 \mu}^+)$. Then the following estimates hold.

\begin{enumerate}
\item Let $(\eta(t))_{t>0}$ be a continuous and uniformly bounded family of functions in $\H^\infty(\S_{\pi - 2 \mu}^+)$. Then for all measurable sets $E, F \subseteq \R^n$, all $t>0$ and all $f \in \L^q \cap \L^2$,
\begin{align*}
 \qquad \|\ind_F \eta(t)(T) \psi(t^2 T) \ind_E f\|_q \lesssim \|\eta\| \|\psi\|_{\sigma,\tau, \mu} \bigg(1 + \frac{\d(E,F)}{t} \bigg)^{-2 \sigma} \|\ind_E f\|_q.
\end{align*}
The norms are $\|\psi\|_{\sigma,\tau,\mu} \coloneqq \sup_{z \in \S_{\pi - 2\mu}^+} \nicefrac{|\psi(z)|}{(|z|^\sigma \wedge |z|^{-\tau})}$ and $\|\eta\| \coloneqq \sup_{t>0} \|\eta(t)\|_\infty$.
\index{$\lVert \, \cdot \, \rVert_{\sigma, \tau, \mu}$ (norm on $\Psi_{\sigma}^\tau(\S_{\pi-2\mu}^+)$)}

\item Furthermore, if $\eta(t)(z) = \varphi(t^2z)$ for some $\varphi \in \Psi_\sigma^0(\S_{\pi - 2\mu}^+)$, then for all $0<r \leq t$ and with the same dependencies,
\begin{align*}
\qquad  \|\ind_F \varphi(r^2 T) \psi(t^2 T) \ind_E f\|_q \lesssim  \|\psi\|_{\sigma,\tau,\mu} \|\varphi\|_{\sigma,0,\mu} \bigg(1 + \frac{\d(E,F)}{r} \bigg)^{-2\sigma} \|\ind_E f\|_q.
\end{align*}

\item Finally for each $\gamma \in [0, \sigma] \cap  [0,\tau)$, it follows for all $r>0$ with the same dependencies,
\begin{align*}
\qquad \qquad \|\varphi(r^2 T) \psi(t^2 T)f \|_q \lesssim  \|\psi\|_{\sigma,\tau,\mu} \|\varphi\|_{\sigma,0,\mu} \bigg(\frac{r^2}{t^2}\bigg)^{\gamma} \|f\|_q.
\end{align*} 
\end{enumerate}
\end{lem}

\begin{proof}
Throughout, let $\|f\|_q = 1$. We pick an angle $\nu \in (\mu, \nicefrac{\theta(\pi - \omega)}{2})$. By Lemma~\ref{lem: Lp bounds for resolvents on sector} we have $\L^q$ off-diagonal estimates of arbitrarily large order for the resolvents for $z \in \cl{\S_{\nu}^+}$. Here, we use the order $2 \sigma +1$.

We begin with the first estimate and put $X \coloneqq \nicefrac{\d(E,F)}{t}$.  Since
\begin{align}
\label{eq1: functional calculus bounds from J(L) abstract}
\begin{split}
\eta(t)(T) \psi(t^2 T)
&= \frac{1}{2 \pi \i} \int_{\bd \S_{\pi - 2\nu}^+} \eta(t)(z) \psi(t^2z) (z-T)^{-1}  \, \d z\\
&= \frac{1}{2 \pi \i} \int_{\bd \S_{\pi - 2\nu}^+}  \eta(t)(z t^{-2}) \psi(z) (1-t^2z^{-1}T)^{-1} \, \frac{\d z}{z},
\end{split}
\end{align}
where $(-t^2 z^{-1})^{1/2} \in \cl{\S_\nu^+}$, we obtain
\begin{align}
\label{eq2: functional calculus bounds from J(L) abstract}
\begin{split}
\|\ind_F \eta(t)(T) &\psi(t^2 T) \ind_E f\|_q \\
&\lesssim \int_{\bd \S_{\pi - 2 \nu}^+} \|\eta\| \frac{|\psi(z)|}{(1+|z|^{1/2}X)^{2\sigma +1}} \, \frac{\d |z|}{|z|} \\
&\leq \|\eta\| \|\psi\|_{\sigma,\tau,\mu} \int_{\bd \S_{\pi - 2 \nu}^+}  \frac{|z|^\sigma \wedge |z|^{-\tau}}{(1+|z|^{1/2}X)^{2\sigma +1}} \, \frac{\d |z|}{|z|}.
\end{split}
\end{align}
In the case $X \leq 1$, we minimize the denominator by $1$ to derive the desirable bound
\begin{align*}
\|\ind_F \eta(t)(T) \psi(t^2 T) \ind_E f\|_q 
\lesssim \|\eta\| \|\psi\|_{\sigma,\tau,\mu}.
\end{align*}
In the case $X \geq 1$, we split the integral at $|z| = X^{-2}$ to give the desirable bound
\begin{align*}
\|\ind_E \eta(t)(T) &\psi(t^2 T) \ind_F f\|_q \\
&\lesssim \|\eta\|  \|\psi\|_{\sigma,\tau,\mu}  \bigg(\int_{0}^{X^{-2}} |z|^\sigma \, \frac{\d |z|}{|z|} 
+ \int_{X^{-2}}^\infty  \frac{|z|^\sigma}{(|z|^{1/2}X)^{2\sigma +1}} \, \frac{\d |z|}{|z|} \bigg)\\
&\lesssim \|\eta\| \|\psi\|_{\sigma,\tau,\mu} X^{-2 \sigma} .
\end{align*}
This completes the proof of (i).

Turning to the second estimate, we take $\eta(r)(z) = \varphi(r^2z)$ in \eqref{eq1: functional calculus bounds from J(L) abstract} and change variables to 
\begin{align*}
\varphi(r^2T) \psi(t^2 T) = \frac{1}{2 \pi \i} \int_{\bd \S_{\pi - 2 \nu}^+} \varphi(z) \psi(t^2 r^{-2}z) (1-r^2z^{-1}T)^{-1} \, \frac{\d z}{z}.
\end{align*}
This time we set $X \coloneqq \nicefrac{\d(E,F)}{r}$ and obtain 
\begin{align}
\label{eq3: functional calculus bounds from J(L) abstract}
\begin{split}
\|\ind_F \varphi(r^2T) \psi(t^2 T) \ind_E f\|_q
\lesssim \int_{\bd \S_{\pi - 2 \nu}^+} \frac{|\varphi(z)||\psi(t^2 r^{-2} z)|}{(1+|z|^{1/2}X)^{2 \sigma +1}} \, \frac{\d |z|}{|z|}.
\end{split}
\end{align}
The important observation is that
\begin{align}
|\varphi(z)| \leq \|\varphi\|_{\sigma,0,\mu}  (|z|^\sigma \wedge 1) 
\end{align}
and, since $r \leq t$,
\begin{align*}
|\psi(t^2 r^{-2} z)| 
\leq \|\psi\|_\infty \wedge (\|\psi\|_{\sigma,\tau,\mu} |t^2 r^{-2} z|^{-\tau})
\leq \|\psi\|_{\sigma,\tau,\mu}  (1 \wedge |z|^{-\tau}),
\end{align*}
so that
\begin{align*}
|\varphi(z)||\psi(t^2 r^{-2} z)| \leq \|\psi\|_{\sigma,\tau,\mu} \|\varphi\|_{\sigma,0,\mu} (|z|^\sigma \wedge |z|^{-\tau}).
\end{align*}
Thus, we can bound the right-hand side in \eqref{eq3: functional calculus bounds from J(L) abstract} by the same parameter integral that already appeared on the far right in \eqref{eq2: functional calculus bounds from J(L) abstract} and get the same bound $(1+X)^{-2\sigma}$ for the integral. Now, (ii) follows.

As for (iii), we first argue as in (ii) with $E = F = \R^n$ and $X=0$ to obtain
\begin{align*}
\|\varphi(r^2T) &\psi(t^2 T) f\|_q \\
&\lesssim \|\psi\|_{\sigma,\tau,\mu} \|\varphi\|_{\sigma,0,\mu} \int_0^\infty (1 \wedge |z|^\sigma) \big(|t^2 r^{-2} z|^\sigma \wedge |t^2 r^{-2} z|^{-\tau} \big) \, \frac{\d |z|}{|z|}.
\end{align*}
Using $(1 \wedge |z|^\sigma) \leq |z|^\gamma$ in order to get a homogeneous estimate and changing variables, we conclude
\begin{align*}
\|\varphi(r^2T) \psi(t^2 T) f\|_q
\lesssim \|\psi\|_{\sigma,\tau,\mu} \|\varphi\|_{\sigma,0,\mu} \bigg(\frac{r^2}{t^2}\bigg)^{\gamma} \int_0^\infty |z|^\gamma \big(|z|^\sigma \wedge |z|^{-\tau} \big) \, \frac{\d |z|}{|z|}
\end{align*}
and the remaining integral is finite since we assume $0 \leq \gamma < \tau$.
\end{proof}

The decay of $\psi$ at the origin can be replaced by the assumption that $\psi(z)$ has a limit as $|z| \to 0$ with order of convergence $\mathcal{O}(|z|^\sigma)$ for some $\sigma > 0$. The exemplary result of this type is as follows. The obtained order of decay is optimal and already attained when $T=-\Delta_{x}$.

\begin{cor}
\label{cor: Poisson semigroup Lp OD}
Let $p \in (1,\infty)$ be such that $((1+t^2T)^{-1})_{t>0}$ is $\L^p$-bounded and let $\theta \in (0,1)$. Then $(\e^{-t T^{1/2}})_{t>0}$ satisfies $\L^{[p,2]_\theta}$ off-diagonal estimates of order $1$. 
\end{cor}

\begin{proof}
This is a consequence of the preceding two lemmata since we can write $\e^{-z^{1/2}} = \psi(z) + (1+z)^{-1}$ with $\psi \in \Psi_{1/2}^1$ on any sector.
\end{proof}
\section{Conservation properties}
\label{sec: conservation}

\noindent In order to extend the operator theory for $L$ to Hardy spaces, we need to guarantee that certain operators $f(L)$ preserve vanishing zeroth moments or have the \emph{conservation property}\index{conservation property} $f(L)c = c$ whenever $c$ is a constant. In absence of integral kernels, the action of such operators on constants is explained via off-diagonal estimates as follows.\index{off-diagonal estimates!operator extensions by} 

\begin{prop}
\label{prop: general OD extension}	
Let $T$ be a bounded linear operator on $\L^2(\R^n; V)$, where $V$ is a finite dimensional Hilbert space. If $T$ satisfies $\L^p$ off-diagonal estimates of order $\gamma > \nicefrac{n}{p}$ for some $p \in [2,\infty)$, then $T$ can be extended to a bounded operator $\L^\infty(\R^n; V) \to \Lloc^p(\R^n; V)$ via 
\begin{align}
\label{eq: OD definition on Linfty}
	T f \coloneqq \sum_{j=1}^\infty T(\ind_{C_j(B(0,1))} f).
\end{align}
Moreover, if $(\eta_j) \subseteq \L^\infty(\R^n;\IC)$ is a  family such that
\begin{align}
	\label{eq: conservation family}
	\begin{minipage}{0.89\linewidth}
	\begin{itemize}
	\item $\displaystyle \sup_j \|\eta_j\|_\infty < \infty$,
	\item $\displaystyle \sum_{j=1}^\infty \eta_j(x) = 1$ for a.e.\ $x \in \R^n$,
	\item $\eta_j$ has compact support, which for some $C,c$ and all sufficiently large $j$ is contained in $B(0, C2^j) \setminus B(0, c 2^j)$,
	\end{itemize}
 	\end{minipage}
\end{align}
then
\begin{align*}
	Tf = \sum_{j=1}^\infty T(\eta_j f),
\end{align*}
where the right-hand side converges in $\Lloc^p(\R^n; V)$ and in particular in $\Lloc^2(\R^n; V)$.
\end{prop}

\begin{rem}
\label{rem: general OD extension}
A particular example for a family with the required properties is $\eta_j = \ind_{C_j(B)}$ for an arbitrary ball (or cube) $B \subseteq \R^n$. 
\end{rem}

\begin{proof}
We put $B \coloneqq B(0,1)$ and fix any compact set $K \subseteq \R^n$. For all large enough $j$ we have $\dist(K, C_j(B)) \geq  2^{j-1}$ and therefore
\begin{align*}
	\|T (\ind_{C_j(B)} f) \|_{\L^p(K)} 
	&\lesssim 2^{-j \gamma} \|f\|_{\L^p(C_j(B))} \\
	& \lesssim 2^{j (\frac{n}{p} - \gamma)} \|f\|_\infty.
\end{align*}
Hence, the series on the right-hand side of \eqref{eq: OD definition on Linfty} converges absolutely in $\L^p(K)$ and the limit satisfies $\|Tf\|_{\L^p(K)} \leq C_K \|f\|_\infty$ for a constant $C_K$ that depends on $K$ but not on $f$.

Next, we pick an integer $j_0 \geq 1$ such that $c2^{j_0} \geq 1$ and therefore $2^J B \subseteq B(0, c 2^{J+j_0})$ for all $J \geq 1$. If $J$ is large enough so that the annular support of $\eta_j$ is granted, then $\sum_{j=1}^J \ind_{C_j(B)} - \sum_{j=1}^{J+ j_0} \eta_j$ vanishes on $2^{J+1}B$, has support in $C' 2^J B$ for some $C'$ that does not depend on  $J$ and is uniformly bounded. The off-diagonal bounds yield again
\begin{align*}
	\bigg\|\sum_{j=1}^J T(\ind_{C_j(B)}f) - \sum_{j=1}^{J+j_0} T(\eta_j f) \bigg\|_{\L^p(K)} 
	\lesssim 2^{J(\frac{n}{p}-\gamma)} \|f\|_\infty,
\end{align*}
which shows that $\sum_{j=1}^\infty T(\eta_j f)$ converges to $Tf$ in $\L^p(K)$. 
\end{proof}

We begin with the conservation property\index{conservation property!for $BD$} for the resolvents of the perturbed Dirac operator $BD$ that has appeared implicitly in several earlier works~\cite{AusSta, R}. The proof relies on the cancellation property $Dc = 0$ for constants $c$ (where $D$ is understood in the sense of distributions).

\begin{prop}
\label{prop: conservation BD}
If $\alpha \in \IN$ and $z \in \S_{\pi/2 - \omega_{BD}}$, then for all $c \in \IC^{m} \times \IC^{mn}$,
\begin{align*}
 (1+\i z BD)^{-\alpha}c = c = (1+z^2 (BD)^2)^{-\alpha}c.
\end{align*}
\end{prop}

\begin{proof}
Let $R>0$ and $(\eta_j)$ be a smooth partition of unity on $\R^n$ subordinate to the sets 
\begin{align*}
	D_1 \coloneqq B(0,4R), \quad D_j \coloneqq B(0,2^{j+1}R)\setminus B(0,2^{j-1}R) \quad (j \geq 2),
\end{align*}
such that $\|\eta_j\|_\infty + 2^j R \|\nabla_x \eta_j\|_\infty \leq C$ for a dimensional constant $C$. 

We begin with the resolvents of $BD$, which satisfy $\L^2$ off-diagonal estimates of arbitrarily large order by Proposition~\ref{prop: OD for Dirac} and composition. According to Proposition~\ref{prop: general OD extension}  we can write
\begin{align*}
	(1+\i z BD)^{-\alpha} c 
	= \sum_{j=1}^\infty (1+\i z BD)^{-\alpha} (\eta_j c),
\end{align*}
so that
\begin{align*}
	(1+ \i z BD)^{-\alpha+1}c - (1+ \i z BD)^{-\alpha}c 
	= \sum_{j=1}^\infty \i z (1+ \i z BD)^{-\alpha} BD(\eta_j c),
\end{align*}
where we set $(1+\i z BD)^{0}c \coloneqq c$ and used $\eta_j c \in \dom(D) = \dom(BD)$. Now, $BD(\eta_jc)$ has support in $B(0,2^{j+1}R)\setminus B(0,2^{j-1}R)$ also for $j = 1$ and satisfies $\|BD(\eta_j c)\|_\infty \leq C |c| \|B\|_\infty (2^jR)^{-1}$. The off-diagonal estimates of order $\nicefrac{n}{2}$ yield
\begin{align*}
	\|(1+ \i z BD)^{-\alpha+1}c -(1+\i z BD)^{-\alpha }c \|_{\L^2(B(0, R/2))} 
	\lesssim R^{\frac{n}{2}-\gamma-1} \sum_{j=1}^\infty 2^{j(\frac{n}{2}-\gamma-1)}
\end{align*}
with an implicit constant that is independent of $R$. Sending $R \to \infty$ gives $(1+ \i z BD)^{-\alpha+1}c = (1+\i z BD)^{-\alpha }c$. Since $(1+\i z BD)^0c=c$, we conclude $(1+\i z BD)^{-\alpha }c = c$ for all $\alpha$.

The argument for the resolvents of $(BD)^2$ is identical and draws upon the identity
\begin{align*}
	(1+ z^2 &(BD)^2)^{-\alpha+1}c - (1+ z^2 (BD)^2)^{-\alpha}c \\
	&= \sum_{j=1}^\infty z^2 BD (1+ z^2 (BD)^2)^{-\alpha} BD(\eta_j c).
\end{align*}
The off-diagonal decay for $z^2 BD (1+ z^2 (BD)^2)^{-\alpha}$ follows again by composition since this operator can be written as
\begin{align*}
	-\frac{ \i z}{2} \Big((1-\i z BD)^{-1} - (1+ \i z BD)^{-1} \Big) (1+ z^2 (BD)^2)^{-\alpha+1}. &\qedhere
\end{align*}
\end{proof}

As a corollary we obtain the conservation property for the second-order operator $L$. The reader can refer to \cite[Sec.~4.4]{Ouhabaz} and references therein for related conservation properties in the realm of semigroups.\index{conservation property!for resolvents of $L$}

\begin{cor}
\label{cor: conservation}
Let $\alpha \in \IN$ and $z \in \S_{(\pi-\omega_L)/2}^+$. Let $c \in \IC^m $ and let $f \in \L^2$ have compact support. Then one has the conservation formula
\begin{align*}
	(1+z^2 L)^{-\alpha} c = c
\end{align*}
and its dual version
\begin{align*}
	\int_{\R^n} a(1+z^2L)^{-\alpha} (a^{-1}f) \, \d x = \int_{\R^n} f \, \d x.
\end{align*}
\end{cor}

\begin{proof}
Since $f$ has compact support, we obtain from off-diagonal estimates of order $\gamma > \nicefrac{n}{2}$ that $a(1+z^2L)^{-\alpha} (a^{-1}f) \in \L^1$, locally uniformly in $z$. The left-hand sides are holomorphic functions of $z$ (valued in $\Lloc^2$ and $\IC^m$, respectively). Hence, it suffices to argue for $z=t \in (0,\infty)$. We have
\begin{align*}
	(1+t^2 (BD)^2)^{-\alpha} = \begin{bmatrix} (1+t^2 L)^{-\alpha} & 0 \\ 0 & (1+t^2 M)^{-\alpha} \end{bmatrix},
\end{align*}
so the fist claim follows from the conservation property for $BD$. As $(a^*)^{-1} L^* a^*$ belongs to the same class as $L$, we also get
\begin{align*}
	\int_{\R^n} a(1+t^2L)^{-\alpha} (a^{-1}f) \cdot \cl{c} \, \d x
	&= \int_{\R^n} f \cdot \cl{(1+t^2 (a^*)^{-1}L^*a^*)^{-\alpha}c} \, \d x\\
	&= \int_{\R^n} f \cdot \cl{c} \, \d x
\end{align*}
and since $c \in \IC^m$ is arbitrary, the second claim follows.
\end{proof}

We turn to more general operators in the functional calculus. In view of Lemma~\ref{lem: functional calculus bounds from J(L) abstract} the decay of the auxiliary function at the origin limits the available off-diagonal decay and hence, in contrast with the case of resolvents, we have to use Proposition~\ref{prop: general OD extension} for exponents $p \neq 2$.

\begin{lem}
\label{lem: cancellation functional calculus}
Let $p \in [2,\infty)$ be such that $((1+t^2 L)^{-1})_{t>0}$ is $\L^p$-bounded. Suppose that $\psi$ is of class $\Psi_\sigma^\tau$ on any sector, where $\tau > 0$ and $\sigma > \nicefrac{n}{(2p)}$. Then 
\begin{align*}
	\psi(t^2 L) c = 0 \quad (c \in \IC^m, \, t >0).
\end{align*}
\end{lem}

\begin{proof}
Let $\theta \in (0,1]$ be such that $q \coloneqq [p,2]_\theta$ satisfies $\sigma > \nicefrac{n}{(2q)}$. According to Lemma~\ref{lem: functional calculus bounds from J(L) abstract} the family $(\psi(t^2 L))_{t>0}$ satisfies $\L^q$ off-diagonal estimates of order $2 \sigma > \nicefrac{n}{q}$. Hence, $\psi(t^2 L) c$ is defined via Proposition~\ref{prop: general OD extension}.

Lemma~\ref{lem: OD extrapolation to sectors} provides $\L^{q}$ off-diagonal decay for the resolvents of $L$ of arbitrarily large order on some sector $\cl{\S_\nu^+}$. We write the definition of $\psi(t^2 L)$ as
\begin{align*}
	\psi(t^2 L) = \frac{1}{2 \pi \i} \int_{\bd \S_{\pi - 2\nu}^+} \psi(t^2 z) (1-z^{-1}L)^{-1} \, \frac{\d z}{z}.
\end{align*}
Setting $B \coloneqq B(0,1) \subseteq \R^n$, we formally have
\begin{align*}
	\sum_{j \geq 1} \psi(t^2 L) (\ind_{C_j(B)} c)  
	&= \frac{1}{2 \pi \i} \int_{\bd \S_{\pi - 2\nu}^+} \psi(t^2 z) \sum_{j \geq 1} (1-z^{-1}L)^{-1}(\ind_{C_j(B)} c) \, \frac{\d z}{z} \\
	&=\frac{1}{2 \pi \i} \int_{\bd \S_{\pi - 2\nu}^+} \psi(t^2 z) c \, \frac{\d z}{z} \\
	& = 0,
\end{align*}
where the second line uses the conservation property and the third one Cauchy's theorem. It remains to justify convergence and interchanging sum and integral sign in the first line. 

To this end, fix any compact set $K \subseteq \R^n$. Using off-diagonal estimates, we obtain for all $j$ large enough to grant for $\d(K, C_j(B)) \geq 2^{j-1}$ that
\begin{align*}
	\|\psi(t^2 z) (1-z^{-1}L)^{-1} &(\ind_{C_j(B)} c) \|_{\L^q(K)} \\
	&\lesssim |\psi(t^2 z)|(1+2^{j-1} |z|^{\frac{1}{2}})^{-\gamma} \|c\|_{\L^q(C_j(B))}  \\
	&\lesssim 2^{j(\frac{n}{q} - \gamma)}
	{\begin{cases}
			t^{-2 \tau} |z|^{-\frac{\gamma}{2}-\tau}  & \text{if } |z| \geq 1 \\
			t^{2 \sigma} |z|^{\sigma-\frac{\gamma}{2}} & \text{if } |z| \leq 1
	\end{cases}},
\end{align*}
where $\gamma >0$ is at our disposal.  We take $\nicefrac{n}{q} < \gamma < 2 \sigma$, in which case the right-hand side takes the form $2^{-j\eps} F_t(z)$ with $\eps > 0$ and $F_t \in \L^1(\bd \S_{\pi - 2\nu}^+, \nicefrac{\d |z|}{|z|})$, locally uniformly in $t$. This justifies at once convergence and interchanging sum and integral sign in $\L^q(K)$.
\end{proof}

Our third conservation property concerns the Poisson semigroup. In line with the previous result we need $\L^p$-boundedness of the resolvents for large $p$ to compensate for the poor decay of $\e^{-\sqrt{z}}-1$ at the origin.

\begin{prop}[Conservation property for the Poisson semigroup\index{conservation property!for the Poisson semigroup}]
\label{prop: conservation Poisson}
If $((1+t^2 L)^{-1})_{t>0}$ is $\L^p$-bounded for some $p>n$, then
\begin{align*}
	\e^{-t {L^{1/2}}} c  = c \quad (c \in \IC^m, \, t >0).
\end{align*}
\end{prop}

\begin{proof}
We have $\e^{-\sqrt{z}} = (1+z)^{-1} + \psi(z)$ with $\psi \in \Psi_{1/2}^1$ on any sector and the claim follows from Corollary~\ref{cor: conservation} and Lemma~\ref{lem: cancellation functional calculus}.
\end{proof}
\section{The four critical numbers}
\label{sec: J and N}

\noindent In this chapter, we define the four numbers that rule the functional calculus properties of our elliptic operators and that will help us to describe the ranges of well-posedness of our boundary value problems. We introduce the sets\index{J@$\cJ(L)$}
\begin{align*}
\cJ(L) &\coloneqq \big\{p \in (1_*, \infty) : ((1+t^2L)^{-1})_{t>0} \text{ is } a^{-1}\H^{p}\text{-bounded}\big\}
\end{align*}
and\index{N@$\cN(L)$}
\begin{align*}
\cN(L) &\coloneqq \big\{p \in (1_*, \infty) : (t\nabla_x(1+t^2L)^{-1})_{t>0} \text{ is }  a^{-1}\H^{p} - \H^{p}\text{-bounded} \big\},
\end{align*}
where we recall that $1_* = \nicefrac{n}{(n+1)}$. These sets contain $p=2$ (Corollary~\ref{cor: OD for second order}) and since the notion of $a_1\H^{p} - a_2\H^{p}$-boundedness interpolates, they are in fact intervals. 

\begin{defn}
\label{def: p+-,q+-}
The lower and upper endpoints of $\cJ(L)$ are denoted by $p_{-}(L)$ and $p_+(L)$, respectively. Similarly, the endpoints of $\cN(L)$ are denoted by $q_{-}(L)$ and $q_+(L)$. 
\end{defn} 

The exponents $p_\pm(L)$ and $q_\pm(L)$ are called \emph{critical numbers}\index{critical numbers} in the following. In this section we study intrinsic relations between these numbers, using the machinery developed in Section~\ref{sec: Hp Hq boundedness}. For the various duality arguments in this section we recall that $L^\sharp = -(a^*)^{-1} \div_x d^* \nabla_x$ is in the same class as $L$ and similar to $L^*$ under conjugation with $a^*$. In particular, we have
\begin{align}
\label{eq1: endpoints for Hodge projector}
\begin{split}
1 \vee p_{-}(L^\sharp) &= p_+(L)', \\
(1 \vee p_{-}(L))' &= p_+(L^\sharp).
\end{split}
\end{align}
\subsection{General facts on critical numbers}
\label{subsec: General facts on crictical numbers}

Here, we prove the following general relations between the four critical numbers. In fact, there are only three of them since $p_-(L)$ and $q_-(L)$ coincide\index{critical numbers!inner relationship}. The two inequalities are best possible in the class of all operators $L$, see Remark~\ref{rem: Frehse} further below.

\begin{thm}
	\label{thm: standard relation J(L) and N(L)}
	The critical numbers satisfy
	\begin{align*}
		p_{-}(L) &= q_{-}(L), \\
		p_{+}(L) &\geq q_+(L)^*, \\
		p_{-}(L) &\leq (q_+(L^\sharp)')_*.
	\end{align*}
\end{thm}

We prepare the proof through a sequence of lemmata that are of independent interest.

\begin{lem}
\label{lem: Lp L2 between Sobolev conjugates n geq 2}
Let $n \geq 2$. Then $(2_*,2^*) \subseteq \cJ(L)$ and $((1+t^2L)^{-1})_{t >0}$ is $\L^2 - \L^q$-bounded and $\L^{q'} - \L^2$-bounded for every $q \in [2, 2^*] \cap [2,\infty)$.
\end{lem}

\begin{proof}
We have $2^* = \infty$ when $n = 2$ and $2^* < \infty$ when $n \geq 3$. The restriction on $q$ is precisely such that we have the Gagliardo--Nirenberg inequality
\begin{align*}
\|u\|_q \lesssim \|\nabla_x u\|_2^\alpha \|u\|_2^{1-\alpha} \quad (u \in \W^{1,2}(\R^n)),
\end{align*}
where $\alpha = \nicefrac{n}{2}-\nicefrac{n}{q}$. We set $u \coloneqq (1+t^2 L)^{-1}f$, $f \in \L^2$, $t >0$, and use the $\L^2$-boundedness of the resolvent and gradient families to give
\begin{align*}
\|(1+t^2 L)^{-1}f\|_q 
&\lesssim t^{-\alpha} \|f\|_2.
\end{align*}
Hence, the resolvents are $\L^2 - \L^q$-bounded. Interpolation with the $\L^2$ off-diagonal estimates by means of Lemma~\ref{lem: OD interpolation} leads to $\L^{2} - \L^{q}$ off-diagonal estimates of arbitrarily large order for any $q \in [2, 2^*)$ and $\L^q$-boundedness follows from Lemma~\ref{lem: OD implies boundedness}.

The  rest follows by duality and similarity by applying the above to $L^\sharp$ in place of $L$. 
\end{proof}

In dimension $n=1$ we have $2_* = \nicefrac{2}{3}$ and by analogy with the previous lemma we expect that $(\nicefrac{2}{3}, \infty) \subseteq \cJ(L)$. However, in the one-dimensional situation we have $\div_x = \nabla_x$ and this allows us to improve the lower bound to the best possible value $1_* = \nicefrac{1}{2}$.

\begin{lem}
\label{lem: Lp L2 between Sobolev conjugates n = 1}
Let $n = 1$. Then $(\nicefrac{1}{2},\infty) \subseteq \cJ(L)$ and $[2,\infty) \subseteq \cN(L)$. Moreover $((1+t^2L)^{-1})_{t >0}$ is $a^{-1} \H^p - \L^2$-bounded for every $p \in (\nicefrac{1}{2},2]$ and $((1+t^2L)^{-1})_{t >0}$ and $(t\tfrac{\d}{\d x} (1+t^2L)^{-1})_{t >0}$ are both $\L^2 - \L^q$-bounded for every $q \in [2,\infty)$.
\end{lem}

\begin{proof}
In the one-dimensional setting the operator $L$ takes the form $L= - a^{-1} \frac{\d}{\d x} (d \frac{\d}{\d x})$ and the space $\cH$ in \eqref{eq: accretivity A} coincides with $\L^2$. In particular, just as $a$, also $d$ is strictly elliptic.

\medskip

\noindent \emph{Step~1: $\L^2 - \L^q$-bound for the gradients.} It suffices to obtain the bound for $t=1$ with an implicit constant that depends on the coefficients only through ellipticity. Indeed, for $t \neq 1$ we can use the change of variable $f_t(x) \coloneqq f(tx)$ in order to write
\begin{align*}
t\tfrac{\d}{\d x}(1+t^2 L)^{-1} f(x) = (\tfrac{\d}{\d x}(1+L_t)^{-1} f_t)(t^{-1} x),
\end{align*}
where $L_t \coloneqq - a_t^{-1} \frac{\d}{\d x} (d_t \frac{\d}{\d x})$ has the same ellipticity constant as $L$.

Let now $f \in \L^2$ and set $u \coloneqq (1+L)^{-1} f$, so that $\frac{\d}{\d x} (d \frac{\d}{\d x} u) = af - au$. In one dimension the Sobolev embedding $\W^{1,2} \subseteq \L^q$ holds for any $q \in [2,\infty)$. Thus, we have
\begin{align*}
\|\tfrac{\d}{\d x} u\|_q
&\simeq \|d \tfrac{\d}{\d x} u\|_q \\
&\lesssim \|d \tfrac{\d}{\d x} u\|_{\W^{1,2}} \\
&\lesssim \|d \tfrac{\d}{\d x} u\|_2 + \|af\|_{2} + \|au\|_2 \\
&\lesssim \|f\|_2,
\end{align*}
where in the final step we have used the $\L^2$-boundedness of the resolvent and gradient families. This is the required $\L^2-\L^q$-bound.

\medskip

\noindent \emph{Step~2: $\L^q$-bound for the gradients}. This follows from Lemma~\ref{lem: OD interpolation} and Lemma~\ref{lem: OD implies boundedness} as in the previous proof. Hence, we have $[2,\infty) \subseteq \cN(L)$.

\medskip

\noindent \emph{Step~3: Bounds for the resolvents}. Let $q \in [2,\infty)$ and define $\varrho \in (\nicefrac{1}{2}, \nicefrac{2}{3}]$ through $1-\nicefrac{1}{q} = \nicefrac{1}{\varrho}-1$. 

For $f \in \L^2$ and $t>0$ we use the Sobolev embedding $\Wdot^{1,q} \subseteq \Lamdot^{1-\nicefrac{1}{q}}$ and the result of Step~1 for $L^\sharp$ to give
\begin{align*}
\|(1+t^2 L^\sharp)^{-1}f\|_{\Lamdot^{1-1/q}}
&\lesssim \|\tfrac{\d}{\d x} (1+t^2 L^\sharp)^{-1}f\|_q \\
&\lesssim t^{-1+\frac{1}{q}-\frac{1}{2}} \|f\|_2.
\end{align*}
Hence, the resolvents of $L^\sharp$ are $\L^2 - \Lamdot^{1-1/q}$-bounded. Since we have $L^\sharp = (a^*)^{-1}L^*a^*$, we obtain by duality that the resolvents of $L$ are $a^{-1} \H^{\varrho} - \L^2$-bounded, see Lemma~\ref{lem: p-q boundedness duality}. They also satisfy $\L^2$ off-diagonal estimates of arbitrarily large order and have the cancellation property $\int_{\R^n} a(1+t^2L)^{-1} a^{-1}f \d x = 0$ if $f \in \L^2$ has compact support and integral zero, see Corollary~\ref{cor: conservation}. Hence, we are in a position to apply Lemma~\ref{lem: OD Hardy implies boundedness} and obtain $a^{-1} \H^p$-boundedness for $p \in (\varrho,1]$. 

Since $q \in [2, \infty)$ was arbitrary, the conclusion is that the resolvents are $a^{-1} \H^{\varrho} - \L^2$-bounded and $a^{-1}\H^p$-bounded for all for all $\varrho \in (\nicefrac{1}{2}, \nicefrac{2}{3}]$ and all $p \in (\nicefrac{1}{2},1]$. By interpolation with the $\L^2$-bound we can allow all $\varrho, p \in (\nicefrac{1}{2},2]$. Finally, the $\L^2 - \L^q$ and $\L^q$-bounds of the resolvents for all $q \in [2,\infty)$ follow again from duality and similarity, by applying the results for $p \in (1,2]$ to $L^\sharp$.
\end{proof}

We also need a result that allows us to switch between powers of the resolvent in $\H^p - \H^q$-estimates.

\begin{lem}
\label{lem: reduce resolvent powers}
Let $1_*<p \leq q < \infty$ with $q>1$ and $\nicefrac{n}{p}-\nicefrac{n}{q} < 1$. Suppose that there exists an integer $\beta \geq 1$ such that $(t \nabla_x (1+t^2 L)^{-\beta-1})_{t>0}$ is $a^{-1} \H^p - \L^q$-bounded. Then also $(t \nabla_x (1+t^2 L)^{-1})_{t>0}$ is $a^{-1} \H^p - \L^q$-bounded.
\end{lem}

\begin{proof}
Let $t>0$ and $f \in \L^2$. The Calder\'on reproducing formula for the injective sectorial operator $T=1+t^2L$ and the auxiliary function $\varphi(z) = z(1+z)^{-\beta-1}$ reads
\begin{align}
\label{eq: resolution of resolvents through high powers}
f =  \beta\int_0^\infty (1+t^2 L)(1+u+t^2u L)^{-\beta-1}f \, \d u.
\end{align}
Applying the bounded operator $t\nabla_x(1+t^2L)^{-1}$ and re-arranging terms gives
\begin{align*}
&t \nabla_x (1+t^2L)^{-1}f \\
&= \beta \int_0^\infty \frac{1}{(1+u)^{\beta+\frac{1}{2}} u^{\frac{1}{2}}} \bigg(\frac{u t^2}{1+u}\bigg)^{\frac{1}{2}} \nabla_x\bigg(1+ \frac{u t^2}{1+u} L \bigg)^{-\beta-1}f \, \d u.
\end{align*}
Now, we let $f \in \H^p \cap \L^2$, apply the formula to $a^{-1}f$, and take $\L^q$ norms on both sides, in order to give
\begin{align*}
\|t &\nabla_x (1+t^2L)^{-1} a^{-1}f \|_q \\
&\lesssim \|f\|_{\H^p} \int_0^\infty \frac{1}{(1+u)^{\beta+\frac{1}{2}} u^{\frac{1}{2}}} \bigg(\frac{u t^2}{1+u}\bigg)^{\frac{1}{2}(\frac{n}{q}-\frac{n}{p})} \,\d u \\
&\leq t^{\frac{n}{q}-\frac{n}{p}} \|f\|_{\H^p} \int_0^\infty \frac{u^{\frac{n}{2q}-\frac{n}{2p} - \frac{1}{2}}}{(1+u)^{\beta+\frac{1}{2} + \frac{n}{2q}-\frac{n}{2p}}} \, \d u.
\end{align*}
The numerical integral in $u$ converges as we have $\nicefrac{n}{q}-\nicefrac{n}{p}  > -1$ by assumption.
\end{proof}

\begin{proof}[Proof of Theorem~\ref{thm: standard relation J(L) and N(L)}]
The argument is in two steps.

\medskip

\noindent \emph{Step 1: Resolvent estimates from gradient bounds.} Here, we show the upper bound in the first line and the second and third lines. In dimension $n = 1$ we have $p_-(L) = 1_*$ and $p_+(L) = \infty$ by Lemma~\ref{lem: Lp L2 between Sobolev conjugates n = 1} and there is nothing to prove. For the rest of the step we assume $n \geq 2$.

Let $\varrho \in \cN(L)$. If $\varrho < n$, then a Sobolev embedding yields for all $f \in \H^\varrho \cap \L^2$ and all $t>0$ that
\begin{align*}
\|(1+t^2 L)^{-1} a^{-1}f\|_{\L^{\varrho^*}} 
&\lesssim t^{-1} \|t\nabla_x (1+t^2 L)^{-1} a^{-1} f\|_{\H^\varrho} \\
&\lesssim t^{-1} \|f\|_{\H^\varrho}.
\end{align*}
Hence, the resolvents of $L$ are $a^{-1}\H^\varrho - \L^{\varrho^*}$-bounded. Likewise, if $\varrho > n$, then we obtain for all $f \in \L^\varrho \cap \L^2$ and all $t>0$ that
\begin{align*}
\|(1+t^2 L)^{-1} f\|_{\Lamdot^{1- n/\varrho}} 
\lesssim t^{-1} \|f\|_{\L^\varrho}.
\end{align*}
By duality the resolvents of $L^\sharp$ are $(a^*)^{-1} \H^{r} - \L^{\varrho'}$-bounded. The exponent $r$ is determined by $1-\nicefrac{n}{\varrho} = n(\nicefrac{1}{r}-1)$, that is, $r= (\varrho')_*$. From these observations, we can infer further mapping properties in each case. 

\medskip

\noindent \emph{Step 1a: The Lebesgue case $1 < \varrho < n$.}  The resolvents of $L$ are $\L^\varrho - \L^{\varrho^*}$-bounded. Lemma~\ref{lem: OD interpolation} yields $\L^{[\varrho,2]_\theta} - \L^{[\varrho^*,2]_\theta}$ off-diagonal estimates of arbitrarily large order, where $\theta \in (0,1)$ is arbitrary. Lemma~\ref{lem: OD implies boundedness} yields both $\L^{[\varrho,2]_\theta}$ and $\L^{[\varrho^*,2]_\theta}$-boundedness. Consequently, we must have $p_-(L) \leq \varrho$ and $p_+(L) \geq \varrho^*$. 
\medskip

\noindent \emph{Step 1b: The Hardy case $\varrho \leq 1$.} Since $\cN(L)$ is an interval, we have $(\varrho,2) \subseteq \cN(L)$. The first part applies to all exponents in $(1,2)$ instead of $\varrho$ and we first get $\L^q$-boundedness of the resolvents of $L$ for all $q \in (1,2^*)$ and then $\L^q$ off-diagonal estimates of arbitrarily large order by interpolation.

If $\varrho = 1$, then $p_-(L) \leq \varrho$ follows directly.

Now, assume $\varrho < 1$. As $\varrho>1_*$, we can take $q \coloneqq \varrho^*$ and have $\L^q$ off-diagonal estimates of arbitrarily large order and $a^{-1}\H^\varrho - \L^{q}$-boundedness. For compactly supported $f \in \L^2$, Corollary~\ref{cor: conservation} yields $\int_{\R^n} a(1+t^2L)^{-1} (a^{-1}f) \d x = 0$. We have verified the assumptions of Lemma~\ref{lem: OD Hardy implies boundedness} and obtain that the resolvents of $L$ are $a^{-1}\H^p$-bounded for every $p \in (\varrho,1]$. Therefore, we have again $p_-(L) \leq \varrho$.

\medskip

\noindent \emph{Step 1c: The H\"older case $\varrho > n$.}  From the preliminary discussion we know that the resolvents of $L^\sharp$ are $(a^*)^{-1} \H^{(\varrho')_*} - \L^{\varrho'}$-bounded. We claim that they satisfy $\L^{\varrho'}$ off-diagonal estimates of arbitrarily large order. Taking the claim for granted, $p_-(L^\sharp) \leq (\varrho')_*$ follows as in the previous step.

For the claim we first prove $(1,2) \subseteq \cJ(L^\sharp)$. In dimension $n=2$ this is due to Lemma~\ref{lem: Lp L2 between Sobolev conjugates n geq 2}. In dimension $n \geq 3$ we have $(2,n) \subseteq \cN(L)$ since the latter is an interval that contains $2$ and $\varrho$. Step~1a applies to all exponents in $(2,n)$ in place of $\varrho$ and yields $(2,\infty) \subseteq \cJ(L)$. By duality, we get again $(1,2) \subseteq \cJ(L^\sharp)$. As we have $\varrho' \in (1,2)$, the $\L^{\varrho'}$ off-diagonal estimates for the resolvents of $L^\sharp$ follow by interpolation with the $\L^2$-result.

\medskip

Let us conclude Step 1. In dimension $n \geq 3$ the set $\cN(L) \cap (1_*,n)$ is non-empty because it contains $2$. Letting $\varrho$ vary over $\cN(L) \cap (1_*,n)$, we conclude $p_-(L) \leq q_-(L)$ and $p_+(L) \geq q_+(L)^*$ from Steps 1a \& 1b. In dimension $n=2$ the same argument applies unless $q_-(L)=2$. But in this case\footnote{In fact this case never occurs as we shall see later on.} the inequalities in question trivially hold because we have $p_-(L) \leq 1$ and $p_+(L) = \infty$ by Lemma~\ref{lem: Lp L2 between Sobolev conjugates n geq 2}. 

As for the third line in the theorem, if $q_+(L^\sharp) \leq n$, then 
\begin{align*}
	p_-(L)
	\leq (p_{+}(L^\sharp))' 
	\leq (q_+(L^\sharp)^*)'
	=(q_+(L^\sharp)')_*
\end{align*}
follows from \eqref{eq1: endpoints for Hodge projector} and the second line. If $q_+(L^\sharp) > n$, then the inequality $p_-(L)\leq (q_+(L^\sharp)')_*$ follows from Step~1c with the roles of $L$ and $L^\sharp$ switched.

\medskip

\noindent \emph{Step 2: Gradient bounds from resolvent estimates.} Let $p \in \cJ(L)$ with $p<2$. Hence, $((1+t^2L)^{-1})_{t>0}$ is $a^{-1}\H^p$-bounded. Lemmata~\ref{lem: Lp L2 between Sobolev conjugates n geq 2} and~\ref{lem: Lp L2 between Sobolev conjugates n = 1} guarantee that this family is $\L^\varrho - \L^2$-bounded for some $\varrho \in (1,2)$. According to Lemma~\ref{lem: extra}, we find for every $q \in (p,2)$ an integer $\beta(q) \geq 1$ such that $((1+t^2 L)^{-\beta(q)})_{t>0}$ is $a^{-1}\H^q - \L^2$-bounded. By composition with the $\L^2$-bounded gradient family, $(t \nabla_x (1+t^2L)^{-\beta(q)-1})_{t>0}$ is $a^{-1} \H^q - \L^2$-bounded. Note that by interpolation the integer $\beta(q)$ can be taken the same on compact subsets of $(p,2)$.

\medskip

\noindent \emph{Step 2a: The Lebesgue case $p \geq 1$.} 
We know that $(t \nabla_x (1+t^2L)^{-\beta(q)-1})_{t>0}$ is $\L^q - \L^2$-bounded. By composition, this family also satisfies $\L^2$ off-diagonal estimates of arbitrarily large order. Since this holds for every $q \in (p,2)$, we can run the usual argument: $\L^q - \L^2$ off-diagonal estimates of arbitrarily larger order follow by interpolation and this implies $\L^q$-boundedness. Thanks to Lemma~\ref{lem: reduce resolvent powers} we get $\L^q$-boundedness also for  $(t \nabla_x (1+t^2L)^{-1})_{t>0}$. Since $q \in (p,2)$ was arbitrary, we have $q_-(L) \leq p$.

\medskip

\noindent \emph{Step 2b: The Hardy case $p < 1$.} We slightly refine the argument in the Lebesgue case by appealing to Lemma~\ref{lem: OD Hardy implies boundedness}. In the following let $q \in (p,1)$ and $s \in (1,2)$ such that $\nicefrac{1}{q} - \nicefrac{1}{s}<\nicefrac{1}{n}$. Such $s$ exists since we have $p> 1_*$. 

First, consider the family $(t\nabla_x (1+t^2 L)^{-\beta(q)-1})_{t>0}$. It is $a^{-1}\H^q - \L^2$-bounded and satisfies $\L^2$ off-diagonal estimates of arbitrarily large order. For compactly supported $f \in \L^2$ we get that $(1+t^2 L)^{-\beta(q)}(a^{-1}f)$ and $\nabla_x (1+t^2 L)^{-\beta(q)}(a^{-1}f) $ are in $\L^1$ from the $\L^2$ off-diagonal decay of order $\gamma > \nicefrac{n}{2}$. The integral of the gradient of a $\W^{1,1}$-function vanishes, so $\int_{\R^n} t \nabla_x (1+t^2 L)^{-\beta(q)}(a^{-1}f) \d x = 0$. We have checked the assumptions of Lemma~\ref{lem: OD Hardy implies boundedness} and obtain $a^{-1}\H^q - \H^q$-boundedness for every $q$. This interpolates with the original $a^{-1}\H^q - \L^2$-boundedness, so that the conclusion is $a^{-1}\H^q - \L^s$-boundedness for all $q$ and $s$.

Now, we consider $(t\nabla_x (1+t^2 L)^{-1})_{t>0}$. Lemma~\ref{lem: reduce resolvent powers} yields $a^{-1}\H^q - \L^s$-boundedness for all $q$ and $s$. Step~2a applies to $s$ and yields $\L^{s}$-boundedness, which implies $\L^s$ off-diagonal decay of arbitrarily large order for every $s$ by interpolation with the $\L^2$-result. As before, we also have $\int_{\R^n} t \nabla_x (1+t^2 L)^{-1} (a^{-1}f) \d x = 0$ for compactly supported $f \in \L^2$. We have again verified the assumptions of Lemma~\ref{lem: OD Hardy implies boundedness} and conclude for $a^{-1} \H^q - \H^q$-boundedness for every $q$. Thus, we have  $q_-(L) \leq p$.

\medskip

As $p \in \cJ(L) \cap (1_*,2)$ was arbitrary, Steps~2a \& 2b yield the missing inequality $q_-(L) \leq p_-(L)$ that completes the proof of the first line in the theorem.
\end{proof}
\subsection{Worst-case estimates for the critical numbers}
\label{subsec: Worst-case estimates for the crictical numbers}

The following extrapolation from the $\L^2$-theory has been proved by an application of \v{S}ne\u{i}berg's stability theorem~\cite{Sneiberg, ABES3}.\index{Theorem!\v{S}ne\u{i}berg's} 

\begin{prop}[{\cite[Prop.~4.5]{AusStaRemarks}}]
\label{prop: Lp boundedness BD}
There exists $\eps > 0$, depending on ellipticity and dimensions, such that whenever $p \in [2-\eps,2+\eps]$, then $((1+\i tDB)^{-1})_{t \in \R}$ is $\L^p$-bounded.
\end{prop}

We use this result to give the following global picture for the critical numbers for the class of all $L$ in all dimensions.\index{critical numbers!general bounds}

\begin{prop}
\label{prop: J(L) contains neighborhood of Sobolev conjugates}
The following relations hold.
\begin{enumerate}
	\item In dimension $n = 1$, 
	\begin{align*}
	p_-(L) = q_-(L) = \frac{1}{2} \quad \& \quad  p_+(L) = q_+(L) = \infty.
	\end{align*}
	\item In dimension $n \geq 2$ there exists $\eps> 0$, depending on ellipticity and dimensions, such that
	\begin{align*}
	\qquad p_-(L) = q_-(L) \leq 2_* - \eps \quad \& \quad p_+(L) \geq 2^* + \eps \quad \& \quad q_+(L) \geq 2 + \eps.
   	\end{align*}
\end{enumerate}
\end{prop}

\begin{proof}
The identification $q_-(L) = p_-(L)$ in any dimension is due to Theorem~\ref{thm: standard relation J(L) and N(L)}. In dimension $n=1$, Lemma~\ref{lem: Lp L2 between Sobolev conjugates n = 1} shows that $p_-(L)=\nicefrac{1}{2}$ and $q_+(L)=\infty=p_+(L)$ take the best possible values. Hence, (i) follows.

As for (ii), we use \eqref{eq: tL and tM} to write, whenever $t>0$,
\begin{align*}
\frac{1}{2} \bigg((1+\i t DB)^{-1} - (1-\i t DB)^{-1} \bigg)
&=-\i t DB (1+t^2(DB)^2)^{-1} \\
&= \begin{bmatrix}
- \i t \div_x d (1+t^2 \tM)^{-1} \\ \i t \nabla_x a^{-1} (1+t^2 \tL)^{-1}
\end{bmatrix}.
\end{align*}
This family is $\L^p$-bounded for $p \in [2-\eps,2+\eps]$ due to Proposition~\ref{prop: Lp boundedness BD}. In particular, the second component is $\L^p$-bounded and since $a$ is strictly elliptic, the same is true for $t \nabla_x  a^{-1} (1+t^2 \tL)^{-1} a = t \nabla_x (1+t^2 L)^{-1}$. Hence, for a possibly different choice of $\eps$ we have $q_+(L) \geq 2+\eps$ and $q_-(L) \leq 2-\eps$. The same thing for $L^\sharp$.  Now, the claim follows from Theorem~\ref{thm: standard relation J(L) and N(L)}.
\end{proof}

\begin{rem}
\label{rem: Frehse}
\begin{enumerate}
	\item 	In the one-dimensional setting the identification of the critical numbers could also be obtained from the kernel estimates in~\cite{AMcT}. They are only stated for $m=1$ but the argument literally applies to systems ($m > 1$) under our ellipticity assumption. In fact, the proof of Lemma~\ref{lem: Lp L2 between Sobolev conjugates n = 1} mimics some intermediate steps in \cite{AMcT}. The value $p_-(L) = \nicefrac{1}{2}$ has appeared in a related context in \cite{AT}.
	\item In higher dimensions the bounds above cannot be improved in general, even when $a=1$ and $m=1$. More precisely, given $\eps>0$, any of $p_-(L) < 2_* -\eps$, $p_+(L) > 2^* + \eps$, $q_+(L) > 2+\eps$ can fail for some $L$.

	Indeed, for $p_\pm$ in dimensions $n \geq 3$, counterexamples rely on Frehse's irregular solution~\cite{Frehse} and can be found in \cite[Prop.~2.10]{HMMc}. In view of Theorem~\ref{thm: standard relation J(L) and N(L)} such counterexamples satisfy  $2^*+\eps \geq p_+(L) \geq q_+(L)^* \geq 2^*$. Hence, they also serve as counterexamples to the general improvement of $q_+$ and show that the inequalities in Theorem~\ref{thm: standard relation J(L) and N(L)} are best possible in the class of all operators $L$.
	
	When $n=2$, the counterexample for $q_+$ due to Kenig comes with $d$ real symmetric~\cite[Sec.~4.2.2]{AT-Asterisque}. The same operator is a counterexample for the general improvement on $p_{-}$, that is, $p_{-}(L)$ can be as close to $2_{*}=1$ as one wants. Hence, the final inequality in Theorem \ref{thm: standard relation J(L) and N(L)} is again best possible. There is no discussion of $p_+$ since $2^* = \infty$.
\end{enumerate}
\end{rem}
\subsection{\texorpdfstring{$\boldsymbol{a}$}{a}-independence of critical numbers}
\label{subsec: a-independence the crictical numbers}

It is tempting to compare the critical numbers for $L$ with those for
\begin{align*}
	L_0 = -\div_x d \nabla_x,
\end{align*}
seeing $L= a^{-1} L_0$ as a multiplicative perturbation of $L_0$.\index{critical numbers!a@$a$-independence} Let us prove that the critical numbers for both operators are indeed the same.

\begin{thm}
\label{thm: critical numbers a-independent}
The critical numbers for $L$ and $L_0$ coincide, that is,
\begin{align*}
	p_\pm(L) = p_\pm(L_0) \quad \& \quad q_\pm(L) = q_\pm(L_0).
\end{align*}
\end{thm}

\begin{proof}
The claim in dimension $n=1$ is an immediate consequence of Proposition~\ref{prop: J(L) contains neighborhood of Sobolev conjugates}.  The proof in dimensions $n \geq 2$ is divided into six steps.

\medskip

\noindent \emph{Step~1: $p_-(L) \leq p_-(L_0)$.} Let $p \in \cJ(L_0) \cap (1_*,2_*]$. This interval is non-empty thanks to Proposition~\ref{prop: J(L) contains neighborhood of Sobolev conjugates}. We set $p_0 \coloneqq p$, define iteratively $p_k \coloneqq p_{k-1}^*$ and stop at the first exponent $k^+ \geq 0$ with $p_{k^+} \in (2_*,2]$. We shall prove by backward induction that $(p_k,2] \subseteq \cJ(L)$ for all $k$. Hence, we eventually find $(p,2] \subseteq \cJ(L)$ and taking the infimum over all $p$ yields $p_-(L) \leq p_-(L_0)$. 

Once again by Proposition~\ref{prop: J(L) contains neighborhood of Sobolev conjugates}, we have $(p_{k^+},2] \subseteq \cJ(L)$. For the inductive step we assume $(p_k,2] \subseteq \cJ(L)$ and pick any $q \in (p_{k-1}, 2_*]$. For all $t>0$ we have
\begin{align*}
	1 = (a + t^2 L_0)(1+t^2 L_0)^{-1} + (1-a) (1+t^2 L_0)^{-1}
\end{align*}
as operators on $\L^2$. Multiplication by $(1+t^2 L)^{-1} a^{-1} = (a+t^2 L_0)^{-1}$ from the left yields the key identity
\begin{align}
\label{eq1: a-independence critical numbers}
\begin{split}
 &(1+t^2 L)^{-1} a^{-1} \\
 &\quad= (1+t^2 L_0)^{-1} + (1+t^2 L)^{-1} a^{-1} (1-a) (1+t^2 L_0)^{-1}.
 \end{split}
\end{align}
On the right-hand side $((1+t^2 L)^{-1})_{t>0}$ is $\L^{q^*}$-bounded by the induction hypothesis. By Theorem~\ref{thm: standard relation J(L) and N(L)} we have $q_-(L_0) = p_-(L_0)$ so that $((t \nabla_x (1+t^2 L_0)^{-1})_{t>0}$ is $\H^q$-bounded. By a Sobolev embedding we have
\begin{align*}
	\|(1+t^2 L_0)^{-1} f\|_{q^*} \lesssim 	\|\nabla_x(1+t^2 L_0)^{-1} f\|_{\H^q} \lesssim t^{-1} \|f\|_{\H^q},
\end{align*}
whenever $f \in \H^q \cap \L^2$ and $t>0$. Hence, $((1+t^2 L_0)^{-1})_{t>0}$ is $\H^q-\L^{q^*}$-bounded. Now, it follows from \eqref{eq1: a-independence critical numbers} that $((1+t^2 L)^{-1})_{t>0}$ is $a^{-1} \H^q - \L^{q^*}$-bounded. This was the key step.

If $q>1$, then we have $\L^q-\L^{q^*}$-boundedness for the resolvents of $L$. Interpolation with the $\L^2$ off-diagonal estimates (Lemma~\ref{lem: OD interpolation}) followed by Lemma~\ref{lem: OD implies boundedness} yields $(q,2] \subseteq \cJ(L)$.

If $q=1$, then $(p_{k-1},2_*]$ also contains exponents that are strictly smaller than $1$ and we can jump right into the following case.

In the remaining case $q < 1$ we have $a^{-1} \H^q - \L^{q^*}$-boundedness for the resolvents of $L$. As $q^*$ is an interior point of $\cJ(L)$ by the induction hypothesis, we get again $\L^{q^*}$ off-diagonal estimates of arbitrarily large order from the ones on $\L^2$ by interpolation. For compactly supported $f \in \L^2$, Corollary~\ref{cor: conservation} yields $\int_{\R^n} a(1+t^2L)^{-1} a^{-1}f \, \d x = 0$. This means that we have verified the assumptions of Lemma~\ref{lem: OD Hardy implies boundedness} and $(q,2] \subseteq \cJ(L)$ follows.

\medskip

\noindent \emph{Step~2: $p_-(L_0) \leq p_-(L)$.} We only need a key identity replacing \eqref{eq1: a-independence critical numbers} and allowing us to deduce $\H^q-\L^{q^*}$-boundedness of $((1+t^2 L_0)^{-1})_{t>0}$ from $\L^{q^*}$-boundedness of $((1+t^2 L_0)^{-1})_{t>0}$ and $a^{-1} \H^q - \L^{q^*}$-boundedness of $((1+t^2 L)^{-1})_{t>0}$. The rest of the proof for $p_-(L) \leq p_-(L_0)$ was symmetric in $L$ and $L_0$.

For the new key identity we split
\begin{align*}
	1 = (1+t^2 L_0) (a+t^2 L_0)^{-1} + (a-1) (a+t^2 L_0)^{-1}
\end{align*}
and multiply by $(1+ t^2 L_0)^{-1}$ from the left in order to get the desirable decomposition
\begin{align*}
	(1+t^2 L_0)^{-1}  = (1+t^2 L)^{-1}  a^{-1} + (1+t^2 L_0)^{-1} (a-1) (1+t^2 L)^{-1}  a^{-1}.
\end{align*}

\noindent \emph{Step~3: $q_-(L) = q_-(L_0)$.} This follows from the first two steps and Theorem \ref{thm: standard relation J(L) and N(L)}.

\noindent \emph{Step~4: $p_+(L) = p_+(L_0)$.} Simply note that by the duality relations \eqref{eq1: endpoints for Hodge projector} and the first two steps we have
\begin{align*}
	p_+(L) = (1 \vee p_-(L^\sharp))' = (1 \vee p_-(L_{0}^\sharp))' = p_+(L_0).
\end{align*}

\noindent \emph{Step~5: $q_+(L_0) \leq q_+(L)$.} Let $2 \leq q < q_+(L_0)$. For $t>0$ we use a new decomposition, namely
\begin{align}
\label{eq2: a-independence critical numbers}
\begin{split}
	&t \nabla_x (1+t^2 L)^{-1} \\
	&\quad= t \nabla_x (1+t^2 L_0)^{-1} ( a + t^2 L_0 + 1 - a) (1+t^2 L)^{-1} \\
	&\quad= t \nabla_x (1+t^2 L_0)^{-1} a + t \nabla_x (1+t^2 L_0)^{-1} (1-a) (1+t^2 L)^{-1}.
\end{split}
\end{align}
On the right-hand side $(t \nabla_x (1+t^2 L_0)^{-1})_{t>0}$ is $\L^q$-bounded by assumption and $((1+t^2 L)^{-1})_{t>0}$ is $\L^q$-bounded since we have $q_+(L_0) \leq p_+(L_0) = p_+(L)$ by Theorem~\ref{thm: standard relation J(L) and N(L)} and Step 4. Thus, $(t \nabla_x (1+t^2 L)^{-1})_{t>0}$ is $\L^q$-bounded. Taking the supremum over all $q$, we obtain $q_+(L_0) \leq q_+(L)$.

\medskip

\noindent \emph{Step~6: $q_+(L) \leq q_+(L_0)$.} The argument follows by reversing the roles of $L$ and $L_0$ in Step~4 and using the identity 
\begin{align*}
	&t \nabla_x (1+t^2 L_0)^{-1} \\
	&\quad= t \nabla_x (a+t^2 L_0)^{-1} ( 1 + t^2 L_0 + a -1) (1+t^2 L_0)^{-1} \\
	&\quad= t \nabla_x (1+t^2 L)^{-1} a^{-1} + t \nabla_x (1+t^2 L)^{-1} a^{-1} (a-1) (1+t^2 L_0)^{-1}
\end{align*}
instead of \eqref{eq2: a-independence critical numbers}.
\end{proof}

As an application of Theorem~\ref{thm: critical numbers a-independent} we determine the critical numbers of multiplicative perturbations of the (coordinatewise acting) Laplacian.\index{critical numbers!for multiplicative perturbations of the Laplacian}

\begin{cor}
\label{cor: McIntosh-Nahmod}
In any dimension it follows that 
\begin{align*}
	p_{-}(-a^{-1}\Delta_x) &= q_{-}(-a^{-1}\Delta_x) = 1_*, \\
	p_{+}(-a^{-1} \Delta_x) &= q_+(-a^{-1} \Delta_x) = \infty.
\end{align*}
\end{cor}

This result is originally due to McIntosh--Nahmod, see Theorem~3.3 and \SS5.(v) in \cite{McN}. Here, we use a rather different and simpler method.\index{Theorem!McIntosh--Nahmod's} In Section~\ref{sec: Critical numbers and kernels}, we shall discuss kernel estimates for these operators.

\begin{proof}[Proof of Corollary~\ref{cor: McIntosh-Nahmod}]
In view of Theorem~\ref{thm: standard relation J(L) and N(L)} we only have to prove that $q_+(-a^{-1} \Delta_x) = \infty$. By Theorem~\ref{thm: critical numbers a-independent} we have $q_+(-a^{-1} \Delta_x) = q_+(-\Delta_x)$ and there are many ways to see that $q_+(-\Delta_x) = \infty$. One is to note that $t\nabla_x(1+t^2 (-\Delta_x))^{-1}$ corresponds to the Fourier multiplier $\xi \mapsto \i t(1+t^2|\xi|^2)^{-1} \xi$, which falls under the scope of the Mihlin multiplier theorem.
\end{proof}
\section{Riesz transform estimates: Part I}
\label{sec: Riesz}

\noindent  In this chapter, we characterize the range of exponents $p$ for $\L^p$-bounded-ness of the Riesz transform $\nabla_x L^{-1/2}$. More generally, we introduce the set\index{I@$\cI(L)$}
\begin{align}
\label{eq: I(L)}
 \cI(L) \coloneqq \big\{ p \in (1_*,\infty) : \nabla_x L^{-1/2} \text{ is $a^{-1} \H^p - \H^p$-bounded}\big\}.
\end{align}
Some clarification on the meaning of $\nabla_x L^{-1/2}$ being $a^{-1} \H^p - \H^p$-bounded is necessary since there are two possible interpretations:
\begin{itemize}
	\item As we have seen in Section~\ref{subsec: Kato and Riesz}, ${L^{1/2}}: \W^{1,2} \to \L^2$ extends to an isomorphism $\Wdot^{1,2} \to \L^2$ that we denote again by ${L^{1/2}}$. In this sense $R_L \coloneqq \nabla_x L^{-1/2}$ is defined as a bounded operator on $\L^2$. The question of $a^{-1} \H^p - \H^p$-boundedness  for $R_L$ fits into the abstract framework of Section~\ref{sec: Hp Hq boundedness} and means that
	\begin{align*}
	\qquad \|R_L a^{-1}f\|_{\H^p} \lesssim \|f\|_{\H^p} \quad (f \in \H^p \cap \L^2)
	\end{align*}
	and when $p>1$ equivalently that 
	\begin{align*}
	\qquad \|R_L f\|_{p} \lesssim \|f\|_{p} \quad (f \in \L^p \cap \L^2).
	\end{align*}
	\item We could also avoid the extension, work directly with $\nabla_x L^{-1/2}$ defined on $\ran({L^{1/2}})$ and ask for $\|\nabla_x L^{-1/2} a^{-1}f\|_{\H^p} \lesssim \|f\|_{\H^p}$ for all $f \in \H^p \cap \ran(a{L^{1/2}})$.
\end{itemize}
We opt for the first interpretation, which is stronger. Then, by interpolation, $\cI(L)$ is an interval and we make the following

\begin{defn}
	\label{def: r+-}
	The lower and upper endpoint of $\cI(L)$ are denoted by $r_{-}(L)$ and $r_+(L)$, respectively.
\end{defn} 

The two interpretations above agree if $\H^p \cap \ran(a{L^{1/2}})$ is dense in $\H^p \cap \L^2$, but \emph{a priori} this information might not be available. It happens for $p \in \cJ(L) \cap (1,\infty)$ though, as the following more general lemma shows.

\begin{lem}
\label{lem: range density in Lp}
If $p \in \cJ(L) \cap (1,\infty)$,  then the spaces $\L^p \cap \dom(L^k) \cap \ran(L^k)$, $k \in \IN$, are all dense in $\L^p \cap \L^2$.
\end{lem}

\begin{proof}
By the Hahn--Banach theorem it suffices to check density for the weak topology. Given $f \in \L^p \cap \L^2$, we consider approximants in $\dom(L^k) \cap \ran(L^k)$ defined by
\begin{align*}
f_j 
&\coloneqq (jL)^k(1+jL)^{-k}(1+j^{-1}L)^{-k} \\
&= (1- (1+jL)^{-1})^k(1+j^{-1}L)^{-k} \quad (j \in \IN).
\end{align*}
By the convergence lemma we have $f_j \to f$ in $\L^2$ as $j \to \infty$. On the other hand, $(f_j)$ is bounded in $\L^p \cap \L^2$ and this space is reflexive since it is isomorphic to a closed subspace of a reflexive space, namely the diagonal in $\L^p \times \L^2$. Hence, it has a weak accumulation point in $\L^p \cap \L^2$, which by $\L^2$-convergence has to be $f$.
\end{proof}

In this section we shall identify $(r_-(L), r_+(L)) \cap (1,\infty)$. Hence, we are studying $\L^p$-boundedness of $\nabla_x L^{-1/2}$. 
Later on, in Section~\ref{sec: Riesz 2}, we will complete the results on the Riesz transform by identifying $\cI(L)$ in the full range of exponents. This will require different methods. 

Here is our main result on the Riesz transform in the $\L^p$-scale.\index{Riesz transform!$\L^p$-boundedness}

\begin{thm}
\label{thm: Riesz}
The endpoints of $\cI(L) \cap (1,\infty)$ can be characterized as follows: 
\begin{align*}
r_-(L) \vee 1 = p_{-}(L) \vee 1 \quad \& \quad r_+(L) = q_{+}(L).
\end{align*}
\end{thm}

Theorem~\ref{thm: Riesz} requires establishing four implications that we shall present in a separate section each. The outline follows \cite[Ch.~5]{A}. To begin with, we need suitable singular integral representations for $R_L$. Let $\alpha \in \IN$. Writing out the Calderón reproducing formula for the auxiliary function $z^{3 \alpha +1/2}(1+z)^{-9\alpha}$ and applying $R_L = \nabla_x L^{-1/2}$ on both sides, we have for all $f \in \L^2$ the representation\index{Riesz transform! singular integral representation} via an improper Riemann integral
\begin{align}
\label{eq: Riesz as singular integral}
R_Lf =  \frac{1}{c_\alpha} \int_0^\infty t \nabla_x (1+t^2L)^{-3\alpha} (t^2L)^{3 \alpha} (1+t^2L)^{-6\alpha}f \, \frac{\d t}{t},
\end{align}
where $c_\alpha$ is a constant depending on $\alpha$. We note that on the right-hand side we do not have to deal with the extension of the square root. More precisely, the truncated Riesz transforms\index{Riesz transform!truncated} defined for $\eps \in (0,1)$ via
\begin{align}
\label{eq: Riesz truncated}
R_L^\eps f 
\coloneqq  \frac{1}{c_\alpha} \nabla_x L^{-1/2}  \int_\eps^{1/\eps} (1+t^2L)^{-3\alpha} (t^2L)^{3 \alpha+1/2} (1+t^2L)^{-6\alpha}f \, \frac{\d t}{t}
\end{align}
converge strongly on $\L^2$ to $R_L$ as $\eps \to 0$. The way to treat the kernel in \eqref{eq: Riesz as singular integral} or \eqref{eq: Riesz truncated} will be through $\L^p - \L^2$ and $\L^2 - \L^p$ off-diagonal bounds that we record in the next lemma.

\begin{lem}
\label{lem: Lp - L2 off diagonal bounds for J(L)}
Let $p \in \cJ(L)$ and let $q$ be between $p$ and $2$. There exists an integer $\beta = \beta(p,q,n)$ with the following property.

\begin{enumerate}
	\item If $p<2$, then $((1+t^2L)^{-\beta})_{t>0}$ and $(t \nabla_x (1+t^2 L)^{-\beta-1})_{t>0}$ satisfy $\L^q - \L^2$ off-diagonal estimates of arbitrarily large order.
	
	\item If $p>2$, then $((1+t^2L)^{-\beta})_{t>0}$ satisfies $\L^2 - \L^q$ off-diagonal estimates of arbitrarily large order.
\end{enumerate}
\end{lem}

\begin{proof}
We begin with (i). The resolvents are $\L^p$-bounded by assumption and $\L^\varrho - \L^2$-bounded for some $\varrho = \varrho(n) \in (1,2)$ due to Lemmata~\ref{lem: Lp L2 between Sobolev conjugates n geq 2} and \ref{lem: Lp L2 between Sobolev conjugates n = 1}. Lemma~\ref{lem: extra} furnishes an integer $\beta = \beta(p,q,n)$ such that $((1+t^2L)^{-\beta})_{t>0}$  is $\L^q-\L^2$-bounded. This holds for all such exponents $q$, so the off-diagonal estimates follow by interpolation with the $\L^2$-result. The claim for the gradients follows by composition since $(t\nabla_x (1+t^2 L)^{-1})$ satisfies $\L^2$ off-diagonal estimates of arbitrarily large order.

As for (ii), we can argue by duality and similarity. Indeed, (i) applies to $L^\sharp = (a^*)^{-1} L^* a^*$ and we have $(2,(p_-(L^\sharp) \vee 1)') = (2, p_+(L))$.
\end{proof}

\subsection{Sufficient condition for \texorpdfstring{$\boldsymbol{1<p<2}$}{1<p<2}}
\label{subsec: Riesz p<2 sufficient}

We prove $(p_{-}(L) \vee 1, 2) \subseteq \cI(L)$. Due to the $\L^2$-bound and the Marcinkiewicz interpolation theorem  it suffices to show that $R_L$ is of weak type $(p,p)$ for every $p \in (p_{-}(L) \vee 1, 2)$. We fix such $p$ and use Blunck--Kunstmann's criterion~\cite{BK}~\index{Theorem!Blunck--Kunstmann's} in its simplified version as stated in \cite[Thm.~1.1]{A}:

\begin{prop}
\label{prop: BK criterion}
Let $p \in [1,2)$. Suppose that $T$ is a sublinear operator of strong type $(2,2)$ and let $A_r$, $r>0$, be a family of bounded linear operators on $\L^2$. Assume for $j \geq 2$ that
\begin{align}
\label{eq: BK1}
\bigg(\frac{1}{|B|}\int_{C_j(B)} |T(1-A_{r(B)})f|^2 \bigg)^{\frac{1}{2}} \leq g(j) \bigg(\frac{1}{|B|} \int_B |f|^p \bigg)^{\frac{1}{p}}
\end{align}
and for $j \geq 1$ that
\begin{align}
\label{eq: BK2}
\bigg(\frac{1}{|B|}\int_{C_j(B)} |A_{r(B)}f|^2 \bigg)^{\frac{1}{2}} \leq g(j) \bigg(\frac{1}{|B|} \int_B |f|^p \bigg)^{\frac{1}{p}},
\end{align}
for all balls $B$ and all $f \in \L^2$ with support in $B$. If $\Sigma \coloneqq \sum_j g(j) 2^{jn/2}$ is finite, then $T$ is of weak type $(p,p)$ with a bound depending on $p$, $\Sigma$ and the strong type $(2,2)$-bound.
\end{prop}

We check \eqref{eq: BK1} and \eqref{eq: BK2} for $T=R_L$ the Riesz transform and 
\begin{align}
\label{eq: BK approximations}
A_r \coloneqq 1-\varphi(r^2L),
\end{align}
where
\begin{align}
\label{eq: phi Riesz p<2}
\varphi(z) \coloneqq (1-(1+z)^{-\beta})^{3\alpha}.
\end{align}
Here, $\alpha \in \IN$ is as in \eqref{eq: Riesz as singular integral}. It will be chosen larger in the further course. Since $p$ is not the lower endpoint of $\cJ(L) \cap (1, 2]$, we can pick $\beta \in \IN$ sufficiently large according to Lemma~\ref{lem: Lp - L2 off diagonal bounds for J(L)} to have $\L^p - \L^2$ off-diagonal estimates of arbitrarily large order for $((1+t^2 L)^{-\beta})_{t>0}$ at our disposal.

\medskip

\noindent \emph{Step 1: Verification of \eqref{eq: BK2} with $g(j) = c2^{-\gamma j}$ and arbitrary $\gamma>0$.} Expanding
\begin{align}
\label{eq: BK approximations expanded}
A_r = -\sum_{k=1}^{3\alpha} \binom{3\alpha}{k} (-1)^k (1+r^2L)^{-\beta k}
\end{align}
and using the $\L^p - \L^2$ off-diagonal decay, we immediately get \eqref{eq: BK2} with $g(j) = c 2^{-j\gamma}$ with $\gamma >0$ as large as we want and $c$ depending on $\alpha, \beta, \gamma$. We take $\gamma > \nicefrac{n}{2}$ to meet the summing condition in Proposition~\ref{prop: BK criterion}.

\medskip

\noindent \emph{Step 2: Verification of \eqref{eq: BK1} with $g(j) = c2^{j(n/2-n/p-\alpha)}$.} Let $B$ be a ball of radius $r>0$ and let $f$ be supported in $B$. We abbreviate $C_j(B)$ by $C_j$ and for $j \geq 2$ we introduce $D_j \coloneqq 2^{j-1}B$. Then 
\begin{align*}
\dist(C_j, D_j) \simeq 2^{j} r \simeq \dist(B, {}^cD_j).
\end{align*}
The representation \eqref{eq: Riesz as singular integral} yields
\begin{align}
\label{eq: Splitting Riesz p<2}
\begin{split}
\|R_L&(1-A_r)f\|_{\L^2(C_j)} \\
&\leq \int_0^\infty \|t\nabla_x (1+t^2L)^{-3\alpha} \ind_{D_j}\psi(t^2L)\varphi(r^2L)f\|_{\L^2(C_j)} \, \frac{\d t}{t} \\
&\quad +\int_0^\infty \|t \nabla_x (1+t^2L)^{-3\alpha} \ind_{{}^c D_j}\psi(t^2L)\varphi(r^2L)f\|_{\L^2(C_j)} \, \frac{\d t}{t},
\end{split}
\end{align}
with an auxiliary function
\begin{align}
\label{eq: psi Riesz p<2}
\psi(z) \coloneqq c_\alpha z^{3\alpha}(1+z)^{-6\alpha}.
\end{align}

From now on we require $3\alpha \geq \beta +1$. Composing $\L^2$ off-diagonal estimates for the resolvents and their gradients and $\L^p - \L^2$ off-diagonal estimates for the $\beta$-th powers of the resolvents, we find that
\begin{align*}
t \nabla_x (1+t^2 L)^{-3\alpha} = t \nabla_x (1+t^2 L)^{-3\alpha+\beta}(1+t^2L)^{-\beta}
\end{align*}
satisfies $\L^p -\L^2$ off-diagonal estimates of arbitrarily large order. Thus,
\begin{align}
\label{eq1: Riesz p<2}
\begin{split}
&\|t \nabla_x (1+t^2L)^{-3\alpha} \ind_{D_j}\psi(t^2 L)\varphi(r^2L)f\|_{\L^2(C_j)} \\
& \quad \lesssim t^{\frac{n}{2}-\frac{n}{p}} \bigg(1+ \frac{2^jr}{t}\bigg)^{-\gamma} \|\psi(t^2L)\varphi(r^2L)f\|_p,
\end{split}
\end{align}
with $\gamma > 0$ at our disposal. From \eqref{eq: phi Riesz p<2} and \eqref{eq: psi Riesz p<2} we can read off the decay properties $\varphi \in \Psi_{3\alpha}^0$ and $\psi \in \Psi_{3\alpha}^{3\alpha}$ on any sector. Thus we find by the third part of Lemma~\ref{lem: functional calculus bounds from J(L) abstract} that
\begin{align}
\label{eq2: Riesz p<2}
\|\psi(t^2L)\varphi(r^2L)f\|_p \lesssim \bigg(\frac{r}{t}\bigg)^{2\alpha} \|f\|_p.
\end{align}
We remark that in applying Lemma~\ref{lem: functional calculus bounds from J(L) abstract} we do not need to switch to an exponent $q \in (p,2]$ since $p$ is not the lower endpoint of $\cJ(L) \cap (1,2]$. The combination of the previous two estimates is
\begin{align*}
\|t \nabla_x (1+t^2L)^{-3\alpha} &\ind_{D_j}\psi(t^2 L)\varphi(r^2L)f\|_{\L^2(C_j)} \\
&\lesssim  \bigg(\frac{r}{t}\bigg)^{2\alpha-\frac{n}{2}+\frac{n}{p}} \bigg(1+ \frac{2^jr}{t}\bigg)^{-\gamma} r^{\frac{n}{2}-\frac{n}{p}} \|f\|_p
\end{align*}
and integrating the resulting bound with respect to $\nicefrac{\d t}{t}$ and changing variables to $s=\nicefrac{2^j r}{t}$ leads us to
\begin{align*}
\int_0^\infty \|t\nabla_x (1+t^2L)^{-3\alpha} \ind_{D_j}\psi(t^2L)\varphi(r^2L)f\|_{\L^2(C_j)} \, \frac{\d t}{t}
\leq g(j) r^{\frac{n}{2}-\frac{n}{p}} \|f\|_p, 
\end{align*}
where
\begin{align*}
g(j) 
\coloneqq 2^{j(\frac{n}{2}-\frac{n}{p}-2\alpha)} \int_0^\infty \frac{s^{2\alpha-\frac{n}{2}+\frac{n}{p}}}{(1+ s)^{\gamma}} \, \frac{\d s}{s}.
\end{align*}
We take $\gamma > 2\alpha -\nicefrac{n}{2}+\nicefrac{n}{p}$ to have a finite integral in $s$ and  $2\alpha > n - \nicefrac{n}{p}$ to take care of the summing condition in Proposition~\ref{prop: BK criterion}. This completes the treatment of the first integral on the right of \eqref{eq: Splitting Riesz p<2}.

For the second integral the roles of uniform boundedness and off-diagonal estimates are reversed. Indeed, as ${}^c D_j$ and $C_j$ intersect, our replacement for \eqref{eq1: Riesz p<2} becomes
\begin{align}
\label{eq3: Riesz p<2}
\begin{split}
&\|t \nabla_x (1+t^2L)^{-3\alpha} \ind_{{}^cD_j}\psi(t^2L)\varphi(r^2L)f\|_{\L^2(C_j)} \\
&\quad \lesssim t^{\frac{n}{2}-\frac{n}{p}} \|\psi(t^2L)\varphi(r^2L)f\|_{\L^p({}^cD_j)}
\end{split}
\end{align}
and from the first and second part of Lemma~\ref{lem: functional calculus bounds from J(L) abstract} we obtain the bound
\begin{align*}
\|\psi(t^2L)\varphi(r^2L)f\|_{\L^p({}^cD_j)} 
\lesssim \begin{cases}
\Big(1+\frac{2^{j-1}r}{t}\Big)^{-6\alpha}\|f\|_p & \text{ if } t \leq r \\
(1+ 2^{j-1})^{-6\alpha}\|f\|_p  & \text{ if } t \geq r.
\end{cases}
\end{align*}
In addition we still have the uniform bound \eqref{eq2: Riesz p<2} and thus, using both estimates raised to the power $\nicefrac{1}{2}$, we have
\begin{align}
\label{eq4: Riesz p<2}
\|\psi(t^2L)\varphi(r^2L)f\|_{\L^p({}^cD_j)}  \lesssim  \bigg(\frac{r}{t}\bigg)^{\alpha}  \bigg(1+ \frac{2^jr}{t}\bigg)^{-3\alpha} \|f\|_p.
\end{align}
We combine the latter estimate with \eqref{eq3: Riesz p<2}, integrate in $t$ and change variables to $s=\nicefrac{2^j r}{t}$ as before in order to obtain
\begin{align*}
\int_0^\infty \|t \nabla_x (1+t^2L)^{-3\alpha} \ind_{{}^cD_j}\psi(t^2L)\varphi(r^2L)f\|_{\L^2(C_j)} \, \frac{\d t}{t}
\lesssim g(j) r^{\frac{n}{2}-\frac{n}{p}} \|f\|_p,
\end{align*}
where this time 
\begin{align*}
g(j) 
\coloneqq 2^{j(\frac{n}{2}-\frac{n}{p}-\alpha)} \int_0^\infty \frac{s^{\alpha-\frac{n}{2}+\frac{n}{p}}}{(1+ s)^{3\alpha}} \, \frac{\d s}{s}.
\end{align*}
We take $2\alpha > \nicefrac{n}{p}-\nicefrac{n}{2}$ to have a finite integral in $s$ and  $\alpha > n - \nicefrac{n}{p}$ to take care of the summing condition in Proposition~\ref{prop: BK criterion}. This completes the treatment of the second integral on the right of \eqref{eq: Splitting Riesz p<2} and also the proof of the weak $(p,p)$-bound for $R_L$ is complete.
\subsection{Sufficient condition for \texorpdfstring{$\boldsymbol{p>2}$}{p>2}}
\label{subsec: Riesz p>2 sufficient}

We prove $(2, q_+(L)) \subseteq \cI(L)$. We let $p \in (2,q_+(L))$ and prove that the Riesz transform $R_L$ is $\L^p$-bounded. We use again the singular integral representation \eqref{eq: Riesz as singular integral} with a parameter $\alpha \in \IN$ to be chosen large in the further course of the proof. 

The kernel of the truncated Riesz transforms $R_L^\eps$ in \eqref{eq: Riesz truncated} given by
\begin{align}
\label{eq: kernel truncated Riesz}
\begin{split}
&t\nabla_x (1+t^2L)^{-3\alpha} (t^2L)^{3 \alpha} (1+t^2L)^{-6\alpha} \\
&\quad = t\nabla_x (1+t^2L)^{-1}\ (1-(1+t^2L)^{-1})^{3\alpha} (1+t^2L)^{-6\alpha + 1}
\end{split}
\end{align}
is $\L^p$-bounded since we have $p < p_+(L)$ by Theorem~\ref{thm: standard relation J(L) and N(L)}. Consequently, each $R_L^\eps$ is $\L^p$-bounded with a bound depending on $\eps$ and it suffices to establish a uniform $\L^p$-bound in order to conclude for $\L^p$-boundedness of $R_L$. To this end we ultimately fix some $p_0 \in (p, q_+(L))$ and employ the following criterion.

\begin{prop}[{\cite[Thm.~1.2]{A}}]
\label{prop: BK replacement}
Let $p_0 \in (2, \infty]$. Suppose that $T$ is a sublinear operator acting on $\L^2$ and let $A_r$, $r>0$, a family of linear operators acting on $\L^2$. Assume
\begin{align}
\label{eq: BK1 replacement}
\bigg(\frac{1}{|B|}\int_B |T(1-A_{r(B)})f|^2 \bigg)^{\frac{1}{2}} \leq C (\Max(|f|^2))^{\frac{1}{2}}(y)
\end{align}
and 
\begin{align}
\label{eq: BK2 replacement}
\bigg(\frac{1}{|B|}\int_B |TA_{r(B)}f|^{p_0} \bigg)^{\frac{1}{p_0}} \leq C (\Max(|Tf|^2))^{\frac{1}{2}}(y)
\end{align}
for all $f \in \L^2$, all balls $B$ and all $y \in B$. If $2<p<p_0$ and $Tf \in \L^p$ whenever $f \in \L^p \cap \L^2$, then 
\begin{align*}
\|Tf\|_p \leq c \|f\|_p,
\end{align*}
where $c$ depends only on $n$, $p$, $p_0$, $C$.
\end{prop}

As $p_0 < p_+(L)$, we can use again Lemma~\ref{lem: Lp - L2 off diagonal bounds for J(L)} to find some large $\beta \in \IN$ for which $((1+t^2L)^{-\beta+1})_{t>0}$ satisfies $\L^2 - \L^{p_0}$ off-diagonal estimates of arbitrarily large order. Then we define the same approximating family $A_r$ as in \eqref{eq: BK approximations} and our task is to verify \eqref{eq: BK1 replacement} and \eqref{eq: BK2 replacement} for $T = R_L^\eps$ and a constant $C$ that does not depend on $\eps$.

We assume right away that $6\alpha \geq \beta$. By composition, this guarantees that the kernel in \eqref{eq: kernel truncated Riesz} is $\L^2 -\L^{p_0}$-bounded and hence that $R_L^\eps$ also maps $\L^2$ into $\L^{p_0}$.

\medskip

\noindent \emph{Step 1: Verification of \eqref{eq: BK1 replacement}.} Let $f \in \L^2$ and $B$ a ball of radius $r$. We claim that
\begin{align}
\label{eq: Claim 1 Riesz p>2}
\|R_L^\eps (1-A_r)f\|_{\L^2(B)} \lesssim r^{\frac{n}{2}} \sum_{j=1}^\infty g(j) \bigg(\barint_{2^{j+1}B} |f|^2 \bigg)^{\frac{1}{2}}
\end{align}
with $g(j) = C 2^{-j(\alpha-\nicefrac{n}{2})}$ and $C$ a constant that does not depend on $\eps$. Since each integral on the right-hand side is bounded by $\Max(|f|^2)(y)$ for every $y \in B$, this bound yields \eqref{eq: BK1 replacement} provided that we take $\alpha > \nicefrac{n}{2}$. 

For the claim we write $f = \sum_{j=1}^\infty f_j$, where $f_j \coloneqq \ind_{C_j}f$ and $C_j \coloneqq C_j(B)$, and obtain
\begin{align*}
\|R_L^\eps (1-A_r)f\|_{\L^2(B)} \leq \sum_{j = 1}^\infty \|R_L^\eps (1-A_r)f_j\|_{\L^2(B)}.
\end{align*}
The term for $j=1$ is readily handled by $\L^2$-boundedness of $R_L^\eps(1-A_r)$:
\begin{align*}
\|R_L^\eps(1-A_r)f_1\|_{\L^2(B)} \lesssim \|f\|_{\L^2(4B)} \simeq r^{\frac{n}{2}} \bigg(\barint_{4B} |f|^2 \bigg)^{\frac{1}{2}}.
\end{align*}
Note that the $\L^2$-bound is independent of $\eps, r$ and depends only on dimensions and ellipticity. This follows from writing $R_L^\eps (1-A_r) = R_L F_{\eps,r}(L)$ as in \eqref{eq: Riesz truncated} and using the functional calculus on $\L^2$. For $j \geq 2$ we re-introduce the auxiliary function $\psi$ from \eqref{eq: psi Riesz p<2} and the sets $D_j = 2^{j-1}B$. In analogy with \eqref{eq: Splitting Riesz p<2} we write
\begin{align*}
&\|R_L^\eps (1-A_r)f_j\|_{\L^2(B)} \\
&\quad\leq \int_0^\infty \|t\nabla_x (1+t^2L)^{-3\alpha} (\ind_{{}^cD_j}+\ind_{D_j})\psi(t^2L)\varphi(r^2L)f_j\|_{\L^2(B)} \, \frac{\d t}{t}.
\intertext{By composition, $t \nabla_x (1+t^2 L)^{-3\alpha}$ satisfies $\L^2$ off-diagonal estimates of arbitrarily large order. Therefore, we continue by}
&\quad\lesssim \int_0^\infty \bigg(1+\frac{2^jr}{t}\bigg)^{-3\alpha} \|\psi(t^2L)\varphi(r^2L)f_j\|_{\L^2({}^c D_j)} \, \frac{\d t}{t}\\ 
&\qquad+ \int_0^\infty \|\psi(t^2L)\varphi(r^2L)f_j\|_{\L^2(D_j)}  \, \frac{\d t}{t}.
\intertext{We can re-use \eqref{eq2: Riesz p<2} with $p=2$ and likewise \eqref{eq4: Riesz p<2} if we replace $({}^cD_j,f)$ with $(D_j, f_j)$ due to the different support properties in the ongoing argument. Indeed, these bounds have been obtained assuming only $p \in (p_{-}(L),2]$. Altogether, we obtain a bound by}
&\quad \lesssim \int_0^\infty \bigg(\frac{r}{t}\bigg)^{2\alpha} \bigg(1+\frac{2^jr}{t}\bigg)^{-3\alpha}  \|f_j\|_2 + \bigg(\frac{r}{t}\bigg)^{\alpha} \bigg(1+\frac{2^jr}{t}\bigg)^{-3\alpha} \|f_j\|_2 \, \frac{\d t}{t} \\
&\quad \leq 2^{-j \alpha} \|f\|_{\L^2(2^{j+1}B)} \int_0^\infty \frac{s^{2\alpha} + s^{\alpha}}{(1+s)^{3\alpha}} \, \frac{\d s}{s},
\end{align*}
where the integral in $s$ is finite. The claim \eqref{eq: Claim 1 Riesz p>2} follows.

\medskip

\noindent \emph{Step 2: Verification of \eqref{eq: BK2 replacement}.} Let $g \in \Wdot^{1,p_0} \cap \W^{1,2}$ and $B$ a ball of radius $r$. We claim that 
\begin{align}
\label{eq: Claim 2 Riesz p>2}
\bigg(\barint_B |\nabla_x A_r g|^{p_0} \bigg)^{\frac{1}{p_0}} \leq C \sum_{j=1}^\infty g(j) \bigg(\barint_{2^{j+1}B} |\nabla_x g|^2 \bigg)^{\frac{1}{2}}
\end{align}
holds with a summable sequence $g(j)$ that does not depend on $\eps$. Taking this for granted, the right-hand side is bounded by $\Max(|\nabla g|^2)(y)^{1/2}$ for every $y \in B$ and, given $f \in \L^2$, the function
\begin{align*}
g 
\coloneqq \frac{1}{c_\alpha} \int_\eps^{1/\eps} t (1+t^2L)^{-3\alpha} (t^2L)^{3 \alpha} (1+t^2L)^{-6\alpha}f \, \frac{\d t}{t}
\end{align*}
verifies $\nabla_x A_r g = R_L^\eps A_r f$ and $\nabla_x g = R_L^\eps f$. At the beginning of the proof we have seen that $R_L^\eps$ maps $\L^2$ into $\L^{p_0}$. Therefore $g \in \Wdot^{1,p_0}$ and we obtain \eqref{eq: BK2 replacement}.

In order to prove \eqref{eq: Claim 2 Riesz p>2}, we perform two more reduction steps. Expanding $A_r$ as in \eqref{eq: BK approximations expanded}, we see that it suffices to establish \eqref{eq: Claim 2 Riesz p>2} with $A_r$ replaced by $(1+r^2L)^{-\beta k}$, $k \geq 1$. Moreover, thanks to the conservation property in Corollary~\ref{cor: conservation} we can replace $g$ by $g-(g)_B$. 

Concerning off-diagonal estimates of arbitrarily large order, we obtain type $\L^2 - \L^{p_0}$ for
\begin{align*}
r\nabla_x(1+r^2L)^{-\beta k} = r\nabla_x(1+r^2L)^{-1}(1+r^2L)^{-\beta k +1}
\end{align*}
by composition: Indeed, for the gradient family we have $\L^{p_0} - \L^{p_0}$ by interpolation of the $\L^2$-result with $\L^q$-boundedness for some $q \in (p_0, q_+(L))$, and $\beta$ was chosen such that already the $(\beta-1)$-th powers of resolvents have type $\L^2 - \L^{p_0}$. As usual, we split $g-g_B = \sum_{j =1}^\infty (g-g_B) \ind_{C_j(B)}$ and obtain
\begin{align*}
\bigg(\barint_B |\nabla_x (1+r^2 L)^{-\beta} (g-g_B)|^{p_0} \bigg)^{\frac{1}{p_0}}
&\lesssim r^{-1-\frac{n}{2}} \sum_{j = 1}^\infty 2^{-j \gamma} \|g-g_B\|_{\L^2(C_j(B))},
\end{align*}
where $\gamma >0$ is at our disposal. Poincaré's inequality \cite[Prop.~7.45]{GT} provides the bound
\begin{align*}
\|g-g_B\|_{\L^2(C_j(B))} \leq \|g-g_B\|_{\L^2(2^{j+1} B)} \lesssim 2^{jn} r \|\nabla_x g\|_{\L^2(2^{j+1} B)}.
\end{align*}
We conclude that
\begin{align*}
\bigg(\barint_B |\nabla_x (1+r^2 L)^{-\beta} (g-g_B)|^{p_0} \bigg)^{1/p_0} \lesssim \sum_{j = 1}^\infty 2^{j(\frac{3n}{2}-\gamma)}  \bigg(\barint_{2^{j+1}B} |\nabla_x g|^2 \bigg)^{1/2}.
\end{align*}
We take $\gamma > \nicefrac{3n}{2}$ to grant summability of $g(j) \coloneqq 2^{j(3n/2-\gamma)}$ and the proof of \eqref{eq: Claim 2 Riesz p>2} is complete.
\subsection{Necessary condition for \texorpdfstring{$\boldsymbol{1<p<2}$}{p<2}.}
\label{subsec: Riesz p<2 necessary}

We suppose that the Riesz transform is $\L^p$-bounded for some $p \in (1,2)$ and prove that $p \geq p_-(L)$. In dimension $n\leq2$ we have $p_-(L) \leq  1$, see Proposition~\ref{prop: J(L) contains neighborhood of Sobolev conjugates}. Hence, we can restrict ourselves to dimensions $n \geq 3$.

We set $p_0 \coloneqq p$, define iteratively $p_k \coloneqq p_{k-1}^*$ and stop at the first exponent $k^+ \geq 0$ with $p_{k^+} \in (2_*,2]$. We shall prove by backward induction that $(p_k,2] \subseteq \cJ(L)$ for all $k$. Hence, we eventually find $(p,2] \subseteq \cJ(L)$, that is to say, $p \geq p_-(L)$.

We have $(p_{k^+},2] \subseteq (2_*,2] \subseteq \cJ(L)$ by Proposition~\ref{prop: J(L) contains neighborhood of Sobolev conjugates}. For the inductive step we assume $(p_k,2] \subseteq \cJ(L)$ and pick any $q \in (p_{k-1},2_*]$. Then $q^*$ is an interior point of $\cJ(L)$ and hence $(tL^{1/2} (1+t^2L)^{-1})_{t>0}$ is $\L^{q^*}$-bounded by Lemma~\ref{lem: functional calculus bounds from J(L) abstract}. For $f \in \L^q \cap \ran({L^{1/2}})$ we can therefore estimate
\begin{align*}
\|(1+t^2 L)^{-1}f \|_{q^*}
&\lesssim t^{-1} \|L^{-1/2}f \|_{q^*} \\
&\lesssim t^{-1} \|\nabla_x L^{-1/2}f\|_q \\
&\lesssim t^{-1} \|f\|_q,
\end{align*}
where the final step uses $q \in (p,2] \subseteq \cI(L)$. We need to make sure that this estimates applies to sufficiently many functions $f$. We stress that Lemma~\ref{lem: range density in Lp} is useless in this regard since $q \in \cJ(L)$ is precisely what we are trying to prove.

\begin{lem}
\label{lem: density for Riesz}
In any dimension $n$, it follows that if $q \in \cI(L)$ satisfies $q < 2_*$, then $\H^q \cap \L^2 \subseteq \ran(a{L^{1/2}})$.
\end{lem}

Momentarily, let us take the lemma for granted. If $q>1$, then multiplication by $a$ is an automorphism of $\L^q \cap \L^2$. Hence, we have $\L^q \cap \L^2 \subseteq \ran({L^{1/2}})$ and the previous bound implies $\L^q - \L^{q^*}$-boundedness of the resolvents. As usual, we can interpolate with the $\L^2$ off-diagonal estimates and then use Lemma~\ref{lem: OD implies boundedness} to obtain $(q,2] \subseteq \cJ(L)$. Since $q \in (p_{k-1},2_*]$ was arbitrary, $(p_{k-1},2] \subseteq \cJ(L)$ follows.

This completes the proof modulo the following:

\begin{proof}[Proof of Lemma~\ref{lem: density for Riesz}]
For clarity we denote by $T_L$ the extension of the bijection ${L^{1/2}}: \W^{1,2} \to \ran({L^{1/2}})$ to an isomorphism $\Wdot^{1,2} \to \L^2$, so that $R_L = \nabla_x T_L^{-1}$. 

Let $f \in \H^q \cap \L^2$. Interpolation yields $f \in \H^{2_*} \cap \L^2$ and $2_* \in \cI(L)$. Hence, $\nabla_x T_L^{-1} a^{-1} f \in \H^{2_*}$. Modulo constants we obtain $T_L^{-1} a^{-1} f \in \L^2$ by the Hardy--Sobolev embedding and consequently $T_L^{-1} a^{-1} f \in \W^{1,2}$. By definition of $T_L$ this means that $a^{-1}f \in \ran({L^{1/2}})$.
\end{proof}
\subsection{Necessary condition for \texorpdfstring{$\boldsymbol{p>2}$}{p>2}}
\label{subsec: Riesz p>2 necessary}

We let $p \in (2, r_+(L))$ and prove that $p \leq q_+(L)$. In fact, it suffices to prove $[2,p) \subseteq \cJ(L)$: For $q \in (2,p)$ we then obtain $\L^q$-boundedness of
\begin{align*}
t \nabla_x (1+t^2L)^{-1} = (\nabla_x L^{-1/2})((t^2L)^{1/2} (1+t^2L)^{-1})
\end{align*}
by composition, applying Lemma~\ref{lem: functional calculus bounds from J(L) abstract} to the second factor. 

The argument for $[2,p) \subseteq \cJ(L)$ is similar to the previous section. We set $p_0 \coloneqq p$, define iteratively $p_k \coloneqq (p_{k-1})_*$ and stop at the first exponent $k^- \geq 0$ with $p_{k^-} \in [2,2^*)$. Then $[2,p_{k-}) \subseteq \cJ(L)$ by Proposition~\ref{prop: J(L) contains neighborhood of Sobolev conjugates}. Now, assume $[2,p_k) \subseteq \cJ(L)$ and pick any $q \in [2^*,p_{k-1})$. Since $q_*$ is an interior point of $\cJ(L)$, the family $(tL^{1/2} (1+t^2L)^{-1})_{t>0}$ is $\L^{q_*}$-bounded by Lemma~\ref{lem: functional calculus bounds from J(L) abstract}. Moreover, $q_* \in [2,p_k) \subseteq \cI(L)$, so for all $f \in \L^{q_*} \cap \L^2$, we get
\begin{align*}
\|(1+t^2 L)^{-1}f \|_{q}
&\lesssim \|\nabla_x (1+t^2 L)^{-1}f \|_{q_*}\\
&\lesssim \|L^{1/2}  (1+t^2 L)^{-1}f\|_{q_*} \\
&\lesssim t^{-1} \|f\|_{q_*},
\end{align*}
which shows that $((1+t^2 L)^{-1})_{t>0}$ is $\L^{q_*} - \L^q$-bounded. Interpolation with the $\L^2$ off-diagonal estimates and then Lemma~\ref{lem: OD implies boundedness} yield $[2,q) \subseteq \cJ(L)$. Since $q \in [2^*, p_{k-1})$ was arbitrary, $[2,p_{k-1}) \subseteq \cJ(L)$ follows. By backward induction we eventually arrive at the desired conclusion $[2,p) \subseteq \cJ(L)$.
\section{Operator-adapted spaces} 
\label{sec: Hardy intro}

\noindent Operator-adapted Hardy--Sobolev spaces are our main tool in this monograph and will be essential for understanding most of the following sections. They have been developed in various references starting with semigroup generators in \cite{AMcR, HM, HMMc, Duong-Li} up to the recent monographs focusing on bisectorial operators~\cite{AA,AusSta}. Still we need some unrevealed features and we take this opportunity to correct some inexact arguments from the literature. 

For general properties of adapted Hardy spaces we closely follow \cite[Sec.~3]{AA}, where the authors develop an abstract framework of two-parameter operator families that provides a unified approach to sectorial and bisectorial operators. The application to bisectorial operators with first-order scaling has been detailed in \cite[Sec.~4]{AA} and we review their results in Section~\ref{subsec: Hardy abstract}. Section~\ref{subsec: Hardy abstract sectorial} provides all necessary details in order to apply the framework to sectorial operators with second-order scaling and we summarize the results that are relevant to us. This will justify using parts of \cite{AA} for sectorial operators in the further course.

The abstract framework allows us to treat operator-adapted Besov spaces simultaneously without any additional effort. These spaces will only be needed in the final Section~\ref{sec: fractional} and the reader might ignore them till then. 
\subsection{Bisectorial operators with first-order scaling}
\label{subsec: Hardy abstract}

To set the stage, we assume that 
\begin{align}
\label{eq: standard assumptions bisectorial}
\begin{minipage}{0.89\linewidth}
\begin{itemize}
\item $T$ is a bisectorial operator in $\L^2 = \L^2(\R^n; V)$ of some angle $\omega \in [0, \frac{\pi}{2})$, where $V$ is a finite-dimensional Hilbert space,
\item $T$ has a bounded $\H^\infty$-calculus on $\cl{\ran(T)}$,
\item $((1+ \i  tT)^{-1})_{t \in \R \setminus \{0\}}$ satisfies $\L^2$ off-diagonal estimates of arbitrarily large order.
\end{itemize}
\end{minipage}
\end{align}
These are called \emph{Standard Assumptions}\index{standard assumptions!for operator-adapted spaces} in \cite[Ch.~4]{AA}. In fact, \cite{AA} requires for all $\nu \in (0, \nicefrac{\pi}{2} - \omega)$ that the family $((1+ \i  zT)^{-1})_{z \in \S_\nu}$ satisfies $\L^2$ off-diagonal estimates of arbitrarily large order but this follows already from the first and third assumption in \eqref{eq: standard assumptions bisectorial} by interpolation, see Lemma~\ref{lem: OD extrapolation to sectors}. The reader may recall from Sections~\ref{subsec: adjoints} and \ref{sec: Hp Hq boundedness} that $T^*$ satisfies the standard assumptions as well. 

In the following we suppress the reference to bisectors from notation of classes of holomorphic functions since we allow any bisector of angle larger than $\omega$. We mimic the extension to the upper half-space by convolutions in the definition of the classical Hardy spaces by associating with each $\psi \in \H^\infty$ on a bisector the \emph{extension operator}\index{$\IQ_{\psi,T}$ (extension operator)}
\begin{align}
\label{eq: def Q extension}
\IQ_{\psi,T} : \cl{\ran(T)} \to \L^\infty(0,\infty; \L^2), \quad (\IQ_{\psi,T}f)(t) = \psi(tT)f.
\end{align}
If in addition $\psi \in \Psi_+^+$, then $\IQ_{\psi,T}$ is defined on all of $\L^2$ and by McIntosh's theorem it maps $\L^2$ boundedly into $\L^2(0,\infty, \frac{\d t}{t}; \L^2) = \T^{0,2} = \Z^{0,2}$. Hence, we can look at the bounded dual operator 
\begin{align*}
\IC_{\psi, T} \coloneqq (\IQ_{\psi^*, T^*})^*: \L^2(0,\infty, \tfrac{\d t}{t}; \L^2) \to \L^2,
\end{align*}
where $\psi^*(z) = \cl{\psi(\cl{z})}$, which is given by the weakly convergent integral
\begin{align}
\label{eq: def C contraction}
\IC_{\psi,T}F = \int_0^\infty \psi(tT)F(t) \, \frac{\d t}{t}.
\end{align}
Of course, the integral converges strongly in  $\L^2$ if $F$ has compact support in $(0,\infty)$. We call $\IC_{\psi, T}$\index{$\IC_{\psi, T}$ (contraction operator)} a \emph{contraction operator}. It is denoted by $\IS_{\psi,T}$ in \cite{AA} and we change notation in order to distinguish it from conical square functions.

\begin{defn}
\label{def: Hardy spaces operator}
Let $\psi \in \H^\infty$, $s \in \R$ and $p\in (0,\infty)$. The sets
\begin{align*}
\IH^{s,p}_{\psi,T} &\coloneqq \{f \in \cl{\ran(T)}: \IQ_{\psi,T}f \in \T^{s,p} \cap \T^{0,2} \},\\
\IB^{s,p}_{\psi,T} &\coloneqq \{f \in \cl{\ran(T)}: \IQ_{\psi,T}f \in \Z^{s,p} \cap \Z^{0,2} \},
\end{align*}
equipped with what will be shown to be quasinorms
\begin{align*}
\|f\|_{\IH^{s,p}_{\psi,T}} \coloneqq \|\IQ_{\psi,T}f\|_{\T^{s,p}}, \quad
\|f\|_{\IB^{s,p}_{\psi,T}} \coloneqq \|\IQ_{\psi,T}f\|_{\Z^{s,p}},
\end{align*}
are called \emph{pre-Hardy--Sobolev}\index{Hardy--Sobolev space!adapted to a bisectorial operator} and \emph{pre-Besov space}\index{Besov space!adapted to a bisectorial operator} space of smoothness $s$ and integrability $p$ adapted to $T$, respectively. The function $\psi$ is called \emph{auxiliary function}. 
\end{defn}

In order to treat pre-Hardy--Sobolev and pre-Besov spaces simultaneously, we introduce the concise notation\index{Y@$(\Y, \IX)$ (either $(\T, \IH)$ or $(\Z, \IB)$)}
\begin{align*}
	\IX^{s,p}_{\psi,T} &\coloneqq \{f \in \cl{\ran(T)}: \IQ_{\psi,T}f \in \Y^{s,p} \cap \Y^{0,2} \},
\end{align*}
where the pair $(\Y, \IX)$ is either $(\T, \IH)$ or $(\Z, \IB)$. These pairs are called $(\X,\IX)$ in \cite{AA} but it will be convenient to keep the symbol $\X$ for a different purpose. For $\psi \in \Psi_+^+$  the condition $\IQ_{\psi,T}f \in \Y^{0,2}$ is redundant and if in addition $\psi$ is non-degenerate, then by McIntosh's theorem we have up to equivalent norms
\begin{align}
\label{eq: Hardy identification for p=2}
\IX^{0,2}_{\psi,T} = \cl{\ran(T)}.
\end{align}
For general values of $s$ and $p$ and auxiliary functions $\psi \in \H^\infty = \Psi_0^0$ we still have that $\IX^{s,p}_{\psi, T}$ is quasinormed~\cite[Prop.~4.3]{AA} and, up to equivalent quasinorms, independent of the auxiliary function in the following classes.

\begin{prop}[{\cite[Prop.~4.4]{AA}}]
\label{prop: auxiliary function bisectorial} Let $s \in \R$ and $p\in (0,\infty)$.
Up to equivalent norms, $\IX^{s,p}_{\psi,T}$ does not depend on the choice of $\psi \in \H^\infty$ as long as it is non-degenerate and of class $\Psi_\sigma^\tau$ with the following technical conditions on the decay parameters\index{admissible auxiliary function}:
\begin{itemize}
	\item $\tau > -s + |\nicefrac{n}{2}-\nicefrac{n}{p}|$ and $\sigma > s$ if $p \leq 2$,
	\item $\tau > -s$ and $\sigma> s + |\nicefrac{n}{2}-\nicefrac{n}{p}|$ and moreover $\psi \in \Psi_+^+$ if $p \geq 2$.
\end{itemize}
\end{prop}

This allows us to drop the dependence on $\psi$. 

\begin{defn}
\label{def: Hardy spaces independent of psi} 
Let $s \in \R$ and $p\in (0,\infty)$. Denote by $\IX^{s,p}_T$ the quasinormed space $\IX^{s,p}_{\psi,T}$ for any $\psi \in \Psi_\sigma^\tau$ as in Proposition~\ref{prop: auxiliary function bisectorial}. When $s=0$, simply write $\IX^p_T \coloneqq \IX^{0,p}_T$.
\end{defn}

Usually, we take $\psi$ with sufficiently large decay to describe these spaces.

\begin{prop}[{\cite[Prop.~4.7]{AA}}]
\label{prop: Hardy via contraction}
Let $s \in \R$ and $p\in (0,\infty)$ and suppose that $\psi \in \Psi_\sigma^\tau \cap \Psi_+^+$ is non-degenerate, where
\begin{itemize}
	\item $\tau > s$ and $\sigma > -s + |\nicefrac{n}{2}-\nicefrac{n}{p}|$ if $p \leq 2$,
	\item $\tau> s + |\nicefrac{n}{2}-\nicefrac{n}{p}|$ and $\sigma > -s$ if $p \geq 2$.
\end{itemize}
Then $\IX^{s,p}_T = \IC_{\psi, T}(\Y^{s,p} \cap \Y^{0,2})$ and 
\begin{align*}
f \mapsto \inf \big\{\|F\|_{\Y^{s,p}} : F \in \Y^{s,p} \cap \Y^{0,2} \, \& \; \, \IC_{\psi, T}F = f \big\}
\end{align*}
is an equivalent quasinorm.
\end{prop}

The spaces $\IX_{T}^{s,p}$ are not complete in general unless $p=2$. This is why we use the subscript `pre' and remove it when taking completions.  As usual, a \emph{completion} of a quasinormed space $Q$ is a linear isometric map $\iota: Q \to \hat{Q}$, where $\hat{Q}$ is a complete quasinormed space and $\iota(Q)$ is dense in $\hat{Q}$. For $Q \coloneqq \IX^{s,p}_{T}$, there are compatible completions of these spaces within the same ambient space $\Lloc^2(\reu)$: the construction in \cite[Prop.~4.20]{AA}, called \emph{canonical completion}, is to take
\begin{align*}
	\iota \coloneqq \IQ_{\psi,T} \text{ with } \psi \in \Psi_\infty^\infty \quad \& \quad \hat{Q} \coloneqq \cl{\IQ_{\psi,T}(\IX^{s,p}_{T})} \subseteq \Y^{s,p}.
\end{align*}

\begin{defn}
	\label{def: Hardy spaces operator completed}
	Let $\psi \in \Psi_\infty^\infty$ be non-degenerate. For $s \in \R$ and $p\in (0,\infty)$ denote by $\psiX^{s,p}_T$ the canonical completion\index{canonical completion ($\psiX^{s,p}_T$)} of the quasinormed space $\IX^{s,p}_{T}$.
\end{defn}

By the Calderón reproducing formula (here for bisectorial operators, see \cite[Prop.~4.2]{AusSta}) the function $\psi$ has a non-degenerate sibling\index{sibling!of an auxiliary function} $\varphi \in \Psi_\infty^\infty$ such that $\IC_{\varphi,T} \IQ_{\psi,T} = 1$ on $\cl{\ran(T)}$. This allows us to summarize the full construction of operator adapted Hardy spaces in one commutative diagram, see Figure~\ref{fig: completions}.

\begin{figure}[ht]
\centering
\begin{center}
	\begin{tikzcd}[column sep=40pt]
		\psiX^{s,p}_{ T} \arrow[hookrightarrow]{r} 
		&\Y^{s,p} \arrow[r, "P"] 
		& \psiX^{s,p}_T  \arrow[hookrightarrow]{r} 
		&\psiX^{s,p}_T 
		\\
		\IQ_{\psi,T}(\IX^{s,p}_{T}) \arrow[hookrightarrow]{u} \arrow[hookrightarrow]{r} \arrow[hookrightarrow, bend right = 30, "\mathrm{Id}"]{rrr}
		& \Y^{s,p} \cap \Y^{0,2} \arrow[hookrightarrow]{u} \arrow[r, "\IC_{\varphi, T}", swap]
		&\IX_{T}^{s,p} \arrow[u, "\IQ_{\psi, T}"] \arrow[r, "\IQ_{\psi, T}", swap]
		& \IQ_{\psi,T}(\IX_{ T}^{s,p}) \arrow[hookrightarrow]{u}
	\end{tikzcd}
\end{center}
\caption{Canonical completion: $\varphi, \psi \in \Psi_\infty^\infty$ are siblings and $P$ is the unique bounded linear map for which the diagram commutes. It follows that $P$ is a projection from $\Y^{s,p}$ onto $\psiX^{s,p}$. By the universal approximation technique for $\Y$-spaces, projections for different choices of admissible spaces are compatible. The bottom part of the diagram also identifies $\psiX^{s,p} \cap \IQ_{\psi,T}(\IX^{0,2}_T) = \IQ_{\psi,T}(\IX^{s,p}_T$).}
\label{fig: completions}
\end{figure}

The canonical completions inherit many properties of tent and $\Z$-spaces via Figure~\ref{fig: completions}. Two important examples are the following approximation results that have been tacitly used in \cite{AA}.\index{universal approximation technique!for $\X^{s,p}_T$} By a slight abuse of notation we allow $\IX \in \{\IB, \IH\}$ to be different in the assumption and the conclusion.

\begin{lem}
\label{lem: universal approximation for canonical completions}
Let $\psi \in \Psi_\infty^\infty$ be non-degenerate. If $F \in \psiX_T^{s_0,p_0}$ for some $s_0 \in \R$, $p_0 \in (0,\infty)$, then there exists $(F_j)_j \subseteq \psiX_T^{0,2}$  with $F_j \to  F$ in every space of type $\psiX_T^{s,p}$ that $F$ belongs to.
\end{lem}

\begin{proof}
This is an immediate consequence of Figure~\ref{fig: completions}. Indeed, since $\ind_{(j^{-1},j) \times B(0,j)} F \in \Y^{0,2}$ is a universal approximation of $F$ with respect to tent and $\Z$-spaces, see Sections~\ref{subsec: tent spaces} and \ref{subsec: Z spaces}, we can take $F_j \coloneqq P (\ind_{(j^{-1},j) \times B(0,j)} F )$.
\end{proof}

\begin{lem}
\label{lem: domain dense in Hardy}
Let $s_0 \in\R$ and $p_0 \in (0,\infty)$. Given $f \in \IX^{s_0,p_0}_T$, there is a sequence $(f_j)_j$ in $\bigcap_{k \in \IZ} \ran(T^k)$ with $f_j \to f$ in every space of type $\IX^{s,p}_T$ that $f$ belongs to. In particular, convergence holds in $\IX^{0,2}_T \subseteq \L^2$.  \index{universal approximation technique!for $\IX_T^{s,p}$}
\end{lem}

\begin{proof}
Again by Figure~\ref{fig: completions} we have $f = \IC_{\varphi,T}F$ with $F \coloneqq \IQ_{\psi,T} f$ and therefore $f_j \coloneqq \IC_{\varphi, T} (\ind_{(j^{-1},j) \times B(0,j)} F)$ have the required universal approximation property. Thanks to $\varphi \in \Psi_\infty^\infty$ we also obtain that 
	\begin{align*}
		f_j
		= T^k \int_{j^{-1}}^j (tT)^{-k}\varphi(tT) (\ind_{B(0,j)} F(t)) \, \frac{\d t}{t^{-k+1}} \in \ran(T^k) \quad (k \in \IZ). &\qedhere
	\end{align*}
\end{proof}

One necessity for the canonical completions is the following interpolation result.\index{interpolation!of operator-adapted spaces $\psiX_T^{s,p}$}

\begin{prop}[{\cite[Thm.~4.28]{AA}}]
\label{prop: interpolation of canonical completions}
Let $\psi \in \Psi_\infty^\infty$ be non-degenerate. Let $0<p_0, p_1 <\infty$, $s_0, s_1 \in \R$, $\theta \in (0,1)$ and set $p \coloneqq [p_0, p_1]_\theta$, $s \coloneqq (1-\theta)s_0 + \theta s_1$. Up to equivalent quasinorms it follows that
\begin{align*}
	[\psiH_T^{s_0,p_0}, \psiH_T^{s_1,p_1}]_\theta &= \psiH_T^{s,p}
\intertext{and if $s_0 \neq s_1$,}	
	(\psiX_T^{s_0,p_0}, \psiX_T^{s_1,p_1})_{\theta,p} &= \psiB_T^{s,p}.
\end{align*}
\end{prop}

When $p \in (1,\infty)$, the spaces $\psiX^{p}_T$ and $\varphiX^{p'}_{T^*}$ are in natural duality with each other as described in \cite[Prop.~4.23]{AA} provided that $\varphi, \psi \in \Psi_\infty^\infty$ are siblings. Since by definition the pre-Hardy--Sobolev and pre-Besov spaces are dense in their completions, we can equivalently state this result as follows.\index{duality!for $\IX^{s,p}_T$}

\begin{prop}
\label{prop: Hardy duality}
Let $p \in (1,\infty)$. Then, whenever $f \in \cl{\ran(T)}$,
\begin{align*}
\sup_{g \in \IX^{p'}_{T^*}} \frac{|\langle f, g \rangle|}{\|g\|_{\IX^{p'}_{T^*}}} \simeq \|f\|_{\IX^p_T},
\end{align*}
where $\langle \cdot\, , \cdot \rangle$ is the inner product on $\L^2$ and the right-hand side is interpreted as $\infty$ if $f \notin \IX_T^p$.
\end{prop}

The `raison d'être' of these spaces is that the $\H^\infty$-calculus of $T$ extends to them in the best possible way.\index{Hinfty calculus@$\H^\infty$-calculus!on $\IX_T^{s,p}$}

\begin{prop}[{\cite[Thm.~4.14]{AA}.}]
\label{prop: FC on pre-Hardy}
Let $s \in \R$, $p \in (0,\infty)$ and $\nu \in (\omega, \nicefrac{\pi}{2})$. Then for all $\eta \in \H^\infty(\S_\nu)$,
\begin{align*}
\|\eta(T) f\|_{\IX^{s,p}_T} \lesssim \|\eta\|_\infty \|f\|_{\IX^{s,p}_T} \quad (f \in \IX^{s,p}_T).
\end{align*}
Moreover, if $\varphi \in \Psi_{-1}^1(\S_\nu)$ and $\psi \in \Psi_{1}^{-1}(\S_\nu)$, then
\begin{align*}
\|\varphi(T) f\|_{\IX^{s+1,p}_T} &\lesssim \|f\|_{\IX^{s,p}_T} \quad (f \in \dom(\varphi(T)) \cap \IX^{s,p}_T )
\intertext{and}
\|\psi(T) f\|_{\IX^{s-1,p}_T} &\lesssim \|f\|_{\IX^{s,p}_T} \quad (f \in \dom(\psi(T)) \cap \IX^{s,p}_T ),
\end{align*}
where the implicit constants also depend on $\varphi$ and $\psi$.
\end{prop}

The second part indicates that the spaces for different smoothness parameters are related through a \emph{lifting property}\index{lifting property!for $\IH_T^{s,p}$}. Indeed, recall that $(\nicefrac{z}{[z]})(T)$ and its inverse are bounded operators on $\cl{\ran(T)}$ since $T$ has a bounded $\H^\infty$-calculus and that therefore $T$ and $[T]$ share the same domain and range. Thus, using $(\varphi,\psi) = (\nicefrac{1}{\psi},\psi)$ with either $\psi(z)=z$ or $\psi(z) = [z]$ in the proposition above, we obtain

\begin{cor}
\label{cor: Lifting pre-Hardy}
The operators $T$ and $(T^2)^{1/2}$ are bijections $\IX^{s+1,p}_T \cap \dom(T) \to \IX^{s,p}_T \cap \ran(T)$ that satisfy
\begin{align*}
\|T f\|_{\IX^{s,p}_T} \simeq \|f\|_{\IX^{s+1,p}_T} \simeq \|[T]f\|_{\IX^{s,p}_T}.
\end{align*}
\end{cor}

From the $\H^\infty$-calculus we immediately obtain that $(\e^{-t[T]})_{t\geq 0}$ is a bounded semigroup on $\IX_T^{s,p}$. In fact, we also have strong continuity and stability.

\begin{prop}[{\cite[Prop.~4.33]{AA}}]
\label{prop: C0 smg on abstract Hardy}
Let $s \in \R$ and $p \in (0,\infty)$. For all $f \in \IX_T^{s,p}$ the following limits hold in $\IX_T^{s,p}$:
\begin{align*}
\lim_{t \to 0} \e^{-t[T]}f = f \quad \text{and} \quad \lim_{t \to \infty} \e^{-t[T]}f = 0.
\end{align*}
\end{prop}
\subsection{Sectorial operators with second-order scaling}
\label{subsec: Hardy abstract sectorial}

In this case our standard assumptions\index{standard assumptions!for operator-adapted spaces} are that
\begin{align}
\label{eq: standard assumptions sectorial}
\begin{minipage}{0.89\linewidth}
\begin{itemize}
\item $T$ is a sectorial operator on $\L^2 = \L^2(\R^n; V)$ of some angle $\omega \in [0,\pi)$, where $V$ is a finite-dimensional Hilbert space,
\item $T$ has a bounded $\H^\infty$-calculus on $\cl{\ran(T)}$,
\item $((1+t^2 T)^{-1})_{t>0}$ satisfies $\L^2$ off-diagonal estimates of arbitrarily large order,
\end{itemize}
\end{minipage}
\end{align}
and we allow holomorphic functions on any sector of angle larger than $\omega$ in the following considerations.

We define the extension\index{$\IQ_{\psi,T}$ (extension operator)} for $\psi \in \H^\infty$ with second-order scaling 
\begin{align*}
\IQ_{\psi,T} : \cl{\ran(T)} \to \L^\infty(0,\infty; \L^2), \quad (\IQ_{\psi,T}f)(t) = \psi(t^2T)f
\end{align*}
and if in addition $\Psi_+^+$, then $\IQ_{\psi,T}$ is again defined on all of $\L^2$, maps into $\L^2(0,\infty, \frac{\d t}{t}; \L^2)$ and we have the dual operator
\begin{align*}
\IC_{\psi, T} \coloneqq (\IQ_{\psi^*, T^*})^*, \quad \IC_{\psi,T}F = \int_0^\infty \psi(t^2T)F(t) \, \frac{\d t}{t},
\end{align*}
where the integral converges weakly in $\L^2$. 

Most of the theory in \cite[Sec.~3 \& 4]{AA} has been written for abstract continuous two-parameter families $(S_{t,\tau})_{t, \tau > 0}$ on $\L^2$ and hence applies \emph{in extenso} to families
\begin{align}
\label{eq: generic second order family}
(\psi(t^2 T) \eta(T) \varphi(\tau^2 T))
\end{align}
with a sectorial operator as above, instead of
\begin{align}
\label{eq: generic first order family}
(\psi(t T) \eta(T) \varphi(\tau T))
\end{align}
with a bisectorial operator. Here, $\psi \in \Psi_{\sigma_1}^{\tau_1}$, $\varphi \in \Psi_{\sigma_2}^{\tau_2}$, $\eta \in \Psi_{\sigma_3}^{\tau_3}$ are auxiliary functions with $\sigma_j, \tau_j \in \R$. The only difference with the results of bisectorial operators lies in how large these parameters have to be in order to arrive at the desired conclusion. \index{Hardy--Sobolev space!adapted to a sectorial operator}\index{Besov space!adapted to a sectorial operator}

The three fundamental mapping properties for families of type \eqref{eq: generic first order family} in \cite{AA} -- Lemma~3.17, Lemma~3.18 and Theorem~3.19 -- remain to hold for families of type \eqref{eq: generic second order family} and then the same conclusion holds already if one replaces $\sigma_j, \tau_j$ by $\nicefrac{\sigma_j}{2}, \nicefrac{\tau_j}{2}$ in the assumptions. Indeed, following the self-contained proofs in \cite{AA}, one readily sees that the assumptions on the auxiliary functions are exclusively determined by \cite[Thm.~3.8]{AA}, which in turn provides the order of $\L^2$ off-diagonal decay that one can get for families of the form $(\eta(t)\psi(t T))_{t>0}$ if $(\eta(t))_{t>0}$ is a continuous bounded family of functions in $\H^\infty$ and $\psi \in \Psi_{\sigma}^\tau$. Precisely, \cite[Thm.~3.8]{AA} allows any order up to $\gamma = \sigma$. On the other hand, in Lemma~\ref{lem: functional calculus bounds from J(L) abstract} we have proved the same conclusion for $(\eta(t) \psi(t^2 T))_{t>0}$ under the mere assumption $\psi \in \Psi_{\sigma/2}^{\tau/2}$. 

From this discussion we conclude that qualitatively the results of Section~\ref{subsec: Hardy abstract} that build on \cite{AA} remain valid for sectorial operators with second-order scaling but there are the following quantitative changes. The technical conditions\index{admissible auxiliary function} of Proposition~\ref{prop: auxiliary function bisectorial} become
\begin{itemize}
\item $\tau > -\nicefrac{s}{2} + |\nicefrac{n}{4}-\nicefrac{n}{(2p)}|$ and $\sigma > \nicefrac{s}{2}$ if $p \leq 2$,
\item $\tau > -\nicefrac{s}{2}$ and $\sigma> \nicefrac{s}{2} + |\nicefrac{n}{4}-\nicefrac{n}{(2p)}|$ and moreover $\psi \in \Psi_+^+$ if $p \geq 2$,
\end{itemize}
with the same type of modification in Proposition~\ref{prop: Hardy via contraction}. In Proposition~\ref{prop: FC on pre-Hardy} the assumption on the angle is again best possible, that is $\nu \in (\omega,\pi)$ and $\eta \in \H^\infty(\S_\nu^+)$, and the second part holds for $\varphi \in \Psi_{-1/2}^{1/2}(\S_\nu^+)$ and $\psi \in \Psi_{1/2}^{-1/2}(\S_\nu^+)$. As a consequence, the lifting property\index{lifting property!for $\IH_T^{s,p}$} of Corollary~\ref{cor: Lifting pre-Hardy} uses $\sqrt{T}$.

Performing only the purely symbolic replacement of $\sqrt{z^2}$ by $\sqrt{z}$ at all occasions in the statement and proof of Proposition~\ref{prop: C0 smg on abstract Hardy}, we immediately obtain the following version for sectorial operators.

\begin{prop}
\label{prop: C0 smg on abstract Hardy sectorial}
Let $s \in \R$ and $p \in (0,\infty)$. For all $f \in \IX_T^{s,p}$ the following limits hold in $\IX_T^{s,p}$:
\begin{align*}
\lim_{t \to 0} \e^{-t\sqrt{T}}f = f \quad \text{and} \quad \lim_{t \to \infty} \e^{-t\sqrt{T}}f = 0.
\end{align*}
\end{prop}
\subsection{Molecular decomposition for adapted Hardy spaces}
\label{subsec: Molecular decomposition}

Molecular decompositions for $\IH^{p}_T$ with $p \in (0,1]$ have been pioneered in \cite{HM, HMMc, Duong-Li} for divergence form operators $T = -\div_x d \nabla_x$. For (bi)sectorial operators satisfying our standard assumptions, the same kind of decomposition has been used in many references including \cite{AA, AusSta} but a proof seems to be missing in the literature. We take the opportunity to close this gap. The construction closely follows \cite{HMMc} but heat semigroup bounds have to be replaced with more technical resolvent bounds.

Throughout this section $T$ is again a (bi)sectorial operator that satisfies the standard assumptions of Section~\ref{subsec: Hardy abstract} or Section~\ref{subsec: Hardy abstract sectorial} and we define $\IH_T^p$ by the abstract theory for first or second-order scaling, respectively. 

\begin{defn}
\label{def: molecule}
Let $p \in (0,1]$, $\eps>0$ and $M \in \IN$. A function $m \in \L^2$ is called \emph{$(\IH^p_T, \eps, M)$-molecule}\index{molecule!$(\IH^p_T, \eps, M)$} if there exists a cube $Q \subseteq \R^n$ and a function $b \in \dom(T^M)$ that satisfies $T^M b = m$ and the following estimates for $j = 1, 2, \ldots$ and $k=0,1,\ldots,M$:
\begin{enumerate}
	\item If $T$ is bisectorial with first-order scaling
	\begin{align*}
	\|(\ell(Q)T)^{-k}m \|_{\L^2(C_j(Q))} \leq (2^j \ell(Q))^{\frac{n}{2}-\frac{n}{p}} 2^{-j \eps}.
	\end{align*}
	
	\item If $T$ is sectorial with second-order scaling
	\begin{align*}
	\|(\ell(Q)^2T)^{-k}m \|_{\L^2(C_j(Q))} \leq (2^j \ell(Q))^{\frac{n}{2}-\frac{n}{p}} 2^{-j \eps}.
	\end{align*}
\end{enumerate}
\end{defn}

\begin{rem}
\label{rem: molecule global bound}
Summing up the bounds in $j$ gives the global $\L^2$-bound $\|((\ell(Q)^\varrho T)^{-k}m \|_2 \leq c \ell(Q)^{n/2-n/p}$, where $\varrho \in \{1,2\}$ is the order of scaling and $c = (2^\eps-1)^{-1}$. If $\eps > n/2$, then we can use H\"older's inequality before summing and obtain $\|(\ell(Q)^\varrho T)^{-k}m \|_1 \leq c \ell(Q)^{n-n/p}$.
\end{rem}

\begin{defn}
\label{def: molecular decomposition}
Let $p \in (0,1]$, $\eps>0$ and $M \in \IN$. A  \emph{molecular $(\IH^p_T, \eps, M)$-representation} of $f \in \cl{\ran(T)}$ is a series $\sum_{i=0}^\infty \lambda_i m_i$ that converges towards $f$ unconditionally in $\L^2$ such that $(\lambda_i) \in \ell^p$ and each $m_i$ is a $(\IH^p_T, \eps, M)$-molecule. The  \emph{molecular Hardy space}
\begin{align*}
\IH^p_{T, \mathrm{mol}, \eps, M} \coloneqq \Big \{ f \in \cl{\ran(T)} : \text{$f$ has a molecular $(\IH^p_T, \eps, M)$-representation}\Big \}
\end{align*}
is equipped with the quasi norm
\begin{align*}
\|f\|_{\IH^p_{T, \mathrm{mol}, \eps, M}} \coloneqq \inf \| (\lambda_i)\|_{\ell^p},
\end{align*}
where the infimum is taken over all admissible representations.
\end{defn}

With these definitions at hand, we establish the following

\begin{thm}
\label{thm: Molecular decomposition adapted Hardy}
Let $p \in (0,1]$, $\eps>0$ and $M \in \IN$ with $M> \nicefrac{n}{p}-\nicefrac{n}{2}$ if $T$ is bisectorial with first-order scaling or $M> \nicefrac{n}{(2p)}-\nicefrac{n}{4}$ if $T$ is sectorial with second-order scaling. Then
\begin{align*}
\IH^p_{T, \mathrm{mol}, \eps, M} = \IH^p_T
\end{align*}
with equivalent quasinorms and the equivalence constants depend on $T$ only through the bounds that are quantified in the standard assumptions.\index{molecular decomposition!for $\IH_T^p$}
\end{thm}

As in many earlier references, the proof relies on the atomic decomposition for tent spaces that we recall beforehand.

\begin{defn}
\label{def: Tp atom}
Let $p \in (0,1]$. A \emph{$\T^p$-atom}\index{atom!for $\T^p$} associated with a cube $Q \subseteq \R^n$ is a measurable function $A: \reu \to V$ with support in $Q \times (0, \ell(Q))$ such that
\begin{align*}
\bigg(\int_0^{\ell(Q)} \int_Q |A(s,y)|^2 \, \frac{\d s \d y}{s}\bigg)^{\frac{1}{2}} \leq \ell(Q)^{\frac{n}{2}-\frac{n}{p}}.
\end{align*}
\end{defn}

\begin{prop}[{\cite[Prop.~5]{CMS}}]
\label{prop: atomic decomposition tent spaces}
Let $p \in (0,1]$. There is a constant $C$ such that every $F \in \T^p$ can be written as $f = \sum_{i=0}^\infty \lambda_i A_i$ with unconditional convergence in $\Lloc^2(\reu)$, where each $A_i$ is a $\T^p$-atom and $\|(\lambda_i)\|_{\ell^p} \leq C \|F\|_{\T^p}$.\index{atomic decomposition!for $\T^p$}
\end{prop}

\begin{rem}
\label{rem: atomic decompsition tent spaces L2 convergence}
The unconditional convergence is not stated explicitly but is immediate from the construction, see \cite[(4.5)]{CMS}. Indeed, we have $\lambda_i A_i = F \ind_{\Delta_i}$, where $(\Delta_i)_i$ is a collection of pairwise disjoint subsets of $\reu$. This also implies that for $f \in \T^p \cap \T^2$ the atomic decomposition converges in $\T^2 = \L^2(\reu, \frac{\d t \d x}{t})$.
\end{rem}

The proof of Theorem~\ref{thm: Molecular decomposition adapted Hardy} relies on two lemmata.

\begin{lem}
\label{lem: molecules uniformly bounded}
Let $p \in (0,1]$ and $\eps>0$. Let $M \in \IN$ and $\psi \in \Psi_+^+$ as follows:
\begin{itemize}
	\item $M> \nicefrac{n}{p}-\nicefrac{n}{2}$ and $\psi(z) = z^{2M} (1+ \i z)^{-4M}$ if $T$ is bisectorial with first-order scaling,
	\item $M> \nicefrac{n}{(2p)}-\nicefrac{n}{4}$ and $\psi(z) = z^{2M}(1+z)^{-4M}$ if $T$ is sectorial with second-order scaling.
\end{itemize}
Then there exists a constant $C$ depending on these parameters and the bounds that are quantified in the standard assumptions such that
\begin{align*}
\|\IQ_{\psi,T} m\|_{\T^p} \leq C
\end{align*}
holds for every $(\IH^p_T, \eps, M)$-molecule $m$.
\end{lem}

\begin{proof}
We give the proof for bisectorial $T$ with first-order scaling. Up to consistently changing the scaling, the argument for sectorial operators is identical. Since
\begin{align*}
\psi(z) = (-\i)^{2M}((1+\i z)^{-1} - (1+\i z)^{-2})^{2M}
\end{align*}
we obtain by composition that $(\psi(tT))_{t>0}$ satisfies $\L^2$ off-diagonal estimates of arbitrarily large order. 

Let $m$ be an $(\IH^p_T, \eps, M)$-molecule associated with a cube $Q$ of sidelength $\ell$. We need a uniform $\L^p$-bound for the square function
\begin{align*}
S_{\psi,T} m(x) \coloneqq \bigg(\iint_{|x-y|<t} |\psi(tT)m(y)|^2 \, \frac{\d t \d y}{t^{1+n}} \bigg)^{1/2}.
\end{align*}
Since $\IH_T^2 = \cl{\ran(T)}$, we have that $\|S_{\psi,T} f\|_2 \lesssim \|f\|_2$ for all $f \in \cl{\ran(T)}$. In particular, we obtain from Hölder's inequality and the molecular decay the local bound
\begin{align*}
\|S_{\psi,T}m\|_{\L^p(16Q)} \leq |16 \ell|^{\frac{n}{p}-\frac{n}{2}} \|S_{\psi,T}m\|_{\L^2(16Q)} \leq C.
\end{align*}
It remains to prove that there is $\alpha > 0$ depending only on $\eps,M,p$ such that for all $j \geq 4$ we have a uniform bound
\begin{align}
\label{eq: goal molecules uniformly bounded}
\|S_{\psi,T}m\|_{\L^2(C_j(Q))} \leq C 2^{-j \alpha} (2^j \ell)^{\frac{n}{2}- \frac{n}{p}}.
\end{align}
Indeed, this implies $\|S_{\psi,T}m\|_{\L^p(C_j(Q))} \leq C 2^{-j \alpha}$ as before and the global $\L^p$-bound for $S_{\psi,T}m$ follows by summing up the $p$-th powers of these estimates.

In order to establish \eqref{eq: goal molecules uniformly bounded}, we split the integral in $t$ at height $2^{\theta (j-1)} \ell$, where $\theta \in (0,1)$ will be fixed later:
\begin{align*}
\|S_{\psi,T}m\|_{\L^2(C_j(Q))}
&= \bigg(\int_{C_j(Q)} \iint_{|x-y|<t}|\psi(tT)m(y)|^2 \, \frac{\d t \d y}{t^{1+n}} \, \d x\bigg)^{1/2}\\
&\lesssim \bigg(\int_0^{2^{\theta (j-1)} \ell} \int_{D_j(Q)}  |\psi(tT)m(y)|^2 \, \frac{\d y \d t}{t} \bigg)^{1/2}\\
&\quad + \bigg(\int_{2^{\theta (j-1)}\ell}^\infty \int_{\R^n} |\psi(tT)m(y)|^2 \, \frac{\d y \d t}{t} \bigg)^{1/2}\\
&\eqqcolon \I + \II,
\end{align*}
where $D_j(Q) \coloneqq 2^{j+2}Q \setminus 2^{j-1}Q$ and we have used Tonelli's theorem to bound the integrals in $x$.
By the molecular properties, we can write $m = T^M b$. Since $\psi \in \Psi_{2M}^{2M}$, we have a uniform $\L^2$-bound for $(tT)^M \psi(tT)$, which together with Remark~\ref{rem: molecule global bound} leads us to
\begin{align*}
\II 
&=  \bigg(\int_{2^{\theta (j-1)}\ell}^\infty \int_{\R^n}  |(tT)^M\psi(tT)b(y)|^2 \, \frac{\d y \d t}{t^{2M+1}} \bigg)^{1/2} \\
&\lesssim (2^{\theta j} \ell)^{-M} \|b\|_2\\
&\leq C 2^{-j(\theta M+\frac{n}{2}-\frac{n}{p})} (2^j \ell)^{\frac{n}{2}-\frac{n}{p}}
\end{align*}
and we can achieve $\alpha \coloneqq \theta M+\nicefrac{n}{2}-\nicefrac{n}{p} > 0$ by taking $\theta$ sufficiently close to $1$. This completes the treatment of $\II$. 

As for $\I$, we decompose further $\I = \I_1 + \I_2$, where $\I_k$ corresponds to replacing $m$ with $m_k$ defined as
\begin{align*}
m_1 \coloneqq \ind_{2^{j+3}Q\setminus 2^{j-2}Q}m, \quad m_2 \coloneqq \ind_{{}^c (2^{j+3}Q\setminus 2^{j-2})}m.
\end{align*}
The $\L^2$-bound for $S_{\psi,T}$ and the molecular estimates yield
\begin{align*}
\I_1 \lesssim \|m_1\|_2 \leq \sum_{k=j-2}^{j+2} \|m\|_{\L^2(C_k(Q))} \leq C 2^{-j\eps} (2^j \ell)^{\frac{n}{2}-\frac{n}{p}}.
\end{align*}
Since the support of $m_2$ is at distance at least $2^{j-2} \ell$ from $D_j(Q)$, we can infer from the off-diagonal decay for $\psi(tT)$ that
\begin{align*}
\I_2 
&\lesssim \|m\|_2 \bigg(\int_0^{2^{\theta (j-1)} \ell} \Big(1+ \frac{2^{j-2} \ell}{t} \Big)^{-2\gamma} \, \frac{\d t}{t} \bigg)^{1/2} \\
&\lesssim \|m\|_2 (2^j \ell)^{-\gamma} \bigg(\int_0^{2^{\theta (j-1)} \ell} t^{2 \gamma} \, \frac{\d t}{t} \bigg)^{1/2} \\
& \leq C(2^j\ell)^{\frac{n}{2}-\frac{n}{p}}  2^{-j((1-\theta)\gamma + \frac{n}{2}-\frac{n}{p})},
\end{align*}
where we have used again Remark~\ref{rem: molecule global bound} in the final step and $\gamma$ is still at our disposal. We have already fixed $\theta \in (0,1)$ and it suffices to take $\gamma$ large enough so that $\alpha \coloneqq (1-\theta)\gamma +\nicefrac{n}{2}-\nicefrac{n}{p} > 0$. This completes the treatment of $\I$ and hence we have established our goal \eqref{eq: goal molecules uniformly bounded}.
\end{proof}

\begin{lem}
\label{lem: coretraction maps atoms to molecules}
Let $p \in (0,1]$. Let $\eps>0$ and $M \in \IN$. Let $\psi(z) = z^{2M} (1+ \i z)^{-4M}$ if $T$ is bisectorial with first-order scaling and $\psi(z) = z^{2M} (1+ z)^{-4M}$ is $T$ is sectorial with second-order scaling. There exists a constant $c$ depending on these parameters and the bounds that are quantified in the standard assumptions, such that $c^{-1} \IC_{\psi,T} A$ is an $(\IH^p_T, \eps, M)$-molecule, whenever $A$ is a $\T^p$-atom.
\end{lem}

\begin{proof}
Again we only do the proof in the bisectorial case and the sectorial case follows line by line up to the usual modifications. 

Let $A$ be a $\T^p$-atom associated with a cube $Q$ of sidelength $\ell$ and set 
\begin{align*}
m \coloneqq \IC_{\psi,T} A = \int_0^\ell (tT)^{2M} (1+ \i tT)^{-4M} A(t) \, \frac{\d t}{t},
\end{align*}
where we have used the support property of $A$. The integral converges weakly in $\L^2$ but as $M \geq 1$, the integral
\begin{align*}
b \coloneqq \int_0^\ell t^{M} (tT)^M (1+ \i tT)^{-4M} A(t) \, \frac{\d t}{t}
\end{align*}
converges strongly and we have $T^M b = m$. We establish the molecular bounds for $m$ up to a generic renorming factor~$c$.

In preparation of the argument, let $g \in \L^2$. For $k=0,\ldots,M$ we bound the $\L^2$ inner product
\begin{align*}
|\langle (\ell T)^{-k} m, g \rangle|
& \leq \ell^{-k} \int_0^\ell |\langle t^k (tT)^{2M-k} (1+ \i tT)^{-4M} A(t), g \rangle| \, \frac{\d t}{t} \\
&= \ell^{-k} \int_0^\ell t^k |\langle A(t), \varphi(tT^*) g \rangle| \, \frac{\d t}{t} \\
&\leq \bigg(\int_0^\ell \|A(t)\|_{\L^2(Q)}^2 \, \frac{\d t}{t}\bigg)^{1/2}  \bigg(\int_0^\ell \|\varphi(tT^*)g\|_{\L^2(Q)}^2 \, \frac{\d t}{t}\bigg)^{1/2} \\
&\leq \ell^{\frac{n}{2}-\frac{n}{p}} \bigg(\int_0^\ell \|\varphi(tT^*)g\|_{\L^2(Q)}^2 \, \frac{\d t}{t}\bigg)^{1/2},
\end{align*}
where $\varphi \in \Psi_M^M$ is given by $\varphi(z) \coloneqq z^{2M-k}(1-\i z)^{-4M}$ and we have used the support and the molecular bound of $A$. Taking the supremum over all $g$ with support in $4Q$ normalized to $\|g\|_2 =1$ and controlling the square function via McIntosh's theorem, we obtain
\begin{align*}
\|(\ell T)^{-k} m\|_{\L^2(4Q)} \leq c \ell^{\frac{n}{2}-\frac{n}{p}},
\end{align*}
which is the required molecular bound for $j=1$. The family $(\varphi(tT))_{t>0}$ satisfies $\L^2$ off-diagonal estimates of arbitrarily large order by composition since we can expand
\begin{align*}
\varphi(z) = \i^{2M-k} (1-(1-\i z)^{-1})^{2M-k}(1-\i z)^{-2M-k}.
\end{align*}
For $j \geq 2$ we take the supremum over all normalized $g$ in $\L^2$ with support in $\L^2(C_j(Q))$ and obtain
\begin{align*}
\|(\ell T)^{-k} m\|_{\L^2(C_j(Q))}
& \lesssim \ell^{\frac{n}{2}-\frac{n}{p}} \bigg(\int_0^\ell \Big(\frac{2^{j-1} \ell}{t} \Big)^{-2\gamma} \, \frac{\d t}{t}\bigg)^{1/2}\\
&\leq c \ell^{\frac{n}{2}-\frac{n}{p}} 2^{-(j-1) \gamma},
\end{align*}
with $\gamma > 0$ at our disposal. We take $\gamma > \nicefrac{n}{p}-\nicefrac{n}{2} + \eps$ to obtain the required molecular decay.
\end{proof}

Putting it all together, we give the

\begin{proof}[Proof of Theorem~\ref{thm: Molecular decomposition adapted Hardy}]
Let $f \in \IH_{T, \mathrm{mol}, \eps, M}^p$ and let $f = \sum_{i=0}^\infty \lambda_i m_i$ be an $\L^2$ convergent molecular representation. We define $\IH_T^p$ via the admissible auxiliary function $\psi$ from Lemma~\ref{lem: molecules uniformly bounded}. Let $\varrho \in \{1,2\}$ be the scaling order. We have
\begin{align}
\label{eq: Hardy norm of molecular representation}
\begin{split}
\|\IQ_{\psi,T} f\|_{\T^p}^p
&=\int_{\R^n} \bigg(\iint_{|x-y| < t} |\psi(t^\varrho T) f(y)|^2 \, \frac{\d t \d y}{t^{1+n}}\bigg)^{p/2} \, \d x\\
&\leq \int_{\R^n} \bigg(\sum_{i=0}^\infty |\lambda_i| \bigg(\iint_{|x-y| < t} |\psi(t^\varrho T)  m_i(y)|^2 \, \frac{\d t \d y}{t^{1+n}} \bigg)^{1/2} \bigg)^p \, \d x\\
&\leq \sum_{i=0}^\infty |\lambda_i|^p  \|\IQ_{\psi,T} m_i\|_{\T^p}^p\\
&\leq C^p \sum_{i=0}^\infty |\lambda_i|^p,
\end{split}
\end{align}
where the first step uses $\L^2$-convergence, the second step is due to $p \leq 1$ and monotone convergence and the third step is by Lemma~\ref{lem: molecules uniformly bounded}. Taking the infimum over all representations yields $\|f \|_{\IH_T^p} \leq C \|f\|_{\IH_{T, \mathrm{mol}, \eps, M}}$.

Conversely, let $f \in \IH_T^p$ and let $\psi$ be the auxiliary function from Lemma~\ref{lem: coretraction maps atoms to molecules}. According to Proposition~\ref{prop: Hardy via contraction}, we can write $f = \IC_{\psi, T}F$ with $F \in \T^p \cap \T^2$ and $\|F\|_{\T^p} \leq 2 \|f\|_{\IH_T^p}$. According to Proposition~\ref{prop: atomic decomposition tent spaces} and the subsequent remark, we can write $F = \sum_{i=0}^\infty \lambda_i A_i$, where the sum converges unconditionally in $\T^2$, each $A_i$ is a $\T^p$-atom and we have $\|(\lambda)_i\|_{\ell^p} \leq C\|F\|_{\T^p}$. Since $\IC_{\psi,T}: \T^2 \to \L^2$ is bounded, we get an unconditionally $\L^2$-convergent representation
\begin{align*}
f = \IC_{\psi,T} F = \sum_{i=0}^\infty \lambda_i \IC_{\psi,T} A_i =  \sum_{i=0}^\infty (c\lambda_i) m_i,
\end{align*}
where $c$ is the constant from Lemma~\ref{lem: coretraction maps atoms to molecules} and the $m_i \coloneqq c^{-1} \IC_{\psi,T} A_i$ are $(\IH_T^p, \eps, M)$-molecules. This proves $\|f\|_{\IH_{T, \mathrm{mol}, \eps, M}} \leq 2Cc \|f\|_{\IH_T^p}$.
\end{proof}	
\subsection{Connection with the non-tangential maximal function}
\label{subsec: Hardy and NT}

We recall the non-tangential maximal function
\begin{align*}
\NT F(x) = \sup_{t>0}  \bigg(\bariint_{W(t,x)} |F(s,y)|^2 \, \d s \d y\bigg)^{1/2},
\end{align*}
where $W(t,x) = (\nicefrac{t}{2},2t)\times B(x,t)$. At this level of generality we do not know whether $\IH_T^p$ could be characterized via $\NT$ as in~\cite{HMMc, Duong-Li} but, using the molecular decompositions, we can give upper bounds for the non-tangential maximal function of resolvent families and Poisson-type semigroups acting on $\IH_T^p$ if $p \leq 1$. Such result can be extended to $p \leq 2$ by interpolation provided the result for $p = 2$ holds, which might be a concern in itself. 

We begin with a simple comparison of the non-tangential maximal function and the uncentered Hardy--Littlewood maximal operator $\cM$ in $\R^n$.
\begin{lem}
\label{lem: Whitney everage controlled by HL-Max}
Let $\psi: (0,\infty) \to \Lop(\L^2)$ be a strongly measurable family that satisfies $\L^2$ off-diagonal estimates of order $\gamma> n/2$.  Then there is a constant $C$ depending on dimensions and the off-diagonal bounds, such that 
\begin{align*}
\bariint_{W(t,x)} |\psi(s)f(y)|^2 \, \d s \d y \leq C \Max(|f|^2)(x)
\end{align*}
for all $f \in \L^2$ and all $(t,x) \in \reu$.
\end{lem}

\begin{proof}
Set $B \coloneqq B(x,t)$ and split $f = \sum_{j =1}^\infty f_j$, where $f_j \coloneqq \ind_{C_j(B)}f$. For $\nicefrac{t}{2} < s < 2t$ we have by assumption
\begin{align*}
\int_B |\psi(s)f_j(y)|^2 \, \d y
& \lesssim \Big(1+ \frac{(2^j - 1)t}{s} \Big)^{-2\gamma} \|f\|_{\L^2(C_j(B))}^2\\
&\lesssim 2^{-2\gamma j} \int_{2^{j+1} B} |f|^2 \\
&\lesssim t^n 2^{-j(2\gamma-n)} \Max(|f|^2)(x).
\end{align*}
The claim follows by summing in $j$ and averaging in $s$.
\end{proof}

We also recall Kolmogorov's lemma for bounding the maximal operator on $\L^\theta$ for $\theta < 1$, see for instance~\cite[Lem.~5.16]{Duo}.

\begin{lem}[Kolmogorov\index{Kolmogorov's lemma}]
\label{lem: Kolmogorov}
Let $\theta \in (0,1)$ and $E \subseteq \R^n$ a set of finite measure. There is a constant $C = C(\theta,n)$ such that
\begin{align*}
\int_E |\Max f(y)|^\theta \, \d y \leq C |E|^{1-\theta} \|f\|_1^\theta \quad (f \in \L^1).
\end{align*}
\end{lem}

With these tools at hand, we establish a first non-tangential maximal bound on $\IH_T^p$.

\begin{prop}
\label{prop: NT bound for p<1}
Let $p \in (0,1]$ and $\eps > 0$. Let $M \in \IN$ and $\psi \in \H^\infty$ as follows:
\begin{itemize}
	\item $M> \nicefrac{n}{p}-\nicefrac{n}{2}$ and $\psi(z) = (1+ \i z)^{-2M}$ if $T$ is bisectorial with first-order scaling,
	\item $M>  \nicefrac{n}{(2p)}-\nicefrac{n}{4}$ and $\psi(z) =(1+z)^{-2M}$ if $T$ is sectorial with second-order scaling.
\end{itemize}
Then there exists a constant $C$ depending on these parameters and the bounds that are quantified in the standard assumptions such that
\begin{align*}
\|\NT(\IQ_{\psi,T} f)\|_p \leq C \|f\|_{\IH_T^p} \quad (f \in \IH_T^p).
\end{align*}
\end{prop}

\begin{proof}
Let $f \in \IH_T^p$ and $f = \sum_{i=0}^\infty \lambda_i m_i$ be an $\L^2$-convergent molecular representation as in Theorem~\ref{thm: Molecular decomposition adapted Hardy}. Then $\IQ_{\psi, T}f = \sum_{i=0}^\infty \lambda_i \IQ_{\psi, T} m_i$ in $\L^\infty(0,\infty; \L^2)$ and by sublinearity of the maximal function we find
\begin{align*}
\|\NT(\IQ_{\psi, T}f)\|_p^p
&\leq \bigg\|\sum_{i=0}^\infty |\lambda_i| \NT(\IQ_{\psi, T}m_i)\bigg \|_p^p
\leq \sum_{i=0}^\infty |\lambda_i|^p \|\NT(\IQ_{\psi, T}m_i) \|_p^p.
\end{align*}
Consequently, it suffices to treat the case that $f=m$ is an $(\IH_T^p, \eps, M)$-molecule (associated with a cube $Q$ of sidelength $\ell$) and derive a uniform bound. We only write out the argument in the bisectorial case. As usual, the proof is identical in the sectorial case upon changing the scaling. 

\medskip

\noindent \emph{Step 1: Local bound for $\NT$}. By composition, the family $(\psi(tT))_{t>0}$ satisfies $\L^2$ off-diagonal estimates of arbitrarily large order. Therefore, Lemma~\ref{lem: Whitney everage controlled by HL-Max} yields $\NT(\IQ_{\psi,T}m) \leq C (\Max(|m|^2))^{1/2}$ a.e.\ on $\R^n$ and by means of Kolmogorov's lemma and Remark~\ref{rem: molecule global bound} we get 
\begin{align*}
\|\NT(\IQ_{\psi, T}m)\|_{\L^p(16Q)}^p
\lesssim \int_{16Q} |\Max(|m|^2)(y)|^{\frac{p}{2}} \, \d y
\leq |16 Q|^{1-\frac{p}{2}} \|m\|_2^p
\leq C^p.
\end{align*}

\medskip

\noindent \emph{Step 2: Decomposition of $\NT$ on annuli}. It remains to show that there is $\alpha>0$ depending only on $\eps, M, p$ such that for all $j \geq 4$ we have a uniform bound
\begin{align}
\label{eq1: NT bound for p<1}
\|\NT(\IQ_{\psi, T}m)\|_{\L^p(C_j(Q))}^p \leq 2^{-j \alpha} C^p.
\end{align}
The claim then follows by summing up in $j$. To this end, we fix $j \geq 4$ and split
\begin{align*}
\NT(\IQ_{\psi, T}m) \leq \NT^{\mathrm{loc}}(\IQ_{\psi, T}m) + \NT^{\mathrm{glob}}(\IQ_{\psi, T}m),
\end{align*}
where the local and global parts correspond to restricting the size of Whitney boxes in the definition of $\NT$ to $t \leq \ell$ and $t\geq \ell$, respectively.

\medskip

\noindent \emph{Step 3: Bound for $\NT^{\mathrm{loc}}$ on $C_j(Q)$}. Let $0< t < \ell$ and $x \in C_j(Q)$. Splitting $m = \sum_{i = 1}^\infty m_i$, where $m_i \coloneqq \ind_{C_i(Q)}m$, we get
\begin{align*}
\bigg(&\bariint_{W(t,x)} |(1+\i s T)^{-2M}m|^2 \, \d s \d y \bigg)^{1/2} \\
&\lesssim \sum_{|i-j| \geq 2} t^{-\frac{n}{2}}\Big(1+ \frac{\dist(B(x,t), C_i(Q))}{t} \Big)^{-\gamma} \|m_i\|_2 
+\sum_{|i-j| \leq 1} \Max(|m_i|^2)(x)^{\frac{1}{2}},
\end{align*}
where we have used $\L^2$ off-diagonal decay of the resolvents whenever $|i-j| \geq 2$ and Lemma~\ref{lem: Whitney everage controlled by HL-Max} whenever $|i-j| \leq 1$. The order $\gamma > 0$ is at our disposal. For any set $E \subseteq \R^n$ we have
\begin{align*}
1+ \frac{\dist(B(x,t), E)}{t} \geq \frac{1}{2} + \frac{\dist(x, E)}{4t}
\end{align*}
as follows by distinguishing whether or not $t \geq \nicefrac{\dist(x,E)}{2}$. Specializing to $E = C_i(Q)$ with $|i-j| \geq 2$, we get
\begin{align*}
1+ \frac{\dist(B(x,t), C_i(Q))}{t} 
\geq \frac{1}{2} + \frac{\dist(x, C_i(Q))}{4t}
\gtrsim \frac{2^{i \vee j} \ell}{4t}.
\end{align*}
We also have
\begin{align*}
\|m_i\|_2 \leq |2^{(i \vee j)+1} Q|^{\frac{1}{2}} \bigg(\barint_{2^{(i \vee j)+1} Q} |m_i|^2 \, \d y\bigg)^{\frac{1}{2}} \lesssim (2^{i \vee j} \ell)^{\frac{n}{2}} \Max(|m_i|^2)(x)^\frac{1}{2}.
\end{align*}
Applying these bounds on the right-hand side of our estimate leads us to
\begin{align}
\label{eq2: NT bound for p<1}
\begin{split}
\bigg(&\bariint_{W(t,x)} |(1+\i s T)^{-2M}m|^2 \, \d s \d y \bigg)^{1/2} \\
&\lesssim \sum_{i \leq j-2} (2^j \ell)^{\frac{n}{2} - \gamma} t^{\gamma- \frac{n}{2}}  \Max(|m_i|^2)(x)^{\frac{1}{2}} \\
&\quad +\sum_{|i-j| \leq 1} \Max(|m_i|^2)(x)^{\frac{1}{2}} \\
&\quad +\sum_{i \geq j+2} (2^i \ell)^{\frac{n}{2} - \gamma} t^{\gamma- \frac{n}{2}}  \Max(|m_i|^2)(x)^{\frac{1}{2}}.
\end{split}
\end{align}
From now on we require $\gamma > \nicefrac{n}{2}$. On the right-hand side $t$ appears with positive exponent and hence the supremum over $0<t \leq \ell$ is attained for $t=\ell$. We conclude that
\begin{align}
\label{eq3: NT bound for p<1}
\begin{split}
\NT^{\mathrm{loc}}(\IQ_{\psi, T}m)(x)
& \lesssim \sum_{i \leq j-2} 2^{j(\frac{n}{2} - \gamma)}  \Max(|m_i|^2)(x)^{\frac{1}{2}} \\
&\quad +\sum_{|i-j| \leq 1} \Max(|m_i|^2)(x)^{\frac{1}{2}} \\
&\quad + \sum_{i \geq j+2} 2^{i(\frac{n}{2} - \gamma)}  \Max(|m_i|^2)(x)^{\frac{1}{2}}.
\end{split}
\end{align}
Kolmogorov's lemma and the molecular bounds for $m$ imply
\begin{align*}
\int_{C_j(Q)} |\Max(|m_i|^2)(x)|^{\frac{p}{2}} \, \d x
&\leq C|C_j(Q)|^{1-\frac{p}{2}} \|m_i\|_2^p \\
&\leq C 2^{j(\frac{n}{p}-\frac{n}{2})p} 2^{i(\frac{n}{2}-\frac{n}{p} - \eps)p},
\end{align*}
so that integrating the $p$-th power of \eqref{eq3: NT bound for p<1} in $x \in C_j(Q)$ yields
\begin{align*}
\|\NT^{\mathrm{loc}}(\IQ_{\psi, T}m)\|_{\L^p(C_j(Q))}^p
&\lesssim \sum_{i \leq j-2}  2^{j(\frac{n}{p} - \gamma )p}2^{i(\frac{n}{2}-\frac{n}{p} - \eps)p} \\
&\quad + \sum_{|i-j| \leq 1} 2^{j(\frac{n}{p}-\frac{n}{2})p} 2^{i(\frac{n}{2}-\frac{n}{p} - \eps)p} \\
&\quad + \sum_{i \geq j-2}  2^{j(\frac{n}{p}-\frac{n}{2})p} 2^{i(n - \frac{n}{p} - \eps - \gamma)p} \\
&\simeq  2^{j(\frac{n}{2} - \eps - \gamma)p} + 2^{-j \eps p} +2^{j(\frac{n}{p} - \gamma )p}.
\end{align*}
This establishes \eqref{eq1: NT bound for p<1} for $\NT^{\mathrm{loc}}$ provided that eventually we take $\gamma > \nicefrac{n}{p}$ (which implies $\gamma > \nicefrac{n}{2}$).

\medskip

\noindent \emph{Step 4: Bound for $\NT^{\mathrm{glob}}$ on $C_j(Q)$}. We write $m = T^M b$ as in Definition~\ref{def: molecule}. We have
\begin{align*}
(1+ \i t T)^{-2M} m
&= (\i t)^{-M} (i t)^M (1+ \i t T)^{-2M} m \\
&= (\i t)^{-M} ((1+ \i t T)^{-1} - (1+ \i t T)^{-2})^M b \\
&\eqqcolon (\i t)^{-M} \varphi(tT)b,
\end{align*}
where $(\varphi(tT))_{t>0}$ satisfies $\L^2$ off-diagonal estimates of arbitrarily large order. Hence, we can repeat the first part of Step~3 with $\varphi, b$ replacing $\psi, m$ and due to the additional factor $(\i t)^{-M}$ our substitute for \eqref{eq2: NT bound for p<1} becomes
\begin{align*}
\bigg(&\bariint_{W(t,x)} |(1+\i s T)^{-2M}m|^2 \, \d s \d y \bigg)^{1/2} \\
&\lesssim \sum_{i \leq j-2} (2^j \ell)^{\frac{n}{2} - \gamma} t^{\gamma- \frac{n}{2}-M}  \Max(|b_i|^2)(x)^{\frac{1}{2}} \\
&\quad +\sum_{|i-j| \leq 1} t^{-M} \Max(|b_i|^2)(x)^{\frac{1}{2}} \\
&\quad +\sum_{i \geq j+2} (2^i \ell)^{\frac{n}{2} - \gamma} t^{\gamma- \frac{n}{2}-M}  \Max(|b_i|^2)(x)^{\frac{1}{2}}
\end{align*}
with $\gamma > 0$ at our disposal and $b_i \coloneqq \ind_{C_i(Q)} b$. We require $\gamma < \nicefrac{n}{2} + M$. Then $t$ appears with negative exponent on the right-hand side and passing to the supremum for all $t \geq \ell$, we get
\begin{align*}
\NT^{\mathrm{glob}}(\IQ_{\psi, T}m)(x)
& \lesssim \sum_{i \leq j-2} 2^{j(\frac{n}{2} - \gamma)}  \Max(|\ell^{-M}b_i|^2)(x)^{\frac{1}{2}} \\
&\quad +\sum_{|i-j| \leq 1} \Max(|\ell^{-M} b_i|^2)(x)^{\frac{1}{2}} \\
&\quad + \sum_{i \geq j+2} 2^{i(\frac{n}{2} - \gamma)}  \Max(|\ell^{-M} b_i|^2)(x)^{\frac{1}{2}}.
\end{align*}
Now, $\ell^{-M} b = (\ell T)^{-M}m$ satisfies the same $\L^2$-bounds on annuli as $m$ and we can repeat the arguments in Step~3 from \eqref{eq3: NT bound for p<1} onward in order to conclude \eqref{eq1: NT bound for p<1} for $\NT^{\mathrm{glob}}$ provided that at the end of the proof we take again $\gamma > \nicefrac{n}{p}$. This requirement is compatible with $\gamma < \nicefrac{n}{2} + M$ since we have $M>\nicefrac{n}{p}-\nicefrac{n}{2}$ by assumption.
\end{proof}

In the context of boundary value problems it will be important to have a statement as above with a Poisson-like semigroup replacing the resolvents. To this end we need the following fact.

\begin{lem}
\label{lem: N < S}
Let $p \in (0,\infty)$. There is a constant $C = C(n,p)$ such that
\begin{align*}
\|\NT(F)\|_p \leq C \|F\|_{\T^{p}} \quad (F \in \Lloc^2(\reu)).
\end{align*}
\end{lem}

We add a proof for convenience.

\begin{proof}
Let $(t,x) \in \reu$. Since $(s,y) \in W(t,x)$ implies $|x-y|<t \leq 2s$ and $t \geq \nicefrac{s}{2}$, we have that
\begin{align*}
\bigg(\bariint_{W(t,x)} |F(s,y)|^2 \, \d s \d y\bigg)^{1/2} \leq C \bigg(\iint_{|x-y| < 2s} |F(s,y)|^2 \, \frac{\d s \d y}{s^{1+n}}\bigg)^{1/2}.
\end{align*}
The right-hand side does not depend on $t$ and its $\L^p$-quasinorm in $x$ is equivalent to $\|F\|_{\T^{p}}$ by a change of aperture. The claim follows by taking the supremum in $t$ and integrating the $p$-th powers in $x$. 
\end{proof}

\begin{prop}
\label{prop: NT bound semigroup for p<1}
Let $p \in (0,1]$. Let $\psi(z) = \e^{-\sqrt{z^2}}$ if $T$ is bisectorial with first-order scaling and $\psi(z) = \e^{-\sqrt{z}}$ if $T$ is sectorial with second-order scaling. Then there exists a constant $C$ depending on $p$ and the bounds that are quantified in the standard assumptions, such that
\begin{align*}
\|\NT(\IQ_{\psi,T} f)\|_p \leq C \|f\|_{\IH_T^p} \quad (f \in \IH_T^p).
\end{align*}
Moreover, the bound continues to hold for $p \in (1,2]$ by interpolation if it holds for $p=2$.
\end{prop}

\begin{proof}
First, let $p \in (0,1]$ and define an auxiliary function $\varphi$ as follows:
\begin{itemize}
	\item If $T$ is bisectorial with first-order scaling, let $M > \nicefrac{n}{p}-\nicefrac{n}{2}$ and $\varphi(z) \coloneqq \psi(z) - (1+\i z)^{-2M}$. Then $\varphi \in \Psi_1^{2M}$, so that the technical condition in Proposition~\ref{prop: auxiliary function bisectorial} holds.
	\item If $T$ is sectorial with second-order scaling, let $M > \nicefrac{n}{(2p)}-\nicefrac{n}{4}$ and $\varphi(z) \coloneqq \psi(z) - (1+ z)^{-2M}$. Then $\varphi \in \Psi_{1/2}^{2M}$ and the corresponding technical condition for sectorial operators (Section~\ref{subsec: Hardy abstract sectorial}) holds.
\end{itemize}
We find for all $f \in \IH_T^p$ that
\begin{align*}
\|\NT(\IQ_{\psi,T} f)\|_p 
&\leq \|\NT(\IQ_{\varphi,T} f)\|_p + \|\NT(\IQ_{\psi-\varphi,T} f)\|_p \\
&\lesssim \|\IQ_{\varphi,T} f\|_{\T^{p}} + \|\NT(\IQ_{\psi-\varphi,T} f)\|_p \\
&\lesssim \|f\|_{\IH_T^p},
\end{align*}
where the second step is due Lemma~\ref{lem: N < S} and the third step uses the definition of the $\IH_T^p$-norm and Proposition~\ref{prop: NT bound for p<1}. 

Suppose in addition that this bound holds for $p=2$. Let $\phi \in \Psi_{\infty}^\infty$ and recall the definition of $\IH_T^p$ via the contraction mapping $\IC_\phi$ (Proposition~\ref{prop: Hardy via contraction}). The claim is then equivalent to $F \mapsto \NT(\IQ_{\psi,T} \IC_\phi F)$ being bounded $\T^{p} \cap \T^2 \to \L^p$ for the respective $p$-norms. By assumption this holds for $p=2$ and from the first part of the proof it follows for $p=1$, so the claim follows by complex interpolation for positive sublinear operators~\cite{Janson-sublinear}. 
\end{proof}
\subsection{\texorpdfstring{$\boldsymbol{D}$}{D}-adapted spaces}
\label{subsec: D-adapted spaces}

The unperturbed Dirac operator $D$ satisfies the standard assumptions of Section~\ref{subsec: Hardy abstract}. In order to fully understand the associated Hardy--Sobolev and Besov spaces, we need the orthogonal projection\index{P@$\IP_D$ (orthogonal projection onto $\cl{\ran(D)}$)} $\IP_D: \L^2 \to \cl{\ran(D)} \subseteq \L^2$. From the specific form of $D^2$ in \eqref{eq: L and M} we see that
\begin{align*}
D^2 f = -\Delta_x f \quad (f \in \dom(D^2) \cap \ran(D))
\end{align*}
and hence that $\IP_D = -\Delta_x^{-1} D^2$ holds on the dense subspace $\dom(D^2)$ of $\L^2$. Now, $-\Delta_x^{-1}D^2$ can also be viewed as a Fourier multiplier with symbol
\begin{align}
\label{eq: symbol of PD}
\begin{bmatrix} 1_{\IC^m} & 0 \\ 0 & (|\xi|^{-2} \xi \otimes \xi) \otimes 1_{\IC^m} \end{bmatrix},
\end{align}
where $\xi \in \R^n$ is the Fourier variable and we think of $\IC^{mn} \simeq (\IC^m)^n$ as $n$-vectors of elements in $\IC^m$ just as in the definition of vector-valued gradient and divergence. This symbol is homogeneous of degree zero and smooth outside of $0$ and hence falls in the scope of the Mihlin multiplier theorem~\cite[Thm.~5.2.2]{Triebel-TheoryOf}\index{Theorem!Mihlin multiplier}. Therefore $-\Delta_x^{-1}D^2$ extends boundedly to $\Xdot^{s,p}$, where $\X \in \{\B, \H\}$, for all $s \in \R$ and $p \in (0,\infty)$. The extension to $\L^2$ is precisely $\IP_D$ and we keep on denoting the extensions to other spaces by the same symbol. From \eqref{eq: symbol of PD} we also obtain the block structure
\begin{align}
\label{eq: matrix for PD}
\IP_D \eqqcolon \begin{bmatrix} 1 & 0 \\ 0 & \IP_{\curl_x} \end{bmatrix}.
\end{align}
Since $\cl{\ran(D)}$ coincides with the space $\cH$ in the ellipticity condition \eqref{eq: accretivity A}, we get that $\IP_{\curl_x}$ is the projection onto the curl-free $\L^2$ vector fields. By \cite[Thm.~5.3]{AA} we have for $s \in \R$ and $p \in (0,\infty)$ that
\begin{align*}
\IX^{s,p}_D = \IP_D(\Xdot^{s,p} \cap \L^2)
\end{align*}
with equivalence of $p$-quasinorms. In particular, $\IP_D(\Xdot^{s,p})$ equipped with the norm of $\Xdot^{s,p}$ is a completion of $\IX^{s,p}_D$ in $\cZ'$. Let now $\psi \in \H^\infty$ for the sectorial functional calculus and put $\varphi(z) \coloneqq \psi(z^2)$. Then \eqref{eq: calculus for L and M from BD} with $B=1$ yields for all $t>0$ that
\begin{align*}
\varphi(tD) = \begin{bmatrix} \psi(-t^2 \Delta_x) & 0 \\ 0 & \psi(-t^2 \nabla_x \div_x) \end{bmatrix},
\end{align*}
that is to say
\begin{align}
\label{eq: splitting of extension operators}
\IQ_{\varphi,D} = \begin{bmatrix} \IQ_{\psi, -\Delta_x} & 0 \\ 0 & \IQ_{\psi, -\nabla_x \div_x} \end{bmatrix}.
\end{align}
On taking $\psi$ with sufficient decay at $0$ and $\infty$, we conclude $\IX^{s,p}_D = \IX^{s,p}_{-\Delta_x} \oplus \IX^{s,p}_{-\nabla_x \div_x}$. Along with \eqref{eq: matrix for PD} we can characterize the $D$-adapted spaces as in Figure~\ref{fig: identification D-adapted}.

\begin{figure}
\centering
\begin{tikzcd}[row sep=0pt, column sep=0pt]
	\IX^{s,p}_{-\Delta_x} \arrow[rrrrrr, dash, bend right, "="] & \oplus & \IX^{s,p}_{-\nabla_x \div_x}  \arrow[rrrrrr, dash, bend left, "="] & = &\IX^{s,p}_D& = &\Xdot^{s,p} \cap \L^2 & \oplus & \IP_{\curl_x} (\Xdot^{s,p} \cap \L^2).
\end{tikzcd}
\caption{Identification of Hardy--Sobolev and Besov spaces up to equivalent quasinorms in the unperturbed case $B=1$.\index{Hardy--Sobolev space!for $D$}}
\label{fig: identification D-adapted}
\end{figure}

As a matter of fact, Theorem~\ref{thm: Molecular decomposition adapted Hardy} for $D$ comprises a molecular decomposition for $\Hdot^{0,p} \cap \L^2 = \H^p \cap \L^2$ when $p \in (0,1]$. In order to illustrate how operator-adapted and standard theory interact for a specific differential operator, we recover an atomic decomposition for $\H^p$ from the molecular decomposition of $\IH_D^p$. 

\begin{prop}
\label{prop: Hp atomic decomposition}
Let $p \in (1_*,1]$. Every $f \in \IH_D^p$ can be written as $f = \sum_{i=0}^\infty \sum_{j=1}^\infty \lambda_i^j a_i^j$ with convergence in $\L^2$, where each $a_i^j$ is an $\L^2$-atom for $\H^p$. Moreover, $\|f\|_{\H^p} \simeq \inf \|(\lambda_i^j)_{i,j}\|_{\ell^p}$, where the infimum is taken over all such representations. \index{atomic decomposition!for $\IH_D^p$}
\end{prop}

\begin{proof}
Let $C$ be such that $\|a\|_{\H^p} \leq C$ for all $\L^2$-atoms $a$ for $\H^p$. We get for any $\L^2$-convergent atomic representation $f = \sum_{i=0}^\infty \sum_{j=1}^\infty \lambda_i^j a_i^j$ that
\begin{align*}
 \|f\|_{\IH_{D}^p}^p \simeq \|f\|_{\H^p}^p \leq C^p \|(\lambda_{i}^j)_{i,j}\|_{\ell^p}^p.
\end{align*}

Conversely, let $f \in \IH_D^p$. Due to Theorem~\ref{thm: Molecular decomposition adapted Hardy} we have $f = \sum_{i=0}^\infty \lambda_i m_i$, where each $m_i$ is an $(\IH_D^p, 1, 1)$-molecule and $\|(\lambda_i)_i\|_{\ell^p} \lesssim \|f\|_{\H^p}$. For the moment fix $m_i$. 

Let $Q_i$ be the associated cube and write $m_i = Db_i$ as in Definition~\ref{def: molecule}. Let $(\chi_i^j)_{j=1}^\infty$ be a smooth partition of unity on $\R^n$ such that 
\begin{align}
\label{eq1: Hp atomic decomposition}
0 \leq \chi_i^j \leq \ind_{C_j(Q_i) \cup C_{j-1}(Q_i)}, \quad \|\nabla_x \chi_i^j\|_\infty \leq c(n) (2^{j} \ell(Q_i))^{-1},
\end{align}
where we set $C_0(Q_i) \coloneqq \emptyset$. Then $b_i = \sum_{j=1}^\infty \chi_i^j b_i$ unconditionally in $\L^2$. Since $D$ is a first-order differential operator, each $D(\chi_i^j b_i)$ is supported in $2^{j+1}Q_i$, has mean value zero and satisfies
\begin{align*}
\|D(\chi_i^j b_i)\|_2 
&\leq \|m_i\|_{\L^2(C_j(Q_i))} + 2^{-j} c\|\ell(Q_i)^{-1} b_i\|_{\L^2(C_j(Q_i))} \\
&\leq c (2^{j+1} \ell(Q_i))^{\frac{n}{2}-\frac{n}{p}} 2^{-j}, 
\end{align*}
where $c$ only depends on dimensions. This means that $a_i^j \coloneqq c^{-1} 2^j D(\chi_i^j b_i)$ is an $\H^p$-atom. Since $D$ is closed, we obtain
\begin{align*}
m_i = Db_i = \sum_{j=1}^\infty D(\chi_i^j b_i) = \sum_{j=1}^\infty c2^{-j} a_i^j,
\end{align*}
and
\begin{align*}
	f = \sum_{i=0}^\infty \sum_{j=1}^\infty c 2^{-j} \lambda_i a_i^j
\end{align*}
is the desired atomic decomposition.
\end{proof}

The proof above showed more.

\begin{cor}
\label{cor: HDp molecules}
Let $p \in (1_*,1]$. There is a constant $C$ that depends on dimensions and $p$ such that every $(\IH_D^p,1,1)$-molecule $m$ satisfies $\|m\|_{\H^p} \leq C$ and hence also $\|m\|_{\IH_D^p} \leq C$.
\end{cor}

Note that Lemma~\ref{lem: molecules uniformly bounded} gives the same result provided that $\nicefrac{n}{p} < 1 + \nicefrac{n}{2}$. We have used the specific structure of $D$ to get the conclusion without this restriction.

We shall also need an atomic decomposition of $\Hdot^{1,p} \cap \W^{1,2}$ as in~\cite{Jawerth-Frazier}, but with $\Wdot^{1,2}$-convergence rather than convergence in $\Hdot^{1,p}$. While this can certainly be inferred from inspection of the proof in \cite{Jawerth-Frazier}, we prefer to give a direct and more transparent argument that relies on the lifting property from Corollary~\ref{cor: Lifting pre-Hardy}.

\begin{defn}
\label{def: H1p atom}
Let $p \in (1_*,1]$. An \emph{$\L^2$-atom for $\Hdot^{1,p}$}\index{atom!for $\Hdot^{1,p}$} is a function $a$ supported in a cube $Q \subseteq \R^n$ such that $\|\nabla_x a\|_2 \leq \ell(Q)^{\frac{n}{2}-\frac{n}{p}}$.
\end{defn}

\begin{prop}
\label{prop: H1p atomic decomposition}
Let $p \in (1_*,1]$. Every $f \in \Hdot^{1,p} \cap \W^{1,2}$ can be written as $f = \sum_{i=0}^\infty \sum_{j=1}^\infty \lambda_i^j a_i^j$ with convergence in $\Wdot^{1,2}$, where each $a_i^j$ is an $\L^2$-atom for $\Hdot^{1,p}$. Moreover, $\|f\|_{\Hdot^{1,p}} \simeq \inf \|(\lambda_i^j)_{i,j}\|_{\ell^p}$, where the infimum is taken over all such representations.
\end{prop}

\begin{proof}
If $f$ is an $\L^2$-atom for $\Hdot^{1,p}$, then $\nabla_x f$ is an $\L^2$-atom for $\H^p$. Hence, if $f = \sum_{i=0}^\infty \sum_{j=1}^\infty\lambda_i^j a_i^j$ is an $\Hdot^{1,p}$ atomic decomposition as above, then 
\begin{align*}
\nabla_x f
= \sum_{i=0}^\infty \sum_{j=1}^\infty \lambda_i^j \nabla_x a_i^j
\end{align*}
is an $\H^p$ atomic decomposition and $\|f\|_{\Hdot^{1,p}} = \|\nabla_x f\|_{\H^p} \lesssim \|(\lambda_i^j)_{i,j}\|_{\ell^p}$ follows. Conversely, let $f \in \Hdot^{1,p} \cap \W^{1,2}$. Then $[f,0]^\top \in \IH^{1,p}_D \cap \dom(D)$, see Figure~\ref{fig: identification D-adapted}, and so $D([f,0]^\top) \in \IH_D^{p}$ by Corollary~\ref{cor: Lifting pre-Hardy}. The atomic decomposition obtained in the proof of Proposition~\ref{prop: Hp atomic decomposition} takes the form 
\begin{align*}
D \begin{bmatrix} f \\ 0 \end{bmatrix} 
=\begin{bmatrix} 0 \\ -\nabla_x f \end{bmatrix} 
= \sum_{i=0}^\infty \sum_{j=1}^\infty c \lambda_i 2^{-j} a_i^j, \quad a_i^j 
= c^{-1} 2^j D(\chi_i^j b_i),
\end{align*}
where each $a_i^j$ is an $\L^2$-atom for $\H^p$ and the $\chi_i^j$ are smooth functions satisfying \eqref{eq1: Hp atomic decomposition}. The function $(c^{-1} 2^j \chi_i^j b_i)_\no$ has support in $2^{j+1} Q_i$ and satisfies $-\nabla_x (c^{-1} 2^j \chi_i^j b_i)_\no = (a_i^j)_\ta$. Hence, it is an $\L^2$-atom for $\Hdot^{1,p}$ and the decomposition we are looking for is 
\begin{align*}
f = \sum_{i=0}^\infty \sum_{j=1}^\infty c \lambda_i 2^{-j} (c^{-1} 2^j \chi_i^j b_i)_\no. &\qedhere
\end{align*}
\end{proof}
\subsection{Spaces adapted to perturbed Dirac operators}
\label{subsec: DB-adapted spaces}

Now, we apply the abstract theory with first-order scaling to the bisectorial operators $BD$ and $DB$ and relate the operator-adapted spaces to those obtained for the sectorial operators $L,M,\tL,\tM$ with second-order scaling. Thanks to the different orders of scaling, the meaning of $s$ as a smoothness parameter is the same for all adapted spaces. 

In analogy with \eqref{eq: splitting of extension operators} we have that whenever $\psi$ is an admissible auxiliary function on a sector for the definition of $\IX_{L}^{s,p}$ and $\IX_{M}^{s,p}$, then $\varphi(z) \coloneqq \psi(z^2)$ is admissible for $\IX_{BD}^{s,p}$ and
\begin{align}
\label{eq: Q-operator splitting L and M}
\IQ_{\varphi,BD} = \begin{bmatrix} \IQ_{\psi, L} & 0 \\ 0 & \IQ_{\psi, M} \end{bmatrix}.
\end{align}
This is again a consequence of \eqref{eq: calculus for L and M from BD}. The same kind of relation holds with $DB$ on the left and $\tL$, $\tM$ on the right and follows from \eqref{eq: tL and tM}. Merely by definition we obtain
\begin{align}
\label{eq: Hardy splitting L and M}
\begin{split}
\IX_{BD}^{s,p} &= \IX_L^{s,p} \oplus \IX_M^{s,p}, \\
\IX_{DB}^{s,p} &= \IX_{\tL}^{s,p} \oplus \IX_{\tM}^{s,p}.
\end{split}
\end{align}
In this sense the theory for the perturbed Dirac operators encompasses the theory of all four second-order operators. Figure~\ref{fig: diagram with X} summarizes their various relations.

\begin{figure}[ht]
\centering
\begin{tikzcd}[column sep=0pt, row sep=40pt]
	\IX^{s,p}_{BD} \cap \ran([BD])  &=&	\IX^{s,p}_{L} \cap \ran({L^{1/2}})  & \oplus & \IX^{s,p}_{M} \cap \ran(M^{1/2})\\
	\IX^{s+1,p}_{BD} \cap \dom(D) \arrow[u, "{[BD]}"] \arrow[d, "D"] &=&	\IX^{s+1,p}_{L} \cap \dom({L^{1/2}}) \arrow[u,"{L^{1/2}}"]\arrow[drr, "\! \! \! \!{-\nabla_x}" near start] & \oplus &\arrow[dll, "\! \! \! \! {\div_x}" near end] \arrow[u, "M^{1/2}"] \IX^{s+1,p}_{M} \cap \dom(M^{1/2}) \\
	\IX^{s,p}_{DB} \cap \ran(D) \arrow[d, "B"] &=&	\IX^{s,p}_{\tL} \cap \ran(\div_x) \arrow[d, "{a^{-1}}"]  &\oplus &  \IX^{s,p}_{\tM} \cap \ran(\nabla_x) \arrow[d, "d"]
	\\
	{\IX}^{s,p}_{BD} \cap \ran(BD)  &=&		\IX^{s,p}_{L} \cap \ran(a^{-1}\div_x) & \oplus & 	\IX^{s,p}_{M}	\cap \ran(d \nabla_x)		
\end{tikzcd}
\caption{Splittings and identifications of pre-Hardy--Sobolev and pre-Besov spaces. Each arrow indicates a bijection that is bounded from below and above for the respective $\IX$-quasinorms. Domains and ranges are taken for the corresponding operators on $\L^2$ with maximal domain. Each appearing space is the intersection of an adapted space $\IX_T^{s,p}$ with one of its dense subsets, where density is with respect to the norm $\|\cdot\|_{\IX_T^{s,p}} + \|\cdot \|_2$, see Lemma~\ref{lem: domain dense in Hardy}.\index{Hardy--Sobolev space!for $DB, BD, L, M$}}
\label{fig: diagram with X}	
\end{figure}

As for the mapping between the second and third row in Figure~\ref{fig: diagram with X}, we first cite the following \emph{regularity shift}\index{regularity shift!for $\IX_{BD}^{s,p}$} from \cite[Prop.~5.6]{AA}: we have that
\begin{align}
\label{eq: regularity shift via D}
D: \IX^{s+1,p}_{BD} \cap \dom(D) \to \IX^{s, p}_{DB} \cap \ran(D) 
\end{align}
is bijective and bounded from below and above for the $\IX$-quasinorms.  In particular, 
\begin{align*}
\|Df\|_{\IX^{s,p}_{DB}} \simeq \|f\|_{\IX^{s+1,p}_{BD}} \quad (f \in \IX^{s+1,p}_{BD} \cap \dom(D)).
\end{align*}
This takes care of the left-hand side. The two ingredients for the proof in \cite{AA} are the intertwining property from Lemma~\ref{lem: intertwining} and the following 

\begin{lem}[Local coercivity inequality\index{local coercivity inequality (for $B$)}, {\cite[Lem.~5.14]{AusSta}}]
\label{lem: local coercivity inequality}
For any $u \in \Lloc^2$ with $Du \in \Lloc^2$ and any ball $B(x,t) \subseteq \R^n$ it follows that
\begin{align*}
\int_{B(x,t)} |Du|^2 \lesssim \int_{B(x,2t)} |BDu|^2 + t^{-2} \int_{B(x,2t)} |u|^2.
\end{align*}
\end{lem}

\begin{rem}
In Lemma~\ref{lem: local coercivity inequality} we understand $Du = [\div_x u_\ta, -\nabla_x u_\no]^\top$ in the sense of distributions. In particular, we can take $u \in \dom(D)$.
\end{rem}

On recalling $\dom(D) = \dom(BD) = \dom([BD])$ from Section~\ref{subsec: Kato and Riesz} and $\ran(D) = \ran(DB)$ from \eqref{eq: ran BD and DB}, we can split the regularity shift \eqref{eq: regularity shift via D} in the spirit of \eqref{eq: Hardy splitting L and M} and obtain the right-hand side between the second and third row.

Similarly, the mappings between the first and second row in Figure~\ref{fig: diagram with X} are due to \eqref{eq: Hardy splitting L and M} and Corollary~\ref{cor: Lifting pre-Hardy}. 

Finally, the mapping from the third to the fourth line follows from the block diagonal structure of $B$ and the following

\begin{lem}
\label{lem: B on HpDB}
Let $s \in \R$ and $p \in (0,\infty)$. The map 
\begin{align*}
	B : \IX^{s,p}_{DB} \cap \ran(D) \to \IX^{s,p}_{BD} \cap \ran(BD)
\end{align*}
is bijective and bounded from below and above for the respective $\IX^{s,p}$-quasinorms.
\end{lem}

\begin{proof}
Let $f \in \IX_{DB}^{s,p} \cap \ran(D)$. We have $Bf \in \ran(BD)$ and for any $\psi \in \Psi_+^+$ we obtain $\IQ_{\psi, BD} Bf = B \IQ_{\psi, DB}f$ from \eqref{eq: similarity BD and DB}. Since $B$ is a bounded multiplication operator, we conclude $\|Bf\|_{\IX^{s,p}_{BD}} \leq \|B\|_\infty \|f\|_{\IX^{s,p}_{DB}}$.

Conversely, let $g \in \IX_{BD}^{s,p} \cap \ran(BD)$ and write $g = Bf$ with $f \in \ran(D)$. In order to bound $f$ in $\IX^{s,p}_{DB}$, we take an auxiliary function $\psi \in \Psi_\infty^\infty$ and define $\varphi \in \Psi_\infty^\infty$ by $\varphi(z) = z \psi(z)$. For fixed $\tau>0$ we have again the intertwining relation $D \psi(\tau BD)g = DB \psi(\tau DB)f$. The local coercivity inequality applied to $u \coloneqq \tau \psi(tBD)g$ can therefore be rewritten as
\begin{align*}
\int_{B(x,t)} |\varphi(\tau DB)f|^2 \lesssim \int_{B(x,2t)} |\varphi(\tau BD)g|^2 + \int_{B(x,2t)} |\psi(\tau BD)g|^2.
\end{align*}
Consequently, 
\begin{align*}
\|\IQ_{\varphi, DB}f\|_{\Y^{s,p}} \lesssim \|\IQ_{\varphi, BD}g\|_{\Y^{s,p}} + \|\IQ_{\psi, BD}g\|_{\Y^{s,p}},
\end{align*}
where in the case $(\IX, \Y) = (\IH, \T)$ we also used a change of angle in the tent space norms. The left-hand side compares to $\|f\|_{\IX^{s,p}_{DB}}$ whereas both terms on the right compare to $\|g\|_{\IX^{s,p}_{BD}}$.
\end{proof}

We could also write down a `completed' version of Figure~\ref{fig: diagram with X} in which all pre-Hardy--Sobolev and pre-Besov spaces are replaced by their canonical completions and all intersections vanish. While conceptually this might seem more satisfactory, the possibility of working with invertible maps in $\L^2$ will have significant advantages for many of our proofs.
\section{Identification of adapted Hardy spaces}
\label{sec: Identification of adapted spaces}

\noindent This section is concerned with identifying three pre-Hardy spaces,  $\IH_L^p$, $\IH_L^{1,p}$ and $\IH_{DB}^p$, that play a crucial role for Dirichlet and regularity problems, with classical smoothness spaces. To this end it will be convenient to have a version of Figure~\ref{fig: diagram with X} around these particular spaces at hand:

\begin{figure}[ht]
\centering
\begin{tikzcd}[column sep=0pt, row sep=40pt]
	\IH^p_{BD} \cap \ran([BD])  &=&	\IH^p_{L} \cap \ran({L^{1/2}})  & \oplus & \IH^p_{M} \cap \ran(M^{1/2})\\
	\IH^{1,p}_{BD} \cap \dom(D) \arrow[u, "{[BD]}"] \arrow[d, "D"] &=&	\IH^{1,p}_{L} \cap \dom({L^{1/2}}) \arrow[u,"{L^{1/2}}"]\arrow[drr, "\! \! \! \!{-\nabla_x}" near start] & \oplus &\arrow[dll, "\! \! \! \! {\div_x}" near end] \arrow[u, "M^{1/2}"] \IH^{1,p}_{M} \cap \dom(M^{1/2}) \\
	\IH^p_{DB} \cap \ran(D) \arrow[d, "B"] &=&	\IH^p_{\tL} \cap \ran(\div_x) \arrow[d, "{a^{-1}}"]  &\oplus &  \IH^p_{\tM} \cap \ran(\nabla_x) \arrow[d, "d"]
	\\
	{\IH}^p_{BD} \cap \ran(BD)  &=&		\IH^p_{L} \cap \ran(a^{-1}\div_x) & \oplus & 	\IH^p_{M}	\cap \ran(d \nabla_x)		
\end{tikzcd}
\caption{Figure~\ref{fig: diagram with X} for $s=0$ and $\IX = \IH$. 
Each appearing space is the intersection of an adapted space with one of its dense subsets.}
\label{fig: diagram}	
\end{figure}

As for the second and third row `identifying'\index{identification!of abstract and concrete spaces} means determining whether the spaces remain the same as sets and with equivalent $p$-quasinorms when $B$ is replaced by the identity matrix. In the fourth row for $\IH_L^p$, we can then expect it is the image of $\H^p \cap \L^2$ under multiplication with $a^{-1}$. If $p>1$, then multiplication by $a^{-1}$ is invertible on $\H^p = \L^p$ and hence the image is the same as  $\L^p \cap \L^2$.
\subsection{Identification regions}
\label{subsec: identification regions}
We introduce three such sets of exponents:\index{H@$\cH(DB)$}\index{identification region! for $\IH_L^p$, $\IH_L^{1,p}$, $\IH_{DB}^p$} \index{H@$\cH(L)$}\index{H1@$\cH^1(L)$}
\begin{align}
\label{eq: H(DB)}
\cH(DB) &\coloneqq \big\{ p \in (1_*, \infty) : \|f\|_{\IH^{p}_{DB}} \simeq \|f\|_{\IH^p_D} \text{ for all } f \in \cl{\ran(D)} \big\}
\intertext{and}
\notag \cH(L) &\coloneqq \big\{ p \in (1_*, \infty) : \|f\|_{\IH^p_L} \simeq \|a f\|_{\H^p} \text{ for all } f \in \L^2 \big\}, \\
\notag \cH^1(L)&\coloneqq \big\{ p \in (1_*, \infty) : \|f\|_{\IH^{1,p}_L} \simeq \|f\|_{\Hdot^{1,p}} \text{ for all } f \in \L^2 \big\}.
\end{align}

The identification region for $DB$ turns out to be the intersection of the two regions associated with $L$. This has nothing to do with the particular Hardy spaces above and follows from Figure~\ref{fig: diagram} for all sorts of adapted spaces. Identification regions for other $DB$ and $L$-adapted spaces will appear much later in the text in Section~\ref{sec: fractional}.

\begin{lem}
\label{lem: intersection of identification regions}
Let $s \in \R$ and $p \in (0,\infty)$. The following are equivalent : 
\begin{enumerate}
	\item $\IX^{s,p}_{DB} = \IX^{s,p}_D$ with equivalent $p$-quasinorms.
	\item $\IX^{s,p}_L = a^{-1}(\Xdot^{s,p} \cap \L^2)$ and $\IX^{s+1,p}_L = \Xdot^{s+1,p} \cap \L^2$, both with equivalent $p$-quasinorms.
\end{enumerate}
\end{lem}

Specializing to $\IX= \IH$ and $s = 0$ in Lemma~\ref{lem: intersection of identification regions}, we obtain

\begin{cor}
\label{cor: easy relations between identification regions} 
It follows that $\cH(DB) = \cH(L) \cap \cH^1(L)$. In particular (by \eqref{eq: Hardy identification for p=2}) all three sets contain $p=2$.
\end{cor}

\begin{proof}[Proof of Lemma~\ref{lem: intersection of identification regions}]
Throughout, equalities of spaces are up to comparable pre-Hardy quasinorms and spaces that arise from multiplication with $a^{-1}$ carry the image topology. 

We start by noting that (i) is equivalent to $\IX^{s,p}_{DB} \cap \ran(D) = \IX^{s,p}_D \cap \ran(D)$ since $\ran(D) = \ran(DB)$ is dense in both adapted spaces. The third row of Figure~\ref{fig: diagram with X} yields equivalence to 
\begin{align*}
	\IX^{s,p}_{\tL} \cap \ran(\div_x) &= \IX^{s,p}_{-\Delta_x} \cap \ran(\div_x), \\
	\IX^{s,p}_{\tM} \cap \ran(\nabla_x) &= \IX^{s,p}_{-\nabla_x \div_x} \cap \ran(\nabla_x).
	\intertext{By moving to the fourth and second row, this is the same as having}
	\IX^{s,p}_{L} \cap \ran(a^{-1} \div_x) &= a^{-1}(\IX^{s,p}_{-\Delta_x} \cap \ran(\div_x)), \\
	\IX^{s+1,p}_{L} \cap \dom(\nabla_x) &= \IX^{s+1,p}_{-\Delta_x} \cap \dom(\nabla_x),
	\intertext{which, by density, is equivalent to having}
	\IX^{s,p}_{L} &= a^{-1} \IX^{s,p}_{-\Delta_x} \\
	\IX^{s+1,p}_{L} &= \IX^{s+1,p}_{-\Delta_x}.
\end{align*}
The spaces associated with the Laplacian have been identified in Figure~\ref{fig: identification D-adapted} and equivalence to (ii) follows.
\end{proof}

\begin{rem}
\label{rem: easy relations between identification regions}
The argument above proves slightly more: it says that we have, all in the sense of continuous inclusions, $\IX^{s,p}_{DB} \subseteq \IX^{s,p}_D$ if and only if we have both $\IX^{s,p}_{L} \subseteq a^{-1} (\Xdot^{s,p} \cap \L^2)$ and $\IX^{s+1,p}_{L} \subseteq \Xdot^{s+1,p} \cap \L^2$, and that the same result holds upon reversing all inclusions.
\end{rem}

In order to show that the identification regions are intervals, we borrow an interpolation argument from \cite[Thm.~4.32]{AA} that uses the canonical completions of adapted Hardy spaces. In fact, for $\cH(DB)$ and $\cH^1(L)$ the result in \cite{AA} would apply `off-the-shelf' but  a slight variant is needed for $\cH(L)$ because of the multiplication by $a$.

\begin{lem}
\label{lem: identification regions are interval}
The sets $\cH(DB)$, $\cH(L)$ and $\cH^1(L)$ are intervals.
\end{lem}

\begin{proof}
We begin with the proof for $\cH(L)$. By definition, we have $p \in \cH(L)$ if and only if the multiplication operators $a: \IH_L^p \to \H^p \cap \L^2$ and $b \coloneqq a^{-1}: \H^p \cap \L^2 \to \IH_L^p$ are well-defined and bounded for the $p$-quasinorms. This is equivalent to saying that these operators have bounded extensions $\hat{a}: \psiH_L^{p} \to \H^{p}$ and $\hat{b}: \H^{p} \to \psiH_L^{p}$ to canonical completions in the sense that the following diagrams commute:
\begin{center}
	\begin{tikzcd}[column sep=40pt]
	\psiH_L^{p} \arrow[r, "\hat{a}"] 
	& \H^{p}
	\\
	\IH_{L}^{p} \arrow[hookrightarrow]{u} \arrow[r, "a"]
	& \H^{p} \cap \L^2 \arrow[hookrightarrow]{u} 
	\end{tikzcd}
\qquad
	\begin{tikzcd}[column sep=40pt]
	\H^p \arrow[r, "\hat{b}"] 
	& \psiH_L^{p} 
	\\
	\H^p \cap \L^2 \arrow[hookrightarrow]{u} \arrow[r, "b"]
	&  \IH_{L}^{p} \arrow[hookrightarrow]{u} 
\end{tikzcd}.
\end{center}
Let now $p_0,p_1 \in \cH(L)$. Since the spaces $\H^p$ and $\psiH_L^p$ have universal approximation techniques, the extensions in the respective diagrams for $p_0$ and $p_1$ are compatible and we can use complex interpolation (Section~\ref{subsec: complex interpolation} and Proposition~\ref{prop: interpolation of canonical completions}) to obtain the same diagrams for all $p$ between $p_0$ and $p_1$. Hence, these exponents are all in $\cH(L)$.

The argument for $\cH^1(L)$ is identical except that we extend the identity operator. The same for $\cH(DB)$ but instead of $\H^p$ we use a canonical completion $\psiH_D^p$.
\end{proof}

\begin{rem}
\label{rem: identification regions are intervals}
Again, the argument above has given a stronger statement about inclusions: The set of exponents for which for instance $\IH_{DB}^p \subseteq \IH_D^p$ holds with continuous inclusion for the $p$-quasinorms is an interval and the same is true for the other five possible inclusions.
\end{rem}

\begin{defn}
\label{def: h+-(DB)}
\begin{enumerate}
	\item The upper and lower endpoints of $\cH(DB)$ are denoted by $h_{-}(DB)$ and $h_{+}(DB)$. 
	\item The upper and lower endpoints of $\cH(L)$ are denoted by $h_{-}(L)$ and $h_{+}(L)$. Likewise $h_\pm^1(L)$ are the endpoints of $\cH^1(L)$.
\end{enumerate}
\end{defn}
\subsection{The identification theorem}
\label{subsec: identification theorem}

We come to the characterization of the identification region's endpoints through the critical numbers $p_-(L)$ and $q_\pm(L)$.

\begin{thm}[Identification Theorem\index{identification Theorem}]
\label{thm: main result Hardy}
The endpoints of $\cH(L)$ and $\cH^1(L)$ can be characterized and controlled as follows:
\begin{align*}
h_\pm(L) &= p_\pm(L),\\
h^1_{-}(L) &\leq (p_{-}(L)_* \vee 1_*),\\
h^1_+(L) &= q_+(L).
\end{align*}
As a consequence, the endpoints of $\cH(DB)$ are $h_-(DB) = p_-(L)$ and $h_+(DB) = q_+(L)$.\index{H@$\cH(L)$! endpoints of}\index{H1@$\cH^1(L)$! endpoints of}\index{H@$\cH(DB)$! endpoints of}
\end{thm}

The relations for $L$ imply those for $DB$ since $\cH(DB) = \cH(L) \cap \cH^1(L)$ and $q_+(L) \leq p_+(L)$ by Theorem~\ref{thm: standard relation J(L) and N(L)}. We later precise this result by showing that these intervals are open at their ends except may be at the lower endpoint of $\cH^1(L)$ for which we cannot even say whether the bound is sharp.

The proof of Theorem~\ref{thm: main result Hardy} is spread over $10$ parts, using different methods for different regimes of parameters. Upper bounds on the size of $\cH(L)$ are easy to obtain (Part 1), whereas lower bounds require establishing two continuous inclusions. Parts 2 - 5 focus on different inclusions of classical and $L$-adapted spaces. Parts 6 - 10 contain the synthesis of these preparatory steps. 

Many arguments are known when $a=1$. However, there are still some new difficulties when $a\ne 1$ that need to be taken care of and for some other parts we can simplify known arguments through the full strength of Figure~\ref{fig: diagram} even when $a=1$. 	
\medskip

\subsection*{Part 1: \texorpdfstring{$\boldsymbol{p_{-}(L) \leq h_{-}(L)}$ and $\boldsymbol{p_{+}(L) \geq h_{+}(L)}$}{p-(L) < h-(L) and p+(L) > h+(L)}}
Being slightly more precise, we show the inclusion $\cH(L) \subseteq \cJ(L)$. Given $p \in (1_*,\infty)$, Proposition~\ref{prop: FC on pre-Hardy} yields
\begin{align*}
\|(1+t^2 L)^{-1} f\|_{\IH^{p}_L} \lesssim \|f\|_{\IH^{p}_L} 
\end{align*}
uniformly for all $f \in \IH^{p}_L$ and all $t>0$. If now $p \in \cH(L)$, then $\IH^p_L = a^{-1} (\H^p \cap \L^2)$ holds with equivalent Hardy norms and $p \in \cJ(L)$ follows.
\subsection*{Part 2: \texorpdfstring{$\boldsymbol{\L^p \cap \L^2 \subseteq \H^p_L}$ for $\boldsymbol{2 \leq p<\infty}$}{Inclusion Lp C IHp for 2<p<oo}}

We are going to prove the continuous inclusion $\L^q \cap \L^2 \subseteq \IH^q_L$ for $q \in [2,\infty)$.

We define $\IH^q_L$ via the auxiliary function $\psi(z) \coloneqq z^\alpha(1+z)^{-2\alpha}$ with an integer $\alpha > \nicefrac{n}{2}$, so that this choice is admissible for all $q$, see Section~\ref{subsec: Hardy abstract sectorial}. We have to establish the bound
\begin{align*}
 \|\IQ_{\psi,L}f\|_{\T^q} \lesssim \|f\|_q \quad (f \in \L^q \cap \L^2).
\end{align*}
For a later purpose, we prove a more general statement. This uses the standard assumptions from Section~\ref{subsec: Hardy abstract sectorial}. For $T=L$ the bound required here follows by simply taking the auxiliary parameters $\theta = 1$ and $p=2$. The further interest in the lemma lies in picking $p$ as large and $\theta$ as small as possible in order to allow for weaker decay assumption of $\psi$ at the origin.

\begin{lem}
\label{lem: Carleson bound Q function}
Let $T$ be a sectorial operator that satisfies the standard assumptions \eqref{eq: standard assumptions sectorial}. Fix $\mu \in (0,\nicefrac{(\pi-\omega)}{2})$ and $\sigma, \tau >0$. Let $\psi \in \Psi_\sigma^\tau(\S_{\pi - 2\mu}^+)$ and consider the square function bound
\begin{align*}
\|\IQ_{\psi, T}f \|_{\T^{q}} \lesssim \|\psi\|_{\sigma, \tau, \mu} \|f\|_{q} \quad (f \in \L^q \cap \L^2),
\end{align*}
where the implicit constant does not depend on $\psi$. Then this bound is valid for every $q \geq 2$ provided that one can find $p \in [2,\infty)$ and $\theta \in (0,1]$ such that $((1+t^2T)^{-1})_{t>0}$ is $\L^p$-bounded and
\begin{align*}
\mu \in \Big(0, \frac{\theta(\pi - \omega)}{2}\Big) \quad \& \quad \sigma > \frac{n}{2[p,2]_\theta}.\end{align*}
\end{lem}

\begin{proof}
In the following implicit constants are allowed to depend on the fixed parameters but not on $\psi$ itself. Via McIntosh's theorem the boundedness for $q=2$ is equivalent to the bounded $\H^\infty$-calculus on $\cl{\ran(T)}$. Hence, we can state
\begin{align*}
\|\IQ_{\psi, T}f \|_{\T^2} \lesssim  \|\psi\|_{\sigma,\tau,\mu} \|f\|_2 \quad (f \in \cl{\ran(T)}).
\end{align*}
Cauchy's theorem yields $\psi(t^2T)f = 0$ for all $t>0$ if $f \in \nul(T)$. Hence, we can state the same bound for all $f \in \L^2$. By complex interpolation it remains to treat the case $q=\infty$, that is to say, to prove for all balls $B \subseteq \R^n$ of radius $r>0$ and all $f \in \L^\infty$ that
\begin{align}
\label{eq: goal H(L) for p>2}
\bigg(\frac{1}{|B|} \int_0^r \int_B |\psi(t^2T)f|^2 \, \frac{\d x \d t}{t} \bigg)^{1/2} \lesssim \|\psi\|_{\tau,\sigma,\mu} \|f\|_\infty.
\end{align}
Here, $\psi(t^2 T)$ is extended to $\L^\infty$ via Proposition~\ref{prop: general OD extension} and we shall see in the further course of the proof that we have the required off-diagonal decay at our disposal.

Having fixed $B$, we write $f = \sum_{j \geq 1} f_j$ with $f_j \coloneqq \ind_{C_j(B)}f$. For $j = 1$ we use that $\T^2 = \L^2(\reu; \frac{\d t \d x}{t})$ and again the $\L^2$-bound to give
\begin{align*}
\frac{1}{|B|} &\int_0^r \int_B |\psi(t^2T)f_1(x)|^2 \, \frac{\d x \d t}{t} \\ 
&\lesssim \frac{1}{|B|} \|\IQ_{\psi,T} f_1\|_{\T^2}^2 \\
&\lesssim \frac{1}{|B|} \|\psi\|_{\sigma,\tau,\mu}^2 \|f_1\|_2^2  \\
&\leq 4^n \|\psi\|_{\sigma,\tau,\mu}^2 \|f\|_\infty^2.
\end{align*}
Next, we let $\varrho \coloneqq [p,2]_\theta$ and obtain from Lemma~\ref{lem: functional calculus bounds from J(L) abstract}.(i) the off-diagonal estimate
\begin{align*}
\|\psi(t^2L)f_j\|_{\L^\varrho(B)}
&\lesssim \|\psi\|_{\tau,\sigma,\mu} \bigg(1+\frac{2^jr}{t} \bigg)^{-2 \sigma} \|f\|_{\L^\varrho(C_j(B))} \\
&\leq \|\psi\|_{\tau,\sigma,\mu}  \frac{t^{2 \sigma}}{r^{2 \sigma}} 2^{-j(2 \sigma - \frac{n}{\varrho})} |B|^{1/\varrho} \|f\|_\infty.
\end{align*}
Since $\varrho\geq 2$, we obtain from H\"older's inequality that
\begin{align*}
\|\psi(t^2L)f_j\|_{\L^2(B)} \lesssim \|\psi\|_{\tau,\sigma,\mu}  \frac{t^{2 \sigma}}{r^{2 \sigma}} 2^{- j (2\sigma - \frac{n}{\varrho})} |B|^{1/2} \|f\|_\infty
\end{align*}
and taking $\L^2$-norms with respect to $\nicefrac{\d t}{t}$, we are led to
\begin{align*}
\bigg(\frac{1}{|B|} &\int_0^r \int_B |\psi(t^2L)f_j(x)|^2 \, \frac{\d x \d t}{t} \bigg)^{1/2} \\
&\lesssim  \|\psi\|_{\tau,\sigma,\mu} 2^{-j(2 \sigma -  \frac{n}{\varrho})} \|f\|_\infty \bigg(\int_0^r \frac{t^{4 \sigma}}{r^{4 \sigma}} \, \frac{\d t}{t}\bigg)^{1/2} \\
&= 2^{-j(2 \sigma -  \frac{n}{\varrho})} \|\psi\|_{\tau,\sigma,\mu} (4 \sigma)^{-\frac{1}{2}} \|f\|_\infty.
\end{align*}
By assumption, we have $2 \sigma > \nicefrac{n}{\varrho}$. Summing up in $j$ yields \eqref{eq: goal H(L) for p>2}.
\end{proof}

\begin{rem}
\label{rem: Carleson bound Q function}
It becomes clear from the proof above that Lemma~\ref{lem: Carleson bound Q function} has very little to do with sectorial operators and could be extended to more general extensions
\begin{align*}
(\IQ_{\psi} f)(t,x) \coloneqq (\psi(t) f) (x)
\end{align*}
where $\psi: (0,\infty) \to \Lop(\L^2)$ is a strongly measurable family of operators. For example, with $p=2$ and $\theta = 1$ the only properties of $(\psi(t))_{t>0}$ that we have used to get for every $q \geq 2$ a bound
\begin{align*}
\|\IQ_{\psi}f\|_{\T^q} \lesssim \|f\|_q \quad (f \in \L^q \cap \L^2),
\end{align*}
is the corresponding $\L^2$-bound and $\L^2$ off-diagonal estimates of order $\gamma > \nicefrac{n}{2}$.
\end{rem}
\subsection*{Part 3: \texorpdfstring{Injection of classical spaces into $\boldsymbol{L}$-adapted spaces for $\boldsymbol{p \in (1,2)}$}{Injection of classical spaces into L-adapted spaces for 1<p<2}}

For this part we work with the auxiliary function $\psi$ defined by
\begin{align}
\label{eq: auxiliary function for weak type bound Hardy}
\psi_\alpha(z) \coloneqq z^{\alpha-1/2} (1+z)^{-3 \alpha},
\end{align}
where $\alpha \in \IN$ will be chosen sufficiently large depending on exponents and dimensions. Throughout this part it will be convenient to write
\begin{align}
\label{eq: SF for weak type bound Hardy}
S_{\psi_\alpha,L} f(x) \coloneqq \bigg(\iint_{|x-y|<t} |\psi_\alpha(t^2L)f(y)|^2 \, \frac{\d t \d y}{t^{1+n}} \bigg)^{1/2},
\end{align}
so that $\|S_{\psi_\alpha,L} \cdot \|_p$ becomes an equivalent norm on $\IH^p_L$ provided that $\alpha > \nicefrac{n}{(2p)}-\nicefrac{n}{4}$, compare with Section~\ref{subsec: Hardy abstract sectorial}.

Our main objective is to establish the following extrapolation result for square functions.

\begin{lem}
\label{lem: weak type bound Hardy}
Suppose for some $q \in (p_{-}(L) \vee 1 , 2]$ and all sufficiently large $\alpha$ (depending on $q, p_-(L), n$) that
\begin{align}
\label{eq: weak type bound hardy assumption}
\|S_{\psi_\alpha,L} ({L^{1/2}} u) \|_q \lesssim \|\nabla_x u\|_q \quad (u \in \Wdot^{1,q} \cap \W^{1,2}).
\end{align}
Then for all $p \in (q_* \vee 1 , q)$ and all sufficiently large $\alpha$ (depending on $p,q, p_-(L),n$) it follows that
\begin{align*}
\|S_{\psi_\alpha,L} ({L^{1/2}} u) \|_p \lesssim \|\nabla_x u\|_p \quad (u \in \Wdot^{1,p} \cap \W^{1,2}).
\end{align*}
\end{lem}

\begin{rem}
\label{rem: weak type bound Hardy}
Assumption \eqref{eq: weak type bound hardy assumption} holds for $q=2$ and any $\alpha \in \IN$. Indeed, this follows from $\IH^2_L = \cl{\ran(L)} = \L^2$ and the solution of the Kato problem. Starting from there, we can iterate Lemma~\ref{lem: weak type bound Hardy} in order to conclude that for every $q \in (p_{-}(L)_* \vee 1, 2]$ the bound \eqref{eq: weak type bound hardy assumption} holds for all sufficiently large $\alpha$.
\end{rem}

Before giving the proof of Lemma~\ref{lem: weak type bound Hardy}, let us state the more important consequences of this lemma for the identification of $L$-adapted Hardy spaces. 

\begin{prop}
\label{prop: Hp in HpL for 1<p<2}
If $p \in (p_{-}(L) \vee 1, 2]$, then $\L^p \cap \L^2 \subseteq \IH^p_L$ with continuous inclusion for the $p$-norms.  
\end{prop}

\begin{proof}
First let us assume $f \in \L^p \cap  \ran({L^{1/2}})$. By Lemma~\ref{lem: range density in Lp} this is a dense subspace of $\L^p \cap \L^2$. We put $u \coloneqq L^{-1/2} f$. Since the Riesz transform is $\L^p$-bounded (Theorem~\ref{thm: Riesz}), we have $u \in \Wdot^{1,p} \cap \W^{1,2}$ with $\|\nabla_x u\|_p \lesssim \|f\|_p$. Remark~\ref{rem: weak type bound Hardy} yields
\begin{align*}
\|S_{\psi_\alpha,L} f \|_p \lesssim \|f \|_p
\end{align*} 
if $\alpha$ is sufficiently large. If in addition $\alpha > \nicefrac{n}{(2p)}-\nicefrac{n}{4}$, then $\psi_\alpha$ is admissible as auxiliary function for $\IH^p_L$ and we obtain
\begin{align*}
\|f\|_{\IH^p_L} \lesssim \|f\|_p
\end{align*}
with an implicit constant independent of $f$. A general $f \in \L^p \cap \L^2$ can be approximated by $(f_j) \subseteq  \L^p \cap \ran({L^{1/2}})$ in $\L^p \cap \L^2$. By $\L^2$-convergence
\begin{align*}
\int_{B(x,t)} |\psi_\alpha(t^2L)f(y)|^2 \, \d y = \lim_{j \to \infty} \int_{B(x,t)} |\psi_\alpha(t^2L)f_j(y)|^2 \, \d y
\end{align*}
holds for all $(t,x) \in \reu$ and we invoke Fatou's lemma to give
\begin{align*}
\|f\|_{\IH^p_L}^p 
&\leq \liminf_{j \to \infty} \int_{\R^n} \bigg(\iint_{|x-y|<t} |\psi_\alpha(t^2L)f_j(y)|^2 \, \frac{\d t \d y}{t^{1+n}} \bigg)^{p/2} \, \d x \\
&= \liminf_{j \to \infty} \|f_j\|_{\IH^p_L}^p.
\end{align*}
On the right-hand side $ \|f_j\|_{\IH^p_L}^p$ is under control by $\|f_j\|_p^p$ thanks to the first part of the proof and $\L^p$-convergence  of $(f_j)$ gives the required bound by $\|f\|_p^p$.
\end{proof}

\begin{prop}
\label{prop: H1p in H1pL for 1<p<2}
If $p \in (p_{-}(L)_* \vee 1, 2]$, then $\Wdot^{1,p} \cap \L^2 \subseteq \IH^{1,p}_L$ with continuous inclusion for the $p$-norms.  
\end{prop}

\begin{proof}
By the universal approximation technique even  $\cZ$ is dense in $\Wdot^{1,p} \cap \L^2$. Hence, the same approximation argument as in the previous proof shows that it suffices to check
\begin{align*}
 \|u\|_{\IH^{1,p}_L} \lesssim \|\nabla_x u\|_p \quad (u \in \Wdot^{1,p} \cap \W^{1,2}).
\end{align*}
We take $\alpha$ large enough so that \eqref{eq: weak type bound hardy assumption} holds at exponent $q=p$ and $\alpha > \nicefrac{n}{(2p)}-\nicefrac{n}{4}$ to make sure that $\IH^{1,p}_L$ can be defined through the auxiliary function $\varphi_\alpha(z) \coloneqq \sqrt{z} \psi_\alpha(z)$. We have $\psi_\alpha(t^2 L) {L^{1/2}} u = t^{-1} \varphi_\alpha(t^2L)u$ for $t>0$ and therefore
\begin{align*}
\|u\|_{\IH^{1,p}_L} 
= \|\IQ_{\varphi_\alpha,L} u\|_{\T^{1,p}} 
&= \|\IQ_{\psi_\alpha,L} ({L^{1/2}} u)\|_{\T^{0,p}} \\
&= \|S_{\psi_\alpha, L} ({L^{1/2}} u) \|_p \\
&\lesssim \|\nabla_x u\|_p. \qedhere
\end{align*}
\end{proof}

We come to the proof of Lemma~\ref{lem: weak type bound Hardy}. We modify the strategy of \cite[pp.42-45]{AusSta}. Henceforth we fix $p,q$ as in the statement and we write $\psi = \psi_\alpha$, where $\alpha$ will be chosen larger from step to step in dependence of $p, q, p_-(L), n$. 

Let $u \in \Wdot^{1,p} \cap \W^{1,2}$ and $\lambda > 0$. It will be enough to obtain the weak-type estimate
\begin{align}
\label{eq: goal weak type bound Hardy}
\Big| \Big \{ x \in \R^n : S_{\psi,L} ({L^{1/2}} u)(x) > 3 \lambda \Big \} \Big| \lesssim \frac{1}{\lambda^p} \|\nabla_x u\|_p^p 
\end{align}
with implicit constant independent of $u$ and $\lambda$. Indeed, consider the positive sublinear operator 
\begin{align*}
T: \cZ \to \L^2, \quad Tu \coloneqq S_{\psi,L}({L^{1/2}}(-\Delta_x)^{-1/2}u)
\end{align*}
and recall that $\cZ$ is dense in all (intersections of) $\L^r$-spaces with $r>1$. Now, $T$ is of strong type $(q,q)$ by \eqref{eq: weak type bound hardy assumption} and of weak type $(p,p)$ by \eqref{eq: goal weak type bound Hardy} since $\nabla_x$ and $(-\Delta_x)^{1/2}$ are comparable in $\L^p$-norm by the Mihlin multiplier theorem. Hence, it is of strong type $(r,r)$ for every $r \in (p,q]$ by the Marcinkiewicz interpolation theorem. As $(-\Delta_x)^{-1/2}$ is invertible on $\cZ$, this means that we have
\begin{align*}
\|S_{\psi,L} ({L^{1/2}} u) \|_r \lesssim \|\nabla_x u\|_r \quad (u \in \cZ).
\end{align*}
This bound extends to $u \in \Wdot^{1,r} \cap \W^{1,2}$ by density as before. Since $p \in (q_* \vee 1, q)$ and $r \in (p,q]$ were arbitrary, the claim follows.

The proof of \eqref{eq: goal weak type bound Hardy} itself comes in 8 steps.

\medskip

\noindent \emph{Step 1: Calder\'{o}n--Zygmund decomposition.\index{Calder\'{o}n--Zygmund decomposition!for Sobolev functions}}  We use the decomposition for Sobolev functions that was introduced in \cite[Lem.~4.12]{A}, see \cite{A-Note} for the correction of an inaccuracy in the original proof. 

Since $u \in \Wdot^{1,p}$, according to this decomposition, there is a countable collection of cubes $(Q_j)_{j \in J}$, measurable functions $g$ and $b_j$ and constants $C$ and $N$ that depend only on dimensions and $p$, such that
\begin{enumerate} 
	\item  $\displaystyle u = g + \sum_{j \in J} b_j$ pointwise almost everywhere, 
	
	\item  $\displaystyle \|\nabla_x g\|_\infty \leq C \lambda$, \\[-5pt]
	
	\item  $b_j$ has support in $Q_j$ and $\|\nabla_x b_j\|_p^p  \leq C \lambda^p |Q_j|$, \\[0pt]
	
	\item  $\displaystyle \sum_{j \in J} |Q_j| \leq C \lambda^{-p} \|\nabla_x u\|_p^p$,
	
	\item  $\displaystyle \sum_{j \in J} \ind_{Q_j} \leq N$.
\end{enumerate}
More precisely, setting $\Omega \coloneqq \{\Max(|\nabla_x u|^p) > \lambda^p\} \subseteq \R^n$, the $b_j$ take the form $b_j = (u-u(x_j))\chi_j$ with $x_j \in2Q_j \cap {}^c \Omega$ and $\chi_j \in \C_0^\infty(Q_j)$ such that $\|\chi_j\|_\infty + \ell(Q_j)\|\nabla_x \chi_j\|_\infty \leq C$. The function $u$ has a representative on ${}^c\Omega$ that satisfies $|u(x) - (u)_Q| \leq C \lambda \ell(Q)$ whenever $Q$ is a cube centered at $x \in {}^c \Omega$ and this is how we understand $u(x_j)$. 

We recall these details on the construction because we need two additional properties in the proof of \eqref{eq: goal weak type bound Hardy}:
\begin{enumerate}
	\item[\quad (i')] If $u \in \Wdot^{1,r}$ for some $r\in (1,\infty)$, then $b_j \in \W^{1,r}$ for all $j$ and $\sum_{j \in J} b_j$ converges unconditionally in $\Wdot^{1,r}$.
	
	\item[\quad (ii')] If $r \in [p, p^*)$, then $\|\nabla_x g\|_r^r \leq C' \lambda^{r-p} \|\nabla_x u\|_p^p$ and $\|b_j\|_r^r  \leq C' \lambda^r |Q_j|^{1 + \frac{r}{n}}$ for all $j$,
    where $C'$ also depends on $r$.
\end{enumerate}

To see property (i'), we let $Q'_j$ be the cube centered at $x_j$ with sidelength $3 \ell(Q_j)$ and write
\begin{align}
\label{eq: i' proof}
b_j 
= (u-(u)_{Q_j})\chi_j + ((u)_{Q_j} - (u)_{Q'_j})\chi_j + ((u)_{Q'_j} - u(x_j))\chi_j.
\end{align}
The special property of $u$ on ${}^c \Omega$ yields $|((u)_{Q'_j} - u(x_j))\nabla_x \chi_j| \leq C \lambda$ on $\R^n$. Next, since $Q_j \subseteq Q'_j$, we obtain from Poincaré's inequality that
\begin{align*}
|(u)_{Q_j} - (u)_{Q'_j}| 
\lesssim \barint_{Q'_j} |u-(u)_{Q'_j}| \, \d x
\lesssim \ell(Q_j) \barint_{Q'_j} |\nabla_x u| \, \d x.
\end{align*}
The right-hand side is bounded by $\lambda$ since $x_j \in {}^c \Omega$ and we obtain $|((u)_{Q_j} - (u)_{Q'_j}) \nabla_x \chi_j| \leq C \lambda$ on $\R^n$. Once again by  Poincaré's inequality we have
\begin{align*}
\int_{\R^n} |\nabla_x((u-(u)_{Q_j}) \chi_j)|^r \, \d x
\lesssim \int_{Q_j} |\nabla_x u|^r \, \d x,
\end{align*}
so that altogether we obtain from \eqref{eq: i' proof} the estimate
\begin{align}
\label{eq1: proof of i'}
\int_{\R^n} |\nabla_x b_j|^r \, \d x \lesssim \lambda^r |Q_j| + \int_{Q_j} |\nabla_x u|^r \, \d x.
\end{align}
Since $b_j$ has compact support, we have $b_j \in \W^{1,r}$ qualitatively. For any partial sum of $j$'s we obtain from (v) and H\"older's inequality that
\begin{align}
\label{eq2: proof of i'}
\int_{\R^n} \bigg|\nabla_x \sum_j b_j  \bigg|^r  \, \d x
\leq N^{r-1} \int_{\R^n} \sum_{j} |\nabla_x b_j|^r \, \d x.
\end{align}
By \eqref{eq1: proof of i'}, (iv) and (v) the right-hand sum has the Cauchy property. Thus, $\sum_{j \in J} \nabla_x b_j$ converges in $\L^r$. The limit is independent of the order of summation since the sum contains at most $N$ non-zero terms at each $x \in \R^n$. 

As for (ii'), the $\L^r$-bound for $b_j$ immediately follows from the Sobolev--Poincar\'e inequality~\cite[Cor.~4.2.3]{Ziemer}\index{inequality!Sobolev--Poincar\'e} and (iii). From \eqref{eq1: proof of i'} and \eqref{eq2: proof of i'} with $r=p$ and then (iv), we obtain $\|\nabla_x \sum_{j \in J} b_j\|_p \leq C \|\nabla_x u\|_p $. We conclude $\|\nabla_x g\|_p \leq C \|\nabla_x u\|_p$ from (i) and the required $\L^r$-bound follows from (ii).

\medskip

\noindent \emph{Step 2: Decomposition of the level set.} For the same $\alpha$ as is the definition of $\psi = \psi_\alpha$ in \eqref{eq: auxiliary function for weak type bound Hardy} we introduce a function $\varphi \in \H^\infty$ through
\begin{align}
\label{eq1: weak type bound Hardy}
\varphi(z) \coloneqq z^\alpha (1+z)^{-\alpha} = (1-(1+z)^{-1})^\alpha
\end{align}
and we decompose $u = g + \wt{g} + b$, using the series
\begin{align*}
\wt g &\coloneqq \sum_{j \in J} (1-\varphi(\ell_j^2 L))b_j,\\
b &\coloneqq \sum_{j \in J} \varphi(\ell_j^2 L)b_j,
\end{align*}
where $\ell_j \coloneqq \ell(Q_j)$. In Step~4 we shall check that the series $\wt{g}$ converges in $\Wdot^{1,q}$, so that by (i') with $r=q$ the same is true for $b$.

Anticipating the convergence of $\wt{g}$, we obtain that the set on the left-hand side of \eqref{eq: goal weak type bound Hardy} is contained in the union of
\begin{align*}
A_1 & \coloneqq \Big \{ x \in \R^n : S_{\psi,L} ({L^{1/2}} g)(x) > \lambda \Big \}, \\
A_2 & \coloneqq \Big \{ x \in \R^n : S_{\psi,L} ({L^{1/2}} \wt{g})(x) > \lambda \Big \}, \\
A_3 & \coloneqq \Big \{ x \in \R^n : S_{\psi,L} ({L^{1/2}} b)(x) > \lambda \Big \},
\end{align*}
where we do not make a notational distinction between $v \mapsto S_{\psi,L} ({L^{1/2}} v)$ and its bounded extension from $\Wdot^{1,q}$ into  $\L^q$. It suffices to bound the measure of each of the three sets by a generic multiple of $\lambda^{-p} \|\nabla_x u\|_p^p$.

\medskip

\noindent \emph{Step 3: Bound of $A_1$.} We use the Markov inequality, the assumption and (ii') to give
\begin{align*}
|A_1| 
\leq \lambda^{-q} \|S_{\psi,L}({L^{1/2}}g)\|_q^q 
\lesssim \lambda^{-q} \|\nabla_x g\|_q^q
\lesssim \lambda^{-p} \|\nabla_x u\|_p^p.
\end{align*}

\medskip

\noindent \emph{Step 4: Convergence and estimate of $\wt{g}$.} For the time being, let $j$ run only through a finite set of $J$. Consider the partial sum of $\wt{g}$ given by
\begin{align}
\label{eq2: weak type bound Hardy}
\sum_j (1-\varphi(\ell_j^2 L))b_j = \sum_{\beta=1}^\alpha \binom{\alpha}{\beta} (-1)^{\beta-1} \bigg(\sum_j (1+\ell_j^2 L)^{-\beta}b_j \bigg),
\end{align}
where we have expanded $\varphi$ from \eqref{eq1: weak type bound Hardy}. We fix $\beta$ and introduce
\begin{align}
\label{eq3: weak type bound Hardy}
f_\beta \coloneqq  \sum_j (1+\ell_j^2 L)^{-\beta}b_j.
\end{align}
Since we have $b_j \in \W^{1,2} = \dom({L^{1/2}})$ by (i'), the same is true for $f_\beta$. We calculate its norm in $\Wdot^{1,q}$ by dualizing $\nabla_x f_\beta$ against $h \in \C_0^\infty$, normalized to $\|h\|_{q'} = 1$:
\begin{align*}
\langle \nabla_x f_\beta, h \rangle 
=  \sum_j \sum_{k = 1}^\infty \langle \nabla_x (1+ \ell_j^2 L)^{-\beta}b_j, h_{j,k} \rangle, 
\end{align*}
where $h_{j,k} \coloneqq \ind_{C_k(Q_j)}h$. We take adjoints, use the support of $b_j$ and then Hölder's inequality to give
\begin{align*}
| \langle \nabla_x f_\beta, h \rangle |
& \leq \sum_j \sum_{k =1} \|b_j\|_{\L^q(Q_j)} \|(\nabla_x (1+\ell_j^2 L)^{-\beta})^* h_{k,j}\|_{\L^{q'}(Q_j)}.
\end{align*}
By (ii') we get
\begin{align}
\label{eq4: weak type bound Hardy}
\begin{split}
|\langle &\nabla_x f_\beta, h \rangle| \\
&\leq \sum_j \sum_{k = 1}^\infty \lambda |Q_j|^{\frac{1}{q}} \|(\ell_j \nabla_x (1+\ell_j^2 L)^{-\beta})^* h_{k,j}\|_{\L^{q'}(Q_j)}.
\end{split}
\end{align}
For $t>0$ the families $((1+t^2L)^{-1})$ and $(t\nabla_x(1+t^2L)^{-1})$ satisfy $\L^2$ off-diagonal estimates of arbitrarily large order. Now $q$ is an inner point of  the interval of resolvent bounds $(p_{-}(L) \vee 1,2)$, which by Theorem~\ref{thm: standard relation J(L) and N(L)} is the same as $(q_{-}(L) \vee 1,2)$ for gradient bounds. By interpolation (Lemma~\ref{lem: OD interpolation}) both families have $\L^q$ off-diagonal bounds of arbitrarily large order. Composition and duality yield $\L^{q'}$ off-diagonal bounds of arbitrarily large order $\gamma >0$ for $((t\nabla_x (1+t^2 L)^{-\beta})^*)$. Consequently, we have 
\begin{align*}
\|(\ell_j \nabla_x (1+\ell_j^2 L)^{-\beta})^*h_{k,j}\|_{\L^{q'}(Q_j)}
&\lesssim 2^{-k \gamma} \|h\|_{\L^{q'}(C_k(Q_j))} \\
&\lesssim  2^{-k \gamma} |2^k Q_j|^{\frac{1}{q'}} (\Max(|h|^{q'})(x))^{\frac{1}{q'}},
\end{align*}
where $x \in Q_j$ is arbitrary. We take $\gamma > \nicefrac{n}{q'}$ so that when substituting this estimate back into \eqref{eq4: weak type bound Hardy}, we obtain a finite sum in $k$:
\begin{align*}
|  \langle \nabla_x f_\beta, h  \rangle |
\lesssim \lambda \sum_j |Q_j| \inf_{x \in Q_j} (\Max(|h|^{q'})(x))^{\frac{1}{q'}}.
\end{align*}
We average in $x \in Q_j$, take into account the finite overlap of the $Q_j$ and apply Kolmogorov's Lemma, in order to conclude that
\begin{align}
\label{eq: weak type bound Hardy Kolmogorov}
\begin{split}
|  \langle \nabla_x f_\beta, h  \rangle |
& \lesssim \lambda \sum_j \int_{Q_j} (\Max(|h|^{q'})(x))^{\frac{1}{q'}} \, \d x \\
&\lesssim \lambda \int_{\cup_j Q_j} (\Max(|h|^{q'})(x))^{\frac{1}{q'}} \, \d x \\
& \lesssim \lambda \bigg|\bigcup_j Q_j \bigg|^{\frac{1}{q}} \|h\|_{q'}.
\end{split}
\end{align}
We recall the definition of $f_\beta$ from \eqref{eq3: weak type bound Hardy} and that $h$ was normalized in $\L^{q'}$. Hence we have shown the estimate
\begin{align*}
\bigg\| \nabla_x \sum_j (1+\ell_j^2 L)^{-\beta}b_j \bigg \|_q \lesssim \lambda \bigg(\sum_j |Q_j| \bigg)^{\frac{1}{q}},
\end{align*}
where $j$ runs over a finite subset of $J$. Property (iv) of the Calder\'on--Zygmund decomposition implies that $\sum_{j \in J} (1+\ell_j^2 L)^{-\beta}b_j $ converges in $\Wdot^{1,q}$ and that its norm is under control by $\lambda^{1-p/q} \|\nabla_x u\|_p^{p/q}$. By definition in \eqref{eq2: weak type bound Hardy}, the series $\wt{g}$ is a finite sum in $\beta$ over series of this type. Hence, it converges in $\Wdot^{1,q}$ as required and is bounded by
\begin{align}
\label{eq: weak type bound Hardy estimate of wtg}
\|\nabla_x \wt{g}\|_q \lesssim \lambda^{1-\frac{p}{q}} \|\nabla_x u\|_p^{\frac{p}{q}}.
\end{align}

\medskip

\noindent \emph{Step 5: Bound of $A_2$.} We argue as in Step~3 and use \eqref{eq: weak type bound Hardy estimate of wtg} instead of (ii') to give
\begin{align*}
|A_2| \lesssim \lambda^{-q} \|\nabla_x \wt{g}\|_q^q \lesssim \lambda^{-p} \|\nabla_x u\|_p^p.
\end{align*}

\medskip

\noindent \emph{Step 6: Preparation of the bound for $A_3$.} By Markov's inequality and the boundedness of $v \mapsto \S_{\psi,L}({L^{1/2}}v)$ from $\Wdot^{1,q}$ into $\L^q$, we have
\begin{align*}
\Big| \Big \{ x \in \R^n : S_{\psi,L} ({L^{1/2}} v)(x) > \lambda \Big \} \Big| \lesssim \lambda^{-q} \|\nabla_x v\|_q^q \quad (v \in \Wdot^{1,q}).
\end{align*}
In particular, the measure of the set on the left tends to $0$ as $v$ tends to $0$ in $\Wdot^{1,q}$. Since the series $b$ converges in $\Wdot^{1,q}$, this argument shows that it suffices to derive the desirable bound $\lambda^{-p} \|\nabla_x u\|_p^p$ for the measure of
\begin{align*}
\wt{A}_3 & \coloneqq \bigg \{ x \in \R^n : S_{\psi,L} ( {L^{1/2}}  \, \wt{b})(x) > \frac{\lambda}{2} \Big \},
\end{align*}
where $\wt{b} \coloneqq \sum_j \varphi(\ell_j^2 L) b_j$ and $j$ runs over a finite subset of $J$.   Again, this reduction bears the advantage that $\wt{b}$ is contained in $\W^{1,2} = \dom({L^{1/2}})$ and hence we can properly work with the functional calculus of $L$. In fact, such type of reduction is necessary since $p$ may lie outside of $\cJ(L)$ and therefore there is no hope for reasonable functional calculus bounds for $L$ on $\L^p$.

First, we can split off $E \coloneqq \bigcup_{j \in J} 6Q_j$ since its measure is under control by property (iii) of the Calder\'on--Zygmund decomposition.
Next, by Markov's inequality and the definition of $S_{\psi, L}$, the measure of the remaining set is at most
\begin{align*}
|\wt{A}_3 \setminus E|
&\leq 4\lambda^{-2} \int_{\wt{A}_3 \setminus E} | (S_{\psi, L} ({L^{1/2}} \, \wt{b}))(x)|^2 \, \d x\\
&\leq 4\lambda^{-2} \iint_{\reu} |(\psi(t^2 L) {L^{1/2}}  \, \wt{b})(y) |^2 \frac{|B(y,t) \setminus E|}{t^n} \, \frac{\d t \d y}{t}. 
\end{align*}
The set $B(y,t) \setminus E$ has of course measure controlled by $t^n$ but if $y$ is contained in the cube $4Q_j$, then this set is empty for all $t< \ell_j$. Hence, introducing the `local' and `global' parts
\begin{align}
\label{eq: weak type bound Hardy local and global}
\begin{split}
f_{\loc}(t,y) &\coloneqq \sum_j \ind_{4Q_j}(y) \ind_{(\ell_j, \infty)}(t) \Big(\psi(t^2 L) {L^{1/2}}  \varphi(\ell_j^2 L) b_j \Big)(y), \\
f_{\glob}(t,y) &\coloneqq \sum_j \ind_{{}^c(4Q_j)} (y) \Big(\psi(t^2 L) {L^{1/2}}  \varphi(\ell_j^2 L) b_j \Big)(y),
\end{split}
\end{align}
we obtain
\begin{align*}
|\wt{A}_3 \setminus E| \lesssim \lambda^{-2} \iint_{\reu} |f_{\loc}(t,y)|^2 + |f_{\glob}(t,y)|^2 \, \frac{\d t \d y}{t}
\end{align*}
and we are left with bounding the two integrals on the right by generic multiples of $\lambda^{2-p} \|\nabla_x u \|_p^p$.

\medskip

\noindent \emph{Step 7: The local part.} Let $h \in \L^2(\reu, \frac{\d t \d x}{t})$ and let $\langle \cdot \,, \cdot \rangle$ be the duality pairing on that space. By the Cauchy--Schwarz inequality we first find
\begin{align*}
|\langle f_{\loc}, h \rangle| \leq \sum_j I_j \bigg(\int_{4Q_j} \int_{ \ell_j}^\infty |h(t,y)|^2 \, \frac{\d t \d y}{t} \bigg)^{1/2} \\
\end{align*}
and then, generously bounding the second integral by a maximal function in $x$, that
\begin{align}
\label{eq5: weak type bound Hardy}
|\langle f_{\loc}, h \rangle| \leq \sum_j  I_j |4Q_j|^{1/2} \inf_{x \in Q_j} (\Max(H^2)(x))^{1/2},
\end{align}
where
\begin{align}
\label{eq: weak type bound Hardy Ij and H}
\begin{split}
I_j &\coloneqq \bigg(\int_{\ell_j}^\infty \int_{4Q_j}  |{L^{1/2}} \psi(t^2 L) \varphi(\ell_j^2 L) b_j (y)|^2 \, \frac{\d y \d t}{t} \bigg)^{1/2}, \\
H(y) &\coloneqq \bigg(\int_0^\infty |h(t,y)|^2 \, \frac{\d t}{t} \bigg)^{1/2}.
\end{split}
\end{align}

At this stage of the proof we introduce a fixed exponent $\varrho \in (p_-(L) \vee 1, q)$ and take the parameter $\alpha$ in \eqref{eq: auxiliary function for weak type bound Hardy} large enough to grant that $(t {L^{1/2}}\psi(t^2L))_{t>0}$ is $\L^\varrho - \L^2$-bounded. This is possible by Lemma~\ref{lem: Lp - L2 off diagonal bounds for J(L)}.(i) since $\varrho$ is not the lower endpoint of $\cJ(L)$ and we can expand
\begin{align*}
t {L^{1/2}}\psi(t^2L) = ((1+t^2L)^{-2} - (1+t^2L)^{-3})^\alpha
\end{align*}
in terms of resolvents of power at least $2 \alpha$. By interpolation with the $\L^2$-bound we then have of course $\L^r - \L^2$-boundedness for all $r \in [\varrho,2]$. Since $\varphi$ from \eqref{eq1: weak type bound Hardy} is bounded, we obtain from the functional calculus on $\L^2$ that
\begin{align}
\label{eq6: weak type bound Hardy}
\|{L^{1/2}} \psi(t^2 L) \varphi(\ell_j^2 L) f\|_2 \lesssim t^{-1} t^{\frac{n}{2}-\frac{n}{r}} \|f\|_r \quad (f \in \L^r \cap \L^2).
\end{align}

In this step we use the above estimate with $r=q$ and $f=b_j$ to bound $I_j$. As we have $n-\nicefrac{2n}{q} \leq 0$, integration in $t$ leads us to
\begin{align*}
I_j \lesssim \|b_j\|_q \bigg(\int_{\ell_j}^\infty t^{n - \frac{2n}{q} -2} \, \frac{\d t}{t}\bigg)^{\frac{1}{2}} \lesssim \ell_j^{\frac{n}{2}-\frac{n}{q} -1} \|b_j\|_q \lesssim \lambda |Q_j|^{\frac{1}{2}},
\end{align*}
where the final step uses (ii'). Going back to \eqref{eq5: weak type bound Hardy}, we have established the bound
\begin{align}
\label{eq7: weak type bound Hardy}
|\langle f_{\loc}, h \rangle|
\lesssim \lambda \sum_j |Q_j|  \inf_{x \in Q_j} (\Max(H^2)(x))^{\frac{1}{2}},
\end{align} 
so that we can bring into play Kolmogorov's lemma as in \eqref{eq: weak type bound Hardy Kolmogorov} and then use property (iv) to conclude 
\begin{align*}
|\langle f_{\loc}, h \rangle|  
\lesssim \lambda \bigg| \bigcup_j Q_j \bigg|^{\frac{1}{2}} \|H^2\|_1^{\frac{1}{2}} 
\lesssim \lambda^{1-\frac{p}{2}} \|\nabla_x u\|_p^{\frac{p}{2}} \|h\|_{\L^2(\frac{\d t \d x}{t})}.
\end{align*}
Since $h$ was arbitrary, we have proved the bound that was required at the end of Step~6:
\begin{align*}
\iint_{\reu} |f_{\loc}(t,y)|^2 \, \frac{\d t \d y}{t} 
\lesssim \lambda^{2-p} \|\nabla_x u\|_p^p.
\end{align*}

\medskip

\noindent \emph{Step 8: The global part.}
We use the same duality argument as in Step~7 except that for $f_{\glob}$ we will have to work on the ${}^c(4Q_j)$, which we split into annuli $C_k(Q_j)$, $k \geq 2$. In this manner, our substitute for \eqref{eq5: weak type bound Hardy} becomes
\begin{align}
\label{eq8: weak type bound Hardy}
|\langle f_{\glob}, h \rangle| 
\leq \sum_j \sum_{k \geq 2}  I_{j,k} |2^{k+1} Q_j|^{\frac{1}{2}} \inf_{x \in Q_j} (\Max(H^2)(x))^{\frac{1}{2}},
\end{align}
where $H$ is still as in \eqref{eq: weak type bound Hardy Ij and H} and
\begin{align*}
I_{j,k} \coloneqq \bigg(\int_{0}^\infty \int_{C_k(Q_j)}  |{L^{1/2}} \psi(t^2 L) \varphi(\ell_j^2 L) b_j (y)|^2 \, \frac{\d y \d t}{t} \bigg)^{\frac{1}{2}}.
\end{align*}

From the definitions in \eqref{eq: auxiliary function for weak type bound Hardy} and \eqref{eq1: weak type bound Hardy} we see that $z \mapsto \sqrt{z}\psi(z)$ and $\varphi$ are of class $\Psi_\alpha^{2\alpha}$ and $\Psi_\alpha^0$, respectively. Lemma~\ref{lem: functional calculus bounds from J(L) abstract}.(i) yields for all $f \in \L^2$ with support in $Q_j$ that
\begin{align*}
\|{L^{1/2}} \psi(t^2L) \varphi(\ell_j^2L) f\|_{\L^2(C_k(Q_j))}
\lesssim t^{-1} \bigg(\frac{2^k \ell_j}{t}\bigg)^{- 2\alpha} \|f\|_{\L^2(Q_j)}.
\end{align*}
For fixed $j$, $k$, $t$, we interpolate this bound with \eqref{eq6: weak type bound Hardy} for $r=\varrho$ by means of the Riesz--Thorin theorem. This results in
\begin{align}
\label{eq9: weak type bound Hardy}
\begin{split}
\|{L^{1/2}} \psi(t^2L) &\varphi(\ell_j^2L)f\|_{\L^2(C_k(Q_j))}\\
&\lesssim t^{-1+\frac{n}{2}-\frac{n}{q}} \bigg(\frac{2^k \ell_j}{t} \bigg)^{-2\theta \alpha} \|f\|_{\L^q(Q_j)},
\end{split}
\end{align}
where $\theta \in (0,1)$ is such that $q = [\varrho,2]_\theta$. In exactly the same manner we can interpolate the assertion of Lemma~\ref{lem: functional calculus bounds from J(L) abstract}.(ii) with \eqref{eq6: weak type bound Hardy} in order to obtain
\begin{align}
\label{eq10: weak type bound Hardy}
\|{L^{1/2}} \psi(t^2L) \varphi(\ell_j^2L) f\|_{\L^2(C_k(Q_j))}
\lesssim t^{-1+\frac{n}{2}-\frac{n}{q}}2^{- 2\theta \alpha k} \|f\|_{\L^q(Q_j)},
\end{align}
provided that $t \geq \ell_j$. 

Now we come back to $I_{j,k}$, split the outer integral at $t=\ell_j$ and use \eqref{eq9: weak type bound Hardy} and \eqref{eq10: weak type bound Hardy} with $f=b_j$ to give
\begin{align*}
I_{j,k}^2
&\lesssim 2^{-4\theta \alpha k} \ell_j^{-4 \theta \alpha} \|b_j\|_q^2 \int_0^{\ell_j} t^{-2+n-\frac{2n}{q} + 4 \theta \alpha}  \, \frac{\d t}{t} \\
&\quad + 2^{-4\theta \alpha k} \|b_j\|_q^2 \int_{\ell_j}^\infty t^{-2+n-\frac{2n}{q}} \, \frac{\d t}{t}.
\end{align*}
There is no issue with convergence of the second integral since we have $q \leq 2$. We pick $\alpha$ large in dependence of $n, q, \theta$ in order to grant convergence of the first integral and get
\begin{align*}
I_{j,k}^2 \leq 2^{-4\theta \alpha k} \ell_j^{-2+n-\frac{2n}{q}} \|b_j\|_q^2 \lesssim \lambda^2 2^{-4\theta \alpha k} |Q_j|,
\end{align*}
where the final step follows from (ii'). We pick $\alpha \geq \nicefrac{n}{(4 \theta)}$ so that when finally going back to \eqref{eq8: weak type bound Hardy}, we find a convergent geometric series in $k$ and obtain
\begin{align*}
|\langle f_{\glob}, h \rangle| \lesssim \lambda \sum_j |Q_j|^{\frac{1}{2}} \inf_{x \in Q_j} (\Max(H^2)(x))^{\frac{1}{2}}.
\end{align*}
At this point, the right-hand side is the same as in the treatment of the local part. We obtain the required bound for the global part by repeating the argument following \eqref{eq7: weak type bound Hardy}. This concludes the proof of Lemma~\ref{lem: weak type bound Hardy}.
\subsection*{Part 4: \texorpdfstring{Injection of $\boldsymbol{L}$-adapted spaces into classical spaces for $\boldsymbol{p \leq 2}$}{Injection of L-adapted spaces into classical spaces for p<2}}

In this section we establish the continuous inclusions
\begin{align}
\label{eq: HpL in Hp}
\IH_L^p &\subseteq a^{-1} (\H^p \cap \L^2) \\
\label{eq: H1pL in H1p}
\IH_L^{1,p} & \subseteq \Hdot^{1,p} \cap \L^2
\end{align}
in the range $1_* < p \leq 2$.

The main observation is the following inclusion for $DB$-adapted spaces. The result appears already in \cite[Sec.~4.4]{AusSta} but for convenience we include a proof.

\begin{lem}
\label{lem: HpDB in Hp}
If $p \in (1_*,2]$, then $\IH_{DB}^p \subseteq \IH_D^p$ and the inclusion is continuous for the $p$-quasinorms.
\end{lem}

\begin{proof}
The claim holds for $p=2$, see \eqref{eq: Hardy identification for p=2}, and by Remark~\ref{rem: identification regions are intervals} the set of exponents for which the claim holds is an interval. Hence, we only have to treat the case $p \leq 1$.

We use the molecular decomposition for $\IH_{DB}^p$ (Theorem~\ref{thm: Molecular decomposition adapted Hardy}) for some admissible $M$ and $\eps = 1$. It suffices to check that there is a constant $c$ such that $\|m\|_{\IH_D^p} \leq c$ for every $(\IH_{DB}^p ,1, M)$-molecule. Writing $m = D (B(DB)^{-1}m)$, we see that $m$ is a generic multiple of an $(\IH_D^p, 1, 1)$-molecule. The required bound follows from Corollary~\ref{cor: HDp molecules}.
\end{proof}

Now, we can use Figure~\ref{fig: diagram} and the identification of $D$-adapted spaces in Figure~\ref{fig: identification D-adapted} as follows to complete Part~4. Moving from the third to the fourth row, we obtain for $f \in \IH_L^p \cap \ran(a^{-1} \div_x)$ that
\begin{align*} 
\|af\|_{\H^p} 
\lesssim \bigg\| \begin{bmatrix} af \\ 0 \end{bmatrix} \bigg\|_{\IH_{DB}^p}
\lesssim \|f\|_{\IH_L^p} .
\end{align*}
The bound extends to $f \in \IH_L^p$ by density, which gives \eqref{eq: HpL in Hp}. Likewise, moving from the third to the second row, we get
\begin{align*}
\|f\|_{\Hdot^{1,p}} 
= \| \nabla_x f \|_{\H^p} 
\lesssim \bigg\| \begin{bmatrix} 0 \\ \nabla_x f \end{bmatrix} \bigg\|_{\IH^p_{DB}} 
\lesssim \|f\|_{\IH^{1,p}_{L}}, 
\end{align*}
first for $f \in \IH^{1,p}_{L} \cap \dom({L^{1/2}})$ and then for all $f \in \IH^{1,p}_{L}$, which gives \eqref{eq: H1pL in H1p}. 

Going one step further to the first row gives an additional Riesz transform bound\index{Riesz transform!$\IH_L^p-\H^p$-bound}, which is of independent interest. It extends \cite[Prop.~5.6]{HMMc} beyond semigroup generators.

\begin{prop}
\label{prop: Riesz transform on Hardy spaces}
If $p \in (1_*,2]$, then 
\begin{align*}
\|\nabla_x L^{-1/2} f \|_{\H^p} \lesssim \|f\|_{\IH_L^p} \quad (f \in \IH_L^p \cap \ran({L^{1/2}})).
\end{align*}
\end{prop}
\subsection*{Part 5: \texorpdfstring{Injection of classical spaces into $\boldsymbol{L}$-adapted spaces for $\boldsymbol{p \leq 1}$}{Injection of classical spaces into L-adapted spaces for p<1}}
\label{subsec: classical into HpL for p<1}

We complement the previous section by proving the reverse continuous inclusions
\begin{align}
\label{eq: Hp in HpL}
a^{-1} (\H^p \cap \L^2) &\subseteq \IH_L^p \qquad (p_{-}(L) < p \leq 1)
\intertext{and}
\label{eq: H1p in H1pL}
\Hdot^{1,p} \cap \L^2 &\subseteq \IH_L^{1,p} \qquad ((p_{-}(L)_* \vee 1_*)< p \leq 1),
\end{align}
if these intervals of exponents are non-empty.

The strategy is the same for both inclusions and relies on the atomic decompositions. We use the auxiliary function $\psi(z) \coloneqq z^{\alpha}(1+z)^{-2\alpha}$, where $\alpha \in \IN$ will be chosen large later on, and introduce the square functions
\begin{align}
\label{eq: SF for Hp in HpL}
S_{\psi,L}^{(0)} f(x) &\coloneqq \bigg(\iint_{|x-y|<t} |\psi(t^2L) f(y)|^2 \, \frac{\d t \d y}{t^{1+n}} \bigg)^{1/2},\\
\label{eq: SF for H1p in H1pL}
S_{\psi,L}^{(1)} f(x) &\coloneqq \bigg(\iint_{|x-y|<t} |t^{-1} \psi(t^2L) f(y)|^2 \, \frac{\d t \d y}{t^{1+n}} \bigg)^{1/2}.
\end{align}
Then $\|S_{\psi,L}^{(0)}(\cdot)\|_p$ and $\|S_{\psi,L}^{(1)}(\cdot)\|_p$ are equivalent norms on $\IH_L^p$ and $\IH_L^{1,p}$ provided that we take at least $\alpha > \nicefrac{n}{(2p)} - \nicefrac{n}{4}$. 

We shall establish the following bounds.

\begin{lem}
\label{lem: atomic H1pL bound}
Let $p \in (p_{-}(L)_* \vee 1_*, 1]$ and $\alpha$ sufficiently large depending on $n, p, p_-(L)$. For all $\L^2$-atoms $m$ for $\Hdot^{1,p}$ it follows that
\begin{align*}
\|S_{\psi,L}^{(1)}(m)\|_p \lesssim 1.
\end{align*}
\end{lem}

\begin{lem}
\label{lem: atomic HpL bound}
Let $p \in (p_{-}(L), 1]$ and $\alpha$ sufficiently large depending on $n, p, p_-(L)$. For all $\L^2$-atoms $m$ for $\H^{p}$ it follows that
\begin{align*}
\|S_{\psi,L}^{(0)}(a^{-1}m)\|_p \lesssim 1.
\end{align*}
\end{lem}

Let us take these estimates for granted and complete the objective of this part first. Given $f \in \L^2$ such that $a f \in \H^p$, we write the latter as an $\L^2$-convergent atomic decomposition $a f = \sum_i \lambda_i m_i$ with $\|(\lambda_i)\|_{\ell^p} \lesssim \|af\|_{\H^p}$. We use Fatou's lemma as in the proof of Proposition~\ref{prop: Hp in HpL for 1<p<2} to obtain
\begin{align*}
S_{\psi,L}^{(0)}f(x) \leq \sum_{i} |\lambda_i| S_{\psi,L}^{(0)} (a^{-1}m_i)(x) \quad (x \in \R^n)
\end{align*}
and we conclude by Lemma~\ref{lem: atomic HpL bound} and as $p \leq 1$,
\begin{align*}
\|S_{\psi,L}^{(0)} (f)\|_p^p 
\leq \sum_{i} |\lambda_i|^p \|S_{\psi,L}^{(0)} (a^{-1}m_i)\|_p^p 
\lesssim \sum_{i} |\lambda_i|^p 
\lesssim \|af\|_{\H^p}^p.
\end{align*}
The left-hand side is equivalent to $\|f\|_{\IH_L^p}^p$ and \eqref{eq: Hp in HpL} follows.

As for \eqref{eq: H1p in H1pL}, it suffices to prove $\|u\|_{\IH_L^{1,p}} \lesssim \|\nabla_x u\|_p$ for all $u \in \Hdot^{1,p} \cap \W^{1,2}$. Indeed, since $\cZ$ is dense in $\Hdot^{1,p} \cap \L^2$, this is yet another application of the Fatou argument above. Now, we can take a $\Wdot^{1,2}$-convergent atomic decomposition $u = \sum_i \sum_j \lambda_i^j m_i^j$ as in Proposition~\ref{prop: H1p atomic decomposition}. By the solution of the Kato problem we have $\L^2$-convergence of 
\begin{align*}
\psi(t^2 L) u  = \sum_{i} \sum_{j} \lambda_i^j L^{-1/2} \psi(t^2 L) L^{1/2} m_i^j
\end{align*}
and the same argument as before applies.

\begin{proof}[Proof of Lemma~\ref{lem: atomic H1pL bound}]
Let $m$ be an $\L^2$-atom for $\Hdot^{1,p}$ associated with a cube $Q$ of sidelength $\ell$ as in Definition~\ref{def: H1p atom}.

We begin with a local bound. By the solution of the Kato problem we have $m \in \dom({L^{1/2}})$. It follows that
\begin{align*}
S_{\psi,L}^{(1)} m(x) 
&= \bigg(\iint_{|x-y|<t} |\varphi(t^2L) {L^{1/2}} m(y)|^2 \, \frac{\d t \d y}{t^{1+n}} \bigg)^{1/2} \\
& \eqqcolon S_{\varphi, L} ({L^{1/2}}m)(x) \quad (x \in \R^n),
\end{align*} 
where $\varphi(z) \coloneqq z^{\alpha-1/2}(1+z)^{-2\alpha}$. H\"older's inequality and the $\L^2$-bound for the square function with $\varphi$ (McIntosh's theorem) yield
\begin{align} 
\label{eq1: atomic H1pL bound}
\begin{split}
\|S_{\psi,L}^{(1)}  (m)\|_{\L^p(4Q)} 
&\leq |4Q|^{\frac{1}{p}-\frac{1}{2}} \|S_{\psi,L}^{(1)} (m)\|_{\L^2(4Q)} \\
&\lesssim |Q|^{\frac{1}{p}-\frac{1}{2}} \|{L^{1/2}} m\|_2\\
&\simeq |Q|^{\frac{1}{p}-\frac{1}{2}} \|\nabla_x m\|_2\\
&\leq 1.
\end{split}
\end{align}

In preparation of the global bound, we pick some $q \in (p_{-}(L), p^*) \cap (1,2]$. This is possible by the assumption on $p$. We also take $\alpha$ large enough in dependence of $q$ and $p_{-}(L)$ in order to have $\L^q - \L^2$ off-diagonal estimates of arbitrarily large order for $(\psi(t^2L))_{t>0}$ at our disposal. This is possible due to Lemma~\ref{lem: Lp - L2 off diagonal bounds for J(L)}.(i) since we can expand
\begin{align}
\label{eq2: atomic H1pL bound}
\begin{split}
\psi(z) 
= z^\alpha(1+z)^{-2\alpha} 
= ((1+z)^{-1}-(1+z)^{-2})^\alpha.
\end{split}
\end{align}
Consequently, we have for all $x \in \R^n$ the estimate
\begin{align}
\label{eq3: atomic H1pL bound}
\begin{split}
\|\psi(t^2 L) m\|_{\L^2(B(x,t))} 
&\lesssim t^{\frac{n}{2}-\frac{n}{q}} \bigg(1+ \frac{\dist(B(x,t),Q)}{t}\bigg)^{-\gamma} \|m\|_q \\
&\simeq t^{\frac{n}{2}-\frac{n}{q}} \bigg(1+ \frac{\dist(x,Q)}{t}\bigg)^{-\gamma} \|m\|_q,
\end{split}
\end{align}
where $\gamma > 0$ is at our disposal and the second step uses $\dist(B(x,t),Q) \geq \nicefrac{\dist(x,Q)}{2}$ for $t\leq \nicefrac{\dist(x,Q)}{2}$ and $2 \geq \nicefrac{\d(x,Q)}{t}$ for $t\geq \nicefrac{\dist(x,Q)}{2}$. Squaring and integrating this bound with respect to $\nicefrac{\d t}{t^{n+3}}$ gives
\begin{align*}
S_{\psi,L}^{(1)}  m(x) 
&\lesssim \bigg(\int_0^\infty t^{-\frac{2n}{q} - 2} \bigg(1+ \frac{\dist(x,Q)}{t}\bigg)^{-2\gamma} \, \frac{\d t}{t} \bigg)^{1/2} \|m\|_q \\
&\simeq \d(x,Q)^{-\frac{n}{q}-1} \|m\|_q,
\end{align*}
where the last step follows by a change of variable $t = s \dist(x,Q)$ and we have taken $2 \gamma > \nicefrac{2n}{q}+2$ in order to have a finite integral in $s$. Thus,
\begin{align*}
\|S_{\psi,L}^{(1)} (m)\|_{\L^p({}^c (4Q))}
 &\lesssim \bigg(\int_{{}^c (4Q)} \d(x,Q)^{-\frac{np}{q}-p} \, \d x \bigg)^{\frac{1}{p}} \|m\|_q \\
 &\lesssim \ell^{\frac{n}{p} - \frac{n}{q} - 1} \|m\|_q,
\end{align*}
where we have used $\nicefrac{np}{q} + p > n$ to calculate the integral in $x$. Since $m$ is supported in $Q$, we obtain from H\"older's and Poincar\'e's inequality that
\begin{align*}
\|S_{\psi,L}^{(1)}  (m)\|_{\L^p({}^c (4Q))} 
\lesssim \ell^{\frac{n}{p}-\frac{n}{2}-1} \|m\|_{\L^2(Q)}
\lesssim \ell^{\frac{n}{p}-\frac{n}{2}}  \|\nabla_x m\|_{\L^2(Q)}
\leq 1,
\end{align*}
which is the required global bound.
\end{proof}

\begin{proof}[Proof of Lemma~\ref{lem: atomic HpL bound}]
Let $m$ be an $\L^2$-atom for $\H^{p}$ associated with a cube $Q$ of sidelength $\ell$, see Definition~\ref{def: Hp atom}.

As before, the local bound $\|S_{\psi,L}^{(0)}  (a^{-1}m)\|_{\L^p(4Q)} \lesssim 1$ follows from H\"older's inequality and the $\L^2$-bound for the square function. 

To prepare the global bound, we pick exponents $p_-(L) < s < r < q < p$. The resolvents of $L$ are $a^{-1}\H^s$-bounded and also $\L^\varrho - \L^2$-bounded for some $\varrho < 2$ thanks to Lemmata~\ref{lem: Lp L2 between Sobolev conjugates n geq 2} and~\ref{lem: Lp L2 between Sobolev conjugates n = 1}. Keeping in mind the expansion \eqref{eq2: atomic H1pL bound}, we take $\alpha$ large and conclude from Lemma~\ref{lem: extra} that $(\psi(t^2L))_{t>0}$ is $a^{-1}\H^r - \L^2$-bounded. Together with the usual $\L^2$ off-diagonal estimates we obtain for all $x \in \R^n$ that
\begin{align*}
\|\psi(t^2L) &(a^{-1}m)\|_{\L^2(B(x,t))} \\
&= \|\psi(t^2L) (a^{-1}m)\|_{\L^2(B(x,t))}^{1-\theta} \|\psi(t^2L) (a^{-1}m)\|_{\L^2(B(x,t))}^{\theta} \\
&\lesssim \bigg(t^{\frac{n}{2}-\frac{n}{r}} \|m\|_{\H^r}\bigg)^{1-\theta} \bigg( \Big(1+\frac{\dist(B(x,t),Q)}{t} \Big)^{-\gamma} \|m\|_2 \bigg)^\theta,
\end{align*}
where $\theta \in (0,1)$ and $\gamma > 0$ are still at our disposal. Since $|Q|^{1/p-1/r}m$ is an $\L^2$-atom for $\H^r$, we have $\|m\|_{\H^r} \lesssim |Q|^{1/r-1/p}$. Picking $\theta$ such that $\nicefrac{(1-\theta)}{r}+\nicefrac{\theta}{2} = \nicefrac{1}{q}$, we obtain
\begin{align*}
\|\psi(t^2L) (a^{-1}m)\|_{\L^2(B(x,t))}
\lesssim t^{\frac{n}{2}-\frac{n}{q}} \Big(1+\frac{\dist(B(x,t),Q)}{t} \Big)^{-\gamma \theta} |Q|^{\frac{1}{q}-\frac{1}{p}}.
\end{align*}
This estimate is of the exact same type as \eqref{eq3: atomic H1pL bound} and we can repeat the previous proof from thereon. Indeed, we integrate the square with respect to $\nicefrac{\d t}{t^{1+n}}$ to obtain
\begin{align*}
S_{\psi,L}^{(0)} (a^{-1}m)(x) 
\lesssim \dist(x,Q)^{-\frac{n}{q}} \ell^{\frac{n}{q}-\frac{n}{p}},
\end{align*}
and then the required global bound
\begin{align*}
\|S_{\psi,L}^{(0)} (a^{-1}m)\|_{\L^p({}^c(4Q))} 
\lesssim \ell^{\frac{n}{p}-\frac{n}{q}} \ell^{\frac{n}{q}-\frac{n}{p}} = 1,
\end{align*}
follows since $\nicefrac{np}{q} > n$.
\end{proof}
\subsection*{Part 6: \texorpdfstring{$\boldsymbol{h_{-}(L) \leq p_{-}(L)}$}{h-(L) < p-(L)}}

Let $p \in (p_{-}(L),2]$. We have to prove that $a^{-1}(\H^p \cap \L^2) = \IH_L^p$ with equivalent $p$-quasinorms.

The inclusion `$\subseteq$' was obtained in Part~5 for $p \in (p_{-}(L),1]$ and in Proposition~\ref{prop: Hp in HpL for 1<p<2} for $p \in (p_{-}(L) \vee 1,2]$. The converse was obtained in Part~4 in the range $p \in (1_*,2]$.
\subsection*{Part 7: \texorpdfstring{$\boldsymbol{h^1_{-}(L) \leq (p_{-}(L)_* \vee 1_*)}$}{Upper bound for h1-}}

Let $p \in (p_{-}(L)_* \vee 1_*,2]$. We have to prove $\Hdot^{1,p} \cap \L^2 = \IH_L^{1,p}$ with equivalent $p$-quasinorms.

We have obtained `$\subseteq$' in Part~5 for $p \in (p_{-}(L)_* \vee 1_*,1]$ and in Proposition~\ref{prop: H1p in H1pL for 1<p<2} for $p \in (p_{-}(L)_* \vee 1,2]$. The converse follows again from Part~4.
\subsection*{Part 8: \texorpdfstring{$\boldsymbol{h_{+}(L) \geq p_{+}(L)}$}{h+(L) < p+(L)}}

Let $p \in (2,p_+(L))$. In Part~2 we have obtained $\L^p \cap \L^2 \subseteq \IH_L^p$ with continuous inclusion for the $p$-norms. It remains to establish the opposite inclusion and this will follow by duality.  

To this end, we recall from Section~\ref{subsec: adjoints} that $L^*$ is an operator in the same class as $\tL$ and similar to an operator $L^\sharp$ in the same class as $L$ under conjugation with $a^*$. By duality and similarity we have $p' \in (p_{-}(L^\sharp) \vee 1,2)$. Replacing systematically $L$ with $L^\sharp$, the result of Part~6 entails $\IH_{L^\sharp}^{p'} = \L^{p'} \cap \L^2$ with equivalent $p'$-norms and from Figure~\ref{fig: diagram} we can read off
\begin{align*}
\IH_{L^*}^{p'} = a^* \IH_{L^\sharp}^{p'} = \L^{p'} \cap \L^2.
\end{align*}
Given $f \in \IH_L^p$, we use Proposition~\ref{prop: Hardy duality} for second-order operators to give
\begin{align*}
|\langle f, g \rangle| \lesssim \|f\|_{\IH_L^p} \|g\|_{\IH_{L^*}^{p'}} \simeq \|f\|_{\IH_L^p} \|g\|_{p'} \quad (g \in \L^{p'} \cap \L^2).
\end{align*}
We conclude $f \in \L^p \cap \L^2$ along with $\|f\|_{p} \lesssim  \|f\|_{\IH_L^p}$ as required.
\subsection*{Part 9: \texorpdfstring{$\boldsymbol{h_{+}^1(L) \geq q_{+}(L)}$}{h1+(L) > q+(L)}}

We have to show that $\Wdot^{1,p} \cap \L^2 = \IH_L^{1,p}$ with equivalent $p$-norms for $p \in (2, q_+(L))$. In fact, we shall establish continuous inclusions for the $p$-Hardy norms
\begin{align}
\label{eq: H1pL in H1p for p>2}
\Wdot^{1,p} \cap \L^2 &\supseteq \IH_L^{1,p} \qquad (2<p<q_+(L))
\intertext{and}
\label{eq: H1p in H1pL for p>2}
\Wdot^{1,p} \cap \L^2 &\subseteq \IH_L^{1,p} \qquad (2<p<p_+(L)),
\end{align}
which is a more general result since by Theorem~\ref{thm: standard relation J(L) and N(L)} we have $p_+(L) \geq q_+(L)^*$.

In the following let $p \in (2, p_+(L))$. Part~8 implies $p< h_{+}(L)$. Hence, we can identify $\IH_L^p = \L^p \cap \L^2$ and the ubiquitous Figure~\ref{fig: diagram} tells us that
\begin{align}
\label{eq1: upper bound q+}
\|f\|_{\IH_L^{1,p}} \simeq \|L^{1/2} f\|_{\IH_L^{p}} \simeq \|L^{1/2} f \|_p \quad (f \in \IH^{1,p}_L \cap \dom(L^{1/2})).
\end{align}

\medskip

\noindent \emph{Proof of \eqref{eq: H1pL in H1p for p>2}}. If even $p< q_+(L)$, then the Riesz transform is $\L^p$-bounded according to Theorem~\ref{thm: Riesz} and we obtain from \eqref{eq1: upper bound q+} that
\begin{align*}
\|f\|_{\IH_L^{1,p}} \gtrsim \|\nabla_x f\|_p \quad (f \in \IH^{1,p}_L \cap \dom(L^{1/2})).
\end{align*}
A general $f \in \IH^{1,p}_L$ can be approximated by a sequence $(f_j) \subseteq \IH^{1,p}_L \cap \dom(L^{1/2})$ simultaneously in $\IH^{1,p}_L$ and $\L^2$, see Section~\ref{subsec: Hardy abstract}. Then $(\nabla_x f_j)$ is a Cauchy sequence in $\L^p$ whose limit coincides with $\nabla_x f$ thanks to $\L^2$-convergence of $(f_j)$. Hence, the previous estimate extends to $f$.

\medskip

\noindent \emph{Proof of \eqref{eq: H1p in H1pL for p>2}}. It suffices to establish the bound
\begin{align}
\label{eq4: upper bound q+}
\|u\|_{\IH_L^{1,p}} \lesssim \|\nabla_x u\|_p \quad (u \in \Wdot^{1,p} \cap \W^{1,2}).
\end{align} 
Indeed, a general $u \in \Wdot^{1,p} \cap \L^{2}$ can be approximated in $\Wdot^{1,p} \cap \L^2$ by a sequence $(u_j) \subseteq \cZ$ and $\L^2$-convergence suffices to infer $\|u\|_{\IH_L^{1,p}} \leq \liminf_{j\to \infty} \|u_j\|_{\IH_L^{1,p}}$, see the proof of Proposition~\ref{prop: Hp in HpL for 1<p<2}. 

We rely on a duality argument using the same notation as in Part~8. Again, we have $p' \in (p_{-}(L^\sharp) \vee 1,2)$ and we obtain from Theorem~\ref{thm: Riesz} that the Riesz transform for $L^\sharp$ is $\L^{p'}$-bounded. For any $g \in \ran(L^*) \cap \dom(L^*) \cap \L^{p'}$  it follows that
\begin{align*}
\big \langle L^{1/2}u, g \big \rangle
&= \big \langle u,  (L^*)^{1/2} g \big \rangle \\
&= \big \langle u,  L^* (L^*)^{-1/2} g \big \rangle\\
&= \big \langle \nabla_x u, d^* \nabla_x (a^*)^{-1} (L^*)^{-1/2} g \big \rangle \\
&= \big \langle d \nabla_x u, \nabla_x (L^\sharp)^{-1/2} (a^*)^{-1}g \big \rangle,
\end{align*}
where the third step is just the definition of $L^*$ and the final step uses that the similarity of operators $L^* = a^* L^\sharp (a^*)^{-1}$ carries over to the functional calculi by construction. Hölder's inequality yields
\begin{align*}
|\langle L^{1/2}u, g \rangle |
\lesssim \|  \nabla_x u\|_p \|\nabla_x \, (L^\sharp)^{-1/2} (a^*)^{-1}g \|_{p'} 
\lesssim \|  \nabla_x u\|_p \, \|g\|_{p'}.
\end{align*}
Since $g$ was taken from a dense subspace of $\L^{p'}$ (as is granted by Lemma~\ref{lem: range density in Lp} applied to $L^\sharp$ and similarity), the bound $\|L^{1/2}u\|_p \lesssim \|\nabla_x u\|_p$
follows. Now, \eqref{eq4: upper bound q+} is a consequence of \eqref{eq1: upper bound q+}.
\subsection*{Part 10: \texorpdfstring{$\boldsymbol{h_{+}^1(L) \leq q_+(L)}$}{h1+(L) < q+(L)}}

Suppose that the interval $\cH^{1}(L)$ contains some exponent $p \geq q_+(L)$. In particular, $q_+(L)$ is finite. 

Since we have $q_+(L) < p_+(L)$ by Theorem~\ref{thm: standard relation J(L) and N(L)}, we can assume $p < p_+(L)$ and by the result of Part~8 this implies $p \in \cH(L)$. Therefore, we have a commutative diagram
\begin{center}
\begin{tikzcd}[column sep=40pt]
		\IH_L^p \cap \ran({L^{1/2}}) \arrow[r, "L^{-1/2}"]	&\IH_L^{1,p} \cap \dom({L^{1/2}}) \arrow[hookrightarrow]{r} &\Wdot^{1,p} \cap \dom({L^{1/2}})\arrow[d, "\nabla_x"]\\
		\L^p \cap \ran({L^{1/2}}) \arrow[hookrightarrow]{u} \arrow[rr, "\nabla_x L^{-1/2}"]  &&\L^p \cap \L^2
		\label{fig: diagram part10},
\end{tikzcd}
\end{center}
where the mapping of $L^{-1/2}$ follows from Figure~\ref{fig: diagram} and the unlabeled arrows indicate continuous inclusions for the $p$-norms. Lemma~\ref{lem: range density in Lp} guarantees that $\L^p \cap \ran({L^{1/2}})$ is dense in $\L^p \cap \L^2$ and we conclude that the Riesz transform is $\L^p$-bounded. But then we must have $p \leq q_+(L)$ according to Theorem~\ref{thm: Riesz} and therefore $p= q_+(L)$. 

This argument has two consequences. First, $q_+(L) \in \cH^1(L)$ is possible only if the Riesz transform is $\L^{q_+(L)}$-bounded. We shall see in the next section that this is never the case. Second, $\cH^1(L)$ cannot contain exponents $p> q_+(L)$ and hence we have $h_+^1(L) \leq q_+(L)$. At this stage the proof of Theorem~\ref{thm: main result Hardy} is complete. \hfill \qedsymbol
\subsection{Consequences for square functions}
\label{subsec: consequences for square functions}

By definition of $\IH_L^p$, the identification theorem (Theorem~\ref{thm: main result Hardy}) can be reformulated in terms of $\L^p$-bounds for conical square functions of type\index{square function!conical ($S_{\psi,L}$)}
\begin{align*}
S_{\psi,L} f(x) \coloneqq S(\IQ_{\psi,L}f)(x) = \bigg(\iint_{|x-y|<t} |\psi(t^2 L) f(y)|^2 \, \frac{\d t \d y}{t^{1+n}} \bigg)^{1/2}.
\end{align*}
Here, we collect and improve these bounds with an emphasis on the decay for the auxiliary function $\psi \in \Psi_+^+$ at $|z|=0$ and $|z|=\infty$ within a sector. This will be important for the applications to boundary value problems. 

When $p\geq 2$, we will use the simple fact that the conical square functions $S$ can be controlled by the vertical square function~\index{square function!vertical ($V$)} defined for $F \in \Lloc^2(\reu)$ as
\begin{align*}
V(F)(x) \coloneqq \bigg(\int_0^\infty |F(t,x)|^2 \, \frac{\d t}{t} \bigg)^{1/2},
\end{align*} 
see for instance \cite[Prop.~2.1]{AHM} for the following lemma.

\begin{lem}
	\label{lem: S less than V}
	Let $p \in [2,\infty)$. There is a constant $c$ depending on $p$ and $n$ such that for all $F \in \Lloc^2(\reu)$,
	\begin{align*}
	\|S(F)\|_{p}  \leq c \|V(F)\|_{p}.
	\end{align*}
\end{lem}

Upper bounds for vertical square functions are provided by an abstract theorem due to Cowling--Doust--McIntosh--Yagi~\cite[Thm.~6.6]{CDMcY}.~\index{Theorem!Cowling--Doust--McIntosh--Yagi's} 
We state the quantitative version found in the textbook \cite{HNVW2}, but inspection of the original argument would yield the same dependence of the constants. We continue to write
\begin{align*}
(\IQ_{\psi,T}f)(t,x) = (\psi(t^2 T)f)(x)
\end{align*}
as in Section~\ref{sec: Hardy intro}, even though $T$ need not act on $\L^2$, and we note that up to a norming factor of $2$ the vertical square function $V(\IQ_{\psi,T}f)$ does not change if instead we use first-order scaling $(\IQ_{\psi,T}f)(t,x) = (\psi(t T)f)(x)$.

\begin{thm}[{\cite[Thm.~10.4.23]{HNVW2}}]
\label{thm: Le Merdy}
Let $p \in (1,\infty)$ and let $T$ be a sectorial operator in $\L^p(\R^n;W)$, where $W$ is a finite-dimensional Hilbert space. Suppose that $T$ has a bounded $\H^\infty$-calculus of angle $\omega \in (0,\pi)$ on $\cl{\ran(T)}$. Let $\mu \in (0, \nicefrac{(\pi-\omega)}{2})$ and choose decay parameters $\sigma, \tau > 0$. Then for all $\psi \in \Psi_\sigma^\tau(\S_{\pi-2 \mu}^+)$ and all $f \in \cl{\ran(T)}$,
\begin{align*}
\|V(\IQ_{\psi,T}f)\|_{p} \lesssim \|\psi\|_{\sigma,\tau, \mu} \|f\|_{p},
\end{align*}
where the implicit constant depends on $T$ through $M_{T,\nu}$ and $M_{T,\nu}^\infty$ for some $\nu \in (\omega, \pi-2 \mu)$.
\end{thm}

\begin{rem}
\label{rem: Le Merdy}
The numbers $M_{T,\nu}$ and $M_{T,\nu}^\infty$ correspond to resolvent and functional calculus bounds, see \eqref{eq: constant for resolvent bound} and \eqref{eq: Hoo bound}. The theorem remains true for all $f \in \L^p(\R^n; W)$ since we have $\psi(t^2T)f = 0$ if $f \in \nul(T)$ and $t>0$.
\end{rem}

With this at hand, we obtain abstract square function bounds. We largely follow the idea for second-order elliptic operators $-\div_x d \nabla _x$ in \cite{AHM}, see also \cite{AusSta, AM}, but with a more direct interpolation argument in tent spaces.

\begin{prop}
\label{prop: upper SFE abstract}
Let $T$ be a sectorial operator that satisfies the standard assumptions \eqref{eq: standard assumptions sectorial}. Let $p \in [2,\infty)$ and suppose that\index{square function!bounds for sectorial operators} 
\begin{align*}
\|f\|_{\IH^p_T} \simeq \|f\|_p \quad (f \in \cl{\ran(T)}).
\end{align*}
Let $\theta \in (0,1]$, fix an angle $\mu \in (0, \nicefrac{\theta(\pi - \omega)}{2})$ and let $\psi \in \Psi_\sigma^\tau(\S_{\pi -  2\mu}^+)$ with $\sigma, \tau > 0$. Consider the square function bound
\begin{align*}
\|\IQ_{\psi, T}f \|_{\T^{q}} \lesssim  \|f\|_{q} \quad (f \in \L^q \cap \L^2),
\end{align*} 
with an implicit constant that depends on $T$ only through \eqref{eq: standard assumptions sectorial} and the comparison constant for the $p$-norms in the assumption. Then this bound is valid provided that
\begin{align*}
q \geq 2 \quad \text{and} \quad \frac{1}{q} > \frac{1}{p} - \frac{[p,2]_\theta}{p} \frac{2\sigma}{n}.
\end{align*}
\end{prop}

\begin{proof}
We organize the proof in four steps. 

\medskip

\noindent \emph{Step~1: $\H^\infty$-calculus for the $\L^p$-realization of $T$}. It follows from Proposition~\ref{prop: FC on pre-Hardy} and the assumption on $p$ that
\begin{align}
\label{eq1: upper SFE abstract}
\|\eta(T) f\|_p \lesssim \|\eta\|_\infty \|f\|_p 
\end{align}
for all $f \in \L^p \cap \cl{\ran(T)}$ and all admissible $\eta \in \H^\infty$. 

Let $\nu \in [0, \frac{\pi-\omega}{2})$ and $\zeta \in \S_{\nu}^+$. For the special choice $\eta(z) \coloneqq (1+ \zeta^2 z)^{-1}$ the operator $\eta(T)$ acts as the identity on $\nul(T)$. Hence, the bound above extends to all $f \in \L^p \cap \L^2$, that is to say, $((1+\zeta^2 T)^{-1})_{\zeta \in \S_\nu^+}$ is $\L^p$-bounded. Hence, $T$ has an $\L^p$-realization described in Proposition~\ref{prop: consistent operator} and this is a sectorial operator in $\L^p$ of the same angle $\omega$ as $T$.

For $\eta \in \Psi_+^+$ the bound \eqref{eq1: upper SFE abstract} also remains true for general $f \in \L^p \cap \L^2$ since $\eta(T)$ vanishes on $\nul(T)$. We have $\eta(T)f = \eta(T_p)f$ since these operators are given by the same Cauchy integral. Since $\L^p \cap \L^2$ is dense in $\L^p$ it follows that $T_p$ has a bounded $\H^\infty$-calculus of angle $\omega$ on the closure of its range. 

The idea of proof is now to interpolate between two square function bounds that we have seen before: Theorem~\ref{thm: Le Merdy} for $T_p$ and Lemma~\ref{lem: Carleson bound Q function} for $T$. 

\medskip

\noindent \emph{Step 2: Definition of an interpolating family}. For $\alpha \in \IC^+ \coloneqq \{z \in \IC: \Re z > 0 \}$ we define
\begin{align}
\label{eq0: upper SFE abstract}
\psi_\alpha: \S_{\pi- 2 \mu}^+ \to \IC, \quad \psi_\alpha(z) \coloneqq \bigg(\frac{z}{1+z} \bigg)^{\alpha-\sigma} \psi(z).
\end{align}
As $\nicefrac{z}{(1+z)} = (1+z^{-1})^{-1} \in \S_{\pi-2\mu}^+$ and $\Re \alpha > 0$, we obtain
\begin{align}
\sup_{z \in \S^+_{\pi- 2 \mu}} \bigg| \bigg(\frac{z}{1+z} \bigg)^{\alpha-\sigma} \bigg| 
\lesssim \e^{(\pi-2\mu) |\Im \alpha|} (|z|^{\Re \alpha - \sigma} \wedge 1),
\end{align}
where the implicit constant is independent of $\Im \alpha$. Consequently, we have $\psi_\alpha \in \Psi_{\Re \alpha}^\tau(\S_{\pi - 2\mu}^+)$ and
\begin{align}
\label{eq2: upper SFE abstract}
\|\psi_\alpha\|_{\Re \alpha,\tau, \mu} \lesssim  \e^{(\pi-2 \mu) |\Im \alpha|} \|\psi\|_{\sigma,\tau, \mu}.
\end{align}
Combining Lemma~\ref{lem: S less than V} and Theorem~\ref{thm: Le Merdy} leads to the following bound for $q \coloneqq p$ and all $f \in \L^{q} \cap \L^2$:
\begin{align*}
\|\IQ_{\psi_\alpha, T}f\|_{\T^q} 
&= \|S(\IQ_{\psi_\alpha,T_p}f)\|_q \\
&\lesssim \|V(\IQ_{\psi_\alpha,T_p}f)\|_q \\
&\lesssim  \e^{(\pi-2\mu) |\Im \alpha|} \|\psi\|_{\sigma,\tau, \mu} \|f\|_q.
\end{align*}
The implicit constant is independent of $\psi$ and $\Im \alpha$. By McIntosh's theorem the same holds for $q = 2$ and hence for all $q \in [2,p]$ by interpolation. If, however, $\Re \alpha > \frac{n}{2[p,2]_\theta}$, then Lemma~\ref{lem: Carleson bound Q function} provides the same bound for \emph{all} $q \in [2,\infty)$, so that in total we obtain
\begin{align}
\label{eq3: upper SFE abstract}
\|\IQ_{\psi_\alpha,T}f \|_{\T^{q}} \lesssim   \e^{(\pi-2\mu) |\Im \alpha|} \|\psi\|_{\sigma,\tau, \mu} \|f\|_{\L^{q}} \quad (f \in \L^{q} \cap \L^2)
\end{align}
if $(\Re \alpha, \nicefrac{1}{q})$ belongs to the interior of the gray shaded region in Figure~\ref{fig: Stein}.

\begin{figure}[ht]
	\begin{center}
		\begin{tikzpicture}[scale=2.4]

		
		\coordinate (P0) at (1,2);
		\coordinate (P1) at (1,1.5);
		\coordinate (P2) at (4.5,2);
		\coordinate (P3) at (4.5,1.5);
		\coordinate (P4) at (4.5,0);
		\coordinate (P5) at (2.7,0);
		\coordinate (P6) at (2.7,1.5);

		\draw [thin] (1,2) -- (4.5,2); 
		\draw [thin] (1,0) -- (4.5,0); 
		\draw [thin] (1,0) -- (1,2); 
		
		\draw [thick,->] (1,-0.5) -- (4.7,-0.5);
		\node [right] at (4.9,-0.5) {$\Re \alpha$};
		
		\draw [thick,->] (0.7,0) -- (0.7,2.2);
		\node [above] at (0.7,2.2) {$\frac{1}{q}$};
		\draw [fill=black] (0.7,2) circle [radius = .5pt];
		\node [left] at (0.7,2) {$\frac{1}{2}$};
		\draw [fill=black] (0.7,0) circle [radius = .5pt];
		\node [left] at (0.7,0) {$0$};
		\draw [fill=black] (0.7,1.5) circle [radius = .5pt];
		\node [left] at (0.7,1.5) {$\frac{1}{p}$};
		\draw [fill=black] (1,-0.5) circle [radius = .5pt];
		\node [below] at (1,-0.5) {$0$};
		\draw [fill=black] (2.7,-0.5) circle [radius = .5pt];
		\node [below] at (2.7,-0.5) {$\frac{n}{2[p,2]_\theta}$};

		\path [fill=lightgray, opacity = 0.6] (P0)--(P2)--(P3)--(P4)--(P5)--(P6)--(P1);
		\path [fill={red!80!black}, opacity = 0.4] (P1)--(P6)--(P5);

		\end{tikzpicture}
	\end{center}
\caption{Visualization of the interpolation in Proposition~\ref{prop: upper SFE abstract}. For $(\Re \alpha,\nicefrac{1}{q})$ in the interior of the grey shaded region, $\IQ_{\psi_\alpha,T}$ is bounded $\L^{q} \to \T^{q}$ with a bound $C \e^{(\pi-2\mu) \Im \alpha}$, where $C$ is independent of $\Im \alpha$. Stein interpolation in Step~3 provides boundedness $\L^{q} \to \T^{q}$ in the interior of the red triangular region, the lower boundary of which is given by $\frac{1}{q} = \frac{1}{p} - \frac{[p,2]_\theta}{p} \frac{2 \Re \alpha}{n}$.}
\label{fig: Stein}
\end{figure}

\medskip

\noindent \emph{Step 3: Abstract Stein interpolation}. For technical reasons it will be more convenient to work with the `truncated' operators
\begin{align*}
\IQ_{\psi_\alpha, T}^{(k)} f = \e^{\alpha^2} \ind_{K_k} (\IQ_{\psi_\alpha,T}f) \quad (k \in \IN),
\end{align*}
where $K_k \coloneqq (k^{-1},k) \times B(0,k) \subseteq \reu$. For fixed $z$ the map $\alpha \mapsto \psi_\alpha(z)$ is holomorphic in the half plane $\IC^+$. Writing out the Cauchy integral for $\psi_\alpha(t^2T)$ and applying the dominated convergence theorem (justified by \eqref{eq2: upper SFE abstract}), we obtain that
\begin{align*}
\IC^+ \to \L^2(K_k), \quad \alpha \mapsto \IQ_{\psi_\alpha,T}^{(k)} f
\end{align*}
is holomorphic, whenever $f \in \L^2$. Moreover, thanks to the factor $\e^{\alpha^2}$ this mapping is \emph{qualitatively} bounded on any strip $\{\alpha \in \IC : c_0 \leq \Re \alpha \leq c_1\} \subseteq \IC^+$ with a bound depending on all parameters at stake. By the choice of $K_k$, the square function $S(\IQ_{\psi_\alpha,T}^{(k)} f)(x)$ vanishes for $x \in {}^cB(0,2k)$. Hence, we get for any $p \in (1,\infty)$ that
\begin{align*}
\|\IQ_{\psi_\alpha,T}^{(k)} f\|_{\T^p} 
\leq |B(0,2k)|^{\frac{1}{p}} k^{\frac{1+n}{2}} \|\IQ_{\psi_\alpha,T}^{(k)} f\|_{\L^2(K_k)},
\end{align*}
which shows that the qualitative mapping properties remain valid if we replace the target space $\L^2(K_k)$ by $\T^p$. 

If in addition $(\Re \alpha, \nicefrac{1}{q})$ belongs to the interior of the gray shaded region in Figure~\ref{fig: Stein} and $f \in \L^q \cap \L^2$, then we obtain the \emph{quantitative} bound
\begin{align*}
\|\IQ_{\psi_{\alpha,T}}^{(k)}f \|_{\T^{q}} 
&\leq |\e^{\alpha^2}| \|\IQ_{\psi_{\alpha,T}}f \|_{\T^{q}} 
\lesssim  \e^{(\Re \alpha)^2} \|\psi\|_{\sigma,\tau, \mu} \|f\|_{\L^{q}},
\end{align*}
where in the second step we have used \eqref{eq3: upper SFE abstract} and the implicit constant is independent of $\psi$, $\Im \alpha$, $k$.

Now, let $(\Re \alpha_j, \nicefrac{1}{q_j})$, $j=0,1$, belong to the interior of the gray shaded region in Figure~\ref{fig: Stein}. We intend to use to Proposition~\ref{prop: Stein interpolation} for
\begin{align*}
T(z) \coloneqq \IQ_{\psi_{\gamma(z),T}}^{(k)}, \quad \gamma(z) \coloneqq (\Re \alpha_0) (1-z) + (\Re \alpha_1) z,
\end{align*}
and the interpolation couples $X_j \coloneqq \L^{q_j}$ and $Y_j \coloneqq \T^{q_j}$. The dense subspace is $Z \coloneqq \L^2 \cap \L^{q_0} \cap \L^{q_1}$. The \emph{qualitative} bounds above yield (i) and the continuity part of (ii) in Proposition~\ref{prop: Stein interpolation}. The \emph{quantitative} bounds determine the constants $M_j$ in (ii). Hence, we get for any $(\Re \alpha, \nicefrac{1}{q})$ on the segment connecting the $(\Re \alpha_j, \nicefrac{1}{q_j})$ a bound
\begin{align*}
\|\IQ_{\psi_{\alpha,T}}^{(k)}f \|_{\T^{q}}
&\lesssim  \|\psi\|_{\sigma,\tau, \mu} \|f\|_{\L^{q}} \qquad (f\in \L^q \cap \L^2),
\end{align*}
where the implicit constant is independent of $\psi$ and $k$. Finally, we can pass to the limit as $k \to \infty$ via Fatou's lemma to obtain the same type of bound with $\IQ_{\psi_\alpha,T} f$ on the left-hand side. We have now completed Figure~\ref{fig: Stein} by adding the triangular region.

\medskip

\noindent \emph{Step 4: Conclusion}. We specialize to $\alpha = \sigma$, so that $\psi_\alpha = \psi$. The corresponding boundedness properties for $\IQ_{\psi,T}$ are dictated by Figure~\ref{fig: Stein}. If $\sigma \leq \frac{n}{2[p,2]_\theta}$, then $\frac{1}{q} > \frac{1}{p} - \frac{[p,2]_\theta}{p} \frac{2\sigma}{n}$ is needed. If $\sigma > \frac{n}{2[p,2]_\theta}$, then every $q \in [2, \infty)$ is admissible and this coincides with the range obtained in the first case.
\end{proof}

We single out the conclusion for the operator $L$ and the most common auxiliary functions $\psi$.\index{square function!bounds for $L$} Note that we can allow any $\psi \in \Psi_+^+$ when $p \geq 2$, which is a significant improvement compared to what is predicted by the abstract theory in Section~\ref{subsec: Hardy abstract sectorial}.

\begin{thm}
\label{thm: SFE from identification}
Let $p_{-}(L) < p < p_{+}(L)$ and let $\sigma, \tau > 0$. Let $\psi$ be of class $\Psi_\sigma^\tau$ on any sector. Then
\begin{align*}
\|S_{\psi,L} f\|_{p} \simeq \|a f\|_{\H^p} \quad (f \in \L^2),
\end{align*}
provided that
\begin{itemize}
	\item $\tau > |\nicefrac{n}{4} - \nicefrac{n}{(2p)}|$ and $\sigma > 0$ if $p \leq 2$,
	\item $\tau > 0$ and $\sigma > 0$ if $p \geq 2$.
\end{itemize}
Moreover, the upper square function bound ` $\lesssim$' remains to hold for $p_+(L) \leq p < \nicefrac{np_{+}(L)}{(n - 2\sigma p_+(L))}$, where the upper exponent bound is interpreted as $\infty$ if $2\sigma p_+(L) >n$.
\end{thm}

\begin{proof}[Proof of Theorem~\ref{thm: SFE from identification}]
If $p \leq 2$, then the assumption means that $\psi$ is an admissible auxiliary function for $\IH^p_L$, see Section~\ref{subsec: Hardy abstract sectorial}. Hence, 
\begin{align*}
\|S_{\psi,L} f\|_{p} = \|\IQ_{\psi,L}f\|_{\T^{p}} \simeq \|f\|_{\IH^p_L} \simeq \|af\|_{\H^p},
\end{align*} 
where the final step is due to Theorem~\ref{thm: main result Hardy}.

If $2 < p < p_+(L)$, then our assumptions on $\psi$ are less restrictive than the ones predicted by the abstract theory. 

We begin with the upper bounds. By Theorem~\ref{thm: main result Hardy} we have $\IH^p_L = \L^p \cap \L^2$ with equivalent $p$-norms. Hence, we can apply Proposition~\ref{prop: upper SFE abstract} for any $p \in (2,p_+(L))$ and by assumption on $\psi$ we may do so for any $\theta \in (0,1)$. This leads to  
\begin{align*}
\|S_{\psi,L} f\|_q \lesssim \|f\|_{q} \quad (f \in \L^q \cap \L^2)
\end{align*}
for any $q \geq 2$ that satisfies $\nicefrac{1}{q} > \nicefrac{1}{p_+(L)} - \nicefrac{2\sigma}{n}$, which is the range stated in the theorem.

For the lower bound we let $f \in \L^2$ with $S_{\psi,L} f \in \L^p$ and take $\varphi \in \Psi_\infty^\infty$ as in Remark~\ref{rem: Calderon reproducing} so that we have the reproducing formula
\begin{align*}
 f = \int_0^\infty \varphi(t^2L) \psi(t^2L) f \, \frac{\d t}{t}.
\end{align*}
Now, we refine the duality argument of Part~8 in the proof of Theorem~\ref{thm: main result Hardy}. We write again $L^* = a^* L^\sharp (a^*)^{-1}$, with $L^\sharp$ an operator in the same class as $L$ and $p' \in (p_{-}(L^\sharp) \vee 1, 2)$. For all $g \in \L^{p'} \cap \L^2$ we get 
\begin{align*}
\langle f, g \rangle
&= \int_0^\infty \langle \psi(t^2L) f, a^* \varphi^*(t^2L^\sharp) (a^*)^{-1}g \rangle \, \frac{\d t}{t} \\
&= \int_0^\infty \int_{\R^n} \IQ_{\psi,L}f \cdot \cl{a^* \IQ_{\varphi^*, L^\sharp} (a^*)^{-1}g} \, \frac{\d x \d t}{t},
\end{align*}
where $\langle \cdot \,, \cdot \rangle$ denotes the $\L^2$ inner product. Thus,
\begin{align*}
|\langle f, g \rangle| 
&\leq \|\IQ_{\psi,L}f\|_{\T^p} \|a^* \IQ_{\varphi^*, L^\sharp} (a^*)^{-1}g\|_{\T^{p'}}\\
&\lesssim \|S_{\psi,L}f\|_{p} \|S_{\varphi^*, L^\sharp} (a^*)^{-1}g\|_{p'} \\
&\lesssim  \|S_{\psi,L}f\|_{p}  \|g\|_{p'},
\end{align*}
where the first step is by the $\T^p - \T^{p'}$ duality and the third step uses the upper square function bound with $\varphi^* \in \Psi_\infty^\infty$ for $L^\sharp$ on $\L^{p'}$. Since $g \in \L^{p'} \cap \L^2$ was arbitrary, the lower bound $\|S_{\psi,L}f\|_{p}  \gtrsim \|f\|_p$ follows.
\end{proof}
\section{A digression: $\H^\infty$-calculus and analyticity}
\label{sec: consequences identification}

\noindent In this short section we present two consequences of the identification theorem for operator-adapted Hardy spaces that are of independent interest. One concerns analyticity, the other one concerns the $\H^\infty$-calculus for $L$. 

Recall that the standard assumptions \eqref{eq: standard assumptions sectorial} that we use to build the $L$-adapted spaces depend only on the configuration on $\L^2$: sectoriality, $\H^\infty$-calculus and off-diagonal estimates for the resolvents $(1+t^2L)^{-1}$ with real $t$.
By the sectorial version of Proposition~\ref{prop: FC on pre-Hardy} discussed in Section~\ref{subsec: Hardy abstract sectorial},
all $L$-adapted spaces inherit the $\H^\infty$-calculus with the same angle as on $\L^2$. 
 
It follows from Theorem~\ref{thm: main result Hardy} that we obtain $\H^\infty$-calculi for $L$ on classical $\H^p$ and $\Hdot^{1,q}$-spaces with the best possible angle. In the range $p \in (1,\infty)$, such results on $\L^p$ could in principle be obtained from Blunck and Kunstmann's theorem~\cite{BK}. 
This is the road taken in \cite[Sec.~5]{A} when $a=1$. We are not aware of 
an analog of the Blunck--Kunstmann result  on Hardy--Sobolev spaces. In fact, we are not even aware of any general results for $p \leq 1$ or $q\le 1$ or even of functional calculus  away from the Banach space range.\index{Hinfty calculus@$\H^\infty$-calculus!on $\H^p$ and $\Hdot^{1,p}$}

We summarize this discussion in the following result. 

\begin{thm}
\label{thm: main result Hoo}
Let $p_-(L) < p < p_+(L)$ and $(p_{-}(L)_* \vee 1_*) < q <  q_+(L)$. For every $\nu \in (\omega, \pi)$ the functional calculus bounds
\begin{align*}
	\|a  \eta(L) a^{-1}f \|_{\H^p} &\lesssim \|\eta\|_\infty \|f\|_{\H^p} \\
	\|\eta(L) g \|_{\Hdot^{1,q}} &\lesssim \|\eta\|_\infty \|g\|_{\Hdot^{1,q}}
\end{align*}
hold for all $\eta \in \H^\infty(\S_\nu^+)$ and all $f \in \H^p \cap \L^2$, $g \in \Hdot^{1,q} \cap \L^2$.
\end{thm}

The open $p$-interval in Theorem~\ref{thm: main result Hoo} is the largest possible one since $\eta(\zeta) = (1+t^2 \zeta)^{-1}$ with real $t$ is admissible.  An example that illustrates the less familiar second inequality is $\|\nabla_{x }(1+t^2 L)^{-1} g \|_{\H^{q}} \lesssim  \|\nabla_{x }g\|_{\H^{q}}$, which is of a different nature than the bounds defining $\cN(L)$ and is valid for $q$  in a   bigger set. 

This also leads us to analyticity, that is, resolvent bounds for parameters in a sector in the complex plane. According to Section~\ref{subsec: bisectorial operators}, $L$ is sectorial in $\L^2$ with angle $\omega_{L}$ not exceeding $2\omega_{DB}<\pi$. We obtain for $\X$ being any of the spaces in the statement above and every $\mu \in (\omega_{L}, \pi)$ that there are extensions by density with operator norm bounds
\begin{align*}
  \sup_{z \in \IC \setminus \cl{\S_\mu^+}}  \|z(z-L)^{-1}||_{\X \to \X} < \infty.
\end{align*}
This means that $\H^p$-boundedness of resolvents $(1+t^2L)^{-1}$ with real $t$ alone self-improves to the same properties for the resolvents $(1+z^2 L)^{-1}$ for $z \in \S_\mu^+$ and $\mu \in (0, \nicefrac{(\pi - \omega_L)}{2})$.

A similar discussion applies to $\L^p$ off-diagonal estimates for $T(z) \coloneqq (1+z^2 L)^{-1}$, $z \in \S_\mu^+$, when $(p_{-}(L)\vee 1)<p<p_{+}(L)$. For a small and $p$-dependent angle they can be obtained from the Stein  interpolation theorem for analytic families of operators, see Lemma~\ref{lem: OD extrapolation to sectors}. Having the $\L^p$-boundedness and the $\L^2$ off-diagonal estimates for the ($p$-independent) optimal angle implies by complex interpolation applied to each single operator $T(z)$
the $\L^p$ off-diagonal estimates for $T(z)$, see Lemma~\ref{lem: OD interpolation}.
If $p_{-}(L)<1$ (resp. $p_{-}(L^\sharp)<1$), we shall see in Section~\ref{sec: Critical numbers and kernels} that we may also include $\L^1$ (resp. $\L^\infty$) off-diagonal estimates here.

In the same manner, we could obtain self-improvements for other families. Of particular interest is the analytic Poisson semigroup generated by $-L^{1/2}$, which has angle $\nicefrac \pi 2- \nicefrac {\omega_{L}} 2$, and when  $\omega_{L}<\nicefrac \pi 2$ --- that is, for instance when $a=1$ --- the analytic heat semigroup $\e^{-zL}$ with angle $\nicefrac \pi 2- \omega_{L}$. 
\section{Riesz transform estimates: Part II}
\label{sec: Riesz 2}

\noindent We come back to the Riesz transform interval
\begin{align*}
\cI(L) \coloneqq \big\{ p \in (1_*,\infty) : R_L \text{ is } \text{ $a^{-1} \H^p - \H^p$-bounded}\big\},
\end{align*}
defined in \eqref{eq: I(L)}, the endpoints of which we have denoted by $r_\pm(L)$. In Section~\ref{sec: Riesz} we have characterized the endpoints of the part of $\cI(L)$ in $(1,\infty)$. The identification theorem for adapted Hardy spaces allows us to complete the discussion through the following theorem. 

\begin{thm}
\label{thm: Riesz complete}
It follows that\index{Riesz transform!$\H^p$-boundedness} \index{I@$\cI(L)$!characterization of}
\begin{align*}
 \cI(L) = (p_{-}(L), q_{+}(L)).
\end{align*}
Moreover, the following hold true:
\begin{enumerate}
	\item The map $aL^{1/2} : \Hdot^{1,p} \cap \Wdot^{1,2} \to \H^p \cap \L^2$ is well-defined and bounded for the $p$-quasinorms if $p_{-}(L)_* \vee 1_* < p < p_+(L)$.
	\item An exponent $p \in (1_*, \infty)$ belongs to $\cI(L)$ if and only if the map in (i) extends by density to an isomorphism $\Hdot^{1,p} \to \H^p$ whose inverse agrees with $L^{-1/2} a^{-1}$ on $\H^p \cap \L^2$. In particular, if $p \in \cI(L)$, then
	\begin{align*}
	\|R_L f\|_{\H^p} \simeq \|af\|_{\H^p} \quad (f \in a^{-1} (\H^p \cap \L^2)).
	\end{align*}
\end{enumerate}
\end{thm}

The reader may wonder if the separate discussion in Section~\ref{sec: Riesz} could have been avoided. The answer is that it can not, since Theorem~\ref{thm: Riesz} was used in proving Theorem~\ref{thm: main result Hardy}.

\begin{proof}
The Hardy space theory yields for $p_{-}(L)_* \vee 1_* < p < p_+(L)$ continuous inclusions for the $p$-quasinorms,
\begin{align}
 \label{eq0: Riesz complete}
 \begin{split}
 \Hdot^{1,p} \cap \L^2  &\subseteq \IH_L^{1,p}, \\
 \IH_L^p &\subseteq a^{-1} (\H^p \cap \L^2).
\end{split}
\end{align}
More precisely, by Theorem~\ref{thm: main result Hardy} the first inclusion is an equality up to equivalent quasinorms if $p < q_+(L)$ and the second one is an equality if $p > p_{-}(L)$. The first inclusion for $q_+(L) \leq p < p_+(L)$ is due to \eqref{eq: H1p in H1pL for p>2} and the second inclusion for $p_-(L)_* \vee 1_*< p \leq p_{-}(L)$ is due to \eqref{eq: HpL in Hp}.

\medskip

\noindent \emph{Step~1: Proof of (i).} As Figure~\ref{fig: diagram} tells us that $L^{1/2}: \IH_L^{1,p} \cap \W^{1,2} \to \IH_L^p$ is bounded for the $p$-quasinorms, we conclude from the inclusions above that $aL^{1/2}: \Hdot^{1,p} \cap \W^{1,2} \to \H^p \cap \L^2$ is well-defined and bounded for the respective $p$-quasinorms. The extension to $\Hdot^{1,p} \cap \Wdot^{1,2}$ follows by density.

\medskip

\noindent \emph{Step~2: Bounds for $R_L$.} Let $p_{-}(L) < p < q_{+}(L)$. Then the inclusions in \eqref{eq0: Riesz complete} become equalities and Figure~\ref{fig: diagram} tells us that
\begin{align*}
\|\nabla_x L^{-1/2} f\|_{\H^p} \simeq \|af\|_{\H^p} \quad (f \in \IH_L^p \cap \ran({L^{1/2}})).
\end{align*}
Since $\IH_L^p \cap \ran({L^{1/2}})$ is dense in $\IH_L^p$ for the norm $\|\cdot\|_{\IH_L^p} + \|\cdot\|_2$, we obtain by approximation and the various quasinorm equivalences that
\begin{align*}
\|R_L f\|_{\H^p} \simeq \|af\|_{\H^p} \quad (f \in a^{-1} (\H^p \cap \L^2)).
\end{align*}
In particular, $R_L$ is $a^{-1} \H^p - \H^p$-bounded.

\medskip

\noindent \emph{Step~3: Identification of the endpoints of $\cI(L)$.} In view of Theorem~\ref{thm: Riesz} it remains to show $p_-(L) = r_-(L)$ in the case that one of these exponents is smaller than $1$. In Step~2 we have already shown $r_{-}(L) \leq p_{-}(L)$ without any such restrictions. The only task remaining is to prove that $r_-(L) < \varrho < 1$ implies $(\varrho,1] \subseteq \cJ(L)$. 

We may assume $ \varrho < 2_*$ since otherwise the claim already follows from Proposition~\ref{prop: J(L) contains neighborhood of Sobolev conjugates}. Since $r_-(L) \vee 1 < \varrho^* < 2$, Theorem~\ref{thm: Riesz} yields $p_-(L) \vee 1 < \varrho^* < 2$ and hence $(tL^{1/2} (1+t^2 L)^{-1})_{t>0}$ is $\L^{\varrho^*}$-bounded, see Lemma~\ref{lem: functional calculus bounds from J(L) abstract}.(i). Now, let $f \in \H^\varrho \cap \L^2$. Then $f \in \ran(a{L^{1/2}})$ thanks to Lemma~\ref{lem: density for Riesz}, so that we can estimate
\begin{align*}
\|(1+t^2 L)^{-1} a^{-1} f \|_{\varrho^*}
&\lesssim t^{-1} \|L^{-1/2} a^{-1}f \|_{\varrho^*} \\
&\lesssim t^{-1} \|\nabla_x L^{-1/2} a^{-1} f\|_{\H^\varrho} \\
&\lesssim t^{-1} \|f\|_{\H^\varrho},
\end{align*}
where we used the assumption $r_{-}(L)<\varrho$ in the last line. 
This means that the resolvents are $a^{-1} \H^\varrho - \L^{\varrho^*}$-bounded. According to Lemma~\ref{lem: Lp bounds for resolvents on sector} they satisfy $\L^{\varrho^*}$ off-diagonal estimates of arbitrarily large order and for compactly supported $f \in \L^2$ with mean value zero we recall from Corollary~\ref{cor: conservation} that $\int_{\R^n} a(1+t^2 L)^{-1} (a^{-1}f) \d x = 0$. With these properties at hand, the required $\H^p$-boundedness of the resolvents for $p \in (\varrho,1]$ follows from Lemma~\ref{lem: OD Hardy implies boundedness}.

\medskip

\noindent \emph{Step~4: Proof of (ii).} If $aL^{1/2}$ extends to an isomorphism with the given property, then
\begin{align*}
\|R_L f\|_{\H^p} = \|L^{-1/2} a^{-1} (af)\|_{\Hdot^{1,p}} \simeq  \|af\|_{\H^p} \quad (f \in a^{-1} (\H^p \cap \L^2))
\end{align*}
as required. 

Conversely, suppose that $p \in \cI(L)$. This means that $R_L = \nabla_x L^{-1/2}$ is $a^{-1}\H^p - \H^p$-bounded and hence $L^{-1/2} a^{-1}: \H^p \cap \L^2 \to \Hdot^{1,p} \cap \Wdot^{1,2}$ is well-defined and bounded for the $p$-quasinorms. According to Step~3 the exponent $p$ must be contained in $[p_{-}(L), q_+(L)] \cap (1_*, \infty)$, which, in view of Theorem~\ref{thm: standard relation J(L) and N(L)}, is a subset of the interval considered in (i). Therefore $aL^{1/2} : \Hdot^{1,p} \cap \Wdot^{1,2} \to \H^p \cap \L^2$ is also bounded for the $p$-quasinorms and hence it extends to an isomorphism with the required properties.

\medskip

\noindent \emph{Step~5: Conclusion.} We already know the endpoints of $\cI(L)$ and it remains to show that this interval is open. The map in (i) is defined and continuous for $p$ in an open interval $I$ that contains $\cI(L)$ and the isomorphism property in (ii) characterizes $\cI(L)$ as a subset of $I$. Since the scales of spaces $(\Hdot^{1,p})_{p \in (1_*,\infty)}$ and $(\H^p)_{p \in (1_*,\infty)}$ interpolate by the complex method, the openness of $\cI(L)$ is a consequence of \v{S}ne\u{i}berg's stability theorem~\cite[Thm.~8.1]{KMM}. See also \cite[Thm.~8.1]{KMM} for the fact that compatibility of the inverses is preserved.\index{Theorem!\v{S}ne\u{i}berg's}
\end{proof}

In Part~10 of the proof of Theorem~\ref{thm: main result Hardy} we have seen that $q_+(L) \in \cH^1(L)$ is possible only if the Riesz transform is $\L^{q_+(L)}$-bounded. Hence, we can note:

\begin{cor}
\label{cor: H1L upper endpoint}
The interval $\cH^1(L)$ is open at the upper endpoint, that is, $q_+(L) \notin \cH^1(L)$.
\end{cor}

The statement (ii) in Theorem~\ref{thm: Riesz complete} can be strengthened to a Riesz transform characterization of abstract and concrete Hardy spaces. For operators of type $- \div_x d \nabla_x$ such results first appeared in \cite[Sec.~5]{HMMc}.\index{Hardy space! Riesz transform characterization via $L$} Interestingly, this observation allows us to strengthen the identification theorem for $\IH_L^p$ itself in that $\cH(L)$ is open and hence identification fails at the endpoints.

\begin{thm}
\label{thm: Riesz characterization HpL}
Let $p \in (p_{-}(L)_* \vee 1_*, q_+(L))$. Then 
\begin{align*}
\IH_L^p = \{f \in \L^2 : R_L f \in \H^p\}
\end{align*}
with equivalent quasinorms $\|\cdot\|_{\IH_L^p} \simeq \|R_L \cdot\|_{\H^p}$. 
In particular, it follows that\index{H@$\cH(L)$!characterization of} 
\begin{align*}
 \cH(L) = (p_-(L), p_+(L)).
\end{align*}
\end{thm}

\begin{proof}
Let $p \in (p_{-}(L)_* \vee 1_*, q_+(L))$. We first prove the quasinorm equivalence for $f \in \IH_L^p$. To this end, we argue as in Step~2 of the proof of Theorem~\ref{thm: Riesz complete}, except that in the given range of exponents only the first inclusion in \eqref{eq0: Riesz complete} is an equality but we cannot identify $\IH_L^p$ unless $p>p_{-}(L)$. This yields
\begin{align*}
\|R_L f\|_{\H^p}  \simeq \|f\|_{\IH_L^p} \quad (f \in \IH_L^p)
\end{align*}
and we can replace $\IH_L^p$ with $a^{-1} (\H^p \cap \L^2)$ if in addition $p>p_{-}(L)$. 

Conversely, let $f \in \L^2$ satisfy $R_L f \in \H^p$. Arguing as in Step~1 of the proof of Theorem~\ref{thm: Riesz complete}, we find that $L^{1/2} : \Hdot^{1,p} \cap \W^{1,2} \to \IH_L^p$ is bounded for the $p$-quasinorms. The only difference is again that we cannot identify $\IH_L^p$. By assumption we have $L^{-1/2} f \in \Hdot^{1,p} \cap \Wdot^{1,2}$. Let $(u_k) \subseteq \cZ$ be a sequence with $u_k \to L^{-1/2}f$ in $\Hdot^{1,p} \cap \Wdot^{1,2}$ as $k \to \infty$ and set $f_k \coloneqq L^{1/2}u_k$. Then $(f_k)$ is a Cauchy sequence in (the possibly non-complete space) $\IH_L^p$ that converges to $f$ in $\L^2$. Let $\IH_L^p$ be defined by the auxiliary function $\psi$. By $\L^2$ convergence
\begin{align*}
\int_{B(x,t)} |\psi(t^2L)f(y)|^2 \, \d y = \lim_{k \to \infty} \int_{B(x,t)} |\psi(t^2L)f_k(y)|^2 \, \d y
\end{align*}
holds for all $(t,x) \in \reu$ and Fatou's lemma yields
\begin{align*}
\|f\|_{\IH^p_L}^p 
&\leq \liminf_{k \to \infty} \int_{\R^n} \bigg(\iint_{|x-y|<t} |\psi(t^2L)f_k(y)|^2 \, \frac{\d y \d t}{t^{1+n}} \bigg)^{p/2} \, \d x \\
&= \liminf_{k \to \infty} \|f_k\|_{\IH^p_L}^p.
\end{align*}
The final expression is finite by the Cauchy property in $\IH_L^p$, which means that $f \in \IH_L^p$.

Concerning the final statement, we recall from Theorem~\ref{thm: main result Hardy} that $p_\pm(L)$ are the endpoints of $\cH(L)$. For the sake of a contradiction, suppose $p \coloneqq p_-(L) \in \cH(L)$. The first part yields $\|R_L f\|_{\H^p} \simeq \|af\|_{\H^p}$ for all $f \in a^{-1}(\H^p \cap \L^2)$, which contradicts Theorem~\ref{thm: Riesz complete}. Likewise, suppose $p_+(L) \in \cH(L)$. Since $\cH(L) \subseteq (1_*,\infty)$, we must have $p_+(L) < \infty$ and therefore $ p_-(L^\sharp) = p_+(L)'>1$ by duality and similarity. Proposition~\ref{prop: Hardy duality} implies with $p=p_-(L^\sharp)$ that  $\IH_{L^\sharp}^{p} = \L^{p} \cap \L^2$ with equivalent $p$-norms, that is $p  \in \cH(L^\sharp)$, which is impossible as we have already seen.
\end{proof}
\section{Critical numbers for Poisson and heat semigroups}
\label{sec: Poisson}

\noindent In this chapter, we show that the critical numbers are intrinsic in the sense that we could have equivalently defined them through other families of functions of $L$ than resolvents. For the applications to boundary value problems we are mainly interested in estimates for the \emph{Poisson semigroup}\index{Poisson semigroup!for $L$} $(\e^{-t L^{1/2}})_{t >0}$. It would have been equally natural to try building the theory in Section~\ref{sec: J and N} from the intervals
\begin{align*}
\cJ^{\mathrm{Pois}}(L) 
&\coloneqq \big\{p \in (1_*, \infty) : (\e^{-t L^{1/2}})_{t>0} \text{ is } a^{-1}\H^{p}\text{-bounded}\big\}, \\
\cN^{\mathrm{Pois}}(L)
&\coloneqq \big\{p \in (1_*, \infty) : (t \nabla_x \e^{-t L^{1/2}})_{t>0} \text{ is } a^{-1}\H^{p} - \H^{p}\text{-bounded}\big\}.
\end{align*}

Note that both intervals contain $p=2$. Indeed, $2 \in \cJ^{\mathrm{Pois}}(L)$ follows from the functional calculus on $\L^2$ and to prove $2 \in \cN^{\mathrm{Pois}}(L)$ we additionally use ellipticity to obtain
\begin{align*}
	\|t \nabla_x \e^{-t L^{1/2}}f\|_2^2 \lesssim \Re \langle a t^2L \e^{-t L^{1/2}}f, \e^{-t L^{1/2}}f \rangle \lesssim \|f\|_2^2
\end{align*} 
for all $t>0$ and all $f \in \L^2$. This gives rise to a definition of critical `Poisson' numbers.\index{critical numbers!via the Poisson semigroup} 
\begin{defn}
\label{def: Poisson critical numbers}
The lower and upper endpoints of $\cJ^{\mathrm{Pois}}(L)$ are denoted by $p_-^{\mathrm{Pois}}(L)$ and $p_+^{\mathrm{Pois}}(L)$, respectively. Likewise, $q_\pm^{\mathrm{Pois}}(L)$ denote the endpoints of $\cN^{\mathrm{Pois}}(L)$.
\end{defn}

The reason why we use $\cJ(L)$ and $\cN(L)$ is that Poisson semigroups offer very limited off-diagonal decay (think of the Poisson kernel for the Laplacian), whereas the resolvents offer exponential decay. One main result in this section is that while the decay properties are strikingly different, the associated critical numbers are the same. 

\begin{thm}
\label{thm: Poisson semigroup interval}
 $p_\pm^{\mathrm{Pois}}(L) = p_\pm(L)$ and $q_\pm^{\mathrm{Pois}}(L) = q_\pm(L)$.
\end{thm}

Aiming in a similar direction, we note that the unperturbed operator $	L_0 =  -\div_x d \nabla_x$ is sectorial of angle $\omega_{L_0} \in (0, \nicefrac{\pi}{2})$ and hence it generates a holomorphic semigroup $(\e^{-t^2 L_0})_{t>0}$ on $\L^2$, called \emph{heat semigroup}. The associated intervals
\begin{align*}
\cJ^{\mathrm{heat}}(L_0)
&\coloneqq \big\{p \in (1_{*}, \infty) : (\e^{-t^2L_0})_{t>0} \text{ is } \H^{p}\text{-bounded}\big\}, \\
\cN^{\mathrm{heat}}(L_0)
&\coloneqq \big\{p \in (1_{*}, \infty) : (t \nabla_x \e^{-t^2L_0})_{t>0} \text{ is } \H^{p}\text{-bounded}\big\}
\end{align*}
contain $p=2$ by the same argument as for the Poisson semigroup and their endpoints are the critical `heat' numbers.\index{critical numbers!via the heat semigroup}

\begin{defn}
\label{def: heat critical numbers}
The lower and upper endpoints of $\cJ^{\mathrm{heat}}(L_0)$ are denoted by $p_-^{\mathrm{heat}}(L_0)$ and $p_+^{\mathrm{heat}}(L_0)$, respectively. Likewise $q_\pm^{\mathrm{heat}}(L_0)$ denote the endpoints of $\cN^{\mathrm{heat}}(L_0)$.
\end{defn}

We refer to \cite{AA} for a systematic treatise of critical heat numbers in the range $p \in (1,\infty)$ and their relation to Riesz transforms, $\H^\infty$-calculus and square function estimates. 

The second main result in this section shows that critical numbers and critical heat numbers are the same in the full interval of exponents. Since also the critical numbers for $L_0$ and $L$ are the same (Theorem~\ref{thm: critical numbers a-independent}), this provides a means of characterizing all intervals of exponents in the monograph through properties of a heat semigroup, even though $L$ itself need not be a generator.

\begin{thm}
\label{thm: heat semigroup interval}
$p_\pm^{\mathrm{heat}}(L_0) = p_\pm(L_0)$ and $q_\pm^{\mathrm{heat}}(L_0) = q_\pm(L_0)$.
\end{thm}

This second result tells us that the theory in \cite{A} relying on the critical heat numbers is in coherence with the one here when restricted to the range $p \in (1,\infty)$. 

The proofs of both theorems follow the same pattern. If we assume resolvent bounds, then semigroup bounds follow immediately from the functional calculus bound in Theorem~\ref{thm: main result Hoo}, whereas in the opposite direction we can represent resolvents via Laplace transforms of the semigroup. For the Poisson semigroup these formul\ae \, become more technical since we want to estimate the resolvents of the square of the semigroup generator $L^{1/2}$. It is for this reason that we shall showcase the strategy for the heat semigroup first although the Poisson semigroup is of greater importance to us.

For both proofs we need part (i) of the following proposition. The extensions in (ii) and (iii) will be needed much later in Section~\ref{sec: single layer}.

\begin{prop}
\label{prop: Poisson bounds}
Let $p_-(L) < p \leq q < p_+(L)$ and consider the families $(a\e^{-t {L^{1/2}}} a^{-1})_{t>0}$ and $(\e^{-t^2 L_0})_{t>0}$.
\begin{enumerate} 
	\item Both families are $\H^p - \H^q$-bounded.
	\item If $p_-(L^\sharp)<1$ and $0<\alpha< n (\nicefrac{1}{p_-(L^\sharp)}-1)$, then the first one is $\H^p - a \Lamdot^{\alpha}$-bounded and the second one is $\H^p - \Lamdot^{\alpha}$-bounded
	\item If $p_-(L^\sharp)<1$, then they are  $\H^p - \L^\infty$-bounded.
\end{enumerate}
\end{prop}

\begin{proof}
We prove the three statements in order.

\medskip

\noindent \emph{Proof of (i).} We recall from Theorem~\ref{thm: critical numbers a-independent} that $p_\pm(L) = p_\pm(L_0)$. Hence, $\H^p$-boundedness follows directly from Theorem~\ref{thm: main result Hoo}. 

As we have $p_-(L) < 2_*$ (Proposition~\ref{prop: J(L) contains neighborhood of Sobolev conjugates}), we can use a Sobolev embedding followed by Theorem~\ref{thm: Riesz complete}.(ii) and Theorem~\ref{thm: main result Hoo} with exponent $p=2_*$ in order to obtain for all $t>0$ and all $f \in \H^{2_*} \cap \L^2$ that
\begin{align}
\label{eq1: Poisson bounds}
\begin{split}
	\|a\e^{-t {L^{1/2}}} a^{-1}f\|_2 
	&\lesssim \|\nabla_x \e^{-t {L^{1/2}}} a^{-1} f\|_{\H^{2_*}} \\
	&\simeq \| a L^{1/2} \e^{-t {L^{1/2}}} a^{-1} f\|_{\H^{2_*}}  \\
	& \lesssim t^{-1} \|f\|_{\H^{2_*}}.
\end{split}
\end{align}
Hence, $(a\e^{-t {L^{1/2}}} a^{-1})_{t>0}$ is $\H^{2_*}-\L^2$-bounded. By the first part and Lemma~\ref{lem: extra} we obtain for each $p \in (p_-(L),2)$ an integer $\beta$ such that $(a \e^{-t \beta {L^{1/2}}}a^{-1})_{t>0}$ is $\H^p - \L^2$-bounded. Interpolation with the first part yields $\H^p - \H^q$-boundedness for all exponents $p_-(L) < p \leq q \leq 2$. 

Applying this result to $L^\sharp$ and using $L^*= {a^*}L^\sharp ({a^*})^{-1}$ yields in particular $\L^{q'} - \L^{p'}$-boundedness for $(\e^{-t {(L^*)^{1/2}}})_{t>0}$ if $(p_-(L^\sharp) \vee 1) < q' \leq p' \leq 2$. Hence, $\L^p - \L^q$ boundedness of $(a\e^{-t {L^{1/2}}} a^{-1})_{t>0}$ follows for $2 \leq p \leq q < p_+(L)$ by duality and ellipticity of $a^*$. In the remaining case that $p$ and $q$ are on opposite sides of $2$ we can use the semigroup property and combine $\H^p - \L^2$ and $\L^2 - \L^q$-boundedness.

The proof for the heat semigroup is \emph{mutadis mutandis} the same since the second-order scaling guarantees that the third step in \eqref{eq1: Poisson bounds} remains valid.

\medskip

\noindent \emph{Proof of (ii).} Let $\alpha = n(\nicefrac{1}{\varrho}-1)$ with $p_-(L^\sharp)<\varrho<1$. Part~(i) yields that $(a^*\e^{-t {(L^\sharp)^{1/2}}}(a^*)^{-1})_{t>0}$ is $\H^\varrho - \L^{p'}$-bounded if $p_-(L) \vee 1 < p < p_+(L)$. By similarity and duality $(\e^{-t {L^{1/2}}})_{t>0}$ is $\L^p-\Lamdot^{\alpha}$-bounded. This is the claim under the additional assumption $p>1$. The full result follows from (i) by the semigroup property. 

The same argument applies to the heat semigroup.

\medskip

\noindent \emph{Proof of (iii).} By the semigroup property and (i) it suffices to treat the case $p >1$. The claim has nothing to do with semigroups and simply follows from (i), (ii) and the following interpolation inequality.
\end{proof}

\begin{lem}
\label{lem: Loo interpolation inequality}
Let $1 \leq p < \infty$ and $0<\alpha<1$. If $g \in \L^p\cap\Lamdot^{\alpha}$, then $g \in \L^\infty$ and
\begin{align*}
	\|g\|_\infty \leq 2 |B(0,1)|^{-\frac{\theta}{p}} \|g\|_p^\theta \|g\|_{\Lamdot^\alpha}^{1-\theta}, \quad \theta = \frac{\alpha}{\alpha + \nicefrac{n}{p}}.
\end{align*}
\end{lem}

\begin{proof}
For $x,y \in \R^n$ we have 
\begin{align*}
	|g(x)| \leq |g(y)| + |x-y|^\alpha \|g\|_{\Lamdot^\alpha}.
\end{align*}
We take the average in $y$ over some ball $B(x,r)$ and use H\"older's inequality to give 
\begin{align*}
	|g(x)| \leq (|B(0,1)| r^n)^{-1/p} \|g\|_p + r^\alpha \|g\|_{\Lamdot^\alpha}.
\end{align*}
We conclude by picking $r$ such that the terms on the right are equal.
\end{proof}

\subsection{Identification of the critical heat numbers}
\label{subsec: identification of the critical heat numbers}

We turn to the proof of the second principal result of this section.

\begin{proof}[Proof of Theorem~\ref{thm: heat semigroup interval}]
We break the argument into three steps.
	
\medskip
	
\noindent \emph{Step~1: From the resolvent to the semigroup}. Proposition~\ref{prop: Poisson bounds}.(i) implies $(p_-(L_0), p_+(L_0)) \subseteq \cJ^{\mathrm{heat}}(L_0)$.  

Next, we let $p \in (q_-(L_0), q_+(L_0))$. Then $p \in (p_-(L_0), p_+(L_0))$ by Theorem~\ref{thm: standard relation J(L) and N(L)}. Combining Theorem~\ref{thm: Riesz complete}.(ii) and Theorem~\ref{thm: main result Hoo}, we get
\begin{align*}
	\|t \nabla_x \e^{-t^2 L_0} f\|_{\H^p}
	\simeq \|t L_0^{1/2} \e^{-t^2 L_0} f\|_{\H^p}
	\lesssim \|f\|_{\H^p},
\end{align*}
which proves $p \in \cN^{\mathrm{heat}}(L_0)$. We conclude that $(q_-(L_0), q_+(L_0)) \subseteq \cN^{\mathrm{heat}}(L_0)$.

\medskip

\noindent \emph{Step~2: From $\cJ^{\mathrm{heat}}(L_0)$ to $\cJ(L_0)$}. For every $t>0$ the operator 
\begin{align*}
	T \coloneqq 1+t^2 L_0
\end{align*}
is invertible and sectorial of angle $\omega_{L_0} < \nicefrac{\pi}{2}$. By the Calderón reproducing formula we have for $f \in \L^2$ as an improper Riemann integral,
\begin{align*}
	f = \int_0^\infty T \e^{-sT}f \, \d s.
\end{align*}
Applying $T^{-1}$ on both sides gives the classical formula
\begin{align}
\label{eq1: heat semigroup interval}
(1+ t^2 L_0)^{-1}f 
&= \int_0^\infty \e^{-s} \e^{-st^2 L_0}f \, \d s
\end{align}
and the integral converges absolutely in $\L^2$ since the heat semigroup is uniformly bounded.

Let now $r \in \cJ^{\mathrm{heat}}(L_0)$ and take any $p$ between $r$ and $2$. We shall show that $((1+t^2 L_0)^{-1})_{t>0}$ is $\H^p$-bounded, that is, $p \in \cJ(L_0)$. Then, together with Step~1, $p_\pm^{\mathrm{heat}}(L_0) = p_\pm(L_0)$ follows.
	
\medskip

\noindent \emph{Step~2a: The Lebesgue case $p>1$}. Since the heat semigroup is $\L^p$-bounded, the integral in \eqref{eq1: heat semigroup interval} converges absolutely in $\L^p$ for all $t>0$ and all $f \in \L^p \cap \L^2$ and we obtain
\begin{align*}
	\|(1+ t^2 L_0)^{-1}  f\|_p \lesssim \|f\|_p
\end{align*}
as required.

\medskip
	
\noindent \emph{Step~2b: The Hardy case $p \leq 1$}. We appeal to Lemma~\ref{lem: OD Hardy implies boundedness} in order to show that the resolvents are $\H^p$-bounded. 

For $f \in \L^2$ with compact support and mean value zero we have $\int_{\R^n} (1+t^2L_0)^{-1} f \d x = 0$ by Corollary~\ref{cor: conservation}. For the other two assumptions in Lemma~\ref{lem: OD Hardy implies boundedness} we use exponents $\varrho \in (r,p)$ and $q \in (1,2)$ with $\nicefrac{n}{\varrho}-\nicefrac{n}{q} <1$. In particular, $\varrho, q$ are interior points of $\cJ^{\mathrm{heat}}(L_0)$.

From Step~2a we obtain $q \in (p_{-}(L_0), 2)$. Hence, $((1+t^2L_0)^{-1})_{t>0}$ satisfies $\L^q$ off-diagonal estimates of arbitrarily large order by interpolation with the $\L^2$ off-diagonal decay.

It remains to show that  $((1+t^2 L_0)^{-1})_{t>0}$ is $\H^\varrho - \L^q$-bounded. The following boundedness properties hold for the heat semigroup: first $\H^r$ and $\H^\varrho$ (by assumption), second $\L^q - \L^2$ (by Proposition~\ref{prop: Poisson bounds}), third $\H^\varrho - \L^2$ (by Lemma~\ref{lem: extra} and the semigroup law), fourth $\H^\varrho - \L^q$ (by interpolation). This allows us to take $\L^q$-norms in \eqref{eq1: heat semigroup interval} and obtain for all $t>0$ and all $f \in \H^\varrho \cap \L^2$,
\begin{align}
\label{eq2: heat semigroup interval}
\begin{split}
	\|(1+ t^2 L_0)^{-1} f\|_q
	&\lesssim \int_0^\infty \e^{- s} (s^{\frac{1}{2}} t)^{\frac{n}{q}-\frac{n}{\varrho}} \|f\|_{\H^\varrho} \, \d s \\
	&\lesssim  t^{\frac{n}{q}-\frac{n}{\varrho}} \|f\|_{\H^\varrho},
\end{split}
\end{align}
where the integral in $s$ is finite by the choice of our exponents. This completes Step~2b.

\medskip

\noindent \emph{Step~3: From $\cN^{\mathrm{heat}}(L_0)$ to $\cN(L_0)$}. Let $r \in \cN^{\mathrm{heat}}(L_0)$ and take any $p$ between $r$ and $2$. We shall show that $(t \nabla_x (1+t^2 L_0)^{-1})_{t>0}$ is $ \H^p$-bounded, that is, $p \in \cN(L_0)$. Then, together with Step~1, $q_\pm^{\mathrm{heat}}(L_0) = q_\pm(L_0)$ follows.
	
\medskip

\noindent \emph{Step~3a: The Lebesgue case $p>1$}. We apply $t \nabla_x$ on both sides of \eqref{eq1: heat semigroup interval} and take $\L^p$-norms in order to get
\begin{align*}
	\|t \nabla_x (1+ t^2 L_0)^{-1}f \|_p
	\lesssim \int_0^\infty s^{-\frac{1}{2}} \e^{- s} \|f\|_{p} \, \d s \lesssim  \|f\|_{p} 
\end{align*}
as required for all $t>0$ and all $f \in \L^p \cap \L^2$. 
\medskip
	
\noindent \emph{Step~3b: The Hardy case $p \leq 1$}. By Theorem~\ref{thm: standard relation J(L) and N(L)} the intervals $\cJ(L_0)$ and $\cN(L_0)$ have the same lower endpoint. Hence, it suffices to prove $p \in \cJ(L_0)$, that is, $\H^p$-boundedness of $((1+t^2 L_0)^{-1})_{t>0}$. Moreover, we can assume $p<2_*$ since otherwise we can directly conclude by Proposition~\ref{prop: J(L) contains neighborhood of Sobolev conjugates}. Once again we appeal to Lemma~\ref{lem: OD Hardy implies boundedness}. We fix any $\varrho \in (r,p)$ and let $q \coloneqq r^* \in (1,2)$. In particular, $\nicefrac{n}{\varrho}-\nicefrac{n}{q} <1$ and $\varrho, q$ are interior points of $\cN^{\mathrm{heat}}(L_0)$.

For $f \in \L^2$ with compact support and mean value zero we have $\int_{\R^n} (1+t^2L_0)^{-1} f \d x = 0$ by Corollary~\ref{cor: conservation}. From Step~3a we obtain that $q$ is an interior point of $\cN(L_0)$, hence of $\cJ(L_0)$. By interpolation with the $\L^2$ off-diagonal decay we find that $((1+t^2L_0)^{-1})_{t>0}$ satisfies $\L^q$ off-diagonal estimates of arbitrarily large order. 

Finally, we obtain by a Sobolev embedding for all $t>0$ and all $f \in \H^r \cap \L^2$,
\begin{align*}
 \|\e^{-t^2 L_0} f\|_{\L^{q}} \lesssim  \| \nabla_x \e^{-t^2 L_0} f\|_{\H^r} \lesssim t^{-1} \|f\|_{\H^r},
\end{align*}
which is $\H^r- \L^{q}$-boundedness of the heat semigroup. Since $q$ is an interior point of  $\cJ(L_0)$, we also have $\L^q$-boundedness of the heat semigroup from Step~1 and hence we obtain $\H^\rho - \L^q$-boundedness by interpolation. This being said, we can take again $\L^q$-norms in \eqref{eq1: heat semigroup interval} and conclude the missing $\H^\rho-\L^q$-boundedness of $((1+t^2 L_0)^{-1})_{t>0}$ as in \eqref{eq2: heat semigroup interval}.
\end{proof}
\subsection{Identification of the critical Poisson numbers}
\label{subsec: identification of the critical Poisson numbers}

We present the proof for the Poisson semigroup \emph{vis-à-vis} and focus on where the argument gets technically more involved.

\begin{proof}[Proof of Theorem~\ref{thm: Poisson semigroup interval}]
We break the argument again in three steps.
\medskip

\noindent \emph{Step~1: From the resolvent to the semigroup}. We get $(p_-(L), p_+(L)) \subseteq \cJ^{\mathrm{Pois}}(L)$ and $(q_-(L), q_+(L)) \subseteq \cN^{\mathrm{Pois}}(L)$ by repeating the argument for the heat semigroup \emph{mutadis mutandis}.


\medskip

\noindent \emph{Step~2: From $\cJ^{\mathrm{Pois}}(L)$ to $\cJ(L)$}. We shall always get from Poisson semigroup bounds to resolvents of $L^{1/2}$ on the imaginary axis and then to resolvents of $L$ via the decomposition
\begin{align*}
	(1+t^2 L)^{-1} = (1- \i t L^{1/2})^{-1} (1+ \i t L^{1/2})^{-1} \quad (t>0).
\end{align*}
As a substitute for \eqref{eq1: heat semigroup interval} we need Laplace transform formul\ae \, on the imaginary axis that we are going to derive next. 

Let $\eps \in (0, \nicefrac{(\pi- \omega_L)}{4})$ and $t>0$. Since $L^{1/2}$ is sectorial of angle $\nicefrac{\omega_L}{2}$, the operator 
\begin{align*}
T \coloneqq \e^{\i(\eps-\frac{\pi}{2})} + t \e^{\i \eps} L^{1/2} = -\i \e^{\i \eps} (1+ \i t L^{1/2})
\end{align*}
is invertible and sectorial of angle $\nicefrac{\pi}{2}-\eps$. By the Calderón reproducing formula we have for $f \in \L^2$ as an improper Riemann integral,
\begin{align*}
f = \int_0^\infty T \e^{-sT}f \, \d s.
\end{align*}
Applying $T^{-1}$ on both sides gives the formula
\begin{align}
\label{eq1: Poisson semigroup interval}
(1+ \i t L^{1/2})^{-1}f 
&= -\i \e^{\i \eps} \int_0^\infty \e^{\i s \e^{\i \eps}} \e^{-st \e^{\i \eps} L^{1/2}}f \, \d s.
\end{align}
The latter integral converges absolutely in $\L^2$ since by the functional calculus on $\L^2$ the Poisson semigroup is uniformly bounded on $\e^{\i \eps} \R^+$ and $\Re(\i \e^{\i \eps}) = -\sin(\eps) < 0$. A similar formula holds for $(1 - \i t L^{1/2})^{-1}f$ upon replacing $\i$ by $-\i$ at each occurrence.

Let now $r \in \cJ^{\mathrm{Pois}}(L)$ and take any $p$ between $r$ and $2$. We shall show that $((1+t^2 L)^{-1})_{t>0}$ is $a^{-1} \H^p$-bounded, that is, $p \in \cJ(L)$. Then, together with Step~1, $p_\pm^{\mathrm{Pois}}(L) = p_\pm(L)$ follows.

\medskip

\noindent \emph{Step~2a: The Lebesgue case $p>1$}. Interpolation (Lemma~\ref{lem: OD extrapolation to sectors}) of the $\L^r$-bound for $(0,\infty)$ and the $\L^2$-bound on some sector  provides us with a smaller $\eps > 0$ such that $\e^{-z L^{1/2}}$ is $\L^p$-bounded for $z \in \cl{\S_\eps^+}$. Hence, the integral on the right-hand side in \eqref{eq1: Poisson semigroup interval} converges absolutely in $\L^p$ for all $t>0$ and all $f \in \L^p \cap \L^2$ and we obtain
\begin{align*}
\|(1+ \i t L^{1/2})^{-1}  f\|_p \lesssim \|f\|_p.
\end{align*}
The same argument applies to $(1 - \i t L^{1/2})^{-1}$ and $\L^p$-boundedness of $(1+t^2 L)^{-1}$ follows by composition.

\medskip

\noindent \emph{Step~2b: The Hardy case $p \leq 1$}. As in the case of the heat semigroup we appeal to Lemma~\ref{lem: OD Hardy implies boundedness} and use exponents $\varrho \in (r,p)$ and $q \in (1,2)$ with $\nicefrac{n}{\varrho}-\nicefrac{n}{q} <1$. In particular, $\varrho, q$ are interior points of $\cJ^{\mathrm{Pois}}(L)$.

The vanishing moments condition and the $\L^q$ off-diagonal estimates of arbitrarily large order for $((1+t^2L)^{-1})_{t>0}$ follow exactly as for the heat semigroup and it remains to show $a^{-1} \H^\varrho - \L^q$-boundedness. As before, we arrive at $a^{-1}\H^\varrho - \L^q$-boundedness for $(\e^{-tL^{1/2}})_{t>0}$ but we need to extend the property to some small sector in order to use \eqref{eq1: Poisson semigroup interval}.
 
Again by Step~2a, we know that there exists a smaller $\eps > 0$ such that the Poisson semigroup $\e^{-z L^{1/2}}$ is $\L^q$-bounded for $z \in \cl{\S_{2 \eps}^+}$. Now, let $z \in \cl{\S_\eps^+}$ and decompose 
\begin{align*}
	z= t + z'  \quad \text{with} \quad t>0, z' \in \cl{\S_{2 \eps}^+},  |z| \simeq |z'| \simeq t.
\end{align*}
By composition, $\e^{-z L^{1/2}} = \e^{-z' L^{1/2}} \e^{-t L^{1/2}}$ is $a^{-1} \H^\varrho - \L^q$-bounded for $z \in \cl{\S_{\eps}^+}$. Taking $\L^q$-norms in \eqref{eq1: Poisson semigroup interval}, we obtain for all $t>0$ and all $f \in \H^\varrho \cap \L^2$,
\begin{align*}
\|(1+ \i t L^{1/2})^{-1} a^{-1}f\|_q
&\lesssim \int_0^\infty \e^{- s \sin(\eps)} (s t)^{\nicefrac{n}{q}-\nicefrac{n}{\varrho}} \|f\|_{\H^\varrho} \, \d s \\
&\lesssim  t^{\nicefrac{n}{q}-\nicefrac{n}{\varrho}} \|f\|_{\H^\varrho},
\end{align*}
where the integral in $s$ is finite by the choice of our exponents. Hence, $((1+ \i t L^{1/2})^{-1})_{t>0}$ is $a^{-1} \H^\varrho - \L^q$-bounded. In Step~2a we have seen that $((1-\i t L^{1/2})^{-1})_{t>0}$ is $\L^q$-bounded. Thus, $((1+t^2 L)^{-1})_{t>0}$ is $a^{-1} \H^\varrho - \L^q$-bounded. This completes Step~2b.

\medskip

\noindent \emph{Step~3: From $\cN^{\mathrm{Pois}}(L)$ to $\cN(L)$}. We cannot work with the representation \eqref{eq1: Poisson semigroup interval}: once the gradient is inside the integral, we would have to deal with a function that behaves like $s^{-1}$ in $\L^2$-norm near $s=0$. For $\eps \in (0, \nicefrac{(\pi- \omega_L)}{4})$, $t>0$, and $T \coloneqq -\i \e^{\i \eps} (1+ \i t L^{1/2})$ as before, we use instead the reproducing formula
\begin{align*}
f = \int_0^\infty s T^2 \e^{-sT}f \, \d s
\end{align*}
for $f \in \L^2$. Applying $T^{-2}$ on both sides, we find the absolutely convergent representation
\begin{align}
\label{eq2: Poisson semigroup interval}
(1+ \i t L^{1/2})^{-2}f 
&= -\e^{2 \i \eps} \int_0^\infty \e^{\i s \e^{\i \eps}}\, s\, \e^{-st \e^{\i \eps} L^{1/2}}f \, \d s
\end{align}
 with an additional factor of $s$. Again, an analogous representation is available for $(1- \i t L^{1/2})^{-2}f$.

Let now $r \in \cN^{\mathrm{Pois}}(L)$ and take any $p$ between $r$ and $2$. We shall show that $(t \nabla_x (1+t^2 L)^{-1})_{t>0}$ is $a^{-1} \H^p-\H^p$-bounded, that is, $p \in \cN(L)$. Then, together with Step~1, $q_\pm^{\mathrm{Pois}}(L) = q_\pm(L)$ follows.

In contrast to the proof for the heat semigroup we also need to distinguish the case $p > 2$ from the rest.

\medskip

\noindent \emph{Step~3a: The case $1 <p \leq 2$}. We can further assume $p < 2_*$ (and hence $n \geq 3$), since otherwise we can directly conclude by Proposition~\ref{prop: J(L) contains neighborhood of Sobolev conjugates}. 

We claim that for every $q \in [p,2]$ there exists a smaller $\eps > 0$ such that the following boundedness properties hold for all $z \in \cl{\S_\eps^+}$: 
\begin{align}
\label{Step3a: prop1}
&\L^{q} - \L^{q} \text{ for } z \nabla_x \e^{-z L^{1/2}} \\
\label{Step3a: prop2}
&\L^q - \L^{q^*} \text{ for } \e^{-z L^{1/2}}.
\end{align}

For \eqref{Step3a: prop1} we use interpolation between the $\L^r$-result on $(0,\infty)$ and the $\L^2$-result on some sector. As for \eqref{Step3a: prop2}, we use the assumption and a Sobolev embedding to give
\begin{align*}
\|\e^{-z L^{1/2}}f\|_{r^*} \leq \|\nabla_x \e^{-z L^{1/2}}f\|_r \lesssim |z|^{-1} \|f\|_r
\end{align*}
for all $z \in (0,\infty)$ and all $f \in \L^r \cap \L^2$. This means $\L^r - \L^{r^*}$-boundedness. The same argument works for $z$ in a sector if we replace the exponent $r$ by $2$ and we can conclude by interpolation as before.

We use \eqref{Step3a: prop1} for $q = p^*$. This choice is admissible since we assume $p<2_*$. Applying $t\nabla_x$ to \eqref{eq2: Poisson semigroup interval} and taking $\L^{p^*}$-norms, we obtain for all $t>0$ and all $f \in \L^{p^*} \cap \L^2$,
\begin{align*}
\|t \nabla_x (1+ \i t L^{1/2})^{-2}f \|_{p^*} 
\lesssim \int_0^\infty \e^{- s \sin(\eps)} \|f\|_{p^*} \, \d s.
\end{align*}
Hence, $(t \nabla_x (1+ \i t L^{1/2})^{-2})_{t>0}$ is $\L^{p^*}$-bounded. In the same manner, \eqref{Step3a: prop2} for $q=p$ implies that $((1- \i t L^{1/2})^{-2})_{t>0}$ is $\L^p - \L^{p^*}$-bounded. By composition, $(t\nabla_x(1+t^2L)^{-2})_{t>0}$ is $\L^p - \L^{p^*}$-bounded. 

Since this works for all $p \in (r, 2_*)$, we get $\L^p$-boundedness of $(t\nabla_x(1+t^2L)^{-2})_{t>0}$. Indeed, it suffices to interpolate with the $\L^2$ off-diagonal estimates and then use Lemma~\ref{lem: OD implies boundedness}. But then we can apply Lemma~\ref{lem: reduce resolvent powers} in order to get $\L^p$-boundedness of $(t\nabla_x(1+t^2L)^{-1})_{t>0}$ as required.

\medskip

\noindent \emph{Step~3b: The case $1_* <p \leq 1$}. As in Step~3b for the heat semigroup we see that it suffices to prove $a^{-1}\H^p$-boundedness of $((1+t^2 L)^{-1})_{t>0}$ and in doing so, we can assume $p<2_*$. As usual, we rely on Lemma~\ref{lem: OD Hardy implies boundedness}. We fix any $\varrho \in (r,p)$ and let $q \coloneqq r^* \in (1,2)$. In particular, $\nicefrac{n}{\varrho}-\nicefrac{n}{q} <1$ and $\varrho, q$ are interior points of $\cN^{\mathrm{Pois}}(L_0)$.

Repeating the argument from Step~3b for the heat semigroup \emph{mutadis mutandis}, we get the vanishing moments condition and the $\L^q$ off-diagonal estimates of arbitrarily large order for $((1+t^2L)^{-1})_{t>0}$ and we get $a^{-1} \H^\varrho - \L^q$-boundedness of the Poisson semigroup $(\e^{-t L^{1/2}})_{t>0}$. From Step~3a we know that $q$ is an interior point of $\cN(L)$, hence of $\cJ(L)$. Theorem~\ref{thm: main result Hoo} yields $\L^q$-boundedness of $(\e^{-zL^{1/2}})_{z \in \cl{\S_\eps^+}}$ for any admissible $\eps>0$. Consequently, we are back in the situation of Step~2b of the ongoing proof and obtain the missing  $a^{-1} \H^\varrho - \L^q$-boundedness of $((1+t^2L)^{-1})_{t>0}$. 

\medskip

\noindent \emph{Step~3c: The case $2 < p < \infty$}. We claim that there exists a smaller $\eps > 0$ and an exponent $q \in (1,p]$ with $\nicefrac{n}{q}-\nicefrac{n}{p} < 1$ such that the following boundedness properties hold for $z \in \cl{\S_\eps^+}$: 
\begin{align}
\label{Step3c: prop1}
\L^{p} - \L^{p} \text{ for } z \nabla_x \e^{-z L^{1/2}},\\
\label{Step3c: prop2}
\L^q - \L^{p} \text{ for } \e^{-z L^{1/2}}.
\end{align}

The first part follows by interpolation between the $\L^r$ and the $\L^2$-result. For the second part we first note that $2^* \leq p^+(L) \leq p_+^{\mathrm{Pois}}(L)$ and $q_{-}^{\mathrm{Pois}}(L) \leq q_{-}(L) < 2_*$ by Step~1 and Proposition~\ref{prop: J(L) contains neighborhood of Sobolev conjugates}. In dimensions $n \leq 2$ we have $2^* = \infty$, hence $2< p< r <  p_+^{\mathrm{Pois}}(L)$. We take $q \coloneqq p$ and obtain the claim by interpolation between the $\L^r$-result on $(0,\infty)$ and the $\L^2$-result on a sector. In dimension $n \geq 3$, we have $r_* \in (2_*,r) \subseteq \cN^{\mathrm{Pois}}(L)$ and we obtain $\L^{r_*} - \L^r$-boundedness on $(0,\infty)$ by the Sobolev embedding as in Step~3a. Now, \eqref{Step3c: prop2} follows by interpolation with the $\L^2$-boundedness on a sector for the choice $q \coloneqq [2,r_*]_\theta$ given that $p = [2,r]_\theta$. Note that $\nicefrac{n}{q}-\nicefrac{n}{p} = \theta < 1$.

Equipped with \eqref{Step3c: prop1} and \eqref{Step3c: prop2}, we can take $\L^p$-norms in \eqref{eq2: Poisson semigroup interval} after having applied the gradient, as well as in the analogous formula for $(1-\i t L^{1/2})^{-2}$. We obtain $\L^p$-boundedness of $(t \nabla_x (1+ \i t L^{1/2})^{-2})_{t>0}$ and $\L^q-\L^p$-boundedness of $((1- \i t L^{1/2})^{-2})_{t>0}$. In the second case the restriction on $q$ guarantees again that the integral in $s$ converges. Hence, $(t\nabla_x (1+ t^2 L)^{-2})_{t>0}$ is $\L^q - \L^p$-bounded. 

At this point we can repeat the argument in the last paragraph of Step~3a to conclude $\L^p$-boundedness.
\end{proof}
\subsection{More on off-diagonal decay for the Poisson semigroup}
\label{subsec: More on off-diagonal decay for the Poisson semigroup}

We include an exemplary result to illustrate the poor off-diagonal decay of the Poisson semigroup. In general, and in stark contrast to the resolvents, there is not enough decay to bridge between $\L^q - \L^2$-estimates and $\L^q - \L^q$-estimates via Lemma~\ref{lem: OD implies boundedness}.

\begin{prop}
\label{prop:sgdecay} 
If $(p_-(L)\vee 1)<q \leq 2$, then $(t{L^{1/2}} \e^{-t {L^{1/2}}})_{t>0}$ satisfies $\L^q - \L^2$ off-diagonal estimates of order $\nicefrac{n}{q}-\nicefrac{n}{2} + 1$.
\end{prop}

\begin{proof}
We pick $p \in (p_{-}(L) \vee 1, q)$ and let $\theta \in (0,1)$ be such that $q = [p,2]_\theta$. For a parameter $\alpha > 1$, to be chosen later on, we consider the family
\begin{align*}
	t^\alpha L^{\alpha/2} \e^{- t {L^{1/2}}} = t^\alpha L^{\alpha/2} \e^{- \frac{t}{2} {L^{1/2}}} \e^{- \frac{t}{2} {L^{1/2}}} \quad (t>0).
\end{align*}
From the left-hand side and Lemma~\ref{lem: functional calculus bounds from J(L) abstract}.(i) we obtain $\L^2$ off-diagonal estimates of order $\alpha$, whereas from the right-hand side and Proposition~\ref{prop: Poisson bounds} we obtain $\L^p - \L^2$-boundedness. This implies $\L^q - \L^2$ off-diagonal estimates of order $\theta \alpha$, see Lemma~\ref{lem: OD interpolation}. 

Now, let $E, F \subseteq \R^n$ be measurable, $f \in \L^q \cap \L^2$ and $t>0$. We use the Calderón reproducing formula 
\begin{align*}
f = c_\alpha \int_0^\infty s^{\alpha -1} L^{\alpha/2-1/2} \e^{-s L^{1/2}} \, \frac{\d s}{s}
\end{align*}
in order to give
\begin{align*}
	\ind_F \big(t {L^{1/2}} &\e^{-t {L^{1/2}}}\big) \ind_E f \\
	&= c_\alpha \int_0^\infty \frac{t s^{\alpha - 1}}{(s+t)^\alpha} \ind_F \big((s+t)^{\alpha} L^{\alpha/2} \e^{-(s+t){L^{1/2}}}\big) \ind_E f \, \frac{\d s}{s}.
\end{align*}
Thus, setting $\gamma \coloneqq \nicefrac{n}{q} - \nicefrac{n}{2} \geq 0$, we get
\begin{align*}
	\|\ind_F &(t {L^{1/2}} \e^{-t {L^{1/2}}}) \ind_E f\|_2 \\
	&\lesssim \|f\|_q \int_0^\infty \frac{t s^{\alpha - 1}}{(s+t)^{\alpha+\gamma}} \bigg(1+ \frac{\dist(E,F)}{s+t} \bigg)^{-\theta \alpha} \, \frac{\d s}{s} \\
	&= \|f\|_q t^{-\gamma} \int_0^\infty \frac{\sigma^{\alpha-1}}{(1+\sigma)^{\alpha + \gamma}} \bigg(1+ \frac{\dist(E,F)/t}{1+\sigma} \bigg)^{-\theta \alpha} \, \frac{\d \sigma}{\sigma}.
\end{align*}
We let $X \coloneqq \nicefrac{\dist(E,F)}{t}$. It remains to show that we can choose $\alpha > 1$ in such a way that with an implicit constant independent of $X$,
\begin{align*}
	\int_0^\infty \frac{\sigma^{\alpha-1}}{(1+\sigma)^{\alpha + \gamma}} \bigg(1+ \frac{X}{1+\sigma} \bigg)^{-\theta \alpha} \, \frac{\d \sigma}{\sigma} \lesssim (1+X)^{-\gamma- 1}.
\end{align*}

In the case $X \leq 1$, we simply bound the left-hand side by
\begin{align*}
	\int_0^\infty \frac{\sigma^{\alpha-1}}{(1+\sigma)^{\alpha + \gamma}} \, \frac{\d \sigma}{\sigma} \lesssim 1 \lesssim (1+X)^{- \gamma - 1}.
\end{align*}
In the case $X > 1$, we split the integral into three pieces and obtain a bound (up to a multiplicative constant depending on $\alpha, \gamma, \theta$) by
\begin{align*}
	\int_0^1 \sigma^{\alpha-1} X^{-\theta \alpha} \, &\frac{\d \sigma}{\sigma} 
	+ \int_1^X \frac{\sigma^{\alpha-1}}{\sigma^{\alpha + \gamma}} \bigg(\frac{X}{\sigma} \bigg)^{-\theta \alpha} \, \frac{\d \sigma}{\sigma}
	+ \int_X^\infty \frac{\sigma^{\alpha-1}}{\sigma^{\alpha + \gamma}} \, \frac{\d \sigma}{\sigma}  \\
 	&\lesssim X^{-\theta \alpha}  + X^{-\theta \alpha}(X^{\theta \alpha - \gamma -1} + 1) + X^{-\gamma-1} \\
 	&\lesssim (1+X)^{-\gamma -1},
\end{align*}
provided that we pick $\alpha > \nicefrac{(\gamma + 1)}{\theta} \vee 1$.
\end{proof}
\section{\texorpdfstring{$\L^p$}{Lp} boundedness of the Hodge projector}
\label{sec: Hodge}

\noindent In this chapter, we discuss $\L^p$-boundedness of the Hodge projector associated to $L_0$ (that is, $L$ in the case when $a = 1$). We obtain a characterization of the range for $p$ in terms of critical numbers. 

Let $p\in (1,\infty)$. The well-known \emph{Leray--Helmholtz decomposition} states that every vector field $f \in \L^p(\R^n; \IC^{mn})$ can be decomposed into a divergence-free part and a gradient field. In order to set the stage for studying operator-adapted counterparts, it will be convenient to reproduce the simple proof. 

\begin{defn}
\label{def: Npdiv and Rpgrad}
For $p\in (1,\infty)$ let 
\begin{align*}
\nul_p(\div_x) &\coloneqq \{g \in \L^p(\R^n; \IC^{nm}): \div_x g = 0\}, \\
\ran_p(\gradx) &\coloneqq \{\nabla_x h: h \in \W^{1,p}(\R^n; \IC^m)\}.
\end{align*}
\end{defn}

\begin{lem}[Leray--Helmholtz decomposition\index{Leray--Helmholtz decomposition}]
\label{lem: Leray--Helmholtz}
Let $p \in (1,\infty)$. There is a topological decomposition 
\begin{align*}
 \L^p(\R^n; \IC^{nm}) = \nul_p(\div_x) \oplus \cl{\ran_p(\gradx)}
\end{align*} 
and the projection onto $\cl{\ran_p(\nabla_x)}$ is given by the $\L^p$-bounded Fourier multiplication operator $-\nabla_x (-\Delta_x^{-1}) \div_x$.
\end{lem}

\begin{proof}
The Fourier symbol $\xi \mapsto |\xi|^{-2} \xi \otimes \xi$ of $-\nabla_x (-\Delta_x^{-1}) \div_x$ is homogeneous of degree $0$ and hence fits into the scope of the Mihlin multiplier theorem. Hence, this operator is defined on $\cZ'$ and restricts to bounded map on $\L^p$ that we call $\IP_p$. As $\IP_p$ is a projection on $\L^p$, it induces the topological decomposition $\L^p = \ran(1 - \IP_p) \oplus \ran(\IP_p)$. By construction, we have $\ran(1 - \IP_p) \subseteq \nul_p(\div_x)$ and $\ran(\IP_p) \subseteq \cl{\ran_p(\nabla_x)}$. Equality in both inclusions follows provided that $\nul_p(\div_x) \cap  \cl{\ran_p(\nabla_x)} = \{0\}$. But if $f$ belongs to this intersection, then $f = \nabla_x h$, where $h \in \Wdot^{1,p}$ satisfies $\Delta_x h = 0$ in $\cZ'$. Therefore $h= 0$ in $\cZ'$, so $h$ must be a polynomial and hence a constant, which in turn means that $f=0$.
\end{proof}

In view of the explicit formula for the projection in Lemma~\ref{lem: Leray--Helmholtz}, the Leray--Helmholtz decomposition is also called \emph{Hodge decomposition}\index{Hodge decomposition} associated with $-\Delta_x$. Following \cite[Sec.~4.5]{A}, we look for similar decompositions adapted to divergence form operators $-\div_x d \nabla_x$. These operators are defined in the sense of distributions modulo constants as bounded operators
\begin{align}
\label{eq: divdgard}
-\div_x d \nabla_x: \Wdot^{1,p}(\R^n; \IC^m) \to \Wdot^{-1,p}(\R^n; \IC^m)
\end{align}
for every $p \in (1,\infty)$. Their action is consistent for different values of $p$ and for $p=2$ we find the operator $\eo$ defined in \eqref{eq: Lax-Milgram operator}. The adjoint to \eqref{eq: divdgard} is given by
\begin{align*}
-\div_x d^* \nabla_x: \Wdot^{1,p'}(\R^n; \IC^m) \to \Wdot^{-1,p'}(\R^n; \IC^m).
\end{align*}
When $p= 2$, it corresponds to the operator $L_{0}^\sharp=L_{0}^*$ in the same way that $-\div_x d \nabla_x$ corresponds to $L_{0}$.

The interval that we are mainly interested in concerns the  bounded extension to $\L^p$ of the $\L^2$-bounded \emph{Hodge projector}\index{Hodge projector!adapted to $\eo$} $\nabla_x \eo^{-1} \div_x$. 
\begin{defn}
\label{def: Hodge interval}
Introduce the interval\index{P@$\cP(L_{0})$}
\begin{align*}
\cP(L_{0}) \coloneqq \big \{p \in (1,\infty): \nabla_x \eo^{-1} \div_x \text{ is $\L^p$-bounded} \big \}.
\end{align*}
\end{defn}
\emph{A priori}, there are two possibilities to incorporate the matrix $d$ into the Leray--Helmholtz decomposition:\index{Hodge decomposition!$d$-adapted}
\begin{align}
\label{eq: Hodge}
\L^p(\R^n; \IC^{nm}) &= \nul_p(\div_x) \oplus \cl{d \ran_p(\gradx)}
\intertext{and}
\label{eq: Hodge-tilde}
\L^p(\R^n; \IC^{nm}) &= \nul_p(\div_x d) \oplus \cl{\ran_p(\gradx)},
\end{align}
where closures are taken in $\L^p$ and $\nul_p(\div_x d) \coloneqq \{f \in \L^p(\R^n; \IC^{nm}) : \div_x(df) = 0\}$. We shall see that these \emph{topological} decompositions always hold when $p=2$ and that this directly relates to \eqref{eq: divdgard} being an isomorphism for $p=2$. We say that such a topological decomposition \emph{compatibly holds}\index{Hodge decomposition!compatible} if in addition for every $f \in \L^p \cap \L^2$ the decomposition in $\L^p$ is the same as in $\L^2$.

Compatibility with the theory for $p=2$ is a key issue here and we take the occasion to clarify some points that had been left unclear in the literature. The central question is whether the set of $p \in (1,\infty)$ for which $-\div_x d \nabla_x: \Wdot^{1,p} \to \Wdot^{-1,p}$ is an isomorphism is an open interval.  While openness turns out to be true in general,   connectedness requires more specific arguments.
	
As a cautionary tale, let us remark that in general the compatibility of the inverse\index{compatibility of the inverse} does not come for free and hence the property of being an isomorphism does not interpolate. To give a simple example, consider the dilation $f \mapsto (t \mapsto f(\frac{t}{2}))$ on the real line. Its restriction $T_p$ to $\L^p(\R)$ is invertible and $\|T_p\|_{p\to p} =2^{1/p} = \|T_p^{-1}\|_{p \to p}^{-1}$. Hence, the spectrum $\sigma(T_p)$ is contained in the circle of radius $2^{1/p}$. Now, pick $\lambda \in \sigma(T_3)$. Then $\lambda-T$ is invertible on $\L^{2}(\R)$ and $\L^{4}(\R)$ but not on $\L^3(\R)$ and therefore the inverses cannot be compatible.

Concerning the isomorphism property for $-\div_x d \nabla_x$, the formulation in \cite[Cor.~4.24]{A} is ambiguous. As far as Hodge decompositions are concerned, a general statement in \cite[Prop.~2.17]{FMcP} asserts (when restricted to our setup) that the set of exponents for which they are valid is an interval, but their proof offers no specific argument. In view of our discussion below, connectedness should still be considered unproved at this stage and compatible invertibility and compatible Hodge decompositions only hold in the connected component that contains $p=2$. The fact that this connected component enters the discussion has previously been noticed in \cite[Section~3]{AusStaRemarks}.
\subsection{Compatible adapted Hodge decompositions}
\label{subsec: compatible Hodge decomposition}

The following discussion extends and streamlines the presentation in \cite[Sec.~4.5]{A}.

\begin{lem}
\label{lem: Hodge interval via isomorphism}
Let $p \in (1,\infty)$. Then $p \in \cP(L_{0})$ if and only if $\eo$ extends by density from $\Wdot^{1,p} \cap \Wdot^{1,2}$ to an isomorphism $\Wdot^{1,p} \to \Wdot^{-1,p}$ whose inverse agrees with $\eo^{-1}$ on $\Wdot^{-1,p} \cap \Wdot^{-1,2}$. In particular, $\cP(L_{0})$ is an open set.
\end{lem}

\begin{proof}
We need some preliminary observations on the Leray--Helmholtz decompositions of $\L^p$ and $\L^2$ in Lemma~\ref{lem: Leray--Helmholtz}. As they are being achieved through projections that coincide on the dense subset $\L^p \cap \L^2$, we also have a direct decomposition
\begin{align}
\label{eq1: Hodge interval via isomorphism}
\L^p \cap \L^2 = \big(\nul_p(\div_x) \cap \nul_2(\div_x)\big) \oplus \big(\cl{\ran_p(\nabla_x)} \cap \cl{\ran_2(\nabla_x)}\big)
\end{align}
that is topological with respect to $\L^p$ and $\L^2$-norms. Moreover, the subspaces on the right are dense in $\nul_p(\div_x)$ and $\cl{\ran_p(\nabla_x)}$ for the $\L^p$-norm, respectively. Now,
\begin{align}
\label{eq2: Hodge interval via isomorphism}
\nabla_x: \Wdot^{1,p} \cap \Wdot^{1,2} \to \cl{\ran_p(\nabla_x)} \cap \cl{\ran_2(\nabla_x)}
\end{align}
is bijective and bounded from above and below for the respective $p$-norms. The same is true for
\begin{align}
\label{eq3: Hodge interval via isomorphism}
\div_x: \cl{\ran_p(\nabla_x)} \cap \cl{\ran_2(\nabla_x)} \to \Wdot^{-1,p} \cap \Wdot^{-1,2}.
\end{align}
Indeed, the upper bound follows right away, injectivity is due to \eqref{eq1: Hodge interval via isomorphism} and surjectivity and the lower bound follow since $\nabla_x \Delta_x^{-1}$ is an explicit right inverse.

We turn to the actual proof. Since $\eo : \Wdot^{1,p} \cap \Wdot^{1,2} \to \Wdot^{-1,p} \cap \Wdot^{-1,2}$ is well-defined and bounded for the $p$-norms, it follows that it extends to an isomorphism as claimed precisely if $\eo ^{-1}: \Wdot^{-1,p} \cap \Wdot^{-1,2} \to \Wdot^{1,p} \cap \Wdot^{1,2}$ is well-defined and bounded for the $p$-norms. Composition with the maps in \eqref{eq2: Hodge interval via isomorphism} and \eqref{eq3: Hodge interval via isomorphism} yields equivalence to well-definedness and boundedness in $p$-norm for
\begin{align*}
\nabla_x \eo^{-1} \div_x : \cl{\ran_p(\nabla_x)} \cap \cl{\ran_2(\nabla_x)} \to \cl{\ran_p(\nabla_x)} \cap \cl{\ran_2(\nabla_x)}.
\end{align*}
Due to \eqref{eq1: Hodge interval via isomorphism} this is the same as saying $p \in \cP(L_{0})$.

Finally, the set of exponents $p \in (1,\infty)$ with the isomorphism property for $\eo$ with compatible inverse is open in $(1,\infty)$ thanks to \v{S}ne\u{i}berg's stability theorem, using that the scales $(\Wdot^{1,p})_{p \in (1,\infty)}$ and $(\Wdot^{-1,p})_{p \in (1,\infty)}$ interpolate by the complex method. See for instance \cite{Sneiberg, ABES3} and also \cite[Thm.~8.1]{KMM} for the compatibility.\index{Theorem!\v{S}ne\u{i}berg's}
\end{proof}

\begin{lem}
\label{lem: Hodge Lp sufficient}
If $p \in \cP(L_{0})$, then the Hodge decompositions \eqref{eq: Hodge} and \eqref{eq: Hodge-tilde} compatibly hold. The projections onto $\cl{d \ran_p(\gradx)}$ and $\cl{\ran_p(\gradx)}$ are the extensions (by density) of $-d\nabla_x \eo^{-1} \div_x$ and $-\nabla_x \eo^{-1} \div_x d$, respectively.
\end{lem}

\begin{proof}
On $\L^2(\R^n; \IC^{nm})$ we consider the bounded projection operators
\begin{align*}
\IP_2 &\coloneqq -d\nabla_x \eo^{-1} \div_x, \\
\wt{\IP}_2 & \coloneqq - \nabla_x \eo^{-1} \div_x d.
\end{align*}
They are $\L^p$-bounded since we assume $p \in \cP(L_{0})$. We call $\IP_p$ and $\wt{\IP}_p$ their extensions by density from $\L^p \cap \L^2$ to bounded projections on $\L^p$, which induce the topological decompositions
\begin{align}
\label{eq1: Hodge Lp sufficient}
\L^p = \ran(1- \IP_p) \oplus \ran(\IP_p), \quad \L^p = \ran(1- \wt{\IP}_p) \oplus \ran(\wt{\IP}_p).
\end{align}
By construction, we have 
\begin{align*}
\ran(\IP_p) \subseteq \cl{d \ran_p(\gradx)}, \quad \ran(\wt{\IP}_p) \subseteq \cl{\ran_p(\gradx)}
\end{align*}
and from $\div_x \IP_p f = \div_x f$ and $\div_x (d \wt{\IP}_p f) = \div_x(df)$ for $f \in \L^p \cap \L^2$ we also conclude
\begin{align*}
\ran(1-\IP_p) \subseteq \nul_p(\div_x), \quad \ran(1- \wt{\IP}_p) \subseteq \nul_p(\div_x d).
\end{align*}
It remains to establish equality in all four inclusions and owing to \eqref{eq1: Hodge Lp sufficient} we only have to show that $\nul_p(\div_x) \cap \cl{d \ran_p(\gradx)}$ and $\nul_p(\div_x d) \cap \cl{\ran_p(\gradx)}$ are trivial.

Let $f \in \nul_p(\div_x) \cap \cl{d \ran_p(\gradx)}$. By density, we find $h_j \in \W^{1,p} \cap \W^{1,2}$ such that $d \nabla_x h_j \to f$ in $\L^p$ as $j \to \infty$. Then $\div_x (d \nabla_x h_j) \to 0$ in $\Wdot^{-1,p}$, whereupon Lemma~\ref{lem: Hodge interval via isomorphism} yields $h_j \to 0$ in $\Wdot^{1,p}$. Consequently, we have $f = 0$. 

Likewise, if $f \in \nul_p(\div_x d) \cap \cl{\ran_p(\gradx)}$, then we pick $h_j \in \W^{1,p} \cap \W^{1,2}$ with $\nabla_x h_j \to f$ in $\L^p$ as $j \to \infty$ and conclude $f=0$ as before.
\end{proof}

We shall see momentarily that $p \in \cP(L_{0})$ also entails the following property.

\begin{defn}
\label{def: p-lower bound} 
Let $p \in (1,\infty)$. Then $d$ is said to satisfy \emph{$p$-lower bounds}\index{plowerbounds@$p$-lower bounds!for $d$} if
\begin{align*}
\|d  f\|_p \gtrsim \|f\|_p \quad (f \in \cl{\ran_p(\nabla_x)}).
\end{align*}
\end{defn}

While this is trivially fulfilled for a strictly elliptic matrix (and probably for that reason has not even been mentioned in \cite{A, FMcP}), in the realm of elliptic systems it imposes a structural condition on $d$.

\begin{lem}
\label{lem: p-lower bound d}
If $p \in \cP(L_{0})$, then $d$ satisfies $p$-lower bounds.
\end{lem}

\begin{proof}
By density it suffices to verify the $p$-lower bound for $f = \nabla_x h$ with $h \in \W^{1,p} \cap \W^{1,2}$. Then $d f \in \L^p \cap \L^2$ and
\begin{align*}
-f = \nabla_x (-\div_x d \nabla_x)^{-1} \div_x d \nabla_x h=  (\nabla_x \eo ^{-1} \div_x) d f.
\end{align*}
The assumption $p \in \cP(L_{0})$ implies $\|f\|_p \lesssim \|d f\|_p$.
\end{proof}

Altogether, we arrive at the following characterization.

\begin{prop}
\label{prop: P(L) vs Hodge}
Let $p \in (1,\infty)$. The following are equivalent:\index{Hodge decomposition!compatible}
\begin{enumerate}
\item $p \in \cP(L_{0})$.
\item $-\div_x d\nabla_x :\Wdot^{1,p} \to \Wdot^{-1,p}$ is an isomorphism whose inverse agrees with $\eo^{-1}$ on $\Wdot^{-1,p} \cap \Wdot^{-1,2}$.
\item $d$ satisfies $p$-lower bounds and \eqref{eq: Hodge} compatibly holds.
\item $d^*$ satisfies $p'$-lower bounds and \eqref{eq: Hodge-tilde} compatibly holds.
\end{enumerate}
\end{prop}

\begin{proof}
We show the following implications.

\medskip

\noindent $(i) \Longleftrightarrow (ii)$. This is Lemma~\ref{lem: Hodge interval via isomorphism}. 

\medskip

\noindent $(i) \Longrightarrow (iii), (iv)$. The compatible Hodge decompositions are due to Lemma~\ref{lem: Hodge Lp sufficient} and the $p$-lower bound for $d$ is due to Lemma~\ref{lem: p-lower bound d}. Moreover, we have $p' \in \cP(L_{0}^*)$ by duality and Lemma~\ref{lem: p-lower bound d} yields the $p'$-lower bound for $d^*$.

\medskip

\noindent $(iii) \Longrightarrow (i)$. We have $2 \in \cP(L_{0})$ and according to Lemma~\ref{lem: Hodge Lp sufficient} the decomposition holds for $p=2$ in virtue of the projection $-d\nabla_x \eo^{-1} \div_x$. The compatibility of the Hodge decomposition implies that this operator is $\L^p$-bounded. Using the $p$-lower bounds, we obtain for all $f \in \L^p \cap \L^2$,
\begin{align*}
\|\nabla_x \eo^{-1} \div_x f\|_p \lesssim \|d \nabla_x \eo^{-1} \div_x f\|_p \lesssim \|f\|_p.
\end{align*}

\medskip

\noindent $(iv) \Longrightarrow (i)$. As in the previous step, we get that $-\nabla_x \eo^{-1} \div_x d$ is $\L^p$-bounded. By duality, $-d^*\nabla_x (\eo^*)^{-1} \div_x$ is $\L^{p'}$-bounded and the $p'$-lower bound implies $p' \in \cP(L_{0}^*)$. Again by duality, $p \in \cP(L_{0})$ follows.
\end{proof}
\subsection{Adapted Hodge decompositions}
\label{subsec: Hodge decomposition}

We drop the compatibility assumption and ask under which conditions the $d$-adapted Hodge decompositions hold.

\begin{prop}
\label{prop: W1p-regularity vs Hodge}
Let $p \in (1,\infty)$. The following are equivalent:\index{Hodge decomposition!$d$-adapted}
\begin{enumerate}	
	\item $-\div_x d\nabla_x :\Wdot^{1,p} \to \Wdot^{-1,p}$ is an isomorphism.
	\item $d$ satisfies $p$-lower bounds and \eqref{eq: Hodge} holds.
	\item $d^*$ satisfies $p'$-lower bounds and \eqref{eq: Hodge-tilde} holds.
\end{enumerate}
\end{prop}

\begin{rem}
 As (i) is equivalent to the  adjoint statement that  $-\div_x d^*\nabla_x :\Wdot^{1,p'} \to \Wdot^{-1,p'}$ is an isomorphism, we could add  to the list three more items. 
\end{rem}

\begin{proof}
We establish the following implications.

\medskip

\noindent \emph{$(i) \Longrightarrow (ii), (iii)$}. Set $\eo_{p}$ the operator in (i). The Hodge decomposition follows by a \emph{verbatim} repetition of the proof of Lemma~\ref{lem: Hodge Lp sufficient}. In fact, it is even easier,  using the operator $\eo_{p}^{-1}:\Wdot^{-1,p} \to \Wdot^{1,p} $ provided by assuming (i).  We can directly define the bounded projections
\begin{align*}
\IP_p &\coloneqq -d\nabla_x \eo_{p}^{-1} \div_x, \\
\wt{\IP}_p & \coloneqq - \nabla_x \eo_{p}^{-1} \div_x d
\end{align*}
on $\L^p$ and use (i) in place of Lemma~\ref{lem: Hodge interval via isomorphism} in the proof. Likewise, for the $p$-lower bound for $d$ we can repeat the proof of Lemma~\ref{lem: p-lower bound d} with $\eo_{p}^{-1}$ in place of $\eo^{-1}$ and working with $f=\nabla_{x}h$, where $h \in \W^{1,p}$. By duality, (i) also implies that $\eo_{p}^* :\Wdot^{1,p'} \to \Wdot^{-1,p'}$ is an isomorphism and hence the $p'$-lower bound for $d^*$ follows as well.

\medskip

\noindent \emph{$(ii) \Longrightarrow (i)$}. The $p$-lower bound implies $\cl{d \ran_p(\nabla_x)} = d \cl{\ran_p(\nabla_x)}$. Hence, $\nul_{p'}(\div_x d^*)$ annihilates $\cl{d \ran_p(\nabla_x)}$ in the $\L^p-\L^{p'}$-duality. In the same duality, $\cl{\ran_{p'}(\nabla_x)}$ annihilates $\nul_p(\div_x)$. The Hodge decomposition \eqref{eq: Hodge} implies $\nul_{p'}(\div_x d^*) \cap \cl{\ran_{p'}(\nabla_x)} = \{0\}$. As we have $\cl{\ran_{p'}(\nabla_x)} = \{\nabla_x h: h \in \Wdot^{1,p'}\}$, injectivity of
\begin{align}
\label{eq1: W1p-regularity vs Hodge}
\eo_{p}^*=-\div_x d^* \nabla_x : \Wdot^{1,p'} \to \Wdot^{-1,p'}
\end{align}
follows. From $\cl{d \ran_p(\nabla_x)} = d \cl{\ran_p(\nabla_x)}$ and \eqref{eq: Hodge} we also obtain directly the injectivity of
\begin{align}
\label{eq2: W1p-regularity vs Hodge}
\eo_{p}=-\div_x d \nabla_x : \Wdot^{1,p} \to \Wdot^{-1,p}.
\end{align}
Hence, both maps have dense range and they become isomorphisms once we have shown that the first map has closed range. To this end, let $h' \in \Wdot^{1,p'}$ and $F \in \L^p$. We decompose $F= G + d \nabla_x f$ according to \eqref{eq: Hodge} and obtain
\begin{align*}
|\langle \nabla_x h', F\rangle|
&=|\langle \nabla_x h', d \nabla_x f\rangle| \\
&=|\langle d^* \nabla_x h', \nabla_x f\rangle| \\
&\lesssim \|\div_x d^* \nabla_xh'\|_{\Wdot^{-1,p'}} \|\nabla_x f\|_{p} \\
&\lesssim \|\div_x d^* \nabla_xh'\|_{\Wdot^{-1,p'}} \|d \nabla_x f\|_{p} \\
&\lesssim \|\div_x d^* \nabla_xh'\|_{\Wdot^{-1,p'}} \|F\|_{p},
\end{align*}
where the third line is just the identification of $\Wdot^{-1,p'}$ with the dual space of $\Wdot^{1,p}$, the fourth is by the $p$-lower bounds and the fifth uses that the splitting \eqref{eq: Hodge} is topological in (ii). Taking the supremum over all $F$ yields $\|h'\|_{\Wdot^{1,p'}} \lesssim  \|\div_x d^* \nabla_xh'\|_{\Wdot^{-1,p'}}$, which implies closed range in \eqref{eq1: W1p-regularity vs Hodge}.

\medskip

\noindent \emph{$(iii) \Longrightarrow (i)$}. The argument is almost identical to the previous step. This time we get $\cl{d^* \ran_{p'}(\nabla_x)} = d^* \cl{\ran_{p'}(\nabla_x)}$, which annihilates $\nul_p(\div_x d)$. By \eqref{eq: Hodge-tilde} we find $\nul_{p'}(\div_x) \cap d^*\cl{\ran_{p'}(\nabla_x)} = \{0\}$ and therefore the map in \eqref{eq1: W1p-regularity vs Hodge} is injective. Injectivity in \eqref{eq2: W1p-regularity vs Hodge} follows directly from  \eqref{eq: Hodge-tilde}. In order to see that we have closed range in \eqref{eq1: W1p-regularity vs Hodge}, we let $h'$ and $F$ as before and decompose $F = G + \nabla_x f$ according to \eqref{eq: Hodge-tilde}. Then
\begin{align*}
|\langle d^*\nabla_x h', F\rangle|
&=|\langle d^* \nabla_x h', \nabla_x f\rangle| \\
&\lesssim \|\div_x d^* \nabla_xh'\|_{\Wdot^{-1,p'}} \|\nabla_x f\|_{p} \\
&\lesssim \|\div_x d^* \nabla_xh'\|_{\Wdot^{-1,p'}} \|F\|_{p},
\end{align*}
which yields $\|d^* \nabla_x h'\|_{p'} \lesssim  \|\div_x d^* \nabla_xh'\|_{\Wdot^{-1,p'}}$.  Using  the $p'$-lower bounds for $d^*$ leads to
\begin{align*}
\|h\|_{\Wdot^{1,p'}}=\|\nabla_x h'\|_{p'} \lesssim \|d^* \nabla_x h'\|_{p'} \lesssim  \|\div_x d^* \nabla_xh'\|_{\Wdot^{-1,p'}}. & \qedhere
\end{align*}
\end{proof}

A comparison between Proposition~\ref{prop: W1p-regularity vs Hodge} and Proposition~\ref{prop: P(L) vs Hodge} shows that compatibility in one of the Hodge decompositions directly relates to compatibility of the inverse of $-\div_x d \nabla_x$ on $\Wdot^{-1,p}$ with the inverse found by the Lax--Milgram lemma on $\Wdot^{-1,2}$. To the best of our knowledge, the question whether incompatibility of the inverses is possible for the operators $-\div_x d \nabla_x$ is still open. A more illuminating comparison between the two results is as follows.

\begin{lem}
\label{lem: comparison P(L) and W1p-regualrity}
Let $\wt{\cP}(L_{0}) \subseteq (1,\infty)$ be the set of exponents $p$ such that $\eo_{p}=-\div_x d \nabla_x: \Wdot^{1,p} \to \Wdot^{-1,p}$ is an isomorphism. Then $\wt{\cP}(L_{0})$ is open and $\cP(L_{0})$ is its connected component that contains $2$.
\end{lem}

\begin{proof}
All relies on the fact that $(\Wdot^{1,p})_{p \in (1,\infty)}$ and $(\Wdot^{-1,p})_{p \in (1,\infty)}$ interpolate by the complex method and have a universal approximation technique. \v{S}ne\u{i}berg's stability theorem yields that $\wt{\cP}(L_{0})$ is open and that locally the inverses agree on the intersection of their domains of definition. If $p_0,p_1 \in \wt{\cP}(L_{0})$ are such that the inverses agree with the one on $\Wdot^{-1,2}$, then by interpolation of the mapping property for the inverses the same is true for all $p \in (p_0, p_1)$. Hence, the subset of exponents with this property is open, closed and connected in $\wt{\cP}(L_0)$, hence is the connected component that contains 2. In Lemma~\ref{lem: Hodge interval via isomorphism} we have identified it to $\cP(L_{0})$.
\end{proof}

\subsection{Characterizations of \texorpdfstring{$\boldsymbol{\cP(L_{0})}$}{P(L)}}
\label{subsec: characterization of P(L)}

\noindent For equations ($m=1$) it has been asserted in \cite[Cor.~4.24]{A} that $\cP(L_{0})$ coincides with the interval $((q_+^{\mathrm{heat}}(L_{0}^*))', q_+^{\mathrm{heat}}(L_{0}))$, albeit being implicit on questions of compatibility. Given Theorem~\ref{thm: heat semigroup interval}, this is the same interval as
$((q_+(L_{0}^*))', q_+(L_{0}))$. We take the opportunity to give the full argument and make compatibilities explicit.

\begin{thm}
\label{thm: endpoints for Hodge projector}
 $\cP(L_{0}) = ((q_+(L_{0}^*))', q_+(L_{0}))$.\index{Hodge projector!$\L^p$-boundedness}\index{P@$\cP(L)$! characterization of}	
\end{thm}

For the proof, we need a particular Sobolev-type inequality and a factorization of $\eo^{-1}$ via Riesz transforms.

\begin{lem}
\label{lem: Sobolev inequality with aL}
If $q \in \cP(L_{0}) \cap (1^*,n)$, then
\begin{align*}
\|\eo^{-1} g\|_{q^*} \lesssim \|g\|_{q_*} \quad (g \in \L^{q_*} \cap \Wdot^{-1,2}).
\end{align*}
\end{lem}

\begin{proof}
We use that $\cZ$ is dense in $\L^{q_*} \cap \Wdot^{-1,2}$, see Section~\ref{subsec: homogeneous smoothness spaces}. Since $\eo^{-1}: \Wdot^{-1,2} \to \Wdot^{1,2}$ is bounded, we may assume $g \in \cZ$. Hence, $f \coloneqq \nabla_x (\Delta_x)^{-1} g$ is defined in $\cZ$ and we have $g = \div_x f$. From the assumption and Sobolev embeddings, we get
\begin{align*}
\|\eo^{-1} g\|_{q^*}
&\lesssim \|\nabla_x \eo^{-1} \div_x f\|_q \\
&\lesssim \|f\|_q \\
&\lesssim \|\nabla_x^2 (\Delta_x)^{-1} g\|_{q_*} \\
&\lesssim \|g\|_{q_*},
\end{align*}
where the final step is due to the Mihlin multiplier theorem.
\end{proof}

\begin{lem}
\label{lem: factorization of aL via Riesz}
Let $R_{L_{0}}=\nabla_{x} L_{0}^{-1/2}$ and $R_{L_{0}^*}=\nabla_{x} (L_{0}^*)^{-1/2}$ be the bounded Riesz transforms on $\L^2$ associated with $L_{0}$ and $L_{0}^*$, respectively. Then
\begin{align}
\label{eq1: factorization of aL via Riesz}
-R_{L_{0}}  (R_{L_{0}^*})^*  = \nabla_x \eo^{-1} \div_x
\end{align}
as bounded operators on $\L^2$. 
\end{lem}

\begin{proof}
The factorization formally follows but some (tedious) density arguments are necessary to make this precise.

Let $f \in \L^2$. The decomposition \eqref{eq: Hodge} with $p=2$ allows us to write $f = f_0+d f_1$, where $f_0 \in \nul(\div_x)$ and $f_1 \in \cl{\ran(\nabla_x)}$. As usual, our notation indicates kernels and ranges of the operators in $\L^2$ with maximal domain. Since $\ran(R_{L_{0}^*}) \subseteq \cl{\ran(\nabla_x)} = \nul(\div_x)^\perp$ by construction, the left-hand side of \eqref{eq1: factorization of aL via Riesz} sends $f_0$ to $0$. Obviously the same is true for the right-hand side. As for the action on $df_1$, we may assume $f_1 = \nabla_x u$ for $u \in \dom(L_0)$. Indeed, the general case follows by density since $\dom(L_{0})$ is dense in $\dom(L_{0}^{1/2}) = \W^{1,2}$. We obtain $\div_x (d f_1) = - \eo u$, so that
\begin{align*}
\nabla_x \eo^{-1} \div_x (df_1) = - \nabla_x u.
\end{align*}
Moreover, for $g \in \ran(L_0^*)$ we get
\begin{align*}
\langle (R_{L_{0}^*})^* (df_1), g \rangle 
&= \langle d f_1, \nabla_x (L_{0}^*)^{-1/2} g \rangle \\
&=  \langle L_{0}u,  (L_{0}^*)^{-1/2}  g \rangle \\
&= \langle L_{0}^{1/2}u,  g\rangle.
\end{align*}
 Since this holds for all $g$ in a dense subspace of $\L^2$, we first obtain $(R_{L_{0}^*})^* (df_1) =  L_{0}^{1/2}u$ and then $-R_{L_{0}}  (R_{L_{0}^*})^* (df_1) = -\nabla_x u$. Altogether, we have justified \eqref{eq1: factorization of aL via Riesz}.
\end{proof}

\begin{proof}[Proof of Theorem~\ref{thm: endpoints for Hodge projector}]
Recall that in the case of $L_{0}$  the duality relations \eqref{eq1: endpoints for Hodge projector} yield $(1 \vee p_{-}(L_{0}^*))= p_+(L_{0})'$ and  
$(1 \vee p_{-}(L_{0}))' = p_+(L_{0}^*)$. The proof of the theorem is organized in 4 Steps. 

\medskip

\noindent \emph{Step~1: Sufficient condition for $\cP(L_{0})$}. Let $(q_+(L_{0}^*))' < p < q_+(L_{0})$. We demonstrate that $p \in \cP(L_{0})$.

Theorem~\ref{thm: standard relation J(L) and N(L)} yields $q_+(L_{0}) \leq p_+(L_{0})$ and $q_+(L_{0}^*) \leq p_+(L_{0}^*)$. Hence, we obtain $p_-(L_{0}) < p < q_+(L_{0})$ and $p_-(L_{0}^*) < p' < q_+(L_{0}^*)$. Theorem~\ref{thm: Riesz} yields that $R_{L_{0}}$ is $\L^p$-bounded and that $R_{L_{0}^*}$ is $\L^{p'}$-bounded. By composition and duality $R_{L_{0}} (R_{L_{0}^*})^*$ is $\L^p$-bounded and the previous lemma yields the claim.

\medskip

\noindent \emph{Step~2: Necessary condition for $\cP(L_{0}) \cap (2,\infty)$}. We let $p \in \cP(L_{0}) \cap (2,\infty)$ and prove that $p \leq q_+(L_{0})$. 

To begin with, we claim that
\begin{align}
\label{eq2: endpoints for Hodge projector}
[2_*,p) \subseteq \cJ(L_{0}).
\end{align}
Thanks to Proposition~\ref{prop: J(L) contains neighborhood of Sobolev conjugates} there is nothing to do if $p \leq 2^*$. Hence, we may assume $p>2^*$ (and therefore $n \geq 3$ implicitly).  We set $p_0 \coloneqq p$, define iteratively $p_k \coloneqq (p_{k-1})_{**} \coloneqq ((p_{k-1})_*)_*$ and stop at the first exponent $k^- \geq 0$ with $p_{k^-} \in [2_*, 2^*)$. Again by Proposition~\ref{prop: J(L) contains neighborhood of Sobolev conjugates} we have $[2_*, p_{k^-}) \subseteq \cJ(L_0)$. Now, suppose $[2_*,p_k) \subseteq \cJ(L_{0})$ and pick any $\varrho \in (p_k \vee 2^*, p_{k-1})$. Let $f \in \L^{\varrho_{**}} \cap \L^2$. The function
\begin{align*}
g \coloneqq \eo (1+t^2L_{0})^{-1}f = L_{0}(1+t^2L_{0})^{-1}f =t^{-2}  (1- (1+t^2L_{0})^{-1})f
\end{align*}
belongs to $\Wdot^{-1,2}$ since it is contained in the range of $\eo$ and it belongs to $\L^{\varrho_{**}}$ since we have $\varrho_{**} \in (2_*, p_k) \subseteq \cJ(L_{0})$ by assumption. We also have $\varrho_* \in (2,p \wedge n) \subseteq \cP(L_{0}) \cap (1^*,n)$, so we can apply Lemma~\ref{lem: Sobolev inequality with aL} with $q = \varrho_*$ in order to obtain
\begin{align*}
\|(1+t^2 L_{0})^{-1}f\|_{\varrho}
= \|\eo^{-1}g\|_{\varrho}
\lesssim \|g\|_{\varrho_{**}}
\lesssim t^{-2} \|f\|_{\varrho_{**}}.
\end{align*}
This means that the resolvents of $L_{0}$ are $\L^{\varrho_{**}} - \L^\varrho$-bounded. Since $\varrho \in (p_k \vee 2^*, p_{k-1})$ was arbitrary, interpolation with the $\L^2$ off-diagonal estimates leads to $\L^\varrho$-boundedness, see Lemma~\ref{lem: OD interpolation} and Lemma~\ref{lem: OD implies boundedness}. Hence, we have $(p_k \vee 2^*, p_{k-1}) \subseteq \cJ(L_{0})$ and since the latter is an interval, we also have $[2_*, p_{k-1}) \subseteq \cJ(L_{0})$. Now, \eqref{eq2: endpoints for Hodge projector} follows by backward induction.

So far, we know that $2<p \leq p_+(L_{0})$ but as $\cP(L_{0})$ is open (see Lemma~\ref{lem: Hodge interval via isomorphism}) we have in fact $2<p<p_+(L_{0})$. By duality we get $(p_-(L_{0}^*) \vee 1) < p' <~2$, so that Theorem~\ref{thm: Riesz characterization HpL} applied to $L_{0}^*$ yields the two-sided estimate
\begin{align*}
\|R_{L_{0}^*} g\|_{p'} \simeq \|g\|_{p'} \quad (g \in \L^2),
\end{align*}
where one (and hence both) sides can be infinite. On the other hand, $p \in \cP(L_{0})$ implies $p' \in \cP(L_{0}^*)$ by duality, that is to say,
\begin{align*}
R_{L_{0}^*}  (R_{L_{0}})^*=- \nabla_x (\eo^*)^{-1} \div_x 
\end{align*}
is $\L^{p'}$-bounded. Here, we used Lemma~\ref{lem: factorization of aL via Riesz} with the roles of $L_{0}$ and $L_{0}^*$ reversed. Altogether, we find for all $f \in \L^{p'} \cap \L^2$ that
\begin{align*}
\|(R_{L_{0}})^*f\|_{p'} 
 \simeq \|R_{L_{0}^*} (R_{L_{0}})^*f\|_{p'} 
\lesssim \|f\|_{p'}.
\end{align*}
This means that $(R_{L_{0}})^*$ is $\L^{p'}$-bounded. By duality, $R_{L_{0}}$ is $\L^p$-bounded and according to Theorem~\ref{thm: Riesz} this can only happen if $p \leq q_+(L_{0})$.
\medskip

\noindent \emph{Step~3: Necessary condition for $\cP(L_{0}) \cap (1,2)$}. Let $p \in \cP(L_{0}) \cap (1,2)$. By duality we get $p' \in \cP(L_{0}^*)$ and Step~2 applied to $L_{0}^*$ gives $p' \leq q_+(L_{0}^*)$. Hence, we have $(q_+(L_{0}^*))' \leq p$.

\medskip

\noindent \emph{Step~4: Conclusion}. Steps~1-3 show that $(q_+(L_{0}^*))'$ and $q_+(L_{0})$ are the endpoints of $\cP(L_{0})$. The latter being an open set by Lemma~\ref{lem: Hodge interval via isomorphism}, we can conclude.
\end{proof}
\section{Critical numbers and kernel bounds}
\label{sec: Critical numbers and kernels}

\noindent In this section, we work out a precise relation between kernel bounds and critical numbers $p_{-}(L)$ \emph{strictly} below $1$. Except for Section~\ref{subsec: kernel estimates for a-1Delta} this is an intermezzo not needed for the application to boundary value problems. However, it nicely illustrates the usefulness of our choice for the interval $\cJ(L)$ compared to \cite{A}  and connects  with the theory of Gaussian estimates in the first chapter of \cite{AT-Asterisque}. In particular, we obtain resolvent kernels from those of high powers of the resolvent without using heat semigroups (which exist only if $\omega_L < \nicefrac \pi 2$). 

It will be convenient to introduce the following notation.\index{$\eta_p$, $p_\eta$ (conversion of $p_-(L)$ and kernel estimates)}

\begin{defn}
\label{def: kernel exponent}
Given $1_{*}<p<1$ and $0<\eta<1$, write $\eta_{p} \coloneqq \nicefrac n p - n$ and conversely $p_{\eta} \coloneqq \frac n {n+\eta}$.
\end{defn}

\subsection{Consequences of \texorpdfstring{$\boldsymbol{p_{-}(L)<1}$}{Consequences of p-(L) < 1}} 
\label{subsec: Consequences of p-(L)<1}

The following result is the core of this section. 
\begin{thm}
 \label{thm: p-<1 equiv high power bounds} 
 The following assertions are equivalent:
 \begin{enumerate}
\item  There exists $p \in (1_*,1)$ such that  $(a(1+t^2L)^{-1}a^{-1})_{t>0}$ is $\H^p$-bounded.
\item There exist $\eta \in (0,1)$ and $\beta(n,\eta) \geq 1$ such that for all integers $\beta\ge \beta(n,\eta)$ the family $((1+t^2L^\sharp)^{-\beta})_{t>0}$ satisfies $\L^2-\L^\infty$ off-diagonal estimates of exponential order and is $\L^2-\Lamdot^\eta$-bounded. 
\end{enumerate}
Moreover, $p_{-}(L)=p_{\eta(L^\sharp)}$, where $\eta(L^\sharp)$ is the supremum of those $\eta$ for which the second property holds.
\end{thm}

For the proof we need an auxiliary result.

\begin{lem}
 \label{lem:L2Linfty OD} 
 Let $(T(t))_{t>0}$ be a family of operators that satisfies $\L^2$ off-diagonal estimates of arbitrarily large (resp.\ exponential) order and that is $\L^2-\Lamdot^\eta$-bounded for some $\eta \in (0,1)$. Then $(T(t))_{t>0}$  satisfies $\L^2-\L^\infty$ off-diagonal estimates of arbitrarily large (resp.\ exponential) order. 
\end{lem}

\begin{proof} 
Lemma~\ref{lem: Loo interpolation inequality} yields $\L^2-\L^\infty$-boundedness. Hence, it suffices to check the off-diagonal estimates when $E,F \subseteq \R^n$ are measurable sets with $\d \coloneqq \d(E,F) \geq t$. 

We let $f \in \L^2$ with support in $E$ and $\|f\|_2 = 1$, set $G \coloneqq \{x \in \R^n: \dist(x,F) \leq \nicefrac{d}{2}\}$, and pick a Lipschitz function $\varphi$ with $\ind_F \leq \varphi \leq \ind_G$ and $\|\nabla \varphi\|_\infty \leq \nicefrac{4}{\dist}$. Lemma~\ref{lem: Loo interpolation inequality}  yields
\begin{align*}
	\|\ind_F T(t) f \|_\infty \lesssim \|\varphi T(t)f \|_2^{\theta} \|\varphi T(t) f \|_{\Lamdot^\eta}^{1-\theta}, 
\end{align*}
where $\theta \in (0,1)$ is such that $(1-\theta)(\nicefrac{n}{2}+\eta) = \nicefrac{n}{2}$. On the right, the first term is bounded by $\|\ind_G T(t) f \|_2^{\theta}$ and as $\dist(G, E) \geq \nicefrac{\d}{2}$, this gives the required off-diagonal decay. The second term is controlled by
\begin{align*}
	 \big(\|\varphi\|_\infty \|T(t) f \|_{\Lamdot^\eta} + \|\varphi\|_{\Lamdot^\eta} \|T(t) f \|_\infty\big)^{1-\theta}
	 &\lesssim \big( t^{-\eta-\frac{n}{2}} + \dist^{-\eta} t^{-\frac{n}{2}}\big)^{1-\theta} \\
	 &\leq t^{-n/2},
\end{align*}
using the $\L^2-\Lamdot^\eta$-bound, the $\L^2-\L^\infty$-bound and $\d \geq t$.
%
\end{proof}

\begin{proof} [Proof of Thm \ref{thm: p-<1 equiv high power bounds}] Let us recall from Corollary~\ref{cor: OD for second order} that the resolvents of $L$ (and hence of $L^\sharp$) satisfy $\L^2$ off-diagonal estimates of exponential order.

\medskip

\noindent $(i) \Longrightarrow (ii)$.
The family $(a(1+t^2L)^{-1}a^{-1})_{t>0}$ is also $\H^\varrho - \L^2$-bounded for some $\varrho < 2$ depending on the dimension $n$, see Lemmata~\ref{lem: Lp L2 between Sobolev conjugates n geq 2} and \ref{lem: Lp L2 between Sobolev conjugates n = 1}. Thus, Lemma \ref{lem: extra} implies that for any $p < q \leq 1$ there exists $\beta(n,q)$ such that for all $\beta \geq \beta(n,q)$ the family $(a(1+t^2L)^{-\beta}a^{-1})_{t>0}$ is $\H^q-\L^2$-bounded.  By duality (Lemma \ref{lem: p-q boundedness duality}) and the fact that $L^\sharp= (a^*)^{-1} L^*a^*$, it follows that $((1+t^2L^\sharp)^{-\beta})_{t>0}$  is $\L^2-\Lamdot^\eta$-bounded with $\eta \coloneqq \eta_{q}$. It remains to apply Lemma \ref{lem:L2Linfty OD}.

\medskip

\noindent $(ii) \Longrightarrow (i)$.
Let $p \coloneqq p_{\eta}$. By duality, $(a(1+t^2L)^{-\beta}a^{-1})_{t>0}$ is $\H^p-\L^2$-bounded. This family satisfies $\L^2$ off-diagonal estimates of exponential order and $\int_{\R^n} a(1+t^2L)^{-\beta}(a^{-1}f)\, \d x=0$  holds for all $f \in \L^2$ with compact support and mean value $0$, see Corollary~\ref{cor: conservation}. Lemma~\ref{lem: OD Hardy implies boundedness} yields $\H^q$-boundedness for any $p<q\le1$, hence for any $p<q\leq2$ by interpolation with the $\L^2$-boundedness. If $\beta = 1$, then we can stop here. From now on, we assume $\beta \geq 2$.

Let $p<q \leq 1$. By interpolation with the original $\H^p-\L^2$-boundedness we obtain $\H^q-\L^r$-boundedness for some $r>1$ with $0\leq \nicefrac{n}{q} - \nicefrac{n}{r} < 1$. Now we apply the formula
\begin{align}
\label{eq:highpower to resolvent}
(1+t^2 L)^{-1}(a^{-1}f) = (\beta-1) \int_0^\infty (1+u+t^2u L)^{-\beta}(a^{-1}f)\, \d u,
\end{align}
which follows for $f\in \H^q \cap \L^2$ by applying $(1+t^2L)^{-1}$ to both sides of \eqref{eq: resolution of resolvents through high powers} and conclude
\begin{align*}
\|(1+t^2 L)^{-1}(a^{-1}f)\|_{r} 
&\lesssim  \int_0^\infty (1+u)^{-\beta}\left(\frac{1+u}{t^2u}\right)^{\frac{n}{2q} - \frac{n}{2r}} \|f\|_{\H^q} \, \d u \\
& \lesssim t^{\frac n r - \frac n q} \|f\|_{\H^q},
\end{align*}
where $0\leq \nicefrac{n}{q} - \nicefrac{n}{r} < 1 $ has guaranteed that the integral in $u$ is finite. This proves that $(a(1+t^2L)^{-1}a^{-1})_{t>0}$ is $\H^q-\L^r$-bounded. 

In the same manner we can start with $\L^r$-boundedness for the higher-order resolvents when $1<r\leq 2$ and obtain first $\L^r$-boundedness of $(a(1+t^2L)^{-1}a^{-1})_{t>0}$ and then $\L^r$ off-diagonal estimates of exponential order by interpolation with the $\L^2$-result.

Now, we apply again Lemma ~\ref{lem: OD Hardy implies boundedness} to conclude $\H^q$-boundedness of $(a(1+t^2L)^{-1}a^{-1})_{t>0}$ whenever $p<q \leq 1$.

Re-examination of the proof shows that the stated relation for $p_{-}(L)$ holds. 
\end{proof}

The following corollary is interesting because $\L^1$ and $\L^\infty$ are not part of our $\cJ(L)$-theory.

\begin{cor}
 \label{cor:Linfty and L1 bounded}
If $p_{-}(L)<1$, then $((1+t^2L)^{-1})_{t>0}$ satisfies $\L^1$ off-diagonal estimates of exponential order and  $((1+t^2L^\sharp)^{-1})_{t>0}$ satisfies $\L^\infty$ off-diagonal estimates of exponential order. In particular, these families are $\L^1$-bounded and $\L^\infty$-bounded, respectively.
\end{cor}

\begin{proof} 
It directly follows from Theorem~\ref{thm: p-<1 equiv high power bounds} and Remark~\ref{rem:OD} that for $\beta \geq 2$ large enough $((1+t^2L^\sharp)^{-\beta})_{t>0}$ satisfies $\L^\infty$ off-diagonal estimates of exponential order. Hence, $((1+t^2L)^{-\beta})_{t>0}$ satisfies $\L^1$ off-diagonal estimates of exponential order.

By the formula \eqref{eq:highpower to resolvent} applied to $af$, where $f \in \L^1\cap \L^2$, we see that $(1+t^2L)^{-1}$ has the desired property. Indeed, if $f$ has support in $E$ and $F$ is another measurable set, then with the change of variable  $v=(\frac {1+u}u)^{1/2}$ in the integral we obtain
\begin{align*}
\|(1+&t^2 L)^{-1}f\|_{\L^1(F)} \\
&\le (\beta-1) \int_0^\infty (1+u)^{-\beta} \Big\|\Big(1+t^2\frac u {1+u} L\Big)^{-\beta}f\Big\|_{\L^1(F)} \, \d u
\\ 
& =  2(\beta-1) \int_{1}^\infty  (v^2-1)^{\beta-2}v^{1-2\beta} \Big\|\Big(1+\frac{t^2}{v^2} L\Big)^{-\beta}f\Big\|_{\L^1(F)}\, \d v
\\
& \lesssim \int_{1}^\infty v^{-3} \e^{-c \frac{\dist (E,F)v}{t}} \|f\|_{1}\, \d v
\\
& \lesssim  \e^{-c \frac{\dist (E,F)}{t}}  \|f\|_{1},
\end{align*}
where we used $\beta \ge 2$ in the third step.

Finally, the claim for $((1+t^2L^\sharp)^{-1})_{t>0}$ follows by duality and similarity.
\end{proof}
\subsection{Equivalence with kernel estimates}
\label{subsec: Equivalence with kernel estimates}

Going one step further, we shall now incorporate pointwise kernel estimates into the machinery.  The convention on the variables for integral kernels is that we always look for representations in the form
\begin{align*}
	(Tf)(x)=\int_{\R^n}K(x,y)f(y)\, \d y.
\end{align*}
We rely on the following lemma.\index{off-diagonal estimates!relation with kernel bounds}

\begin{lem}
 \label{lem: half-equivalence}
Let $(T(t))_{t>0}$ be a family of bounded operators on $\L^2$ and denote by $K_{t}(x,y)$ their distribution kernels.  For every $\eta \in (0,1)$ the following assertions are equivalent:
 \begin{enumerate}
\item   $(T(t))_{t>0}$ satisfies $\L^2-\L^\infty$ off-diagonal estimates of exponential order and is $\L^2-\Lamdot^\eta$-bounded. 
\item For each $t>0$, $K_{t}(x,y)$ agrees with a measurable function and there are constants $C, c > 0$ that do not depend on $t$ such that for all $x, h\in \R^n$ and all measurable sets $E$,
\begin{align}
\label{eq:ktL2Linfty}
\int_{E}|K_{t}(x,y)|^2\, \d y \le C t^{-n} \e^{-c \frac{\d(x,E) }{t}},
\end{align}
\begin{align}
\label{eq:ktL2holder}
\int_{\R^n}|K_{t}(x+h,y)-K_{t}(x,y)|^2\, \d y \le C |h|^{2\eta}t^{-n-2\eta}.
\end{align}
\end{enumerate}

\end{lem}

\begin{proof}
The implication (ii) $\Longrightarrow$ (i) is a direct consequence of the Cauchy--Schwarz inequality.

Next, assume that (i) holds. Fix $t>0$. As pointed out in \cite[Thm.~1.3]{Arendt-Bukhvalov}, any linear operator $T(t)$ that is bounded from $\L^2$ to $\L^\infty$ has an integral representation 
\begin{align*}
	T(t)f(x)=\int_{\R^n}K_{t}(x,y)f(y)\, \d y \quad (f \in \L^2, \text{ a.e. } x \in \R^n)
\end{align*}
with a measurable kernel that belongs to $\L_x^\infty( \L_y^2)$ with norm equal to the operator norm. Hence, $K_{t}(x,y)$ can indeed be identified to a measurable function that satisfies \eqref{eq:ktL2Linfty}. For $h \in \R^n$ let $\tau_h$ be the translation operator $f \mapsto f(\, \cdot \, + h)$. Since $(T(t))_{t>0}$ is $\L^2 - \Lamdot^\eta$-bounded, the family $((\nicefrac{t}{|h|})^{\eta}(1-\tau_{h})T(t))_{t>0}$ is  $\L^2-\L^\infty$-bounded, uniformly in $h$, and we may apply the above result again to obtain \eqref{eq:ktL2holder}.  
\end{proof}
 
We introduce two auxiliary functions that naturally appear in kernel estimates for the resolvents.
 
\begin{defn}
\label{def: kernel auxiliary functions}
Define functions $\omega_n, \wt{\omega}_n: (0,\infty) \to (0,\infty)$ by\index{$\omega_n, \wt{\omega}_n$ (kernel estimates)}
\begin{align*}
	\omega_n(s) \coloneqq \begin{cases} 1 & \text{if $n=1$} \\ |\ln s\,|+1 & \text{if $n=2$} \\ s^{2-n} & \text{if $n\geq3$} \end{cases} \quad \text{and} \quad 
	\wt{\omega}_n(s) \coloneqq \begin{cases} 1 & \text{if $n=1,2$}  \\ s^{2-n} & \text{if $n\geq3$} \end{cases}.
\end{align*}
\end{defn}

Combining Theorem~\ref{thm: p-<1 equiv high power bounds} with Lemma~\ref{lem: half-equivalence} allows us to characterize the property $p_-(L) < 1$ through $\L^2$ kernel bounds of a large power of the resolvent. What is missing to get to pointwise kernel bounds is dual information on $L^\sharp$. 

 \begin{thm}
\label{thm:pointwiseboundsequivalence} 
The following assertions are equivalent:

\begin{enumerate}
 \item There exists $p\in (1_{*},1)$ such that $(a(1+t^2L)^{-1}a^{-1})_{t>0}$ and $(a^*(1+t^2L^\sharp)^{-1}(a^*)^{-1})_{t>0}$ are $\H^p$-bounded.
 
  \item There exists $\eta\in (0,1)$ such that for all $t>0$ the operator $(1+t^2 L)^{-1}a^{-1}$ is given by a measurable kernel $G_t(x,y)$ that satisfies, for some constants $C, c > 0$,  the following bounds:
	\begin{align}
	\label{eq:pointwiseGt}
	|G_t(x,y)| &\leq C t^{-n} \omega_{n}\bigg(\frac{|x-y|}{t}\bigg) \e^{-c \frac{|x-y|}{t}},
	\\
	\label{eq:holderGty}
	\begin{split}
	  |G_t(x,y+h)-&G_t(x,y)|  \\
	&\leq C t^{-n} \bigg(\frac{|h|}{|x-y|}\bigg)^{\eta}  \wt{\omega}_{n}\bigg(\frac{|x-y|}{t}\bigg) \e^{-c \frac{|x-y|}{t}},
	\end{split}
	\\
	\label{eq:holderGtx}
	\begin{split}
	  |G_t(x+h,y)-&G_t(x,y)| \\
	&\leq C t^{-n} \bigg(\frac{|h|}{|x-y|}\bigg)^{\eta}  \wt{\omega}_{n}\bigg(\frac{|x-y|}{t}\bigg) \e^{-c \frac{|x-y|}{t}},
	\end{split}
	\end{align}
	provided that $2|h| \leq |x-y|$.
\end{enumerate}
Moreover, if either condition holds, then
\begin{align*}
	p_{-}(L)=p_{\eta(L^\sharp)} \quad \& \quad p_{-}(L^\sharp)=p_{\eta(L)},
\end{align*}
where $\eta(L^\sharp)$ and $\eta(L)$ are the suprema of those $\eta$ for which  \eqref{eq:holderGty} and \eqref{eq:holderGtx} hold, respectively. \index{critical numbers!relation to kernel bounds}
\end{thm}

\begin{proof}
We argue in three steps.

\medskip

\noindent \emph{Step 1: $(i) \Longrightarrow (ii)$.}
We apply Theorem~\ref{thm: p-<1 equiv high power bounds} and Lemma~\ref{lem: half-equivalence} to both $L$ and $L^\sharp$. Hence, there is an even integer $\beta$ such that $(1+t^2L)^{-\beta/2}$ and $(1+t^2 L^\sharp)^{-\beta/2} (a^*)^{-1}$ are given by measurable kernels $K_t^{(L)}(x,y)$ and $K_t^{(L^\sharp)}(x,y)$, respectively, and both kernels satisfy \eqref{eq:ktL2Linfty} and \eqref{eq:ktL2holder}. By duality and composition, we see that $(1+t^2 L)^{-\beta} a^{-1}$ is an integral operator given by the kernel
\begin{align*}
	G_t^\beta(x,y) \coloneqq \int_{\R^n} K_t^{(L)}(x,z) K_t^{(L^\sharp)}(y,z) \, \d z.
\end{align*}
We claim that there are constants $C,c \in (0,\infty)$, $\eta \in (0,1)$ such that for all $x,y,h$,   
\begin{align}
\label{eq:pointwiseGt beta}
|G_t^\beta(x,y)| &\leq C t^{-n}  \e^{-c \frac{|x-y|}{t}},
\\
\label{eq:holderGty beta no decay}
|G_t^\beta(x,y+h)-G_t^\beta(x,y)|  
&\leq C t^{-n} \bigg(\frac{|h|}{t}\bigg)^{\eta}   ,
\\
\label{eq:holderGtx beta no decay}
|G_t^\beta(x+h,y)-G_t^\beta(x,y)| 
&\leq C t^{-n} \bigg(\frac{|h|}{t}\bigg)^{\eta}.
\end{align}
Indeed, \eqref{eq:holderGty beta no decay} and \eqref{eq:holderGtx beta no decay} follow directly from \eqref{eq:ktL2Linfty}, \eqref{eq:ktL2holder} and the Cauchy--Schwarz inequality. The same argument yields the first estimate if we split integration in $z$ into the parts where $|x-z| \geq \nicefrac{|x-y|}{2}$ and $|y-z| \geq \nicefrac{|x-y|}{2}$ beforehand. Note that $\eta$ in \eqref{eq:holderGty beta no decay} and \eqref{eq:holderGtx beta no decay} can be any exponent such that $p_\eta > p_-(L)$ and $p_\eta > p_-(L^\sharp)$, respectively. 

Taking logarithmic convex combinations of \eqref{eq:pointwiseGt beta} with \eqref{eq:holderGty beta no decay} and \eqref{eq:holderGtx beta no decay}, we obtain in the same ranges of $\eta$ but with different constants $C,c>0$ the following H\"older estimates with exponential decay when $2|h| \leq |x-y|$:
\begin{align}
\label{eq:holderGty beta}
|G_t^\beta(x,y+h)-G_t^\beta(x,y)|  
&\leq C t^{-n} \bigg(\frac{|h|}{t}\bigg)^{\eta}   \e^{-c \frac{|x-y|}{t}},
\\
\label{eq:holderGtx beta}
|G_t^\beta(x+h,y)-G_t^\beta(x,y)| 
&\leq C t^{-n} \bigg(\frac{|h|}{t}\bigg)^{\eta}   \e^{-c \frac{|x-y|}{t}}.
\end{align}

From there, it suffices to use again the formula \eqref{eq:highpower to resolvent} for all $f\in \L^1 \cap \L^2$ and \eqref{eq:pointwiseGt beta}
to see that $(1+t^2L)^{-1} a^{-1}$ is given by a kernel with bound
\begin{align*}
|G_t(x,y)| 
&\leq C t^{-n} \int_0^\infty (1+u)^{-\beta} (1+\tfrac{1}{u})^{\frac{n}{2}} \e^{-c\frac{|x-y|}{t}(1+\frac{1}{u})^{\nicefrac{1}{2}}} \, \d u \\
&=2 C t^{-n} \int_1^\infty v^{n-2\beta+1} (v^2-1)^{\beta-2} \e^{-c \frac{|x-y|}{t}v} \, {\d v} \\
&\le 2C t^{-n}  \int_{1}^\infty v^{n-3} \e^{-c \frac{|x-y|}{t}v} \, \d v,
\end{align*}
where we used the change of variable $v=(1+\frac{1}{u})^{\nicefrac{1}{2}}$ and $\beta\ge 2$. The latter integral is controlled by $\omega_n(\nicefrac{|x-y|}{t}) \e^{-c\nicefrac{|x-y|}{2t}}$ and  \eqref{eq:pointwiseGt} follows.

Next, we use the same strategy starting from \eqref{eq:holderGtx beta}. In that case, we assume $2|h|\le |x-y|$ and we obtain  
\begin{align*}
|G_t(x+h ,y)-G_{t}(x,y)|  \leq 2C |h|^\eta t^{-n-\eta}
  \int_{1}^\infty v^{n+\eta -3} \e^{-c \frac{|x-y|}{t}v} \, \d v,
\end{align*} 
from which we conclude readily for \eqref{eq:holderGtx}. The argument to obtain  \eqref{eq:holderGty} from \eqref{eq:holderGty beta} is the same. 

\medskip

\noindent \emph{Step~2: $(ii) \Longrightarrow (i)$.}
For the converse, let $\eta$ be given as in the estimates and let $1>p>p_{\eta}$. It is enough to  show that $(a(1+t^2L)^{-1}a^{-1})_{t>0}$ is $\H^p$-bounded. The argument for the adjoint is the same.

To this end, it suffices to establish for some $C, \eps > 0$ the molecular bounds 
\begin{align}
\label{eq1: molecular}
\|a(1+t^2L)^{-1} a^{-1}m\|_{\L^2(C_j(B))} \leq C (2^{j} r(B))^{\frac{n}{2}-\frac{n}{p}} 2^{-\eps j} \quad (j \geq 1),
\end{align}
whenever $t>0$ and $m$ is an $\L^2$-atom for $\H^p$ associated with a ball $B$. Indeed, since $a(1+t^2L)^{-1}a^{-1}m$ has integral zero by Corollary~\ref{cor: conservation}, we can first use Lemma~\ref{lem: Taibleson-Weiss} to get a uniform $\H^p$-bound and then conclude by the ($\L^2$-convergent) atomic decomposition for functions in $\H^p \cap \L^2$. 
 
For $j = 1$ we have as required
\begin{align*}
\|a(1+t^2L)^{-1} a^{-1}m\|_{\L^2(C_1(B))}
\leq  C\|m\|_2 \leq C r(B)^{\frac{n}{2}-\frac{n}{p}}.
\end{align*}
For $j \geq 2$ we use the mean value property of $m$ to write
\begin{align*}
a (1+t^2 L)^{-1} a^{-1} m(x) = \int_{B} a(x)(G_t(x,y)-G_t(x,y_B))m(y) \, \d y
\end{align*}
with $y_B$ the center of $B$ and obtain for $x \in C_j(B)$ that
\begin{align*}
|a(1+t^2 L_0)^{-1}a^{-1} m(x)| 
&\leq C t^{-n} 2^{-j \eta} \wt{\omega}_{n}\bigg(\frac{2^{j-1} r(B)}{t}\bigg) \e^{-c \frac{2^{j-1} r(B)}{t}} \|m\|_1 \\
&\leq C C' (2^{j-1}r(B))^{-n} 2^{-j \eta} r(B)^{n-\frac{n}{p}},
\end{align*}
where $C' \coloneqq \sup_{s>0} s^n\wt{\omega}_{n}(s) \e^{-cs}$. Integrating  the square of this inequality on $C_j(B)$ 
and sorting powers of $2^j$ and $r(B)$  gives us \eqref{eq1: molecular} with $\eps = -\nicefrac{n}{p}+n+\eta$. Now $\eps > 0$ is equivalent to $p> p_{\eta}$, which we have assumed. 

\medskip

\noindent \emph{Step~3: The formul\ae \, for the critical numbers.} In Step~1 we have obtained \eqref{eq:holderGty} if $p_\eta > p_-(L)$, whereas in Step~2 we have obtained $\H^p$-boundedness if $p > p_{\eta(L^\sharp)} $. Thus, we have $p_-(L) = p_{\eta(L^\sharp)}$. We have also seen the same conclusions with the roles of $L$ and $L^\sharp$ interchanged.
\end{proof}

\begin{rem}
In dimension $n=1$ it is shown in \cite{AT} that the first-order derivatives of $G_{t}(x,y)$ in $x$ and $y$ exist and have an exponentially decaying pointwise bound in $|x-y|$. In particular,  $\eta(L)=1$ is attained.  
 \end{rem}
 
\begin{rem} 
Under one of the conditions of Theorem~\ref{thm:pointwiseboundsequivalence} one can also obtain pointwise and H\"older bounds for the kernel $G(x,y)$ of $L^{-1}a^{-1}$ when $n\ge 3$. Since $L^{-1}a^{-1}= L_{0}^{-1}$, this kernel $G$ is just the Green kernel of $L_{0}$ and does not depend on $a$. To see the estimates,  it suffices to replace  the formula \eqref{eq:highpower to resolvent} by the Calderón reproducing formula $L_{0}^{-1}f=(\beta -1) \int_0^\infty (1+u L_{0})^{-\beta}f\, \d u$ that is valid for $f \in \ran(L_0)$ and to plug in the estimates \eqref{eq:pointwiseGt beta}, \eqref{eq:holderGty beta} and \eqref{eq:holderGtx beta}. This does not work for $n=1,2$.
\end{rem}
\subsection{Dirichlet property, stability and examples}
\label{subsec: Stability and examples}

Having made the link between critical numbers strictly below one and kernel estimates for the resolvent, opens the door to further characterizations of either property in terms of regularity theory for the corresponding elliptic system in $\R^n$. 

We shall use the notion of weak solutions and Caccioppoli's inequality. A reader who is not familiar with these tools will find all necessary background material (written for systems in $\ree$) in Section~\ref{sec: Basic properties of weak solutions} below. The Dirichlet property for $L_{0} = -\div_x d \nabla_x$ is the following quantitative regularity property.

\begin{defn}
\label{def: Dirichlet property}
The operator $L_0$ satisfies the \emph{Dirichlet property}\index{Dirichlet property} if there are  $\mu \in (0,1)$ and  $C_{0} \in (0,\infty)$ such that for all $R >0$ and all $x_0\in \R^n$ it follows that any weak solution $v \in \W^{1,2}(B(x_0,R))$ to $\div_x d \nabla_x v = 0$ in $B(x_0,R)$ satisfies
\begin{align}
\label{eq:Dirichet property}
\int_{B(x_0,\rho)} |\nabla v|^2 \, \d x \le C_0
\left( \frac{\rho}
{R}\right)^{n-2+2\mu}\int_{B(x_0,R)} |\nabla v|^2 \, \d x,
\end{align}
when $0<\rho\le R$. The supremum of those $\mu$ for which this property holds is denoted by $\mu(L_{0})$.\index{$\mu(L_{0})$ (Dirichlet property)}
\end{defn}

\begin{rem}
 \label{rem:Dirichlet  property}
 The Dirichlet property has been discussed in detail in \cite[Sec.~1]{AT-Asterisque} for elliptic equations ($m=1$). Among others, it was shown that it is stable under small $\L^\infty$-perturbations of the coefficients $d$ and   that it holds  when $n=1$ with $\mu(L_{0})=1$, when $n=2$ with $\mu(L_{0})>0$ and when $n\ge 3$ for real-valued $d$ with $\mu(L_{0})>0$ or with $\mu(L_{0})=1$ when $d$ has small enough $\BMO$-norm. The latter example includes in particular the case of constant coefficients. More exotic examples are given by coefficients $d$ that depend only on one coordinate. In this case $\mu(L_{0})=1$, see \cite[App.~B]{AT-Asterisque}. For systems ($m\ge 2$) all examples but the case of real-valued coefficients can be adapted.
\end{rem}

Let us prove that that critical numbers below $1$ are also characterized through the Dirichlet property for the adjoint.

\begin{thm}
\label{thm:Dirichlet property and p-} 
$p_{-}(L_{0})<1$ if and only if $L_{0}^*$ satisfies the Dirichlet property. Moreover, $p_{-}(L_{0})=  p_{\mu(L_{0}^*)}$.
\end{thm}

By Theorem~\ref{thm: critical numbers a-independent} the critical numbers for $L$ and $L_0$ are the same. Hence, we immediately obtain

\begin{cor}
 \label{cor:stability} The condition $p_{-}(L)<1$ is satisfied exactly when $L_{0}^*$ has the Dirichlet property. \index{critical numbers!relation to Dirichlet property}
\end{cor}

\begin{proof}[Proof of Theorem~\ref{thm:Dirichlet property and p-}]
By Theorem \ref{thm: p-<1 equiv high power bounds} we can replace the assertion $p_-(L_0) < 1$ by the existence of $0<\eta<1$ and $\beta(n,\eta) \geq 1$ such that for all integers $\beta\ge \beta(n,\eta)$ the family $((1+t^2L_{0}^*)^{-\beta})_{t>0}$ is  $\L^2-\Lamdot^\eta$-bounded and satisfies $\L^2-\L^\infty$ off-diagonal estimates of exponential order. 

We shall prove that under this assumption the Dirichlet property for $L_0^*$ holds for any $\mu \in (0,\eta)$ and that conversely the Dirichlet property for $L_0^*$ with exponent $\mu$ implies the above for any $\eta \in (0,\mu)$. Once this is done, also $p_{-}(L_{0})=  p_{\mu(L_{0}^*)}$ follows from Theorem~\ref{thm: p-<1 equiv high power bounds}. 

\medskip

\noindent \emph{Step~1: From $p_{-}(L)<1$ to property $(H)$.} Let $0<\mu<\eta$. We prove that $L_{0}^*$ has the so-called property $(H)$ with exponent $\mu$: There is a constant
$C$ depending on $L_{0}^*$  such that for any ball  $B$ of radius $R>0$ and any $u\in \W^{1,2}(B)$ with $\div_x d^* \nabla_x u=0$ on $B$ in the weak sense it follows that
\begin{align*}
\sup_{\frac{1}{4} B}|u| +  R^\mu 
\sup_{  
{(x,y) \in \frac{1}{4} B, \, x \neq y}}
\frac{
|u(x) - u(y)|}{|x-y|^\mu} \le C \Big( \barint_{B}  |u|^2 \, \d x \Big)^{1/2}. 
\end{align*} 
The proof is a modification of an argument in \cite[Sec.~1.4.2]{AT-Asterisque}.

Let $u \in \W^{1,2}(B)$ be a weak solution to $\div_x d^* \nabla_x u=0$ in $B$. Let  $\chi \in \C^\infty_0(\R^n)$ be supported in $\frac{8}{9}B$ with $\chi= 1$ on
$\frac{7}{8} B$ and $\|\nabla \chi\|_{\infty}\leq c R^{-1}$ for a dimensional constant $c$. Let $v \coloneqq u\chi$. Since $v=u$ on $\frac{7}{8}B$, it
suffices to show that for any $\varphi \in \C^\infty_0(\frac{1}{4}B)$ and any $h \in \frac{1}{2} B$ we have
\begin{align}
\label{eq: H1}
\bigg|\int_{\R^n} v(x) \cl{\varphi(x)} \, \d x\bigg| \le CR^{-\frac{n}{2}} \|\varphi\|_1 \|u\|_{\L^2(B)} 
\end{align}
and 
\begin{align}
\label{eq: H2}
\bigg|\int_{\R^n} (v(x+h) - v(x)) \cl{\varphi(x)} \, \d x\bigg| \le C |h|^\mu
R^{-\mu-\frac{n}{2}} \|\varphi\|_1 \|u\|_{\L^2(B)}. 
\end{align}

We abbreviate inner products in $\L^2(\R^n)$ by $\angle v \varphi \coloneqq \int_{\R^n} v \cl{\varphi} \, \d x$ and set $T(t) \coloneqq
(1+t^2L_{0})^{-1}$. Since $T(t)^\beta\varphi \in  \W^{1,2}(\R^n)$ and $v \in \W^{1,2}(\R^n)$, we can write
\begin{align}
\label{eq: p- to H goal}
\begin{split}
\angle v \varphi &=
  \angle {u\chi} {T(R)^\beta \varphi} - \int_0^{R} \Big \langle u \chi, \frac{\d}{\d t} T(t)^\beta \varphi \Big \rangle  \, \d t
 \\
&= \angle {u\chi} {T(R)^\beta \varphi} + 2\beta \int_0^{R}
\angle {\nabla_{x}(u\chi)} {d\nabla_{x} T(t)^{\beta+1}\varphi} \,t  \d t.
\end{split}
\end{align}
By duality the family $(T(t)^\beta)_{t>0}$ satisfies $\L^1-\L^2$ off-diagonal estimates of exponential order. In particular, it is $\L^1-\L^2$-bounded and we obtain
\begin{align}
\label{eq1: p- to H goal}
|\angle {u\chi} {T(R)^\beta \varphi}|  \le \|{u\chi}\|_2\|T(R)^\beta \varphi\|_2
\lesssim \|u\|_{\L^2(B)} R^{-\frac{n}{2}}\|\varphi\|_1.
\end{align}
Next, we rewrite the inner product inside the integral in \eqref{eq: p- to H goal} as
\begin{align*}
&\angle {d^* \nabla_{x} u  } {\nabla_{x} (\chi T(t)^{\beta+1}\varphi)} \\
&\quad +\angle{  d^*(\nabla_{x} \chi \otimes u)} {\nabla_{x} T(t)^{\beta+1}\varphi} \\
&\quad - \angle{  d^*\nabla_{x} u} {\nabla_{x} \chi \otimes T(t)^{\beta+1}\varphi} \\
&\quad \eqqcolon \I +\II - \III,  
\end{align*}
where $\nabla_{x} \chi \otimes u$ is short for the vector in $(\IC^m)^n$ that comes from the product rule when calculating $\nabla_x(\chi u)$. 

The term $\I$ vanishes thanks to the equation for $u$.

For the term $\II$ we note that  $(t\nabla_{x} T(t)^{\beta+1})_{t>0}$ satisfies $\L^1-\L^2$ off-diagonal estimates of exponential order by composing the $\L^1-\L^2$-estimates for $(T(t)^\beta)_{t>0}$ and the $\L^2$-estimates for $(t\nabla_x T(t))_{t>0}$ from Corollary~\ref{cor: OD for second order}. As the supports of $\varphi$ and $\nabla_x \chi$ have distance at least $\frac{5}{8} R$, we obtain for  some $\alpha>0$ that
\begin{align*}
| t \II| 
\lesssim  R^{-1} t^{- \frac{n}{2}} \e^{- \frac{\alpha R}{t}} \|u\|_{\L^2(B)} \|\varphi\|_1.
\end{align*}
Similarly, we get
\begin{align*}
|t \III|
\lesssim R^{-1} \|\nabla u\|_{\L^2(\frac{8}{9}B)} t^{1-\frac{n}{2}} \e^{- \frac{\alpha R}{t}}
\|\varphi\|_1
\end{align*}
and hence by the Caccioppoli inequality
\begin{align*}
|t\III| \lesssim R^{-2} t^{1-\frac{n}{2}} \e^{- \frac{\alpha R}{t}} \|u\|_{\L^2(B)}\|\varphi\|_1.
\end{align*}
Going back to \eqref{eq: p- to H goal}, we obtain by integration that
\begin{align*}
\bigg|\int_0^{R} \angle {\nabla_{x}(u\chi)} {d\nabla_{x} T(t)^{\beta+1}\varphi}\, t 
\d t\bigg| \lesssim R^{-\frac{n}{2}}  \|u\|_{\L^2(B)}\|\varphi\|_1
\end{align*}
as desired. Together with \eqref{eq1: p- to H goal} this proves \eqref{eq: H1}.

The integral in \eqref{eq: H2} can be interpreted as $\angle v {\varphi_{h}}$, where $\varphi_{h} \coloneqq {(1-\tau_{-h})\varphi}$ and $\tau_h$ is the translation operator $f \mapsto f(\, \cdot \, + h)$ as before. We replace $T(t)^\beta\varphi $ by $T(t)^\beta\varphi_{h} $ and can run the same argument provided we still have the necessary bounds, namely:
\begin{enumerate}
	\item $((\nicefrac{t}{|h|})^{\mu} T(t)^\beta (1-\tau_{-h}))_{t>0}$ is $\L^1-\L^2$-bounded, uniformly in $h$.
	\item When $|h|\le \nicefrac{R}{2}$ and $S(t)$ is one of $T(t)^{\beta+1}$ or $t\nabla_{x}T(t)^{\beta+1}$, then $(\nicefrac{t}{|h}|)^{\mu} S(t)(1-\tau_{-h})$ is bounded from $\L^1(\frac{1}{4}B)$ into $\L^2(\frac{8}{9}B \setminus \frac{7}{8}B)$ with norm controlled by $t^{-n/2}\e^{-cR/t}$.
\end{enumerate}
To obtain (i), we note that by duality $((\nicefrac{t}{|h|})^{\eta} T(t)^\beta (1-\tau_{-h}))_{t>0}$ is $\L^1-\L^2$-bounded, uniformly in $h$, and the same is true for $(T(t)^\beta (1-\tau_{-h}))_{t>0}$, using the uniform bound of the translations on $\L^1$. A logarithmic convex combination of the two estimates gives (i). Then (ii) follows in the same way but we use the $\L^1-\L^2$ off-diagonal decay for $S(t)$ and that $(1-\tau_{-h})\varphi$ is supported in $\frac{3}{4} B$ whenever $|h| \leq \nicefrac{R}{2}$ and $\varphi$ is supported in $\frac{1}{4} B$.

This completes the proof of property $(H)$.

\medskip

\noindent \emph{Step~2: From property $(H)$ to the Dirichlet property.} 
Condition $(H)$ for $L_{0}^*$ implies the Dirichlet property for $L_{0}^*$ with the same $\mu$. This argument is done in \cite[p.45]{AT-Asterisque}.

\medskip

\noindent \emph{Step~3: From the Dirichlet property to resolvent bounds.} Assuming the Dirichlet property for $L_{0}^*$ with exponent $\mu$, it suffices to follow line by line the argument in \cite[Sec.~1.4.3]{AT-Asterisque} up until the intermediate result of formula (38) which, in particular, states that for large enough integer $k_{0}$ the family $((1+t^2L_{0}^*)^{-k_{0}})_{t>0}$ is 
$\L^2-\Lamdot^\eta$-bounded for any $\eta<\mu$. Then we conclude using Lemma~\ref{lem:L2Linfty OD}.
\end{proof}
\subsection{Remarks on  multiplicative perturbations}
\label{subsec: multiplicative perturbations}

It is instructive to put our results in perspective with Theorem~\ref{thm: critical numbers a-independent}, which states that the numbers $p_{-}(L)$ are $a$-independent, that is $p_{-}(L)=p_{-}(L_{0})$ if we write $ L$ as a multiplicative perturbation $L=a^{-1}L_{0}$. There is no other condition on $a$ than the standing ellipticity condition from Section~\ref{subsec: definition elliptic operators}. This implies that the set of estimates on the kernel  for $(1+t^2L_{0})^{-1}$ in Theorem~\ref{thm:pointwiseboundsequivalence} is equivalent to the similar one for the kernel of $(1+t^2L)^{-1} a^{-1}$. 
 
Prior to that there were works on multiplicative perturbations involving semigroups.  Duong and Ouhabaz~\cite{Ouhabaz, Duong-Ouhabaz} proved that semigroup kernel estimates for  $\e^{-tL_{0}}$ imply semigroup kernel estimates for $\e^{-ta^{-1}L_{0}}a^{-1}$ if $d$ is an $n\times n$  matrix with real valued coefficients (so $m=1$) under the additional assumption that $d$ is symmetric or, more generally, that the sectoriality angle of $a^{-1}L_0$ does not exceed $\nicefrac{\pi}{2}$. This condition is of course necessary to define a holomorphic semigroup and allows one to use contour integrals. 

Before that, work of McIntosh--Nahmod dealt with the specific case of $L=-a^{-1}\Delta_x$, see \cite{McN}. It was shown in \cite{AQ} that the only restriction to transfer a set of estimates called condition $(G)$ on the semigroup kernel of  $\e^{-tL_{0}}$ to the corresponding one for $\e^{-tL}a^{-1}$  is the sectoriality of $L$. 

The conclusion is that if estimates on the resolvent kernels or their high powers suffice for an application, then the existence of the semigroup generated by $-L$  can be removed. Besides, the arguments  are somewhat less involved than those passing through semigroups. 
\subsection{Kernel estimates for \texorpdfstring{$\boldsymbol{L=-a^{-1}\Delta_x}$}{a^-1 Delta}}
\label{subsec: kernel estimates for a-1Delta}

We close this section with kernel estimates in the special case of $L=-a^{-1}\Delta_x$ that are used later in this monograph. Some of them are due to \cite{McN}. Interestingly, we use a much simpler method than the original proof and we obtain further estimates, notably those on mixed second-order derivatives. Corollary~\ref{cor: McIntosh-Nahmod} yields $p_{-}(L)=1_{*}$ and so we could try to apply the previous theory. However, we wish to give a complete argument with the minimal tools.

\begin{prop} 
\label{prop: pointwise estimates} 
For all integers $\beta> \nicefrac{n}{2}+2$ the following properties hold for the kernel $H_t^\beta(x,y)$ of the higher-order resolvents $(1-t^2a^{-1}\Delta_x)^{-\beta} a^{-1}$. \index{Theorem!McIntosh--Nahmod's}\index{kernel bounds!for perturbations of the Laplacian} 
\begin{enumerate}
\item There are $C,c>0$, depending on ellipticity, dimensions and $\beta$, such that for for all $t>0$ and $x,y\in \R^n$, 
\begin{align*}
\qquad |H_t^\beta(x,y)|+ |t\nabla_x H_t^\beta(x,y)|+ |t\nabla_y H_t^\beta(x,y)|
\\
+  |t^2\nabla_x \nabla_y H_t^\beta(x,y)| \le Ct^{-n}\e^{-\frac{c|x-y|}{t}}.
\end{align*}
\item For all $\eta \in (0,1)$, the kernels
\begin{align*}
\qquad t\nabla_x H_t^\beta(x,y), \quad  t\nabla_y H_t^\beta(x,y), \quad t^2\nabla_x \nabla_y H_t^\beta(x,y)
\end{align*}
are H\"older continuous in both variables with exponent $\eta$ and norms in this space of the order of $t^{-\eta-n}$.  In particular, 
$H_t^\beta \in \C^{1, \eta}(\R^n \times \R^n)$, the space of $\C^1$-functions having H\"older continuous first-order derivatives of exponent $\eta$.
\end{enumerate}
\end{prop}

\begin{proof} We set $L \coloneqq -a^{-1}\Delta_x$ and $L_{0} \coloneqq -\Delta_{x}$ is acting componentwise on $\IC^m$-valued functions. It suffices to prove the properties of $H_{t}^\beta$ for $t=1$ with implicit constants that depend on dimensions and ellipticity. Indeed, a change of variables yields that $H_{t}^\beta(x,y)= t^{-n}\wt{H}_{1}^\beta (\nicefrac{x}{t}, \nicefrac{y}{t})$, where $\wt{H}^\beta_1$ corresponds to the coefficients $a_t(x) \coloneqq a(tx)$, which have the same ellipticity constant as $a$. We split the proof into four steps.

\medskip

\noindent \emph{Step~1: Pointwise estimates for $H^\beta$.} Let $s>0$. When  $m=1$,  $(1-s\Delta_{x})^{-1}$ is given by convolution with a classical Bessel potential, that is, a positive function with integral $1$ that is in $\L^r$ whenever $\nicefrac{1}{r} > 1 - \nicefrac{2}{n}$, see for instance \cite[Sec.V.3]{Stein-Singular}. When $m \ge2$, $(1+sL_0)^{-1}$ is given by componentwise convolution with the same potential.

By positivity, we get for $f \in \L^2$ and $s>0$ the pointwise bound 
\begin{align}
\label{eq:pointwisebound}
|(1+ sL_{0})^{-1}f| \le (1-s\Delta_{x})^{-1}|f|,
\end{align}
where $|\cdot|$ is the $\IC^m$-norm and the resolvent on the right-hand side is scalar-valued. In particular, $(1+ sL_{0})^{-1}$ is a contraction on $\L^2$.

We can write 
\begin{align}
\label{eq:accr-representation}
	a= \tau(1-b)
\end{align} 
for some $\tau \ > 0$ and $b \in \L^\infty(\R^n;\Lop(\IC^m))$ with $\|b\|_\infty  < 1$. We shall give the well-known argument in the final step of the proof. 

Using the above decomposition of $a$,  we find
\begin{align*}
(a+L_{0})^{-1} 
= \tfrac{1}{\tau} (1- (1+ \tfrac{1}{\tau}L_0)^{-1}b)^{-1}(1+ \tfrac{1}{\tau}L_0)^{-1}
\end{align*}
as operators on $\L^2$, and the first term on the right can be computed by a Neumann series. Expanding this series explicitly and applying \eqref{eq:pointwisebound} inductively with $s= \nicefrac{1}{\tau}$, we have 
\begin{align*}
	|(a+L_{0})^{-1} f |  \leq \tfrac{1}{\tau} \sum_{k=0}^\infty ((1 -\tfrac{1}{\tau}\Delta_{x})^{-1} \|b\|_{\infty})^k(1 -\tfrac{1}{\tau}\Delta_{x})^{-1}|f|,
\end{align*}
so that summing backward, we obtain the pointwise bound
\begin{align*}
	|(a+L_{0})^{-1} f |  \leq (\alpha -\Delta_{x})^{-1} |f|, \quad \text{where}\quad \alpha=\tau(1-\|b\|_{\infty}).  
\end{align*}
Applying this estimate to $af$ in place of $f$, we get
\begin{align}
\label{eq:a pointwise bound}
	|(1+L)^{-1} f | 
\leq  \|a\|_\infty (\alpha-\Delta_x)^{-1}  |f|
\end{align}
and obtain extensions by density with operator-norm bounds
\begin{align}
\label{eq1: a pointwise bound}
\| (1+L)^{-1} \|_{\L^p\to\L^q} 
\le  \|a\|_\infty \| (\alpha-\Delta_{x})^{-1} \|_{\L^p\to\L^q},
\end{align}
whenever $1 \leq p \leq q \leq \infty$, not both infinite, and the right-hand side is finite. By Young's inequality for convolutions, this is the case if $\nicefrac{1}{p}-\nicefrac{1}{q}<\nicefrac{2}{n}$. This gives $\L^1-\L^\infty$-boundedness of $(1+L)^{-\beta} a^{-1}$ provided that $\beta> \nicefrac{n}{2}$. 

By the Dunford--Pettis theorem~\cite[Thm.~1.3]{Arendt-Bukhvalov}  we obtain that $(1+L)^{-\beta} a^{-1}$ is given by a bounded kernel $H^\beta(x,y)$ and the bound depends only on dimensions and ellipticity. Iterating \eqref{eq:a pointwise bound}, we see that $|H^\beta(x,y)|$ is dominated by the kernel of $(\alpha-   \Delta_x)^{-\beta}$ up to a factor $\|a\|_{\infty}^\beta \|a^{-1}\|_{\infty}$. The latter operator is given by convolution with a higher-order Bessel potential. Since $\beta > \nicefrac{n}{2}$, we get exponential decay and no singularity at $x=y$ as required in (i), see \cite[Sec.V.3]{Stein-Singular}.

\medskip

\noindent \emph{Step~2: Proof of (ii) and the other bounds in (i) with $c=0$.}
Write
\begin{align*}
	\Delta_x(1+L)^{-1}a^{-1} = - 1+a(1+L)^{-1}a^{-1},
\end{align*}
the Laplacian acting componentwise, so that 
\begin{align}
\label{eq2: a pointwise bound}
\begin{split}
\| \Delta_x(1+L)^{-1}&a^{-1} \|_{\L^p\to\L^p} \\
&\le 1+\|a\|_\infty \| (\alpha-\Delta_{x})^{-1} \|_{\L^p\to\L^p}.
\end{split}
\end{align} 

The operator norm in the line above is controlled for all $p\in [1,\infty]$. If $1<p<\infty$, then by the Mihlin multiplier theorem \eqref{eq1: a pointwise bound} and \eqref{eq2: a pointwise bound} imply that $(1+L)^{-1} a^{-1}$ and  $\nabla_x(1+L)^{-1} a^{-1}$ are bounded from $\L^p$ to $\W^{1,p}$. In particular, for $p>n$, we have the inhomogeneous Sobolev embedding $\W^{1,p}\subseteq \Lamdot^\eta \cap \L^\infty$, $\eta=1-\nicefrac{n}{p}$. The same applies with $a^*$ in place of $a$ and by duality $(1+L)^{-1} a^{-1} \div_x$  is bounded from $\L^1$ to $\L^{p'}$. By composition, we obtain that for $\beta > \nicefrac{n}{2} +2$ the operators
\begin{align*}
&\nabla_x(1+L)^{-\beta}a^{-1}, \quad  -(1+L)^{-\beta}a^{-1}\div_x, \quad -\nabla_x(1+L)^{-\beta} a^{-1}\div_x  
\end{align*}
are bounded from $\L^1$ into $\Lamdot^\eta \cap \L^\infty$. In particular they are bounded from $\L^1$ into $\L^\infty$ and, invoking again the Dunford--Pettis theorem, they correspond to the kernels $\nabla_x H^\beta(x,y)$, $\nabla_y H^\beta(x,y)$, $\nabla_x \nabla_y H^\beta(x,y)$, which therefore are bounded measurable functions. 

We can then use the mapping properties from $\L^1$ into $\Lamdot^\eta$ and once more the Dunford--Pettis theorem, in order to obtain first H\"older continuity of the kernels in $x$ (with any exponent $\eta \in (0,1)$), uniformly in $y$, and then by duality the same with the roles of $x$ and $y$ reversed. Such type of argument has appeared earlier on in the proof of Lemma~\ref{lem: half-equivalence}. This proves (ii) and finishes the proof of (i) with $c=0$.

\medskip

\noindent{\emph{Step~3: Exponential decay for the other kernels.}} We begin with $\partial_{x_i} H^\beta$, where $1 \leq i \leq n$. Let $e_i \in \R^n$ be the $i$-th standard unit vector and let $h>0$. By the fundamental theorem of calculus we have
\begin{align*}
	&\frac{1}{h} (H^\beta(x+he_i,y) - H^\beta(x,y))  \\
	&=\quad \barint_0^h \partial_{x_i} H^\beta(x+s e_i,y) \, \d s \\
	&=\quad \partial_{x_i}H^\beta(x,y) + \barint_0^h \partial_{x_i} H^\beta(x+s e_i,y) - \partial_{x_i}H^\beta(x,y) \, \d s,
\end{align*}
where $x,y \in \R^n$. If  $2|h| \leq |x-y|$, then $|x+h e_i - y| \geq \nicefrac{|x-y|}{2}$ and we get from (i) for $H^\beta$ and (ii) for $\nabla_x H^\beta$ that
\begin{align*}
	|\partial_{x_i}H^\beta(x,y)| \leq \frac{C}{h}\big(\e^{-\frac{c}{2}|x-y|} + \e^{-c|x-y|}\Big) + h^\eta \|\nabla_x H^\beta(\cdot, y)\|_{\Lamdot^\eta}. 
\end{align*}
Since in Step~2 we have already obtained a uniform bound for $\partial_{x_i} H^\beta$, it suffices to prove the decay for $|x-y|$ large, say $|x-y| \e^{\frac{c}{4}|x-y|} \geq 2$. This restriction is manufactured such that we can take $h \coloneqq \e^{-\frac{c}{4}|x-y|}$, resulting in the desirable estimate
\begin{align*}
		|\partial_{x_i}H^\beta(x,y)| \leq C' \e^{-\frac{\eta c}{4}|x-y|}
\end{align*}
for some new constant $C'$ that depends on $a$ only through ellipticity. This completes the proof for $\nabla_x H^\beta$.

The argument above has only used the exponential decay for $H^\beta$, the uniform boundedness of $\nabla_x H^\beta$ and the $\Lamdot^\eta$-estimate for $\nabla_x H^\beta$ in the $x$-variable uniformly in the $y$-variable, in order to give exponential decay for $\nabla_x H^\beta$. Thus, it can be repeated \emph{verbatim} for the decay of $\nabla_y H^\beta$. Then, replacing $H^\beta$ by $\nabla_y H^\beta$ gives decay of $\nabla_x \nabla_y H^\beta$.

\medskip

\noindent \emph{Step~4: Proof of \eqref{eq:accr-representation}.} We let $\tau \coloneqq \lambda^{-1} \|a\|_\infty^{2}$  and  $b \coloneqq 1- \tau^{-1}a$. If $\xi  \in  \IC^m$ is normalized to $|\xi| = 1$, then 
\begin{align*}
	|b(x)\xi|^2 
	&= 1 + \tau^{-2}|a(x) \xi |^2 - 2 \tau^{-1} \Re \langle a(x) \xi, \xi \rangle \\
	&\leq 1 + \tau^{-2} \|a\|_\infty^2 - 2 \tau^{-1} \lambda \\
	&= 1 - \lambda^{2} \|a\|_\infty^{-2}<1. \qedhere
\end{align*}
\end{proof}
\section{Comparison with the Auscher--Stahlhut interval}
\label{sec: Comparison with Auscher-Stahlhut interval}

\noindent The identification of adapted Hardy spaces as a key tool to treating boundary value problems has appeared first in  \cite{AusSta}. Although we argue independently of this reference concerning this particular issue, we need to make the bridge and the results of this section are explicitly used in Section~\ref{sec: Neumann problems} on Neumann problems.

In \cite[Thm.~5.1]{AusSta} an interval of values of $p$ is constructed, where one has the identification $\IH^p_{DB}=\IH^p_{D}$ even for more general operators $DB$. (The matrix $B$ need not be block-diagonal.) Its upper endpoint is denoted by $p_{+}(DB)$ and the lower endpoint is at most the lower Sobolev conjugate of another exponent $p_{-}(DB)$. To avoid confusion, we denote these exponents by $p_{\pm}^{\mathrm{AS}}(DB)$ here.
They have a precise meaning that we recall next. The following material is all taken from \cite[Sec.~3.2]{AusSta}. 

Let
\begin{align*}
\dom_p(D) \coloneqq \{f \in \L^p: Df \in \L^p\},
\end{align*}
where $\L^p = \L^p(\R^n; \IC^m \times \IC^{mn})$ and the action of $D$ is in the sense of distributions. Then $BD$ is defined as an unbounded operator in $\L^p$ with domain $\dom_p(BD) = \dom_p(D)$, null space $\nul_p(BD)$ and range $\ran_p(BD)$. Similar to Definition~\ref{def: p-lower bound}, one introduces the set of exponents with $p$-lower bounds\index{plowerbounds@$p$-lower bounds!for $B$}
\begin{align*}
\cI(BD) &\coloneqq \{p \in (1,\infty) : \|Bf\|_p \gtrsim \|f\|_p \text{ for all } f \in \ran_p(D) \}
\end{align*}
and the analogous set with $B^*$ replacing $B$. They are open but possibly non-connected and $\cI_2$ denotes the connected component of $\cI(BD) \cap \cI(B^*D)'$ that contains $p=2$. Here, $I' = \{p' : p \in I\}$ is the dual set of a given $I \subseteq (1,\infty)$. In passing, we point out that the use of $\cI(B^*D)$ instead of its dual set in \cite[Rem.~3.5]{AusSta} is a typo that does not appear in the original reference \cite[Sec.~5]{AusStaRemarks}.

Then  $(p_{-}^{\mathrm{AS}}(DB), p_{+}^{\mathrm{AS}}(DB))$ is the interval of exponents $q \in \cI_2$ such that for all $p$ between $2$ and $q$ there is a topological decomposition
\begin{align}
\label{eq: AusSta kernel/range splitting}
\L^p = \nul_p(BD) \oplus \cl{\ran_p(BD)},
\end{align}
see \cite[Thm.~3.6]{AusSta}. It is proved in \cite[Thm.~5.1]{AusSta} that for $(p_{-}^{\mathrm{AS}}(DB))_{*}<p< p_{+}^{\mathrm{AS}}(DB)$ one has $\IH^p_{DB}=\IH^p_{D}$. It was not proved that this interval is optimal for the identification in the class of $DB$-operators there and for some examples it was shown that this is not the case, especially for the lower endpoint. Hence, \cite{AusSta} does not provide the whole identification interval, yet \cite[Prop.~6.4 \& 6.5]{AusSta} there describe it as an open interval.

Using the same framework as \cite{AusSta}, it became clear  in the classification theorems of \cite{AM} as well as in the uniqueness statements of \cite{AE} that the full interval of identification is the object of interest. Both references introduce the set of exponents $p \in (1_*, p_{+}^{\mathrm{AS}}(DB))$ for which $\IH^p_{DB}=\IH^p_{D}$ holds with equivalent $p$-quasinorms. It is called $I_L$ in \cite{AM} and $\cH_{\cL}$ in \cite{AE}. Hence, either of these intervals is of the form\index{interval!of Auscher--Stahlhut}\index{interval!of Auscher--Mourgoglou}
\begin{align*}
	(a^{\mathrm{AS}}(DB), p_{+}^{\mathrm{AS}}(DB)) 
\end{align*}
for some number $a^{\mathrm{AS}}(DB) \geq 1_*$ which could be in particular less than $(p_{-}^{\mathrm{AS}}(DB))_{*}$.

In the block situation of this monograph, we proceeded differently and introduced the set of identification $\cH(DB)$ in \eqref{eq: H(DB)} directly as the largest set of exponents $p \in (1_*, \infty)$ for which $\IH^p_{DB}=\IH^p_{D}$ holds with equivalent $p$-quasinorms. Then we proved that it is an open interval and characterized its endpoints as $h_{-}(DB) = p_{-}(L)$ and $h_+(DB) = q_{+}(L)$, see Theorem~\ref{thm: main result Hardy}, Corollary~\ref{cor: H1L upper endpoint} and Theorem~\ref{thm: Riesz characterization HpL}. Hence, in order to be able to apply the results in \cite{AM, AE} within the interval of identification $\cH(DB)$, we need to connect both approaches.  

The discussion above already shows that $a^{\mathrm{AS}}(DB) = p_-(L)$ and $q_+(L) = h_+(DB) \geq p_{+}^{\mathrm{AS}}(DB)$. Identifying the upper endpoints requires a specific argument.

\begin{prop} 
\label{prop: AS number}
In the block case setting of this monograph the number $p_{+}(DB) = p_{+}^{\mathrm{AS}}(DB)$ of \cite{AusSta} coincides with $h_{+}(DB) = q_+(L)$.
\end{prop}

\begin{proof}
As said, it remains to prove $q_+(L) \leq p_{+}^{\mathrm{AS}}(DB)$.

Let $2 \leq p < q_+(L)$. First, we recall that $q_+(L)=q_+(L_{0})$ from Theorem~\ref{thm: critical numbers a-independent}, so that by Theorem~\ref{thm: endpoints for Hodge projector} we have $p \in \cP(L_{0})$. Hence, Proposition~\ref{prop: P(L) vs Hodge} implies  $p$-lower bounds for $d$ and $p'$-lower bounds for $d^*$, as well as  the topological Hodge decomposition \eqref{eq: Hodge}.

To reinterpret this, we recall that
\begin{align*}
B = \begin{bmatrix} a^{-1} & 0  \\ 
0 & d \end{bmatrix},
\quad 
D = \begin{bmatrix} 0 & \div_x  \\ 
-\nabla_x & 0 \end{bmatrix}.
\end{align*}
Using the notation of Section~\ref{sec: Hodge}, we have
\begin{align*}
\nul_p(D) &= \{0\} \times \nul_p(\div_x),\quad \cl{\ran_p(D)} = \L^p \times \cl{\ran_p(\nabla_x)}.
\end{align*}
Since $a^{-1}$ is strictly elliptic, we see that the conditions $p \in \cI(BD)$ and $p \in \cI(B^*D)'$ are equivalent to $p$-lower bounds for $d$ and $p'$-lower bounds for $d^*$, respectively. Moreover, using the $p$-lower bounds to determine the null space,  we find
\begin{align*}
\nul_p(BD) &= \{0\} \times \nul_p(\div_x),\quad \cl{\ran_p(BD)} = \L^p \times \cl{d\ran_p(\nabla_x)}.
\end{align*}
In turn, this shows that \eqref{eq: AusSta kernel/range splitting} is equivalent to the Hodge decomposition \eqref{eq: Hodge}. 

Altogether, we have shown that $p \in \cI(BD) \cap \cI(B^*D)'$ as well as the Hodge decomposition \eqref{eq: AusSta kernel/range splitting}. As we have done this for all $p$ in the interval $ [2,q_{+}(L))$, this proves that $q_{+}(L) \leq p_{+}^{\mathrm{AS}}(DB)$.
\end{proof}

Summarizing, we have obtained

\begin{cor}
\label{cor: AE interval and AM interval}
In the block case setting of this monograph the open intervals $\cH_\Le$ from \cite{AE} and $I_L$ from \cite{AM} both equal $(p_-(L), q_+(L))$.\index{interval!of Auscher--Stahlhut}\index{interval!of Auscher--Mourgoglou}
\end{cor}
\section{Basic properties of weak solutions}
\label{sec: Basic properties of weak solutions}

\noindent At this point in the monograph we begin to slightly change our perspective from Hardy spaces adapted to $L = -a^{-1} \div_x d \nabla_x$ to weak solutions to the elliptic system
\begin{align}
\label{eq: elliptic equation in weak sol section}
\Le u = -\div(A \nabla u) = -\partial_t (a \partial_t u) - \div_x d \nabla_x u = 0 
\end{align}
in $\ree$, where as before we write
\begin{align*}
A = \begin{bmatrix} a & 0 \\ 0 & d \end{bmatrix}
\end{align*}
for the coefficient matrix in block form. In this section, we gather well-known properties of weak solutions that will frequently be used in the further course.

As usual, a \emph{weak solution}\index{solution!weak} to the equation
\begin{align*}
\Le u = g \in \Lloc^2(O)
\end{align*}
in an open set $O \subseteq \ree$ is a function $u \in \Wloc^{1,2}(O; \IC^m)$ that satisfies for all $\phi \in \C_0^\infty(O; \IC^m)$,
\begin{align*}
\iint_{O}  A \nabla u \cdot \cl{\nabla \phi} \, \d t \d x = \iint_O g \cdot \cl{\phi} \, \d t \d x.
\end{align*}
By density we can allow any test function $\phi \in \W^{1,2}(O; \IC^m)$ with compact support in $O$.
\subsection{Energy solutions}
\label{subsec: Energy solutions}

The most common construction of weak solutions is by the Lax--Milgram lemma, using the \emph{energy class}\index{energy class!$\Wdot^{1,2}(\reu)$}
\begin{align*}
	\Wdot^{1,2}(\reu)
	\coloneqq \{v \in \Lloc^2(\reu): \nabla v \in \L^2(\reu) \} / \IC^m.
\end{align*}
This is a Hilbert space for the inner product $\langle \nabla \cdot \, , \nabla \cdot \rangle$ and it contains the restrictions of $\C_0^\infty(\ree)$-functions to $\reu$ as a dense subspace, see for instance~\cite[Lem.~3.1]{AMcM}. By Poincar\'{e}'s inequality it embeds continuously into $\Lloc^2(\reu)$.

We recall the well-known trace and extension results. For convenience and a later use we include elementary proofs in our homogeneous Sobolev setting.\index{energy class!trace}
\begin{lem}
\label{lem: energy trace}
Every equivalence class $v \in \Wdot^{1,2}(\reu)$ has a representative that is continuous on $[0,\infty)$ with values in $\Lloc^2$. In this sense $v \in \C_0([0,\infty); \Hdot^{\nicefrac{1}{2},2})$ and
\begin{align*}
 \sup_{t \geq 0} \|v(t,\cdot)\|_{\Hdot^{\nicefrac{1}{2},2}} \lesssim \|\nabla v\|_2.
\end{align*}
Conversely, every $f \in \Hdot^{\nicefrac{1}{2},2}$ can be extended to a function $v \in  \Wdot^{1,2}(\reu)$ with $v(0,\cdot) = f$ and $\|\nabla v\|_2 \simeq  \|f\|_{\Hdot^{\nicefrac{1}{2},2}}$.
\end{lem}

\begin{proof}
That $v$ has a representative that is continuous on $[0,\infty)$ valued in $\Lloc^2$ is just the one-dimensional Sobolev embedding in the $t$-variable. This property is not affected by adding constants to $v$ and amounts to re-defining $v$ a.e.\ on $\reu$.

By density it suffices to prove the embedding into $\C_0([0,\infty); \Hdot^{\nicefrac{1}{2},2})$ in the case that $v$ is the restriction of a function in $\C_0^\infty(\ree)$. For all $t \geq 0$ we have
\begin{align*}
	\frac{\d }{\d t} \| (-\Delta_x)^{1/4} v(t,\cdot)\|_2^2
	&= 2 \Re \langle (-\Delta_x)^{1/4} v(t,\cdot), (-\Delta_x)^{1/4} \partial_t v (t,\cdot) \rangle \\
	& \leq 2 \| (-\Delta_x)^{1/2} v(t,\cdot) \|_2 \| \partial_t v(t,\cdot) \|_2 \\
	&\lesssim  \| \nabla_x v(t,\cdot)\|_2 ^2 +  \| \partial_t  v(t,\cdot)\|_2 ^2,
\end{align*}
where the final step is by the solution of the Kato problem. Integration in $t$ up to $\infty$ gives 
\begin{align*}
\| (-\Delta_x)^{1/4} v(t,\cdot)\|_2^2 \lesssim \|\nabla v\|_2^2 \quad (t \geq 0)
\end{align*}
and the left-hand side is comparable to $\|v(t,\cdot)\|_{\Hdot^{1/2,2}}^2$ by Corollary~\ref{cor: Kato}.
	
Again by density it suffices to prove the extension part for $f \in \Hdot^{\nicefrac{1}{2},2} \cap \L^2$. We set $v(t,\cdot) \coloneqq \e^{-t (-\Delta_x)^{1/2}}f$. Clearly $v$ is continuous on $[0,\infty)$ valued in $\L^2$ with $v(0,\cdot) = f$. Moreover, we have
\begin{align}
\label{eq: energy solution with Laplace}
\begin{split}
\|\nabla v\|_2^2
&= \int_0^\infty \|\partial_t v(t,\cdot)\|_2^2 + \|\nabla_x v(t,\cdot)\|_2^2 \, \d t \\
&= \int_0^\infty  \|(-\Delta_x)^{1/2}  \e^{-t (-\Delta_x)^{1/2}} f\|_2^2 + \|\nabla_x  \e^{-t (-\Delta_x)^{1/2}} f\|_2^2 \, \d t\\
& \simeq \int_0^\infty \|(-\Delta_x)^{1/2}  \e^{-t (-\Delta_x)^{1/2}} f\|_2^2 \, \d t \\
&= \int_0^\infty \|(-t^2 \Delta_x)^{1/4}  \e^{-(-t^2 \Delta_x)^{1/2}} (-\Delta_x)^{1/4} f\|_2^2 \, \frac{\d t}{t} \\
&\simeq \|(-\Delta_x)^{1/4} f\|_2^2 \\
&\simeq \|f\|_{\Hdot^{1/2,2}}^2,
\end{split}
\end{align}
where the fourth step is by McIntosh's theorem.
\end{proof}

We also obtain the usual characterization of the subspace with trace zero at the boundary.

\begin{lem}
\label{lem: energy zero trace}
The subspace\index{energy class!with trace zero ($\Wdot_0^{1,2}(\reu)$)}
\begin{align*}
	\Wdot_0^{1,2}(\reu) \coloneqq \{u\in \Wdot^{1,2}(\reu) : u(0,\cdot) = 0 \text{ in } \Hdot^{\nicefrac{1}{2},2} \}
\end{align*}
coincides with the closure of $\C_0^\infty(\reu)$ in $\Wdot^{1,2}(\reu)$.
\end{lem}

\begin{proof}
Since the restriction $R: \Wdot^{1,2}(\reu) \to \Hdot^{\nicefrac{1}{2},2}$ to $t = 0$ is bounded, $\Wdot_0^{1,2}(\reu)$ is a closed subspace and it contains $\C_0^\infty(\reu)$. 

Conversely, let $u \in \Wdot_0^{1,2}(\reu)$. Let $E: \Hdot^{\nicefrac{1}{2},2}\to  \Wdot^{1,2}(\reu)$ be the extension operator from the proof of Lemma~\ref{lem: energy trace}. We pick a sequence $(u_k) \subseteq \C_0^\infty(\ree)$ with $u_k \to u$ in $\Wdot^{1,2}(\reu)$ as $k \to \infty$ and set $v_k \coloneqq (1-ER)u_k$. Then $Rv_k = 0$ and $v_k \to u$ in $\Wdot^{1,2}(\reu)$. Therefore it suffices to approximate each $v_k$ by $\C_0^\infty(\reu)$-functions. In fact, it suffices to find approximants with compact support in $\reu$ since then we can conclude via convolution with smooth kernels.

To this end, we note that $R u_k \in \L^2$ together with the explicit construction of $E$ implies $v_k \in \C_0([0,\infty); \L^2)$ with $v_k(0, \cdot) = 0$. Extending $v_k$ to $\ree$ by $0$ and using the $\L^2$-continuity of the translation in the $t$-direction, we obtain approximants $w_k$ with the same properties that have their support in $\reu$. Now, we take $\eta \in \C_0^\infty(\ree)$ with $\eta(0,0) = 1$ and set $\eta_\eps(t,x) \coloneqq \eta(\eps t, \eps x)$. We can bound
\begin{align*}
\|\nabla (\eta_\eps w_k) - \nabla w_k\|_{2}
\lesssim  \|(1-\eta_\eps) \nabla w_k\|_{2} + \eps^{\frac{1}{2}} \|w_k\|_{\L^\infty((0,\infty); \L^2)}.
\end{align*}
In the limit as $\eps \to 0$ the first term on the right vanishes by dominated convergence, whereas the second one vanishes thanks to the additional information $w_k \in \L^\infty((0,\infty); \L^2)$.
\end{proof}

We can now use the Lax--Milgram lemma to prove the following well-posedness result. Neither the block structure of $A$ nor its $t$-independence are needed in the argument. We call $u$ the \emph{energy solution}\index{solution!energy with Dirichlet datum} to $\Le u = 0$ in $\reu$ with Dirichlet data $f$.

\begin{prop}
\label{prop: existence of energy solution}
For all $f \in \Hdot^{\nicefrac{1}{2},2}$ there exists a unique solution $u$ (modulo constants) to the problem
\begin{equation*}
\begin{cases}
\Le u=0   & (\text{in }\reu), \\
\nabla u \in \L^2(\reu),   \\
u(0,\cdot) = f & (\text{in } \Hdot^{\nicefrac{1}{2},2}).
\end{cases}
\end{equation*}
Moreover, $\|\nabla u\|_2 \lesssim  \|f\|_{\Hdot^{1/2,2}}$ and $\lim_{t \to \infty} u(t,\cdot) = 0$ in $ \Hdot^{\nicefrac{1}{2},2}$.\index{solution!operator for the energy class}
\end{prop}

\begin{proof}
If $u$ is any solution, then we obtain by density and Lemma~\ref{lem: energy zero trace} that
\begin{align*}
\iint_{\reu} A \nabla u \cdot \cl{\nabla \phi} \,  \d t \d x = 0 \quad (\phi \in \W^{1,2}_0(\reu)).
\end{align*}	
Since $A$ is elliptic, $u \in \W^{1,2}_0(\reu)$ implies $\nabla u = 0$. Hence, solutions are unique modulo constants. In order to construct a solution, let $v \in \Wdot^{1,2}(\reu)$ be an extension of $f$ as in Lemma~\ref{lem: energy trace}. By the Lax--Milgram lemma, there exists $w \in \W^{1,2}_0(\reu)$ solving
\begin{align*}
 \iint_{\reu} A \nabla w \cdot \cl{\nabla \phi} \,  \d t \d x = - \iint_{\reu} A \nabla v \cdot \cl{\nabla \phi} \,  \d t \d x \quad (\phi \in \W^{1,2}_0(\reu)).
\end{align*}
Hence, $u \coloneqq v+w$ is a solution to the given problem and Lemma~\ref{lem: energy trace} yields the limit at $t=\infty$ as well as the bound
\begin{align*}
\|\nabla u \|_2 \leq \|\nabla v\|_2 + \|\nabla w\|_2 \lesssim \|\nabla v\|_2 \simeq \|f\|_{\Hdot^{1/2,2}}. &\qedhere
\end{align*}
\end{proof}
\subsection{Semigroup solutions}
\label{subsec: semigroup solutions}

In the specific situation of coefficients in block form, we can also use the Poisson semigroup for $L$ to construct weak solutions. Here, the natural boundary space is $\L^2$ rather than $\Hdot^{\nicefrac{1}{2},2}$.

\begin{prop}
\label{prop: Poisson smg is weak solution}
Let $f \in \L^2$. Then $u(t,x) \coloneqq \e^{-t {L^{1/2}}}f(x)$ is a weak solution to $\Le u = 0$ in $\reu$ of class $\C_0([0,\infty); \L^2) \cap \C^\infty((0,\infty); \L^2)$ with $u(0,\cdot) = f$.\index{solution!semigroup}
\end{prop} 

\begin{proof}
The regularity in $t$ follows directly from the functional calculus. In particular, $u(t,\cdot)$ is in the domain of $L$ for every $t>0$ and $\frac{\d^2}{\d t^2} u(t,\cdot)= Lu(t,\cdot)$. Since $a$ is bounded and independent of $t$, the function $au$ has the same properties and we have $\frac{\d^2}{\d t^2} (a u(t,\cdot))= a Lu(t,\cdot)$. Let now $\phi \in \C_0^\infty(\reu)$. For any $t>0$, the Lax--Milgram interpretation of $aL$ in \eqref{eq: Lax-Milgram operator} yields
\begin{align*}
\int_{\R^n} \frac{\d^2}{\d t^2} (a u(t, \cdot)) \cdot \cl{\phi(t,\cdot)} \, \d x
= \int_{\R^n} d \nabla_x u(t,\cdot) \cdot \cl{\nabla_x \phi(t,\cdot)} \, \d x
\end{align*} 
and the claim follows by integrating both sides in $t$ and then integrating by parts in $t$ on the left-hand side.
\end{proof}

We have the following \emph{compatibility}\index{solution!compatible} between semigroup and energy solutions. This could be deduced from more general results in \cite{AMcM} but in the block situation there is a particularly simple proof. 

\begin{prop}
\label{prop: smg solution is compatible}
If $f \in \Hdot^{\nicefrac{1}{2},2} \cap \L^2$, then $u(t,x) \coloneqq \e^{-t L^{1/2}}f(x)$ is the energy solution to $\Le u = 0$ in $\reu$ with Dirichlet data $f$.
\end{prop}

\begin{proof}
We already know that $u$ is a weak solution to $\Le u = 0$ in $\reu$ with $u(0,\cdot) = f$ in the sense of $\C_0([0,\infty); \L^2)$. Furthermore, $\nabla u \in \L^2(\reu)$ follows by a literal repetition of the argument in \eqref{eq: energy solution with Laplace}, replacing $-\Delta_x$ by $L$ at each occurrence. In fact, this is why we have justified \eqref{eq: energy solution with Laplace} by abstract arguments instead of using the Fourier transform.
\end{proof}
\subsection{Interior estimates}
\label{subsec: interior estimates}

We continue with the standard interior estimates. All this is well-known but precise references for systems with our ellipticity assumption are hard to find. One is \cite{Barton}, where even systems of higher order are treated, but for the reader's convenience we include the simple arguments in the second-order case. Again the block structure of $A$ and its $t$-independence are not needed for this part.

We call $W \subseteq \reu$ a \emph{cylinder of radius $r$} if $W = I \times B$, where $I \subseteq (0,\infty)$ is an interval of length $r$ and $B \subseteq \R^n$ a ball of radius $r$ (or a cube of sidelength $r$). As usual, we write $\alpha W$ for the concentrically scaled version of $W$. 

\begin{lem}[Caccioppoli\index{inequality!Caccioppoli}]
\label{lem: Caccioppoli}
Let $O \subseteq \ree$ be open, $g \in \Lloc^2(O)$ and $u$ a weak solution to $\Le u = g$ in $O$. Let $W \subseteq \ree$ a cylinder of radius $r$ and $\alpha > 1$ be such that $\alpha W \Subset O$. Then there is a constant $C$ depending on dimensions, ellipticity and $\alpha$, such that
\begin{align*}
\iint_{W} |\nabla u|^2 \, \d s \d y \leq C \iint_{\alpha W} r^{-2} |u|^2 + r^2 |g|^2 \, \d s \d y.
\end{align*}
\end{lem}

\begin{proof}
Fix $\eta \in \C_0^\infty(\ree)$ with $\ind_W \leq \eta \leq \ind_{\alpha W}$ and $|\nabla \eta| \leq c r^{-1}$ for a constant $c$ depending on $n$ and $\alpha$. We write $\langle \cdot \,, \cdot \rangle$ for the inner product on $\L^2(\ree)$. By ellipticity and multiple applications of the product rule, we have
\begin{align*}
\lambda \|\nabla (\eta u)\|_2^2
&\leq |\langle A \nabla (\eta u), \nabla (\eta u) \rangle| \\
&\leq |\langle A \nabla u, \nabla (\eta^2 u) \rangle| + |\langle \eta A \nabla u, u \otimes \nabla \eta \rangle| \\
&\quad + |\langle A (u \otimes \nabla \eta), \nabla (\eta u) \rangle|\\
&\eqqcolon \I_1 + \I_2 + \I_3,
\end{align*}
where our notation is $\nabla(\eta u) \coloneqq \eta \nabla u + u \otimes \nabla \eta$ in the sense prescribed by the product rule. By the equation for $u$, the Cauchy--Schwarz inequality and the elementary bound $xy \leq \frac{\eps}{2} x^2 + \frac{1}{2\eps} y^2$ for positive numbers $x,y,\eps$, we have
\begin{align*}
|\I_1| = |\langle g, \eta^2 u \rangle| \leq \frac{1}{2r^2} \|\eta u\|_2^2 + \frac{r^2}{2} \|\eta g\|_2^2.
\end{align*} 
Similarly, we get
\begin{align*}
\I_2 + \I_3 
&\leq \frac{\lambda}{4} \|\eta \nabla u\|_2^2 + \frac{\lambda}{4} \|\nabla (\eta u)\|_2^2 + C \|u \otimes \nabla \eta\|_2^2 \\
&\leq \frac{3\lambda}{4} \|\nabla (\eta u)\|_2^2 + \bigg(C + \frac{\lambda}{2}\bigg) \|u \otimes \nabla \eta\|_2^2,
\end{align*}
where $C$ depends on dimensions and ellipticity. Rearranging terms leads to 
\begin{align*}
\frac{\lambda}{4} \|\nabla (\eta u)\|_2^2 \leq \Big(C+\frac{\lambda}{2}\Big) \|u \otimes \nabla \eta\|_2^2 + \frac{1}{2r^2} \|\eta u\|_2^2 + \frac{r^2}{2} \|\eta g\|_2^2
\end{align*}
and by choice of $\eta$ we are done.
\end{proof}

\begin{lem}[Reverse H\"older\index{inequality!reverse H\"older}]
\label{lem: RH}
Let $u$ be a weak solution to $\Le u = 0$ in an open set $O \subseteq \ree$ and let $\alpha > 1$. There is a constant $C$ depending on dimensions, ellipticity and $\alpha$, such that for all cylinders $W$ with $\alpha W \Subset O$ it follows that
\begin{align*}
\bigg(\bariint_{W} |\nabla u|^2 \, \d s \d y \bigg)^{1/2} \leq C \bariint_{\alpha W} |\nabla u| \, \d s \d y.
\end{align*}
Moreover, with $q \coloneqq \frac{2(n+1)}{n-1}$ in dimension $n \geq 2$ and $q \in (2, \frac{2(n+1)}{n-1})$ arbitrary in dimension $n=1$, it follows that
\begin{align*}
\bigg(\bariint_{W} |u|^q \, \d s \d y \bigg)^{1/q} \leq C \bariint_{\alpha W} |u| \, \d s \d y,
\end{align*}
where $C$ also depends on $q$.
\end{lem}

\begin{proof}
We begin with the first inequality. Let $c \coloneqq \bariint_{\alpha W} u$ and $p \coloneqq \frac{2(n+1)}{n+3}$, the lower Sobolev conjugate of $2$ in dimension $n+1$. We apply the Caccioppoli inequality to $u-c$ and bound the right-hand side by the Sobolev--Poincaré inequality in order to give:
\begin{align*}
\bigg(\bariint_{W} |\nabla u|^2 \, \d s \d y \bigg)^{1/2} \leq C \bigg(\bariint_{\alpha W} |\nabla u|^p \, \d s \d y \bigg)^{1/p}.
\end{align*}
As we have $p<2$, this is a reverse H\"older inequality for $\nabla u$. It remains to lower the exponent to $p=1$ but this is always possible by
a general feature of such inequalities, see \cite[Thm.~2]{Iwaniec-Nolder}. Strictly speaking, this reference is for $W=I \times B$ with $B$ a cube and the case of a ball then follows by a straightforward covering argument.

For the second inequality we let $c \coloneqq \bariint_W u$ and note that $\frac{2(n+1)}{n-1}$ is the upper Sobolev conjugate of $2$ in dimension $n+1$. It follows that
\begin{align*}
\bigg(\bariint_{W} |u|^q \, \d s \d y \bigg)^{1/q}
&\leq \bigg(\bariint_{W} |u-c|^q \, \d s \d y \bigg)^{1/q} + \bariint_{W} |u| \, \d s \d y \\
&\leq C \bigg(\bariint_{\alpha W} |u|^2 \, \d s \d y \bigg)^{1/2},
\end{align*}
where the second step follows again by combining the Sobolev--Poincaré inequality with the Caccioppoli inequality. The exponent on the right-hand side can be lowered as before.
\end{proof}

We close with a simple but important approximation result for weak solutions.

\begin{lem}
\label{lem: weak solution L1loc convergence}
Let $(u_k)$ be a sequence of weak solutions to $\Le u_k = 0$ in an open set $O \subseteq \ree$ that converges to $u$ in $\Lloc^1(O)$. Then $u$ is a weak solution to \eqref{eq: elliptic equation in weak sol section} in $O$ and $(u_k)$ tends to $u$ in $\Wloc^{1,2}(O)$.
\end{lem}

\begin{proof}
The Cauchy property in $\Wloc^{1,2}(O)$ follows by applying the reverse H\"older and the Caccioppoli inequality to $u_k - u_j$ on arbitrary admissible cylinders. Hence, we can pass to the limit in $k$ in the weak formulation of the equation for $u_k$.
\end{proof}

\begin{cor}
\label{cor: weak solution smooth in t}
If u is a weak solution to $\Le u = 0$ in $\reu$, then so is $\partial_t u$. In particular, $u$ is of class $\C^\infty((0,\infty); \Lloc^2)$.
\end{cor}

\begin{proof}
For $\eps > 0$ and $h \in (-\eps,\eps)$ define $v_h(t,x) \coloneqq \frac{1}{h}(u(t+h,x)-u(t,x))$ in $\R_{+,\eps}^{1+n} \coloneqq \{(s,y) \in \ree : s > \eps \}$. All $v_h$ are weak solutions in $\R_{+,\eps}^{1+n}$ since the coefficients of $\Le$ are independent of $t$ and we have $v_h \to \partial_t u$ in $\Lloc^2(\R_{+,\eps}^{1+n})$ as $h \to 0$. By the preceding lemma, $\partial_t u$ is a weak solution in $\reu$, so that in particular $\partial_t^2 u \in \Lloc^2((0,\infty); \Lloc^2)$. By iteration the same applies to $\partial_t^k u$ for any $k \in \IN$ and the claimed regularity follows by (one-dimensional) Sobolev embeddings.
\end{proof}
\section{Existence in \texorpdfstring{$\H^p$}{Hp} Dirichlet and Regularity problems}
\label{sec: estimates semigroup}

\noindent In this section we establish the existence part in our main results on the Dirichlet and Regularity problems with $\H^p$-data, Theorems~\ref{thm: blockdir} and \ref{thm: blockreg}. When the data $f$ additionally belongs to $\L^2$, the (eventually unique) solution is given by the Poisson semigroup. Hence, we proceed in two steps: First, we establish the required semigroup estimates for data $f \in a^{-1}(\H^p \cap \L^2)$ and $f \in \Hdot^{1,p} \cap \W^{1,2}$, respectively. Second, we obtain existence of a solution by a density argument for the full class of data.
\subsection{Estimates towards the Dirichlet problem}
\label{subsec: estimates Dir} 

We begin with the square function bound. 

\begin{prop}
\label{prop: SFE Dirichlet}
Let $p_{-}(L) < p < p_+(L)^*$. If $f \in \L^2$, then  $u(t,x) \coloneqq \e^{-t {L^{1/2}}}f(x)$ satisfies
\begin{align*}
\|S(t \nabla u)\|_p \simeq \|af\|_{\H^p}.
\end{align*}
\end{prop}

\begin{proof}
We organize the argument in three steps. For $p \leq 2$ we will be able to use Hardy space theory `off-the-shelf' but for $p>2$ different arguments on the level of the second-order equation for $u$ are needed since $p$ might lie outside of $\cH(L)$.

\medskip

\noindent \emph{Step~1: The case $p_-(L) < p \leq 2$}. We have
\begin{align*}
t \partial_t u
= -t L^{1/2} \e^{-t L^{1/2}} f \eqqcolon\psi(t^2 L)f,
\end{align*}
and, recalling \eqref{eq: BD and BD} - \eqref{eq: tL and tM},
\begin{align*}
t \nabla_x u 
= t \nabla_x a^{-1} \e^{-t \tL^{1/2}} (af) 
= (- tDB \e^{-t [DB]} g)_\ta \eqqcolon (\varphi(tDB)g)_\ta
\end{align*}
where $g = [af, 0]^\top$. We recall from Proposition~\ref{prop: auxiliary function bisectorial} and the corresponding result for sectorial operators in Section~\ref{subsec: Hardy abstract sectorial} that $\psi \in \Psi_{1/2}^\infty$ and $\varphi \in \Psi_1^\infty$ are admissible auxiliary functions for $\IH_L^p$ and $\IH_{DB}^p$, respectively. By Theorem~\ref{thm: main result Hardy} we have $p \in \cH(L) \cap \cH(DB)$ and hence we get as required
\begin{align*}
\|af\|_{\H^p} 
&\simeq \|f\|_{\IH_L^p} \\
&\simeq \|S(t \partial_t u)\|_p \\
&\leq \|S(t\nabla u)\|_p \\
&\lesssim \|f\|_{\IH_L^p} + \|g\|_{\IH_{DB}^p} \\
&\simeq \|af\|_{\H^p} + \|g\|_{\H^p} \\
&\simeq \|af\|_{\H^p}.
\end{align*}

\medskip

\noindent \emph{Step~2: Upper bound for $2 < p < p_+(L)^*$}. Consider the auxiliary function $\phi(z) \coloneqq \e^{-\sqrt{z}} - (1+z)^{-2}$. Then $\phi \in \Psi_{1/2}^2$ on any sector. Differentiating the resolvent twice, we find that $v \coloneqq \phi(t^2L)f$ solves the following equation in $\reu$ in the weak sense:
\begin{align*}
&(a \partial_t^2 + \div_x d \nabla_x) v \\
&\quad = 4a L (1+t^2L)^{-3}f - 24 at^2 L^2 (1+t^2L)^{-4}f - aL (1+t^2L)^{-2}f \\
&\quad \eqqcolon t^{-2} a \psi(t^2 L)f,
\end{align*}
where $\psi\in \Psi_1^1$ on any sector. For $x \in \R^n$ and $t>0$ consider Whitney boxes $W(t,x) \coloneqq (t,2t) \times B(x,2t)$ and $\wt{W}(t,x) \coloneqq (\nicefrac{t}{2},4t) \times B(x,4t)$. The Caccioppoli inequality yields
\begin{align}
\label{eq0: SFE Dirichlet}
\bariint_{W(t,x)} |s \nabla v|^2 \, \d s  \d y \lesssim \bariint_{\wt{W}(t,x)} |v|^2 + |\psi(s^2 L)f|^2 \, \d s \d y.
\end{align}
Summing up these estimates for $t = 2^{-k}$, $k \in \IZ$, leads to\index{inequality! Caccioppoli on a cone}
\begin{align*}
\iint_{|x-y|<s} |s \nabla v|^2 \, \frac{\d s  \d y}{s^{1+n}} \lesssim \iint_{|x-y|<8s} |v|^2 + |\psi(s^2 L)f|^2 \, \frac{\d s  \d y}{s^{1+n}},
\end{align*}
where we have used that at most $3$ of the enlarged boxes $\wt{W}(2^{-k},x)$ overlap in order to get the term on the right. By definition of $v$ we conclude
\begin{align}
\label{eq1: SFE Dirichlet}
\|S(t \nabla u)\|_p
\lesssim \|S_{\phi,L}f\|_p + \|S_{\psi,L}f\|_p + \|S(t\nabla(1+t^2L)^{-2}f)\|_p,
\end{align}
where as usual $S_{\phi,L}f$ denotes the square function of $\phi(t^2L)f(x)$. 

Since $\phi \in \Psi_{1/2}^2$ and $\psi \in \Psi_1^1$, Theorem~\ref{thm: SFE from identification} applies in our range of exponents and yields
\begin{align*}
\|S_{\phi,L}f\|_p + \|S_{\psi,L}f\|_p \lesssim \|f\|_p.
\end{align*}

The analogous bound for the third square function in \eqref{eq1: SFE Dirichlet} is a consequence of Remark~\ref{rem: Carleson bound Q function}.  Indeed, the family
\begin{align*}
t\nabla(1+t^2L)^{-2} 
&= \begin{bmatrix}
-2 t^2 L(1+t^2 L)^{-3} \\ t \nabla_x (1+t^2L)^{-2}
\end{bmatrix} \\
&= \begin{bmatrix}
2 ((1+t^2 L)^{-3} - (1+t^2L)^{-2}) \\ t \nabla_x (1+t^2L)^{-2}
\end{bmatrix}
\end{align*}
satisfies $\L^2$ off diagonal estimates of arbitrary large order by composition and we have for all $t>0$
\begin{align*}
\|t\nabla(1+t^2L)^{-2}f \|_2^2 \simeq \|t^2 L(1+t^2 L)^{-3}f\|_2^2 + \|t {L^{1/2}} (1+t^2L)^{-2}f\|_2^2
\end{align*}
by the solution of the Kato problem, so that the theorems of Fubini and McIntosh yield the $\L^2$-bound
\begin{align*}
\|S(t\nabla(1+t^2L)^{-2}f)\|_2^2
\simeq \int_0^\infty \|t\nabla(1+t^2L)^{-2}f\|_2^2 \, \frac{\d t}{t} 
\simeq \|f\|_2^2.
\end{align*}
\medskip

\noindent \emph{Step~3: Lower bound for $2 < p < p_+(L)^*$}. Introduce the adapted Laplacian $H \coloneqq -(a^*)^{-1} \Delta_x$ and for $f \in \L^2$ with $S(t\nabla u) \in \L^p$ and $g \in \L^{p'} \cap \L^2$ set
\begin{align*}
\Phi: (0,\infty) \to \IC, \quad \Phi(t) \coloneqq \langle a \e^{-t {L^{1/2}}}f, \e^{-t {H^{1/2}}} g \rangle,
\end{align*}
where $\langle \cdot \,, \cdot \rangle$ is the $\L^2$ inner product. By the functional calculus on $\L^2$, this is a smooth function and we have
\begin{align*}
\Phi'(t) 
&= - \langle a {L^{1/2}} \e^{-t {L^{1/2}}}f, \e^{-t {H^{1/2}}} g \rangle - \langle a \e^{-t {L^{1/2}}}f, {H^{1/2}} \e^{-t {H^{1/2}}} g \rangle \\
\Phi''(t)
&= \langle a L \e^{-t {L^{1/2}}}f, \e^{-t {H^{1/2}}} g \rangle + 2 \langle a {L^{1/2}} \e^{-t {L^{1/2}}}f, {H^{1/2}} \e^{-t {H^{1/2}}} g \rangle  \\ &\quad +\langle a \e^{-t {L^{1/2}}}f, H \e^{-t {H^{1/2}}} g \rangle \\
&= \langle d \nabla_x \e^{-t {L^{1/2}}}f, \nabla_x \e^{-t {H^{1/2}}} g \rangle + 2 \langle a {L^{1/2}} \e^{-t {L^{1/2}}}f, {H^{1/2}} \e^{-t {H^{1/2}}} g \rangle  \\ &\quad + \langle \nabla_x \e^{-t {L^{1/2}}}f, \nabla_x \e^{-t {H^{1/2}}} g \rangle,
\end{align*}
as well as
\begin{align*}
\lim_{t \to \infty} \Phi(t) = \lim_{t \to \infty} t \Phi'(t) = \lim_{t \to 0} t \Phi'(t)= 0, \quad \lim_{t \to 0} \Phi(t) = \langle af, g \rangle.
\end{align*}
Putting all together and integrating by parts twice in $t$, we obtain
\begin{align*}
\langle af, g \rangle
&= \int_0^\infty t^2 \Phi''(t) \, \frac{\d t}{t} \\
&= \int_0^\infty \langle d t \nabla_x \e^{-t {L^{1/2}}}f, t \nabla_x \e^{-t {H^{1/2}}} g \rangle \, \frac{\d t}{t} \\
&\quad + 2 \int_0^\infty \langle a t {L^{1/2}} \e^{-t {L^{1/2}}}f, t {H^{1/2}} \e^{-t {H^{1/2}}} g \rangle \, \frac{\d t}{t} \\ 
&\quad + \int_0^\infty \langle t \nabla_x \e^{-t {L^{1/2}}}f, t \nabla_x \e^{-t {H^{1/2}}} g \rangle \, \frac{\d t}{t}.
\end{align*}
We regard the right-hand side as $\T^p - \T^{p'}$ duality pairings in order to give
\begin{align*}
|\langle af, g \rangle|
&\lesssim
\|S(t \nabla_x \e^{-t {L^{1/2}}}f)\|_p \|S(t \nabla_x \e^{-t {H^{1/2}}}g) \|_{p'} \\
&\quad +\|S(t {L^{1/2}} \e^{-t {L^{1/2}}}f)\|_p \|S(t {H^{1/2}} \e^{-t {H^{1/2}}}g) \|_{p'} \\
&\leq 2 \|S(t \nabla \e^{-t {L^{1/2}}}f)\|_p \|S(t \nabla \e^{-t {H^{1/2}}}g) \|_{p'}.
\end{align*}
We know that $p_-(H) = 1_*$ from Corollary~\ref{cor: McIntosh-Nahmod}. Hence, Step~1 for $H$ on $\L^{p'}$ yields $\|S(t \nabla \e^{-t {H^{1/2}}}g) \|_{p'} \lesssim \|g\|_{p'}$ and since $g \in \L^{p'} \cap \L^2$ was arbitrary, we conclude
\begin{align*}
 \|af\|_p \lesssim \|S(t \nabla \e^{-t {L^{1/2}}}f)\|_p. & \qedhere
\end{align*}
\end{proof} 

We turn to bounds for the non-tangential maximal function and begin by recalling the respective $\L^2$-bound for our perturbed Dirac operators. 
\begin{thm}[{\cite[Thm.~9.9]{AusSta}}]
	\label{thm: NT and Fatou in L2}
	Let $T$ be one of $DB$ or $BD$. Then
	\begin{align*}
		\|\NT(\e^{-t[T]}f)\|_2 \simeq \|f\|_2 \quad (f \in \cl{\ran(T)})
	\end{align*}
	and for every $f \in \L^2$ the Whitney averages converge in the $\L^2$-sense
	\begin{align*}
		\lim_{t \to 0} \bariint_{W(t,x)} |\e^{-t [T]}f - f(x)|^2 \, \d s \d y= 0 \quad (\text{a.e. } x \in \R^n).
	\end{align*}
\end{thm}

We remark that the result above for $T=BD$ is originally due to Rosén~\cite[Thm.~5.1]{R}.

\begin{prop}
\label{prop: NT for Dirichlet}
Let $p_{-}(L) < p < p_+(L)^*$. If $f \in \L^2$, then  $u(t,x) \coloneqq \e^{-t {L^{1/2}}}f(x)$ satisfies
\begin{align*}
\|\NT(u)\|_p \simeq \|af\|_{\H^p}
\end{align*}
and
\begin{align*}
\lim_{t \to 0} \bariint_{W(t,x)} |u(s,y) - f(x)|^2 \, \d s \d y = 0 \quad (\text{a.e. } x \in \R^n).
\end{align*}
\end{prop}

\begin{proof}
We recall from \eqref{eq: L and M} that $L$ is incorporated in the matrix operator $(BD)^2$. Hence, we have
\begin{align*}
\e^{-t [BD]} \begin{bmatrix} f \\ 0 \end{bmatrix} 
= \begin{bmatrix} u \\ 0 \end{bmatrix} 
\end{align*}	
and the claim for $p=2$ as well as the convergence of averages follows from Theorem~\ref{thm: NT and Fatou in L2}.

\medskip

\noindent \emph{Step 1: Upper bound}. If $p \in (p_{-}(L),2)$, then according to Proposition~\ref{prop: NT bound semigroup for p<1} and Theorem~\ref{thm: main result Hardy} we have
\begin{align*}
\|\NT(u)\|_p \lesssim \|f\|_{\IH_L^p} \simeq \|af\|_{\H^p}.
\end{align*}
If $p \in (2, p_+(L)^*)$, we first introduce $\psi(z) \coloneqq \e^{-\sqrt{z}} - (1+z)^{-1}$ and split 
\begin{align*}
u = v + w \coloneqq \psi(t^2L)f + (1+ t^2 L)^{-1}f.
\end{align*}
We have $\psi \in \Psi_{1/2}^1$ on any sector. Combining Lemma~\ref{lem: N < S} and Theorem~\ref{thm: SFE from identification}, we find that
\begin{align*}
\|\NT(v)\|_p  \lesssim \|S_{\psi,L}f \|_p \lesssim \|f\|_{p}.
\end{align*}
As for $w$, we use that the resolvents satisfy off-diagonal estimates of arbitrarily large order. Consequently, Lemma~\ref{lem: Whitney everage controlled by HL-Max} and the $\L^{p/2}$-bound for the Hardy--Littlewood maximal operator yield
\begin{align*}
\|\NT(w)\|_p \leq \|(\Max(|f|^2))^{1/2}\|_p \lesssim \|f\|_p.
\end{align*}

\medskip

\noindent \emph{Step 2: Lower bound for $p >1$}. The convergence of Whitney averages implies $\NT(u) \geq f$ a.e.\ on $\R^n$ and $\|\NT(u)\|_p \geq \|f\|_p$ follows.

\medskip

\noindent \emph{Step 3: Lower bound for $p_-(L) < p \leq 1$}.  We calculate the $\H^p$-norm of $af$ using the Fefferman--Stein characterization of $\H^p$.\index{Fefferman--Stein characterization (of $\H^p$)} This argument works for all $p \in (\nicefrac{n}{(n+1)},1]$, not only $p \in (p_{-}(L),1]$.

Fix $\phi \in \C_0^\infty(\R^n; \R)$ with support in $B(0,1)$ and $\int_{\R^n} \phi \d y = 1$ and let $\phi_t(y) \coloneqq t^{-n} \phi(\nicefrac{y}{t})$. Then a function $h \in \L^2$ belongs to $\H^p$ if and only if the maximal function
\begin{align*}
(\Max_\phi h)(x) \coloneqq \sup_{t>0} |h \ast \phi_t|(x) \quad (x \in \R^n)
\end{align*}
is in $\L^p$ and in this case $\|h\|_{\H^p} \simeq \|\Max_\phi h\|_p$, see e.g.\ \cite[Thm.~6.4.4]{Grafakos}.

Temporarily fix $t>0$ and $x \in \R^n$. Let $\chi : [0,\infty) \to [0,1]$ be smooth with $\ind_{[0,\nicefrac{1}{2}]} \leq \chi \leq \ind_{[0,2)}$, set $\chi_t(s) \coloneqq \chi(\nicefrac{s}{t})$ and introduce $\Phi(s,y) \coloneqq \phi_t(x-y) \chi_t(s)$. The functional calculus on $\L^2$ and the compact support of $\Phi$ justify writing
\begin{align*}
(af \ast \phi_t)(x)
&=\int_{\R^n} (af)(y) \phi_t(x-y) \, \d y \\
&= \lim_{\eps \to 0} \int_{\R^n} -\bigg(\int_\eps^\infty \partial_{s}(\Phi au) \, \d s \bigg) \d y.
\end{align*}
For $\eps < \nicefrac{t}{2}$ we expand, integrate by parts and use $a \partial_s^2 u = L u$, to give
\begin{align*}
\int_{\R^n} &\int_\eps^\infty \partial_{s}(\Phi au) \, \d s \d y \\
& = \int_{\R^n} \int_\eps^\infty (\partial_{s} \Phi) au + \Phi a \partial_s u\, \d s \d y \\
& = \int_{\R^n} \int_\eps^\infty (\partial_{s} \Phi) au - (\partial_s \Phi) s a \partial_s u - \Phi s a \partial_s^2 u \, \d s \d y \\
&\quad + \int_{\R^n} \Phi(\eps,y) \eps a(y) \partial_s u(\eps,y) \, \d y \\
& = \int_{\R^n} \int_\eps^\infty (\partial_{s} \Phi) au - (\partial_s \Phi) a s \partial_s u + \nabla_y \Phi \cdot s d \nabla_y u \, \d s \d y \\
&\quad + \int_{\R^n} \phi_t(x-y) \eps a(y) \partial_s u(\eps,y) \, \d y.
\end{align*}
By the functional calculus for $L$ we have as a limit in $\L^2$,
\begin{align*}
\lim_{\eps \to 0} \eps a \partial_s u(\eps,\cdot ) = - \lim_{\eps \to 0} \eps a L^{1/2} \e^{-\eps L^{1/2}} f = 0.
\end{align*}
By Young's convolution inequality we get $\phi_t \ast (\eps a \partial_s u(\eps,\cdot )) \to 0$ uniformly on $\R^n$ as $\eps \to0$. Altogether,
\begin{align}
\label{eq1: NT for Dirichlet}
\begin{split}
|(af &\ast \phi_t)(x)| \\
&\leq \iint_{\reu} |(\partial_{s} \Phi) au| \, \d s \d y + \iint_{\reu} |(\partial_s \Phi) a s \partial_s u| \, \d s \d y \\
&\quad + \iint_{\reu} |\nabla_y \Phi \cdot s d \nabla_y u| \, \d s \d y \\
&\eqqcolon \I + \II + \III.
\end{split}
\end{align}
Since $\partial_s \Phi$ is bounded by $t^{-1-n}$ and supported in $W(t,x)$, we get
\begin{align*}
|\I| + |\II| \lesssim \NT(u)(x) +  \NT(t \partial_t u)(x).
\end{align*}
As for $\III$, we get
\begin{align*}
|\III| \lesssim t^{-1-n} \iint_{\reu} |\ind_{(0,2t) \times B(x,t)}s \nabla_x u| \, \d s \d y,
\end{align*}
so that Lemma~\ref{lem:AM} applied to $F \coloneqq |\ind_{(0,2t) \times B(x,t)}s \nabla_x u|$ with $r=1$ and $p = \nicefrac{n}{(n+1)}$ yields
\begin{align*}
|\III| \lesssim t^{-1-n} \|\NT( F)\|_{\frac{n}{n+1}}.
\end{align*}
If a Whitney ball $W(r,z)$ intersects the support of $F$ at some $(s,y) \in \reu$, then 
\begin{align*}
|x-z| \leq |x-y| + |y-z| \leq t + r \leq t + 2s \leq 5t,
\end{align*}
which means that $\NT (F)$ has support in $B(x,5t)$. Thus, we have
\begin{align*}
|\III| 
\lesssim \bigg( t^{-n} \int_{B(x,5t)} |\NT (F)|^{\frac{n}{n+1}} \bigg)^{\frac{n+1}{n}} \\
\lesssim \Max( |\NT(t \nabla_x u)|^{\frac{n}{n+1}})^{\frac{n+1}{n}}(x). 
\end{align*}
Going back to \eqref{eq1: NT for Dirichlet} and taking the supremum in $t$, leads us to
\begin{align*}
\Max_\phi (a f) \lesssim \NT(u) + \NT(t \partial_t u) + \Max( |\NT(t \nabla_x u)|^{\frac{n}{n+1}})^{\frac{n+1}{n}}.
\end{align*}
By assumption we have $p> \nicefrac{n}{(n+1)}$. Hence, $\Max$ is bounded on $\L^{p(n+1)/n}$ and it follows that
\begin{align*}
\|a f\|_{\H^p} 
\simeq \|\Max_\phi(af)\|_p
\lesssim \|\NT(u)\|_p +  \|\NT(t \nabla u)\|_p 
\lesssim \|\NT(u)\|_p,
\end{align*}
where the final step is due to Caccioppoli's inequality.
\end{proof}

Finally, we establish uniform bounds and strong continuity at $t=0$.

\begin{prop}
\label{prop: sup estimates Dirichlet}
Let $p_{-}(L) < p < p_+(L)$. If $f \in a^{-1} (\H^p \cap \L^2)$ and  $u(t,x) \coloneqq \e^{-t {L^{1/2}}}f(x)$, then $au$ is of class 
\begin{align*}
\C_0([0,\infty);  \H^p) \cap \C^\infty((0,\infty); \H^p)
\end{align*}
and satisfies
\begin{align*}
\sup_{t>0} \|a u(t,\cdot) \|_{\H^p} \simeq \|af\|_{\H^p}
\end{align*}
and for all $k \in \IN$,
\begin{align*}
\sup_{t>0}  \|t^{\frac{k}{2}} \partial_t^k (a u(t,\cdot))\|_{\H^p} \lesssim (\tfrac{k}{2})^{\frac{k}{2}} \e^{-k} \|af\|_{\H^p}.
\end{align*}
\end{prop}

\begin{proof}
According to Theorem~\ref{thm: main result Hardy} we have $a^{-1} (\H^p \cap \L^2) = \IH_L^p$ with equivalent $p$-quasinorms $\|f\|_{\IH_L^p} \simeq \|af\|_{\H^p}$.

The upper bounds for $u$ and $\partial_t^k u$ now follow immediately from the bounded $\H^\infty$-calculus on $\IH_L^p$, see Section~\ref{subsec: Hardy abstract sectorial}. Likewise, Proposition~\ref{prop: C0 smg on abstract Hardy sectorial} provides the limits $a u(t,\cdot) \to af$ as $t \to 0$ and $a u(t,\cdot) \to 0$ as $t \to \infty$ in $\H^p$ and the limit at $t = 0$ implies lower bound for $u$.
\end{proof}

For exponents $p \ge p_+(L)$, the space $\IH_L^p$ does not equal $a^{-1} (\H^p \cap \L^2)$ and the previous argument breaks down.   However, using off-diagonal estimates, we can still obtain the continuity at the boundary $t=0$ with values in $\Lloc^2$ if $p_+(L) \leq p < p_+(L)^*$.

\begin{lem}
\label{lem: sup L2loc estimates Dirichlet}
If $p_+(L) \leq p < p_+(L)^*$, then for all $f \in \L^p \cap \L^2$, all balls $B \subseteq \R^n$ and all $t>0$,
\begin{align*}
\|\e^{-t {L^{1/2}}}f-f\|_{\L^2(B)} 
\lesssim r(B)^{\frac{n}{2}-\frac{n}{p}-1}(r(B)+  t) \|f\|_p.
\end{align*}
\end{lem}

\begin{proof}
We can pick $q$ such that $2 \leq q < p_+(L)$ and $\nicefrac{1}{q}-\nicefrac{1}{p} < \nicefrac{1}{n}$. We split $f = \sum_{j \geq 1} f_j$, where $f_j \coloneqq \ind_{C_j(B)} f$, and obtain from H\"older's inequality that
\begin{align*}
\|\e^{-t {L^{1/2}}}f-f\|_{\L^2(B)}
&\leq \|\e^{-t {L^{1/2}}}f_1- f_1\|_{\L^2(B)} \\
&\quad+  r(B)^{\frac{n}{2}-\frac{n}{q}}  \sum_{j \geq 2} \|\e^{-t {L^{1/2}}} f_j\|_{\L^q(B)} \\
&\leq r(B)^{\frac{n}{2}-\frac{n}{p}} \|f\|_p \\
&\quad+  r(B)^{\frac{n}{2}-\frac{n}{q}}  \sum_{j \geq 2} \|\e^{-t {L^{1/2}}} f_j\|_{\L^q(B)}.
\end{align*}
Since the Poisson semigroup satisfies $\L^q$ off-diagonal estimates of order $1$, see Corollary~\ref{cor: Poisson semigroup Lp OD}, we can bound the sum in $j$ by
\begin{align*}
\sum_{j \geq 2} t 2^{-j} r(B)^{-1}  \|f_j\|_{\L^q(B)}  \lesssim  t r(B)^{\frac{n}{q}-\frac{n}{p}-1}   \sum_{j \geq 2} 2^{j(\frac{n}{q}-\frac{n}{p}-1)} \|f\|_p,
\end{align*}
where the right-hand side is finite by choice of $q$. The claim follows.
\end{proof}
\subsection{Estimates towards the Regularity problem}
\label{subsec: estimates Reg}

We begin again with the square function bounds.

\begin{prop}
\label{prop: SFE Regularity}
Let $(p_{-}(L)_* \vee 1_*) < p < q_+(L)$. If $f \in \W^{1,2}$, then  $u(t,x) \coloneqq \e^{-t {L^{1/2}}}f(x)$ satisfies
\begin{align*}
\|S(t \nabla \partial_t u )\|_p \simeq \|\nabla_x f\|_{\H^p}.
\end{align*}
\end{prop}

\begin{proof}
Let us first interpret the exponents. The identification Theorem~\ref{thm: main result Hardy} tells us that we have $\IH^{1,p}_L = \Hdot^{1,p} \cap \L^2$ with equivalent $p$-quasinorms and then $\|g\|_{\IH_{\tM}^p} \simeq \|g\|_{\H^p}$ for all $g \in \cl{\ran(\nabla_x)}$ follows from Figure~\ref{fig: diagram}. The square function we have to control contains
\begin{align*}
t \nabla_x \partial_t u
&= -t \nabla_x {L^{1/2}} \e^{-t {L^{1/2}}}f 
= -t \tM^{1/2} \e^{-t \tM^{1/2}} \nabla_x f 
\eqqcolon \psi(t^2 \tM) \nabla_x f,
\intertext{where $\psi \in \Psi_{1/2}^\infty$ on any sector and we used an intertwining relation for the functional calculus on $\L^2$, as well as}
t \partial_t^2 u 
&= - \psi(t^2 L) {L^{1/2}} f
=t^{-1} (t^2 L \e^{-t {L^{1/2}}}) f 
\eqqcolon t^{-1} \phi(t^2L)f,
\end{align*}
where $\phi \in \Psi_1^\infty$ on any sector.

If $p \leq 2$, then $\phi$ and $\psi$ are admissible auxiliary functions for defining $\IH_L^{1,p}$ and $\IH_{\tM}^p$, respectively. Thus, we get 
\begin{align*}
\|S(t \nabla \partial_t u)\|_p
&\simeq \|S(t^{-1} \phi(t^2L)f)\|_p + \|S(\psi(t^2 \tM) \nabla_x f)\|_p \\
&\simeq \|f\|_{\IH_L^{1,p}} + \|\nabla_x f\|_{\IH_{\tM}^p} \\
&\simeq \|\nabla_x f\|_{\H^p}
\end{align*}
right away. 

If $p \geq 2$, then Proposition~\ref{prop: upper SFE abstract} applies to $\tM$ with auxiliary function $\psi$ and $q=p$. The same holds for $L$ since from Theorem~\ref{thm: main result Hardy} and the general bound $q_+(L) < p_+(L)$ in Theorem~\ref{thm: standard relation J(L) and N(L)} we obtain $\IH_L^p = \L^p \cap \L^2$ with equivalent $p$-norms. Consequently, we get the upper bound
\begin{align*}
\|S(t \nabla \partial_t u)\|_p
&\simeq \|S(\psi(t^2 L) {L^{1/2}} f)\|_p + \|S(\psi(t^2 \tM) \nabla_x f)\|_p \\
&\lesssim \|{L^{1/2}}f\|_p + \|\nabla_x f\|_p \\
&\simeq \|\nabla_x f\|_p,
\end{align*}
where the final estimate is due to Theorem~\ref{thm: Riesz complete}. Instead of Proposition~\ref{prop: upper SFE abstract}, we can also use the square function estimate for the Dirichlet problem in Proposition~\ref{prop: SFE Dirichlet} directly. This argument also gives the lower bound: It suffices to note that $\partial_t  u = v$, where $v \coloneqq \e^{-t {L^{1/2}}} (-{L^{1/2}}f)$ with $(-{L^{1/2}}f) \in \L^2$, so that 
\begin{align*}
\|S(t \nabla \partial_t u)\|_p
= \|S(t \nabla v)\|_p
\simeq \|{L^{1/2}}f\|_p 
\simeq \|\nabla_x f\|_p
\end{align*}
follows as required.
\end{proof}

We continue with the non-tangential maximal function bounds. 

\begin{prop}
\label{prop: NT for regularity}
Let $(p_{-}(L)_* \vee 1_*) < p < q_+(L)$. If $f \in \W^{1,2}$, then $u(t,x) \coloneqq \e^{-t {L^{1/2}}}f(x)$ satisfies
\begin{align*}
\|\NT(\nabla u)\|_p \simeq \|\nabla_x f\|_{\H^p}
\end{align*}
and
\begin{align*}
\lim_{t \to 0} \bariint_{W(t,x)} \bigg|\begin{bmatrix} a \partial_t u \\ \nabla_x u \end{bmatrix} - \begin{bmatrix} -a{L^{1/2}}f(x) \\ \nabla_x f(x) \end{bmatrix}\bigg|^2 \, \d s \d y = 0 \quad (\text{a.e. } x \in \R^n).
\end{align*}
\end{prop}

\begin{proof}
We use the intertwining property to write $\nabla_x u = \e^{- t \tM^{1/2}} \nabla_x f$. Moreover, we have $\partial_t u = \e^{-t {L^{1/2}}} (-{L^{1/2}}f)$, so that by similarity $a \partial_t u = \e^{-t {\tL^{1/2}}} (-a{L^{1/2}}f)$. We recall from \eqref{eq: tL and tM} that $\tM$ and $\tL$ are incorporated in the matrix operator $(DB)^2$. Hence, we have
\begin{align*}
\e^{-t [DB]} \begin{bmatrix} -a{L^{1/2}}f \\ \nabla_x f \end{bmatrix} 
= \begin{bmatrix} a\partial_t u  \\ \nabla_x u \end{bmatrix}.
\end{align*}
The claim for $p=2$ as well as the convergence of averages now follows from Theorem~\ref{thm: NT and Fatou in L2} and the comparison $\|a{L^{1/2}} f\|_2 \simeq \|\nabla_x f\|_2$.

\medskip

\noindent \emph{Step 1: Upper bound for $p \neq 2$}. As in the proof of Proposition~\ref{prop: SFE Regularity} we have $\IH^{1,p}_L = \Hdot^{1,p} \cap \L^2$ with equivalent $p$-quasinorms.

If $p \in (p_{-}(L)_* \vee 1_*,2]$, then Proposition~\ref{prop: NT bound semigroup for p<1} applied to $\tM$ and $L$ directly yields
\begin{align*}
\|\NT(\nabla u)\|_p 
&\leq \|\NT(\nabla_x u)\|_p + \|\NT(\partial_t u)\|_p \\
&\lesssim \|\nabla_x f\|_{\IH_{\tM}^p} +  \|{L^{1/2}}f\|_{\IH_L^p}
\intertext{and the ubiquitous Figure~\ref{fig: diagram} allows us to compare with}
&\simeq \|f\|_{\IH_L^{1,p}}\\
&\simeq \|\nabla_x f\|_{\H^p}
\end{align*}
as required. If $p \in (2, q_+(L))$, we first introduce $\psi(z) \coloneqq \e^{-\sqrt{z}} - (1+z)^{-1}$ and split 
\begin{align*}
\nabla u 
= v+w
\coloneqq \begin{bmatrix} -\psi(t^2 L) {L^{1/2}} f \\ \psi(t^2 \tM) \nabla_x f \end{bmatrix} + \begin{bmatrix} -(1+t^2 L)^{-1} {L^{1/2}} f \\ (1+ t^2\tM)^{-1} \nabla_x f  \end{bmatrix}.
\end{align*}
We have $\psi \in \Psi_{1/2}^1$ on any sector. As in the preceding proof,  Proposition~\ref{prop: upper SFE abstract} with $q=p$ and auxiliary function $\psi$ applies to both $\tM$ and $L$ in our range of exponents. Along with Lemma~\ref{lem: N < S}, we find that
\begin{align*}
\|\NT(v)\|_p 
&\leq \|\NT(\IQ_{\psi,\tM} \nabla_x f)\|_p + \|\NT(\IQ_{\psi,L} {L^{1/2}}f)\|_p\\
&\lesssim \|S_{\psi, \tM} (\nabla_x f) \|_p + \|S_{\psi, L} (L^{1/2} f) \|_p \\
&\lesssim \|\nabla_x f\|_p + \|{L^{1/2}} f\|_p\\
&\lesssim \|\nabla_x f\|_p,
\end{align*}
where the final estimate is due to Theorem~\ref{thm: Riesz complete}. As for $w$, we use that the resolvents of $L$ and $\tM$ satisfy off-diagonal estimates of arbitrarily large order. Consequently, Lemma~\ref{lem: Whitney everage controlled by HL-Max} and the $\L^{p/2}$-bound for the Hardy--Littlewood maximal operator yield
\begin{align*}
\|\NT(w)\|_p 
&\leq \|(\Max(|\nabla_xf|^2)^{1/2}\|_p + \|(\Max(|{L^{1/2}}f|^2)^{1/2}\|_p \\
&\lesssim \|\nabla_x f\|_p + \|{L^{1/2}}f\|_p 
\end{align*}
and we conclude as before. Combining these estimates gives the required upper bound for $\NT(\nabla u)$.

\medskip

\noindent \emph{Step 2: Lower bound for $p > 1$}. Since $\nabla_x f \in \L^2$, we obtain from the convergence of Whitney averages that $\NT(\nabla u) \geq |\nabla_x f|$ a.e.\ on $\R^n$ and $\|\NT(\nabla u)\|_p  \geq \|\nabla_x f\|_p$ follows.

\medskip

\noindent \emph{Step 3: Lower bound for $p\leq 1$}. As in Step~3 of the proof of Proposition~\ref{prop: NT for Dirichlet} we calculate the $\H^p$-norm of $\nabla_x f$ through the Fefferman--Stein characterization of $\H^p$. The argument works again for all $p \in (\frac{n}{n+1},1]$.

Fix $\phi \in \C_0^\infty(\R^n; \R)$ with support in $B(0,1)$ and $\int_{\R^n} \phi \d y = 1$ and let $\phi_t(y) \coloneqq t^{-n} \phi(\nicefrac{y}{t})$. We need to control the $\L^p$-norm of
\begin{align*}
\Max_\phi (\nabla_x f)(x) \coloneqq \sup_{t>0} |\nabla_x f \ast \phi_t|(x) \quad (x \in \R^n).
\end{align*}
Temporarily fix $t>0$ and $x \in \R^n$. Let $\chi : [0,\infty) \to [0,1]$ be smooth with $\ind_{[0,\nicefrac{1}{2}]} \leq \chi \leq \ind_{[0,2)}$, set $\chi_t(s) \coloneqq \chi(\nicefrac{s}{t})$ and introduce $\Phi(s,y) \coloneqq \phi_t(x-y) \chi_t(s)$. As $\nabla_x u(s,y) = \e^{- s \tM^{1/2}} \nabla_x f(y)$, the functional calculus on $\L^2$ and the compact support of $\Phi$ justify writing
\begin{align*}
(\nabla_x f \ast \phi_t)(x)
&=\int_{\R^n} \nabla_xf(y) \phi_t(x-y) \, \d y \\
&= \lim_{\eps \to 0} \int_{\R^n} -\bigg(\int_\eps^\infty \partial_{s}(\Phi \nabla_x u) \, \d s \bigg) \d y \\
&= \lim_{\eps \to 0} \int_\eps^\infty \int_{\R^n} - \partial_{s} \Phi \nabla_x u - \Phi \partial_s \nabla_x u \, \d y \d s,
\end{align*} 
so that
\begin{align*}
|(\nabla_x f \ast \phi_t)(x)|
& \leq\iint_{\reu} |\partial_{s} \Phi \nabla_x u| \, \d s \d y + \iint_{\reu} |\nabla_x \Phi \otimes \partial_s u| \, \d s \d y \\
&\eqqcolon \I + \II,
\end{align*}
where $\nabla_x \Phi \otimes \partial_s u$ is the vector in $(\IC^m)^n$ coming from integration by parts in $x$. Now, we can literally repeat the arguments in Step~3 of the proof of Proposition~\ref{prop: NT for Dirichlet} and arrive at
\begin{align*}
\I \lesssim \NT(\nabla_x u)(x)
\end{align*}
and
\begin{align*}
\II 
\lesssim \Max( |\NT(\partial_t u)|^{\frac{n}{n+1}})^{\frac{n+1}{n}}(x)
\end{align*}
for all $x \in \R^n$. Consequently, we have a pointwise bound
\begin{align*}
\Max_\phi (\nabla_x f) \lesssim \NT(\nabla_x u) + \Max( |\NT(\partial_t u)|^{\frac{n}{n+1}})^{\frac{n+1}{n}}
\end{align*}
and since $\Max$ is bounded on $\L^{p(n+1)/n}$ we get $\|\nabla_x f\|_{\H^p} \lesssim \|\NT(\nabla u)\|_p$ as required.
\end{proof}

Uniform boundedness and strong continuity follow again by abstract semigroup theory.

\begin{prop}
\label{prop: sup estimates regularity}
Let $(p_{-}(L)_* \vee 1_*) < p < q_+(L)$. If $f \in \Hdot^{1,p} \cap \W^{1,2}$, then  $u(t,x) \coloneqq \e^{-t {L^{1/2}}}f(x)$ satisfies the following:
\begin{enumerate}
\item $\nabla_x u \in \C_0([0,\infty);  \H^p) \cap \C^\infty((0,\infty); \H^p)$ with
\begin{align*}
\sup_{t>0} \|\nabla_x u(t,\cdot)\|_{\H^p} \simeq \|\nabla_x f\|_{\H^p}
\end{align*}
and, for every $k \in \IN$,
\begin{align*}
\sup_{t>0}  \|t^{\frac{k}{2}} \partial_t^k \nabla_x u(t,\cdot)\|_{\H^p} \lesssim (\tfrac{k}{2})^{\frac{k}{2}} \e^{-k} \|\nabla_x f\|_{\H^p}.
\end{align*}
\item If $p<n$, then $u \in \C_0([0,\infty);  \L^{p^*}) \cap \C^\infty((0,\infty); \L^{p^*})$ with
\begin{align*}
\|f\|_{p^*} \leq \sup_{t > 0} \|u(t,\cdot)\|_{p^*} \lesssim \|\nabla_x f\|_{\H^p} + \|f\|_{p^*}.
\end{align*}
\end{enumerate}
\end{prop}

\begin{proof}
From the proofs of Propositions~\ref{prop: NT for regularity} and \ref{prop: SFE Regularity} we know $\nabla_x u = \e^{-t \tM^{1/2}}\nabla_x f$ and that in the given range of exponents $\|g\|_{\IH_{\tM}^p} \simeq \|g\|_{\H^p}$ holds for all $g \in \cl{\ran(\nabla_x)}$. Hence, (i) follows \emph{verbatim} as for the Dirichlet problem in Proposition~\ref{prop: sup estimates Dirichlet} by appealing to the abstract theory for $\tM$ instead of $L$.

For $p<n$ we have the Sobolev embedding $\Hdot^{1,p} \subseteq \L^{p^*}/\IC^m$ but since $\L^{p^*}+ \L^2$ does not contain any constants but $0$ we also have $\Hdot^{1,p} \cap \L^2 \subseteq \L^{p^*}$. This yields the regularity statement in (ii) and the upper bound, whereas the lower bound follows again from the continuity at $t=0$.
\end{proof}
\subsection{Conclusion of the existence part}
\label{subsec: existence for Dp and Rp} 

We now guide the reader through collecting and extending by density the respective estimates in order to obtain the existence part in our main results. 

\medskip

\noindent\emph{Existence of a solution with the properties in Theorem~\ref{thm: blockdir}}. First, let $f \in \L^p\cap \L^2$ if $p>1$ and $f \in a^{-1}(\H^1 \cap \L^2)$ if $p=1$. Then $u(t,x) \coloneqq \e^{-t {L^{1/2}}}f(x)$ is a weak solution to $\Le u = 0$ in $\reu$, see Proposition~\ref{prop: Poisson smg is weak solution}, and parts of (i) - (iv) are contained in the previous sections:
\begin{center}
\renewcommand{\arraystretch}{1.5}
\begin{longtable}{|c|c|}
\hline
Part & Obtained in \\
\hhline{|=|=|}
(i) & Propositions~\ref{prop: NT for Dirichlet} \& \ref{prop: SFE Dirichlet} \\
\hline
(ii) & Proposition~\ref{prop: NT for Dirichlet} \\
\hline
(iii) & Proposition~\ref{prop: sup estimates Dirichlet} (including quantitative  \\[-7pt]  
& bounds on the $t$-derivatives)\\
\hline
(iv) & Lemma~\ref{lem: sup L2loc estimates Dirichlet} \& qualitative continuity with values in $\L^2$ \\[-7pt]  &(by the functional calculus on $\L^2$)\\
\hline
\end{longtable}
\vspace{-25pt}
\end{center}
The non-tangential convergence with $\L^2$-averages in (ii) is stronger than what is asked for in $(D)_p^\Le$. Hence, $u$ solves $(D)_p^\Le$ with data $f$.

Now, consider general data $f \in \L^p$ if $p>1$ and $f \in a^{-1} \H^1$ if $p=1$. Take any sequence of data $(f_k) \subseteq \L^2$ that approximates $f$ in the data space as $k \to \infty$. Here, $a^{-1}\H^1$ is considered as a subspace of $\L^1$ with natural norm $\|a \cdot\|_{\H^1}$. Denote the corresponding solutions by $u_k$. 

By (i), we have that $(u_k)$ is a Cauchy sequence in $\T^{0,p}_\infty$ and that $(t \nabla u_k)$ is a Cauchy sequence in $\T^p$. Both topologies are stronger than $\Lloc^2(\reu)$. Hence, $(u_k)$ has a limit $u$ in $\Lloc^2(\reu)$ that satisfies (i) and it follows from Lemma~\ref{lem: weak solution L1loc convergence} that $u$ is a weak solution to  $\Le u = 0$. Note that this construction is independent of the choice of the $(f_k)$. In the same way we obtain (iii) and (iv) for $u$ since we can identify limits for the respective topologies in $\Lloc^1(\reu)$. 

Property (ii) for $u$ can be obtained by a well-known argument for maximal functions. More precisely, we obtain from (ii) for the $u_k$ that for a.e.\ $x \in \R^n$,
\begin{align}
\label{eq: extension of Whitney convergence}
\begin{split}
\limsup_{t \to 0} \bigg(\bariint_{W(t,x)} &|u-f(x)|^2 \, \d s \d y\bigg)^{\frac{1}{2}} \\
&\leq \NT(u-u_k)(x) + |f(x)-f_k(x)|.
\end{split}
\end{align}
If the left-hand side exceeds a fixed threshold $\eps>0$, then at least one of the terms on the right exceeds $\nicefrac{\eps}{2}$. By (i) applied to $u-u_k$ and Markov's inequality, this can only happen on a set of measure 
\begin{align*}
	C\eps^{-p} (\|a(f-f_k)\|_{\H^p}+\|f-f_k\|_{\L^p}),
\end{align*}
which tends to $0$ as $k \to \infty$ since $\H^p \subseteq \L^p$ with continuous inclusion for $p\geq1$. Hence, the left-hand side of \eqref{eq: extension of Whitney convergence} vanishes for a.e.\ $x \in \R^n$.

Finally, suppose that $f$ is also an admissible datum for energy solutions. In the case $p>1$ this means that we assume $f \in \L^p \cap \Hdot^{\nicefrac{1}{2},2}$ and by the universal approximation technique in Hardy--Sobolev spaces we can take the $f_k$ above in such a way that $f_k \to f$ also in $\Hdot^{\nicefrac{1}{2},2}$. We know from Proposition~\ref{prop: smg solution is compatible} that $u_k$ is the energy solution with Dirichlet data $f_k$ and it follows from Proposition~\ref{prop: existence of energy solution} that $u$ is the energy solution with Dirichlet data $f$.

In the case $p = 1$ we assume $f \in (a^{-1} \H^1) \cap \Hdot^{\nicefrac{1}{2},2}$. We claim that this is a subspace of $\L^2$. Taking the claim for granted, no approximation is necessary to construct the solution $u(t,x) = \e^{-t L^{1/2}}f(x)$ and by Proposition~\ref{prop: smg solution is compatible} this is the energy solution with data $f$. The easiest way to see the claim is to note that $f \in \L^1 \cap \Hdot^{\nicefrac{1}{2},2}$ and hence its Fourier transform satisfies
\begin{align*}
	\int_{\R^n} |\cF f(\xi)|^2 \, \d \xi 
	&\leq \int_{B(0,1)} \|f\|_1^2  \, \d \xi + \int_{{}^cB(0,1)} |\xi| |\cF f(\xi)|^2 \, \d \xi \\
	&\leq C \|f\|_1^2 + \|f\|_{\Hdot^{1/2,2}}^2 .
\end{align*}

\medskip

\noindent\emph{Existence of a solution with the properties in Theorem~\ref{thm: blockreg}}. First, recall from Theorem \ref{thm: standard relation J(L) and N(L)} that $p_-(L)=q_-(L)$. Let $f \in \Hdot^{1,p}\cap \W^{1,2}$. As before, $u(t,x) \coloneqq \e^{-t {L^{1/2}}}f(x)$ is a weak solution to $\Le u = 0$ in $\reu$ from Proposition~\ref{prop: Poisson smg is weak solution}  and (i) as well as (iii) - (v) are contained in the previous sections. Part (ii) will mostly follow from a general trace theorem that we comment on below:
\begin{center}
\renewcommand{\arraystretch}{1.5}
\begin{longtable}{|c|c|}
\hline
Part & Obtained in \\
\hhline{|=|=|}
(i) & Propositions~\ref{prop: NT for regularity} \& \ref{prop: SFE Regularity} \& Theorem \ref{thm: Riesz complete}.(i)\\
\hline
(iii) & Proposition~\ref{prop: NT for regularity}
\\
\hline
(iv) & Proposition~\ref{prop: sup estimates regularity} \\
\hline
(v) & Proposition~\ref{prop: NT for Dirichlet} \& Proposition~\ref{prop: sup estimates Dirichlet} \& Theorem~\ref{thm: Riesz complete} \\[-7pt]
&since $\partial_t u = -\e^{-t {L^{1/2}}} ({L^{1/2}} f)$ with ${L^{1/2}} f \in a^{-1} (\H^p \cap \L^2)$\\
\hline
\end{longtable}
\vspace{-25pt}
\end{center}

As for the extension to general data $f \in \Hdot^{1,p}$, we first treat the case $p<n$. We can assume $f \in \L^{p^*}$ since the general case follows by modifying data and solution by the same additive constant. 

Take any sequence $(f_k) \subseteq \Hdot^{1,p} \cap \W^{1,2}$ with $f_k \to f$ as $k \to \infty$ in $\Hdot^{1,p} \cap \L^{p^*}$. It follows from (iv) that $(u_k)$ is a Cauchy sequence in $\C([0,\infty); \L^{p^*})$, hence in $\Lloc^1$. Lemma~\ref{lem: weak solution L1loc convergence} asserts that $(u_k)$ converges in $\Wloc^{1,2}$ to a weak solution to $\Le u = 0$. The properties (i), (iv), (v) for $u$ follow by identifying limits as before and for (iii) we rely on the same type of density argument as in \eqref{eq: extension of Whitney convergence}. In particular, (iv) implies $\lim_{t\to 0} u(t,\cdot) = f$ in $\cD'$ as claimed in (ii). This being said, the non-tangential limit in (ii) follows from the Kenig--Pipher trace theorem (Proposition~\ref{prop: KP}). In conclusion, $u$ solves $(R)_p^\Le$ and has all required properties.

In the case $p\geq n$ we can only take a sequence $(f_k) \subseteq \Hdot^{1,p} \cap \W^{1,2}$ with $f_k \to f$ in $\Hdot^{1,p}$ as $k \to \infty$ . We use (i) to infer that for the corresponding solutions $(\nabla u_k)$ converges in $\T_\infty^{0,p}$, hence in $\Lloc^2$. Define the averages $c_k \coloneqq (u_k)_W$ with $W \subseteq \reu$ a fixed cube. By Poincar\'e's inequality $(u_k - c_k)$ is bounded in $\Wloc^{1,2}$. By compactness, we can define, up to passing to a subsequence,
\begin{align*}
u \coloneqq \lim_{k \to \infty} u_k - c_k \quad (\text{in } \Lloc^2).
\end{align*}
Lemma~\ref{lem: weak solution L1loc convergence} asserts again that $u$ is a weak solution to $\Le u = 0$ and modulo constants the construction of $u$ is independent of the particular choice of the $(f_k)$. With this definition all properties but (ii) follow as before. For the latter we fix the representative for $f$. Since $n > p_-(L)$, see  Proposition~\ref{prop: J(L) contains neighborhood of Sobolev conjugates}, we obtain from (v) that $\partial_t u \in \C_0([0,\infty); \L^p)$. Hence, $u(t,\cdot)$ has a limit in $\cD'$ as $t\to 0$. By (iv) we can fix the free constant for $u$ such that this limit is $f$ and the non-tangential convergence follows again from Proposition~\ref{prop: KP}.

Finally, if $f \in \Hdot^{1,p} \cap \Hdot^{\nicefrac{1}{2},2}$, then the same argument as for the Dirichlet problem yields that modulo constants $u$ is the energy solution with Dirichlet datum $f$.
\section{Existence in the Dirichlet problems with \texorpdfstring{$\Lamdot^\alpha$}{Hölder}-data}
\label{sec: existence of Holder-dir}

\noindent Here, we establish the existence part of Theorem~\ref{thm: Holder-dir}, our main result on the Dirichlet problems $(D)_{\Lamdot^\alpha}^\Le$ and $(\wtD)_{\Lamdot^\alpha}^\Le$ with boundary data in $\Lamdot^\alpha$.  Let us stress that in accordance with the formulation of these problems the data space is \emph{not} considered modulo constants.

Since $\Lamdot^\alpha \cap \L^2$ is not dense in $\Lamdot^\alpha$ for the strong topology, we cannot proceed in two well-separated steps as in the previous section. Instead, given $f \in \Lamdot^\alpha$, we directly define
\begin{align}
\label{eq: Holder-dir solution}
u(t,\cdot) \coloneqq \sum_{j = 1}^\infty \e^{-tL^{1/2}}(\ind_{C_j(Q)} f) \quad (t>0), 
\end{align}
where $Q \subseteq \R^n$ is any cube, and check that this is a solution with all required properties for both Dirichlet problems. More concisely, we can write
\begin{align*}
	u(t,\cdot) = \lim_{j \to \infty} \e^{-t L^{1/2}} (\ind_{2^{j+1}Q}f) \quad (t>0),
\end{align*}
but the representation as a series will be advantageous for most considerations. In fact, the assumptions of Theorem~\ref{thm: Holder-dir} are already required to prove convergence in $\Lloc^2$ via off-diagonal estimates. More precisely, we work with the following exponents for most of the section:
\begin{align}
\label{eq: assumptions Holder-Dir}
\begin{minipage}{0.88\linewidth}
	\begin{itemize}
		\item $p_+(L) > n$ and  $0 \leq \alpha < 1 - \nicefrac{n}{p_+(L)}$.
		\item When $\alpha$ is fixed, $p$ denotes a fixed exponent with $2 \leq p < p_+(L)$ and $\alpha < 1 - \nicefrac{n}{p}$.
	\end{itemize}
\end{minipage}
\end{align}
We break the argument into six parts.
\subsection*{Part 1: Well-definedness of the solution}
We begin with an elementary oscillation estimate.

\begin{lem}
\label{lem: averages of Lamdot-alpha}
Let $\alpha \in [0,1)$ and $p \in [1,\infty)$. For all $f \in \Lamdot^\alpha$, all cubes $Q \subseteq \R^n$ and all $j \geq 1$, it follows that
\begin{align*}
\bigg(\barint_{2^j Q} |f - (f)_Q|^p \, \d y\bigg)^{\frac{1}{p}} \lesssim \gamma_j \ell(Q)^\alpha \|f\|_{\Lamdot^\alpha},
\end{align*}
where $\gamma_j \coloneqq j +1$ if $\alpha = 0$ and $\gamma_j \coloneqq 2^{\alpha j}$ if $\alpha>0$.\index{$\gamma_j$ (oscillation estimate)}
\end{lem}

\begin{proof}
If $\alpha = 0$, then $\Lamdot^\alpha = \BMO$ and hence for all cubes $Q \subseteq \R^n$,
\begin{align*}
\bigg(\barint_{2Q} |f - (f)_Q|^p \, \d y\bigg)^{\frac{1}{p}} \lesssim \|f\|_{\Lamdot^0}.
\end{align*}
A telescopic sum of the estimates for $Q, 2Q,\ldots,2^{j-1}Q$ yields the claim. If $\alpha > 0$, then $|f(x)-f(y)| 
\lesssim (2^j \ell(Q))^\alpha \|f\|_{\Lamdot^\alpha}$ for $x \in Q$ and $y \in 2^j Q$ and the claim follows immediately.
\end{proof}

The oscillation estimate allows us to prove convergence of the right-hand side in \eqref{eq: Holder-dir solution} and obtain further useful representations of $u$.

\begin{lem}
\label{lem: well-definedness for solution Holder-dir}
Assume \eqref{eq: assumptions Holder-Dir}. Then the following hold true.
\begin{enumerate}
	\item The sum defining $u$ converges absolutely in $\Lloc^p(\R^n)$, locally uniformly in $t$. In particular, $u$ is a weak solution to $\Le u = 0$ in $\reu$. 
	\item If a family $(\eta_j) \subseteq \L^\infty(\R^n; \IC)$ satisfies \eqref{eq: conservation family}, then $u(t,\cdot) = \sum_{j = 1}^\infty \e^{-tL^{1/2}}(\eta_j f)$
 	with absolute convergence in $\Lloc^p(\R^n)$, locally uniformly in $t$. In particular, $u$ is independent of $Q$.
 	\item If $f=c$ is constant, then $u=c$ almost everywhere.
\end{enumerate}
\end{lem}

\begin{proof}
By Corollary~\ref{cor: Poisson semigroup Lp OD} the Poisson semigroup satisfies $\L^p$ off-diagonal estimates of order $1$. Let $K \subseteq \R^n$ be any compact set and set  $\ell \coloneqq \ell(Q)$. For $j$ large enough we have $\dist(K, C_j(Q)) \geq 2^{j-1} \ell$ and hence
\begin{align}
\label{eq1: well-definedness for solution Holder-dir}
\begin{split}
	\|\e^{-tL^{1/2}}&(\ind_{C_j(Q)} f)\|_{\L^p(K)} \\
	&\lesssim t (2^{j} \ell)^{-1} \|f\|_{\L^p(C_j(Q))} \\
	&\lesssim t (2^{j} \ell)^{-1} \Big(\|f - (f)_Q\|_{\L^p(2^{j+1}Q)} + (2^{j} \ell)^{\frac{n}{p}} |(f)_Q| \Big) \\
	& \lesssim t (2^{j} \ell)^{\frac{n}{p}-1} 
	\Big(\ell^\alpha \gamma_j \|f\|_{\Lamdot^\alpha} +|(f)_Q|\Big),
\end{split}
\end{align}
where we have used Lemma~\ref{lem: averages of Lamdot-alpha} in the final step. The right-hand side is summable in $j$ since $\alpha < 1 - \nicefrac{n}{p}$, which proves absolute convergence of the series in \eqref{eq: Holder-dir solution} in  $\Lloc^p$, locally uniformly in $t$. Since all partial sums are weak solutions to the equation for $\Le$ in $\reu$, the same is true for $u$, see Proposition~\ref{prop: smg solution is compatible} and Lemma~\ref{lem: weak solution L1loc convergence}. This completes the proof of (i).

Now, (ii) follows by repeating the proof of Proposition~\ref{prop: general OD extension} word by word  up to incorporating the new off-diagonal estimate above. Finally, (iii) is due to the conservation property for Poisson semigroups (Proposition~\ref{prop: conservation Poisson}).
\end{proof}
\subsection*{Part 2: Proof of (ii)}

We start by proving continuity and convergence towards the boundary data in $\Lloc^2$.

\begin{lem}
\label{lem: Holder-Dir solution Lloc2 continuous}
The solution $u$ is of class $\C([0,T]; \Lloc^2)$ with $u(0, \cdot) = f$ for every $T>0$.
\end{lem}

\begin{proof}
Continuity on $(0,T]$ is a general property of weak solutions, see Corollary~\ref{cor: weak solution smooth in t}. We fix an arbitrary cube $Q$ of sidelength $\ell$ and prove the limit at $t=0$ in $\L^2(Q)$. 

Set $f_j \coloneqq (f-(f)_Q) \ind_{C_j(Q)}$. By Lemma~\ref{lem: well-definedness for solution Holder-dir} we have, whenever $y \in Q$ and $s>0$,
\begin{align}
	\label{eq1: Holder-Dir solution Lloc2 continuous}
	\begin{split}
		u(s,y) - f(y)
		&= \sum_{j=1}^\infty \e^{-sL^{1/2}}f_j(y) + (f)_Q - f(y) \\
		&= \sum_{j=2}^\infty \e^{-sL^{1/2}}f_j(y) + (\e^{-sL^{1/2}}f_1(y) - f_1(y)).
	\end{split}
\end{align}
For the error terms with $j \geq 2$ we use again that the Poisson semigroup satisfies $\L^p$ off-diagonal estimates of order $1$, see Corollary~\ref{cor: Poisson semigroup Lp OD}. Here, $p$ is as in \eqref{eq: assumptions Holder-Dir}. Together with Lemma~\ref{lem: averages of Lamdot-alpha}, we obtain
\begin{align}
\label{eq2: Holder-Dir solution Lloc2 continuous}
	\begin{split}
		\bigg \|\sum_{j=2}^\infty \e^{-sL^{1/2}}f_j \bigg\|_{\L^p(Q)}
		&\leq \sum_{j=2}^\infty \|\e^{-s ^{1/2}}f_j\|_{\L^p(Q)} \\
		&\lesssim \sum_{j=2}^\infty \frac{s}{2^j \ell} \|f_j\|_{\L^p(Q)} \\
		&\lesssim \frac{s}{\ell^{1-\frac{n}{p}- \alpha}} \sum_{j=2}^\infty 2^{j(\frac{n}{p}-1)} \gamma_j \|f\|_{\Lamdot^\alpha},
	\end{split}
\end{align}
where the sum in $j$ is finite by the choice of $p$. In particular, we have by H\"older's inequality that
\begin{align*}
		\bigg \|\sum_{j=2}^\infty \e^{-sL^{1/2}}f_j \bigg\|_{\L^2(Q)}
		&\lesssim \frac{s}{\ell^{1-\frac{n}{2}- \alpha}} \|f\|_{\Lamdot^\alpha},
\end{align*}
which in combination with \eqref{eq1: Holder-Dir solution Lloc2 continuous} leads us to
\begin{align*}
	\|u(s,\cdot) - f\|_{\L^2(Q)} \lesssim  \frac{s}{\ell^{1-\frac{n}{2}- \alpha}}  \|f\|_{\Lamdot^\alpha} +   \|\e^{-sL^{1/2}}f_1 - f_1\|_2.
\end{align*}
The right-hand side tends to $0$ in the limit as $s \to 0$ since we have $f \in \Lamdot^\alpha$ and $f_1 \in \L^2$.
\end{proof}

We turn to non-tangential convergence towards the boundary data and control of the corresponding sharp functional on Whitney averages. In the case $\alpha > 0$ this would come for free from Proposition~\ref{prop: NT trace Y} once we have established the upper bound for the Carleson functional as stated in (i) but the following direct argument also works for $\alpha = 0$.

\begin{lem}
\label{lem: Holder-Dir solution NT}
The solution $u$ satisfies 
\begin{align*}
	 \lim_{t \to 0} \bariint_{W(t,x)} |u(s,y)-f(x)|^2 \, \d s \d y = 0 \quad (\text{a.e. } x\in \R^n)
\end{align*}
and
\begin{align*}
	\| \NTsharpalpha(u-f)\|_\infty \lesssim \|f\|_{\Lamdot^\alpha}.
\end{align*}
\end{lem}

\begin{proof}
We only need a slight refinement of the previous argument. To this end let $x \in \R^n$, $\ell \geq t$ and let $Q$ be the axis-parallel cube of sidelength $\ell$ centered at $x$.  

For any $(s,y) \in W(t,x) = (\nicefrac{t}{2}, 2t) \times B(x,t)$ we can use \eqref{eq1: Holder-Dir solution Lloc2 continuous} and \eqref{eq2: Holder-Dir solution Lloc2 continuous} with this choice of $Q$ and the same definition of $f_j$, $j \geq 1$, in order to obtain
\begin{align*}
	&\|u(s,\cdot) - f \|_{\L^2(B(x,t))}\\
	&\quad \lesssim t^{\frac{n}{2}-\frac{n}{p}} \bigg \|\sum_{j=2}^\infty \e^{-sL^{1/2}}f_j \bigg\|_{\L^p(Q)} + \|\e^{-sL^{1/2}}f_1 - f_1\|_{\L^2(B(x,t))}\\
	&\quad \lesssim \frac{s t^{\frac{n}{2}- \frac{n}{p}}}{\ell^{1-\frac{n}{p}- \alpha}} \|f\|_{\Lamdot^\alpha} + \|\e^{-sL^{1/2}}f_1 - f_1\|_{\L^2(B(x,t))}.
\end{align*}
Thus, we get our key estimate
\begin{align}
\label{eq1: Holder-Dir solution NT}
\begin{split}
	&\bigg(\bariint_{W(t,x)} |u(s,\cdot)- f|^2 \, \d s \d y \bigg)^{1/2} \\
	& \quad \lesssim \frac{t^{1-\frac{n}{p}}}{\ell^{1-\frac{n}{p}-\alpha}} \|f\|_{\Lamdot^\alpha} + \bigg(\bariint_{W(t,x)} |\e^{-sL^{1/2}}f_1 - f_1|^2 \, \d s \d y \bigg)^{1/2}.
\end{split}
\end{align}

For the first claim it suffices (by the Lebesgue differentiation theorem) to prove that the left-hand side in \eqref{eq1: Holder-Dir solution NT} vanishes in the limit as $t \to 0$ for a.e.\ $x \in \R^n$. But passing to the limit on the right-hand side, the first term vanishes since we have $p>n$ by \eqref{eq: assumptions Holder-Dir} and the second term vanishes for a.e.\ $x\in \R^n$ thanks to the Lebesgue differentiation theorem and Proposition~\ref{prop: NT for Dirichlet} applied to $f_1 \in \L^2$.

In order to bound the sharp functional, we use \eqref{eq1: Holder-Dir solution NT} with $\ell = t$. This yields for all $t>0$ and all $x \in \R^n$ the required uniform bound
\begin{align*}
	\frac{1}{t^\alpha} &\bigg(\bariint_{W(t,x)} |u(s,\cdot)- f|^2 \, \d s \d y \bigg)^{1/2}  \\
	&\lesssim \|f\|_{\Lamdot^\alpha} + \frac{1}{t^\alpha} \bigg(\bariint_{W(t,x)} |\e^{-sL^{1/2}}f_1 - f_1|^2 \, \d s \d y \bigg)^{1/2} \\
	&\lesssim \|f\|_{\Lamdot^\alpha} +  \frac{1}{t^{\alpha+\frac{n}{2}}} \sup_{s>0} \|\e^{-sL^{1/2}}f_1 - f_1\|_2 \\
	&\lesssim \|f\|_{\Lamdot^\alpha} + \frac{1}{t^{\alpha+\frac{n}{2}}}  \|f_1\|_2 \\
	&\lesssim \|f\|_{\Lamdot^\alpha},
\end{align*}
where the final step is due to Lemma~\ref{lem: averages of Lamdot-alpha}, keeping in mind that by definition $f_1 = (f-(f)_Q) \ind_{4Q}$ and that $t$ is the sidelength of $Q$. 
\end{proof}
\subsection*{Part 3: The upper bound for the Carleson functional}

In this part we prove the upper bound $\|C_\alpha(t \nabla u)\|_\infty \lesssim \|f\|_{\Lamdot^\alpha}$. It will be convenient to use cubes instead of balls for the Carleson functional and to show that for all cubes $Q \subseteq \R^n$ of sidelength $\ell$ we have
\begin{align}
\label{eq: Goal Part 3 Holder-dir}
\bigg(\int_0^\ell \barint_Q |s \nabla u|^2 \; \frac{\d y \d s}{s} \bigg)^{1/2} \lesssim \ell^\alpha \|f\|_{\Lamdot^\alpha}.
\end{align}

From now on $Q$ is fixed. Since both sides stay the same under adding constants to $u$ and $f$, we can assume $(f)_Q = 0$. For $j\geq1$ we introduce
\begin{align*}
f_j \coloneqq \ind_{C_j(Q)}f, \quad u_j(t,\cdot) \coloneqq \e^{-tL^{1/2}}f_j.
\end{align*}

\noindent \emph{Step~1: The local bound}. By Lemma~\ref{lem: averages of Lamdot-alpha} we have $\|f_1\|_2^2 \lesssim |Q| \ell^{2 \alpha}\|f\|_{\Lamdot^\alpha}^2$. Hence, the local term $u_1$ can be handled via the $\L^2$-bound for the square function in Proposition~\ref{prop: SFE Dirichlet}:
\begin{align*}
\int_0^\ell \barint_Q |s \nabla u_1|^2 \; \frac{\d y \d s}{s} 
&\leq |Q|^{-1} \iint_{\reu} |s \nabla \e^{-tL^{1/2}}f_1|^2 \; \frac{\d s \d y}{s}  \\
&\lesssim \ell^{2 \alpha} \|f\|_{\Lamdot^\alpha}^2.
\end{align*}

\noindent \emph{Step~2: Decomposition of the non-local terms}. Set $W(t,x) \coloneqq (t,2t) \times Q(x,t)$ and $\wt{W}(t,x) \coloneqq (\nicefrac{t}{2},4t) \times Q(x,2t)$. Let $\phi(z) \coloneqq \e^{-\sqrt{z}} - (1+z)^{-2}$ and recall from \eqref{eq0: SFE Dirichlet} the Caccioppoli estimate
\begin{align}
\label{eq0: Carleson bound Dirichlet}
\begin{split}
\iint_{W(t,x)} |s \nabla &\phi(s^2L)f_j|^2 \, \frac{\d s  \d y}{s} \\
&\lesssim \iint_{\wt{W}(t,x)} |\phi(s^2L)f_j|^2 + |\psi(s^2 L)f_j|^2 \, \frac{\d s \d y}{s},
\end{split}
\end{align}
where $\psi \in \Psi_1^1$ on any sector. Let the regions $(W(t_k, x_k))_{k}$ cover $(0,\ell) \times Q$ modulo a set of measure zero such that the $(\wt{W}(t_k,x_k))_{k}$ are contained in $(0,2\ell) \times 2Q$ and at most $4\cdot 2^{n}$ of them overlap at each point. Summing up in $k$ yields\index{inequality!Caccioppoli on a Carleson box}
\begin{align*}
\int_0^\ell \int_{Q} &|s\nabla \phi(s^2L)f_j|^2 \, \frac{\d y  \d s}{s} \\
&\lesssim \int_0^{2\ell} \int_{2Q} |\phi(s^2L)f_j|^2 + |\psi(s^2 L)f_j|^2 \, \frac{\d y  \d s}{s},
\end{align*}
so that in total
\begin{align*}
&\bigg(\int_0^\ell \barint_{Q} |s \nabla u_j|^2 \, \frac{\d y  \d s}{s}\bigg)^{\frac{1}{2}}\\
&\lesssim \bigg(\int_0^{2\ell} \barint_{2Q} |\phi(s^2L)f_j|^2 + |\psi(s^2 L)f_j|^2 + |s \nabla (1+s^2 L)^{-2}f_j|^2 \, \frac{\d y  \d s}{s} \bigg)^{\frac{1}{2}}.
\end{align*}
From Lemma~\ref{lem: well-definedness for solution Holder-dir} and Caccioppoli's inequality we obtain that $u =\sum_{j=1}^\infty u_j$ converges in $\Wloc^{1,2}(\reu)$. We can use Fatou's lemma to conclude
\begin{align}
\label{eq1: Carleson bound Dirichlet}
\bigg(\int_0^\ell \barint_{Q} |s \nabla u|^2 \, \frac{\d y  \d s}{s}\bigg)^{\frac{1}{2}}
\lesssim \ell^\alpha \|f\|_{\Lamdot^\alpha} + \sum_{j=2}^\infty \I_j + \II_j + \III_j,
\end{align}
where
\begin{align*}
\I_j &\coloneqq \bigg(\int_0^{2\ell} \barint_{2Q} |\phi(s^2L)f_j|^2 \frac{\d y  \d s}{s}\bigg)^{\frac{1}{2}},\\
\II_j &\coloneqq \bigg(\int_0^{2\ell} \barint_{2Q} |\psi(s^2L)f_j|^2 \frac{\d y  \d s}{s}\bigg)^{\frac{1}{2}},\\
\III_j &\coloneqq \bigg(\int_0^{2\ell} \barint_{2Q} |s \nabla (1+s^2 L)^{-2}f_j|^2 \frac{\d y  \d s}{s}\bigg)^{\frac{1}{2}}.\\
\end{align*}

\noindent \emph{Step~3: Bounds for the off-diagonal pieces}. We begin with the bound for $\I_j$. The family $(\phi(t^2 L))_{t > 0}$ satisfies $\L^p$ off-diagonal estimates of order $1$. This is due to Lemma~\ref{lem: functional calculus bounds from J(L) abstract} since $\phi \in \Psi_{1/2}^2$ on any sector. Hence,
\begin{align*}
\bigg(\barint_{2Q} |\phi(s^2L)f_j|^2 \, \d y\bigg)^{\frac{1}{2}}
&\leq \bigg(\barint_{2Q} |\phi(s^2L)f_j|^p \, \d y\bigg)^{\frac{1}{p}}\\
&\lesssim \Big(\frac{2^j \ell }{s} \Big)^{-1} 2^{j \frac{n}{p}}\bigg(\barint_{2^j Q} |f|^p \, \d y\bigg)^{\frac{1}{p}} \\
&\lesssim s \ell^{\alpha-1} \gamma_j 2^{j (\frac{n}{p} -1)} \|f\|_{\Lamdot^\alpha},
\end{align*}
where the final step is again due to Lemma~\ref{lem: averages of Lamdot-alpha}. We take $\L^2$-norms with respect to $\frac{\d s}{s}$ on both sides to give
\begin{align}
\label{eq2: Carleson bound Dirichlet}
\I_j &\leq \ell^\alpha \|f\|_{\Lamdot^\alpha}  \gamma_j 2^{j (\frac{n}{p} -1)}.
\end{align}
Summing these estimates in $j$ leads to a desirable bound in \eqref{eq1: Carleson bound Dirichlet}. 

In estimating $\I_j$ we have only used $\phi \in \Psi_{1/2}^{\tau}$ on any sector for some $\tau >0$. Hence, we can use the same strategy for $\II_j$ and the first component of
\begin{align*}
s \nabla (1+s^2 L)^{-2}f_j
=\begin{bmatrix} -4s^2 L (1+s^2 L)^{-3}f_j \\ s \nabla_x (1+s^2 L)^{-2}f_j \end{bmatrix}
\end{align*}
in $\III_j$. As for the second component, we have $\L^2$ off-diagonal estimates of arbitrarily large order $\gamma>0$ for $(t \nabla_x (1+t^2 L)^{-2})_{t>0}$ by composition. Therefore, we can run the same argument as before but with $p=2$ in Lemma~\ref{lem: averages of Lamdot-alpha} and obtain
\begin{align*}
\bigg(\int_0^{2\ell} \barint_{2Q} |s \nabla_x (1+s^2 L)^{-2}f_j|^2 \, \frac{\d y \d s}{s} \bigg)^{\frac{1}{2}}
&\leq \ell^\alpha \|f\|_{\Lamdot^\alpha} \gamma_j 2^{j (\frac{n}{2} - \gamma)}.
\end{align*}
We take $\gamma \coloneqq \nicefrac{n}{2} - \nicefrac{n}{p} +1$ and conclude a desirable bound for $\III_j$ in \eqref{eq1: Carleson bound Dirichlet}. This completes the proof of \eqref{eq: Goal Part 3 Holder-dir}.
\subsection*{Part 4: Compatibility}

In this section we work with $\Lamdot^\alpha$ as a homogeneous smoothness space modulo constants. In view of Lemma~\ref{lem: well-definedness for solution Holder-dir} this determines $u$ modulo constants. 

Our goal is to establish compatibility of $u$ with the energy class, that is, we assume $f \in \Lamdot^\alpha \cap \Hdot^{\nicefrac{1}{2},2}$ and have to show that modulo constants $u$ is the energy solution with Dirichlet data $f$. This is a delicate matter since no density argument can help us here. We shall rely on the following two lemmata.

\begin{lem}
\label{lem: Carleson plus L2 is constant}
Let $g_1 \in \L^2$ and $g_2 \in \T^{-1,\infty;\alpha}$ for some $\alpha \in [0,1)$ be such that $g_1-g_2$ is constant on $\reu$. Then $g_1=g_2$ almost everywhere.
\end{lem}

\begin{proof}
Let $g_1-g_2 = c$ almost everywhere. We obtain for all $r>0$ that
\begin{align*}
|c|^2 
&\simeq r^{-1-n} \int_r^{2r} \int_{B(0, 2r)} |g_1-g_2|^2 \, \d x \d t \\
&\lesssim r^{-1-n} \|g_1\|_2^2 + r^{2 \alpha - 2} \|g_2\|_{\T^{-1,\infty;\alpha}}^2.
\end{align*}
As $\alpha < 1$, sending $r \to \infty$ yields $c=0$.
\end{proof}

\begin{lem}
\label{lem: LP decomposition of boundary data}
Let $\alpha \in [0,1)$. Each $f \in \Lamdot^\alpha \cap \Hdot^{\nicefrac{1}{2},2}$ can be decomposed in $\Lamdot^\alpha \cap \Hdot^{\nicefrac{1}{2},2}$ as $f = f_{\loc} + f_{\glob}$, where $f_{\loc} \in \Wdot^{1,2}$ and $f_{\glob} \in \L^2$.
\end{lem}

\begin{proof}
We pick $\varphi \in \C_0^\infty(\R^n; \R)$ such that $\ind_{B(0,1)} \leq \varphi \leq \ind_{B(0,2)}$ and set
\begin{align*}
f_{\loc} \coloneqq \cF^{-1} (\varphi \cF f), \quad f_{\glob} \coloneqq \cF^{-1} ((1-\varphi) \cF f).
\end{align*}
Then obviously $f=f_{\loc} + f_{\glob}$ and since $\varphi$ and $1-\varphi$ are smooth Fourier multipliers in the scope of the Mihlin multiplier theorem, both $f_{\loc}$ and $f_{\glob}$ remain in $\Lamdot^\alpha \cap \Hdot^{\nicefrac{1}{2},2}$. Moreover, $m_{\loc}(\xi) \coloneqq |\xi|^{1/2} \varphi(\xi)$ and $m_{\glob}(\xi) \coloneqq |\xi|^{-1/2} (1-\varphi(\xi))$ are bounded and since we have $g \coloneqq \cF^{-1}(|\xi|^{1/2} \cF f) \in \L^2$ by assumption, we obtain that
\begin{align*}
\cF^{-1} (|\xi| \cF f_{\loc}) 
= \cF^{-1} (m_{\loc} \cF g) \in \L^2, \quad 
f_{\glob}= \cF^{-1} (m_{\glob} \cF g) \in \L^2
\end{align*}
as required.
\end{proof}
	
As we are dealing with a linear problem, the benefit from Lemma~\ref{lem: LP decomposition of boundary data} is that it suffices to prove compatibility under the additional assumption that either $f \in \L^2$ or $f \in \Wdot^{1,2}$. 

If additionally $f \in \L^2$, then $\sum_{j=1}^\infty \ind_{C_j(Q)} f$ converges to $f$ in $\L^2$ and from \eqref{eq: Holder-dir solution} we get back 
\begin{align*}
u(t,\cdot) = \e^{-t L^{1/2}}f \quad (t>0).
\end{align*}
According to Proposition~\ref{prop: smg solution is compatible} this is the energy solution with Dirichlet data $f$.
	
Now, suppose that additionally $f \in \Wdot^{1,2}$ and let $\wt{u}$ be the energy solution with Dirichlet data $f$. We claim that it suffices to show that for all $g \in \C_0^\infty$ with $\int_{\R^n} g \d x = 0$ and all $t>0$ we have
\begin{align}
\label{eq: Goal Holder-dir compatibility}
\langle u(t,\cdot), g \rangle = \langle \wt{u}(t,\cdot), g \rangle,
\end{align}
where the angular brackets denote the (extended) inner product on $\L^2$. Indeed, the claim implies that $u-\wt{u}$ is independent of the $x$-variable but looking at the equation $\Le (u - \wt{u}) = 0$ in $\reu$, we also obtain $a \partial_t^2 (u-\wt{u}) = 0$, so $\partial_t u- \partial_t \wt{u}$ is constant. By definition we have $\partial_t \wt{u} \in \L^2$ and by the Carleson bound in Part~3 we have $\partial_t u \in \T^{-1,\infty; \alpha}$. Lemma~\ref{lem: Carleson plus L2 is constant} yields $\partial_t u- \partial_t \wt{u} = 0$ and the desired compatibility $u = \wt{u}$ (modulo constants) follows.

In order to prove \eqref{eq: Goal Holder-dir compatibility}, we pick a cube $Q$ that contains the support of $g$ and use Lemma~\ref{lem: well-definedness for solution Holder-dir} to write
\begin{align}
\label{eq0: Holder-dir compatible}
u(t,x) = \sum_{j=1}^\infty \e^{-t L^{1/2}} (\eta_j f)(x) \quad ((t,x) \in \reu),
\end{align}
with $(\eta_j)_j$ a smooth partition of unity on $\R^n$ subordinate to the sets $D_1 \coloneqq 4Q$ and $D_j \coloneqq 2^{j+1}Q \setminus 2^{j-1}Q$, $j \geq 2$, such that $\|\eta_j\|_\infty + 2^j \ell(Q) \|\nabla_x \eta_j\|_\infty \leq C$ for a dimensional constant $C$. 
	
Since $g$ has integral $0$, we can write $g = \div_x G$ with $G \in \C_0^\infty(Q)$. Indeed, in dimension $n=1$ it suffices to take a suitable primitive of $g$ and in dimension $n \geq 2$ this is Bogovski\u{i}'s lemma~\cite[Lemma~III.3.1]{Galdi}\index{Bogovski\u{i}'s lemma}. By duality and the intertwining relations, we obtain 
\begin{align}
\label{eq1: Holder-dir compatible}
\begin{split}
\langle \e^{-t L^{1/2}} (\eta_j f), g \rangle
&= \langle \eta_j f,  \e^{-t (L^*)^{1/2}} \div_x G \rangle \\
&= \langle \eta_j f,  \div_x \e^{-t (M^\sharp)^{1/2}} G \rangle \\
&= -\langle \eta_j \nabla_x f, \e^{-t (M^\sharp)^{1/2}} G \rangle \\
&\quad - \langle \nabla_x \eta_j \otimes f, \e^{-t (M^\sharp)^{1/2}} G \rangle \\
& \eqqcolon - \I_j - \II_j,
\end{split}
\end{align}
where $M^\sharp \coloneqq -d^* \nabla_x (a^*)^{-1} \div_x$ intertwines with $L^*$ in the same ways as $M$ intertwines with $\tL$. Our notation is $\nabla_x (\eta_j f) = \eta_j \nabla_x f + \nabla_x \eta_j \otimes f$ as predicted by the product rule. The assumption $\nabla_x f \in \L^2$ and the fact that $\e^{-t (M^\sharp)^{1/2}} G \in \L^2$ allow us to sum up
\begin{align}
\label{eq2: Holder-dir compatible}
\sum_{j=1}^\infty \I_j = \langle \nabla_x f, \e^{-t (M^\sharp)^{1/2}} G \rangle.
\end{align}
	
As for the error terms $\II_j$, we shall need the qualitative information 
\begin{align}
\label{eq3: Holder-dir compatible}
e^{-t (M^\sharp)^{1/2}} G \in \L^q \quad (\text{for some } q<2).
\end{align}
In each of the following steps we take $q$ as close to $2$ as necessary for the respective result to apply. First, we write $G = G_1 + G_2$ with $G_1 \in \nul(\div_x)$ and $G_2 \in \cl{\ran(d^* \nabla_x)}$ as in the Hodge decomposition \eqref{eq: Hodge} with $d^*$ replacing $d$. By Proposition~\ref{prop: P(L) vs Hodge} and Lemma~\ref{lem: Hodge interval via isomorphism}, this decomposition can be taken topological in $\L^q$. The identification Theorem~\ref{thm: main result Hardy} tells us that we can have $\IH^{1,q}_{L^\sharp} = \Wdot^{1,q} \cap \L^2$ with equivalent $q$-norms and then $\IH_{M^\sharp}^q = \L^q \cap \cl{\ran(d^* \nabla_x)}$ follows by moving from the second to the fourth row in Figure~\ref{fig: diagram}. Proposition~\ref{prop: FC on pre-Hardy} yields $\e^{-t (M^\sharp)^{1/2}} G_2 \in \L^q$ and from $G_1 \in \nul(M^\sharp)$ we obtain by the functional calculus in $\L^2$ that  $\e^{-t (M^\sharp)^{1/2}} G_1 = G_1$, which also belongs to $\L^q$. Hence, \eqref{eq3: Holder-dir compatible} follows.
	
Now, we go back to \eqref{eq1: Holder-dir compatible}. We pick exponents $r,s \in (1,\infty)$ such that $\nicefrac{1}{q}+ \nicefrac{1}{r} + \nicefrac{1}{s} = 1$ and obtain for all $J \geq 1$ that
\begin{align*}
\bigg|\sum_{j=1}^J \II_j \bigg|
&\lesssim \bigg\|\sum_{j=1}^J \nabla_x \eta_j \bigg\|_r \|f\|_{\L^s(2^{J+1}Q)} \|\e^{-t (M^\sharp)^{1/2}} G\|_q \\
&\lesssim 2^{J(\frac{n}{r}-1)} \|f\|_{\L^s(2^{J+1}Q)},
\end{align*}
where we have used that $\sum_{j=1}^J \nabla_x \eta_j$ has support in $2^{J+1} Q$ and is controlled in $\L^\infty$-norm by $2^{-J} \ell(Q)$. The implicit constant depends on all variables but  $J$. The choice of $s$ depends on Sobolev embeddings. In dimension $n \geq 3$ we can assume $f \in \L^{2^*}$ up to modifying $f$ (and hence $u$) by a constant. Then we pick $s \coloneqq 2^*$ and obtain
\begin{align*}
\bigg|\sum_{j=1}^J \II_j \bigg| 
\lesssim 2^{J(\frac{n}{r}-1)} = 2^{J(\frac{n}{2}-\frac{n}{q})},
\end{align*}
which tends to $0$ as $J \to \infty$ since $q<2$. In dimension $n \leq 2$ we can assume $f \in \Lamdot^{1-n/2}$ and also change $f$ to $f-(f)_Q$ in \eqref{eq0: Holder-dir compatible}, which changes $u$ by a constant. With this modification, we  obtain together with Lemma~\ref{lem: averages of Lamdot-alpha} that
\begin{align*}
\bigg|\sum_{j=1}^J \II_j \bigg| 
\lesssim 
2^{J(\frac{n}{r}-1)} 2^{J\frac{n}{s}} \gamma_J 
= \begin{cases}
2^{J(\frac{1}{2}- \frac{1}{q})} &\text{if } n=1 \\
2^{J(1-\frac{2}{q})} (1+J) &\text{if } n=2 \\
\end{cases},
\end{align*}
which also tends to $0$ as $J \to \infty$. Together with \eqref{eq0: Holder-dir compatible} - \eqref{eq2: Holder-dir compatible}, we arrive at
\begin{align*}
\langle  u(t,\cdot), g \rangle = - \langle \nabla_x f, \e^{-t (M^\sharp)^{1/2}} G \rangle.
\end{align*}

Since $f \in \Wdot^{1,2} \cap \Hdot^{\nicefrac{1}{2},2}$, the universal approximation technique lets us pick a sequence $(f_k) \subseteq \Wdot^{1,2} \cap \Hdot^{\nicefrac{1}{2},2} \cap \L^2$ with $f_k \to f$ in both $\Wdot^{1,2}$ and $\Hdot^{\nicefrac{1}{2},2}$. We let $u_k$ be the energy solution with Dirichlet data $f_k$.
Then $(u_k)$ tends to the energy solution $\wt{u}$ with data $f$ in $\Wdot^{1,2}(\reu)$. By Lemma~\ref{lem: energy trace}, this implies $u_k(t,\cdot) \to u(t, \cdot)$ in the sense of distributions modulo constants. On the other hand, we know from Proposition~\ref{prop: smg solution is compatible} that $u_k(t, \cdot) = \e^{- t L^{1/2}} f_k$ and we can undo the duality and intertwining in order to give
\begin{align*}
\langle  u(t,\cdot), g \rangle 
&= - \lim_{k \to \infty} \langle \nabla_x f_k, \e^{-t (M^\sharp)^{1/2}} G \rangle \\
&= \lim_{k \to \infty} \langle \e^{-t L^{1/2}} f_k, g \rangle \\
&= \lim_{k \to \infty} \langle u_k(t,\cdot), g \rangle \\
&= \langle  \wt{u}(t,\cdot), g \rangle.
\end{align*}
This establishes the remaining claim \eqref{eq: Goal Holder-dir compatibility} and the proof is complete.
\subsection*{Part 5: The lower bound for the Carleson functional}

Our goal is to show that for all $g \in \C_0^\infty$ with $\int_{\R^n} g \d x = 0$ the solution $u$ in \eqref{eq: Holder-dir solution} satisfies
\begin{align}
\label{eq: Holder-Dir goal}
|\langle f, g \rangle| \lesssim \|C_\alpha(t \nabla u)\|_\infty \|g\|_{\H^{\varrho}},
\end{align}
where $\varrho \in (1_*,1]$ is such that $n(\nicefrac{1}{\varrho}-1) = \alpha$ and $\langle \cdot\,, \cdot \rangle$ is the extended $\L^2$-duality pairing. Indeed, then density and duality yield the lower bound 
\begin{align*}
\|f\|_{\Lamdot^\alpha} \lesssim \|C_\alpha(t \nabla u)\|_\infty.
\end{align*}
We suggest that the reader recalls the argument of Step~3 of Proposition~\ref{prop: SFE Dirichlet} beforehand. The proof here follows the same line of thought but since $u(t,\cdot)$ and $f$ may not be globally in $\L^2$, we cannot as directly rely on the functional calculus in $\L^2$ as before. This is the major technical challenge. 

From now on we fix $g$ and pick a cube $Q$ that contains its support. Since both sides in \eqref{eq: Holder-Dir goal} do not change when adding constants to $f$ or $u$, we can assume $(f)_Q = 0$ and write $u$ as
\begin{align*}
u(t,\cdot) = \sum_{j=1}^\infty \e^{-t L^{1/2}}f_j, \quad f_j \coloneqq \ind_{C_j(Q)}f.
\end{align*}
Next, we introduce again $H\coloneqq -(a^*)^{-1} \Delta_x$ and set
\begin{align*}
v(t,\cdot) \coloneqq (1+t^2 H)^{-\beta} ((a^*)^{-1} g)
\end{align*}
for an integer $\beta > \nicefrac{n}{2} + 2$. Then the kernel estimates in Proposition~\ref{prop: pointwise estimates} become available and this is why we use the resolvents of $H$ and not the Poisson semigroup as in the proof of Proposition~\ref{prop: SFE Dirichlet}. The auxiliary function we are working with is
\begin{align}
\label{eq: Holder-dir Phi function}
\begin{split}
\Phi: (0, \infty) \to \IC, \quad \Phi(t) \coloneqq \sum_{j=1}^\infty \langle \e^{-t L^{1/2}}f_j, a^*v(t,\cdot) \rangle.
\end{split}
\end{align}
This turns out to be the appropriate way of defining $\langle u(t,\cdot), a^*v(t,\cdot) \rangle$ as we shall see momentarily. We divide the proof of \eqref{eq: Holder-Dir goal} into eight steps.

\medskip

\noindent \emph{Step~1: Qualitative growth bounds for $v$}. We claim that there are $c>0$ and $C>0$ depending also on $\beta$, $g$ and $Q$ such that 
\begin{align}
\label{eq: Step 1 Holder-dir lower bound}
|v(t,x)| + |t \nabla v(t,x)| \leq C (1 \wedge t^{-n-1}) \e^{-c \frac{\dist(x,Q)}{t}} \quad ((t,x) \in \reu).
\end{align}

To this end, we recall that $(1+t^2H)^{-\beta}(a^*)^{-1}$ is given by an integral kernel denoted by $H_{t}^\beta(x,y)$ with bounds
\begin{align*}
|H_{t}^\beta(x,y)|+ |t\nabla_x H_{t}^\beta(x,y)|+ |t\nabla_y H_{t}^\beta(x,y)|
\\
+  |t^2\nabla_x \nabla_y H_{t}^\beta(x,y)| \le Ct^{-n}\e^{-c\frac{|x-y|}{t}},
\end{align*}
see Proposition~\ref{prop: pointwise estimates}. Hence, by the support of $g$ and an $\L^1$-bound on the kernel,
\begin{align*}
|v(t,x)| 
&\leq \e^{-\frac{c}{2} \frac{\dist(x,Q)}{t}} \int_{\R^n} \e^{\frac{c}{2} \frac{|x-y|}{t}}|H_{t}^\beta(x,y)| |g(y)| \, \d y\\
&\lesssim \e^{-\frac{c}{2} \frac{\dist(x,Q)}{t}} \| g\|_\infty.
\end{align*}
Since $g \in \C_0^\infty(Q)$ has mean value $0$, we can also write $g = \div_x F$ with $F \in \C_0^\infty(Q)$, using a suitable primitive in dimension $n=1$ and Bogovski\u{i}'s lemma in higher dimensions. Thus,
\begin{align*}
v(t,x) = - \int_{\R^n} \nabla_y H_{t}^\beta(x,y) \cdot F(y) \, \d y
\end{align*}
and the $\L^\infty$-bound for the kernel yields
\begin{align*}
|v(t,x)| \leq Ct^{-n-1}\e^{-c\frac{\dist(x,Q)}{t}} \|F\|_1.
\end{align*}
This completes the estimate in \eqref{eq: Step 1 Holder-dir lower bound} for $v$. The bounds for $t \nabla_x v$ follow \emph{mutadis mutandis}, using the kernel bounds for $t \nabla_x H_{t}^\beta$ and $t^2 \nabla_x \nabla_y H_{t}^\beta$. Eventually,
\begin{align*}
t \partial_t v(t,\cdot) 
&= - 2 \beta t^2 H(1+t^2 H)^{-\beta - 1} ((a^*)^{-1}g) \\
&= - 2 \beta ((1+t^2 H)^{-\beta } - (1+t^2 H)^{-\beta - 1})((a^*)^{-1}g) 
\end{align*}
is a linear combination of two functions of the same type as $v$.

\medskip

\noindent \emph{Step~2: $\Phi$ is well-defined}. More precisely, we shall show the qualitative bound
\begin{align}
\label{eq: Step 2 Holder-dir lower bound}
\sum_{j=2}^\infty \| |\e^{-t L^{1/2}}f_j| |a^*v(t, \cdot)| \|_1 \leq C (t \wedge t^{-\frac{n}{p}}) < \infty,
\end{align}
where $C$ is independent of $t>0$ but may depend on all other parameters.

By H\"older's inequality, we have
\begin{align}
\label{eq: Holder-dir 3}
\begin{split}
\| |\e^{-t L^{1/2}}f_j| |a^*v(t, \cdot)| \|_1
&\leq \|\ind_{2^{j-1/2}Q} \e^{-t L^{1/2}}f_j\|_p \|a^*v(t, \cdot)\|_{p'}  \\
&\quad + \|\e^{-t L^{1/2}}f_j\|_p \|\ind_{{}^c(2^{j-1/2}Q)}a^*v(t, \cdot)\|_{p'}.
\end{split}
\end{align}
Since $p \in [2, p_+(L))$, the Poisson semigroup satisfies $\L^p$ off-diagonal estimates of order $1$, see Corollary~\ref{cor: Poisson semigroup Lp OD}. From the support of $f_j$ and Lemma~\ref{lem: averages of Lamdot-alpha} we obtain for $j \geq 2$ with implicit constants independent of $j$ and $t$,
\begin{align}
\label{eq: Holder-dir 3bis}
\|\ind_{2^{j-1/2}Q} \e^{-t L^{1/2}}f_j\|_p 
\lesssim t 2^{-j} \|f_j\|_p 
\lesssim t \gamma_j 2^{j(\frac{n}{p}-1)}
\end{align}
and
\begin{align*}
\|\e^{-t L^{1/2}}f_j\|_p 
\lesssim \|f_j\|_p 
\lesssim \gamma_j 2^{j\frac{n}{p}}.
\end{align*}
Likewise, integrating the $p'$-th powers of both sides of \eqref{eq: Step 1 Holder-dir lower bound} gives 
\begin{align*}
\|a^*v(t, \cdot)\|_{p'} 
\lesssim (1 \wedge t^{-n-1}) (1+t^{\frac{n}{p'}}) 
\lesssim 1 \wedge t^{-1-\frac{n}{p}}
\end{align*}
and, with a smaller constant $c$ then in \eqref{eq: Step 1 Holder-dir lower bound},
\begin{align*}
\|\ind_{{}^c(2^{j-1/2}Q)} a^*v(t, \cdot)\|_{p'} 
\lesssim (1 \wedge t^{-n-1}) t^{\frac{n}{p'}} \e^{-c \frac{2^j}{t}}  
\lesssim (t^{n+1-\frac{n}{p}} \wedge t^{-\frac{n}{p}}) 2^{-j},
\end{align*}
where in the final step we have used the crude bound $\e^{-s} \lesssim s^{-1}$ for $s>0$ in order to restore the right homogeneity in $t$. Using these bounds on the right-hand side of \eqref{eq: Holder-dir 3}, leads us to
\begin{align}
\label{eq: Holder-dir 4}
\| |\e^{-t L^{1/2}}f_j| |a^*v(t, \cdot)| \|_1
\lesssim  (t \wedge t^{-\frac{n}{p}}) \gamma_j 2^{j\frac{n}{p}-1}.
\end{align}
Since $\alpha < 1 - \nicefrac{n}{p}$, we can sum in $j$ and conclude \eqref{eq: Step 2 Holder-dir lower bound}. 

As a matter of fact, the same estimate holds if in the definition of $\Phi$ we replace $v(t,\cdot)$ by $t \nabla v(t,\cdot)$, which satisfies the same pointwise bounds. We can also replace $\e^{-t L^{1/2}}$ by $(t L^{1/2})^k \e^{-t L^{1/2}}$ for an integer $k \geq 1$ since the latter satisfies again $\L^p$ off-diagonal estimates of order at least $1$, see Lemma~\ref{lem: functional calculus bounds from J(L) abstract}. All such sums are called \emph{of $\Phi$-type}.\index{sum of $\Phi$-type} We also remark that it was only the bound \eqref{eq: Holder-dir 3bis} that required $j \geq 2$. All other estimates in this step also work for $j=1$.

\medskip

\noindent \emph{Step~3: Integration by parts in $t$}.
Since we have left out the term for $j=1$ in Step~2, the full estimate for $\Phi(t)$ is
\begin{align*}
|\Phi(t) - \langle \e^{-t L^{1/2}}f_1, a^*v(t,\cdot) \rangle| \lesssim t \wedge t^{-\frac{n}{p}}.
\end{align*}
By the functional calculus on $\L^2$ we have
\begin{align*}
\lim_{t \to 0} \langle \e^{-t L^{1/2}}f_1, a^*v(t,\cdot) \rangle = \langle f_1, g \rangle = \langle f, g \rangle,
\end{align*}
where in the final step we used the support of $f_1$, and likewise
\begin{align*}
\lim_{t \to \infty} \langle \e^{-t L^{1/2}}f_1, a^*v(t,\cdot) \rangle = 0.
\end{align*}
We conclude $\lim_{t \to 0} \Phi(t) = \langle f, g \rangle$ and $\lim_{t \to \infty} \Phi(t) = 0$. Next,
\begin{align*}
\frac{\d}{\d t} \langle \e^{-t L^{1/2}}f_j, a^*v(t,\cdot) \rangle 
&= -\langle L^{1/2} \e^{-t L^{1/2}}f_j, a^*v(t,\cdot) \rangle \\
&\quad + \langle \e^{-t L^{1/2}}f_j, a^* \partial_t v(t,\cdot) \rangle
\end{align*}
gives rise to two sums of $\Phi$-type (times a factor of $t^{-1}$), which converge locally uniformly in $t$ by Step~2. Hence, we can differentiate $\Phi$ term by term. The upshot is that we can integrate $\Phi$ by parts to obtain
\begin{align}
\label{eq: Step 3 Holder-dir lower bound}
\langle f, g \rangle = -\int_0^\infty \Phi'(t) \, \d t = \int_0^\infty \Phi^{(1)}(t) \, \d t - \int_0^\infty \Psi^{(1)}(t) \, \d t,
\end{align}
where $\Phi^{(1)}(t), \Psi^{(1)}(t): (0,\infty) \to \IC$ are given by
\begin{align*}
\Phi^{(1)}(t) &\coloneqq \sum_{j = 1}^\infty \langle L^{1/2} \e^{-t L^{1/2}}f_j, a^*v(t,\cdot) \rangle, \\
\Psi^{(1)}(t) &\coloneqq \sum_{j = 1}^\infty \langle \e^{-t L^{1/2}}f_j, a^* \partial_t v(t,\cdot) \rangle
\end{align*}
and $t \Phi^{(1)}$ and $t \Psi^{(1)}$ are of $\Phi$-type. The idea here is that $\Phi^{(1)}$ is the bad term that we have to keep, whereas the part involving $\Psi^{(1)}$ can be treated directly.

\medskip

\noindent \emph{Step~4: Integral estimate for $\Psi^{(1)}$}. We introduce 
\begin{align*}
	\wt{v}(t,\cdot) \coloneqq 2 \beta (1+t^2 H)^{-\beta-1}((a^*)^{-1} g),
\end{align*}
which is of the same type as $v$ but with a higher resolvent power. The objective in this step is to establish the bound
\begin{align}
\label{eq: Step 3bis Holder-dir lower bound}
\int_0^\infty |\Psi^{(1)}(t)| \, \d t 
\leq \iint_{\reu} |t\nabla_x u| \cdot |t\nabla_x \wt{v}| \, \frac{\d t  \d x}{t}.
\end{align}

Let $\eta \in \C_0^\infty(\R^n; \R)$ be such that $\ind_Q \leq \eta \leq \ind_{2Q}$ and for $R>0$ set $\eta_R(x) \coloneqq \eta(\nicefrac{x}{R})$. We note that
\begin{align*}
a^* \partial_t v(t,\cdot) 
&= -2\beta a^* tH(1+t^2 H)^{-\beta-1}((a^*)^{-1} g) \\
& = -2 \beta t \Delta_x (1+t^2 H)^{- \beta-1}((a^*)^{-1} g) \\
&\eqqcolon -t \Delta_x \wt{v}(t,\cdot)
\end{align*}
and, having split
\begin{align*}
\Psi^{(1)}(t) 
&= \sum_{j = 1}^\infty \langle \eta_R \e^{-t L^{1/2}}f_j, -t \Delta_x \wt{v}(t,\cdot) \rangle + \langle (1-\eta_R) \e^{-t L^{1/2}}f_j, a^*\partial_t v(t,\cdot) \rangle,
\end{align*}
we can integrate by parts the term with $\eta_R$ to give
\begin{align}
\label{eq: Holder-dir IPP}
\begin{split}
\Psi^{(1)}(t) 
&=\frac{1}{t} \langle \eta_R t \nabla_x u(t,\cdot), t \nabla_x \wt{v}(t,\cdot) \rangle \\
&\quad+\frac{1}{t} \sum_{j = 1}^\infty \langle (t \nabla_x \eta_R) \otimes \e^{-t L^{1/2}}f_j, t \nabla_x \wt{v}(t,\cdot) \rangle\\
&\quad +\frac{1}{t}\sum_{j = 1}^\infty \langle (1-\eta_R) \e^{-t L^{1/2}}f_j, t a^*\partial_t v(t,\cdot) \rangle.
\end{split}
\end{align}
Our notation is 
\begin{align*}
\nabla_x ( \eta_R \e^{-t L^{1/2}}f_j) = \eta_R \nabla_x \e^{-t L^{1/2}}f_j + (\nabla_x \eta_R) \otimes \e^{-t L^{1/2}}f_j
\end{align*}
as predicted by the product rule and for the sum with $\eta_R \nabla_x \e^{-t L^{1/2}}f_j$ we have used that the series that defines $u(t,\cdot)$ converges in $\Wloc^{1,2}$ as a consequence of $\Lloc^2$-convergence and the Caccioppoli inequality.

So far, \eqref{eq: Holder-dir IPP} holds for any $t>0$ and any $R>0$. We let now $k\geq 2$, set $R \coloneqq (1 \vee t)k$ and integrate in $t$ to obtain
\begin{align}
\label{eq: Holder-dir 5}
\begin{split}
\int_0^\infty &|\Psi^{(1)}(t)| \, \d t\\
&\leq \iint_{\reu} |t\nabla_x u| \cdot |t\nabla_x \wt{v}| \, \frac{\d t  \d x}{t} \\
&+ c_n \int_0^\infty \int_{{}^c(2Q)} |t \nabla_x \eta_{(1 \vee t)k}| \sum_{j = 1}^\infty |\e^{-t L^{1/2}}f_j| |t \nabla_x \wt{v}| \, \frac{\d x  \d t}{t} \\
&+ c_n \|a^*\|_\infty \int_0^\infty \int_{{}^c(2Q)} |1-\eta_{(1 \vee t)k}| \sum_{j = 1}^\infty |\e^{-t L^{1/2}}f_j| |t \partial_t v| \, \frac{\d x  \d t}{t},
\end{split}
\end{align}
where $c_n$ only depends on $n$. We also used that the terms with $\eta$ vanish on $2Q$ and interchanged the sum with the integral in $x$ using the monotone convergence theorem. 

The sums in $j$ are of $\Phi$-type and when using the bounds from Step~2 for such sums only on ${}^c(2Q)$, we can allow $j=1$ and pick up the same behavior in $t$. Indeed, on the right-hand side of \eqref{eq: Holder-dir 3} we would only get the second term when $j=1$, so that we do not need \eqref{eq: Holder-dir 3bis}. It follows that
\begin{align*}
\sum_{j = 1}^\infty \||\e^{-t L^{1/2}}f_j| |t \nabla_x \wt{v}|\|_{\L^1({}^c(2Q))} \lesssim t \wedge t^{-\frac{n}{p}}
\end{align*}
and likewise with $\partial_t v$ replacing $t \nabla_x \wt{v}$. Hence, in \eqref{eq: Holder-dir 5} the sums in $j$ are of class $\L^1((0,\infty) \times {}^c(2Q); \frac{\d x \d t}{t})$. Since $|t \nabla_x \eta_{(1 \vee t)k}|$ and $|1-\eta_{(1 \vee t)k}|$ are uniformly bounded in $t$ and tend to $0$ pointwise as $k \to \infty$, we can use the dominated convergence theorem in \eqref{eq: Holder-dir 5} to conclude \eqref{eq: Step 3bis Holder-dir lower bound}.

\medskip

\noindent \emph{Step~5: Completing the treatment of $\Psi^{(1)}$ by duality}. We can interpret the right-hand side in \eqref{eq: Step 3bis Holder-dir lower bound} as a $\T^{0,\infty;\alpha} - \T^\varrho$ duality pairing, where $\varrho \in (1_*,1]$ is such that $\alpha = n(\nicefrac{1}{\varrho}-1)$, see Section~\ref{subsec: tent spaces}. Consequently, we have
\begin{align*}
\int_0^\infty |\Psi^{(1)}(t)| \, \d t 
\lesssim \|C_\alpha (t\nabla_x u)\|_\infty \|S(t\nabla_x \wt{v})\|_\varrho.
\end{align*}
In order to bound the square function, let
\begin{align*}
B_H \coloneqq \begin{bmatrix} (a^*)^{-1} & 0 \\ 0 & 1 \end{bmatrix}
\end{align*} 
be the matrix that corresponds to $H$ in the same way as $B$ corresponds to $L$. Recalling \eqref{eq: BD and BD} and the intertwining relation \eqref{eq: similarity BD and DB}, we write
\begin{align}
\label{eq: Holder dir v-tilde via BH}
\begin{split}
\begin{bmatrix} 0 \\ t \nabla_x \wt{v} \end{bmatrix}
&= - 2 \beta tD(1+(tB_H D)^2)^{-\beta-1} B_H \begin{bmatrix} g \\ 0 \end{bmatrix} \\
&\eqqcolon \psi(tDB_H)\begin{bmatrix} g \\ 0 \end{bmatrix},
\end{split}
\end{align}
where $\psi(z) = -2\beta z(1+z^2)^{-\beta-1}$ is of class $\Psi_1^{2\beta +1}$ on any sector. As $\beta > \nicefrac{n}{2} +2$, this is an admissible auxiliary function for $\IH_{DB_H}^\varrho$. From $p_{-}(H) = 1_*$ (Corollary~\ref{cor: McIntosh-Nahmod}) and the identification theorem (Theorem~\ref{thm: main result Hardy}) we obtain
\begin{align*}
\|S(t\nabla_x \wt{v})\|_\varrho 
\simeq \bigg\|\begin{bmatrix} g \\ 0 \end{bmatrix} \bigg \|_{\IH_{DB_H}^\varrho} 
\simeq \|g\|_{\H^\varrho}.
\end{align*}
Thus, we have found
\begin{align}
\label{eq: Step 4 Holder-dir lower bound}
\int_0^\infty |\Psi^{(1)}(t)| \, \d t 
\lesssim \|C_\alpha (t\nabla_x u)\|_\infty \|g\|_{\H^\varrho},
\end{align}
which is a desirable bound for the second term in \eqref{eq: Step 3 Holder-dir lower bound}.

\medskip

\noindent \emph{Step~6: Setting up an iteration on $\Psi^{(1)}$.} At this point we are left with proving
\begin{align*}
\bigg|\int_0^\infty \Phi^{(1)}(t) \, \d t \bigg| \lesssim \|C_\alpha (t\nabla u)\|_\infty \|g\|_{\H^\varrho},
\end{align*}
where
\begin{align*}
\Phi^{(1)}(t) &= \sum_{j = 1}^\infty \langle L^{1/2} \e^{-t L^{1/2}}f_j, a^*v(t,\cdot) \rangle.
\end{align*}
Since $t\Phi^{(1)}(t)$ is of $\Phi$-type, we can repeat Step~3 with this function replacing $\Phi(t)$. The only difference is that now $\lim_{t\to 0} t\Phi^{(1)}(t) = 0$ and we can integrate by parts without boundary terms to give
\begin{align*}
\int_0^\infty  \Phi^{(1)}(t) \, \d t =  \int_0^\infty  t\Phi^{(2)}(t) \, \d t - \int_0^\infty  t\Psi^{(2)}(t) \, \d t,
\end{align*}
where
\begin{align*}
\Phi^{(2)}(t) &\coloneqq \sum_{j = 1}^\infty \langle L \e^{-t L^{1/2}}f_j, a^*v(t,\cdot) \rangle, \\
\Psi^{(2)}(t) &\coloneqq \sum_{j = 1}^\infty \langle L^{1/2} \e^{-t L^{1/2}}f_j, a^* \partial_t v(t,\cdot) \rangle.
\end{align*}
Now, $t^2 \Phi^{(2)}$ and $t^2 \Psi^{(2)}$ are of $\Phi$-type and $t \Psi^{(2)}$ is of the same structure as $\Psi^{(1)}$ except for an extra $t \partial_t$-derivative on the Poisson semigroup. Hence, we can repeat Step~4 and Step~5 \emph{mutadis mutandis} for $\Psi^{(2)}$ and arrive at
\begin{align*}
\int_0^\infty |t \Psi^{(2)}(t)| \, \d t 
\lesssim \|C_\alpha (t \nabla_x (t \partial_t u))\|_\infty \|g\|_{\H^\varrho}
\end{align*}
as replacement for \eqref{eq: Step 4 Holder-dir lower bound}. But since $\partial_t u$ is a weak solution to the same equation, we can use Caccioppoli's inequality on Carleson boxes $ (0, \ell(Q))  \times Q$ as in Part 3 to bound 
\begin{align*}
\|C_\alpha (t^2 \nabla_x \partial_t u)\|_\infty  \lesssim \|C_\alpha (t\partial_t u)\|_\infty
\end{align*} and conclude with a desirable bound.

The upshot is that we can iterate this scheme until for some large $N$, depending on the dimension, we can control
\begin{align}
\label{eq: Step 5 Holder-dir lower bound}
\bigg|\int_0^\infty t^{N-1} \Phi^{(N)}(t) \, \d t \bigg| \lesssim \iint_{\reu} |t^{N-1} \nabla_x \partial_t^{N-2} u| \cdot |t\nabla_x v| \, \frac{\d t  \d x}{t},
\end{align}
where
\begin{align*}
\Phi^{(N)}(t) &\coloneqq \sum_{j = 1}^\infty \langle L^{N/2} \e^{-t L^{1/2}}f_j, a^*v(t,\cdot) \rangle.
\end{align*}
Indeed, a desirable bound for the right-hand side of \eqref{eq: Step 5 Holder-dir lower bound} follows by $\T^{0,\infty;\alpha} - \T^\varrho$-duality and Caccioppoli's inequality as before.
\medskip

\noindent \emph{Step~7: Reduction to a final estimate of $\Phi$-type}. We shall establish \eqref{eq: Step 5 Holder-dir lower bound} for the first integer that satisfies $N>\nicefrac{n}{2}+3$. As
\begin{align*}
 \langle L^{N/2} \e^{-t L^{1/2}}f_j, a^*v(t,\cdot) \rangle = \langle -\div_x d \nabla_x L^{N/2-1} \e^{-t L^{1/2}}f_j, v(t,\cdot) \rangle,
\end{align*}
we can integrate by parts as in Step~4 but in the opposite direction. Using the same notation, the replacement for \eqref{eq: Holder-dir IPP} is
\begin{align*}
t^{N-1} &\Phi^{(N)}(t) \\
&=\frac{(-1)^{N-2}}{t} \langle d t^{N-1} \nabla_x \partial_t^{N-2} u(t,\cdot), \eta_R t \nabla_x v(t,\cdot) \rangle \\
&\quad+\frac{1}{t} \sum_{j = 1}^\infty \langle d t \nabla_x (t L^{1/2})^{N-2} \e^{-t L^{1/2}}f_j, (t\nabla_x \eta_R) \otimes v(t,\cdot) \rangle\\
&\quad +\frac{1}{t}\sum_{j = 1}^\infty \langle (t L^{1/2})^{N} \e^{-t L^{1/2}}f_j, (1-\eta_R)  a^* v(t,\cdot) \rangle
\end{align*}
and the replacement for \eqref{eq: Holder-dir 5} is
\begin{align*}
\int_0^\infty &|t^{N-1} \Phi^{(N)}(t)| \, \d t\\
&\leq \|d\|_\infty \iint_{\reu} |t^{N-1} \nabla_x \partial_t^{N-2} u| \cdot |t\nabla_x v| \, \frac{\d t  \d x}{t} \\
&+ c_n \int_0^\infty \int_{{}^c(2Q)} |t \nabla_x \eta_{(1 \vee t)k}| \sum_{j = 1}^\infty |t d\nabla_x (t L^{1/2})^{N-2} \e^{-t L^{1/2}}f_j| |v| \, \frac{\d x  \d t}{t} \\
&+ c_n \|a^*\|_\infty \int_0^\infty \int_{{}^c(2Q)} |1-\eta_{(1 \vee t)k}| \sum_{j = 1}^\infty |(t L^{1/2})^N \e^{-t L^{1/2}}f_j| |v| \, \frac{\d x  \d t}{t},
\end{align*}
where $c_n$ only depends on $n$. Thus, we have to prove that the second and third term on the right vanish in the limit as $k \to \infty$. The third term contains a sum of $\Phi$-type, so that we can use dominated convergence as in Step~4. The middle term is not of $\Phi$-type since we do not have $\L^p$ off-diagonal estimates for the gradient of the Poisson semigroup. We claim that nonetheless there are $\sigma,\tau > 0$ such that we have the qualitative bound
\begin{align}
\label{eq: Holder-dir ultimate goal}
\sum_{j=1}^\infty \||t \nabla_x (t L^{1/2})^{N-2} \e^{-t L^{1/2}}f_j| |v|\|_{\L^1({}^c (2Q))} \lesssim t^{\sigma} \wedge t^{- \tau}
\end{align}
with implicit constant independent of $t$. Taking this estimate for granted, dominated convergence also applies to the middle term and \eqref{eq: Step 5 Holder-dir lower bound} follows.

\medskip

\noindent \emph{Step~8: Conclusion}. In order to prove the final missing bound \eqref{eq: Holder-dir ultimate goal}, we argue as in Step~2 but with $p=2$. To simplify notation, let 
\begin{align*}
T(t) \coloneqq t \nabla_x (t L^{1/2})^{N-2} \e^{-t L^{1/2}} \quad (t>0).
\end{align*}
This family satisfies $\L^2$ off-diagonal estimates of order $N-2$ by composition and Lemma~\ref{lem: functional calculus bounds from J(L) abstract} since we can write
\begin{align*}
	T(t) = \Big(t \nabla_x (1+t^2 L)^{-1} \Big) \Big( (t L^{1/2})^{N-2} (1+t^2 L) \e^{-t L^{1/2}} \Big).
\end{align*}

By H\"older's inequality, we have
\begin{align}
\label{eq: Step 6 Holder-dir lower bound}
\begin{split}
\| |T(t) f_j| |v(t, \cdot)| &\|_{\L^1({}^c (2Q))} \\
&\leq \|\ind_{2^{j-1/2}Q} T(t) f_j\|_{\L^2({}^c (2Q))} \|v(t, \cdot)\|_{\L^2({}^c (2Q))}  \\
&\quad + \|T(t)f_j\|_2 \|\ind_{{}^c(2^{j-1/2}Q)}v(t, \cdot)\|_2,
\end{split}
\end{align}
where the first term on the right vanishes for $j=1$. From the support of $f_j$ and Lemma~\ref{lem: averages of Lamdot-alpha} we obtain for $j \geq 2$ that 
\begin{align*}
\|\ind_{2^{j-1/2}Q} T(t) f_j\|_{\L^2({}^c (2Q))}
\lesssim t^\gamma 2^{-j \gamma} \|f_j\|_2 
\lesssim t^\gamma \gamma_j 2^{j(\frac{n}{2}-\gamma)}
\end{align*}
and for $j \geq 1$ that
\begin{align*}
\|T(t)f_j\|_2
\lesssim \|f_j\|_2
\lesssim \gamma_j 2^{j\frac{n}{2}}
\end{align*}
with $\gamma \in (0, N-2]$ at our disposal and implicit constants independent of $j$ and $t$. The bounds for $v$ are obtained by squaring both sides of \eqref{eq: Step 1 Holder-dir lower bound} and integrating. They take the form
\begin{align*}
\|v(t, \cdot)\|_{\L^2({}^c (2Q))} 
\lesssim (1 \wedge t^{-n-1}) t^{\frac{n}{2}}
= t^{\frac{n}{2}} \wedge t^{-1-\frac{n}{2}}
\end{align*}
and
\begin{align*}
\|\ind_{{}^c(2^{j-1/2}Q)} v(t, \cdot)\|_{2} 
\lesssim (1 \wedge t^{-n-1}) t^{\frac{n}{2}} \e^{-\frac{c}{2} \frac{2^j}{t}}  
\lesssim (t^{\frac{n}{2} + \gamma} \wedge t^{\gamma-\frac{n}{2}-1}) 2^{-j \gamma},
\end{align*}
where in the final step we have used the crude bound $\e^{-s} \lesssim s^{-\gamma}$ for $s>0$ to restore the homogeneity in $t$. Using these bounds on the right-hand side of \eqref{eq: Step 6 Holder-dir lower bound}, we find
\begin{align*}
\| |T(t) f_j| |v(t, \cdot)| \|_{\L^1({}^c (2Q))} \lesssim  (t^{\frac{n}{2} + \gamma} \wedge t^{\gamma-\frac{n}{2}-1}) \gamma_j 2^{j(\frac{n}{2}-\gamma)}.
\end{align*}

We need $\gamma > \nicefrac{n}{2} + \alpha$ to be able to sum in $j$ and $\gamma <  \nicefrac{n}{2} + 1$ to pick up decay at $t=\infty$. Such $\gamma$ exists since $\alpha < 1$ and the choice is admissible because we have assumed $N> \nicefrac{n}{2} + 3$. It is only at this point where we need the size of $N$. Now,  \eqref{eq: Holder-dir ultimate goal} follows from \eqref{eq: Step 6 Holder-dir lower bound} and the proof is complete.
\subsection*{Part 6: Proof of (iii)}

Instead of \eqref{eq: assumptions Holder-Dir} we work with the following exponents in this part:
\begin{align}
\begin{minipage}{0.86\linewidth}
\label{eq: assumptions Holder-Dir strenghtened}
\begin{itemize}
\item $p_-(L^\sharp) < 1$ and  $0 \leq \alpha < n(\nicefrac{1}{p_-(L^\sharp)} -1)$.
\item When $\alpha$ is fixed, $p_-(L^\sharp) < p \leq 1$ is such that $\alpha = n(\nicefrac{1}{p} -1)$.
\end{itemize}
\end{minipage}
\end{align}
This is a stronger assumption than in the previous parts since $p_-(L^\sharp) < 1$ implies $p_+(L) = \infty$ by duality and similarity. 

In particular, $(\e^{-t (L^\sharp)^{1/2}})_{t>0}$ is $(a^*)^{-1} \H^p$-bounded by Theorem~\ref{thm: Poisson semigroup interval} and we can define $(\e^{-t L^{1/2}})_{t>0}$  as a bounded semigroup on $\Lamdot^\alpha$ via duality and similarity:
\begin{align*}
\langle \e^{-t L^{1/2}} f, g \rangle \coloneqq \langle f, a^* \e^{-t {(L^\sharp)^{1/2}}}(a^*)^{-1}g \rangle \quad (f \in \Lamdot^\alpha, g \in \H^p \cap \L^2).
\end{align*}
Next, we identify the solution $u$ from \eqref{eq: Holder-dir solution} with such a semigroup extension.

\begin{lem}
\label{lem: identification u dual semigroup}
Assume \eqref{eq: assumptions Holder-Dir strenghtened}. If $g \in \C_0^\infty$ with $\int_{\R^n} g \d x = 0$, then 
\begin{align*}
\langle u(t,\cdot), g \rangle = \langle f, a^* \e^{-t {(L^\sharp)^{1/2}}}(a^*)^{-1}g \rangle \quad (t>0),
\end{align*}
where the left-hand side is the (extended) $\L^2-\L^2$-duality and the right-hand side is the $\Lamdot^\alpha - \H^p$-duality.
\end{lem}

\begin{proof}
We fix $t$ and $g$ and let $Q$ be a cube that contains the support of $g$. As $a^* \e^{-t {(L^\sharp)^{1/2}}} (a^*)^{-1}g \in \H^p \cap \L^2 \subseteq \H^1$ we have in particular that $\int_{\R^n} a^* \e^{-t {(L^\sharp)^{1/2}}} (a^*)^{-1}g \d x = 0$. Therefore we can assume $(f)_Q = 0$. In the following, $C$ denotes a constant that may depend on all parameters but on $j \geq 1$ used for the annuli $C_j(Q)$. 

Since $p \in (1_*,1)$, we can fix $q \in (1,2)$ such that $\eps \coloneqq \nicefrac{n}{q}-\nicefrac{n}{p} + 1 > 0$. Then $(\e^{-t (L^\sharp)^{1/2}})_{t>0}$ is $\L^q$-bounded and satisfies $\L^q$ off-diagonal estimates of order $1$, see Corollary~\ref{cor: Poisson semigroup Lp OD}. We conclude that
\begin{align*}
\|a^* \e^{-t {(L^\sharp)^{1/2}}}(a^*)^{-1}g\|_{\L^q(C_j(Q))} 
\leq C2^{-j}
= C 2^{j(\frac{n}{q}-\frac{n}{p})} 2^{-\eps j}.
\end{align*}
Hence, we can use Lemma~\ref{lem: Taibleson-Weiss} in order to write 
\begin{align*}
a^* \e^{-t {(L^\sharp)^{1/2}}}(a^*)^{-1}g = \sum_{j=1}^\infty C 2^{-\eps j} a_j \quad (\text{in } \H^p \cap \Lloc^1),
\end{align*}
where the $a_j$ are $\L^q$-atoms for $\H^p$ with support in $C_{j+1}(Q) \cup C_{j}(Q)$. Using Lemma~\ref{lem: averages of Lamdot-alpha} with exponent $q'$ and the atomic bounds, we obtain
\begin{align}
\label{eq2: identification u dual semigroup}
|\langle f, C 2^{-\eps j} a_j \rangle| 
\leq C \gamma_j 2^{j \frac{n}{q'}} 2^{-\eps j} \|a_j\|_q
\leq  C \gamma_j 2^{- \alpha j} 2^{-\eps j}.
\end{align}
Now, we use the definition of $u$, duality for the semigroups on $\L^2$ and absolute convergence of the series following from \eqref{eq2: identification u dual semigroup} in order to write, setting $a_0 \coloneqq 0$,
\begin{align*}
\langle u(t,\cdot), g \rangle 
&= \sum_{j=1}^\infty \langle \e^{-tL^{1/2}} \ind_{C_j(Q)} f, g \rangle 
\\
&= \sum_{j=1}^\infty \langle \ind_{C_j(Q)} f, a^* \e^{-t {(L^\sharp)^{1/2}}}(a^*)^{-1}g  \rangle 
\\
&= \sum_{j=1}^\infty C \langle  f, \ind_{C_j(Q)} (2^{-\eps j} a_j+ 2^{-\eps (j-1)} a_{j-1}) \rangle 
\\
&= \sum_{j=1}^\infty C \langle  f, (\ind_{C_{j}(Q)} + \ind_{C_{j+1}(Q)}) 2^{-\eps j} a_j \rangle
\\
&= \sum_{j=1}^\infty \langle f, C 2^{-\eps j} a_j \rangle = \langle f, a^* \e^{-t {(L^\sharp)^{1/2}}}(a^*)^{-1}g \rangle. &\qedhere
\end{align*}
\end{proof}

Since $\C_0^\infty$-functions with integral zero are dense in $\H^p$,  we obtain from the lemma and Proposition~\ref{prop: sup estimates Dirichlet} applied to $L^\sharp$ that  $u$ is of class
\begin{align*}
\C_0([0,\infty);  \Lamdot^\alpha_{\text{weak}^*}) \cap \C^\infty((0,\infty); \Lamdot^\alpha_{\text{weak}^*}),
\end{align*}
where the subscript indicates that $\Lamdot^\alpha$ carries the $\text{weak}^*$ topology as the dual of $\H^p$, with bound
\begin{align*}
\sup_{t>0}  \|u(t,\cdot)\|_{\Lamdot^\alpha} \lesssim \|f\|_{\Lamdot^\alpha}.
\end{align*}
In the opposite direction, Part~2 implies for all $g \in \C_0^\infty$ with integral zero that
\begin{align*}
|\langle f,g \rangle| = \lim_{t \to 0} |\langle u(t,\cdot),g \rangle| \leq \sup_{t>0}  \|u(t,\cdot)\|_{\Lamdot^\alpha} \|g\|_{\H^p}
\end{align*}
and $\|f\|_{\Lamdot^\alpha} \leq \sup_{t>0}  \|u(t,\cdot)\|_{\Lamdot^\alpha}$ follows. Hence, we have
\begin{align}
\label{eq: sup Holder bound}
\sup_{t>0}  \|u(t,\cdot)\|_{\Lamdot^\alpha} \simeq \|f\|_{\Lamdot^\alpha}.
\end{align}

For the global upper bound in $\Lamdot^\alpha(\overline{\reu})$ we need a variant of the Poincar\'e inequality that we prove at the end of the section.\index{inequality!Poincar\'e on a Carleson box}

\begin{lem}
\label{lem:Poincare} 
Let $v\in \L^2_{\loc}({\reu})$ with $t^{1/2} \nabla v \in \L^2_{\loc}(\cl{\reu})$. There is a dimensional constant $c$ such that for all cubes $Q \subseteq \R^n$, 
\begin{align*}
\bariint_{T(Q)}  |v-(v)_{T(Q)}|^2 \, \d s \d y \le c \int_{0}^{\ell(Q)}  \barint_{Q}s |\nabla v|^2 \, \d y \d s,
\end{align*}
where $T(Q) \coloneqq (0, \ell(Q)) \times Q$. In particular, $v \in \L^2_{\loc}({\overline{\reu}})$. The same inequality holds with balls instead of cubes.
\end{lem}

Together with the upper Carleson bound of Part~3 we now obtain, for all cubes $Q \subseteq \R^n$,
\begin{align*}
\bigg(\bariint_{T(Q)}  |u-(u)_{T(Q)}|^2 \, \d s \d y\bigg)^{\frac{1}{2}}
&\lesssim \bigg(\int_{0}^{\ell(Q)} \barint_{Q}  s |\nabla u|^2 \, \d y \d s \bigg)^{\frac{1}{2}} \\
&\leq \ell(Q)^\alpha \|C_\alpha(t \nabla u)\|_\infty \\
&\lesssim \ell(Q)^\alpha \|f\|_{\Lamdot^\alpha}.
\end{align*}
This is an oscillation estimate at the boundary of $\reu$. In order to replace $T(Q)$ by an arbitrary cube $T(Q) + (t_0,0)$ with $t_0 > 0$, we use that according to Lemma~\ref{lem: identification u dual semigroup} we have the semigroup property $u(t+t_0,\cdot)= \e^{-tL^{1/2}}f_{t_0}=: u_{t_0}(t,\cdot)$, where $f_{t_0} \coloneqq u(t_0,\cdot) \in \Lamdot^\alpha$. The previous estimate with $u_{t_0}$ in place of $u$ becomes
\begin{align*}
\bigg(\bariint_{T(Q) + (t_0,0)}  |u-(u)_{T(Q) + (t_0,0)}|^2 \, \d s \d y\bigg)^{\frac{1}{2}}
&\lesssim \ell(Q)^\alpha \|f_{t_0}\|_{\Lamdot^\alpha} \\
&\lesssim \ell(Q)^\alpha \|f\|_{\Lamdot^\alpha},
\end{align*}
where the final step is due to \eqref{eq: sup Holder bound}. By definition of the $\BMO$-norm if $\alpha = 0$ and by the Morrey--Campanato characterization of H\"older continuity if $\alpha \in (0,1)$, see~\cite{Meyers-Holder}, we conclude
\begin{align*}
\|u\|_{\Lamdot^\alpha(\overline{\reu})} \lesssim \|f\|_{\Lamdot^\alpha}.
\end{align*}
The proof of (iii) is complete, modulo the

\begin{proof}[Proof of Lemma~\ref{lem:Poincare}]
We can assume that $Q$ is the unit cube centered at the origin, as a scaling argument gives the general result. Let $T_\eps(Q) \coloneqq (\eps,1)\times (1-\eps)Q$ for $\eps \in (0,1)$.  We apply first the Hardy--Poincaré inequality of Boas--Straube~\cite{Boas-Straube}:
\begin{align*}
\iint_{T_\eps(Q)}   |v-(v)_{T_\eps(Q)}|^2 \, \d s \d y 
\le c_\eps \iint_{T_\eps(Q)}  \dist((s,y), \partial T_\eps(Q)) |\nabla v|^2 \, \d s \d y.
\end{align*}
\emph{A priori}, the constant $c_\eps$ depends on $T_\eps(Q)$ but scaling and translation to $(1,2) \times Q$ reveals that we can take $c_\eps = (1-\eps) c$, where $c$ is dimensional. We conclude
\begin{align}
\label{eq1: Poincare}
\iint_{T_\eps(Q)} |v-(v)_{T_\eps(Q)}|^2 \, \d s \d y 
\le c \iint_{T(Q)}   s |\nabla v|^2 \, \d s \d y,
\end{align}
where the right-hand side is assumed to be finite. 

Now, consider a decreasing sequence of values $\eps \in (0, \nicefrac{1}{2})$ with $\eps \to 0$. Since $T_{1/2}(Q) \subseteq T_\eps(Q)$ and $v \in \L^2(T_{1/2}(Q))$, it follows from \eqref{eq1: Poincare} that the numerical sequence $((v)_{T_\eps(Q)})_\eps$ is bounded. Let $C$ be one of its accumulation points. Via Fatou's lemma we can pass to the limit in \eqref{eq1: Poincare} along a subsequence of $\eps$ to give
\begin{align*}
\iint_{T(Q)}   |v-C|^2 \, \d s \d y
\le c \iint_{T(Q)}   s |\nabla v|^2 \, \d s \d y.
\end{align*}
This implies that $v$ is (square) integrable on $T(Q)$ and therefore we have $C = (v)_{T(Q)}$ by dominated convergence.

The argument for balls instead of cubes is the same.
\end{proof}
\section{Existence in Dirichlet problems with fractional regularity data}
\label{sec: fractional}

\noindent In this section we prove the compatible existence on Dirichlet problems with data in homogeneous Hardy--Sobolev and Besov spaces of fractional smoothness that have been announced in Section~\ref{intro:fractionalDP}.  We also compare them to what can be obtained by the general first-order approach~\cite{AA} when specialized to elliptic systems in block form. We recall the color code for our various exponent regions and segments:
\begin{itemize}
	\item Gray  corresponds to what can be obtained from the theory of $DB$-adapted spaces in \cite{AA} and our identification of the interval from \cite{AE, AM} in Corollary~\ref{cor: AE interval and AM interval}.
	\item Blue shows additional results obtained from the theory of $L$-adapted spaces.
	\item Red  indicates results outside of the theory of operator-adapted spaces. 
\end{itemize}  
When we speak of `colored' points or regions, we always mean points or regions that are displayed in one of these three colors.
\subsection{Fractional identification regions}
\label{subsec: fractional identification}

As in Section~\ref{sec: Hardy intro} we treat adapted Hardy--Sobolev and Besov spaces simultaneously by letting $\X$ denote one of $\B$ or $\H$.\index{identification region!for $\IX_L^{s,p}$} As before, `identification'\index{identification!of abstract and concrete spaces} means `equality of sets with equivalent $p$-quasinorms'.

\begin{prop}
\label{prop: blue region}
Identification $\IX^{s,p}_L = \Xdot^{s,p} \cap \L^2$ holds for all exponents corresponding to the interior of the colored trapezoidal region in Figure~\ref{fig: diagram general proof}.
\end{prop}

\begin{figure}[ht]
\begin{center}
\begin{tikzpicture}[scale=2.4]
	
	\newcommand\fracspace{\vphantom{\frac{1}{1}}};
	\def\dimension{6};
	\def\xlength{3.5};
	\pgfmathsetmacro\xstretch{\xlength/(1+1/\dimension)}; 
	
	\pgfmathsetmacro\UpSob{\xstretch*(1/2-1/\dimension)-0.2};
	\pgfmathsetmacro\UpSobStar{\xstretch*(\UpSob/\xstretch-1/\dimension)};
	\pgfmathsetmacro\LowSob{\xstretch*(1/2+1/\dimension)+0.5};
	\pgfmathsetmacro\LowSobStar{\xstretch*(\LowSob/\xstretch+1/\dimension)};
	\pgfmathsetmacro\UpGradient{\xstretch*(1/2)-0.2};
	\pgfmathsetmacro\UpGradientDual{\xstretch*(1/2)+0.3};
	\pgfmathsetmacro\Half{\xstretch*0.5};
	\pgfmathsetmacro\One{\xstretch*1};
	\coordinate (P00) at (1+\UpGradient,2);
	\coordinate (P01) at (1+\LowSob,2);
	\coordinate (P10) at (1+\UpSob,0);
	\coordinate (P11) at (1+\UpGradientDual,0);
	\coordinate (P12) at (1+\LowSob,0);
	\coordinate (DirExtra) at (1+\UpSobStar,0);
	\coordinate (RegExtra) at (1+\LowSobStar,2);
	
	\draw [thin] (1,2) -- (1+\xlength,2); 
	\draw [thin] (1,0) -- (1+\xlength,0); 
	
	\draw [thick,->] (1,-0.5) -- (1+\xlength+0.2,-0.5);
	\draw [fill=black] (1+\Half,-0.5) circle [radius = .5pt];
	\node [below] at (1+\Half,-0.5) {$\frac{1}{2 \fracspace}$};
	\draw [fill=black] (1+\UpGradient,-0.5) circle [radius = .5pt];
	\node [below] at (1+\UpGradient,-0.5) {$\frac{1}{q_+^L\fracspace}$};
	\draw [fill=black] (1+\UpGradientDual,-0.5) circle [radius = .5pt];
	\node [below] at (1+\UpGradientDual,-0.5) {$\frac{1}{(q_+^{L^\sharp})'\fracspace}$};
	\draw [fill=black] (1+\UpSob,-0.5) circle [radius = .5pt];
	\node [below] at (1+\UpSob,-0.5) {$\frac{1}{p_+^L\fracspace}$};
	\draw [fill=black] (1+\LowSob,-0.5) circle [radius = .5pt];
	\node [below] at (1+\LowSob,-0.5) {$\frac{1}{p_-^L \vee 1 \fracspace}$};
	\draw [fill=black] (1+\LowSobStar,-0.5) circle [radius = .5pt];
	\node [below] at (1+\LowSobStar,-0.5) {$\frac{1}{{(p_-^L)_* \vee 1_*\fracspace}}$};
	\draw [fill=black] (1+\xlength,-0.5) circle [radius = .5pt];
	\node [below] at (1+\xlength,-0.5) {$\frac{n+1}{n \fracspace}$};
	\node [right] at (1+\xlength+0.2,-0.5) {$\frac{1}{p \fracspace}$};
	\draw [fill=black] (1,-0.5) circle [radius = .5pt];
	\node [below] at (1,-0.5) {$0$};
	
	\draw [thick,->] (0.7,0) -- (0.7,2.2);
	\node [above] at (0.7,2.2) {$s$};
	\draw [fill=black] (0.7,2) circle [radius = .5pt];
	\node [left] at (0.7,2) {$1$};
	\draw [fill=black] (0.7,0) circle [radius = .5pt];
	\node [left] at (0.7,0) {$0$};
	
	\draw [thin,dotted] (1,2) -- (1,0); 
	\draw [thin,dotted] (1+\LowSob,2) -- (1+\LowSob,0); 
	\draw [thin,dotted] (1+\LowSobStar,2) -- (1+\LowSobStar,0); 
	\draw [thin, dotted] (1+\Half, 2) -- (1+\Half,0);
	\draw [thin,dotted] (1+\UpGradient,2) -- (1+\UpGradient,0);
	\draw [thin, dotted] (1+\xlength,2) -- (1+\xlength,0); 
	
	\path [fill=lightgray, opacity = 0.6] (P00)--(P01)--(P11)--(P10)--(P00);
	\path [fill=blue!80!black, opacity = 0.4] (P11)--(P12)--(RegExtra)--(P01)--(P11);
	\draw [ultra thick, gray] (P00) -- (P01); 
	\draw [ultra thick, blue!80!black] (P01) -- (RegExtra); 
	\draw [ultra thick, gray] (P10) -- (P11); 
	\draw [ultra thick, blue!80!black] (P11) -- (P12);
	\draw [fill=white] (P00) circle [radius = .75pt];
	\draw [fill=blue!80!black] (P01) circle [radius = .75pt];
	\draw [fill=blue!80!black] (P11) circle [radius = .75pt];
	\draw [fill=white] (P12) circle [radius = .75pt];
	\draw [fill=white] (P10) circle [radius = .75pt];
	\draw [fill=white] (RegExtra) circle [radius = .75pt];
\end{tikzpicture}
\end{center}
\caption{Identification $\IX_L^{s,p} = \Xdot^{s,p} \cap \L^2$ up to equivalent $p$-quasinorms holds for all exponents corresponding to the interior of the colored trapezoidal region. The picture is up to scale when $p_-(L) \ge1$. When $p_-(L) < 1$, the top blue point is situated at $(\nicefrac{1}{p_-(L)}, 1)$.}
\label{fig: diagram general proof}
\end{figure}

\begin{proof} Theorem~\ref{thm: main result Hardy} yields $\IH^{1,p}_L = \Hdot^{1,p} \cap \L^2$	and $\IH_L^p = a^{-1}(\H^p\cap \L^2)=\L^p \cap \L^2$ if $(\nicefrac{1}{p},s)$ belongs to the open segments that join $(\nicefrac{1}{q_+(L)}, 1)$ to $(\nicefrac{1}{(p_-(L)_*\vee1_*)},1)$ and $(\nicefrac{1}{p_+(L)}, 0)$ to $(\nicefrac{1}{(p_-(L)\vee 1)},0)$, respectively. Both cases can be summarized as saying $\IH_L^{s,p} = \IH_{-\Delta_x}^{s,p}$, see Figure~\ref{fig: identification D-adapted}. By real and complex interpolation (\cite[Thm.4.32]{AA} or equivalently the argument in the proof of Lemma~\ref{lem: identification regions are interval}) we conclude $\IX_L^{s,p} = \IX_{-\Delta_x}^{s,p}$ in the interior of the convex hull of the two segments and the claim follows by using Figure~\ref{fig: identification D-adapted} again.
\end{proof}

\begin{rem}
\label{rem: blue region}
If $p_-(L) < 1$, then we could also combine Theorem~\ref{thm: main result Hardy} with Corollary~\ref{cor: McIntosh-Nahmod} and write $\IH_L^{s,p} = \IH_{-a^{-1}\Delta_x}^{s,p}$ on the top segment, which in this case joins $(\nicefrac{1}{q_+(L)}, 1)$ to $(\nicefrac{1}{1_*},1)$, and the full bottom segment joining $(\nicefrac{1}{p_+(L)}, 0)$ to $(\nicefrac{1}{p_-(L)},0)$.  Extending Figure~\ref{fig: diagram general proof} to the right by the triangle with vertices $(\nicefrac{1}{1_*},1), (1,0), (\nicefrac{1}{p_-(L)},0)$, the same interpolation argument yields identification $\IX^{s,p}_L = \IX^{s,p}_{-a^{-1} \Delta_x}$ in the interior of that extended region. The reason why we do not use this extension is that we do not know whether  $\IX^{s,p}_{-a^{-1} \Delta_x}=\IX^{s,p}_{-\Delta_x}$ and not even whether a completion of  $\IX^{s,p}_{-a^{-1} \Delta_x}$ can be realized as a space of distributions, except if $a=1$ of course. In the first-order $DB$-theory this phenomenon does not appear as $B$ is applied first. As a cautionary tale we remark that even when $a=1$ not all of our arguments for solvability of Dirichlet problems would go through in the extended gray region, notably the non-tangential trace used in Proposition~\ref{prop: solvability identification region}.
\end{rem}

Identification in the interior of the gray region in Figure~\ref{fig: diagram general proof} has previously been obtained (implicitly) in \cite[Sec.~7.2.4]{AA} and we shall next explain why. 

Let us first recall that in Theorem~\ref{thm: main result Hardy} we have identified $\IH_{DB}^{0,p} = \IH_D^{0,p}$ for $p \in (p_-(L), q_+(L))$. For $p \in (1,\infty)$, the  \emph{$\heartsuit$-duality}\index{duality!$\heartsuit$-duality} from \cite[Cor.~5.14]{AA} states that $\IH_{DB^*}^{0,p'} = \IH_D^{0,p'}$ implies $\IH_{DB}^{-1,p} = \IH_D^{-1,p}$. Thus, the latter follows for $p \in (q_+(L^\sharp)' , p_+(L^\sharp))$ by duality and similarity. As before, interpolation leads to the identification region that is shown in Figure~\ref{fig: diagram DB}\index{identification region!for $\IX_{DB}^{s,p}$}. Lemma~\ref{lem: intersection of identification regions} `maps' the gray region in Figure~\ref{fig: diagram DB} onto the gray region in Figure~\ref{fig: diagram general proof} since $\IX^{s,p}_{DB} = \IX_D^{s,p}$ implies $\IX^{s+1,p}_L = \Xdot^{s+1,p} \cap \L^2$.

\begin{figure}[ht]
\begin{center}
\begin{tikzpicture}[scale=2.4]
	
	\newcommand\fracspace{\vphantom{\frac{1}{1}}};
	\def\dimension{6};
	\def\xlength{3.5};
	\pgfmathsetmacro\xstretch{\xlength/(1+1/\dimension)}; 
	
	\pgfmathsetmacro\UpSob{\xstretch*(1/2-1/\dimension)-0.2};
	\pgfmathsetmacro\UpSobStar{\xstretch*(\UpSob/\xstretch-1/\dimension)};
	\pgfmathsetmacro\LowSob{\xstretch*(1/2+1/\dimension)+0.5};
	\pgfmathsetmacro\LowSobStar{\xstretch*(\LowSob/\xstretch+1/\dimension)};
	\pgfmathsetmacro\UpGradient{\xstretch*(1/2)-0.2};
	\pgfmathsetmacro\UpGradientDual{\xstretch*(1/2)+0.3};
	\pgfmathsetmacro\Half{\xstretch*0.5};
	\pgfmathsetmacro\One{\xstretch*1};
	\coordinate (P00) at (1+\UpGradient,2);
	\coordinate (P01) at (1+\LowSob,2);
	\coordinate (P10) at (1+\UpSob,0);
	\coordinate (P11) at (1+\UpGradientDual,0);
	\coordinate (P12) at (1+\LowSob,0);
	\coordinate (DirExtra) at (1+\UpSobStar,0);
	\coordinate (RegExtra) at (1+\LowSobStar,2);
	
	\draw [thin] (1,2) -- (1+\xlength,2); 
	\draw [thin] (1,0) -- (1+\xlength,0); 
	
	\draw [thick,->] (1,-0.5) -- (1+\xlength+0.2,-0.5);
	\draw [fill=black] (1+\Half,-0.5) circle [radius = .5pt];
	\node [below] at (1+\Half,-0.5) {$\frac{1}{2 \fracspace}$};
	\draw [fill=black] (1+\UpGradient,-0.5) circle [radius = .5pt];
	\node [below] at (1+\UpGradient,-0.5) {$\frac{1}{q_+^L\fracspace}$};
	\draw [fill=black] (1+\UpGradientDual,-0.5) circle [radius = .5pt];
	\node [below] at (1+\UpGradientDual,-0.5) {$\frac{1}{(q_+^{L^\sharp})'\fracspace}$};
	\draw [fill=black] (1+\UpSob,-0.5) circle [radius = .5pt];
	\node [below] at (1+\UpSob,-0.5) {$\frac{1}{p_+^L\fracspace}$};
	\draw [fill=black] (1+\LowSob,-0.5) circle [radius = .5pt];
	\node [below] at (1+\LowSob,-0.5) {$\frac{1}{p_-^L\fracspace}$};
	\draw [fill=black] (1+\xlength,-0.5) circle [radius = .5pt];
	\node [below] at (1+\xlength,-0.5) {$\frac{n+1}{n \fracspace}$};
	\node [right] at (1+\xlength+0.2,-0.5) {$\frac{1}{p \fracspace}$};
	\draw [fill=black] (1,-0.5) circle [radius = .5pt];
	\node [below] at (1,-0.5) {$0$};
	
	\draw [thick,->] (0.7,0) -- (0.7,2.2);
	\node [above] at (0.7,2.2) {$s$};
	\draw [fill=black] (0.7,2) circle [radius = .5pt];
	\node [left] at (0.7,2) {$0$};
	\draw [fill=black] (0.7,0) circle [radius = .5pt];
	\node [left] at (0.7,0) {$-1$};
	
	\draw [thin,dotted] (1,2) -- (1,0); 
	\draw [thin,dotted] (1+\LowSob,2) -- (1+\LowSob,0); 
	\draw [thin,dotted] (1+\UpGradient,2) -- (1+\UpGradient,0);
	\draw [thin,dotted] (1+\xlength,2) -- (1+\xlength,0); 
	
	\path [fill=lightgray, opacity = 0.6] (P00)--(P01)--(P11)--(P10)--(P00);
	\draw [ultra thick, gray] (P00) -- (P01); 
	\draw [ultra thick, gray] (P10) -- (P11); 
	\draw [fill=white] (P00) circle [radius = .75pt];
	\draw [fill=white] (P01) circle [radius = .75pt];
	\draw [fill=white] (P11) circle [radius = .75pt];
	\draw [fill=white] (P10) circle [radius = .75pt];
\end{tikzpicture}
\end{center}
\caption{In the interior of the gray region $\IX_{DB}^{s,p} = \IX_D^{s,p}$ holds (up to equivalent $p$-quasinorms). By $\heartsuit$-duality~\cite[Cor.~5.14]{AA} this is equivalent to $\IX^{-s-1,p'}_{DB^*} = \IX^{-s-1,p'}_D$.}
\label{fig: diagram DB}
\end{figure}

In the particular case $p_-(L^\sharp)< 1$ we have $p_+(L) = \infty$ by duality and similarity and hence the left lower vertex of the identification regions is situated at the origin. However, results can be improved further as follows. We reproduce the argument from \cite[Sec.~7.2.1]{AA} for the sake of clarity.

\begin{prop}
\label{prop: blue region large}
If $p_-(L^\sharp)< 1$, then identification $\IX_{DB}^{s,p} = \IX_D^{s,p}$ and $\IX^{s,p}_L = \Xdot^{s,p} \cap \L^2$ hold in the interior of the extended gray regions of Figure~\ref{fig: diagram DB large} and Figure~\ref{fig: diagram general proof2}, respectively.
\end{prop}

\begin{figure}[ht]
	\begin{center}
		\begin{tikzpicture}[scale=2.4]
			
			\newcommand\fracspace{\vphantom{\frac{1}{1}}};
			\def\dimension{6};
			\def\xlength{3.5};
			\pgfmathsetmacro\xstretch{\xlength/(1+1/\dimension)}; 
			
			\pgfmathsetmacro\UpSob{\xstretch*(1/2-1/\dimension)-0.2};
			\pgfmathsetmacro\UpSobStar{\xstretch*(\UpSob/\xstretch-1/\dimension)};
			\pgfmathsetmacro\LowSob{\xstretch*(1/2+1/\dimension)+0.5};
			\pgfmathsetmacro\LowSobStar{\xstretch*(\LowSob/\xstretch+1/\dimension)};
			\pgfmathsetmacro\UpGradient{\xstretch*(1/2)-0.2};
			\pgfmathsetmacro\UpGradientDual{\xstretch*(1/2)+0.3};
			\pgfmathsetmacro\Half{\xstretch*0.5};
			\pgfmathsetmacro\One{\xstretch*1};
			\coordinate (P00) at (1+\UpGradient,2);
			\coordinate (P01) at (1+\LowSob,2);
			\coordinate (P10) at (1,0);
			\coordinate (P11) at (1+\UpGradientDual,0);
			\coordinate (P12) at (1+\LowSob,0);
			\coordinate (HolderAlex) at (1, 0.3);
			\coordinate (RegExtra) at (1+\LowSobStar,2);
			
			\draw [thin] (1,2) -- (1+\xlength,2); 
			\draw [thin] (1,0) -- (1+\xlength,0); 
			
			\draw [thick,->] (1,-0.5) -- (1+\xlength+0.2,-0.5);
			\draw [fill=black] (1+\Half,-0.5) circle [radius = .5pt];
			\node [below] at (1+\Half,-0.5) {$\frac{1}{2 \fracspace}$};
			\draw [fill=black] (1+\UpGradient,-0.5) circle [radius = .5pt];
			\node [below] at (1+\UpGradient,-0.5) {$\frac{1}{q_+^L\fracspace}$};
			\draw [fill=black] (1+\UpGradientDual,-0.5) circle [radius = .5pt];
			\node [below] at (1+\UpGradientDual,-0.5) {$\frac{1}{(q_+^{L^\sharp})'\fracspace}$};
			\draw [fill=black] (1+\LowSob,-0.5) circle [radius = .5pt];
			\node [below] at (1+\LowSob,-0.5) {$\frac{1}{p_-^L\fracspace}$};
			\draw [fill=black] (1+\xlength,-0.5) circle [radius = .5pt];
			\node [below] at (1+\xlength,-0.5) {$\frac{n+1}{n \fracspace}$};
			\node [right] at (1+\xlength+0.2,-0.5) {$\frac{1}{p \fracspace}$};
			\draw [fill=black] (1,-0.5) circle [radius = .5pt];
			\node [below] at (1,-0.5) {$0=\frac{1}{p_+^L\fracspace}$};
			
			\draw [thick,->] (0.7,0) -- (0.7,2.2);
			\node [above] at (0.7,2.2) {$s$};
			\draw [fill=black] (0.7,2) circle [radius = .5pt];
			\node [left] at (0.7,2) {$0$};
			\draw [fill=black] (0.7,0) circle [radius = .5pt];
			\node [left] at (0.7,0) {$-1$};
			\node [above] at (HolderAlex) {$x_{A^*}^\heartsuit$};
			
			\draw [thin,dotted] (1,2) -- (1,0); 
			\draw [thin,dotted] (1+\LowSob,2) -- (1+\LowSob,0); 
			\draw [thin,dotted] (1+\UpGradient,2) -- (1+\UpGradient,0);
			\draw [thin,dotted] (1+\xlength,2) -- (1+\xlength,0); 
			
			\path [fill=lightgray, opacity = 0.6] (P00)--(P01)--(P11)--(P10)--(HolderAlex)--(P00);
			\draw [ultra thick, gray] (P00) -- (P01); 
			\draw [ultra thick, gray] (P10) -- (P11); 
			\draw [fill=white] (P00) circle [radius = .75pt];
			\draw [fill=white] (P01) circle [radius = .75pt];
			\draw [fill=white] (P11) circle [radius = .75pt];
			\draw [fill=white] (P10) circle [radius = .75pt];
			\draw [fill=white] (HolderAlex) circle [radius = .75pt];
		\end{tikzpicture}
	\end{center}
	\caption{Extension of Figure~\ref{fig: diagram DB} to the left in the case $p_-(L^\sharp)<1$. The extension only concerns exponents with $p \geq q_+(L)>2$. The length of the vertical segment on the left is at most $\nicefrac{n}{p_-(L^\sharp)}-n$.}
	\label{fig: diagram DB large}
\end{figure}
\begin{figure}[ht]
	\begin{center}
		\begin{tikzpicture}[scale=2.4]
			
			\newcommand\fracspace{\vphantom{\frac{1}{1}}};
			\def\dimension{6};
			\def\xlength{3.5};
			\pgfmathsetmacro\xstretch{\xlength/(1+1/\dimension)}; 
			
			\pgfmathsetmacro\UpSob{\xstretch*(1/2-1/\dimension)-0.2};
			\pgfmathsetmacro\UpSobStar{\xstretch*(\UpSob/\xstretch-1/\dimension)};
			\pgfmathsetmacro\LowSob{\xstretch*(1/2+1/\dimension)+0.5};
			\pgfmathsetmacro\LowSobStar{\xstretch*(\LowSob/\xstretch+1/\dimension)};
			\pgfmathsetmacro\UpGradient{\xstretch*(1/2)-0.2};
			\pgfmathsetmacro\UpGradientDual{\xstretch*(1/2)+0.3};
			\pgfmathsetmacro\Half{\xstretch*0.5};
			\pgfmathsetmacro\One{\xstretch*1};
			\coordinate (P00) at (1+\UpGradient,2);
			\coordinate (P01) at (1+\LowSob,2);
			\coordinate (P10) at (1,0);
			\coordinate (P11) at (1+\UpGradientDual,0);
			\coordinate (P12) at (1+\LowSob,0);
			\coordinate (HolderAlex) at (1, 0.3);
			\coordinate (RegExtra) at (1+\LowSobStar,2);
			
			\draw [thin] (1,2) -- (1+\xlength,2); 
			\draw [thin] (1,0) -- (1+\xlength,0); 
			
			\draw [thick,->] (1,-0.5) -- (1+\xlength+0.2,-0.5);
			\draw [fill=black] (1+\Half,-0.5) circle [radius = .5pt];
			\node [below] at (1+\Half,-0.5) {$\frac{1}{2 \fracspace}$};
			\draw [fill=black] (1+\UpGradient,-0.5) circle [radius = .5pt];
			\node [below] at (1+\UpGradient,-0.5) {$\frac{1}{q_+^L\fracspace}$};
			\draw [fill=black] (1+\UpGradientDual,-0.5) circle [radius = .5pt];
			\node [below] at (1+\UpGradientDual,-0.5) {$\frac{1}{(q_+^{L^\sharp})'\fracspace}$};
			\node [below] at (1+\LowSobStar,-0.5) {$\frac{1}{{(p_-^L)_* \vee 1_*\fracspace}}$};
			\draw [fill=black] (1+\LowSob,-0.5) circle [radius = .5pt];
			\node [below] at (1+\LowSob,-0.5) {$\frac{1}{p_-^L\fracspace \vee 1}$};
			\draw [fill=black] (1+\xlength,-0.5) circle [radius = .5pt];
			\node [below] at (1+\xlength,-0.5) {$\frac{n+1}{n \fracspace}$};
			\node [right] at (1+\xlength+0.2,-0.5) {$\frac{1}{p \fracspace}$};
			\draw [fill=black] (1,-0.5) circle [radius = .5pt];
			\node [below] at (1,-0.5) {$0=\frac{1}{p_+^L\fracspace}$};
			
			\draw [thick,->] (0.7,0) -- (0.7,2.2);
			\node [above] at (0.7,2.2) {$s$};
			\draw [fill=black] (0.7,2) circle [radius = .5pt];
			\node [left] at (0.7,2) {$1$};
			\draw [fill=black] (0.7,0) circle [radius = .5pt];
			\node [left] at (0.7,0) {$0$};
			
			\draw [thin,dotted] (1,2) -- (1,0); 
			\draw [thin,dotted] (1+\LowSob,2) -- (1+\LowSob,0); 
			\draw [thin,dotted] (1+\LowSobStar,2) -- (1+\LowSobStar,0); 
			\draw [thin,dotted] (1+\UpGradient,2) -- (1+\UpGradient,0);
			\draw [thin,dotted] (1+\xlength,2) -- (1+\xlength,0); 
			
			\path [fill=lightgray, opacity = 0.6] (P00)--(P01)--(P11)--(P10)--(HolderAlex)--(P00);
			\path [fill=blue!80!black, opacity = 0.4] (P11)--(P12)--(RegExtra)--(P01)--(P11);
			\draw [ultra thick, gray] (P00) -- (P01); 
			\draw [ultra thick, blue!80!black] (P01) -- (RegExtra);
			\draw [ultra thick, gray] (P10) -- (P11); 
			\draw [ultra thick, blue!80!black] (P11) -- (P12); 
			\draw [fill=white] (P00) circle [radius = .75pt];
			\draw [fill=blue!80!black] (P01) circle [radius = .75pt];
			\draw [fill=blue!80!black] (P11) circle [radius = .75pt];
			\draw [fill=white] (P12) circle [radius = .75pt];
			\draw [fill=white] (P10) circle [radius = .75pt];
			\draw [fill=white] (RegExtra) circle [radius = .75pt];
			\draw [fill=white] (HolderAlex) circle [radius = .75pt];
			\node [above] at (1.1,0.3) {$x_{A^*}^\heartsuit \!+\!\begin{bsmallmatrix} 0 \\ 1 \end{bsmallmatrix}$};
		\end{tikzpicture}
	\end{center}
	\caption{Extension of Figure~\ref{fig: diagram general proof} to the left in the case $p_-(L^\sharp)<1$. The extension only concerns exponents with $p \geq q_+(L)>2$. The length of the vertical segment on the left is at most $\nicefrac{n}{p_-(L^\sharp)}-n$.}
	\label{fig: diagram general proof2}
\end{figure}

\begin{proof} 
It suffices to argue for Figure~\ref{fig: diagram DB large} since the extension in Figure~\ref{fig: diagram general proof2} follows from Lemma~\ref{lem: intersection of identification regions} as before.

Consider the analog of Figure~\ref{fig: diagram DB} but for $B^*$. Since we assume $p_-(L^\sharp) < 1$, the right-hand segment of the gray trapezoid described by
\begin{align*}
 \frac{1}{p} = \frac{-s}{q_+(L)'} + \frac{s+1}{p_-(L^\sharp)}
\end{align*}
intersects the vertical line $\nicefrac{1}{p} = 1$ at a point that is called $x_{A^*}$ in~\cite{AA}. Let $x_{A^*}^\heartsuit$ be the symmetric point with respect to $(\nicefrac{1}{2}, - \nicefrac{1}{2})$. By  $\heartsuit$-duality this is a boundary point of the identification region for $\IX_{DB}^{s,p} = \IX_D^{s,p}$. Interpolation with the exponents that have already been obtained in Figure~\ref{fig: diagram DB} yields the extension that is displayed in Figure~\ref{fig: diagram DB large}.

The length of the vertical segment that we have been able to add on the line $\nicefrac{1}{p} = 0$ is given by $\sigma$, where 
\begin{align*}
	\sigma \bigg(\frac{1}{p_-(L^\sharp)} - \frac{1}{q_+(L)'} \bigg) = \frac{1}{p_-(L^\sharp)} - 1.
\end{align*}
Since Theorem~\ref{thm: standard relation J(L) and N(L)} for $L^\sharp$ yields $p_-(L^\sharp) \leq (q_+(L)')_*$, the left-hand side is bounded from below by $\nicefrac{\sigma}{n}$ and we obtain 
\begin{align*}
	\sigma \leq \frac{n}{p_-(L^\sharp)} - n
\end{align*}
as we have claimed.
\end{proof}

Let us illustrate  these diagrams in special cases. When $m=1$, $n\ge 3$ and $d$ is real-valued, we know that $p_-(L)=q_-(L)<1$ and $p_+(L)=\infty$ (Remark~\ref{rem:Dirichlet  property}). Thus, we are in the case of Figure~\ref{fig: diagram general proof2} for the blue and gray identification regions. This is also the generic situation in dimension $n=2$ for any $L$ (Proposition~\ref{prop: J(L) contains neighborhood of Sobolev conjugates}).

In dimension $n=1$, Proposition~\ref{prop: J(L) contains neighborhood of Sobolev conjugates} yields $p_-(L)=q_-(L)=1_* (= \nicefrac{1}{2})$ and $p_+(L)=q_+(L)=\infty$. {The same holds for $L^\sharp$ in place of $L$ and therefore $x_{A^*}^\heartsuit = [0,0]^\top$}. Consequently, we already have the largest possible gray region shown in Figure~\ref{fig: diagram n=1} and there is no additional blue region. In any dimension, the same situation occurs for operators of type $-a^{-1}\Delta_x$ (Corollary~\ref{cor: McIntosh-Nahmod}) or more generally when $d$ depends only on one coordinate (Remark~\ref{rem:Dirichlet  property}).

\begin{figure}[h]
\begin{center}
		\begin{tikzpicture}[scale=2.4]
	
	\newcommand\fracspace{\vphantom{\frac{1}{1}}};
	\def\dimension{6};
	\def\xlength{3.5};
	\pgfmathsetmacro\xstretch{\xlength/(1+1/\dimension)}; 
	
	\pgfmathsetmacro\UpSob{\xstretch*(1/2-1/\dimension)-0.2};
	\pgfmathsetmacro\UpSobStar{\xstretch*(\UpSob/\xstretch-1/\dimension)};
	\pgfmathsetmacro\LowSob{\xstretch*(1/2+1/\dimension)+0.5};
	\pgfmathsetmacro\LowSobStar{\xstretch*(\LowSob/\xstretch+1/\dimension)};
	\pgfmathsetmacro\UpGradient{\xstretch*(1/2)-0.2};
	\pgfmathsetmacro\UpGradientDual{\xstretch*(1/2)+0.3};
	\pgfmathsetmacro\Half{\xstretch*0.5};
	\pgfmathsetmacro\One{\xstretch*1};
	\coordinate (P00) at (2,2);
	\coordinate (P10) at (1,0);
	\coordinate (HolderAlex) at (1, 2);
	\coordinate (RegExtra) at (1+\xlength,2);
	\coordinate (P11) at (1+\LowSob,0);
	
	\draw [thin] (1,2) -- (1+\xlength,2); 
	\draw [thin] (1,0) -- (1+\xlength,0); 
	
	\draw [thick,->] (1,-0.5) -- (1+\xlength+0.2,-0.5);
	\draw [fill=black] (1+\Half,-0.5) circle [radius = .5pt];
	\node [below] at (1+\Half,-0.5) {$\frac{1}{2 \fracspace}$};
	\draw [fill=black] (1+\LowSob,-0.5) circle [radius = .5pt];
	\node [below] at (1+\LowSob,-0.5) {$1$};
	\draw [fill=black] (1+\xlength,-0.5) circle [radius = .5pt];
	\node [below] at (1+\xlength,-0.5) {$\frac{n+1}{n \fracspace}$};
	\node [right] at (1+\xlength+0.2,-0.5) {$\frac{1}{p \fracspace}$};
	\draw [fill=black] (1,-0.5) circle [radius = .5pt];
	\node [below] at (1,-0.5) {$0$};
	
	\draw [thick,->] (0.7,0) -- (0.7,2.2);
	\node [above] at (0.7,2.2) {$s$};
	\draw [fill=black] (0.7,2) circle [radius = .5pt];
	\node [left] at (0.7,2) {$1$};
	\draw [fill=black] (0.7,0) circle [radius = .5pt];
	\node [left] at (0.7,0) {$0$};
	
	\draw [thin,dotted] (1,2) -- (1,0); 
	\draw [thin,dotted] (1+\LowSob,2) -- (1+\LowSob,0); 
	\draw [thin,dotted] (1+\xlength,2) -- (1+\xlength,0); 
	
	\path [fill=lightgray, opacity = 0.6] (HolderAlex) -- (RegExtra) -- (P11)--(P10);
	\draw [ultra thick, gray] (HolderAlex) -- (RegExtra);
	\draw [ultra thick, gray] (P10) -- (P11);
	\draw [fill=white] (P11) circle [radius = .75pt];
	\draw [fill=white] (P10) circle [radius = .75pt];
	\draw [fill=white] (RegExtra) circle [radius = .75pt];
	\draw [fill=white] (HolderAlex) circle [radius = .75pt];
\end{tikzpicture}
\end{center}
\caption{Figure~\ref{fig: diagram general proof2} in dimension $n=1$ and in any dimension for the special case $L= -a^{-1} \Delta_x$ or more generally when $d$ depends only on one coordinate.}
\label{fig: diagram n=1}
\end{figure}
\subsection{Solvability for fractional regularity data}
\label{subsec: fractional solvability}

We turn to solvability of the Dirichlet problems $(D)_{\Hdot^{s,p}}^\Le$ and $(D)_{\Bdot^{s,p}}^\Le$ when $0<s<1$ and $0 < p \leq \infty$ satisfy $\nicefrac{1}{p} < 1 + \nicefrac{s}{n}$. The restrictions on $s$ and $p$ guarantee that all distributions in $\Hdot^{s,p}$ and $\Bdot^{s,p}$ are locally integrable functions. Indeed, for $p=\infty$ we have $\H^{s,\infty} \subseteq \B^{s,\infty} = \Lamdot^s$, whereas for $p<\infty$ both are interpolation spaces between $\Hdot^{0,p_0} = \L^{p_0}$ and $\Hdot^{1,p_1}$ for some exponents $p_0>1, p_1 > 1_*$. 

In the formulation of the Dirichlet problems for fractional regularity data we consider the data spaces as classes of measurable functions and do not factor out constants. We use  the pair $(\Y, \Xdot)$ to denote either $(\Z, \Bdot)$ or $(\T, \Hdot)$. By definition of tent and $\Z$-spaces, all problems that appear in Section~\ref{intro:fractionalDP} can simultaneously be phrased as asking for given $f \in \Xdot^{s,p}$ to find a solution  to
\begin{equation*}
	(D)_{\Xdot^{s,p}}^\Le  \quad\quad
	\begin{cases}
		\Le u=0   & (\text{in } \reu), \\
		\nabla u \in \Y^{s-1,p},   \\
		\lim_{t \to 0} \bariint_{W(t,x)} |u(s,y)-f(x)| \, \d s \d y = 0 & (\text{a.e. } x\in \R^n).
			\end{cases}
\end{equation*}
Let us mention that another way of formulating the boundary condition is $\lim_{t \to 0} u(t,\cdot) = f$ in  $\cD'(\R^n)/\IC^m$, see \cite{BM, AA}. In all cases, we recover this condition in the construction of our solutions. We do not impose a condition at $t=\infty$, contrarily to \cite{AA}. 

\begin{rem}
\label{rem: energy solution}
For $(\nicefrac{1}{p},s) = (\nicefrac{1}{2},\nicefrac{1}{2})$ we obtain $\Xdot^{\nicefrac{1}{2},2} = \Hdot^{\nicefrac{1}{2},2}$ by Fubini's theorem and $\Y^{-\nicefrac{1}{2},2} = \L^2$ by the averaging trick, so that $(D)_{\Xdot^{\nicefrac{1}{2},2}}^\Le$ is a Dirichlet problem for the energy class. The energy solution given by Proposition~\ref{prop: existence of energy solution} is (modulo a constant) a solution to this problem. Indeed, consider $f\in \Hdot^{\nicefrac{1}{2},2}$ and let $u$ be  the energy solution. It converges to $f$ as $t \to 0$ in $\Xdot^{\nicefrac{1}{2},2}$. By Proposition~\ref{prop: NT trace Y}, there exists a non-tangential trace $u_0$ and the Ces\`aro means of $u(t, \cdot)$ converge in $\cD'$ to $u_0$ as $t \to 0$. It follows that $f = u_0 +c$ for some $c \in \IC^m$. From now on, we call $u + c$ the \emph{energy solution with Dirichlet datum $f$.}\index{energy solution!with Dirichlet datum} 
\end{rem}

\emph{Solvability}\index{solvability} of $(D)_{\Xdot^{s,p}}^\Le$ means that for any given data there exists a solution. \emph{Compatible solvability}\index{solvability!compatible} means that the energy solution is a solution if the data is also in $\Hdot^{\nicefrac{1}{2},2}$. This notion of (compatible) solvability differs from parts of the literature in that we do not require an \emph{a priori} estimate for solutions by the data, compare with \cite[Section~2.4]{BM}. Such estimate usually holds since a specific method was used to construct solutions. We find it natural to separate these two aspects of solvability theory by using the concept of solution operators. This notion is manufactured in a way that is amenable to interpolation, independently of any uniqueness result.

\begin{defn}
\label{defn: solution operator}
Let $s \in (0,1)$ and $p \in (1_*,\infty]$ satisfy $\nicefrac{1}{p} < 1 + \nicefrac{s}{n}$. Consider $\Xdot^{s,p}$ as a (quasi-)Banach space modulo constants. A \emph{solution operator}\index{solution!operator for $(D)_{\Xdot^{s,p}}^\Le$} for $(D)_{\Xdot^{s,p}}^\Le$ is a linear map $\sol: \Xdot^{s,p} \to \cD'(\reu)/\IC^m$ such that for all $f \in  \Xdot^{s,p}$ the function $u \coloneqq \sol f$ satisfies
\begin{align}
\label{eq: solution operator}
\begin{cases}
	\Le u=0   & (\text{in } \reu), \\
	\|\nabla u\|_{\Y^{s-1,p}} \lesssim \|f\|_{\Xdot^{s,p}},   \\
	\lim_{t \to 0} u(t,\cdot) = f & (\text{in } \cD'(\R^n)/\IC^m),
\end{cases}
\end{align}
where the implicit constant in the second line is independent of $f$.
The solution operator is \emph{compatible} if it agrees on $\Xdot^{s,p} \cap \Hdot^{\nicefrac{1}{2},2}$ with the solution operator for the energy class (Proposition~\ref{prop: existence of energy solution}).
\end{defn}

Recall that a weak solution of $\Le u=0$ in $\reu$ is in $\Wloc^{1,2}$ by definition and of class $\C^\infty((0,\infty); \Lloc^2)$ by Corollary~\ref{cor: weak solution smooth in t}. Hence, all conditions in our definition make sense. The second line implies that $\sol: \Xdot^{s,p} \to \cD'(\reu)/\IC^m$ is continuous. In passing, we note that in the existence parts of both Theorem~\ref{thm: blockdir} (Section~\ref{subsec: existence for Dp and Rp}) and Theorem~\ref{thm: Holder-dir} (Section~\ref{sec: existence of Holder-dir}) we have already encountered such operators for different classes of data without using the terminology. Proposition~\ref{prop: existence of energy solution} provides a solution operator for $(D)_{\Hdot^{\nicefrac{1}{2},2}}^\Le$.

\begin{lem}
\label{lem: solution operator implies solvability}
Let $s \in (0,1)$ and $p \in (1_*,\infty]$ satisfy $\nicefrac{1}{p} < 1 + \nicefrac{s}{n}$. If there is a (compatible) solution operator for $(D)_{\Xdot^{s,p}}^\Le$, then $(D)_{\Xdot^{s,p}}^\Le$ is (compatibly) solvable.
\end{lem}

\begin{proof}
For solvability we do not consider $\Xdot^{s,p}$ modulo constants. Given $f \in \Xdot^{s,p}$, the assumption yields a solution $u$ to \eqref{eq: solution operator}. Now, $u$ has a non-tangential trace $u_0$ and the Ces\`aro means of $u(t,\cdot)$ converge to $u_0$ in $\cD'$ as $t\to 0$, see  Proposition~\ref{prop: NT trace Y}. Thus, $f = u_0 +c$ for some $c \in \IC^m$ and $u+c$ is a solution to $(D)_{\Xdot^{s,p}}^\Le$ with data $f$. 

If the solution operator is compatible and $f$ also belongs to $\Hdot^{\nicefrac{1}{2},2}$, then $u+c$ is the energy solution, see Remark~\ref{rem: energy solution}.
\end{proof}

We shall now construct solution operators in a series of results, enlarging the range of boundary spaces step by step.

We begin with exponents in the blue and gray identification regions from the previous section.  Note that the $\H^p$ regularity problem $(R)_p^\Le$ does not fit into the scheme of problems $(D)_{\Xdot^{s,p}}^\Le$ because of the missing square function control for $\nabla u$. Hence, no interpolation argument between the existence parts of Theorem~\ref{thm: blockdir} and \ref{thm: blockreg} can help us here. Instead, we rely on the first-order theory and adapted Hardy spaces as in Section~\ref{sec: estimates semigroup}.

\begin{prop}
\label{prop: solvability identification region}
Suppose that $(\nicefrac{1}{p}, s)$ is contained in the interior of the colored region described in Figure~\ref{fig: diagram general proof} and Figure~\ref{fig: diagram general proof2} in the particular case $p_-(L^\sharp) < 1$. Then $(D)_{\Xdot^{s,p}}^\Le$ is solvable. There is a compatible solution operator that assigns to each $f \in \Xdot^{s,p}$ a solution of class $\C_0([0, \infty); \Xdot^{s,p}) \cap \C^\infty((0,\infty); \Xdot^{s,p})$ with $u(0,\cdot) = f$ and comparability
	\begin{align*}
		\sup_{t > 0} \|u(t,\cdot)\|_{\Xdot^{s,p}} \simeq \|f\|_{\Xdot^{s,p}} \simeq \|\nabla u\|_{\Y^{s-1,p}}.
	\end{align*}
\end{prop}

\begin{proof}
In view of Lemma~\ref{lem: solution operator implies solvability} it suffices to construct the solution operator. We first consider $f \in \Xdot^{s,p} \cap \W^{1,2}$. In this case, we set of course $u(t,x) \coloneqq \e^{-t L^{1/2}}f(x)$.

\medskip

\noindent \emph{Step 1: Regularity and the first comparability}. Since we have $\IX_L^{s,p} = \Xdot^{s,p} \cap \L^2$ with equivalent $p$-quasinorms, the regularity for $u$ and the first comparability immediately follow from the bounded $\H^\infty$-calculus and the semigroup properties on $\IX_L^{s,p}$, see Section~\ref{subsec: Hardy abstract sectorial}. This argument also yields quantitative bounds for $\|t^{k/2} \partial_t u(t,\cdot)\|_{\Xdot^{s,p}}$ that will be needed to carry the $\C^\infty$-property over to general data $f \in \Xdot^{s,p}$ in Step~4.

\medskip

\noindent \emph{Step~2: The second comparability when $p \leq 2$}. By means of the intertwining property we find
\begin{align*}
	\|\nabla u \|_{\Y^{s-1,p}} 
	&\simeq \|L^{1/2} \e^{-tL^{1/2}} f\|_{\Y^{s-1,p}}  + \| \nabla_x \e^{-tL^{1/2}} f\|_{\Y^{s-1,p}} \\
	&\simeq \|t L^{1/2} \e^{-tL^{1/2}} f\|_{\Y^{s,p}}  + \| \e^{-t\tM^{1/2}} \nabla_x f\|_{\Y^{s-1,p}} \\
	&\eqqcolon \|\phi(t^2 L)f \|_{\Y^{s,p}} + \|\psi(t^2 \tM) \nabla_x f \|_{\Y^{s-1,p}}
\intertext{and the auxiliary functions are of class $\phi \in \Psi_{1/2}^\infty$ and $\psi \in \Psi_0^\infty$. They are admissible for defining $\IX_L^{s,p}$ and $\IX_{\tM}^{s-1,p}$, respectively, since we have $p\leq 2$ and $s<1$. Hence, we can continue with}
	&\simeq \|f \|_{\IX_L^{s,p}} + \|\nabla_x f \|_{\IX_{\tM}^{s-1,p}}\\
	&\simeq \|f\|_{\IX_L^{s,p}}\\
	&\simeq \|f\|_{\Xdot^{s,p}},
\end{align*}
where we used Figure~\ref{fig: diagram with X} in the second step.

\medskip

\noindent \emph{Step~3: The second comparability when $p>2$}. In this case we are in the gray identification region. We know from Figure~\ref{fig: diagram DB} (or Figure~\ref{fig: diagram DB large}) that we can identify $\IX_{DB}^{s-1,p} = \IX_D^{s-1,p}$ and therefore the Cauchy characterization of adapted spaces in \cite[Thm.~5.26]{AA} and \cite[Rem.~5.28]{AA} yields
\begin{align}
\label{eq1: solvability identification region}
	\|\e^{-t[DB]} \ind_{\IC^+}(DB) g\|_{\Y^{s-1,p}} \simeq \|g\|_{\IX_D^{s-1,p}} \quad (g \in \cl{\ran(DB)}).
\end{align}
We pick 
\begin{align*}
	g 
	\coloneqq  \begin{bmatrix} 0 \\ \nabla_x f \end{bmatrix} 
	= DB \begin{bmatrix} - af  \\ 0 \end{bmatrix} .
\end{align*}
As for the right-hand side in \eqref{eq1: solvability identification region}, Figure~\ref{fig: identification D-adapted} yields $\|g\|_{\IX_D^{s-1,p}}  \simeq \|f\|_{\Xdot^{s,p}}$. Next, we use the identity $2 (\ind_{\IC^+}(z)) = 1 + \nicefrac{\sqrt{z^2}}{z}$ to write
\begin{align*}
		2 (\ind_{\IC^+}(DB) g)
		= \begin{bmatrix} 0 \vphantom{\tL} \\ \nabla_x f \vphantom{\tM} \end{bmatrix} 
		+ [DB] \begin{bmatrix} - af \vphantom{\tL} \\ 0 \vphantom{\tM} \end{bmatrix}
		= \begin{bmatrix}  -\tL^{1/2} af \vphantom{\tL} \\ \nabla_x f \vphantom{\tM} \end{bmatrix}.
\end{align*}
The intertwining relation and the similarity of $L$ and $\tL$ lead to
\begin{align*}
	2\e^{-t[DB]} \ind_{\IC^+}(DB) g 
	= \begin{bmatrix} - \e^{-t \tL^{1/2}} \tL^{1/2} af  \vphantom{\tL} \\ \e^{- t \tM^{1/2}} \nabla_x \vphantom{\tM} \end{bmatrix}
	= \begin{bmatrix} - a L^{1/2} \e^{-t L^{1/2}} f  \vphantom{\tL} \\ \nabla_x \e^{- t L^{1/2}} f \vphantom{\tM} \end{bmatrix}
	= \begin{bmatrix} a \partial_t u  \vphantom{\tL} \\ \nabla_x u \vphantom{\tM} \end{bmatrix}.
\end{align*}
Thus, the left-hand side in \eqref{eq1: solvability identification region} is comparable to $\|\nabla u\|_{\Y^{s-1,p}}$.

\medskip

\noindent \emph{Step~4: Extension to a solution operator.} By the same density argument as for the regularity problem in Section~\ref{subsec: existence for Dp and Rp} when $p\geq n$, we can construct for general $f \in \Xdot^{s,p}$ a weak solution to $\Le u = 0$ in $\reu$ that has all the properties stated in the proposition. The construction depends linearly on the data and since $u(t,\cdot) \to f$ in $\Xdot^{s,p} \subseteq \cD'/\IC^m$, we see that $u$ solves \eqref{eq: solution operator}. This means that we have constructed a compatible solution operator.
\end{proof}

If $p_+(L) < \infty$, then (the existence part of) Theorem~\ref{thm: blockdir} contains existence of a solution to the Dirichlet problem $(D)_p^\Le$ in a range of exponents that exceeds the identification region for $\IH_L^{0,p}$ by up to one Sobolev conjugate. This leads to the following improvement of the previous result in that case.\index{Dirichlet problem!solvability for $\Hdot^{s,p}$/$\Bdot^{s,p}$-data}

\begin{figure}[h]
\begin{center}
\begin{tikzpicture}[scale=2.4]

	\newcommand\fracspace{\vphantom{\frac{1}{1}}};
	\def\dimension{6};
	\def\xlength{3.5};
	\pgfmathsetmacro\xstretch{\xlength/(1+1/\dimension)}; 
	
	\pgfmathsetmacro\UpSob{\xstretch*(1/2-1/\dimension)-0.2};
	\pgfmathsetmacro\UpSobStar{\xstretch*(\UpSob/\xstretch-1/\dimension)};
	\pgfmathsetmacro\LowSob{\xstretch*(1/2+1/\dimension)+0.5};
	\pgfmathsetmacro\LowSobStar{\xstretch*(\LowSob/\xstretch+1/\dimension)};
	\pgfmathsetmacro\UpGradient{\xstretch*(1/2)-0.2};
	\pgfmathsetmacro\UpGradientDual{\xstretch*(1/2)+0.3};
	\pgfmathsetmacro\Half{\xstretch*0.5};
	\pgfmathsetmacro\One{\xstretch*1};
	\coordinate (P00) at (1+\UpGradient,2);
	\coordinate (P01) at (1+\LowSob,2);
	\coordinate (P10) at (1+\UpSob,0);
	\coordinate (P11) at (1+\UpGradientDual,0);
	\coordinate (P12) at (1+\LowSob,0);
	\coordinate (DirExtra) at (1+\UpSobStar,0);
	\coordinate (RegExtra) at (1+\LowSobStar,2);

	\draw [thin] (1,2) -- (1+\xlength,2); 
	\draw [thin] (1,0) -- (1+\xlength,0); 
	
	\draw [thick,->] (1,-0.5) -- (1+\xlength+0.2,-0.5);
	\draw [fill=black] (1+\Half,-0.5) circle [radius = .5pt];
	\node [below] at (1+\Half,-0.5) {$\frac{1}{2 \fracspace}$};
	\draw [fill=black] (1+\UpGradient,-0.5) circle [radius = .5pt];
	\node [below] at (1+\UpGradient,-0.5) {$\frac{1}{q_+^L\fracspace}$};
	\draw [fill=black] (1+\UpGradientDual,-0.5) circle [radius = .5pt];
	\node [below] at (1+\UpGradientDual,-0.5) {$\frac{1}{(q_+^{L^\sharp})'\fracspace}$};
	\draw [fill=black] (1+\UpSob,-0.5) circle [radius = .5pt];
	\node [below] at (1+\UpSob,-0.5) {$\frac{1}{p_+^L\fracspace}$};
	\draw [fill=black] (1+\LowSob,-0.5) circle [radius = .5pt];
	\node [below] at (1+\LowSob,-0.5) {$\frac{1}{p_-^L \vee 1 \fracspace}$};
	\draw [fill=black] (1+\UpSobStar,-0.5) circle [radius = .5pt];
	\node [below] at (1+\UpSobStar,-0.5) {$\frac{1}{(p_+^L)^*\fracspace}$};
	\draw [fill=black] (1+\LowSobStar,-0.5) circle [radius = .5pt];
	\node [below] at (1+\LowSobStar,-0.5) {$\frac{1}{{(p_-^L)_* \vee 1_*\fracspace}}$};
	\draw [fill=black] (1+\xlength,-0.5) circle [radius = .5pt];
	\node [below] at (1+\xlength,-0.5) {$\frac{n+1}{n \fracspace}$};
	\node [right] at (1+\xlength+0.2,-0.5) {$\frac{1}{p \fracspace}$};
	\draw [fill=black] (1,-0.5) circle [radius = .5pt];
	\node [below] at (1,-0.5) {$0$};
	
	\draw [thick,->] (0.7,0) -- (0.7,2.2);
	\node [above] at (0.7,2.2) {$s$};
	\draw [fill=black] (0.7,2) circle [radius = .5pt];
	\node [left] at (0.7,2) {$1$};
	\draw [fill=black] (0.7,0) circle [radius = .5pt];
	\node [left] at (0.7,0) {$0$};

	\draw [thin,dotted] (1,2) -- (1,0); 
	\draw [thin,dotted] (1+\LowSob,2) -- (1+\LowSob,0); 
	\draw [thin,dotted] (1+\LowSobStar,2) -- (1+\LowSobStar,0); 
	\draw [thin, dotted] (1+\Half, 2) -- (1+\Half,0);
	\draw [thin,dotted] (1+\UpGradient,2) -- (1+\UpGradient,0);
	\draw [thin,dotted] (1+\xlength,2) -- (1+\xlength,0); 
	
	\path [fill=lightgray, opacity = 0.6] (P00)--(P01)--(P11)--(P10)--(P00);
	\path [fill=blue!80!black, opacity = 0.4] (P11)--(P12)--(RegExtra)--(P01)--(P11);
	\path [fill=red!80!black, opacity = 0.4] (DirExtra)--(P10)--(P00)--(DirExtra);
	\draw [ultra thick, gray] (P00) -- (P01); 
	\draw [ultra thick, blue!80!black] (P01) -- (RegExtra); 
	\draw [ultra thick, gray] (P10) -- (P11); 
	\draw [ultra thick, red!80!black] (DirExtra) -- (P10); 
	\draw [ultra thick, blue!80!black] (P11) -- (P12); 
	\draw [fill=white] (P00) circle [radius = .75pt];
	\draw [fill=blue!80!black] (P01) circle [radius = .75pt];
	\draw [fill=red!80!black] (P10) circle [radius = .75pt];
	\draw [fill=blue!80!black] (P11) circle [radius = .75pt];
	\draw [fill=white] (P12) circle [radius = .75pt];
	\draw [fill=white] (DirExtra) circle [radius = .75pt];
	\draw [fill=white] (RegExtra) circle [radius = .75pt];
\end{tikzpicture}
\end{center}
\caption{Extended region for compatible solvability of $(D)_{\Xdot^{s,p}}^\Le$ when $p_+(L) < \infty$. Recall that $p_+(L)^* = \infty$ if $p_+(L) \geq n$.}
\label{fig: diagram general proof3}
\end{figure}

\begin{prop}
\label{prop: solvability fractional}
Suppose that $p_+(L) < \infty$. If $(\nicefrac{1}{p}, s)$ is contained in the interior of the colored region in Figure~\ref{fig: diagram general proof3}, then there is a compatible solution operator for $(D)_{\Xdot^{s,p}}^\Le$. In particular, the problem is compatibly solvable.
\end{prop}

\begin{proof}
The blue and gray regions have been treated in Proposition~\ref{prop: solvability identification region}. We need to add the red triangle to the picture.  It suffices to show for any $P_{0} \coloneqq (\nicefrac{1}{p_0}, 0)$ with $p_+(L) \leq p_0 < p_+(L)^*$ (bottom red segment) and any $P_1 \coloneqq (\nicefrac{1}{p_1}, s_1)$ in the interior of the gray region that a compatible solution operator exists for all points on the open segment $\overline{P_0P_1}$. Compatible solvability then follows by Lemma~\ref{lem: solution operator implies solvability}.

We argue by interpolation and consider the data classes as Banach spaces embedded into $\cD' / \IC^m$. In Section~\ref{subsec: existence for Dp and Rp} we have established existence of a solution with the properties (i) and (iv) of Theorem~\ref{thm: blockdir}. This construction furnishes a continuous linear solution operator $\sol_{0}: \Hdot^{0,p_0} \to \cD'(\reu)/\IC^m$ such that $u =\sol_{0} f$ solves
\begin{equation*}
	\begin{cases}
		\Le u=0   & (\text{in } \reu), \\
		\|\nabla u\|_{\T^{-1,p_0}} \lesssim \|f\|_{\Hdot^{0,p_0}},   \\
		\lim_{t \to 0} u(t,\cdot) = f & (\text{in } \cD'(\R^n)/\IC^m),
	\end{cases}
\end{equation*}
whereas Proposition~\ref{prop: solvability identification region} furnishes a continuous linear solution operator $\sol_{1}: \Hdot^{s_1, p_1} \to \cD'(\reu)/\IC^m$ such that $u =\sol_{1} f$ solves
\begin{equation*}
	\begin{cases}
		\Le u=0   & (\text{in } \reu), \\
		\|\nabla u\|_{\T^{s_1-1,p_1}} \lesssim \|f\|_{\Hdot^{s_1, p_1}},   \\
		\lim_{t \to 0} u(t,\cdot) = f & (\text{in } \cD'(\R^n)/\IC^m).
	\end{cases}
\end{equation*}
Since both operators produce compatible solutions, the universal approximation technique implies that they coincide on $\Hdot^{0, p_0} \cap \Hdot^{s_1, p_1}$. Hence, we have a well-defined continuous linear operator
\begin{align*}
	\sol: \Hdot^{0, p_0} + \Hdot^{s_1, p_1} \to \cD'(\reu)/\IC^m
\end{align*}
such that $u = \sol f$ solves $\Le u = 0$ in $\reu$ and satisfies $u(t,\cdot) \to f$ as $t \to 0$ in $\cD'/\IC^m$. 

Pick any point $(\nicefrac{1}{p},s)$ on the open segment $\overline{P_0P_1}$. Since the real and complex interpolation spaces of an interpolation couple continuously embed into the sum space, we obtain that $\sol: \Xdot^{s,p} \to \cD'(\reu)/\IC^m$ is continuous. The map $\sol$ and the continuous solution map for energy solutions from Proposition~\ref{prop: existence of energy solution} agree on $\Xdot^{s,p} \cap \Hdot^{\nicefrac{1}{2},2} \cap \Hdot^{s_1,p_1}$ and hence on $\Xdot^{s,p} \cap \Hdot^{\nicefrac{1}{2},2}$. Since the maps $\nabla \sol: \Hdot^{0, p_0} \to \T^{-1,p_0}$ and $\nabla \sol: \Hdot^{s_1, p_1} \to \T^{s_1-1,p_0}$ are bounded, we obtain by real and complex interpolation that $\nabla \sol: \Xdot^{s,p} \to \Y^{s-1,p}$ is bounded. This means that we have constructed a solution operator for $(D)_{\Xdot^{s,p}}^\Le$.
\end{proof}

In the case $p_+(L) > n$ we can go one step further and study endpoint problems $(D)_{\Xdot^{\alpha, \infty}}^\Le$ for $0 < \alpha < 1 - \nicefrac{n}{p_+(L)}$. We have $\Bdot^{\alpha,\infty} = \Lamdot^\alpha$ with equivalent norms, so that $(D)_{\Bdot^{\alpha,\infty}}^\Le$ is a third way of posing a Dirichlet problem with H\"older continuous data. The other endpoint problem uses the data space $\Hdot^{\alpha, \infty} = \BMOdot^\alpha$, which is continuously embedded into $\Lamdot^\alpha$ and carries the equivalent norm \eqref{eq: Strichartz BMO norm}. The upshot is that, given $f \in \Xdot^{\alpha,\infty}$, the existence part of Theorem~\ref{thm: Holder-dir} already shows that $u$ defined in \eqref{eq: Holder-dir solution} is a compatible solution that converges to $f$ at the boundary in the non-tangential sense. The following addendum guarantees that this solution also solves the new endpoint problem\index{Dirichlet problem!solvability for $\Lamdot^\alpha/\BMOdot^\alpha$-data} and that \eqref{eq: Holder-dir solution} defines a compatible solution operator to $(D)_{\Xdot^{\alpha, \infty}}^\Le$ .

\begin{prop}
\label{prop: BMO-Besov-dir}
Suppose that $p_+(L) > n$ and that $0 < \alpha < 1 - \nicefrac{n}{p_+(L)}$. Then the Dirichlet problem $(D)_{\Xdot^{\alpha,\infty}}^{\Le}$ is compatibly solvable. More precisely, given $f \in \Xdot^{\alpha,\infty}$, the same solution $u$ that was defined in \eqref{eq: Holder-dir solution} and solves $(D)_{\Lamdot^\alpha}^{\Le}$ and $(\wtD)_{\Lamdot^\alpha}^{\Le}$, also solves $(D)_{\Xdot^{\alpha,\infty}}^{\Le}$ and satisfies
\begin{align*}
 \|\nabla u\|_{\Y^{\alpha -1, \infty}} \simeq \|f\|_{\Xdot^{\alpha,\infty}}.
\end{align*}
\end{prop}

\begin{rem}
\label{rem: BMO-Besov-dir}
Combining Theorem~\ref{thm: Holder-dir} with Proposition~\ref{prop: BMO-Besov-dir} yields comparability
\begin{align*}
 \|\nabla u\|_{\T^{-1, \infty; \alpha}} =	\|C_\alpha(t \nabla u)\|_\infty \simeq \|W(t^{1-\alpha} \nabla u) \|_\infty = 	\|\nabla u \|_{\Z^{\alpha-1,\infty}},
\end{align*}
whenever $u$ is a solution to $(D)_{\Lamdot^\alpha}$. A simple comparison of the two functionals shows that the estimate `\,$\gtrsim$' holds for any $\Lloc^2$-function $F$ in place of $\nabla u$. The converse is a special property of  weak solutions to $(D)_{\Lamdot^\alpha}^\Le$.
\end{rem}

\begin{rem}
\label{rem: BMO-Besov-dir-2}
If $p_-(L^\sharp) < 1$ and $\alpha < n(\nicefrac{1}{p_-(L^\sharp)} -1)$, then $u$ is given by a weak$^*$-continuous semigroup on $\Lamdot^\alpha$ as the dual of $\H^p$, $\alpha = n(\nicefrac{1}{p} -1)$, see Lemma~\ref{lem: identification u dual semigroup}. In essence, this followed from the identification $\IH_{L^\sharp}^p = (a^*)^{-1}(\H^p \cap \L^2)$. By interpolation one can obtain a subregion of the red region where the (unique) solution to $(D)_{\Bdot^{s,p}}^{\Le}$ is given by a $\C_0$-semigroup. 

An analogous result for $\BMOdot^\alpha$ would require boundedness of the Poisson semigroup for $L^\sharp$ on $(a^*)^{-1}(\Hdot^{-\alpha,1} \cap \L^2)$, which we do not know when $\alpha>0$. One can use the first-order approach to obtain the semigroup property of the solution to $(D)_{\Hdot^{\alpha,\infty}}^{\Le}$ for $0<\alpha<\theta$, where $\theta$ appears in Figure~\ref{fig: diagram p-small} or equivalently as the upper endpoint of the vertical boundary segment of the gray region in Figure~\ref{fig: diagram general proof2}. The semigroup property for $\theta \le \alpha < n(\nicefrac{1}{p_-(L^\sharp)} -1)$ is unclear. These observations will not be needed in the further course, so we do not detail them.
\end{rem}

\begin{proof}[Proof of Proposition~\ref{prop: BMO-Besov-dir}]
We fix (a representative for) $f \in \Xdot^{\alpha,\infty}$ and let $u$ be the solution to both $(D)_{\Lamdot^\alpha}^{\Le}$ and $(\wtD)_{\Lamdot^\alpha}^{\Le}$ defined in \eqref{eq: Holder-dir solution}. Since we are working within the same or even a smaller class of boundary data, we have at our disposal all properties for $u$ from Section~\ref{sec: existence of Holder-dir} and only at distinguished places we have to intervene in order to obtain the additional features that we claimed above. More precisely, we have to modify Part~3 for the upper bound of $\|\nabla u\|_{\Y^{\alpha -1, \infty}}$ and Part~5 for the converse.

\medskip

\noindent \emph{Modification of Part~3: The bound `\,$\lesssim$'}. 
In the case $\X = \B$ it suffices to combine the observation from Remark~\ref{rem: BMO-Besov-dir} and the existence part of Theorem~\ref{thm: Holder-dir} in order to obtain
\begin{align*}
	\|\nabla u \|_{\Z^{\alpha-1,\infty}}  \lesssim  \|C_\alpha(t \nabla u)\|_\infty \lesssim \|f\|_{\Lamdot^{\alpha}}.
\end{align*}

We turn to the case $\X = \H$. We have to prove that for all cubes $Q \subseteq \R^n$ of sidelength $\ell$ we have
\begin{align}
\label{eq: Goal Part 3 Holder-dir BMO}
	\bigg(\int_0^\ell \int_Q |s^{1-\alpha} \nabla u|^2 \; \frac{\d y \d s}{s} \bigg)^{1/2} \lesssim |Q| \|f\|_{\BMOdot^\alpha}.
\end{align}
From now on $Q$ is fixed. Since both sides stay the same under adding constants to $u$ and $f$, we can assume $(f)_{Q} = 0$. 

In contrast to Section~\ref{sec: existence of Holder-dir} we use a smooth resolution for $f$ in order to represent $u$. We let $(\eta_j)_j$ be a smooth partition of unity on $\R^n$ subordinate to the sets $D_1 \coloneqq 4Q$ and $D_j \coloneqq 2^{j+1}Q \setminus 2^{j-1}Q$, $j \geq 2$, such that $\|\eta_j\|_\infty + 2^j \ell(Q) \|\nabla_x \eta_j\|_\infty \leq C$ for a dimensional constant $C$. For $j\geq1$ we introduce
\begin{align*}
	f_j \coloneqq \eta_j f, \quad u_j(t,\cdot) \coloneqq \e^{-tL^{1/2}}f_j.
\end{align*}

The main difficulty is to handle the local term for $j=1$. For the moment, let us take for granted the estimate
\begin{align}
\label{eq: Local bound Holder-dir BMO}
  \|f_1 \|_{\Hdot^{\alpha,2}}^2 \lesssim |Q| \|f\|_{\BMOdot^\alpha}^2.
\end{align}
This is where the smoothness of $\eta_1$ is needed and we include the argument at the end. Thus, it suffices to prove the local bound
\begin{align}
 \int_0^\ell \int_Q |s^{1-\alpha}  \nabla u_1|^2 \; \frac{\d y \d s}{s} \lesssim  \|f_1\|_{\Hdot^{\alpha,2}}^2.
\end{align}
In doing so, we can work under the qualitative assumption $f_1 \in \W^{1,2}$, which can be removed afterwards via density of $\W^{1,2}$ in $\Hdot^{\alpha,2} \cap \L^2$ and Fatou's lemma. We use the intertwining property to write
\begin{align*}
	s^{1-\alpha} \nabla u_1(s,y)
	= \begin{bmatrix} -s^{1-\alpha} L^{1/2} \e^{-s L^{1/2}} f_1 \\ s^{1-\alpha} \e^{-s \tM^{1/2}} \nabla_x f_1 \end{bmatrix} 
	\eqqcolon \begin{bmatrix} s^{-\alpha} \phi(s^2 L)  f_1\\ s^{1-\alpha} \psi(s^2 \tM) \nabla_x f_1 \end{bmatrix},
\end{align*}
where $\phi \in \Psi_{1/2}^\infty$ and $\psi \in \Psi_0^\infty$. These auxiliary functions are admissible for $\IH_L^{\alpha,2}$ and $\IH_{\tM}^{\alpha - 1,2}$, respectively. Hence, we get as required
\begin{align*}
\int_0^\ell \int_Q |s^{1-\alpha}  \nabla u_1|^2 \; \frac{\d y \d s}{s}
&\leq \iint_{\reu} \bigg| \begin{bmatrix} s^{-\alpha} \phi(s^2 L) f_1 \\ s^{1-\alpha} \psi(s^2 \tM) \nabla_x f_1 \end{bmatrix} \bigg|^2 \; \frac{\d y \d s}{s} \\
&\simeq \|f_1\|_{\IH_L^{\alpha,2}}^2 + \|\nabla _x f_1\|_{\IH_{\tM}^{\alpha-1,2}}^2 \\
&\simeq \|f_1\|_{\IH_L^{\alpha,2}}^2\\
&\simeq  \|f_1\|_{\Hdot^{\alpha,2}}^2,
\end{align*}
where the third step is due to Figure~\ref{fig: diagram with X} and the final step uses that $(\nicefrac{1}{2},\alpha)$ belongs to the identification region of Figure~\ref{fig: diagram general proof}.

For the non-local terms with $j \geq 2$ we can now follow Steps~2 and 3 \emph{verbatim}, the only modification being that we multiply the Caccioppoli estimate \eqref{eq0: Carleson bound Dirichlet} by $s^{-2\alpha}$ before summing. This leads to \eqref{eq1: Carleson bound Dirichlet} with the local bound $\ell^\alpha \|f\|_{\Lamdot^{\alpha}}$ replaced by $\|f\|_{\BMOdot^\alpha}$ and additional powers $s^{-2\alpha}$ in each of the off-diagonal pieces, so that the power $\ell^\alpha$ in \eqref{eq2: Carleson bound Dirichlet} disappears. Thus, we control the sum of the off-diagonal pieces by $\|f\|_{\Lamdot^{\alpha}} \lesssim \|f\|_{\BMOdot^\alpha}$.

The proof is complete modulo the argument for \eqref{eq: Local bound Holder-dir BMO} that we give now. By translation we can assume that $Q$ is centered at the origin. A classical argument using the Fourier transform of $f_1 \in \L^2$ yields
\begin{align*}
	\|f_1\|_{\Hdot^{\alpha,2}}^2 \simeq \int_{\R^n} \int_{\R^n} \frac{|f_1(y) - f_1(z)|^2}{|y-z|^{n+2 \alpha}} \, \d z \d y \eqqcolon \I,
\end{align*}
see for example~\cite[Prop.~1.3.7]{Danchin}. According to \eqref{eq: Strichartz BMO norm} it suffices to prove $\I \lesssim \mathrm{A}$, where
\begin{align*}
  \mathrm{A} \coloneqq \int_{4Q} \int_{4Q} \frac{|f(y) - f(z)|^2}{|y-z|^{n+2 \alpha}} \, \d z \d y.
\end{align*}
By symmetry, we have
\begin{align*}
	\I = 2	\iint_{|y|_\infty \geq |z|_\infty} \frac{|f_1(y) - f_1(z)|^2}{|y-z|^{n+2 \alpha}} \,\d z \d y,
\end{align*}
where $|\, \cdot\,|_\infty$ is the $\ell^\infty$-norm on $\R^n$. We write the numerator as $\eta_1(y)(f(z)-f(y)) + f(z)(\eta_1(z)-\eta_1(y))$. The first term vanishes unless $y \in 4Q$ and in that case $z \in 4Q$ follows from $|y|_\infty \geq |z|_\infty$. Hence, the integral of this part is controlled by $\mathrm{A}$. Likewise, $z \in {}^c(4Q)$ implies $y \in {}^c(4Q)$ and the second term vanishes. Altogether, we obtain
\begin{align*}
	\I 
	&\leq  \mathrm{A} + \int_{\R^n} \int_{4Q} \frac{|f(z)|^2 |\eta_1(z)-\eta_1(y)|^2}{|y-z|^{n+2 \alpha}}  \, \d z \d y,
	\intertext{where we bound $|\eta_1(z)-\eta_1(y)|$ via the mean value theorem if $|z-y| \leq \ell(Q)$ and in $\L^\infty$-norm if not, in order to get}
	&\leq \mathrm{A} + \frac{C^2}{\ell(Q)^2} \int_{4Q} \int_{|y-z| \leq \ell(Q)} \frac{|f(z)|^2}{|y-z|^{n+2 \alpha -2}}  \, \d y \d z \\
	&\quad + 4C^2 \int_{4Q} \int_{|y-z| \geq \ell(Q)} \frac{|f(z)|^2}{|y-z|^{n+2 \alpha }}  \, \d y \d z \\
	&\lesssim \mathrm{A} + \frac{1}{\ell(Q)^{2 \alpha}} \int_{4Q} |f(z)|^2 \, \d z \\
	&= \mathrm{A} + \frac{1}{\ell(Q)^{2 \alpha}} \int_{4Q} \bigg|\barint_Q f(z) - f(y) \, \d y \bigg|^2 \d z \\
	&\lesssim \mathrm{A}.
\end{align*}
We used $(f)_Q = 0$ in the second to last step and Jensen's inequality and $|y-z| \lesssim \ell(Q)$ in the final step.

\medskip

\noindent \emph{Modification of Part~5: The bound `\,$\gtrsim$'}. We fix $g \in \C_0^\infty$ with $\int_{\R^n} g \d x = 0$ and consider the extended $\L^2$-duality pairing $\langle f, g \rangle$. We use the same notation as in Part~5 of Section~\ref{sec: existence of Holder-dir}. The only difference in the argument appears in Step~5, where we have to handle
\begin{align}
\label{eq: BMO-Besov-dir lower bound}
 \iint_{\reu} |t\nabla_x u| \cdot |t\nabla_x \wt{v}| \, \frac{\d t  \d x}{t}
\end{align}
by a duality. The argument is repeated twice in Step~6 for $t$-derivatives of $u$. The control of these integrals determines the bound for $|\langle f,g \rangle|$. We recall from \eqref{eq: Holder dir v-tilde via BH} the notation
\begin{align*}
	\begin{bmatrix} 0 \\ t \nabla_x \wt{v} \end{bmatrix} = \psi(tDB_H)\begin{bmatrix} g \\ 0 \end{bmatrix},
\end{align*}
where $\psi \in \Psi_1^{2 \beta +1}$ with $\beta > \nicefrac{n}{2} +2$ and $D B_H$ correspond to $H = -(a^*)^{-1} \Delta_x$ in the same way as $DB$ corresponds to $L$. In Section~\ref{sec: existence of Holder-dir} we have interpreted \eqref{eq: BMO-Besov-dir lower bound} as a $\T^{0,\infty;\alpha} - \T^{\varrho}$ duality pairing, where $\varrho \in (1_*,1]$ was such that $\alpha = n(\nicefrac{1}{\varrho}-1)$, in order to bring $C_\alpha(t \nabla u)$ into play. 

Now, we use the $\Y^{\alpha, \infty} - \Y^{-\alpha, 1}$ pairing, see Sections~\ref{subsec: tent spaces} and \ref{subsec: Z spaces}, in order to give
\begin{align*}
	\iint_{\reu} |t\nabla_x u| \cdot |t\nabla_x \wt{v}| \, \frac{\d t  \d x}{t} 
	&\lesssim \|\nabla_x u\|_{\Y^{\alpha-1, \infty}} \|t\nabla_x \wt{v}\|_{\Y^{-\alpha,1}}.
\end{align*}
Since $\beta > \nicefrac{n}{2} + 2$, the function $\psi$ is admissible for defining $\IX^{-\alpha,1}_{DB_H}$. We have $p_-(H) = 1_*$ and $q_+(H) = \infty$ (Corollary~\ref{cor: McIntosh-Nahmod}) and consequently the identification region for $DB_H$ in Figure~\ref{fig: diagram DB} contains the full open segment that joins $(1,-1)$ to $(1,0)$, see also Figure~\ref{fig: diagram n=1}. In particular, $\IX^{-\alpha,1}_{DB_H} = \IX^{-\alpha,1}_{D}$ and together with Figure~\ref{fig: identification D-adapted} we obtain
\begin{align*}
 \|t\nabla_x \wt{v}\|_{\Y^{-\alpha,1}}
 \simeq \bigg\|\begin{bmatrix} g \\ 0 \end{bmatrix} \bigg \|_{\IX_{DB_H}^{-\alpha,1}} 
 \simeq \|g\|_{\Xdot^{-\alpha,1}}.
\end{align*}
Thus, we control \eqref{eq: BMO-Besov-dir lower bound} by $\|\nabla u\|_{\Y^{\alpha-1, \infty}}  \|g\|_{\Xdot^{-\alpha,1}}$ and we conclude for all $g \in \C_0^\infty$ with $\int_{\R^n} g \d x = 0$ that
\begin{align*}
	|\langle f,g \rangle| \lesssim \|\nabla u\|_{\Y^{\alpha-1, \infty}}  \|g\|_{\Xdot^{-\alpha,1}}.
\end{align*}
These $g$ form a dense subclass of $\Xdot^{-\alpha,1}$. There are probably many ways to see this -- one is to use the smooth atomic decomposition for $\Xdot^{-\alpha,1}$ in \cite[Thm.~5.11 \& 5.18]{Frazier-Jawerth-Weiss}. By duality, we obtain the lower bound
\begin{align*}
	\|f\|_{\Xdot^{\alpha,\infty}} \lesssim \|\nabla_x u\|_{\Y^{\alpha-1, \infty}}. &\qedhere
\end{align*}
\end{proof}

Let us come back to Figure~\ref{fig: diagram general proof3} but for $p_+(L) > n$, so that the left lower vertex of the red triangle is situated at the origin. Proposition~\ref{prop: BMO-Besov-dir} allows us to add a segment on the line $\nicefrac{1}{p} = 0$ and we can try to interpolate again to enlarge the region for compatible solvability as illustrated in Figure~\ref{fig: diagram p+large proof}. This is the content of the final result in this section.

\begin{figure}[h]
\begin{center}
\begin{tikzpicture}[scale=2.4]
	\newcommand\fracspace{\vphantom{\frac{1}{1}}};
	\def\dimension{6};
	\def\xlength{3.5};
	\pgfmathsetmacro\xstretch{\xlength/(1+1/\dimension)}; 
	
	\pgfmathsetmacro\UpSob{\xstretch*(1/\dimension)-0.25};
	\pgfmathsetmacro\UpSobStar{\xstretch*(\UpSob/\xstretch-1/\dimension)};
	\pgfmathsetmacro\UpSobStarBesov{(-2)*\UpSobStar*\dimension/\xlength};
	\pgfmathsetmacro\LowSob{\xstretch*(1/2+1/\dimension)+0.5};
	\pgfmathsetmacro\LowSobStar{\xstretch*(\LowSob/\xstretch+1/\dimension)};
	\pgfmathsetmacro\UpGradient{\xstretch*(1/2)-0.2};
	\pgfmathsetmacro\UpGradientDual{\xstretch*(1/2)+0.3};
	\pgfmathsetmacro\Half{\xstretch*0.5};
	\pgfmathsetmacro\One{\xstretch*1};
	\coordinate (P00) at (1+\UpGradient,2);
	\coordinate (P01) at (1+\LowSob,2);
	\coordinate (P10) at (1+\UpSob,0);
	\coordinate (P11) at (1+\UpGradientDual,0);
	\coordinate (P12) at (1+\LowSob,0);
	\coordinate (DirExtra) at (1,0);
	\coordinate (DirHolder) at (1+\UpSobStar,0);
	\coordinate (DirHolderBesov) at (1, \UpSobStarBesov);
	\coordinate (RegExtra) at (1+\LowSobStar,2);
	\coordinate (P000) at (1-\xstretch/\dimension,0);
	\coordinate (Energy) at (1+\Half,1);

	\draw [thin] (0.3,2) -- (1+\xlength,2); 
	\draw [thin] (0.3,0) -- (1+\xlength,0); 
	
	\draw [thick,->] (0.3,-0.5) -- (1+\xlength+0.2,-0.5);
	\draw [fill=black] (1+\Half,-0.5) circle [radius = .5pt];
	\node [below] at (1+\Half,-0.5) {$\frac{1}{2 \fracspace}$};
	\draw [fill=black] (1+\UpGradient,-0.5) circle [radius = .5pt];
	\node [below] at (1+\UpGradient,-0.5) {$\frac{1}{q_+^L\fracspace}$};
	\draw [fill=black] (1+\UpGradientDual,-0.5) circle [radius = .5pt];
	\node [below] at (1+\UpGradientDual,-0.5) {$\frac{1}{(q_+^{L^\sharp})'\fracspace}$};
	\draw [fill=black] (1+\UpSob,-0.5) circle [radius = .5pt];
	\node [below] at (1+\UpSob,-0.5) {$\frac{1}{p_+^L\fracspace}$};
	\draw [fill=black] (1+\LowSob,-0.5) circle [radius = .5pt];
	\node [below] at (1+\LowSob,-0.5) {$\frac{1}{p_-^L \vee 1 \fracspace}$};
	\draw [fill=black] (1+\LowSobStar,-0.5) circle [radius = .5pt];
	\node [below] at (1+\LowSobStar,-0.5) {$\frac{1}{{(p_-^L)_* \vee 1_*\fracspace}}$};
	\draw [fill=black] (1+\xlength,-0.5) circle [radius = .5pt];
	\node [below] at (1+\xlength,-0.5) {$\frac{n+1}{n \fracspace}$};
	\node [right] at (1+\xlength+0.2,-0.5) {$\frac{1}{p \fracspace}$};
	\draw [fill=black] (1,-0.5) circle [radius = .5pt];
	\node [below] at (1,-0.5) {$0$};
	\draw [fill=black] (0.3,-0.5) circle [radius = .5pt];
	\node [below] at (0.2,-0.5) {$-\frac{1}{n \fracspace}$};
	
	\draw [thick,->] (-0,0) -- (0,2.2);
	\node [above] at (0,2.2) {$s$};
	\draw [fill=black] (0,2) circle [radius = .5pt];
	\node [left] at (0, \UpSobStarBesov) {$\hspace{-30pt}1\!-\! \frac{n}{p_+^L}$};
	\draw [fill=black] (0,  \UpSobStarBesov) circle [radius = .5pt];	
	\node [left] at (0,2) {$1$};
	\draw [fill=black] (0,0) circle [radius = .5pt];
	\node [left] at (0,0) {$0$};

	\draw [thin,dotted] (1,2) -- (1,0); 
	\draw [thin, dotted] (1+\Half, 2) -- (1+\Half,0);
	\draw [thin,dotted] (0.3,2) -- (0.3,0); 
	\draw [thin,dotted] (1+\UpGradient,2) -- (1+\UpGradient,0);
	\draw [thin,dotted] (1+\LowSob,2) -- (1+\LowSob,0);
	\draw [thin,dotted] (1+\xlength,2) -- (1+\xlength,0); 
	\draw [thin,dotted] (1+\LowSobStar,2) -- (1+\LowSobStar,0);
	
	\draw [thin,dotted] (0.3, \UpSobStarBesov) -- (1,\UpSobStarBesov); 
	\draw [thin,dotted] (0.3, 1) -- (1+\xlength,1); 
	
	\path [fill=lightgray, opacity = 0.6] (P00)--(P01)--(P11)--(P10)--(P00);
	\path [fill=blue!80!black, opacity = 0.4] (P11)--(P12)--(RegExtra)--(P01)--(P11);
	\path [fill=red!80!black, opacity = 0.4] (DirExtra)--(P10)--(P00)--(DirExtra);
	\path [fill=red!80!black, opacity = 0.4]
	(DirExtra) -- (DirHolderBesov) -- (P00) -- (DirExtra);
	\draw [ultra thick, gray] (P00) -- (P01); 
	\draw [ultra thick, blue!80!black] (P01) -- (RegExtra); 
	\draw [ultra thick, gray] (P10) -- (P11); 
	\draw [ultra thick, red!80!black] (DirExtra) -- (P10); 
	\draw [ultra thick, red!80!black] (DirExtra) -- (DirHolderBesov); 
	\draw [ultra thick, blue!80!black] (P11) -- (P12); 
	\draw [fill=white] (P00) circle [radius = .75pt];
	\draw [fill=blue!80!black] (P01) circle [radius = .75pt];
	\draw [fill=red!80!black] (P10) circle [radius = .75pt];
	\draw [fill=blue!80!black] (P11) circle [radius = .75pt];	
	\draw [fill=white] (P12) circle [radius = .75pt];
	\draw [fill=red!80!black] (DirExtra) circle [radius = .75pt];
	\draw [fill=white] (RegExtra) circle [radius = .75pt];
	\draw [fill=white] (DirHolderBesov) circle [radius = .75pt];
	\draw [fill=black] (Energy) circle [radius = .75pt];
	\draw [thick, dotted] (DirExtra) -- (P00); 
	\draw [thick, dotted] (DirHolderBesov) -- (Energy) -- (DirExtra); 
	\node [below] at (DirExtra) {$O$};
	\node [left] at (DirHolderBesov) {$X$};
	\node [above] at (Energy) {$E$};
	\node [above] at (P00) {$R$};
	\node [below] at (P10) {$D$};
\end{tikzpicture}
\end{center}
\caption{Extended region for compatible solvability of $(D)_{\Xdot^{s,p}}^\Le$ if $p_+(L) > n$ via a two-step interpolation argument. The picture is up to scale when $p_+(L) < \infty$. If $p_+(L) = \infty$, then $D=O$ and $X=(0,1)$. In case this happens while $p_-(L^\sharp)$ is still larger or equal to $1$, then the red region becomes the triangle $ORX$. If even $p_-(L^\sharp) < 1$, then also the color code changes and parts of the red region should turn into gray since now they belong to the identification region of Figure~\ref{fig: diagram general proof2}. This special situation is showcased in Figure~\ref{fig: diagram p-small} in the introduction, to which we refer.}
\label{fig: diagram p+large proof}
\end{figure}

\begin{prop}
\label{prop: fractional solvability p+>n}
Suppose that $p_+(L) > n$. If $(\nicefrac{1}{p}, s)$ is contained in the interior of the colored region in Figure~\ref{fig: diagram p+large proof}, then there is a compatible solution operator for $(D)_{\Xdot^{s,p}}^\Le$. In particular, the problem is compatibly solvable.
\end{prop}

\begin{proof} 
As before, it suffices to construct the compatible solution operator. In view of Proposition~\ref{prop: solvability fractional} it remains to consider points in the interior of the triangle $ORX$ and on the open segment $\cl{OR}$. Our starting point is that by Proposition~\ref{prop: BMO-Besov-dir} there is a compatible solution operator for the problems corresponding to the open segment $\cl{OX}$ and that the constructed solution has all the properties listed in Theorem~\ref{thm: Holder-dir}.

Fix any $P = (0,\alpha) \in \cl{OX}$. At $E \coloneqq (\nicefrac{1}{2},\nicefrac{1}{2})$ the corresponding problem is the Dirichlet problem for the energy class and we have the solution operator $\sol_E: \Xdot^{\nicefrac{1}{2}, 2} \to \cD'(\reu)/\IC^m$ from Proposition~\ref{prop: existence of energy solution}, which is compatible with the solution operator $\sol_P: \Xdot^{\alpha,\infty} \to \cD'(\reu)/\IC^m$ at $P$. Hence, we obtain a well-defined linear operator
\begin{align*}
	\sol: \Xdot^{\alpha, \infty} + \Xdot^{\nicefrac{1}{2},2} \to \cD'(\reu)/\IC^m
\end{align*}
such that $u = \sol f$ solves $\Le u = 0$ in $\reu$ and satisfies $u(t,\cdot) \to f$ as $t \to 0$ in $\cD'(\R^n)/\IC^m$. This time the compatibility with $\sol_E$ already holds by construction and no density argument is needed. Real and complex interpolation of the mapping properties at the endpoints yields that $\nabla \sol: \Xdot^{s,p} \to \Y^{s-1,p}$ is bounded provided $(\nicefrac{1}{p},s)$ belongs to the open segment $\cl{PE}$. This yields the required solution operator for $(D)_{\Xdot^{s,p}}^\Le$ and we can add the interior of the triangle $OEX$ in Figure~\ref{fig: diagram p+large proof} to the region of compatible solvability.

Now that we have successfully moved away from the line $\nicefrac{1}{p} = 0$ of infinite exponents, we can repeat the argument in the proof of Proposition~\ref{prop: solvability fractional} once more for any $P_0$ in the interior of $OEX$ and any $P_1$ in the interior of the gray region. In particular, we reach any point in the interior of $ORX$ and on the open segment $\cl{OR}$.
\end{proof}
\section{Single layer operators for \texorpdfstring{$\Le$}{L} and estimates for \texorpdfstring{$\Le^{-1}$}{L}}
\label{sec: single layer}

\noindent This section is needed to prepare the next section on uniqueness. We consider the divergence form operator 
\begin{align*}
\Le u = -\div A \nabla u = -\partial_t (a \partial_t u) - \div_x d \nabla_x u
\end{align*}
on $\ree$. It is of the same class as $\eo$ in \eqref{eq: Lax-Milgram operator} but in one dimension higher. Hence, $\Le $ is defined on $\Wdot^{1,2}(\R^{1+n})$ via the Lax-Milgram lemma and invertible onto $\Wdot^{-1,2}(\R^{1+n})$. It turns out that the inverse $\Le^{-1}$ on particular test functions can explicitly be constructed using abstract single layer operators $\cS_{t}^{\Le}$. All this relies on the fundamental observation of Rosén~\cite{R}\index{Theorem!Ros\'en's} that what is called \emph{single layer potential} in the classical context of elliptic operators with real coefficients can abstractly be defined using the $\H^\infty$-calculus for the perturbed Dirac operator $DB$ in \eqref{eq: BD and BD}. Here, we cite the equivalent formulation from \cite[Sec.~7]{AE} or \cite[pp.~100-101]{AusSta}, which is somewhat closer to our terminology. 

We define the \emph{conormal gradient}\index{conormal gradient} $\nabla_A \coloneqq [a \partial_t, \nabla_x]^\top$. For all $f \in \L^2$ and $t>0$ there is a unique distribution (up to a constant) that we denote by $\cS_{t}^\Le f$ such that
\begin{equation} 
\label{eq:st}
\nabla_A \cS_{t}^\Le f
\coloneqq \begin{cases}
+\e^{-tDB}\ind_{\IC^+}(D B) \begin{bmatrix} f \\ 0\end{bmatrix}   & \text{if } t>0, \\[20pt]
-\e^{-tD B}\ind_{\IC^-}(D B) \begin{bmatrix} f \\ 0\end{bmatrix}     &  \text{if } t<0.
\end{cases}
\end{equation}
Note that $[f, 0]^\top \in \cH = \cl{\ran(D)} = \cl{\ran(DB)}$, so that the right-hand side is defined in the same space via the bounded $\H^\infty$-calculus. Then, we have the following result.

\begin{prop}[{\cite[Prop.~4.5]{AE}}]
\label{prop: slp DB}
Assume $\wt G= \div_{x}G^\sharp$ with $G^\sharp \in \C_{0}^\infty(\ree; \IC^{mn})$. Let $H \coloneqq \Le^{-1}(\partial_t \wt G)$. Then $H$ is given for all $t\in \R$ as an $\L^2$-valued Bochner integral
\begin{align*} 
H(t,\cdot) = \int_{\R} \partial_t \cS_{t-s}^\Le \wt G(s,\cdot) \, \d s.
\end{align*}
\end{prop}

The reader may be surprised that the representation by convolution with the single layer is not a singular integral. This is due to a hidden integration by parts because we represent $\Le^{-1}(\partial_t \wt G)$ and not $ \Le^{-1}(\wt G)$, see \cite[Rem.~4.6]{AE}. We also note that $\partial_t \wt G \in \Wdot^{-1,2}(\ree)$ because it is a derivative of a test function.

For our purpose it will be more convenient to write the single layer operators in terms of the second-order operator $L$. This is the content of the following proposition.\index{single layer operator}

\begin{prop}
\label{prop: slo} Let $t\in \R, t\ne 0$ and $f\in \L^2$. Then
\begin{align*}
 \partial_t \cS_{t}^{\Le}f = \frac{1}{2} \sgn(t) \e^{-|t| {L^{1/2}}}(a^{-1}f).
\end{align*}
\end{prop}

\begin{proof}
We have $[z] = \sqrt{z^2} = \pm z$ in the complex half-planes $z \in \IC_{\pm}$. Hence, we can write the $\no$-component of \eqref{eq:st} as
\begin{align}
\label{eq1: slo}
\begin{split}
a \partial_t \cS_{t}^\Le f
\coloneqq \sgn(t) \begin{cases}
\bigg(\e^{-|t| [DB]} \ind_{\IC^+}(D B) \begin{bmatrix} f \\ 0\end{bmatrix} \bigg)_\no   & \text{if } t>0, \\[20pt]
\bigg(\e^{-|t| [DB]} \ind_{\IC^-}(D B)\begin{bmatrix} f \\ 0\end{bmatrix} \bigg)_\no  &  \text{if } t<0.
\end{cases}
\end{split}
\end{align}
If $[g_\no, g_\ta]^\top$ is in the range of $[DB]$, then the functional calculus on $\cl{\ran(D)}$ translates the identity of functions $\ind_{\IC^\pm}(z) = \nicefrac{1}{2}(1 \pm \nicefrac{z}{\sqrt{z^2}})$ into
\begin{align}
\label{eq2: slo}
\begin{split}
&\ind_{\IC^\pm}(DB)\begin{bmatrix} g_\no \vphantom{\tL} \\ g_\ta \vphantom{\tM} \end{bmatrix} \\
&\quad= 
\frac{1}{2} \left(
\begin{bmatrix} g_\no \vphantom{\tL} \\ g_\ta \vphantom{\tM} \end{bmatrix} \pm
\begin{bmatrix} 0\vphantom{\tL} & \div_x d \\ -\nabla_x a^{-1} & 0\vphantom{\tM} \end{bmatrix}
\begin{bmatrix} (\tL)^{-1/2} g_\no \\ (\tM)^{-1/2} g_\ta \end{bmatrix}
\right) \\
&\quad= \frac{1}{2} \left(
\begin{bmatrix} g_\no \vphantom{\tL} \\ g_\ta \vphantom{\tM} \end{bmatrix} \pm
\begin{bmatrix} \div_x d (\tM)^{-1/2} g_\ta \\ -\nabla_x L^{-1/2} a^{-1} g_\no \vphantom{\tL} \end{bmatrix}
\right),
\end{split}
\end{align}
compare with the matrix representations in \eqref{eq: BD and BD} and \eqref{eq: tL and tM}. We set $g_\ta = 0$ and apply the $[DB]$-semigroup to give
\begin{align*}
\bigg(\e^{-|t| [DB]} \ind_{\IC^\pm}(D B) \begin{bmatrix} g_\no \\ 0\end{bmatrix} \bigg)_\no 
= \frac{1}{2} \e^{-|t|\tL^{1/2}} g_\no
= \frac{1}{2} a \e^{-|t|L^{1/2}} a^{-1} g_\no.
\end{align*}
This identity extends to general $f \in \L^2$ in place of $g_\no$ since $\tL$ has dense range in $\L^2$ and the claim follows from \eqref{eq1: slo}. 
\end{proof}

Combining the previous two results gives us the following representation.\index{single layer operator!representation by}

\begin{cor}
\label{cor: H representation}
Assume $G=\pd_{t}\div_{x}G^\sharp$ with $G^\sharp \in \C_{0}^\infty(\ree; \IC^{mn})$ and set $\wt G \coloneqq \div_{x}G^\sharp$. Let $H \coloneqq \Le^{-1}(G)$. Then for all $t\in \R$, $H(t,\cdot)$ is given as an $\L^2$-valued Bochner integral by
\begin{equation}
\label{eq:H}
\begin{split}
H(t,\cdot) &= \frac{1}{2} \int_{\R}  \sgn(t-s) \e^{-|t-s| {L^{1/2}}}(a^{-1}\wt G(s,\cdot)) \, \d s.
\end{split}
\end{equation}
\end{cor}

As the formula for $H$ only uses the Poisson semigroup, we can use the range where the semigroup enjoys $\L^p$-estimates.  This leads to additional estimates as compared to \cite{AE} in the non-block case.  

\begin{lem}
\label{lem:H} Let $G, H$ be as in Corollary~\ref{cor: H representation} and suppose in addition that $\supp G^\sharp \subseteq [\nicefrac{1}{\beta},\beta]\times \R^n$ for some $\beta > 1$. Let $h: =H(0,\cdot)$. Then if $r\in (p_{-}(L)\vee 1, p_{+}(L))$, there is some $\gamma >0$ such that for all $t>0$, 
\begin{align}
\label{eq:Hr}
\begin{split}
\| H(t,\cdot)- \e^{-t L^{1/2}}h\|_{r} & \lesssim  t\wedge t^{-\gamma}, \\
\| \pd_{t}(H(t,\cdot)- \e^{-t L^{1/2}}h)\|_{r} &\lesssim  1\wedge t^{-1-\gamma}. 
\end{split}
\end{align}
If, in addition, $p_-(L^\sharp)<1$, then this also holds for $r=\infty$.
\end{lem}

\begin{proof} 
We remark that $a^{-1}\wt G(s,\cdot)$ belongs to any $\L^q$-space, uniformly in $s\in [\nicefrac{1}{\beta},\beta]$. We will choose $q$ at our convenience. 

We treat the case $r<\infty$ first. For the exponent $r$ we have at hand the estimates for the Poisson semigroup from Proposition~\ref{prop: Poisson bounds} and the $\H^\infty$-calculus on $\L^r \cap \L^2$, see Theorem~\ref{thm: main result Hoo}.

\medskip

\noindent \emph{Step~1}. We begin with the estimate for $H- \e^{-t L^{1/2}}h$ using \eqref{eq:H}. For $0< t\le \nicefrac{1}{4\beta}$, we have
\begin{align}
\label{eq:H-H1, tsmall}
\begin{split}
H(t, \cdot)- &\e^{-t L^{1/2}}h 
\\
&= - \frac{1}{2}   \int_{\nicefrac{1}{\beta}}^\beta   \Big(\e^{-(s-t)  L^{1/2}}-\e^{-(s+t)  L^{1/2}}\Big) (a^{-1}\wt G(s,\cdot)) \, \d s
\end{split}
\end{align}
and the operator in brackets is $\L^r$-bounded with bound $C t$ by the $\H^\infty$-calculus.

For $\nicefrac{1}{4\beta}<t< 4\beta$, we see that 
\begin{align}
\label{eq:H-H1, tmiddle}
\begin{split}
&H(t, \cdot)- \e^{-t L^{1/2}}h \\
&=\frac{1}{2}  \int_{\nicefrac{1}{\beta}}^\beta  \Big(\sgn(t-s) \e^{-|t-s| {L^{1/2}}} + \e^{-(s+t)  L^{1/2}} \Big) (a^{-1}\wt G(s,\cdot)) \, \d s,
\end{split}
\end{align}
and we get a uniform $\L^r$-bound.

Finally for $t\geq 4 \beta$,  we have
\begin{align}
\label{eq:H-H1, tlarge}
\begin{split}
&H(t, \cdot)- \e^{-t L^{1/2}}h \\
&=\frac{1}{2}  \int_{\nicefrac{1}{\beta}}^\beta  \Big( \e^{-(t-s) {L^{1/2}}} + \e^{-(s+t)  L^{1/2}} \Big) (a^{-1}\wt G(s,\cdot)) \, \d s
\\
&=\frac{1}{2}  \int_{\nicefrac{1}{\beta}}^\beta   \e^{-(t-2\beta) {L^{1/2}}} \Big( \e^{-(2\beta-s) {L^{1/2}}} + \e^{-(2\beta+s)  L^{1/2}} \Big) (a^{-1}\wt G(s,\cdot)) \, \d s.
\end{split}
\end{align}
We pick any $q \in (p_-(L)\vee 1, r)$. In the last line, the operator in brackets is $\L^q$-bounded, uniformly, and the operator to its left is $\L^q-\L^r$-bounded with norm controlled by $t^{-n/q+n/r}$. Hence, we get the required estimate with $\gamma = \nicefrac{n}{q}-\nicefrac{n}{r}$.

\medskip

\noindent \emph{Step~2}. We turn to estimates for $\pd_{t}(H- \e^{-t L^{1/2}}h)$ on differentiating  \eqref{eq:H}. For $t > 4\beta$ we have
\begin{align*}
\begin{split}
&\pd_{t}(H(t, \cdot)- \e^{-t L^{1/2}}h) \\
&=  -  \frac{1}{2}   \int_{\nicefrac{1}{\beta}}^\beta   L^{1/2}\Big(\e^{-(t-s)  L^{1/2}} + \e^{-(s+t)  L^{1/2}}\Big) (a^{-1}\wt G(s,\cdot)) \, \d s.
\end{split}
\end{align*}
We expand the kernel as
\begin{align*}
L^{1/2}&\Big(\e^{-(t-s)  L^{1/2}} + \e^{-(s+t)  L^{1/2}}\Big) \\
&=  \Big( \e^{-(\frac{t}{2}-\beta) L^{1/2}}\Big) \Big(L^{1/2}   \e^{-(\frac{t}{2}-\beta) L^{1/2}}\Big)  \Big(\e^{-(2\beta-s)  L^{1/2}} + \e^{-(2\beta+s)  L^{1/2}}\Big)
\end{align*}
and pick again any $q \in (p_-(L)\vee 1, r)$. On the right-hand side the third operator is uniformly $\L^q$-bounded, the second one is $\L^q$-bounded with bound controlled by $t^{-1}$ and the first one is $\L^q -\L^r$-bounded with bound controlled by $t^{-n/q+n/r}$.  We get again the required bound with $\gamma = \nicefrac{n}{q}-\nicefrac{n}{r}$.

For $0 < t \leq 4\beta$ we need a uniform $\L^r$-bound. We are integrating over the singularity at $t=s$ in \eqref{eq:H} but using the convolution structure, we can compute with $G=\pd_{t}\wt G$, 
\begin{align} 
\label{eq:dt(H-H1), tsmall}
\begin{split}
&\pd_{t}(H(t, \cdot)- \e^{-t L^{1/2}}h) 
\\&= \frac{1}{2}   \int_{\nicefrac{1}{\beta}}^\beta  \sgn(t-s)  \e^{-|t-s|  L^{1/2}} (a^{-1} G(s,\cdot)) \, \d s \\
&\quad -\frac{1}{2}   \int_{\nicefrac{1}{\beta}}^\beta L^{1/2} \e^{-(s+t)  L^{1/2}} (a^{-1}\wt G(s,\cdot)) \, \d s.
\end{split}
\end{align}
The operators inside the integrals are $\L^r$-bounded, uniformly for $s,t$ in the prescribed range. 

\medskip

Finally, we establish the $\L^\infty$-bounds under the additional assumption  $p_-(L^\sharp)<1$. This implies $p_+(L) = \infty$ by duality and similarity.

\medskip

\noindent \emph{Step~3}. We modify Step~1 as follows. 

If $t\le \nicefrac{1}{4\beta}$, then we pick any $r \in (p_-(L) \vee 1, \infty)$ and write the operator inside the integral in \eqref{eq:H-H1, tsmall} as 
\begin{align*}
	\e^{-(s-t)L^{1/2}} - \e^{-(s+t)L^{1/2}} = \e^{-\frac{s}{2} L^{1/2}} \Big(\e^{-(\frac{s}{2}-t)L^{1/2}} - \e^{-(\frac{s}{2}+t)L^{1/2}}\Big).
\end{align*}
On the right we use the $\L^r-\L^\infty$-bound for $\e^{- \nicefrac{s}{2} L^{1/2}}$, see Proposition \ref{prop: Poisson bounds}.(iii), which is uniform since $s \in [\nicefrac{1}{\beta}, \beta]$, and the $\L^r$-bound with bound $Ct$ that follows from the $\H^\infty$-calculus as before.

If $\nicefrac{1}{4\beta}<t< 4\beta$, then for any $q \in (p_-(L) \vee 1, \infty)$ the operator inside the integral in \eqref{eq:H-H1, tmiddle} is $\L^q-\L^\infty$-bounded with norm controlled by $|s-t|^{-n/q}$. We can pick $q>n$ and this bound becomes integrable on $[\nicefrac{1}{\beta}, \beta]$. 

If $t\geq 4 \beta$, then thanks to Proposition~\ref{prop: Poisson bounds}.(iii) the same argument as before applies with $r = \infty$.

\medskip

\noindent \emph{Step~4}. We modify Step~2 as follows.

If  $t > 4\beta$, then thanks to Proposition \ref{prop: Poisson bounds}.(iii) the same argument as before applies with $r=\infty$.

If $0 < t \leq 4\beta$, then for any $q \in (p_-(L) \vee 1, \infty)$ the operator inside the first integral in \eqref{eq:dt(H-H1), tsmall} is $\L^q-\L^\infty$-bounded with norm controlled by $|s-t|^{-n/q}$ and choosing 
$q>n$ gives an integrable singularity. In the second integral we write
\begin{align*}
L^{1/2} \e^{-(s+t)L^{1/2}} = \e^{-(\frac{s}{2} + \frac{t}{2} )L^{1/2}}  \Big(L^{1/2} \e^{-(\frac{s}{2} + \frac{t}{2} )L^{1/2}}\Big).
\end{align*}
The operator on the far right is $\L^q$-bounded and the one to its left is $\L^q-\L^\infty$-bounded, both with uniform bounds since $\nicefrac{s}{2} + \nicefrac{t}{2} \in [\nicefrac{1}{2\beta}, \nicefrac{5 \beta}{2}]$.

\end{proof}

\section{Uniqueness in regularity and Dirichlet problems}
\label{sec: uniqueness}

\noindent This section complements Sections~\ref{sec: estimates semigroup}, \ref{sec: existence of Holder-dir} and \ref{sec: fractional}. We shall prove the uniqueness parts in Theorems~\ref{thm: blockdir}, \ref{thm: blockreg}, \ref{thm: Holder-dir} and \ref{thm: uniqueness fractional}.

In \cite{AE}, we developed a strategy to prove uniqueness for elliptic systems without regularity assumptions and with coefficients not necessarily in block form. We streamline the strategy in the case of the block system $\Le u = 0$ to obtain uniqueness of solutions in much greater generality.
\subsection{The strategy of proof}
\label{subsec: review of strategy uniqueness}

Throughout, we denote by $\dual \cdot \cdot$ the sesquilinear duality pairing between distributions and test functions in $\reu$. Since we are dealing with a linear equation, it suffices to assume that $u$ solves one of 
\begin{align*}
	(R)_p^\Le, (D)_p^\Le, (D)_{\Lamdot^\alpha}^\Le, (\wtD)_{\Lamdot^\alpha}^\Le, (D)_{\Xdot^{s,p}}^\Le
\end{align*}
with boundary data $0$ and show that this forces $u$ to vanish almost everywhere.

It begins with the following lemma in order to restrict the class of necessary testing conditions for $u$. The possible combinations of  an interior control with a boundary limit cover all cases that can occur in our BVPs.

\begin{lem}
\label{lem:prep} 
Let $u$ be a weak solution to $\Le u=-\div A \nabla u =0$ on $\reu$.  Let $\alpha \in [0,1)$ and $p \in (0,\infty)$. Assume one of the interior controls
\begin{multicols}{2}
\begin{itemize}
	\item $\NT( u) \in \L^p$, 
	\item $\NT (\nabla u) \in \L^p$,
	\item $W(t^{1-\alpha} \nabla u) \in \{\L^p, \L^\infty\}$,
	\item $S(t^{1-\alpha} \nabla u) \in \L^p$,
	\item $\NTsharpalpha(u) \in \L^\infty$,
	\item $C_\alpha(t \nabla u) \in \L^\infty$,
	\item $C_0(t^{1-\alpha} \nabla u) \in \L^\infty$,
\end{itemize}
\end{multicols}
\noindent and one of the boundary limits
\begin{align}
\label{eq:aecv0} \lim_{t \to 0} \bariint_{W(t,x)} |u(s,y)| \, \d s \d y &= 0 \quad  (\text{a.e. } x\in \R^n),\\
\label{eq:L2cv} \lim_{t \to 0} \barint_{t/2}^{2t} |u(s,\cdot)| \, \d s &= 0 \quad (\text{in } \L_{\loc}^2).
\end{align}
If $\dual u G =0$ for all test functions of the form $G=\pd_{t}\div_{x}G^\sharp$ with $G^\sharp \in \C_{0}^\infty(\reu; \IC^{mn})$, then $u=0$ almost everywhere.
\end{lem}

\begin{proof} We have  $\dual {\nabla_{x} \partial_{t}u} { G^\sharp}=0$, where $G^\sharp$ is an arbitrary test function in $\reu$. Hence, $\partial_t u \in \Lloc^2$ is independent of $x$ and we obtain 
\begin{align*}
	u(t,x)=g(t)+f(x)
\end{align*}
with $f \in \Lloc^2$ and $g: (0, \infty) \to \IC^m$ smooth (Corollary~\ref{cor: weak solution smooth in t}). If \eqref{eq:aecv0} holds, then we write
\begin{align*}
 \bariint_{W(t,x)} u(s,y)\, \d s \d y = \barint_{t/2}^{2t} g(s) \, \d s + \barint_{B(x,t)} f(y) \, \d y,
\end{align*}
where in the limit as $t \to 0$ the left-hand side tends to $0$ for a.e.\ $x \in \R^n$ by assumption and the second term on the right-hand side tends to $f(x)$ by  Lebesgue's differentiation theorem. Hence, $\barint_{t/2}^{2t} g(s) \d s$ has a limit as $t\to 0$ that we call $\beta \in \IC^m$ and we have $f(x) = -\beta$ almost everywhere.  The same conclusion holds under the assumption \eqref{eq:L2cv} since then
\begin{align*}
	\barint_{t/2}^{2t} u(s,\cdot) \, \d s = \barint_{t/2}^{2t} g(s) \, \d s + f(\cdot)
\end{align*}
tends to $0$ in $\Lloc^2$ as $t \to 0$.

So far we know that $u(t,x) = g(t) - \beta$. The equation for $u$ yields $a \partial_t^2 g = -\Le u =0$. Consequently, $g$ is a linear function. By definition of $\beta$ we get $g(t) = \gamma t + \beta$ for some $\gamma \in \IC^m$ and therefore $u(t,x) = \gamma t$. If $\gamma \neq 0$, then we get for all $x \in \R^n$ and all $t>0$ that 
\begin{multicols}{2}
	\begin{itemize}
		\item $\NT(u)(x) = \infty$,
		\item $\NT (\nabla u)(x) = |\gamma|$,
		\item $W(t^{1-\alpha} \nabla u)(t,x)~\simeq~t^{1-\alpha} |\gamma|$,
		\item $S(t^{1-\alpha} \nabla u) (x) = \infty$,
		\item $\NTsharpalpha(u)(x) = \infty$,
		\item $C_\alpha(t \nabla u)(x) = \infty$,
		\item $C_0(t^{1-\alpha} \nabla u)(x) = \infty$,
	\end{itemize}
\end{multicols}
\noindent and none of the interior controls is satisfied. Thus, $\gamma = 0$.
\end{proof}

Now, let $u$ be a solution to $\Le u=-\div A \nabla u =0$ on $\reu$.  We take $G$ as above. To compute $\dual u G$, we then pick a second function $\theta$, compactly supported in $\reu$, real-valued, Lipschitz continuous  and equal to $1$ on the support of $G$. Finally, we let $H \coloneqq (\Le^*)^{-1} G$, which is a weak solution to the adjoint equation
\begin{align*}
\Le ^*H= - \div A^* \nabla H = G \quad (\text{on } \ree).
\end{align*}
As $u\theta$ is a test function for this equation, we have
\begin{align*}
\dual u G= \dual {u\theta} G = \dual {A\nabla (u\theta)} {\nabla H}.
\end{align*}
Next,
\begin{align*}
\dual{A\nabla (u\theta)}{&\nabla H}   \\
& = \dual {A(u \otimes \nabla \theta)}{\nabla H} +  \dual {A (\theta \nabla u) }{\nabla H}  \\
& =  \dual {A(u \otimes \nabla \theta)}{\nabla H} - \dual {A\nabla u }{H \otimes \nabla \theta  } + \dual {A\nabla u }{\nabla ( \theta  H)},
\end{align*}
and the last term vanishes because $\theta H$ is a test function for $\Le u=0$. Our notation $\nabla(u \theta) = u \otimes \nabla \theta + \theta \nabla u$ is as predicted by the product rule.

We let $h\coloneqq H(0,\cdot) \in \L^2$ (see Corollary~\ref{cor: H representation}) and take
\begin{align}
\label{eq:H1}
H_1 (t,x) \coloneqq \e^{-t (L^\sharp)^{1/2}} h(x),
\end{align}
where  we recall that $L^\sharp$ corresponds to $\Le^*$ in the same way as $L$ corresponds to $\Le$. In particular, $H_{1}$ a solution to the adjoint problem $\Le ^*H_{1}=0 $ on $\reu$  with  boundary condition $h$, see Proposition~\ref{prop: Poisson smg is weak solution}. We can apply the same decomposition to $\dual {A\nabla (u\theta)} {\nabla H_{1}}$ and remark that this term vanishes since $u\theta$ is a test function for $\Le ^*H_{1}=0 $. Hence, we obtain
\begin{align}
\label{eq:uG}
\dual u G = \dual {A(u \otimes \nabla \theta)}{\nabla (H-H_{1})} - \dual {A\nabla u }{  (H-H_{1}) \otimes \nabla \theta}.
\end{align}
We remark that $u$ and $H-H_{1}$ vanish at the boundary in some sense. In fact, the reason to introduce $H_{1}$ is to help convergence near the boundary. 

Lemma~\ref{lem:prep} implies the following reformulation of uniqueness in the five BVPs.

\begin{prop} 
Suppose that $u$ solves one of the problems $(R)_p^\Le$, $(D)_p^\Le$, $(D)_{\Lamdot^\alpha}^\Le$, $(\wt D)_{\Lamdot^\alpha}^\Le$, $(D)_{\Xdot^{s,p}}^\Le$ with boundary data $f = 0$. If the right-hand side of \eqref{eq:uG} converges to $0$ as $\theta\to 1$ everywhere on $\reu$, then $u = 0$ almost everywhere.
\end{prop}

We prepare the limit procedure by picking $\theta$ in a more explicit way. For the rest of the section the following parameters will be used:\index{test functions!for uniqueness proofs}
\begin{align}
\label{eq:theta}
\begin{minipage}{0.88\linewidth}
\begin{itemize}
\item $G^\sharp \in \C_0^\infty(\reu; \IC^{mn})$ with support in $[\nicefrac{1}{\beta}, \beta] \times B(0, \beta) \subseteq \reu$ and $G \coloneqq \partial_t \div_x G^\sharp$,
\item $\chi\in \C_{0}^\infty(\R^n; \R)$ with $\ind_{B(0,1)} \leq \chi \leq \ind_{B(0,2)}$,
\item $\eta$ the continuous piecewise linear function, which is equal to $0$ on $[0,\frac{2}{3}]$, equal to $1$ on $ [\frac{3}{2}, \infty)$, and linear in between,
\item $M> 8\beta$ and  $0<\eps<\nicefrac{1}{4\beta}$ and $8\beta< R<\infty$ to finally set 
$$\qquad \theta(t,x) \coloneqq \chi\Big(\frac{x}{M}\Big)\eta\Big(\frac{t}{\eps} \Big) \Big(1-\eta \Big(\frac{t}{R}\Big)\Big).$$
\end{itemize}
\end{minipage}
\end{align}
We also use the block structure of $A$ to write
\begin{align*}
A(u \otimes \nabla \theta) 
= \begin{bmatrix} a u \partial_t \theta\\ d (u \otimes \nabla_x \theta) \end{bmatrix}, \quad 
A \nabla u
= \begin{bmatrix} a \partial_t u \\ d \nabla_x u  \end{bmatrix}.
\end{align*}
Due to the explicit form of $\theta$, we obtain for the first term on the right-hand side of \eqref{eq:uG} that
\begin{align}
\label{eq:firstuG}
|\dual {A(u \otimes \nabla \theta)}{\nabla (H-H_{1})}| \lesssim  I_{M,\varepsilon,R}+ J_{\varepsilon, M}+J_{R,M},
\end{align}
with  
\begin{align}
\label{eq:IM}
I_{M,\varepsilon,R} \coloneqq \frac{1}{M} \int_{M\le |y| \le 2M} \int_{2\varepsilon/3}^{3R/2} |u||\nabla_{x}(H-H_{1})| \, \d s \d y
\end{align}
and 
\begin{equation}
\label{eq:Jalpha}
J_{\tau, M} \coloneqq \int_{|y|\le 2M} \barint_{2\tau/3}^{3\tau/2} |u||\pd_{t} (H-H_{1})| \, \d s \d y.
\end{equation}
For the second term, we have
\begin{equation}
\label{eq:seconduG}
|\dual {A\nabla u}{(H-H_{1}) \otimes \nabla \theta}| \lesssim  \wt I_{M,\varepsilon,R}+ \wt J_{\varepsilon,M}+\wt J_{R,M},
\end{equation}
with  
\begin{equation}
\label{eq:tildeIM}
\wt I_{M,\varepsilon,R} \coloneqq \frac{1}{M} \int_{M\le |y| \le 2M} \int_{2\varepsilon/3}^{3R/2} |\nabla_{x} u||H-H_{1}| \, \d s \d y
\end{equation}
and 
\begin{equation}
\label{eq:tildeJalpha}
\wt J_{\tau, M} \coloneqq \int_{|y|\le 2M} \barint_{2\tau/3}^{3\tau/2} |\pd_{t}u||H-H_{1}| \, \d s \d y.
\end{equation}
Implicit constants depend only on dimensions and ellipticity. We need to specify the way how the parameters $M,R$ tend to $\infty$ and $\varepsilon $ tends to $0$ in order to make the terms on the right of \eqref{eq:firstuG} and \eqref{eq:seconduG} all tend to $0$.
\subsection{Uniqueness for \texorpdfstring{$\boldsymbol{(R)_p^\Le}$}{(R)_p} -- conclusion of the proof of Theorem~\ref{thm: blockreg}}
\label{subsec: uniqueness regularity}

We shall obtain uniqueness of solutions to $(R)_p^\Le$ in the range
\begin{align*}
	p_{-}(L)_{*}\vee  1_* <p<p_{+}(L).
\end{align*}
By Theorem~\ref{thm: standard relation J(L) and N(L)} we have $q_+(L) \leq p_+(L)$, so that this is even a larger range than for existence of a solution in Theorem~\ref{thm: blockreg}. We assume the interior control $\NT(\nabla u)\in \L^p$ and the convergence at the boundary \eqref{eq:aecv0} for almost every $x$. Then we distinguish two cases:
\begin{itemize}
	\item $(p_{-}(L)_{*} \vee 1_*) < p \leq (p_-(L) \vee 1)$,
	\item $(p_-(L) \vee 1) < p < p_+(L)$.
\end{itemize}
\subsection*{Case 1: \texorpdfstring{$\boldsymbol{{p_{-}(L)_* \vee 1_* <p \leq (p_-(L) \vee 1)}}$}{p- v 1* < p < p-v1}}

To implement the strategy in Section~\ref{subsec: review of strategy uniqueness}, we begin with the following lemma.

\begin{lem} 
\label{lem: exponent regularity proof}
If $0<p<r \leq 2$, then for any weak solution $u$ to $\Le u=0$ on $\reu$,  
\begin{equation}
\label{eq:pdtup*}
\bigg(\iint_{\reu} |\nabla u|^{r} t^{n(\frac{r}{p} - 1)} \, \frac{\d t \d x}{t}\bigg)^{\frac{1}{r}} \lesssim \|\NT(\nabla u)\|_{p}.
\end{equation}
Moreover, if \eqref{eq:aecv0} holds and in addition $p>1_*$, then
\begin{equation}
\label{eq:u/tp*}
\bigg(\iint_{\reu} | {u}|^{r} t^{n(\frac{r}{p} - 1)-r} \, \frac{\d t \d x}{t} \bigg)^{\frac{1}{r}} \lesssim \|\NT(\nabla u)\|_{p}.
\end{equation}
\end{lem}

\begin{proof}
The first inequality is due to Lemma~\ref{lem:AM} applied to $F \coloneqq |\nabla u|$. For the second inequality, Proposition~\ref{prop: KP}.(iii) yields  $\|\NTone(\nicefrac{u}{t})\|_{p}\lesssim \|\NT(\nabla u)\|_{p}$, where $\NTone$ is a non-tangential maximal function that uses $\L^1$-averages instead of $\L^2$-averages. But as $u$ is a weak solution to $\Le u=0$, it satisfies reverse H\"older inequalities. Hence, $\|\NT(\nicefrac{u}{t})\|_{p} \lesssim  \|\NT(\nabla u)\|_{p}$, where we also used a change of parameters in non-tangential maximal functions (Lemma~\ref{lem: independence of Whitney parameters}). Applying Lemma~\ref{lem:AM} to $F \coloneqq \nicefrac{u}{t}$ concludes the proof.
\end{proof}

We fix an exponent $r$ such that $(p_-(L) \vee 1) < r \leq 2$. Then the assumption $p< r$ in Lemma~\ref{lem: exponent regularity proof} holds automatically and we have $2 \leq r' < p_+(L^\sharp)$ by duality and similarity. Next, we recall that $H_{1}(t,\cdot)= \e^{-t(L^{\sharp})^{1/2}} h$, where $h$ is the trace of $H$ at $t=0$. We have at our disposal the estimates of Lemma~\ref{lem:H} with $L^{\sharp}$ replacing $L$. In particular, we obtain for some $\gamma>0$ and all $t>0$,
\begin{align}
\label{eq:Hr'}
\begin{split}
\| H(t,\cdot)- H_1(t,\cdot)\|_{r'} & \lesssim  t\wedge t^{-\gamma}, \\
\| \pd_{t}(H(t,\cdot)- H_1(t,\cdot))\|_{r'} &\lesssim  1\wedge t^{-1-\gamma}. 
\end{split}
\end{align}

We come to taking limits in \eqref{eq:firstuG} and \eqref{eq:seconduG}. We shall send $M \to \infty$ for $\eps, R$ fixed and then send $\eps \to 0$ and $R \to \infty$. We start with the terms on the right-hand side of \eqref{eq:firstuG}.

\medskip
\noindent{\emph{The term $I_{M,\varepsilon,R}$}}. We can bound $MI_{M,\varepsilon,R}$ by a finite number (depending on $\varepsilon,R$) of integrals
\begin{align*}
K_{\tau, M} \coloneqq \int_{|y| \ge M} \int_{2\tau/3}^{3\tau/2}  |u||\nabla_x(H-H_{1})| \, \d s \d y \quad (\eps \leq \tau \leq R)
\end{align*}
and it suffices to bound each of them uniformly for $M$ large, say $M\ge 10R$. 

Because we do not have global bounds on $\nabla_x (H-H_{1})$, we argue as follows. We let $w(\tau,x)\coloneqq (\nicefrac{2 \tau}{3}, \nicefrac{3 \tau}{2}) \times B(x, \nicefrac{\tau}{2})$ denote slightly smaller Whitney boxes and use an averaging trick to give
\begin{align*}
K_{\tau, M}    
&\lesssim  \int_{|x| \ge M/2}  \bigg(\bariint_{w(\tau,x)} |u||\tau\nabla_x (H-H_{1})|\bigg) \, \d x  \\
&\lesssim  \int_{|x| \ge M/2} \bigg(\bariint_{w(\tau,x)} | u|^2 \bigg)^{\frac{1}{2}}\bigg(\bariint_{w(\tau,x)} |\tau \nabla_x (H-H_{1})|^2 \bigg)^{\frac{1}{2}}\ \d x \\
&\lesssim  \int_{|x| \ge M/2} \bigg(\bariint_{W(\tau,x)}| u|^r \bigg)^{\frac{1}{r}}\bigg(\bariint_{ W(\tau,x)} | H-H_{1}|^{r'} \bigg)^{\frac{1}{r'}} \, \d x,
\end{align*} 
where for the last line we used reverse H\"older estimates for $u$ and the Caccioppoli estimate followed by H\"older's inequality for $H-H_{1}$, which is a weak solution to $\Le^*(H-H_1) = 0$ in a neighborhood of each $W(\tau,x)$. Indeed, $\Le^*(H-H_1) = G$ on $\reu$ but as $W(\tau,x) \subseteq \{(t,y): |y| \geq 4R\}$, all Whitney boxes are outside the support of $G$, see \eqref{eq:theta}. By H\"older's inequality in $x$, we have
\begin{align*}
K_{\tau, M} 
&\lesssim   
\bigg(\!\int_{\R^n} \bariint_{W(\tau,x)} \!\!\!|u|^r \d y \d t \d x \bigg)^{\frac{1}{r}} \bigg(\!\int_{\R^n}\bariint_{ W(\tau,x)} \!\!\!| H-H_{1}|^{r'} \d y \d t \d x \bigg)^{\frac{1}{r'}} \\
& \lesssim  \bigg(\int_{\tau/2}^{2 \tau} \int_{\R^n} |u|^r \, \frac{\d y \d t}{t} \bigg)^{\frac{1}{r}} \bigg(\barint_{\tau/2}^{2\tau} \int_{\R^n} | H-H_{1}|^{r'} \, \d y \d t \bigg)^{\frac{1}{r'}} 
\\ &\lesssim \tau^{1-n(\frac{1}{p} - \frac{1}{r})} \|\NT(\nabla u)\|_{p} \ (\tau\wedge \tau^{-\gamma})
\end{align*}
using  \eqref{eq:u/tp*} and \eqref{eq:Hr'}. Summing up in $\tau$, we conclude $MI_{M,\varepsilon,R}\lesssim  \|\NT(\nabla u)\|_{p}$ with an implicit constant that depends on $\eps, R$ but not on $M$. Thus, $I_{M,\varepsilon,R} \to 0$ in the limit as $M \to \infty$.

\medskip
\noindent \emph{The term $J_{\varepsilon, M}$}. When $M\to \infty$, we have to take the $\d x$-integral on all of $\R^n$ and we can use H\"older's inequality directly to obtain a bound for the limit by
\begin{align}
\label{eq: JepsM bound regularity}
\begin{split}
&\bigg(\barint_{2\varepsilon/3}^{3\varepsilon/2} \int_{\R^n}  | u  |^r \, \d y \d t\bigg)^{\frac{1}{r}} \bigg(\barint_{2\varepsilon/3}^{3\varepsilon/2} \int_{\R^n}  | \pd_{t}(H-H_{1})|^{r'} \, \d y \d t \bigg)^{\frac{1}{r'}}
\\
&\lesssim \eps^{1-n(\frac{1}{p} - \frac{1}{r})} \bigg(\int_{2\varepsilon/3}^{3\varepsilon/2} \int_{\R^n}  | u  |^r t^{n(\frac{r}{p} - 1)-r} \, \frac{\d y \d t}{t}\bigg)^{\frac{1}{r}}
\end{split}
\end{align}
using the estimate  \eqref{eq:Hr'} when $\eps<1$. At this point we have to discuss the choice of $r$.

In dimension $n \geq 2$ we set $r \coloneqq p^*$. In order to see that this choice is admissible, we first note that Proposition~\ref{prop: J(L) contains neighborhood of Sobolev conjugates} yields $(p_-(L) \vee 1) \leq 2_*$ and therefore $r \leq 2$ follows from the upper bound on $p$. Likewise, the lower bound on $p$ implies $r> (p_-(L) \vee 1)$. For this choice of $r$ the exponent of $\eps$ in \eqref{eq: JepsM bound regularity} vanishes and we conclude from \eqref{eq:u/tp*} that the remaining integral converges to $0$ as $\eps \to 0$.

In dimension $n=1$ we have $1^* = \infty$ and hence we must argue differently. Proposition~\ref{prop: J(L) contains neighborhood of Sobolev conjugates} yield $p_-(L) = \nicefrac{1}{2} = 1_*$. Hence, our assumption on $p$ is $\nicefrac{1}{2} < p \leq 1$ and this allows us to pick $r>1$ sufficiently close to $1$ such that $\nicefrac{1}{r} \geq \nicefrac{1}{p} - 1$. Consequently, the exponent for $\eps$ in \eqref{eq: JepsM bound regularity} is non-negative and we conclude as before.

\medskip
\noindent \emph{The term $J_{R, M}$}. Similarly, we have a bound for the limit as $M \to \infty$ by
\begin{align}
\label{eq: JRM bound regularity}
\begin{split}
&\bigg(\barint_{2R/3}^{3R/2} \int_{\R^n}  | u  |^r \, \d y \d t\bigg)^{\frac{1}{r}} \bigg(\barint_{2R/3}^{3R/2} \int_{\R^n} | \pd_{t}(H-H_{1})|^{r'} \, \d y \d t \bigg)^{\frac{1}{r'}}
\\
&\lesssim  R^{-n(\frac{1}{p} - \frac{1}{r})-\gamma} \bigg(\int_{2R/3}^{3R/2} \int_{\R^n}  | u  |^r t^{n(\frac{r}{p} - 1)-r} \, \frac{\d y \d t}{t}\bigg)^{\frac{1}{r}},
\end{split}
\end{align}
using \eqref{eq:Hr'} when $R>1$. Since we have $r>p$ in any case, we get a negative power of $R$ in front of the integral and in view of \eqref{eq:u/tp*} this term tends to $0$ as $R\to \infty$.

\medskip

We next consider the terms on the right-hand side of \eqref{eq:seconduG}.

\medskip 
\noindent \emph{The term $\wt I_{M,\varepsilon,R}$}. H\"older's inequality yields that $M\wt I_{M,\varepsilon,R}$ is bounded by
\begin{align*}
\bigg(\int_{2\varepsilon/3}^{3R/2} \! \!\int_{\R^n} |t \nabla_x u|^r \, \d y \d t \bigg)^{\frac{1}{r}} \bigg(\int_{2\varepsilon/3}^{3R/2} \!\! \int_{\R^n} \bigg|\frac{H-H_{1}}{t}\bigg|^{r'} \, \d y \d t \bigg)^{\frac{1}{r'}}.
\end{align*}
Using \eqref{eq:pdtup*} and \eqref{eq:Hr'}, we thus obtain $M\wt I_{M,\varepsilon,R} \lesssim \|\NT(\nabla u)\|_{p}$ with an implicit constant that depends on $\eps, R$ but not on $M$. Hence, we have $\wt I_{M,\varepsilon,R} \to 0$ in the limit as $M \to \infty$.

\medskip

\noindent \emph{The term $\wt J_{\varepsilon, M}$}. We have again the following bound for the limit as $M \to \infty$ by taking the $\d x$-integral on $\R^n$ and using H\"older's inequality directly:
\begin{align}
\label{eq: JepsM-tilde bound regularity}
\begin{split}
&\bigg(\barint_{2\varepsilon/3}^{3\varepsilon/2} \int_{\R^n}  |t\pd_{t}u  |^r \, \d y \d t\bigg)^{\frac{1}{r}} \bigg(\barint_{2\varepsilon/3}^{3\varepsilon/2} \int_{\R^n} \bigg| \frac{H-H_{1}}{t}\bigg|^{r'} \, \d y \d t \bigg)^{\frac{1}{r'}}
\\
&\lesssim \eps^{1-n(\frac{1}{p} - \frac{1}{r})} \bigg(\int_{2\varepsilon/3}^{3\varepsilon/2} \int_{\R^n}  | \partial_t u  |^r t^{n(\frac{r}{p} - 1)} \, \frac{\d y \d t}{t}\bigg)^{\frac{1}{r}},
\end{split}
\end{align}
where the second step is due to \eqref{eq:Hr'}. The exponent for $\eps$ is the same as in \eqref{eq: JepsM bound regularity} and thus becomes non-negative for the same choice of $r$ as before. It follows from \eqref{eq:pdtup*} that the remaining integral tends to $0$ as $\varepsilon\to 0$.
\medskip

\noindent \emph{The term $\wt J_{R, M}$}. Similarly, for the limit of $\wt J_{R,M}$ as $M \to \infty$, we have the bound
\begin{align}
\label{eq: JRM-tilde bound regularity}
\begin{split}
&\bigg(\barint_{2R/3}^{3R/2} \int_{\R^n} |t\pd_{t}u  |^r \, \d y \d t\bigg)^{\frac{1}{r}} \bigg(\barint_{2R/3}^{3R/2} \int_{\R^n} \bigg| \frac{H-H_{1}}{t}\bigg|^{r'} \, \d y \d t \bigg)^{\frac{1}{r'}}
\\
&\lesssim R^{-n(\frac{1}{p} - \frac{1}{r})-\gamma}  \bigg(\int_{2R/3}^{3R/2} \int_{\R^n}  |\pd_{t}u  |^r t^{n(\frac{r}{p} - 1)} \, \frac{\d y \d t}{t}\bigg)^{\frac{1}{r}},
\end{split}
\end{align}
using \eqref{eq:Hr'} when $R>1$. The exponent for $R$ is negative and in the limit as $R\to \infty$, the right-hand side tends to $0$, taking into account \eqref{eq:pdtup*}. The argument is complete.
\subsection*{Case 2: \texorpdfstring{$\boldsymbol{(p_-(L) \vee 1)<p<p_+(L)}$}{p- v 1 < p < p+}}

For this case we organize the limit procedure differently. We set $R = M$ and first send $\eps \to 0$ and then $M \to \infty$ in \eqref{eq:firstuG} and \eqref{eq:seconduG}. 

The interior control $\NT(\nabla u) \in \L^p$ and the boundary limit \eqref{eq:aecv0} enter the calculations in a particularly concise form via the trace estimate
\begin{align*}
 \bariint_{W(t,x)} |u| \, \d s \d y \lesssim t \NT(\nabla u)(x) \quad ((t,x) \in \reu)
\end{align*}
from Proposition~\ref{prop: KP}. The non-tangential maximal function $\NT(\nabla u)$ has no further meaning to our argument and we can proceed without any additional effort under the following general assumption: Besides \eqref{eq:aecv0} we assume that there exists $ \Theta \in \L^p$ and $\alpha \in [0,1]$ such that
\begin{align}
\label{eq: NTeq 1}
	U_t (x) \coloneqq  \bariint_{W(t,x)} |u| \, \d s \d y
\end{align}
is controlled via
\begin{align}
\label{eq: G-alpha-assumption}
	U_t(x) \leq t^\alpha \Theta(x) \quad ((t,x) \in \reu).
\end{align}
This generalization will have fruitful implications for some of the other boundary value problems.

We begin with the terms in \eqref{eq:firstuG}.

\medskip

\noindent \emph{The term $J_{\varepsilon, M}$}. We let $w(\tau,x)\coloneqq (\nicefrac{2 \tau}{3}, \nicefrac{3 \tau}{2}) \times B(x, \nicefrac{\tau}{2})$ denote slightly smaller Whitney boxes and use an averaging trick to give
\begin{align}
\label{eq: JepsM generic}
\begin{split}
	J_{\varepsilon, M}    
	& \lesssim  \int_{|x| \le 3M}  \bigg(\bariint_{w(\varepsilon,x)} |u||\pd_{t} (H-H_{1})|\bigg) \, \d x   \\
	& \lesssim  \int_{|x| \le 3M} \bigg(\bariint_{w(\varepsilon,x)} | u|^2 \bigg)^{\frac{1}{2}}\bigg(\bariint_{w(\varepsilon,x)} | \pd_t (H-H_{1})|^2 \bigg)^{\frac{1}{2}} \, \d x \\
	& \lesssim  \int_{|x| \le 3M} \bigg(\bariint_{W(\varepsilon,x)} | u| \bigg)\bigg(\bariint_{W(\varepsilon,x)} \bigg| \frac{ H-H_{1}}{\eps}\bigg| \bigg) \, \d x \\
	& =  \int_{|x| \le 3M} U_\eps(x) \bigg(\bariint_{ W(\varepsilon,x)} \bigg| \frac{ H-H_{1}}{\eps}\bigg| \, \d t \d y\bigg) \, \d x.
\end{split}
\end{align}
The third line is the combination of Caccioppoli's estimate and the reverse H\"older inequality for $H-H_{1}$, which solves $\Le^*(H-H_1) = G$ on $\reu$ but $W(\eps,x) \subseteq \{(t,y): 0 < t < \nicefrac{1}{2\beta}\}$ is outside the support of $G$, see \eqref{eq:theta}. 

Changing the order of integration brings into play the maximal function $\Max^\varepsilon$ restricted to balls with radii not exceeding $\varepsilon$, acting on $U_\eps$:
\begin{align*}
	J_{\varepsilon, M}    
	& \lesssim \int_{|y| \le 4M} \Max^\varepsilon(U_\eps)(y) \barint_{\varepsilon/2}^{2\varepsilon}  \bigg| \frac{ H-H_{1}}{t}\bigg| \, \d t \d y\\
	& \lesssim \bigg(\int_{|y|\le 4M} \Max^\varepsilon(U_\eps)^{p} \, \d y \bigg)^{\frac{1}{p}} \bigg(\int_{\R^n} \bigg(\barint_{\varepsilon/2}^{2\varepsilon}  \bigg| \frac{ H-H_{1}}{t}\bigg| \, \d t \bigg)^{p'} \d y \bigg)^{\frac{1}{p'}} \\
	& \lesssim \bigg(\int_{|y|\le 4M} \Max^\varepsilon(U_\eps)^{p} \, \d y\bigg)^{\frac{1}{p}} \bigg(\int_{\R^n} \barint_{\varepsilon/2}^{2\varepsilon}  \bigg| \frac{ H-H_{1}}{t}\bigg|^{p'}\, \d t \d y \bigg)^{\frac{1}{p'}} \\
	& \lesssim \bigg(\int_{|y|\le 5M} U_\eps^{p} \, \d y\bigg)^{\frac{1}{p}} \bigg(\barint_{\varepsilon/2}^{2\varepsilon} \int_{\R^n}  \bigg| \frac{ H-H_{1}}{t}\bigg|^{p'}\, \d y \d t \bigg)^{\frac{1}{p'}},
\end{align*} 
where we have used $\Max^\eps(U_\eps) \leq \Max(\ind_{B(0,5M)}U_\eps)$ on $B(0,4M)$ and the maximal theorem in the last line. 

The assumption on $p$ implies $(p_-(L^\sharp) \vee 1) <p'<p_+(L^\sharp)$ by duality and similarity. Thus, we may use Lemma~\ref{lem:H} for $H-H_1$ with $r=p'$ and obtain  for some $\gamma>0$ and all $t>0$ the bound
\begin{align}
	\label {eq:Hp'}
	\|H(t,\cdot)-H_1(t,\cdot)\|_{p'} \lesssim t \wedge t^{-\gamma}.
\end{align}
Thus, the second integral on the right in the estimate above is uniformly bounded in $\eps \leq 1$ and we are left with
\begin{align*}
	J_{\varepsilon, M} \lesssim \bigg(\int_{|y|\le 5M} U_\eps^{p} \, \d y\bigg)^{\frac{1}{p}} \qquad (\eps \leq 1).
\end{align*}
According to \eqref{eq: G-alpha-assumption} we have $U_\eps \leq \Theta \in \L^p$ for all $\eps \leq 1$, so that we can use the dominated convergence theorem when passing to the limit as $\eps \to 0$. By assumption \eqref{eq:aecv0} we have $U_\eps(x) \to 0$ for a.e.\ $x \in \R^n$ and $J_{\varepsilon, M}  \to 0$ follows. This completes the treatment of this term.

\medskip

\noindent \emph{The terms $I_{M,\eps,M}$ and $J_{M,M}$}. Having sent $\varepsilon\to 0$, we have to estimate 
\begin{align*}
	\lim_{\eps \to 0} I_{M,\eps,M}+ J_{M,M} \eqqcolon I_M + J_M,
\end{align*}
where 
\begin{align}
	\label{eq: IM and JM}
	\begin{split}
		I_{M} & \coloneqq \frac{1}{M} \int_{M\le |y| \le 2M} \int_{0}^{3M/2}  |u||\nabla_x(H-H_{1})| \, \d s \d y, \\
		J_{M} & \coloneqq \int_{|y|\le 2M} \barint_{2M/3}^{3M/2} |u||\pd_t(H-H_{1})|\, \d s \d y.
	\end{split}
\end{align}

We begin with $I_M$. In the following we use small Whitney regions $w(\tau,x) = (\nicefrac{2 \tau}{3}, \nicefrac{3 \tau}{2}) \times B(x, \nicefrac{2 \tau}{9})$. Let $\tau_{j} \coloneqq (\nicefrac{9}{4})^j$ for $j\in \IZ$ and let $j_{M}$ be the unique integer with $\tau_{j_{M} -1}\le M<\tau_{j_{M}}$. Then $\tau_{j_M} \leq \nicefrac{9M}{4}$ and 
\begin{align}
	\label{eq: MIM bound using Ktau}
	MI_{M}\le \sum_{j=-\infty}^{j_{M}}  K_{\tau_{j},M}
\end{align}
with 
\begin{align}
	\label{eq: Ktau definition}
	K_{\tau, M} \coloneqq \int_{M\le |y| \le 2M} \int_{2\tau/3}^{3\tau/2}  |u||\nabla_x(H-H_{1})| \, \d s \d y.
\end{align}
Applying Caccioppoli and reverse H\"older inequalities as usual, we obtain for $\tau \leq \nicefrac{9M}{4}$ that
\begin{align*}
	K_{\tau, M}
	& \lesssim  \int_{\frac{M}{2} \le |x|\le \frac{5M}{2}}  \bigg(\bariint_{w(\tau,x)} |u||\tau\nabla_x (H-H_{1})|\bigg) \, \d x   \\
	&   \lesssim  \int_{\frac{M}{2} \le |x|\le \frac{5M}{2}}  \bigg(\bariint_{w(\tau,x)} | u|^2 \bigg)^{\frac{1}{2}}\bigg(\bariint_{w(\tau,x)} | \tau\nabla_x (H-H_{1})|^2 \bigg)^{\frac{1}{2}} \, \d x \\
	&  \lesssim  \int_{\frac{M}{2} \le |x|\le \frac{5M}{2}}  \bigg(\bariint_{W(\tau,x)}| u| \bigg)\bigg(\bariint_{ W(\tau,x)} | {H-H_{1}}| \bigg) \, \d x \\
	&= \int_{\frac{M}{2} \le |x|\le \frac{5M}{2}} U_\tau(x) \bigg(\bariint_{ W(\tau,x)} | {H-H_{1}}| \bigg) \, \d x.
\end{align*}
To justify the interior estimates for $H-H_1$ on $w(\tau,x)$, we remark that $W(\tau,x)$ lies outside the support of $G$. Indeed, if $\tau \leq \nicefrac{M}{4}$, then $|y| \geq \nicefrac{M}{4} = \nicefrac{R}{4} > 2 \beta$ for all $(s,y) \in W(\tau,x)$, and if $\tau \geq \nicefrac{M}{4}$, then $s \geq  \nicefrac{M}{8} > \beta$, see \eqref{eq:theta}.

Now, we change again the order of integration to bring the maximal function into play and use \eqref{eq: G-alpha-assumption} to give
\begin{align}
\label{eq: Ktau estimate generic}
	K_{\tau, M}  \lesssim \tau^\alpha \int_{|y|\le 5 M} \Max(\Theta)(y) \bigg(\barint_{\tau/2}^{2\tau}  | H-H_{1}| \, \d t \bigg) \d y.
\end{align} 
We continue by
\begin{align*}
K_{\tau, M}  
	&\lesssim \tau^\alpha \bigg(\int_{\R^n} \Max(\Theta)(y)^p \, \d y\bigg)^{\frac{1}{p}}  \bigg(\int_{\R^n} \bigg|\barint_{\tau/2}^{2\tau}  | H-H_{1}| \, \d t\bigg|^{p'} \d y\bigg)^{\frac{1}{p'}} \\
	&\lesssim \tau^\alpha \bigg(\int_{\R^n} \Theta(y)^p \, \d y\bigg)^{\frac{1}{p}}  \bigg(\barint_{\tau/2}^{2\tau} \int_{\R^n} |H-H_{1}|^{p'} \, \d y  \d t \bigg)^{\frac{1}{p'}}\\
	&\lesssim \tau^{\alpha}(\tau \wedge \tau^{-\gamma}) \|\Theta\|_p,
\end{align*}
where we have used H\"older's inequality in the first line, the maximal theorem and Jensen's inequality in the second one and \eqref{eq:Hp'} in the third one. At this point we can go back to \eqref{eq: MIM bound using Ktau} and sum up the estimates for $\tau = \tau_j$ in order to obtain
\begin{align*}
	M I_M \lesssim (1+M^{\alpha - \gamma})\|\Theta\|_p.
\end{align*}
By assumption we have $\alpha \leq 1$ and $\gamma > 0$. Hence, $M$ appears with exponent smaller than $1$ on the right-hand side and we conclude $I_M \to 0$ in the limit as $M \to \infty$.

For $J_{M}$ in \eqref{eq: IM and JM} we can argue just as for $K_{\tau,M}$ with $\tau = M$ since we have not used the lower bound on $|y|$ to justify the interior estimates in \eqref{eq: Ktau estimate generic} when $\tau \geq \nicefrac{M}{4}$. This leads to
\begin{align}
\label{eq: JM estimate generic}
\begin{split}
	M J_M 
	& \lesssim \int_{|x|\le \frac{5M}{2}}  U_M(x) \bigg(\bariint_{ W(M,x)} | {H-H_{1}}| \bigg) \, \d x \\
	& \lesssim M^\alpha \int_{|y|\le 5 M} \Max(\Theta)(y) \bigg(\barint_{\tau/2}^{2\tau}  | H-H_{1}| \, \d t \bigg) \d y
\end{split}
\end{align}
and repeating the argument from \eqref{eq: Ktau estimate generic} onward yields the same bound
\begin{align*}
	M J_M \lesssim (1+M^{\alpha - \gamma})\|\Theta\|_p.
\end{align*}
As before, we conclude $J_M \to 0$ as $M\to \infty$.

\medskip At this point we have handled the terms in \eqref{eq:firstuG}. The argument for the terms in \eqref{eq:seconduG} is \emph{verbatim} the same. Indeed, all of our estimates concerning \eqref{eq:firstuG} have used reverse H\"older estimates on $u$ and $H-H_1$ and the Caccioppoli inequality to replace $\nabla (H-H_1)$ by $\frac{H-H_1}{t}$ in the $\L^2$-averages. Now, we simply use Caccioppoli inequalities to replace $t \nabla u $ by $u$ and obtain the same bounds.  The proof of Theorem~\ref{thm: blockreg} is complete.
\subsection{Uniqueness for \texorpdfstring{$\boldsymbol{(D)_p^\Le}$}{(D)_p} -- conclusion of the proof of Theorem~\ref{thm: blockdir}}
\label{subsec: uniqueness Lp Dirichlet}

We shall implement again the formalism of Section~\ref{subsec: review of strategy uniqueness}. The interval of allowable exponents is $p_{-}(L) <p<p_{+}(L)^*$ if $p_-(L)\geq 1$ and $1\le p<p_{+}(L)^*$ if $p_-(L)< 1$. Hence, we assume $\NT (u)\in \L^p$ and that \eqref{eq:aecv0} holds. We distinguish three cases: 
\begin{itemize}
	\item $(p_{-}(L)\vee 1)< p < p_{+}(L)$, 
	\item $p=1$ if $p_-(L)<1$,
	\item $p_{+}(L) \le p < p_{+}(L)^*$.
\end{itemize}
\subsection*{Case 1: \texorpdfstring{$\boldsymbol{{(p_{-}(L) \vee 1)<p<p_{+}(L)}}$}{(p- v 1)< p < p+}}

This is the range of exponents for the generic argument under the assumptions \eqref{eq:aecv0} and \eqref{eq: G-alpha-assumption}. In our concrete setting the latter holds with $\alpha = 0$ and $\Theta = \NT(u)$ and there is nothing more to do.
\subsection*{Case 2: \texorpdfstring{$\boldsymbol{{p_{-}(L) <1= p}}$}{p- < 1=p }}
We basically follow the generic argument in Case~2 for the regularity problem with $\alpha = 0$ and $\Theta = \NT(u) \in \L^1$. In addition, we incorporate the following estimate for $H-H_1$ that comes from Lemma \ref{lem:H} in the case $r=\infty$: for some $\gamma>0$ and all $t>0$, 
\begin{align}
\label{eq:Hinfty}
 \|H(t,\cdot)-H_1(t,\cdot)\|_\infty & \lesssim t \wedge t^{-\gamma}.
\end{align}
This uniform bound will allow us to avoid the maximal operator, which would not be $\L^p$-bounded in this case.

\medskip

\noindent \emph{The term $J_{\varepsilon, M}$}. By \eqref{eq: JepsM generic} we have
\begin{align*}
J_{\varepsilon, M} \lesssim  \int_{|x| \le 3M} U_\eps(x) \bigg(\bariint_{ W(\varepsilon,x)} \bigg| \frac{ H-H_{1}}{\eps}\bigg| \, \d t \d y\bigg) \, \d x
\end{align*}
and thanks to \eqref{eq:Hinfty} we get for $\eps \leq 1$ 
\begin{align*}
	J_{\varepsilon, M}  \lesssim  \int_{|x| \le 3M} U_\eps(x) \, \d x.
\end{align*}
The assumption \eqref{eq:aecv0} together with the pointwise bound $U_\eps \leq \NT(u)$ and the dominated convergence theorem yield again $J_{\varepsilon, M}  \to 0$ as $\eps \to 0$.

\medskip

\noindent \emph{The terms $I_{M}$ and $J_{M}$}. We have to estimate $I_M$ and $J_M$ in \eqref{eq: IM and JM}. Once again, we intervene before introducing the maximal operator in \eqref{eq: Ktau estimate generic} and simply use \eqref{eq:Hinfty}. In this way we obtain
\begin{align*}
K_{\tau, M}  \lesssim (\tau \wedge \tau^{-\gamma}) \|\NT(u)\|_1.
\end{align*}
Since the right-hand side is summable for $\tau = \tau_j $, $j \in \IZ$, we conclude $	M I_M  \lesssim \|\NT(u)\|_1$. Thus, we have $I_M \to 0$ in the limit as $M \to \infty$. For $M J_M$ we obtain the same type of bound by arguing as for $K_{\tau,M}$ with $\tau = M$.

At this point we have handled the terms in \eqref{eq:firstuG} and the argument at the end of Case~2 for the regularity problem explains why our proof automatically covers the terms in \eqref{eq:seconduG}.
\subsection*{Case 3: \texorpdfstring{$\boldsymbol{{p_{+}(L) \leq p<p_+(L)^*}}$}{p+ < p <p+*}}

We fine-tune the strategy in Case~2 for the regularity problem. Once again, working under the general assumptions \eqref{eq:aecv0} and \eqref{eq: G-alpha-assumption} does not pose any additional difficulty. However, the range of admissible exponents now changes with the parameter $\alpha$ in \eqref{eq: G-alpha-assumption} and we need to assume
\begin{align}
\label{eq: G-alpha exponents}
0 < \frac{1}{p_+(L)} - \frac{1}{p} < \frac{1-\alpha}{n}.
\end{align}
For the Dirichlet problem $(D)_p^\Le$ we have $\Theta = \NT(u)$ and $\alpha = 0$, so that this is the range that we are aiming at.

\medskip

\noindent \emph{The term $J_{\varepsilon, M}$}. By \eqref{eq: JepsM generic} we have
\begin{align*}
J_{\eps,M} \lesssim  \int_{|x| \le 3M} U_\eps(x) \bigg(\bariint_{ W(\varepsilon,x)} \bigg| \frac{ H-H_{1}}{\eps}\bigg| \, \d t \d y\bigg) \, \d x.
\end{align*}
We introduce the maximal function $\Max^\eps$ restricted to balls with radii not exceeding $\varepsilon$ as before and use the Cauchy--Schwarz inequality to give
\begin{align*}
J_{\eps,M} 	
& \lesssim \bigg(\int_{|y|\le 4M} \Max^\varepsilon(U_\eps)^{2} \, \d y \bigg)^{\frac{1}{2}} \bigg(\int_{\R^n} \bigg(\barint_{\varepsilon/2}^{2\varepsilon}  \bigg| \frac{ H-H_{1}}{t}\bigg| \, \d t \bigg)^{2} \d y \bigg)^{\frac{1}{2}} \\
& \lesssim \bigg(\int_{|y|\le 4M} \Max^\varepsilon(U_\eps)^{2} \, \d y\bigg)^{\frac{1}{2}} \bigg(\int_{\R^n} \barint_{\varepsilon/2}^{2\varepsilon}  \bigg| \frac{ H-H_{1}}{t}\bigg|^{2}\, \d t \d y \bigg)^{\frac{1}{2}} \\
& \lesssim \bigg(\int_{|y|\le 5M} U_\eps^{2} \, \d y\bigg)^{\frac{1}{2}} \bigg(\barint_{\varepsilon/2}^{2\varepsilon} \int_{\R^n}  \bigg| \frac{ H-H_{1}}{t}\bigg|^{2}\, \d y \d t \bigg)^{\frac{1}{2}},
\end{align*} 
where we have used $\Max^\eps(U_\eps) \leq \Max(\ind_{B(0,5M)}U_\eps)$ on $B(0,4M)$ and the maximal theorem in the last line. The second integral on the right is uniformly bounded in $\eps \leq 1$ by Lemma~\ref{lem:H} applied with $r=2$ and we are left with
\begin{align}
\label{eq:uniqCase4generic}
	J_{\varepsilon, M} \lesssim \bigg(\int_{|y|\le 5M} U_\eps^{2} \, \d y\bigg)^{\frac{1}{2}} \qquad (\eps \leq 1),
\end{align}
so far under the mere assumption that $u$ is a weak solution to $\Le u = 0$ in $\reu$. H\"older's inequality yields
\begin{align*}
	J_{\varepsilon, M} \lesssim M^{\frac{n}{2}-\frac{n}{p}} \bigg(\int_{\R^n} U_\eps^{p} \, \d y\bigg)^{\frac{1}{p}},
\end{align*}
which goes to $0$ as $\eps \to 0$, using \eqref{eq:aecv0}, the pointwise bound $U_\eps \leq \Theta$ and the dominated convergence theorem.

\medskip

\noindent \emph{The terms $I_{M}$ and $J_{M}$}. We are left with treating the terms $I_M$ and $J_M$ in \eqref{eq: IM and JM}. To this end, we recall the generic decomposition from \eqref{eq: MIM bound using Ktau} and \eqref{eq: Ktau estimate generic}: 
\begin{align*}
	MI_{M}\le \sum_{j=-\infty}^{j_{M}}  K_{\tau_{j},M},
\end{align*}
where $\tau_{j} = (\nicefrac{9}{4})^j$, $j_{M}$ is the unique integer with $\tau_{j_{M} -1}\le M<\tau_{j_{M}}$ and for $\tau \leq \nicefrac{9M}{4}$,
\begin{align*}
	K_{\tau, M}  \lesssim \tau^\alpha \int_{|y|\le 5 M} \Max(\Theta)(y) \bigg(\barint_{\tau/2}^{2\tau}  | H-H_{1}| \, \d t \bigg) \d y.
\end{align*} 
With a choice of $(p_-(L) \vee 1) < r < p_+(L)$ that we will specified later on, we obtain
\begin{align*}
	K_{\tau, M}  
	&\lesssim \tau^\alpha \bigg(\int_{\R^n} \Max(\Theta)(y)^p \, \d y\bigg)^{\frac{1}{p}}  \bigg(\int_{|y|\le 5 M} \bigg|\barint_{\tau/2}^{2\tau}  | H-H_{1}| \, \d t\bigg|^{p'} \d y\bigg)^{\frac{1}{p'}} \\
	&\lesssim \tau^\alpha M^{\frac{n}{p'}-\frac{n}{r'}} \bigg(\int_{\R^n} \Theta(y)^p \, \d y\bigg)^{\frac{1}{p}}  \bigg(\int_{|y|\le 5 M} \bigg|\barint_{\tau/2}^{2\tau}  | H-H_{1}| \, \d t\bigg|^{r'} \d y\bigg)^{\frac{1}{r'}} \\
	&\leq \tau^\alpha M^{\frac{n}{p'}-\frac{n}{r'}} \|\Theta\|_p \bigg(\barint_{\tau/2}^{2\tau} \int_{\R^n} |H-H_{1}|^{r'} \, \d y  \d t \bigg)^{\frac{1}{r'}}\\
	&\lesssim \tau^{\alpha}(\tau \wedge \tau^{-\gamma}) M^{\frac{n}{p'}-\frac{n}{r'}} \|\Theta\|_p \\
	&= \tau^{\alpha}(\tau \wedge \tau^{-\gamma}) M^{\frac{n}{r}-\frac{n}{p}} \|\Theta\|_p.
\end{align*}
We have used H\"older's inequality in the first line, the maximal theorem and again H\"older's inequality in the second one, Jensen's inequality in the third one and Lemma~\ref{lem:H} with exponent $r'$ in the fourth one. The exponent $\gamma$ is positive and depends on $r$. Summing up the estimates for $\tau = \tau_j$ yields
\begin{align*}
	M I_M \lesssim M^{\frac{n}{r}-\frac{n}{p}} (1+M^{\alpha - \gamma})\|\Theta\|_p.
\end{align*}
The assumption \eqref{eq: G-alpha exponents} guarantees that we can pick $r$ such that $\nicefrac{1}{r} - \nicefrac{1}{p} < \nicefrac{(1-\alpha)}{n}$. In this case $M$ appears with exponent smaller than $1$ on the right-hand side and we conclude $I_M \to 0$ in the limit as $M \to \infty$.

For $J_{M}$ we have the bound
\begin{align*}
M J_M \lesssim M^\alpha \int_{|y|\le 5 M} \Max(\Theta)(y) \bigg(\barint_{\tau/2}^{2\tau}  | H-H_{1}| \, \d t \bigg) \d y,
\end{align*}
see \eqref{eq: JM estimate generic}. The steps above with $\tau = M$ yield the same bound
\begin{align*}
M J_{M}\lesssim M^{\frac{n}{r}-\frac{n}{p}} (1+M^{\alpha - \gamma})\|\Theta\|_p,
\end{align*}
from which we conclude $J_M \to 0$ as $M\to \infty$.

\medskip We have handled the terms in \eqref{eq:firstuG} and once again the discussion at the the end of Case~2 for the regularity problem explains why our proof automatically covers the terms in \eqref{eq:seconduG}. This completes the proof of Theorem~\ref{thm: blockdir}.
\subsection{Uniqueness for \texorpdfstring{$\boldsymbol{(\wt D)_{\Lamdot^\alpha}^\Le}$}{(D)_a}}
\label{subsec: uniqueness Holder Dirichlet tilde}

We turn to the situation when $p_{+}(L)>n$ and prove that solutions to $(\wt D)_{\Lamdot^\alpha}^\Le$ are unique in the range of exponents $0\le \alpha< 1-\nicefrac{n}{p_{+}(L)}$. Hence, we assume \eqref{eq:aecv0} and $\NTsharpalpha(u) \in \L^\infty$. The control of the sharp functional means that we have
\begin{align}
\label{eq: G-alpha-assumption infinity}
	U_t(x) \leq t^\alpha \NTsharpalpha(u)(x) \quad ((t,x) \in \reu),
\end{align}
which is an assumption of the same type as \eqref{eq: G-alpha-assumption} but for $p = \infty$. Fortunately, this only requires a slight modification of the generic argument in the previous section.

\medskip

\noindent \emph{The term $J_{\varepsilon,M}$}. According to \eqref{eq:uniqCase4generic} we have
\begin{align*}
	J_{\varepsilon, M} 
	\lesssim \bigg(\int_{|x|\le 5M} U_\eps^{2} \, \d y\bigg)^{\frac{1}{2}}
\end{align*} 
and \eqref{eq: G-alpha-assumption infinity} still allows us to use the dominated convergence theorem when passing to the limit as $\eps \to 0$. Hence, $J_{\eps,M} \to 0$ follows.

\medskip

\noindent \emph{The terms $I_{M}$ and $J_{M}$}. For the terms in \eqref{eq: IM and JM} we start out with the usual decomposition from \eqref{eq: MIM bound using Ktau} and the estimate before \eqref{eq: Ktau estimate generic}: 
\begin{align*}
	MI_{M}\le \sum_{j=-\infty}^{j_{M}}  K_{\tau_{j},M},
\end{align*}
where $\tau_{j} = (\nicefrac{9}{4})^j$, $j_{M}$ is the unique integer with $\tau_{j_{M} -1}\le M<\tau_{j_{M}}$ and for $\tau \leq \nicefrac{9M}{4}$,
\begin{align}
\label{eq: Ktau generic 2}
 	K_{\tau, M}  \lesssim  \int_{\frac{M}{2} \le |x|\le \frac{5M}{2}} U_\tau(x) \bigg(\bariint_{ W(\tau,x)} | {H-H_{1}}| \bigg) \, \d x.
\end{align}
We use \eqref{eq: G-alpha-assumption infinity}, H\"older's inequality and Lemma~\ref{lem:H} with an exponent $(p_-(L^\sharp) \vee 1) <  r' < p_+(L^\sharp)$ to be specified yet in order to give
\begin{align*}
	K_{\tau, M}    
	& \lesssim \tau^\alpha \|\NTsharpalpha(u)\|_\infty \int_{|y|\le 5 M} \bigg(\barint_{\tau/2}^{2\tau}  | H-H_{1}| \, \d t \bigg) \d y \\
	&\lesssim \tau^\alpha \|\NTsharpalpha(u)\|_\infty M^{\frac{n}{r}} \barint_{\tau/2}^{2\tau}  \bigg( \int_{\R^n} |H-H_{1}|^{r'} \, \d y \bigg)^{\frac{1}{r'}}  \, \d t\\
	& \lesssim \tau^\alpha (\tau \wedge \tau^{-\gamma}) M^{\frac{n}{r}} \|\NTsharpalpha(u)\|_\infty,
\end{align*}
where $\gamma > 0$ depends on $r$. Summing up the estimates for $\tau = \tau_j$ leads to
\begin{align*}
	M I_M \lesssim M^{\frac{n}{r}} (1+M^{\alpha-\gamma}) \|\NTsharpalpha(u)\|_\infty.
\end{align*}
By assumption on $\alpha$ we can pick $n < r < p_+(L)$ such that $\alpha < 1 -\nicefrac{n}{r}$. Then the exponent for $M$ on the right-hand side becomes smaller than $1$ and $I_M \to 0$ in the limit as $M \to \infty$ follows.
 
For $J_{M}$ we recall from \eqref{eq: JM estimate generic} the bound
\begin{align*}
	MJ_M \lesssim \int_{|x|\le \frac{5M}{2}}  U_M(x) \bigg(\bariint_{ W(M,x)} | {H-H_{1}}| \bigg) \, \d x
\end{align*}
and the previous argument for $\tau = M$ yields $J_M \to 0$ in the limit as $M \to \infty$.

\medskip At this point we have handled the terms in \eqref{eq:firstuG} and as in the earlier steps the limits for the terms in \eqref{eq:seconduG} come for free.
\subsection{Uniqueness for \texorpdfstring{$\boldsymbol{(D)_{\Lamdot^\alpha}^\Le}$}{(D)_a} -- conclusion of the proof of Theorem~\ref{thm: Holder-dir}}
\label{subsec: uniqueness Holder Dirichlet}

We turn to uniqueness of solutions to the Dirichlet problem $(D)_{\Lamdot^\alpha}^{\Le}$ with interior Carleson control. We work under the same assumptions
$p_{+}(L)>n$ and $0\le \alpha< 1-\nicefrac{n}{p_{+}(L)}$ as in the previous section.

The case $\alpha > 0$ is particularly simple. We merely need the following general lemma to compare several functionals that all measure smoothness of order $\alpha-1$. 

\begin{lem}
\label{lem: Calpha controls Zalpha-1}
Let $\alpha \in (0,1)$. There is a dimensional constant $\omega_n$ such that for all $u \in \Wloc^{1,2}(\reu)$,
\begin{align*}
	\|W(t^{1-\alpha} \nabla u) \|_\infty \leq \omega_n 2^\alpha \|C_\alpha(t \nabla u)\|_\infty \leq  \omega_n 2^\alpha \|C_0(t^{1-\alpha} \nabla u)\|_\infty.
\end{align*}
Moreover, if \eqref{eq:aecv0} holds, then
\begin{align*}
	\|\NTsharpalpha(u)\|_\infty \lesssim \|W(t^{1-\alpha} \nabla u) \|_\infty.
\end{align*}
\end{lem}

\begin{proof}
For the first claim we simply note that for any $F \in \Lloc^2(\reu)$,
\begin{align*}
	\bigg(\bariint_{W(t,x)} &|s^{1-\alpha} F(s,y)|^2 \, \d s \d y\bigg)^{\frac{1}{2}} \\
	&\leq \frac{\omega_n 2^\alpha}{t^\alpha} \bigg(\int_0^{2t} \barint_{B(x,2t)} |sF(s,y)|^2 \, \frac{\d y \d s}{s}\bigg)^{\frac{1}{2}} \\
	&\leq \omega_n 2^\alpha \bigg(\int_0^{2t} \barint_{B(x,2t)} |s^{1-\alpha} F(s,y)|^2 \, \frac{\d y \d s}{s}\bigg)^{\frac{1}{2}},
\end{align*}
and taking suprema in $t$ and $x$ on both sides yields the claim. Under assumption \eqref{eq:aecv0} the second claim follows from the trace theorem in Proposition~\ref{prop: NT trace Y}.(iii).
\end{proof}

To prove uniqueness of solutions to $(D)_{\Lamdot^\alpha}^\Le$ when $\alpha>0$, we assume $C_\alpha(t\nabla u) \in \L^\infty$ and that \eqref{eq:aecv0} holds.  Lemma~\ref{lem: Calpha controls Zalpha-1} yields $\NTsharpalpha(u) \in \L^\infty$ and under this weaker assumption we have already shown $u=0$ in the previous section. 

It remains to treat the $\BMO$ Dirichlet problem $(D)_{\Lamdot^0}^\Le$. We assume therefore $C_0(t \nabla u) \in \L^\infty$ and for the first time \eqref{eq:L2cv}. We implement the strategy of Section~\ref{subsec: review of strategy uniqueness} with $R = M$ and first send $\eps \to 0$ and then $M \to \infty$ in \eqref{eq:firstuG} and \eqref{eq:seconduG}. 

\medskip

\noindent \emph{The terms $J_{\varepsilon,M}$ and $\wt J_{\eps, M}$}. The Cauchy--Schwarz inequality yields
\begin{align}
\label{eq: Jeps Holder}
J_{\varepsilon,M}
\lesssim  \bigg(\int_{|x|\le 2M} \barint_{2\varepsilon/3}^{3\varepsilon/2} |u|^2 \, \d t \d x \bigg)^{\frac{1}{2}} \bigg(\int_{\R^n} \barint_{2\varepsilon/3}^{3\varepsilon/2}  |\pd_{t} (H-H_{1})|^2\, \d t \d x\bigg)^{\frac{1}{2}}.
\end{align} 
By covering $B(0,2M)$ up to a set of measure zero by pairwise disjoint cubes $Q_k$ of sidelength $\varepsilon$ with $2Q_k \subseteq B(0, 2M+1)$ and using reverse H\"older inequalities for $u$, we obtain
\begin{align*}
\int_{B(0,2M)} \barint_{2\varepsilon/3}^{3\varepsilon/2} |u|^2 \, \d t \d x
& \lesssim \sum_{k} |Q_{k}|   \barint_{Q_{k}}\barint_{2\varepsilon/3}^{3\varepsilon/2} |u|^2 \, \d t \d x \\
&   \lesssim  \sum_{k} |Q_{k}|   \bigg(\barint_{2 Q_{k}}\barint_{\varepsilon/2}^{2\varepsilon} |u| \, \d t \d x\bigg)^2
\\
&   \lesssim  \sum_{k} |Q_{k}|  \barint_{2 Q_{k}}  \bigg(\barint_{\varepsilon/2}^{2\varepsilon} |u|\, \d t \bigg)^2 \d x
\\
&   \leq    \int_{B(0,2M+1)}  \bigg(\barint_{\varepsilon/2}^{2\varepsilon} |u| \, \d t \bigg)^2 \d x.
\end{align*}
By assumption \eqref{eq:L2cv}, this integral tends to $0$ as $\varepsilon\to 0$.  As for the term with $H-H_{1}$ in \eqref{eq: Jeps Holder}, we use Lemma~\ref{lem:H} with $r =2$ to deduce a uniform bound in $\eps \in (0,1)$.  

The estimate for $\wt J_{\varepsilon,M}$ is very similar. Indeed, $t\pd_{t}u$ is handled via the same argument and incorporating the Caccioppoli inequality, whereas for $\nicefrac{(H-H_{1})}{t}$ we use Lemma~\ref{lem:H} again. 

\medskip

\noindent \emph{The terms $I_{M}$, $\wt{I}_{M}$ and $J_{M}$, $\wt{J}_{M}$}. We estimate the terms in \eqref{eq: IM and JM}. Only one change to the corresponding argument for $(\wt D)_{\Lamdot^0}^\Le$ in Section~\ref{subsec: uniqueness Holder Dirichlet tilde} will be necessary. In particular, the estimates for the tilde terms that correspond to \eqref{eq:seconduG} come again for free.

The argument for $I_M$ with $\alpha =0$ in the previous section uses the interior control only once, namely to bound $U_\tau(x)$ in \eqref{eq: Ktau generic 2} uniformly by $ \|\NTsharpalpha(u)\|_\infty$. This bound is not available under our current assumption but the following lemma provides a substitute that still suits our purpose.

\begin{lem}
\label{lem: Whitney average bound under Carleson}
If $v \in \Wloc^{1,2}(\reu)$ is such that $C_0(t \nabla v) \in \L^\infty(\R^n)$, then
\begin{align*}
	\bariint_{W(t,x)} |v| \, \d s \d y \lesssim 1+|\ln(t)| + \ln(1+|x|) \quad ((t,x) \in \reu),
\end{align*}
where the implicit constant also depends on $v$.
\end{lem}

We defer the proof and use Lemma~\ref{lem: Whitney average bound under Carleson} to bound $U_\tau(x)$ in \eqref{eq: Ktau generic 2}. This yields an additional factor $(1 + |\ln(\tau)| + \ln(M))$ compared to the estimates in the previous section and hence we obtain
\begin{align*}
	K_{\tau, M}  \lesssim (1 + |\ln(\tau)| + \ln(M)) M^{\frac{n}{r}} (\tau \wedge \tau^{-\gamma})
\end{align*}
with $r>n$ and $\gamma > 0$. Summing up the estimates for $\tau = \tau_j$ yields
\begin{align*}
	MI_{M}\le \sum_{j=-\infty}^{j_{M}}  K_{\tau_{j},M}  \lesssim \ln(M) M^{\frac{n}{r}},
\end{align*}
which still implies that $I_M$ tends to $0$ as $M \to \infty$.

Likewise, using Lemma~\ref{lem: Whitney average bound under Carleson} to control $U_M(x)$ in \eqref{eq: JM estimate generic} leads to $M J_{M} \lesssim \ln(M) M^{\frac{n}{r}}$ and we conclude as before. The proof of Theorem~\ref{thm: Holder-dir} is complete modulo the following:

\begin{proof}[Proof of Lemma~\ref{lem: Whitney average bound under Carleson}]
Set $w \coloneqq |v|$, which satisfies the same assumptions. Suppose that $W_j = W(t_j,x_j)$ and $W_k = W(t_k,x_k)$ are two Whitney regions with non-empty intersection and suppose that $t_j \leq t_k$. Then $W_j$ and $W_k$ are comparable in measure and the cylinder $W \coloneqq (\nicefrac{t_j}{2}, 8 t_j) \times B(x_k,8t_j)$  contains both $W_j$ and $W_k$. Hence, we can use Poincaré's inequality in order to give
\begin{align*}
	|(w)_{W_j} - (w)_{W_k}|
	&\lesssim \bariint_{W} |w - (w)_{W_k}| \, \d t \d x \\
	&\lesssim \bariint_{W} |t \nabla w| \, \d t \d x \\
	&\lesssim \|C_0(t \nabla w)\|_\infty
\end{align*}
with a implicit constant that depends only on $n$. If $W_1,\ldots, W_k$ is a chain of Whitney regions with the property that each region intersects its successor, then a telescopic sum yields
\begin{align*}
	|(w)_{W_1} - (w)_{W_k}| \lesssim k \|C_0(t \nabla w)\|_\infty.
\end{align*}
We write $W_1 \to W_k$ in that case.

Now, we fix $(t,x) \in \R^n$. Since $w$ is locally integrable, it suffices to construct a chain $W(t,x) \to W(1,0)$ of length controlled by $1 + |\ln(t)| + \ln(1+|x|)$. One possible construction is as follows. Successively halving or doubling $t$, we obtain a chain $W(t,x) \to W(1,x)$ of length comparable to $1 + |\ln(t)|$. If $|x| < 1$, then $W(1,x)$ and $W(1,0)$ intersect and we are done. If $|x| \geq 1$, then in the same manner we obtain chains $W(1,x) \to W(2|x|,x)$ and $W(1,0) \to W(2|x|,0)$ of length comparable to $\ln(1+|x|)$. Moreover, $W(2|x|,x)$ and $W(2|x|,0)$ intersect.
\end{proof}
\subsection{Uniqueness for \texorpdfstring{$\boldsymbol{(D)_{\Xdot^{s,p}}^\Le}$}{(D)_Xsp} -- conclusion of the proof of Theorem~\ref{thm: uniqueness fractional}}
\label{subsec: uniqueness fractional}

The last uniqueness result concerns the problems $(D)_{\Xdot^{s,p}}$ with fractional regularity data. As usual, $\X$ denotes $\B$ or $\H$ and $\Y$ is the corresponding solution space of type $\Z$ or $\T$. Figure~\ref{fig: diagram-uniqueness-general} and Figure~\ref{fig: diagram-uniqueness-p+large} show the regions of exponents that we are aiming at in an $(\nicefrac{1}{p},s)$-plane. In the previous sections we have already obtained uniqueness on the bottom and top segments. 

\begin{figure}[h]
\begin{center}
\begin{tikzpicture}[scale=2.2]

\newcommand\fracspace{\vphantom{\frac{1}{1}}};
\def\dimension{6};
\def\xlength{3.5};
\pgfmathsetmacro\xstretch{\xlength/(1+1/\dimension)}; 

\pgfmathsetmacro\UpSob{\xstretch*(1/2-1/\dimension)-0.2};
\pgfmathsetmacro\UpSobStar{\xstretch*(\UpSob/\xstretch-1/\dimension)};
\pgfmathsetmacro\LowSob{\xstretch*(1/2+1/\dimension)+0.5};
\pgfmathsetmacro\LowSobStar{\xstretch*(\LowSob/\xstretch+1/\dimension)};
\pgfmathsetmacro\UpGradient{\xstretch*(1/2)-0.2};
\pgfmathsetmacro\UpGradientDual{\xstretch*(1/2)+0.3};
\pgfmathsetmacro\Half{\xstretch*0.5};
\pgfmathsetmacro\One{\xstretch*1};
\coordinate (P00) at (1+\UpGradient,2);
\coordinate (P01) at (1+\LowSob,2);
\coordinate (P10) at (1+\UpSob,0);
\coordinate (P10shift) at (1+\UpSob,2);
\coordinate (P11) at (1+\UpGradientDual,0);
\coordinate (P12) at (1+\LowSob,0);
\coordinate (DirExtra) at (1+\UpSobStar,0);
\coordinate (RegExtra) at (1+\LowSobStar,2);

\draw [thin] (1,2) -- (1+\xlength,2); 
\draw [thin] (1,0) -- (1+\xlength,0); 

\draw [thick,->] (1,-0.5) -- (1+\xlength+0.2,-0.5);
\draw [fill=black] (1+\Half,-0.5) circle [radius = .5pt];
\node [below] at (1+\Half,-0.5) {$\frac{1}{2 \fracspace}$};
\draw [fill=black] (1+\UpGradient,-0.5) circle [radius = .5pt];
\node [below] at (1+\UpGradient,-0.5) {$\frac{1}{q_+^L\fracspace}$};
\draw [fill=black] (1+\UpGradientDual,-0.5) circle [radius = .5pt];
\node [below] at (1+\UpGradientDual,-0.5) {$\frac{1}{(q_+^{L^\sharp})'\fracspace}$};
\draw [fill=black] (1+\UpSob,-0.5) circle [radius = .5pt];
\node [below] at (1+\UpSob,-0.5) {$\frac{1}{p_+^L\fracspace}$};
\draw [fill=black] (1+\LowSob,-0.5) circle [radius = .5pt];
\node [below] at (1+\LowSob,-0.5) {$\frac{1}{p_-^L \vee 1 \fracspace}$};
\draw [fill=black] (1+\UpSobStar,-0.5) circle [radius = .5pt];
\node [below] at (1+\UpSobStar,-0.5) {$\frac{1}{(p_+^L)^*\fracspace}$};
\draw [fill=black] (1+\LowSobStar,-0.5) circle [radius = .5pt];
\node [below] at (1+\LowSobStar,-0.5) {$\frac{1}{{(p_-^L)_* \vee 1_*\fracspace}}$};
\draw [fill=black] (1+\xlength,-0.5) circle [radius = .5pt];
\node [below] at (1+\xlength,-0.5) {$\frac{n+1}{n \fracspace}$};
\node [right] at (1+\xlength+0.2,-0.5) {$\frac{1}{p \fracspace}$};
\draw [fill=black] (1,-0.5) circle [radius = .5pt];
\node [below] at (1,-0.5) {$0$};

\draw [thick,->] (0.7,0) -- (0.7,2.2);
\node [above] at (0.7,2.2) {$s$};
\draw [fill=black] (0.7,2) circle [radius = .5pt];
\node [left] at (0.7,2) {$1$};
\draw [fill=black] (0.7,0) circle [radius = .5pt];
\node [left] at (0.7,0) {$0$};

\draw [thin,dotted] (1,2) -- (1,0); 
\draw [thin,dotted] (1+\LowSobStar,2) -- (1+\LowSobStar,0); 
\draw [thin,dotted] (1+\UpSobStar,2) -- (1+\UpSobStar,0);
\draw [thin,dotted] (1+\LowSob,2) -- (1+\LowSob,0); 
\draw [thin,dotted] (1+\UpSob,2) -- (1+\UpSob,0);
\draw [thin,dotted] (1+\xlength,2) -- (1+\xlength,0); 
	
\path [fill=green!40!black, opacity = 0.3] (P10shift)--(RegExtra)--(P12)--(DirExtra)--(P10shift);	
\draw[ultra thick, green!40!black] (P10shift)--(RegExtra);
\draw[ultra thick, green!40!black] (P12)--(DirExtra);
\draw [fill=white] (P10shift) circle [radius = .75pt];
\draw [fill=white] (RegExtra) circle [radius = .75pt];
\draw [fill=white] (P12) circle [radius = .75pt];
\draw [fill=white] (DirExtra) circle [radius = .75pt];
\end{tikzpicture}
\end{center}
\caption{Exponents for uniqueness in Dirichlet and regularity problems in the case $p_+(L) \leq n$. Uniqueness holds on the open bottom and top segments (Sections~\ref{subsec: uniqueness Lp Dirichlet} and \ref{subsec: uniqueness regularity}) and the interior of the trapezoidal region (Section~\ref{subsec: uniqueness fractional}).
}
\label{fig: diagram-uniqueness-general}
\end{figure}
\begin{figure}[h]
\begin{center}
\begin{tikzpicture}[scale=2.2]
	\newcommand\fracspace{\vphantom{\frac{1}{1}}};
	\def\dimension{6};
	\def\xlength{3.5};
	\pgfmathsetmacro\xstretch{\xlength/(1+1/\dimension)}; 
	
	\pgfmathsetmacro\UpSob{\xstretch*(1/\dimension)-0.3};
	\pgfmathsetmacro\UpSobStar{\xstretch*(\UpSob/\xstretch-1/\dimension)};
	\pgfmathsetmacro\UpSobStarBesov{(-2)*\UpSobStar*\dimension/\xlength};
	\pgfmathsetmacro\LowSob{\xstretch*(1/2+1/\dimension)+0.5};
	\pgfmathsetmacro\LowSobStar{\xstretch*(\LowSob/\xstretch+1/\dimension)};
	\pgfmathsetmacro\UpGradient{\xstretch*(1/2)-0.2};
	\pgfmathsetmacro\UpGradientDual{\xstretch*(1/2)+0.3};
	\pgfmathsetmacro\Half{\xstretch*0.5};
	\pgfmathsetmacro\One{\xstretch*1};
	\coordinate (P00) at (1+\UpGradient,2);
	\coordinate (P01) at (1+\LowSob,2);
	\coordinate (P10) at (1+\UpSob,0);
	\coordinate (P10shift) at (1+\UpSob,2);
	\coordinate (P11) at (1+\UpGradientDual,0);
	\coordinate (P12) at (1+\LowSob,0);
	\coordinate (DirExtra) at (1,0);
	\coordinate (DirHolder) at (1.07+\UpSobStar,0); 
	\coordinate (DirHolderBesov) at (1, \UpSobStarBesov);
	\coordinate (RegExtra) at (1+\LowSobStar,2);
	\coordinate (P000) at (1-\xstretch/\dimension,0);

	\draw [thin] (0.3,2) -- (1+\xlength,2); 
	\draw [thin] (0.3,0) -- (1+\xlength,0); 
	
	\draw [thick,->] (0.3,-0.5) -- (1+\xlength+0.2,-0.5);
	\draw [fill=black] (1+\Half,-0.5) circle [radius = .5pt];
	\node [below] at (1+\Half,-0.5) {$\frac{1}{2 \fracspace}$};
	\draw [fill=black] (1+\UpGradient,-0.5) circle [radius = .5pt];
	\node [below] at (1+\UpGradient,-0.5) {$\frac{1}{q_+^L\fracspace}$};
	\draw [fill=black] (1+\UpGradientDual,-0.5) circle [radius = .5pt];
	\node [below] at (1+\UpGradientDual,-0.5) {$\frac{1}{(q_+^{L^\sharp})'\fracspace}$};
	\draw [fill=black] (1+\UpSob,-0.5) circle [radius = .5pt];
	\node [below] at (1+\UpSob,-0.5) {$\frac{1}{p_+^L\fracspace}$};
	\draw [fill=black] (1+\LowSob,-0.5) circle [radius = .5pt];
	\node [below] at (1+\LowSob,-0.5) {$\frac{1}{p_-^L \vee 1 \fracspace}$};
	\draw [fill=black] (1+\UpSobStar+0.07,-0.5) circle [radius = .5pt];
	\node [below] at (1+\UpSobStar,-0.5) {$\frac{1}{p_+^L\fracspace} \!- \!\frac{1}{n\fracspace}$};
	\draw [fill=black] (1+\LowSobStar,-0.5) circle [radius = .5pt];
	\node [below] at (1+\LowSobStar,-0.5) {$\frac{1}{{(p_-^L)_* \vee 1_*\fracspace}}$};
	\draw [fill=black] (1+\xlength,-0.5) circle [radius = .5pt];
	\node [below] at (1+\xlength,-0.5) {$\frac{n+1}{n \fracspace}$};
	\node [right] at (1+\xlength+0.2,-0.5) {$\frac{1}{p \fracspace}$};
	\draw [fill=black] (1,-0.5) circle [radius = .5pt];
	\node [below] at (1,-0.5) {$0$};
	\draw [fill=black] (0.3,-0.5) circle [radius = .5pt];
	\node [below] at (0.2,-0.5) {$-\frac{1}{n \fracspace}$};
	
	\draw [thick,->] (-0,0) -- (0,2.2);
	\node [above] at (0,2.2) {$s$};
	\draw [fill=black] (0,2) circle [radius = .5pt];
	\node [left] at (0, \UpSobStarBesov) {$\hspace{-30pt}1\!-\! \frac{n}{p_+^L}$};
	\draw [fill=black] (0,  \UpSobStarBesov) circle [radius = .5pt];	
	\node [left] at (0,2) {$1$};
	\draw [fill=black] (0,0) circle [radius = .5pt];
	\node [left] at (0,0) {$0$};

	\draw [thin,dotted] (1,2) -- (1,0); 
	\draw [thin,dotted] (0.3,2) -- (0.3,0); 
	\draw [thin,dotted] (1+\UpGradient,2) -- (1+\UpGradient,0);
	\draw [thin,dotted] (1+0.07+\UpSobStar,2) -- (1+0.07+\UpSobStar,0);
	\draw [thin,dotted] (P10shift) -- (DirHolder);
	\draw [thin,dotted] (1+\LowSob,2) -- (1+\LowSob,0);
	\draw [thin,dotted] (1+\xlength,2) -- (1+\xlength,0); 
	\draw [thin,dotted] (1+\UpSob,2) -- (1+\UpSob,0);
	
	\draw [thin,dotted] (0.3, \UpSobStarBesov) -- (1,\UpSobStarBesov); 
	
	\path [fill=green!40!black, opacity = 0.3] (P10shift)--(RegExtra)--(P12)--(DirExtra)--(DirHolderBesov)--(P10shift);	
	\draw[ultra thick, green!40!black] (P10shift)--(RegExtra);
	\draw[ultra thick, green!40!black] (P12)--(DirHolder);
	\draw[ultra thick, green!40!black] (DirExtra)--(DirHolderBesov);
	\draw [fill=white] (P10shift) circle [radius = .75pt];
	\draw [fill=white] (RegExtra) circle [radius = .75pt];
	\draw [fill=white] (P12) circle [radius = .75pt];
	\draw [fill=white] (DirHolder) circle [radius = .75pt];
	\draw [fill=green!40!black] (DirExtra) circle [radius = .75pt];
	\draw [fill=white] (DirHolderBesov) circle [radius = .75pt];
\end{tikzpicture}
\end{center}
\caption{Exponents for uniqueness in Dirichlet and regularity problems in the case $p_+(L) > n$. Uniqueness holds on the open bottom and top segments (Sections~\ref{subsec: uniqueness Lp Dirichlet}, \ref{subsec: uniqueness regularity} and \ref{subsec: uniqueness Holder Dirichlet}), the open vertical segment at $\nicefrac{1}{p}=0$ and the interior of the trapezoidal region (Section~\ref{subsec: uniqueness fractional}). Exponents with $\nicefrac{1}{p} < 0$ correspond to the spaces $\Lamdot^\alpha$ with $\alpha = -\nicefrac{n}{p}$ as usual.}
\label{fig: diagram-uniqueness-p+large}
\end{figure}

We distinguish four cases.

\begin{itemize}
	\item The rectangle $$(p_-(L) \vee 1) < p < p_+(L) \quad \& \quad  s \in (0,1),$$
	\item the left-hand triangle ($p_+(L) \leq n$) or trapezoid ($p_+(L) > n$) $$p_+(L) \leq p < p_+(L)^* \quad \& \quad \frac{1}{p_+(L)} - \frac{1}{p} < \frac{1-s}{n},$$
	\item $p_+(L) > n$ and the vertical segment $$p = \infty \quad \& \quad 0<s<1- \frac{n}{p_+(L)},$$
	\item the right-hand triangle $$(p_-(L)_* \vee 1_*) < p \leq (p_-(L) \vee 1) \quad \& \quad \frac{1}{p} - \frac{1}{p_-(L) \vee 1} < \frac{s}{n}.$$
\end{itemize}

In any case we assume \eqref{eq:aecv0} and $\|\nabla u\|_{\Y^{s-1,p}} < \infty$, which by definition of tent and $\Z$-spaces corresponds to one of the interior controls in Lemma~\ref{lem:prep}.
\subsection*{Case 1: The rectangle}

According to the trace theorem from Proposition~\ref{prop: NT trace Y}, there exists a function $\Theta \in \L^p$ such that
\begin{align*}
	U_t(x) \leq t^s\Theta(x) \quad ((t,x) \in \reu).
\end{align*}
Hence, \eqref{eq: G-alpha-assumption} holds with $\alpha = s$ and the general result from Case~2 for the regularity problem applies directly.
\subsection*{Case 2: The left-hand triangle or trapezoid}

Since we still work with finite exponents, assumption \eqref{eq: G-alpha-assumption} holds as in Case~1 with exponent $\alpha = s$. Thus, we can apply the general result from Case~3 for the Dirichlet problem provided that the exponents satisfy the respective assumption \eqref{eq: G-alpha exponents}. But this is exactly the restriction that defines this region.
\subsection*{Case 3: \texorpdfstring{$\boldsymbol{p_+(L) > n}$}{p+ > n} and the vertical segment}

Let $0< \alpha < 1- \nicefrac{n}{p_+(L)}$. We assume one of $C_0(t^{1-\alpha} \nabla u) \in \L^\infty$ or $W(t^{1-\alpha} \nabla u) \in \L^\infty$ and in any case that \eqref{eq:aecv0} holds at the boundary. Lemma~\ref{lem: Calpha controls Zalpha-1} yields $\NTsharpalpha(u) \in \L^\infty$ and under this weaker assumption we have already shown $u=0$ in Section~\ref{subsec: uniqueness Holder Dirichlet tilde}.
\subsection*{Case 4: The right-hand triangle}

The argument in Case~1 for the regularity problem implicitly contains a more general result that applies here. In view of the technicalities concerning the choice of exponents in that argument we have decided to stick with the version at regularity $s=1$ earlier on and here we provide the required generalization.  

We begin with the substitute for Lemma~\ref{lem: exponent regularity proof}.

\begin{lem} 
\label{lem: exponent fractional proof}
If $s \in (0,1)$ and $0<p<r \leq 2$, then for any weak solution $u$ to $\Le u=0$ on $\reu$,  
\begin{equation}
	\label{eq:pdtup* fractional}
	\bigg(\iint_{\reu} |\nabla u|^{r} t^{n(\frac{r}{p} - 1)+(1-s)r} \, \frac{\d t \d x}{t}\bigg)^{\frac{1}{r}} \lesssim \|\nabla u\|_{\Y^{s-1,p}}.
\end{equation}
Moreover, if \eqref{eq:aecv0} holds and in addition $p > \nicefrac{n}{n+s}$, then
\begin{equation}
	\label{eq:u/tp* fractional}
	\bigg(\iint_{\reu} | {u}|^{r} t^{n(\frac{r}{p} - 1)- sr} \, \frac{\d t \d x}{t} \bigg)^{\frac{1}{r}} \lesssim \|\nabla u\|_{\Y^{s-1,p}}.
\end{equation}
\end{lem}

\begin{proof}
Since $p<r$, we can use the mixed embedding for tent and $\Z$-spaces from \cite[Thm.~2.34]{AA} to the effect that $\Y^{s,p}\subseteq \Z^{\alpha,r}$ if $\alpha - s = n(\nicefrac{1}{r} - \nicefrac{1}{p})$. This means that
\begin{align*}
	\bigg( \iint_\reu W(t^{-\alpha} F)(\tau,y)^r\, \frac{\, \d \tau \d y}{\tau}\bigg)^{\frac{1}{r}} \lesssim \|F\|_{\Y^{s,p}}.
\end{align*}
As $r \leq 2$,  H\"older's inequality implies 
\begin{align}
\label{eq1: exponent fractional proof}
\bariint_{W(\tau,y)} |t^{-\alpha} F(t,x)|^r\, \d t \d x \leq W(t^{-\alpha} F)(\tau,y)^r
\end{align}
and applying the averaging trick backwards yields
\begin{align*}
	\bigg( \iint_\reu |t^{-\alpha} F(t,x)|^r\, \frac{\, \d t \d x}{t}\bigg)^{\frac{1}{r}} \lesssim \|F\|_{\Y^{s,p}}.
\end{align*}
If $F \coloneqq |t\nabla u|$, then $\|F\|_{\Y^{s,p}}= \|\nabla u\|_{\Y^{s-1,p}}$ and sorting out the exponent for $t$ on the left-hand side yields \eqref{eq:pdtup* fractional}. 

We have $\nicefrac{1}{r} - \nicefrac{\alpha}{n} = \nicefrac{1}{p} - \nicefrac{s}{n} < 1$ by assumption, which means that $r> \nicefrac{n}{n+\alpha}$. Since also $r \leq 2$, we can use part (i) of the trace theorem in Proposition~\ref{prop: NT trace Y} for $\nabla u \in \Z^{\alpha-1,r}$ with the same exponent $r$. Owing to \eqref{eq:aecv0}, we obtain
\begin{align*}
	\bariint_{W(\tau,y)} |t^{-\alpha} u(t,x)|^r \, \d t \d x \lesssim \Theta(y)^r
\end{align*}
for some function $\Theta$ with $\|\Theta\|_r \lesssim \|\nabla u\|_{\Z^{\alpha-1,r}}$. Integrating in $y$ and applying the averaging trick backwards yields
\begin{align*}
	\bigg(\iint_{\reu} |t^{-\alpha} u(t,x)|^r \, \frac{\, \d t \d x}{t} \bigg)^{\frac{1}{r}} \lesssim \|\nabla u\|_{\Z^{\alpha-1,r}} \lesssim \|\nabla u\|_{\Y^{s-1,p}}
\end{align*}
and as before $-\alpha r$ reveals itself as the same exponent than in the claim.
\end{proof}

Our standing assumption implies $\nicefrac{1}{p} < 1 + \nicefrac{s}{n}$. Hence, we have both parts of Lemma~\ref{lem: exponent fractional proof} at our disposal. This lemma allows us to control $\nabla u$ and $u$ in certain Lebesgue norms exactly as it was the case with Lemma~\ref{lem: exponent regularity proof}, except that we have different powers of $t$ to compensate: $t^{n(\frac{r}{p} - 1)+(1-s)r}$ and $t^{n(\frac{r}{p} - 1)-sr}$ replace $t^{n(\frac{r}{p} - 1)}$ and $t^{n(\frac{r}{p} - 1)-r}$, respectively, that is to say, we have an additional power $t^{(1-s)r}$. Armed with this observation, we pick again $(p_-(L) \vee 1) < r \leq 2$  and follow the proof in Case~1 of Section~\ref{subsec: uniqueness regularity} \emph{verbatim}. We only have to check that the additional power of $t$ still allows us to pass to the limits.

As for $I_{M,\varepsilon,R}$ and $\wt I_{M,\varepsilon,R}$, the different power of $t$ only changes the implicit constant that depends on $\eps, R$. Hence, these terms vanish when sending $M \to \infty$ as before. 

The estimates for $J_{\eps,M}$ and $\wt{J}_{\eps, M}$ are more delicate since now we obtain $\eps^{s-n(\frac{1}{p} - \frac{1}{r})}$ as factor in \eqref{eq: JepsM bound regularity} and \eqref{eq: JepsM-tilde bound regularity} if we want to control the respective integral on the right via Lemma~\ref{lem: exponent fractional proof}. We need to pick an admissible $r$ such that the exponent is non-negative. 

In dimension $n \geq 2$ we pick $r \coloneqq \frac {np}{n-ps}$ since then the exponent of $\eps$ vanishes. In particular, using also the restriction on $p$, we have
\begin{align*}
 \frac{1}{p} -\frac{1}{r} = \frac{s}{n} > \frac{1}{p} - \frac{1}{p_-(L) \vee 1},
\end{align*}
which in turn implies that $r > (p_-(L) \vee 1)$. On the other hand, $s \leq 1$ implies $r \leq p^*$ and at the same time we have $p \leq (p_-(L) \vee 1) \leq 2_*$ by Proposition~\ref{prop: J(L) contains neighborhood of Sobolev conjugates}. Thus, $r \leq 2$ and we conclude that $r$ is admissible.

In dimension $n=1$, Proposition~\ref{prop: J(L) contains neighborhood of Sobolev conjugates} yields $p_-(L) = \nicefrac{1}{2} = 1_*$. Hence, our assumption on $p$ is $\nicefrac{1}{(s+1)} < p \leq 1$ and this allows us to pick $r>1$ sufficiently close to $1$ such that $\nicefrac{1}{r} \geq \nicefrac{1}{p} - s$. Consequently, the exponent for $\eps$ in \eqref{eq: JepsM bound regularity} is non-negative and we conclude as before.

Likewise, we obtain for $J_{R,M}$ and $\wt{J}_{R, M}$ the new factor $R^{(s-1)-\gamma - n(\frac{1}{p} - \frac{1}{r})}$ in \eqref{eq: JRM bound regularity} and \eqref{eq: JRM-tilde bound regularity}. The exponent is negative since we have $s<1$, $\gamma > 0$ and $r>p$. This completes the proof.
\section{The Neumann problem}
\label{sec: Neumann problems}

\noindent In this final chapter we are concerned with the Neumann problem. In particular, we shall give the proof of Theorem~\ref{thm: blockneu}. We begin by recalling the construction of energy solutions to the Neumann problem. We use again the energy space $\Wdot^{1,2}(\reu)$ from Section~\ref{subsec: Energy solutions}. 

If $u \in \Wdot^{1,2}(\reu)$ is a weak solution to $\Le u = 0$ in $\reu$, then there exists a unique element $\dnuA u(0, \cdot) \in \Hdot^{-\nicefrac{1}{2},2}$ such that
\begin{align*}
\iint_{\reu} A \nabla u \cdot \cl{\nabla \phi} \, \d t \d x = - \langle \dnuA u(0, \cdot), \phi(0,\cdot) \rangle \quad (\phi \in \Wdot^{1,2}(\reu)),
\end{align*}
where on the right-hand side we use the duality pairing between $\Hdot^{-\nicefrac{1}{2},2}$ and $\Hdot^{\nicefrac{1}{2},2}$. Indeed, by assumption on $u$ and Lemma~\ref{lem: energy zero trace}, the left-hand side is a bounded anti-linear functional on $\Wdot^{1,2}(\reu)$ that vanishes whenever $\phi(0,\cdot) = 0$ and therefore it defines a bounded anti-linear functional on the trace space $\Hdot^{\nicefrac{1}{2},2}$. We call $\dnuA u(0, \cdot)$ the (inward pointing) \emph{conormal derivative}\index{conormal derivative} of $u$ at the boundary. 

\begin{prop}
\label{prop: existence of energy solution Neumann}
For all $f \in \Hdot^{-\nicefrac{1}{2},2}$ there exists a unique solution $u$ (modulo constants) to the problem
\begin{equation*}
\begin{cases}
\Le u=0   & (\text{in }\reu), \\
\nabla u \in \L^2(\reu),   \\
\dnuA u(0,\cdot) = f & (\text{in } \Hdot^{-\nicefrac{1}{2},2}).
\end{cases}
\end{equation*}
Moreover, $\|\nabla u\|_2 \lesssim  \|f\|_{\Hdot^{-1/2,2}}$ and $\lim_{t \to \infty} u(t,\cdot) = 0$ in $ \Hdot^{\nicefrac{1}{2},2}$.
\end{prop}

\begin{proof}
This is just the Lax--Milgram lemma applied in $\Wdot^{1,2}(\reu)$. The limit at $t=\infty$ follows from Lemma~\ref{lem: energy trace}.
\end{proof}

In the situation above we call $u$ the \emph{energy solution}\index{solution!energy with Neumann datum} to $\Le u = 0$ in $\reu$ with Neumann data $f$. Much alike to Section~\ref{subsec: Energy solutions} the energy solution coincides with the Poisson semigroup extension for suitable data. Throughout, we use the (extension to an) isomorphism $L^{1/2}: \Wdot^{1,2} \to \L^2$ with inverse $L^{-1/2}$. By duality and similarity we also obtain an (extension to an) isomorphism $a L^{1/2}: \L^2 \to \Wdot^{-1,2}$.

\begin{prop}
\label{prop: Neumann weak solution}
If $f \in \L^2 \cap  \Wdot^{-1,2}$, then the energy solution with Neumann data $f$ is given by $u(t,x)  = -\e^{-t {L^{1/2}}}(aL^{1/2})^{-1}f(x)$.\index{solution!semigroup}\index{solution!compatible}  
\end{prop}

\begin{proof} 
Set $g\coloneqq - (aL^{1/2})^{-1}f$. Then $g \in \Wdot^{1,2} \cap \L^2$ and, by interpolation, $g \in \Hdot^{\nicefrac{1}{2},2}$. It follows from Proposition~\ref{prop: smg solution is compatible} that  $u(t,x) \coloneqq \e^{-t {L^{1/2}}}g(x)$ is an energy solution to $\Le u = 0$ in $\reu$. In order to determine its Neumann datum, we let $\phi \in \C_0^\infty(\ree)$. By the functional calculus on $\L^2$ we have $u \in \C^1([0,\infty); \L^2)$ with $a \partial_t u(0,\cdot) = f$. Hence, we can integrate by parts in $t$ and use the definition of $L$ to give 
\begin{align*}
\iint_{\reu} A \nabla u \cdot \cl{\nabla \phi} \, \d t \d x
& = - \int_{\R^n} f \cdot\cl{\phi(0,\cdot)} \, \d x - \iint_{\reu} a \partial_t^2 u \cdot \cl{\phi} \, \d t \d x \\
&\quad + \iint_{\reu} d \nabla_x u \cdot \cl{\nabla_x \phi} \, \d t \d x \\
& = - \int_{\R^n} f \cdot\cl{\phi(0,\cdot)} \, \d x.
\end{align*}
The $\L^2$-pairing on the right-hand side can also be viewed as the $\Hdot^{-\nicefrac{1}{2},2} -\Hdot^{\nicefrac{1}{2},2}$-duality. Then the identity can be extended to all $\phi \in \Wdot^{1,2}(\reu)$ by density and we conclude $\dnuA u(0,\cdot) = f$.
\end{proof}

The semigroup construction provides solutions to the Neumann problem $(N)_{p}^\Le$ in an appropriate range of exponents.\index{energy solution!with Neumann datum}

\begin{prop}
\label{prop: NT for Neumann}
Let $q_{-}(L) < p < q_+(L)$. If $f \in \L^2 \cap \Wdot^{-1,2}$, then the energy solution $u$ with Neumann data $f$ satisfies
\begin{align*}
\|\NT(\nabla u)\|_p \simeq \|f\|_{\H^p}.
\end{align*}
\end{prop}

\begin{proof}
We have $q_-(L) = p_-(L)$ and $q_+(L) < p_+(L)$, see Theorem~\ref{thm: standard relation J(L) and N(L)}. Letting $g \coloneqq -(a L^{1/2})^{-1}f \in \W^{1,2}$ as before, we obtain from Proposition~\ref{prop: NT for regularity} and Theorem~\ref{thm: Riesz complete} that
\begin{align*}
\|\NT(\nabla u)\|_p 
\simeq \|\nabla_x g\|_{\H^p} 
\simeq \|a L^{1/2} g\|_{\H^p}
= \|f\|_{\H^p}. &\qedhere
\end{align*}
\end{proof}

\begin{proof}[Proof of Theorem~\ref{thm: blockneu}]	
Let $q_{-}(L) < p < q_+(L)$. According to Corollary~\ref{cor: AE interval and AM interval} this range is the same as what is called $I_L$ in \cite{AM}. We have seen in the introduction (Section~\ref{intro:Neumann}) that it suffices to prove the bound $\|\NT(\nabla u)\|_p \lesssim \|f\|_{\H^p}$, whenever $f \in \H^p \cap \Hdot^{-\nicefrac{1}{2},2}$ and $u$ is the energy solution with Neumann data $f$. 

By the universal approximation technique we an pick for any given $f$ a sequence $(f_k) \subseteq \H^p \cap \L^2 \cap \Wdot^{-1,2}$ with $f_k \to f$ as $k \to \infty$ in both $\H^p$ and $\Hdot^{-\nicefrac{1}{2},2}$. It follows from Proposition~\ref{prop: existence of energy solution Neumann} that the corresponding energy solutions $u_k$ tend to $u$ in $\Wdot^{1,2}(\reu)$, whereas Proposition~\ref{prop: NT for Neumann} implies that $(\nabla u_k)_k$ is a Cauchy sequence in $\T^{0,p}_\infty$. The limits for the gradients can be identified in $\Lloc^2(\reu)$ and the conclusion follows.
\end{proof}

Let us conclude with an additional uniqueness result for the Neumann problem. We remark that in our formulation of the Neumann problem the convergence of the conormal derivative to its trace is in the sense of distributions. By \cite[Cor.~1.2]{AM}, the Whitney averages convergence
\begin{align*}
\lim_{t \to 0} \bariint_{W(t,x)} a \partial_t u \, \d s \d y = g(x) \quad (\text{a.e. } x \in \R^n)
\end{align*}
of the conormal derivative of the unique solution to its trace comes as a bonus if $p \geq 1$ with $q_-(L) < p<q_+(L)$. In the case of block systems, one can reverse these interpretations of the boundary behavior and still obtain uniqueness, hence compatible well-posedness.  

\begin{thm} If $p \geq 1$ with $q_-(L) < p<q_+(L)$, then the following Neumann problem with non-tangential boundary trace is compatibly well-posed (modulo constants).\index{Neumann problem!with non-tangential trace} Given $g \in \H^p(\R^n; \IC^m)$, solve
\begin{equation*}
(\wt{N})_{p}^\Le  \quad\quad
\begin{cases}
\Le u=0   & (\text{in } \reu), \\
\NT(\nabla u)\in \L^p(\R^n),    \\
\lim_{t \to 0} \bariint_{W(t,x)} a \partial_tu \, \d s \d y = g(x) & (\text{a.e. } x\in \R^n).\end{cases}
\end{equation*}
\end{thm}

\begin{proof} 
In view of the preceding discussion we only need to establish uniqueness. 

According to \cite[Thm.~1.1]{AM} and our identification of $I_L$, the condition $\NT(\nabla u)\in \L^p(\R^n)$ implies (is equivalent to, in fact) the representation of the conormal gradient of $u$ via the $[DB]$-semigroup: 
\begin{align*}
\nabla_Au(t,\cdot)=\e^{-t[DB]}\nabla_Au\! \mid_{t=0} \quad (t>0),
\end{align*}
where $F_0 \coloneqq \nabla_Au\! \mid_{t=0}\in \H^p_D$ is characterized by $\ind_{\IC^+}(D B)F_0=F_0$ and the functional calculus is extended from $\IH_D^p$ to its completion $\H_D^p$ for the $\H^p$-norm as $\IH_{DB}^p=\IH_D^p$. It follows from \eqref{eq2: slo} that $F_0= [g, -\nabla_xL^{-1/2}(a^{-1}g)]^\top$ for some $g\in \H^p$.
Assume now that the Whitney averages of $a \partial_{t}u$ converge to $0$ almost everywhere at the boundary. By \cite[Cor.~1.2]{AM}, we know that this limit agrees with $g$ almost everywhere. Thus, $g=0$. We conclude that $F_0 $ vanishes identically and it follows that $u$ is constant in $\reu$. 
\end{proof}

\appendix
\section{Non-tangential maximal functions and traces}
\label{sec:technical}

\noindent In this appendix we collect some technical results involving non-tangential maximal functions with a focus on non-tangential trace theorems. 

Throughout, we consider the Whitney parameters $c_0 > 1$ and $c_1 > 0$ fixed, write $W(t,x) \coloneqq (c_0^{-1}t,c_0t) \times B(x,c_1 t)$ for $(t,x) \in \reu$ and for $q>0$ we use the $q$-adapted non-tangential maximal functions\index{non-tangential maximal function!$q$-adapted}
\begin{align*}
\NTq (F)(x) \coloneqq \sup_{t>0} \bigg(\bariint_{W(t,x)} |F(s,y)|^q \, \d s \d y \bigg)^{1/q} \quad (x \in \R^n)
\end{align*}
defined for measurable functions on $\reu$. In the case $q=2$ we simply write $\NT$ as before. Implicit constants always depend only on the Whitney parameters, dimensions and the exponents at stake. We shall not mention this at each occurrence.

It is common knowledge that different choices of Whitney parameters yield maximal functions with comparable $\L^p$-norms. For the reader's convenience we include a proof.\index{change of Whitney parameters!for $\NT$}

\begin{lem}[Change of Whitney parameters]
\label{lem: independence of Whitney parameters}
Let $0<p,q<\infty$. Let $c_0,c_1$ and $d_0,d_1$ be two pairs of Whitney parameters and let $\NTq^{(c)}$ and $\NTq^{(d)}$ be the corresponding maximal functions. Then,
\begin{align*}
\|\NTq^{(d)} (F)\|_p \simeq \|\NTq^{(c)}(F)\|_p
\end{align*}
for all measurable functions $F: \reu \to \R$.
\end{lem}

\begin{proof}
By symmetry it suffices to prove the estimate `$\lesssim$'. We write $W_{c_0,c_1}(t,x) \coloneqq (c_0^{-1}t, c_0t) \times B(x, c_1t)$. By compactness, we find points $(t_i, x_i) \in W_{d_0,d_1}(1,0)$, $i=1,\ldots,N$, such that the sets $W_{c_0, c_1/2}(t_i,x_i)$ cover $W_{d_0,d_1}(1,0)$. Using the affine transformation $(s,y) \mapsto (ts, x+ty)$, we obtain
\begin{align}
\label{eq: Whitney region covering different parameters}
W_{d_0,d_1}(t,x) \subseteq \bigcup_{i=1}^N W_{c_0, c_1/2}(t_i t,x+tx_i)
\end{align}
for any $(t,x) \in \reu$. Since $|x-(x+tx_i)| \leq t d_1 \leq d_0 d_1 t_i t$, we get
\begin{align*}
\bigg(\bariint_{W_{d_0,d_1}(t,x)} |F|^q \bigg)^{1/q} \leq C \sum_{i=1}^N \sup_{|x-y| < d_0 d_1 t_i t} \bigg(\bariint_{W_{c_0, c_1/2}(t_i t,y)} |F|^q \bigg)^{1/q}
\end{align*}
for an admissible constant. For measurable $H: \reu \to \R$ let $H_{*,\eta}(x) \coloneqq \sup_{|x-y|<\eta t} |H(t,y)|$ be the pointwise non-tangential maximal function with aperture $\eta$. With $H(t,y) \coloneqq (\bariint_{W_{c_0, c_1/2}(t,y)} |F|^q )^{1/q}$ the previous bound yields
\begin{align*}
\NTq^{(d)}(F)(x) \leq C H_{*,d_0 d_1}(x) \quad (x \in \R^n).
\end{align*}
On the other hand, $|y-x| < \nicefrac{tc_1}{2}$ implies $B(y,\nicefrac{tc_1}{2}) \subseteq B(x, c_1 t)$, so that
\begin{align*}
H_{*,c_1/2}(x) \leq C \NTq^{(c)}(F)(x)  \quad (x \in \R^n).
\end{align*}
For the classical pointwise non-tangential maximal functions we can change the aperture~ \cite[Sec.~II.2.5.1]{Stein93}: There is $C=C(n,c_0,c_1, d_0,d_1)$ such that
\begin{align*}
|\{x: \R^n : H_{*,d_0 d_1}(x) > \alpha \}| \leq C |\{x: \R^n : H_{*,c_1/2}(x) > \alpha \}| \quad (\alpha >0).
\end{align*}
The claim follows from the previous three bounds and the layer cake formula.
\end{proof}

\begin{rem}
\label{rem: independence of Whitney parameters}
The covering argument in \eqref{eq: Whitney region covering different parameters} implies directly that different choices of Whitney parameters for the Whitney average functionals yield equivalent $\Z$-space norms.\index{change of Whitney parameters!for $\Z$-spaces} 
\end{rem}

We continue with a useful non-tangential embedding.

\begin{lem}[{\cite[Lem.~2.2]{AM} \& \cite[Lem.~A.2]{HMiMo}}]
\label{lem:AM} 
If $0<p<r\le 2$, then there is a constant $C$ such that for all measurable functions $F: \reu \to \R$,
\begin{align*}
\bigg(\iint_{\reu} |F(t,x)|^r\,  t^{n(\frac{r}{p}-1)}\, \frac{\d t \d x}{t}\bigg)^{1/r} \leq C\|\NT (F)\|_{p}.
\end{align*}
\end{lem}

We turn our attention to non-tangential trace theorems.

\begin{defn}\label{def: NT trace}
A locally integrable function $u$ on $\reu$ is said to have a \emph{non-tangential trace} (in the sense of Whitney averages)\index{non-tangential trace} if there exists a function $u_0$ on $\R^n$ such that for almost every $x\in \R^n$, 
\begin{align*}
\lim_{t\to 0} \bariint_{W(t,x)} u(s,y)  \, \d s \d y =  u_0(x).
\end{align*}
\end{defn}
As a pointwise limit of measurable functions, such a trace is necessarily measurable. The following is a variation of Kenig--Pipher's trace theorem~\cite[Thm.~3.2]{KP}\index{Theorem!Kenig--Pipher trace} that covers exponents $p \leq 1$ and applies to averaged non-tangential maximal functions. This has appeared (without proof) in many earlier works and we take the opportunity to close the gap.

\begin{prop}
\label{prop: KP}
Let $p \in (\nicefrac{n}{(n+1)},\infty)$ and $q \in [1,\infty)$. Let $u \in \Wloc^{1,q}(\reu)$ satisfy $\|\NTq(\nabla u)\|_p < \infty$. Then there exists a non-tangential trace $u_0$  with the following properties.
\index{Theorem!Kenig--Pipher trace}
\begin{enumerate}
	\item Let $r \in (0, \infty)$ and assume $r \leq \frac{(n+1)q}{n+1-q}$ if $q < n+1$. For almost every  $x \in \R^n$ and all $t>0$,
	\begin{align*}
	\qquad \bigg(\bariint_{W(t,x)} |u(s,y) - u_0(x)|^r \, \d s \d y \bigg)^{\frac{1}{r}} \leq C t \NTq(\nabla u)(x).
	\end{align*}
	In particular, the left-hand side tends to $0$ as $t \to 0$ and $u_0$ does not depend (in the almost everywhere sense) on the choice  of the Whitney parameters.
	
	\item $u_0$ is of class $\Hdot^{1,p}(\R^n)$ with $\|\nabla_x u_0\|_{\H^p} \leq C \|\NTq(\nabla u)\|_p$.
	\item Let $r$ be as in (i) and suppose in addition that $r< \frac{np}{n-p}$ if $p<n$. Then,
	\begin{align*}
	\bigg\| \NTr \bigg( \frac{u-u_0}{t} \bigg) \bigg\|_p \leq C \|\NTq(\nabla u)\|_p.
	\end{align*}
	\item Suppose that either $p \geq 1$ or that $p<1$ and that there exists $\eps > 0$ such that $\sup_{0<t<\eps} \|u(t,\cdot)\|_{\frac{np}{n-p}} < \infty$. Then,
	\begin{align*}
	\lim_{t\to 0} \barint_{(c_0)^{-1}t}^{c_0 t} u(s,\cdot) \, \d s = u_0 \quad (\text{in } \cD'(\R^n)).
	\end{align*}
\end{enumerate}
\end{prop}

\begin{rem}
\label{rem: KP}
In applications we usually have $q=2$ and $r \in (0,2]$, which is admissible in (i). Also $r \in (0,1]$ is always admissible in (iii). {Identification of the non-tangential trace with a distributional limit seems to be far from obvious in the case $p<1$. We got the idea to impose the additional condition on $u$ in (iv) from \cite[Lem.~5.2]{HMiMo}. In our applications to the regularity problem $(R)_p^\Le$ it follows from Sobolev embeddings and strong continuity of the Poisson semigroup.}
\end{rem}

For the proof we need a simple lemma on real functions.

\begin{lem}
\label{lem: Whitney trace preparation}
Let $h: (0,\infty) \to \R$ be a function for which there are constants $\theta > 1$, $\alpha > 0$ and $C \geq 0$ such that $|h(t) - h(\tau)| \leq C t^\alpha$, whenever $\tau \in [\theta^{-1}t, t]$. Then $h(0) \coloneqq \lim_{s \to 0} h(s)$ exists and satisfies
\begin{align*}
	|h(t) - h(0)| \leq \frac{Ct^\alpha}{1-\theta^{-\alpha}} \quad (t>0).
\end{align*}
\end{lem}

\begin{proof}
Given $0 < \tau \leq t$, let $k$ be the largest integer with $\tau \leq \theta^{-k} t$. By a telescopic sum we find
\begin{align*}
|h(t)- h(\tau)|
&\leq |h(\theta^{-k}t)- h(\tau)| + \sum_{j=1}^k |h(\theta^{-j+1}t) - h(\theta^{-j}t)| \\
&\leq \sum_{j=1}^{k+1} C \theta^{\alpha(-j+1)}t^\alpha \leq \frac{C t^\alpha}{1-\theta^{-\alpha}}.
\end{align*}
This proves the Cauchy property for $h$ at $0$. Hence, $h(0)$ is defined and the estimate follows by sending $\tau \to 0$.
\end{proof}

\begin{proof}[Proof of Proposition~\ref{prop: KP}]
Throughout the proof we write $\NTq^{\mathrm{large}}$ for a non-tangential maximal function with Whitney parameters $c_0^{\mathrm{large}} > c_0$ and $c_1^{\mathrm{large}} \geq c_1$  that will be further specified if needed. We denote the associated Whitney regions by $W^{\mathrm{large}}(t,x)$.

\medskip

\noindent \emph{Proof of} (i). Let $\theta > 1$ be such that $c_0^{\mathrm{large}} = \theta c_0$. If $\tau \in [\theta^{-1}t, t]$, then both $W(\tau,x)$ and $W(t,x)$ are contained in $W^{\mathrm{large}}(t,x)$ and we can estimate
\begin{align}
\label{eq0: KP}
\begin{split}
|(u)_{W(\tau,x)}& - (u)_{W(t,x)}|
\\
&\leq  \bariint_{W(\tau, x)} |u - (u)_{W(t,x)}| \, \d s \d y \\
&\lesssim \bigg(\bariint_{W^{\mathrm{large}}(t,x)} |u - (u)_{W(t,x)}|^q \, \d s \d y \bigg)^{\frac{1}{q}}\\
&\lesssim  t \bigg( \bariint_{W^{\mathrm{large}}(t,x)} |\nabla u|^q \, \d s \d y \bigg)^{\frac{1}{q}}\\
&\leq t \NTq^{\mathrm{large}}(\nabla u)(x),
\end{split}
\end{align}
where the third step is due to the Poincaré inequality on cylinders. From the assumption on $u$ and Lemma~\ref{lem: independence of Whitney parameters} we obtain that $\NTq^{\mathrm{large}}(\nabla u)(x)$ is finite for a.e.\ $x \in \R^n$. In this case Lemma~\ref{lem: Whitney trace preparation} yields the existence of a non-tangential trace $u_0(x)$ with control
\begin{align}
\label{eq1: prop KP}
|(u)_{W(t,x)} - u_0(x)| \leq C t \NTq^{\mathrm{large}}(\nabla u)(x).
\end{align} 
This argument works for any choice of Whitney parameters. In order to see that $u_0$ is always the same, it suffices (by transitivity) to verify that the trace $u_0^{\mathrm{large}}$ corresponding to the regions $W^{\mathrm{large}}(t,x)$ agrees with $u_0$. By the argument in \eqref{eq0: KP} we also have
\begin{align*}
	|(u)_{W(t,x)} - (u)_{W^{\mathrm{large}}(t,x)}|  \lesssim t \NTq^{\mathrm{large}}(\nabla u)(x)
\end{align*}
and hence the limits as $t \to 0$ are the same almost everywhere.

As for the estimate in (i) we pick some smaller Whitney parameters with associated regions $w(t,x)$ such that $W(t,x) = w^{\mathrm{large}}(t,x)$. In this scenario \eqref{eq1: prop KP} becomes
\begin{align*}
	|(u)_{w(t,x)} - u_0(x)| \leq C t \NTq(\nabla u)(x)
\end{align*}
and the restriction on $r$ allows us to use the Sobolev--Poincaré inequality in order to give
\begin{align*}
\bigg(\bariint_{W(t,x)} &|u - u_0(x)|^r\, \d s \d y \bigg)^{\frac{1}{r}} \\
&\lesssim \bigg(\bariint_{W(t,x)} |u - (u)_{w(t,x)}|^r \, \d s \d y \bigg)^{\frac{1}{r}} + |(u)_{w(t,x)} - u_0(x)|\\
&\lesssim t \NTq(\nabla u)(x).
\end{align*}

\noindent \emph{Proof of} (ii). We use the following result: If there is $g \in \L^p(\R^n)$ such that for almost every $x,y \in \R^n$,
\begin{align}
\label{eq2: prop KP}
|u_0(x) - u_0(y)| \leq |x-y| (g(x) + g(y)),
\end{align}
then $u_0 \in \Hdot^{1,p}$ with $\|\nabla_x u_0\|_{\H^p} \lesssim \|g\|_p$. For $p>1$ this is Hajlasz's Sobolev space characterization~\cite[Thm.~1]{Hajlasz} and the result for exponents $\nicefrac{n}{(n+1)} < p \leq 1$ has been obtained in \cite[Thm.~1 \& Prop.~5]{Koskela-Saksman}. 

Now, let $x, y \in \R^n$ and set $t \coloneqq |x-y|$. We take $c_1^{\mathrm{large}} \geq 1+c_1$. Since $B(y, c_1t) \subseteq B(x,(1+c_1)t)$, we have $W(t,y) \subseteq W^{\mathrm{large}}(t,x)$ and Poincaré's inequality yields again
\begin{align*}
|(u)_{W(t,y)} - (u)_{W(t,x)}| \leq C t \NTq^{\mathrm{large}} (\nabla u)(x).
\end{align*}
Together with \eqref{eq1: prop KP}, we see that we can take $g \coloneqq 2 C \NTq^{\mathrm{large}} (\nabla u)$. Note that $\|g\|_p \simeq \|\NTq(\nabla u)\|_p$ by Lemma~\ref{lem: independence of Whitney parameters}.

\medskip

\noindent \emph{Proof of} (iii). It suffices to find a function $h$ with $\|h\|_p \leq C \|\NTq(\nabla u)\|_p$ such that for a.e.\ $x \in \R^n$ and all $t>0$,
\begin{align}
\label{eq3: prop KP}
\bigg(\bariint_{W(t,x)} |u-u_0|^r \, \d y \d s \bigg)^{\frac{1}{r}} \leq h(x).
\end{align}
Indeed, since we are integrating $s$ on $(c_0^{-1}t,c_0 t)$ on the left-hand side, the bound required in (iii) follows immediately. The argument slightly differs depending on whether or not we have $p > 1$. Let us first assume that this is the case. 

The additional restriction on $r$ makes sure that we can find some $\varrho \in (1, p \wedge n)$ such that $\nicefrac{n \varrho}{(n-\varrho)} \geq r$. Hence, by H\"older's inequality followed by the Sobolev--Poincar\'e inequality, we have
\begin{align}
\label{eq4: prop KP}
\begin{split}
\bigg(\barint_{B(x,c_1t)} |u_0-&(u_0)_{B(x,c_1 t)}|^r \, \d y \bigg)^{\frac{1}{r}} \\
&\leq \bigg(\barint_{B(x,c_1t)} |u_0-(u_0)_{B(x,c_1 t)}|^{\frac{n \varrho}{n-\varrho}} \, \d y \bigg)^{\frac{1}{\varrho} - \frac{1}{n}} \\
&\leq C t \bigg(\barint_{B(x,c_1t)} |\nabla_x u_0|^\varrho \, \d y \bigg)^{\frac{1}{\varrho}}.
\end{split}
\end{align}
Since $\nicefrac{n \varrho}{(n-\varrho)} \geq 1$, we can also argue as in \eqref{eq0: KP}, using the Sobolev--Poincar\'e inequality in the second step, to get whenever $\tau \in [\nicefrac{t}{2},t]$,
\begin{align*}
|(u_0)_{B(x,c_1t)} - (u_0)_{B(x,c_1\tau)}| 
\lesssim t \Max(|\nabla_x u_0|^\varrho)(x)^{\frac{1}{\varrho}}.
\end{align*}
By Lebesgue differentiation and Lemma~\ref{lem: Whitney trace preparation} we get 
\begin{align}
\label{eq5: prop KP}
|(u_0)_{B(x,c_1t)} - u_0(x)| 
\lesssim t \Max(|\nabla_x u_0|^\varrho)(x)^{\frac{1}{\varrho}}
\end{align}
for a.e.\ $x \in \R^n$. 
Using the decomposition
\begin{equation*}
u(s,y) - u_0(y)= u(s,y) - u_0(x)+ u_0(x)-(u_0)_{B(x,c_1t)}+ (u_0)_{B(x,c_1t)}- u_0(y)
\end{equation*}
and combining (i), \eqref{eq4: prop KP} and \eqref{eq5: prop KP}, we arrive at
\begin{align*}
\bigg(\bariint_{W(t,x)} |u(s,y) - &u_0(y) |^r \, \d s \d y\bigg)^{\frac{1}{r}} \\
&\lesssim t\Big(\NTq(\nabla u)(x)  + \Max(|\nabla_x u_0|^\varrho)(x)^{\frac{1}{\varrho}}\Big)
\end{align*}
for a.e.\ $x \in \R^n$ and all $t>0$. The right-hand side is admissible for \eqref{eq3: prop KP} by assumption on $u$, the $\L^{p/\varrho}$-boundedness of the maximal function and the result of (ii).

We turn to the case $p \leq 1$. Since $p > \nicefrac{n}{(n+1)}$, we can pick $\varrho \in (\nicefrac{n}{(n+1)}, p)$ with $\nicefrac{n \varrho}{(n-\varrho)} \geq (r \vee 1)$. Since the function $g$ in \eqref{eq2: prop KP} is locally $\varrho$-integrable, we have Hajlasz's Sobolev--Poincaré inequality\index{inequality!Hajlasz's Sobolev--Poincaré}
\begin{align*}
\bigg( \barint_{B(x,c_1t)} |u_0 - (u_0)_{B(x,c_1t)} |^{\frac{n \varrho}{n-\varrho}} \, \d y \bigg)^{\frac{1}{\varrho} - \frac{1}{n}} \lesssim t \bigg( \barint_{B(x,2c_1t)} g^\varrho \, \d y\bigg)^{\frac{1}{\varrho}},
\end{align*}
see~\cite[Thm.~8.7]{Hajlasz-Survey}. Hence, except for replacing $\nabla_x u_0$ by $g$, the argument stays the same. 

\medskip

\noindent \emph{Proof of} (iv). Let $B \subseteq \R^n$ be a ball and let $\phi \in \C_0^\infty(B)$. We use the averaging trick to write
\begin{align}
\label{eq6: KP}
\begin{split}
\int_{\R^n} \bigg(\barint_{(c_0)^{-1}t}^{c_0 t} &u(s, y) \, \d s - u_0(y) \bigg) \phi(y) \, \d y \\
&=\int_{\R^n} \bigg(\bariint_{W(t,x)} (u -u_0) \phi \, \d s \d y\bigg)  \d x \\
&\eqqcolon \int_{\R^n} F_{t}^\phi(x) \, \d x.
\end{split}
\end{align}
We have to show that the right-hand side tends to $0$ as $t \to 0$. From now on, we require $t< \nicefrac{r(B)}{c_1}$, so that all functions $F_t^\phi$ have support in $2 B$. 

If $p \geq 1$, then  \eqref{eq3: prop KP} for the admissible choice $r=1$ gives us $|F_{t}^\phi(x)| \leq \|\phi\|_\infty t h(x)$ and $h$ is locally integrable, so we are done.

In the case $p<1$ we need a different argument and this is where the additional assumption $C\coloneqq \sup_{0<t<\eps} \|u(t,\cdot)\|_{np/(n-p)} < \infty$ comes into play. We abbreviate $p^* \coloneqq \nicefrac{np}{(n-p)} > 1$.  We can restrict ourselves to $t < (\nicefrac{r(B)}{c_1} \wedge \nicefrac{\eps}{c_0})$ and $x \in 2 B$. In this case, $B(x,c_1t) \subseteq 3 B$ and by H\"older's inequality we can crudely bound
\begin{align*}
\bariint_{W(t,x)} |u-u_0| \, \d s \d y
&\lesssim t^{1-\frac{n}{p}} \barint_{c_0^{-1}t}^{c_0 t} \|u(s,\cdot) - u_0\|_{\L^{p_*}(3 B)} \, \d s \\
&\leq (C+\|u_0\|_{\L^{p^*}(3B)}) t^{1-\frac{n}{p}}.
\end{align*}
We have  $u_0 \in \Lloc^{p^*} $ from  the Hardy--Sobolev embedding or by the following direct argument. We have, for $t$ small enough, using H\"older's inequality and averaging, 
\begin{equation*}
\int_{3B} | (u)_{W(t,x)}|^{p^*}\, \d x \leq \barint_{c_0^{-1}t}^{c_0 t} \int_{4B} |u(s, y)|^{p^*}\, \d y\d s \le C^{p^*}.
\end{equation*}
Thus, by Fatou's lemma, we obtain $\int_{3B} | u_0|^{p^*} \, \d x \le C^{p^*}$.
Now, we use the $p$-th power of \eqref{eq3: prop KP} (with $r=1$) and the $(1-p)$-th power of the crude bound in order to get for a.e.\ $x \in 2B$ that
\begin{align*}
 |F_{t}^\phi(x)|\leq \|\phi\|_\infty \bariint_{W(t,x)} |u-u_0| \, \d s \d y
\leq \|\phi\|_\infty C' t^{1+ n - \frac{n}{p}} h(x)^p.
\end{align*}
On the right the power of $t$ is positive since $p>\nicefrac{n}{(n+1)}$ and we have $h^p \in \L^1$.  Thus, we get the desired convergence in \eqref{eq6: KP} when passing to the limit as $t \to 0$.
\end{proof}

Next, we present variants of the non-tangential trace theorem for tent and $\Z$-spaces. In our applications we shall only encounter functionals based on $\L^2$-averages such as $S$ and $W$ used to define tent and $\Z$-spaces, respectively. For simplicity we stick to that case. The following results have appeared in~\cite[Thm.~6.3]{BM} ($p=\infty$) and \cite[Sec.~6.6]{AA} ($p< \infty$). For the sake of self-containedness we include a proof that follows the same pattern as before.\index{Theorem! non-tangential trace ($\Y^{\alpha-1,p}$)} The lower bound on $p$, notably to identify the non-tangential trace with a distributional limit, is now related to fractional Sobolev embeddings and the argument turns out to be conceptually simpler than in Proposition~\ref{prop: KP}. 

As usual, we treat both scales of spaces simultaneously and let $\Y$ denote one of $\T$ or $\Z$.

\begin{prop}
\label{prop: NT trace Y}
Let $\alpha \in (0,1)$ and $\nicefrac{n}{(n+\alpha)}<p<\infty$. Let $u \in \Wloc^{1,2}(\reu)$ satisfy $\|\nabla u\|_{\Y^{\alpha-1,p}} < \infty$. Then there exists a non-tangential trace $u_0$  with  the following properties.
\begin{enumerate}
	\item Let $r \in (0, \infty)$ and assume $r \leq \frac{2(n+1)}{n-1}$ if $n>1$. 
	For all $x \in \R^n$ and all $t>0$,
	\begin{align*}
		\bigg(\bariint_{W(t,x)} |u(s,y) - u_0(x)|^r \, \d s \d y\bigg)^{\frac{1}{r}} \leq C t^\alpha \Theta(x)
	\end{align*}
	with $\|\Theta\|_p \leq C \|\nabla u\|_{\Y^{\alpha-1,p}}$.  In particular, the left-hand side tends to $0$ almost everywhere as $t \to 0$ and $u_0$ does not depend on the choice  of the Whitney parameters.
	\item There is convergence
	\begin{align*}
		\lim_{t\to 0} \barint_{(c_0)^{-1}t}^{c_0 t} u(s,\cdot) \, \d s = u_0 \quad (\text{in } \cD'(\R^n)).
	\end{align*}
	\item The results above continue to hold for $p=\infty$ and $\nabla u \in \Z^{\alpha-1,\infty}$. In that case $\Theta(x) = \| \nabla u\|_{\Z^{\alpha-1,\infty}}$ and $u_0$ is of class $\Lamdot^\alpha$ with $\|u_0\|_{\Lamdot^\alpha} \leq C \|\nabla u\|_{\Z^{\alpha-1,\infty}}$.
\end{enumerate}
\end{prop}

The following lemma contains the construction of the function $\Theta$ in part (i) for finite $p$.

\begin{lem}
\label{lem: NT trace Theta function}
Let $\alpha \in \R$, $p \in (0,\infty)$ and $F \in \Y^{\alpha-1,p}$. There exists a measurable function $\Theta: \R^n \to [0,\infty)$ with $\|\Theta\|_p \leq C \|F\|_{\Y^{\alpha-1,p}}$ such that
\begin{align*}
	\bigg( \bariint_{W(t,x)} |s^{1-\alpha}F|^2 \, \d s \d y \bigg)^{\frac{1}{2}} \leq C \Theta(x) \quad ((t,x) \in \reu).
\end{align*}
\end{lem}

\begin{proof}
We begin with the case $\Y = \Z$ and set
\begin{equation*}
	\Theta(x) \coloneqq  \bigg(\int_0^\infty \bigg( \bariint_{W^{\mathrm{large}}(t,x)} |s^{1-\alpha} F|^2 \, \d s \d y \bigg)^{\frac{p}{2}}  \frac{\d t}{t} \bigg)^{\frac 1 p},
\end{equation*}
where $W^{\mathrm{large}}(t,x)$ are Whitney regions with Whitney parameter $c_0^{\mathrm{large}} \coloneqq 2 c_0$. Since
\begin{equation*}
	W(t,x) \subseteq W^{\mathrm{large}}(\tau,x)  \quad ( \tau \in [\nicefrac{t}{2}, t]),
\end{equation*}
we can infer that
\begin{align*}
	\bigg( \bariint_{W(t,x)} &|s^{1-\alpha}F|^2 \, \d s \d y \bigg)^{\frac{p}{2}} \\
	&\lesssim \int_{t/2}^t \bigg( \bariint_{W^{\mathrm{large}}(\tau,x)} |s^{1-\alpha}F|^2 \, \d s \d y \bigg)^{\frac{p}{2}}
	\frac{\d \tau}{\tau}
\end{align*}
and the right-hand side is bounded by $\Theta(x)^p$. Moreover, a change of Whitney parameters for $\Z$-space quasinorms yields $\|\Theta\|_p \simeq \|F\|_{\Z^{\alpha-1,p}}$.

In the case $\Y = \T$ we can simply set
\begin{equation*}
	\Theta(x) \coloneqq \bigg(\iint_{|x-y|<2c _1 s}  |s^{1-\alpha}F|^2 \, \frac{\d s \d y}{s^{n+1}} \bigg)^{\frac{1}{2}}  
\end{equation*}
since $W(t,x)$ is contained in the cone appearing in the integral. By a change of aperture in tent space quasinorms we  conclude that $\|\Theta\|_p \simeq  \|F\|_{\T^{\alpha-1,p}}$.
\end{proof}

\begin{proof}[Proof of Proposition~\ref{prop: NT trace Y}]

We use the same notation as in the proof of Proposition~\ref{prop: KP} and follow the same line of thoughts.

\medskip

\noindent \emph{Proof of} (i). Let $c_0^{\mathrm{large}} \coloneqq 2 c_0$. If $\tau \in [\nicefrac{t}{2}, t]$, then both $W(\tau,x)$ and $W(t,x)$ are contained in $W^{\mathrm{large}}(t,x)$ and using the Poincaré inequality with $q= 2$ as in \eqref{eq0: KP}, we obtain 
\begin{align*}
	|(u)_{W(\tau,x)} - (u)_{W(t,x)}| 
	&\lesssim t \bigg( \bariint_{W^{\mathrm{large}}(t,x)} |\nabla u|^2 \, \d s \d y \bigg)^{\frac{1}{2}}
	\\
	& \lesssim t^\alpha \bigg( \bariint_{W^{\mathrm{large}}(t,x)} |s^{1-\alpha}\nabla u|^2 \, \d s \d y 
	\bigg)^{\frac{1}{2}}.
\end{align*}
Lemma~\ref{lem: NT trace Theta function} applied to the `large' Whitney regions yields a function $\Theta$ with $\|\Theta\|_p \lesssim \|\nabla u\|_{\Y^{\alpha-1,p}}$ such that
\begin{align*}
	|(u)_{W(\tau,x)} - (u)_{W(t,x)}| \lesssim t^\alpha \Theta(x).
\end{align*}
Now, we can apply Lemma~\ref{lem: Whitney trace preparation} to obtain a non-tangential trace $u_0(x)$ with control
\begin{align}
	\label{eq2: NT trace Z finite exposants}
	|(u)_{W(t,x)} - u_0(x)| \lesssim t^\alpha \Theta(x),
\end{align}
whenever $\Theta(x)<\infty$, that is, almost everywhere. That $u_0$ is independent of the choice of Whitney parameters follows as in the proof of Proposition~\ref{prop: KP} and the restriction on $r$ allows us to use the Sobolev--Poincaré inequality again in order to conclude 
\begin{align*}
	\bigg(\bariint_{W(t,x)}& |u - u_0(x)|^r\, \d s \d y \bigg)^{\frac{1}{r}} \\
	&\lesssim \bigg(\bariint_{W(t,x)} |u - (u)_{W(t,x)}|^r \, \d s \d y \bigg)^{\frac{1}{r}} + |(u)_{W(t,x)} - u_0(x)|\\
	&\lesssim t^\alpha \bigg(\bariint_{W(t,x)} |s^{1-\alpha} \nabla u|^2 \, \d s \d y \bigg)^{\frac{1}{2}} + |(u)_{W(t,x)} - u_0(x)|\\
	&\lesssim t^\alpha \Theta(x).
\end{align*}

\medskip

\noindent \emph{Proof of (ii)}. 
We begin with the case $p>1$. With the notation of the proof of Proposition~\ref{prop: KP}.(iv) we have to show that
\begin{align}
\label{eq3: NT trace Z finite exposants}
\lim_{t \to 0} \int_{\R^n} |F_{t}^\phi(x)|\, \d x,
\end{align}
where $F_t^\phi(x) = \bariint_{W(t,x)} (u - u_0) \phi \, \d s \d y$ is supported in $2B$ if the support of $\phi$ is contained in $B$ and $t\le \nicefrac{r(B)}{c_1}$. We record two elementary observations.
\begin{itemize}
\item If $y$ belongs to a ball $B(x,c_1t)$, then
\begin{align*}
	|(u)_{W(t,y)} - (u)_{W(t,x)}| \lesssim  t^\alpha \Theta(x).
\end{align*}
Indeed, we take $c_1^{\mathrm{large}} \geq 2c_1$. Since $B(y, c_1t) \subseteq B(x,2 c_1 t)$, we have $W(t,y) \subseteq W^{\mathrm{large}}(t,x)$ and   Poincaré's inequality yields this inequality as before. 
\item  For almost every $y\in B(x,c_1t)$ the first observation together with \eqref{eq2: NT trace Z finite exposants} yields
\begin{align*}
	|u_0(y) - &u_0(x)| \\
	&\leq |u_0(y) - (u)_{W(t,y)}| + |(u)_{W(t,y)}-(u)_{W(t,x)}| \\
	&\quad + |(u)_{W(t,x)} - u_0(x)|\\
	& \lesssim  t^\alpha (\Theta(x) + \Theta(y)).
\end{align*}
\end{itemize}
The second observation implies
\begin{align*}
|u(s,y) -u_0(y)| \lesssim |u(s,y) -u_0(x)|  + t^\alpha (\Theta(x) + \Theta(y))
\end{align*}
and taking into account (i) with $r=1$, we are left with
\begin{align*}
|F_t^\phi(x)| \lesssim \|\phi\|_\infty  t^\alpha  \Max(\Theta)(x).
\end{align*}
The maximal theorem ensures that $\Max(\Theta) \in \L^p$ and since $F_t^\phi$ is supported in $2B$ we conclude \eqref{eq3: NT trace Z finite exposants}.

In the case $p \leq 1$ we use the embedding $\Y^{\alpha-1,p} \subseteq \Y^{\beta-1,q}$ for $0<p<q<\infty$ and $\alpha- \beta = n(\nicefrac{1}{p}-\nicefrac{1}{q})$, see \cite[Thm.~2.34]{AA}. We have $\alpha > n(\nicefrac{1}{p} - 1)$ by assumption, which allows us to pick $q>1$ and $0<\beta<\alpha$. Hence, we are back in the case of integrability above $1$. 

\medskip

\noindent \emph{Proof of (iii)}. If $p=\infty$ and $\nabla u \in \Z^{\alpha - 1,\infty}$, then the constant function $\Theta(x) \coloneqq C\|\nabla u\|_{\Z^{\alpha-1,\infty}}$ has the properties stated in Lemma~\ref{lem: NT trace Theta function} by definition of the $\Z^{\alpha-1,\infty}$-norm. Hence, we can repeat the first two steps and the second observation in the proof of (ii) yields $\|u_0\|_{\Lamdot^\alpha} \leq C \|\nabla u\|_{\Z^{\alpha-1,\infty}}$.
\end{proof}
\section{The \texorpdfstring{$\L^p$}{Lp}-realization of a sectorial operator in \texorpdfstring{$\L^2$}{L2}}
\label{sec: The Lp realization}

\noindent The following result is folklore but we could not find a precise statement in the literature.\index{Lp-realization@$\L^p$-realization (of an operator)}

\begin{prop}
\label{prop: consistent operator}
Let $T$ be a sectorial operator in $\L^2$ and let $p \in (1,\infty)$. Suppose that there exists $\mu \in (\omega_T, \pi)$ such that
\begin{align}
\label{eq0: consistent operator}
\|z (z - T)^{-1}f\|_p \lesssim \|f\|_p \quad (f \in \L^p \cap \L^2,\, z \in \IC \setminus \cl{\S_{\mu}^+}).
\end{align}
The case $\mu=\pi$ with the convention that $\IC \setminus \cl{\S_\pi^+} \coloneqq (-\infty,0)$ is also permitted.
Then there is a (unique) sectorial operator $T_p$ in $\L^p$ of angle smaller than $\mu$ that satisfies
\begin{align}
\label{eq1: consistent operator}
(z - T_p)^{-1}f = (z - T)^{-1}f \quad (f \in \L^p \cap \L^2, \, z \in \IC \setminus \cl{\S_{\mu}^+}).
\end{align}
Moreover, $T_p f = Tf$ for $f \in \dom(T_p) \cap \dom(T)$ and if $T$ is injective, then so is $T_p$. The corresponding statement for bisectorial operators also holds.
\end{prop}

\begin{rem}
\label{rem: consistent operator}
The assumption with $\mu = \pi$ simply means that $T$ satisfies $\|(1+t^2 T)^{-1}f\|_p \lesssim \|f\|_p$ for all $f \in \L^p \cap \L^2$ and all $t>0$.
\end{rem}

The operator $T_p$ is usually called \emph{$\L^p$-realization} of $T$. We have tried to avoid passing to an $\L^p$-realization whenever possible, but knowing that we always can, turns out helpful when dealing with abstract results that do not need a distinguished space such as $\L^2$ to start with. One such example is Theorem~\ref{thm: Le Merdy}.

Condition~\eqref{eq0: consistent operator} is obviously necessary for the existence of a $\L^p$-realization with consistent resolvents as in \eqref{eq1: consistent operator} and the latter uniquely determines $T$. We also obtain consistency of $T_p$ and $T$, whereas consistency of general invertible operators does not imply consistency of their inverses, compare with the discussion below Definition~\ref{def: Hodge interval}.   

\begin{proof}
By \eqref{eq0: consistent operator} we can define $R(z)$ as the extension by density of $(z-T)^{-1}$ to $\L^p$. Then $(z R(z))_{z \in \IC \setminus \cl{\S_{\mu}^+}}$ is a uniformly bounded family in $\L^p$ with the property
\begin{align}
\label{eq2: consistent operator}
R(z) - R(z') = (z' - z)R(z)R(z') \quad (z,z' \in \IC \setminus \cl{\S_\mu^+}).
\end{align}
We claim that for $f \in \L^p$ we have
\begin{align}
\label{eq3: consistent operator}
\lim_{ z\in (-\infty,0), z \to -\infty} z R(z)f &= f \quad (\text{weakly in } \L^p)
\intertext{and if in addition $T$ is injective, also that}
\label{eq4: consistent operator}
\lim_{z \to 0} z R(z)f &= 0 \quad (\text{weakly in } \L^p).
\end{align}
Indeed, since $T$ is sectorial in $\L^2$, the limits exist strongly in $\L^2$ if $f \in \L^p \cap \L^2$, see \cite[Prop.~2.1.1(a)]{Haase}. The extension then follows by uniform boundedness and density. 

By \eqref{eq2: consistent operator}, $R(-1)f=0$ implies $R(z)f=0$ for all $z$. Then $f=0$ follows from \eqref{eq3: consistent operator}, so $R(-1)$ is injective. We show that $T_p \coloneqq - R(-1)^{-1} -1 $ has the required properties. 

For $f \in \dom(T_p)$ we have
\begin{align*}
R(z)(z - T_p) f = R(z)((z+1) R(-1) + 1) R(-1)^{-1}f = f,
\end{align*}
where the final step uses \eqref{eq2: consistent operator}. Likewise, for $g \in \L^p$ we have
\begin{align*}
R(z)g = R(-1)(g-(z+1)R(z)g) \in \dom(T_p)
\end{align*}
and
\begin{align*}
(z - T_p)R(z)g &= (z+1+R(-1)^{-1})R(z)g = g.
\end{align*}
This proves $(z - T_p)^{-1} = R(z)$, so  \eqref{eq1: consistent operator} holds. By a Neumann series, the uniform boundedness of the family $(z R(z))$ implies that $T_p$ is a sectorial operator of angle smaller than $\mu$. Now, suppose that $f \in \dom(T_p) \cap \dom(T)$. Then
\begin{align}
\label{eq5: consistent operator}
z (z - T_p)^{-1} T_p f 
= z (z - T)^{-1} T f
\end{align}
since both terms can be expanded in terms of $R(z)$. When $z \in (-\infty,0)$ tends to $-\infty$, the left-hand side tends to $T_pf$ weakly in $\L^p$ and the right-hand side tends to $Tf$ strongly in $\L^2$, see \eqref{eq3: consistent operator}. This proves $T_p f = Tf$. Finally, if $f \in \nul(T_p)$, then $f = z R(z)f$ for all $z$ and if $T$ is injective, then $f=0$ follows from \eqref{eq4: consistent operator}.

The argument for a bisectorial operator is exactly the same, using $z \in \i(0,\infty)$ instead of $z \in (-\infty,0)$ for the limits. In this case we can allow $\mu = \nicefrac{\pi}{2}$ with the convention that $\IC \setminus \cl{\S_{\pi/2}} \coloneqq \i \R$.
\end{proof}

\bibliography{Refs}
\bibliographystyle{abbrv}
\printindex
\end{document}